\documentclass[a4paper]{amsart}

\RequirePackage{amsmath} 
\RequirePackage{amssymb}
\usepackage{amscd,latexsym,amsthm,amsfonts,amssymb,amsmath,amsxtra}
\usepackage[colorlinks=true,urlcolor=blue,citecolor=blue,pdfusetitle]{hyperref}
\usepackage{color}
\usepackage[all]{xy}
\usepackage[OT2,T1]{fontenc}
\usepackage{bm}
\usepackage{mathtools}
\usepackage{mathrsfs}
\usepackage{xcolor}
\usepackage{comment}
\usepackage{marginnote}
\usepackage{enumitem}
\usepackage{thmtools, thm-restate}

\let\Re\undefined
\let\Im\undefined

\DeclareMathOperator{\Re}{Re}
\DeclareMathOperator{\Im}{Im}
\DeclareMathOperator{\Tr}{Tr}

\newcommand{\Res}{\operatorname{Res}}

\newcommand{\sgn}{\operatorname{sgn}}

\newcommand{\Ad}{\operatorname{Ad}}

\newcommand{\reg}{\operatorname{reg}}

\newcommand{\fin}{\operatorname{fin}}
\newcommand{\diag}{\operatorname{diag}}
\newcommand{\Vol}{\operatorname{Vol}}

\newcommand{\Ind}{\operatorname{Ind}}

\newcommand{\RNum}[1]{\uppercase\expandafter{\romannumeral #1\relax}}

\makeatletter
\@ifundefined{subjclassname@2020}{%
  \@namedef{subjclassname@2020}{\textup{2020} Mathematics Subject Classification}%
}{}
\makeatother

\begin{document}
\theoremstyle{plain}
\newtheorem{thm}{Theorem}[section]
	
\newtheorem{cor}[thm]{Corollary}

\newtheorem{thmx}{Theorem}
\renewcommand{\thethmx}{\Alph{thmx}}

\newtheorem{lemma}[thm]{Lemma}
\newtheorem{prop}[thm]{Proposition}
	
\theoremstyle{definition}
\newtheorem{defn}[thm]{Definition}
\theoremstyle{remark}
\newtheorem{remark}[thm]{Remark}	
\numberwithin{equation}{section}
	
\title[]{Spectral Reciprocity: A Fourier--Analytic Approach}%
\makeatletter
\hypersetup{pdftitle={\@title}}
\makeatother
\author{Liyang Yang}
\address{Department of Mathematics, Texas A\&M University\\
College Station, TX 77843, USA}
\email{liyangy@tamu.edu}

\begin{abstract}
We develop a unified Fourier--analytic theory of spectral reciprocity over number fields.
Both Motohashi-type (type~\RNum{1}) and Blomer--Khan-type (type~\RNum{2}) formulas follow from partial Fourier identities for the same automorphic form on $\mathrm{GL}_{3}$.
Type~\RNum{1} compares the $\mathrm{GL}_{2}$ spectral and Mellin expansions of a partial Whittaker coefficient, while type~\RNum{2} uses a Weyl element to relate the Fourier expansions for two embeddings of $\mathrm{GL}_{2}$.
The formulas apply to cuspidal and Eisenstein data and include the degenerate and residual terms.
We construct local test vectors and bound their dual transforms.
The explicit twisted fourth moment and newform sieve give a fifth-moment bound for $\mathrm{PGL}_{2}$ cusp forms spherical at infinity over totally real fields. It covers composite squarefree and arbitrary prime-power levels with polynomial dependence on the Archimedean cutoff.
For type~\RNum{1}, local bounds for the dual transforms of analytic newvector weights give a $\mathrm{GL}_{3}\times\mathrm{GL}_{2}$ first-moment estimate with explicit conductor dependence for cuspidal and maximal Eisenstein data of orthogonal type and trivial central character.
This estimate is used in Rankin--Selberg subconvexity, including for adjoint lifts with a $2+1$ decomposition.
Mixed first-moment estimates with polynomial Archimedean dependence, including for noncuspidal adjoint lifts, are used to prove Weyl subconvexity for $\mathrm{GL}_{2}\times\mathrm{GL}_{2}$ $L$-functions.
Further applications include triple product subconvexity and quantitative simultaneous nonvanishing.
\end{abstract}

\date{November 26, 2026}
\maketitle
\tableofcontents

\section{Introduction}
Motohashi-type (type~\RNum{1}) and Blomer--Khan-type (type~\RNum{2})
reciprocity both follow from partial Fourier identities for the same
$\mathrm{GL}_{3}$ automorphic form. Type~\RNum{1} compares the
$\mathrm{GL}_{2}$ spectral and Mellin expansions of a partial Whittaker
coefficient. Type~\RNum{2} uses a Weyl element to relate Fourier expansions
for two $\mathrm{GL}_{2}$ embeddings.

We outline both derivations for cuspidal data in \textsection\ref{sec1.3}.
The starting point for type~\RNum{1} is Motohashi's formula \cite{Mot93}, which gives
the schematic relation
\begin{equation}\label{1.1} 
\sum_{\phi\ \mathrm{GL}_{2}\text{-form}} L(1/2,\phi)^{3}\, h_{\phi} \ \leftrightsquigarrow\ \int_{\mathbb{R}} |\zeta(1/2+it)|^{4}\, h^{\vee}(t)\, dt, 
\end{equation} 
which connects the cubic moment of $\mathrm{GL}_{2}$ $L$-functions with the
fourth moment of the Riemann zeta function.  Many variants and
generalizations have since appeared; see, for example,
\cite{CI00,Ivi01,Ivi02,Pet15,You17,Fro20,Nel20,Wu22,PK25}.  Two developments
are especially relevant here.
\begin{itemize} 

\item[\textbf{(1)}] \textbf{Motohashi-type formulas.}  
For a cusp form $\varphi$ on $\mathrm{SL}_{3}(\mathbb{Z})$, Kwan
\cite{Kwa24,Kwa25} established refinements of \eqref{1.1}, whose twisted
form can be summarized schematically as
\begin{multline}\label{eq1.1}
\sum_{\text{$\phi$: level 1}} \lambda_{\phi}(n)L(1/2, \varphi \times \phi)h_{\varphi,\phi}\\
\rightsquigarrow \sum_{\chi\bmod n}\int_{\mathbb{R}} L(1/2-it, \widetilde{\varphi}\otimes\chi)\, L(1/2+it,\overline{\chi})\, h_{\varphi,n}^{\vee}(\chi,t)\, dt,
\end{multline}
for arbitrary $n\geq1$, with suitable weights incorporating local factors.
Here we display only the continuous dual moment. The untwisted identity also holds
for minimal Eisenstein data \cite{Kwa23}.
Kwan also indicates the specialization to
maximal Eisenstein data $\tau\boxplus\mathbf{1}$ in
\cite[\textsection 1.6, (1.19)]{Kwa23}. In the cuspidal case,
Wu \cite{Wu25} obtained a distributional extension over general number
fields.

\item[\textbf{(2)}] \textbf{Blomer--Khan reciprocity.}  
For a cusp form or a minimal Eisenstein series $\varphi$ on
$\mathrm{SL}_{3}(\mathbb{Z})$, Blomer and Khan \cite[Theorem~1]{BK19} proved the
reciprocity
\begin{equation}\label{e1.3}
\sum_{\phi\ \text{level } q} \lambda_{\phi}(n)
L(1/2,\varphi \times \phi) L(1/2,\phi)
\rightsquigarrow
\sum_{\phi\ \text{level } n} \lambda_{\phi}(q)
L(1/2,\varphi \times \phi)L(1/2,\phi).
\end{equation}
For an asymmetric version with a primitive central character, see \cite{BK19a}.
For minimal Eisenstein data $\mathbf{1}\boxplus\mathbf{1}\boxplus\mathbf{1}$,
\eqref{e1.3} becomes a twisted fourth-moment formula for $\mathrm{GL}_{2}$ $L$-functions and leads
to fifth-moment estimates and uniform subconvexity.  For cuspidal
$\mathrm{GL}_{3}$ data, Nunes \cite{Nun23} extended \eqref{e1.3} to number
fields.
\end{itemize}

This Fourier perspective subsequently led to the construction in
\cite{YangCanonical26}, which fixes the embedding into the ambient group
and compares different subgroups of $\mathrm{GL}_{2}$. Here we compare
different embeddings of $\mathrm{GL}_{2}$ in $\mathrm{GL}_{3}$.

The principal contributions of this paper are as follows.
\begin{itemize}
\item Our Fourier method gives exact type~\RNum{1} and
type~\RNum{2} reciprocity over number fields for cuspidal and Eisenstein
$\mathrm{GL}_{3}$ representations and makes the direct connection between
the two types explicit at the Fourier level. The formulas include all degenerate and
residual contributions for minimal and maximal Eisenstein data, extending
the cuspidal settings of \cite{Wu25,Nun23}. In particular, both types give
mixed-moment identities for $2+1$ data, including ramification.
\item For type~\RNum{1}, we control the essential support and size of the
dual transforms of our analytic newvector weights at both real and complex
places. The spectral weights are nonnegative and bounded below on the
prescribed conductor ranges, including for nonspherical $\mathrm{GL}_{2}$
components at infinity.
For type~\RNum{2}, we construct spherical weights and obtain explicit
formulas and quantitative estimates for their dual transforms.
\item The type~\RNum{1} first-moment estimate of
Theorem~\ref{thm1.3} has explicit conductor dependence for cuspidal and
maximal Eisenstein $\mathrm{GL}_{3}$ representations of orthogonal type and
trivial central character. It is the first-moment estimate used in the
Rankin--Selberg subconvexity argument of \cite[Theorem~14.4]{HY26}, including adjoint
lifts with a $2+1$ decomposition.
\item Over totally real fields, type~\RNum{2} gives the explicit twisted
fourth-moment identity of Theorem~\ref{thmC}. By isolating newforms of each
conductor and bounding the Archimedean dual transform, we extend the
fifth-moment bound beyond prime level to composite squarefree and arbitrary
prime-power levels, with polynomial dependence on the Archimedean cutoff
(Theorem~\ref{thm1.5}).
\item Further type~\RNum{2} estimates give weighted mixed first-moment
bounds in the level aspect, with polynomial Archimedean dependence,
including for noncuspidal adjoint lifts. See \cite[Remark~1.11]{Yan26Weyl}
for the moment estimates needed in the $\mathrm{GL}_{2}\times\mathrm{GL}_{2}$
Weyl argument.
Other applications include triple product subconvexity and a power lower
bound for simultaneous nonvanishing.
\end{itemize}

\subsection{The Spectral Reciprocity Formulas}
Let $F$ be a number field, and put $G=\mathrm{GL}_{3}$, $G'=\mathrm{GL}_{2}$, and $C_{F}:=F^{\times}\backslash\mathbb{A}_{F}^{(1)}$, with $\rho_{F}:=\Vol(C_{F})$. Fix unitary Hecke characters $\omega$, $\omega'$, and $\eta$ of $F$, and a unitary generic automorphic representation $\pi$ of $G/F$ with central character $\omega$. For $\varphi\in\pi$, write $W_{\varphi}$ for its Whittaker function. For Eisenstein data, the inducing cuspidal representations are assumed unitary.

\subsubsection{Spectral reciprocity of type \RNum{1}}

Let $|\Re(s)|<1/4$. Define
\begin{align*}
J_{\mathrm{cusp}}^{\heartsuit}(s,\varphi,\omega')
   := \sum_{\sigma \in \mathcal{A}_{0}([G'],\omega')}
      \ \sum_{\phi \in \mathfrak{B}(\sigma)}
      \widetilde{\Psi}(1/2+s, W_{\varphi}, W_{\phi})\,
      \overline{W_{\phi}}(I_{2}).
\end{align*}
Here $\widetilde{\Psi}$ is the meromorphically continued Rankin--Selberg period, and $J_{\mathrm{Eis}}^{\heartsuit}$ denotes the corresponding continuous contribution, defined in Proposition~\ref{prop4.1}. On the dual side, put
\begin{align*}
J_{\mathrm{dual}}^{\heartsuit}(s,\varphi,\omega')
   := \frac{1}{2\pi i\rho_{F}}\sum_{\xi \in \widehat{C_{F}}}
      \int_{i\mathbb{R}}
      \widetilde{\mathcal{J}}(s,\lambda;\varphi,\xi,\omega')\, d\lambda.
\end{align*}
The continuation $\widetilde{\mathcal{J}}$ has the completed $L$-factors in
\eqref{e4.10}. The following identity generalizes \eqref{eq1.1}.

\begin{restatable}[]{thmx}{typei} \label{thmA}
Let $|\Re(s)|<1/4$. Then the following identity of meromorphic functions holds:  
\begin{multline}\label{f8.1}
J_{\mathrm{cusp}}^{\heartsuit}(s,\varphi,\omega')
+J_{\mathrm{Eis}}^{\heartsuit}(s,\varphi,\omega')
=J_{\mathrm{sing}}^{\heartsuit}(s,\varphi,\omega')
+J_{\mathrm{dual}}^{\heartsuit}(s,\varphi,\omega')\\
-J_{\mathrm{degen}}^{\heartsuit}(s,\varphi,\omega')
+J_{\mathrm{degen}}^{\dag,\heartsuit}(s,\varphi,\omega')
+R_{\mathrm{dual}}^{+}(s,\varphi,\omega')
-R_{\mathrm{dual}}^{-}(s,\varphi,\omega')-2R_{\RNum{1}}^{+}(s,\varphi,\omega').
\end{multline}
\end{restatable}

See \textsection\ref{sec8.1} for the proof and the terms in \eqref{f8.1}.

\subsubsection{Spectral reciprocity of type \RNum{2}}
Let $\mathbf{s} = (s_{1}, s_{2})$ and $\mathbf{s}^{\vee}:=(\frac{s_{2}-s_{1}}{2},\frac{3s_{1}+s_{2}}{2})$, with $|\Re(s_{1})|<1/2$ and $|\Re(s_{2})|<1/2$. Define 
\begin{align*}
I_{\mathrm{cusp}}^{\heartsuit}(\mathbf{s},\varphi,\omega',\eta)
:=\sum_{\pi'\in \mathcal{A}_{0}([G'],\omega')}
\sum_{\phi\in\mathfrak{B}(\pi')}
\widetilde{\Psi}(1/2+s_{1},W_{\varphi},W_{\phi})
\overline{\widetilde{\Psi}(1/2+\overline{s_{2}},W_{\phi},\overline{\eta})}.
\end{align*}
Let $I_{\mathrm{Eis}}^{\heartsuit}(\mathbf{s},\varphi,\omega',\eta)$ denote the corresponding contribution from the continuous spectrum; see Proposition \ref{prop4.4} for the precise definition. 

The following identity generalizes \eqref{e1.3}.

\begin{restatable}[]{thmx}{typeii} \label{thmii}
Suppose $|\Re(s_{j})|<1/2$ and $|\Re(s_{j}^{\prime})|<1/2$ for $j=1,2$, where $(s_{1}^{\prime},s_{2}^{\prime})=\mathbf{s}^{\vee}$. Then the following identity of meromorphic functions holds:  
\begin{multline}\label{8.1}
I_{\mathrm{cusp}}^{\heartsuit}(\mathbf{s},\varphi,\omega',\eta)
+I_{\mathrm{Eis}}^{\heartsuit}(\mathbf{s},\varphi,\omega',\eta)
=I_{\mathrm{cusp}}^{\heartsuit}(\mathbf{s}^{\vee},\pi(w_{2})\varphi,
\overline{\omega}\overline{\omega}'\eta,\eta)\\
+I_{\mathrm{Eis}}^{\heartsuit}(\mathbf{s}^{\vee},\pi(w_{2})\varphi,
\overline{\omega}\overline{\omega}'\eta,\eta)+I_{\mathrm{degen}}^{\heartsuit}(\mathbf{s}^{\vee},\pi(w_{2})\varphi,\overline{\omega}\overline{\omega}'\eta,\eta)\\
-I_{\mathrm{gen}}^{\heartsuit}(\mathbf{s}^{\vee},\pi(w_{2})\varphi,
\overline{\omega}\overline{\omega}'\eta,\eta)
-I_{\mathrm{degen}}^{\heartsuit}(\mathbf{s},\varphi,\omega',\eta)
+I_{\mathrm{gen}}^{\heartsuit}(\mathbf{s},\varphi,\omega',\eta)\\
+2\sum_{j=1}^{2}
\big[R_{\RNum{2}}^{j,+}(\mathbf{s}^{\vee},\pi(w_{2})\varphi,
\overline{\omega}\overline{\omega}'\eta,\eta)
-R_{\RNum{2}}^{j,+}(\mathbf{s},\varphi,\omega',\eta)\big].
\end{multline}
\end{restatable}

The restrictions on $\mathbf{s}$ and $\mathbf{s}^{\vee}$ give a common
domain for the two spectral expansions in Proposition~\ref{prop4.4}.
See \textsection\ref{sec8.2} for the proof and the terms in \eqref{8.1}.

Proposition~\ref{prop:residues-initial-data} expresses every residual
term in Theorems~\ref{thmA} and~\ref{thmii} as an integral
of the initial automorphic data. The formulas for the degenerate
terms complete the explicit description of both identities.

The type~\RNum{2} formula pairs a $\mathrm{GL}_{3}\times\mathrm{GL}_{2}$
Rankin--Selberg factor with a standard $\mathrm{GL}_{2}$ factor. For
symmetric reciprocity formulas, see \cite{HY26,Yan25b}. Minimal
Eisenstein data give the twisted fourth moment. Maximal Eisenstein data
$\pi=\tau\boxplus\chi$ give a mixed moment of a
$\mathrm{GL}_{2}\times\mathrm{GL}_{2}$ factor and two standard
$\mathrm{GL}_{2}$ factors.

\subsection{Applications of the Spectral Reciprocity Formulas}
\subsubsection{\texorpdfstring{The first moment of $\mathrm{GL}_{3}\times\mathrm{GL}_{2}$ $L$-functions}{The first moment of GL3 x GL2 L-functions}}
Let $\mathcal{F}(\mathfrak{q};\omega')$ denote the set of unitary generic automorphic representations of $G'/F$ whose arithmetic conductor divides $\mathfrak{q}$ and whose central character is $\omega'$.
Use the spectral measure $d\mu_{\sigma}$ in \textsection\ref{def:spectral-measure}.
The type~\RNum{1} formula gives the following estimate.
\begin{thmx}\label{thm1.3}
Let $F$ be a number field and let $\mathfrak{q}$ be an integral ideal.
For each Archimedean place $v\mid\infty$, let $C_{v}\geq10$, and put
$\mathbf{C}_{\infty}=\prod_{v\mid\infty}C_{v}$.
Let $\pi$ be a unitary isobaric automorphic representation of
$\mathrm{GL}_{3}/F$ of orthogonal type with trivial central character,
either cuspidal or of the form $\tau\boxplus\chi$, where $\tau$ is
unitary cuspidal on $\mathrm{GL}_{2}/F$ and $\chi$ is a unitary Hecke character.
Write $\mathfrak{m}$ for its finite conductor ideal and put
$Q=N_{F}(\mathfrak{q}\cap\mathfrak{m})$ and
$A_{v}=C_{v}^{-1}C_{v}(\pi)$. Then, for every
$\varepsilon>0$,
\begin{multline}\label{fc1.1}
\int_{\substack{\sigma\in\mathcal{F}(\mathfrak{q};\mathbf{1})\\
 C_{v}(\sigma)\leq C_{v},\ v\mid\infty}}
 L(1/2,\pi\times\sigma)\,d\mu_{\sigma}
 \ll_{F,\varepsilon}C(\pi)^{\varepsilon}
 \Big\{Q^{1+\varepsilon}\mathbf{C}_{\infty}^{1+\varepsilon}
       \prod_{v\mid\infty}(1+A_{v})\\
 {}+Q^{1/2+\varepsilon}\mathbf{C}_{\infty}^{1/2+\varepsilon}
 C(\pi)^{1/4+\varepsilon}\prod_{v\mid\infty}(1+A_{v}^{2})
 +Q^{\varepsilon}\mathbf{C}_{\infty}^{1/2+\varepsilon}
 C(\pi)^{1/2+\varepsilon}\prod_{v\mid\infty}(1+A_{v})\Big\}.
\end{multline}
If $\pi$ is cuspidal, the product $\prod_{v\mid\infty}(1+A_{v})$
in the first term may be omitted.
\end{thmx}
See \textsection\ref{sec13.2} and \textsection\ref{sec:proof-maximal-first-moment} for the two cases.

\begin{remark}
The Rankin--Selberg subconvexity argument of \cite[Theorem~14.4]{HY26} uses this
estimate with $\pi=\Ad\pi'$.
The cuspidal case applies when $\pi'$ is non-dihedral; the maximal
Eisenstein case covers an adjoint lift with a $2+1$ decomposition.
The additional product in the first term of \eqref{fc1.1} comes from
$R_{\mathrm{dual}}^{-}$; see
Corollary~\ref{cor:maximal-negative-residue-covered-range}. It is
bounded in that application, where $C_{v}(\pi)\ll C_{v}$.
\end{remark}

When $F=\mathbb{Q}$ and $\pi$ is a fixed self-dual Hecke--Maass form on
$\mathrm{SL}_{3}(\mathbb{Z})$, standard Kuznetsov and Petersson arguments are
expected to recover \eqref{fc1.1}. Theorem~\ref{thm1.3} allows
nonspherical Archimedean components of $\sigma$ and gives explicit
dependence on $\pi$. We adapt the analytic newvector construction of
\cite{MV10} to obtain weights that are nonnegative and bounded below in the prescribed conductor ranges. Theorem~\ref{thm9.5}
bounds the size and decay of the dual transforms of these test vectors
at real and complex places. The finite-place bounds and the moment
estimate are in \textsection\ref{sec10} and \textsection\ref{sec13}.
The maximal Eisenstein case gives a mixed first moment of
$L(1/2,\tau\times\sigma)$ and $L(1/2,\sigma\otimes\chi)$,
with explicit conductor dependence as both $\tau$ and $\chi$ vary.
The additional degenerate terms and residues are bounded in
\textsection\ref{sec13.maximal}. For $\pi=\tau\boxplus\chi$, orthogonal type means that
$\tau$ is of orthogonal type and $\chi^{2}=\mathbf{1}$; trivial central
character further gives $\omega_{\tau}=\chi$.
These hypotheses ensure nonnegative central values on the whole
generic $\mathrm{GL}_{2}$ spectrum.

\subsubsection{\texorpdfstring{Explicit twisted fourth moment of $\mathrm{GL}_{2}$ $L$-functions}{Explicit twisted fourth moment of GL2 L-functions}}\label{sec:explicit-fourth-moment-statement}
Let $F$ be a totally real field, and let $\mathfrak{q}, \mathfrak{n} \subseteq \mathcal{O}_{F}$ be ideals such that $\mathfrak{n}+\mathfrak{q}=\mathcal{O}_{F}$. Define 
\begin{align*}
\mathcal{F}^{\circ}(\mathfrak{q};\mathbf{1}):=\big\{\sigma\in \mathcal{F}(\mathfrak{q};\mathbf{1}):\ \text{the Archimedean component of $\sigma$ is unramified}\big\},
\end{align*}
and let $\mathcal{F}_{\mathrm{cusp}}^{\circ}(\mathfrak{q};\mathbf{1})$ be its cuspidal subfamily.

Let $\mathbf{s}=(s_{1},s_{2})$ and $\mathbf{s}^{\vee}=(s_{1}',s_{2}')=(\frac{s_{2}-s_{1}}{2},\frac{3s_{1}+s_{2}}{2})$ with $|\Re(s_{1})|+|\Re(s_{2})|\leq 10^{-1}$. Let $\boldsymbol{\nu}=(\nu_{1},\nu_{2},\nu_{3})\in \mathbb{C}^{3}$ satisfy $|\Re(\nu_{i})|\leq 10^{-2}$, $1\leq i\leq 3$, and $\nu_{1}+\nu_{2}+\nu_{3}=0$. Put $\pi=|\cdot|^{\nu_{1}}\boxplus |\cdot|^{\nu_{2}}\boxplus|\cdot|^{\nu_{3}}$.

\phantomsection\label{def:finite-fourth-moment-factor}
For a unitary generic automorphic representation $\sigma$ of $\mathrm{PGL}_{2}/F$, put
\begin{align*}
\mathcal{L}_{\fin}(\mathbf{s},\pi,\sigma):=L(1/2+s_{1},\pi\times\sigma)
L(1/2+s_{2},\sigma)L(1,\sigma,\Ad)^{-1}.
\end{align*} 
On the continuous spectrum, use the denominator in \eqref{eq:regularized-adjoint}.

For $\mathcal{T}=(\mathcal{T}_{v})_{v\mid\infty}$, with each $\mathcal{T}_{v}=\mathcal{S} h_{v}$ from \eqref{eq10.4}, define
\begin{align*}
&\widetilde{\mathcal{I}}_{\mathcal{T}}
(\mathbf{s},\boldsymbol{\nu};\mathfrak{q},\mathfrak{n})
:=\int_{\sigma\in\mathcal{F}^{\circ}(\mathfrak{q};\mathbf{1})}
\lambda_{\sigma}(\mathfrak{n})\mathcal{L}_{\fin}(\mathbf{s},\pi,\sigma)
\mathcal{H}_{\sigma,\mathcal{T}}
(\mathbf{s},\boldsymbol{\nu};\mathfrak{q},\mathfrak{n})\,d\mu_{\sigma},\\
&\widetilde{\mathcal{I}}_{\mathcal{T}}^{\vee}
(\mathbf{s}^{\vee},\boldsymbol{\nu};\mathfrak{q},\mathfrak{n})
:=\frac{N_{F}(\mathfrak{q})^{1/2-s_{2}'}}
{N_{F}(\mathfrak{n})^{1/2+s_{2}}}
\int_{\sigma\in\mathcal{F}(\mathfrak{n};\mathbf{1})}
\widetilde{\lambda}_{\sigma,\mathbf{s}}(\mathfrak{q})
\mathcal{L}_{\fin}(\mathbf{s}^{\vee},\pi,\sigma)
\mathcal{H}_{\sigma,\mathcal{T}}^{\vee}
(\mathbf{s}^{\vee},\boldsymbol{\nu};\mathfrak{q},\mathfrak{n})\,d\mu_{\sigma}.
\end{align*}
Here $\lambda_{\sigma}(\mathfrak{n})$ and $\widetilde{\lambda}_{\sigma,\mathbf{s}}(\mathfrak{q})$ are defined in \eqref{eq12.3} and \eqref{eq14.123}.
The spectral weights are defined in \textsection\ref{def:fourth-spectral-weights}.

For this $\pi$ and $\omega'=\eta=\mathbf{1}$, Theorem~\ref{thmii} gives the following twisted fourth-moment formula.
\begin{restatable}[]{thmx}{typeC} \label{thmC}
For the minimal Eisenstein data and test functions of
\textsection\ref{sec14},
\begin{multline}\label{fc1.3}
\widetilde{\mathcal{I}}_{\mathcal{T}}(\mathbf{s},\boldsymbol{\nu};\mathfrak{q},\mathfrak{n})
=\widetilde{\mathcal{I}}_{\mathrm{gen},\mathcal{T}}(\mathbf{s},\boldsymbol{\nu};\mathfrak{q},\mathfrak{n})+\widetilde{\mathcal{I}}_{\mathcal{T}}^{\vee}(\mathbf{s}^{\vee},\boldsymbol{\nu};\mathfrak{q},\mathfrak{n})\\
+\widetilde{\mathcal{I}}_{\mathrm{degen},\mathcal{T}}(\mathbf{s},\boldsymbol{\nu};\mathfrak{q},\mathfrak{n})
+\widetilde{\mathcal{I}}_{\mathrm{res},\mathcal{T}}(\mathbf{s},\boldsymbol{\nu};\mathfrak{q},\mathfrak{n}).
\end{multline}
The identity holds meromorphically in the parameters. At $\nu_{i}=\nu_{j}$
or another pole of an individual term, combine the full right-hand side
before continuation. The generic, degenerate, and residual terms and
their local factors are defined in \textsection\ref{sec14.8}.
\end{restatable}
See \textsection\ref{sec:proof-explicit-fourth} for the proof.

\begin{itemize}
\item The \emph{generic term} is defined by
\begin{multline*}
\widetilde{\mathcal{I}}_{\mathrm{gen},\mathcal{T}}(\mathbf{s},\boldsymbol{\nu};\mathfrak{q},\mathfrak{n}):=2D_{F}^{1/2}\big[\mathcal{L}_{\fin}(\mathbf{s},\pi)
N_{F}(\mathfrak{q})N_{F}(\mathfrak{n})^{-1/2-s_{2}}
\mathcal{H}_{\mathrm{gen},\mathcal{T}}(\mathbf{s},\boldsymbol{\nu};\mathfrak{q},\mathfrak{n})\\
-\mathcal{L}_{\fin}(\mathbf{s}^{\vee},\pi)
N_{F}(\mathfrak{q})^{-s_{1}-s_{2}}N_{F}(\mathfrak{n})^{1/2+s_{2}}
\mathcal{H}_{\mathrm{gen},\mathcal{T}}^{\vee}(\mathbf{s}^{\vee},\boldsymbol{\nu};\mathfrak{q},\mathfrak{n})\big],
\end{multline*}
where $\mathcal{L}_{\fin}(\mathbf{s},\pi)
=L(1+s_{1}+s_{2},\pi)L(1+2s_{1},\widetilde{\pi})/\zeta_{F}(2+3s_{1}+s_{2})$.
The generic weights are defined in \textsection\ref{def:fourth-generic-weights}.

\item The \emph{degenerate term} is defined by
\begin{multline*}
\widetilde{\mathcal{I}}_{\mathrm{degen},\mathcal{T}}(\mathbf{s},\boldsymbol{\nu};\mathfrak{q},\mathfrak{n}):=2D_{F}^{1/2}\sum_{w\in \mathbb{W}}\Big[\mathcal{L}_{\fin}(\mathbf{s}^{\vee},\boldsymbol{\nu};w)\mathcal{H}_{\mathrm{degen},\mathcal{T}}^{\vee}(\mathbf{s}^{\vee},\boldsymbol{\nu};\mathfrak{q},\mathfrak{n};w)\\
-\mathcal{L}_{\fin}(\mathbf{s},\boldsymbol{\nu};w)\mathcal{H}_{\mathrm{degen},\mathcal{T}}(\mathbf{s},\boldsymbol{\nu};\mathfrak{q},\mathfrak{n};w)\Big],
\end{multline*}
where $\mathbb{W}=\{w_{2}, w_{1}w_{2}, w_{1}w_{2}w_{1}\}$. The factors
$\mathcal{L}_{\fin}(\cdot,\boldsymbol{\nu};w)$ and the weights are defined
in \textsection\ref{sec14.8.2}.

\item The \emph{residual term} is defined by
\begin{align*}
\widetilde{\mathcal{I}}_{\mathrm{res},\mathcal{T}}(\mathbf{s},\boldsymbol{\nu};\mathfrak{q},\mathfrak{n}):=\mathscr{R}_{\mathbf{s}^{\vee},\boldsymbol{\nu}}
[D_{\mathcal{T}}^{\vee}(\mathbf{s}^{\vee},\cdot)]
-\mathscr{R}_{\mathbf{s},\boldsymbol{\nu}}
[D_{\mathcal{T}}(\mathbf{s},\cdot)],
\end{align*}
with densities $D_{\mathcal{T}},D_{\mathcal{T}}^{\vee}$ in \textsection\ref{def:fourth-residue-densities}
and residue sum $\mathscr{R}$ in \eqref{eq:residue-operator}.

\end{itemize}

The test function $\mathcal{T}_{v}$ determines the original local weight in
\textsection\ref{def:fourth-spectral-weights}. Right translation by $w_{2}$
gives the dual weight by \eqref{fc10.2}.
Thus the formula exchanges the level and Hecke twist with explicit
Archimedean weights.

\subsubsection{Fifth moments at squarefree and prime-power levels}

For the remaining applications, $F$ is totally real and
$0\leq\vartheta<1/2$ denotes an admissible exponent toward the Ramanujan
conjecture at the finite places for $\mathrm{PGL}_{2}/F$.

In \textsection\ref{sec15.2} we subtract the oldforms from the twisted
fourth moment before applying the approximate functional equation to
the fifth factor at each conductor.
\begin{thmx}\label{thm1.5}
Let $\mathfrak{q}$ be an integral ideal that is squarefree or a prime
power. For $T>10$ and $0<\varepsilon<10^{-3}$, we have
\begin{align*}
 \sum_{\substack{\sigma\in\mathcal{F}_{\mathrm{cusp}}^{\circ}
 (\mathfrak{q};\mathbf{1})\\ C_{\infty}(\sigma)\leq T}}
 L(1/2,\sigma)^{5}
 \ll_{F,\varepsilon}N_{F}(\mathfrak{q})^{1+\vartheta+\varepsilon}
 T^{5/2+\varepsilon}.
\end{align*}
\end{thmx}
For fixed $T$, finite-place Ramanujan gives level exponent $1+\varepsilon$.
No temperedness at infinity is needed. In \eqref{eq:fifth-moment-dyadic-bound},
the dual fourth moment and the length of the approximate functional
equation account for the exponent $5/2$.

In the nonvanishing argument of \textsection\ref{sec:simultaneous-nonvanishing},
the weight decays faster than any fixed power of $T$, so polynomial
dependence on $T$ does not change the level exponent. The level-aspect $\mathrm{GL}_{2}\times\mathrm{GL}_{2}$ Weyl argument
outlined in \cite[Remark~1.11]{Yan26Weyl} over $\mathbb{Q}$ with trivial
central character needs only polynomial dependence on the spectral cutoff
in the fifth-moment estimate.
In Nelson's level-aspect Weyl argument \cite[(4.5)--(4.6)]{Nel20},
the rapid decay of the dual weight likewise makes polynomial dependence
on the Archimedean cutoffs of the fourth moment of Hecke characters sufficient.

\begin{remark}
For squarefree $\mathfrak{q}$, the estimate is uniform in the number of
prime divisors, and for $\mathfrak{q}=\mathfrak{p}^{r}$ it is uniform in
$\mathfrak{p}$ and $r$. Over $\mathbb{Q}$ at prime level and fixed $T$, this
is the level bound of \cite[Theorem~3]{BK19}. The composite squarefree
case and prime powers with $r\geq2$ appear to be new even over $\mathbb{Q}$.

K{\i}ral and Young \cite{KY21} work at prime level in a small range of
holomorphic weights where the relevant cusp forms are automatically new;
they point to a newform Petersson formula as a way to relax this
restriction. Our newform sieve separates representations by conductor before
applying the functional equation and requires no such oldform restriction.
\end{remark}

\begin{remark}
Lemma~\ref{lem:fixed-level-fifth-moment} gives the fifth-moment bound
at level $\mathfrak{f}\mathfrak{q}$ for fixed $\mathfrak{f}$ and
squarefree $\mathfrak{q}$ coprime to $\mathfrak{f}$, with an implied
constant depending on $\mathfrak{f}$.
\end{remark}

\begin{cor}
Let $\sigma\in\mathcal{F}_{\mathrm{cusp}}^{\circ}(\mathfrak{q};\mathbf{1})$
have conductor $\mathfrak{q}$, where $\mathfrak{q}$ is squarefree or a
prime power. Then
\begin{align*}
L(1/2,\sigma)\ll
C_{\fin}(\sigma)^{\frac{1+\vartheta}{5}+\varepsilon},
\end{align*}
with an implied constant depending on its Archimedean component.
\end{cor}
This follows from Theorem~\ref{thm1.5} and nonnegativity.

\subsubsection{\texorpdfstring{Mixed moments of $\mathrm{GL}_{3}\times\mathrm{GL}_{2}$ and $\mathrm{GL}_{2}$ $L$-functions}{Mixed moments of GL3 x GL2 and GL2 L-functions}}
The cuspidal type~\RNum{2} formula also gives mixed first moments.
Let $\pi$ be a fixed unitary self-dual cuspidal representation of
$\mathrm{PGL}_{3}/F$, unramified at all places. For the local weights $\mathcal{V}_{v}^{\circ}$ of
Lemma~\ref{lem:mixed-spherical-arch},
Lemma~\ref{lem:mixed-full-spectrum-mass} gives, for every integral ideal
$\mathfrak{a}$,
\begin{align*}
 \int_{\sigma\in\mathcal{F}(\mathfrak{a};\mathbf{1})}
 \frac{L(1/2,\pi\times\sigma)L(1/2,\sigma)}{L(1,\sigma,\Ad)}
 \prod_{v\mid\infty}\mathcal{V}_{v}^{\circ}(\sigma_{v})\,d\mu_{\sigma}
 \ll_{F,\pi,\varepsilon}N_{F}(\mathfrak{a})^{1+\varepsilon}.
\end{align*}
This includes the continuous spectrum. Lemma~\ref{lem:mixed-full-spectrum-mass}
also covers $\pi=\Ad\pi_{0}$, with $\pi_{0}$ fixed cuspidal, possibly
dihedral, and $\Ad\pi_{0}$ unramified at all places; it gives polynomial
dependence on a growing Archimedean cutoff. The fixed-weight estimate
supplies a moment input for the $\mathrm{GL}_{2}\times\mathrm{GL}_{2}$
Weyl argument outlined in \cite[Remark~1.11]{Yan26Weyl}.

An additional approximate functional equation gives the following moment
with a second standard $L$-factor and a growing Archimedean cutoff.
For an unramified finite component $\pi_{v}$, let
$\alpha_{1,v},\alpha_{2,v},\alpha_{3,v}$ be its Satake parameters and put
\begin{align*}
 \alpha_{\pi}(\mathfrak{q}):=
 \prod_{v\mid\mathfrak{q}}
 \prod_{j=1}^{3}(1+|\alpha_{j,v}|)(1+|\alpha_{j,v}|^{-1}).
\end{align*}
This factor bounds the finite weights in the newform sieve. It depends
only on the prime divisors of $\mathfrak{q}$, and not on their exponents.

\begin{cor}\label{thm1.7}
Fix a unitary self-dual cuspidal representation $\pi$ of
$\mathrm{PGL}_{3}/F$ unramified at all places. Let $T>10$ and
$0<\varepsilon<10^{-3}$. Suppose that $\mathfrak{q}$ is squarefree, or
that $\mathfrak{q}=\mathfrak{p}^{r}$ with $\pi_{\mathfrak{p}}$ tempered.
Then
\begin{equation}\label{eq1.3}
\sum_{\substack{\sigma\in \mathcal{F}_{\mathrm{cusp}}^{\circ}(\mathfrak{q};\mathbf{1})\\
C_{\infty}(\sigma)\leq T}}
L(1/2,\pi\times \sigma)L(1/2,\sigma)^{2}
\ll_{F,\pi,\varepsilon}\alpha_{\pi}(\mathfrak{q})
N_{F}(\mathfrak{q})^{1+\vartheta+\varepsilon}T^{O_{\pi,F}(1)}.
\end{equation}
\end{cor}
See \textsection\ref{sec:proof-mixed-moment} for the proof.

\begin{remark}
If $\pi_{v}$ is tempered at every finite place $v\mid\mathfrak{q}$, then
$\alpha_{\pi}(\mathfrak{q})\leq64^{\omega(\mathfrak{q})}\ll
N_{F}(\mathfrak{q})^{\varepsilon}$. For fixed $T$,
\eqref{eq1.3} then gives $O(N_{F}(\mathfrak{q})^{1+\vartheta+\varepsilon})$,
which is Lindel\"of on average in the level aspect under Ramanujan. The exponent of $T$ is not optimized.
\end{remark}

\subsubsection{\texorpdfstring{Subconvexity for triple product $L$-functions}{Subconvexity for triple product L-functions}}
Amplifying the mixed first moment with $\pi=\Ad\pi_{1}$ gives the following triple product bound.

\begin{thm}\label{thm1.6}
Let $\pi_{1}$ be a fixed unitary cuspidal representation of
$\mathrm{GL}_{2}/F$, unramified at all places, with $\Ad\pi_{1}$ cuspidal.
Let $\pi_{2}\in\mathcal{F}_{\mathrm{cusp}}^{\circ}(\mathfrak{q};\mathbf{1})$
have conductor $\mathfrak{q}\neq\mathcal{O}_{F}$, where $\mathfrak{q}$ is
squarefree or a prime power. Suppose that $\pi_{1,v}$ is tempered for
every $v\mid\mathfrak{q}$. Then
\begin{align*}
L(1/2,\pi_{1}\times\widetilde{\pi}_{1}\times\pi_{2})
\ll C_{\fin}(\pi_{2})^{1-\frac{1-2\vartheta}{6}+\varepsilon},
\end{align*}
where the implied constant depends on $\pi_{1}$, $F$, $\varepsilon$, and
the Archimedean components of $\pi_{2}$.
\end{thm}
See \textsection\ref{sec:proof-triple-product} for the proof.

\begin{remark}
By \cite{Yan23a}, the representation $\pi$ in Corollary~\ref{thm1.7}
has infinitely many unramified tempered
places. The same holds for $\Ad\pi_{1}$ when it is cuspidal, and
temperedness of $(\Ad\pi_{1})_{\mathfrak{p}}$ is equivalent to that of
$\pi_{1,\mathfrak{p}}$ at those places. Thus the prime-power estimates
in Corollary~\ref{thm1.7} and Theorem~\ref{thm1.6} hold for infinitely
many $\mathfrak{p}$, uniformly for $r\geq1$.
\end{remark}

\subsubsection{Simultaneous nonvanishing}

The cuspidal type~\RNum{2} formula also gives simultaneous nonvanishing
for its two central factors.

\begin{thm}\label{thm1.9}
Let $\pi$ be as in Corollary~\ref{thm1.7}, and assume $\vartheta<1/4$.
There is a smooth weight $0\leq \alpha_{\infty}\leq1$ on the spherical
Archimedean unitary dual, depending only on $F$ and $\pi$, with
$\alpha_{\infty}(\sigma_{\infty})\ll_{\ell} C_{\infty}(\sigma)^{-\ell}$ for every $\ell>0$,
such that, for $0<\varepsilon<10^{-3}$ and squarefree $\mathfrak{q}$
of sufficiently large norm,
\begin{equation}\label{1.5}
 \sum_{\substack{\sigma\in
 \mathcal{F}_{\mathrm{cusp}}^{\circ}(\mathfrak{q};\mathbf{1})\\
 L(1/2,\pi\times\sigma)L(1/2,\sigma)\neq0}}
 \alpha_{\infty}(\sigma_{\infty})
 \gg_{F,\pi,\varepsilon}
 N_{F}(\mathfrak{q})^{\frac{1-4\vartheta}{20}-\varepsilon}.
\end{equation}
\end{thm}
See \textsection\ref{sec:proof-simultaneous-nonvanishing} for the proof.

Since $\alpha_{\infty}\leq1$, \eqref{1.5} gives the same lower
bound for the number of distinct representations in the displayed
family for which both central values are nonzero.
For prime conductor, the first-moment asymptotic of
\cite[Corollary~1.3]{Nun23} and $\mathrm{GL}_{2}$ subconvexity \cite{MV10}
already imply a power lower bound. Theorem~\ref{thm1.9} gives the explicit
exponent in \eqref{1.5} for conductors dividing any squarefree $\mathfrak{q}$.

For $\mathrm{GL}_{4}\times\mathrm{GL}_{1}$, Radziwi{\l}{\l} and
Yang \cite{RY23} proved nonvanishing for infinitely many primitive
Dirichlet twists of a fixed cuspidal $\mathrm{GL}_{4}$ representation. For
$\mathrm{GL}_{4}\times\mathrm{GL}_{2}$, Blomer, Li, and Miller
\cite[Theorem~4]{BLM19} proved nonvanishing for infinitely many odd
Maass cusp forms for $\mathrm{SL}_{2}(\mathbb{Z})$ against a fixed
self-dual cuspidal $\mathrm{GL}_{4}$ representation over $\mathbb{Q}$
unramified at all finite places. Their quantitative lower bound is
$c_{0}\log\log T$ nonvanishing forms with spectral parameter at most $T$,
for every fixed $c_{0}>1$ and all sufficiently large $T$. Equation~\eqref{1.5}
gives a power lower bound in the level aspect for the noncuspidal isobaric $3+1$ representation
$\pi\boxplus\mathbf{1}$.

\subsection{Reciprocity via Partial Fourier Expansion}\label{sec1.3}

Write $n_{ij}(u)=I_{3}+uE_{ij}$ for $i\ne j$, where $E_{ij}$ is the
$3\times3$ matrix unit in row $i$, column $j$.
Assume throughout this subsection that $\pi$ is cuspidal. For $\varphi\in\pi$,
the Fourier--Whittaker expansion is
\begin{align*}
\varphi(I_{3})
   =\sum_{\delta \in N'(F)\backslash G'(F)}
      W_{\varphi}\left(
      \diag(\delta,1)
\right),
\end{align*}
where $W_{\varphi}$ denotes the Whittaker function associated to $\varphi$; see \eqref{w2.5}.

\subsubsection{Type \RNum{1}: partial Fourier expansion}

Automorphy and the Bruhat decomposition reorganize the Fourier--Whittaker
expansion into the partial expansion

\begin{align*}
\int_{F\backslash \mathbb{A}_{F}}
\varphi(n_{12}(a))
\overline{\psi(a)}\, da
=
\sum_{\alpha \in F^{\times}}
\sum_{\beta \in F}
W_{\varphi}\bigl(n_{32}(\beta)\diag(\alpha,\alpha,1)\bigr),
\end{align*}
where $\psi(\cdot)$ is an additive character of $F\backslash\mathbb{A}_{F}$; see Proposition \ref{prop2.2}. 

We apply this expansion to
\begin{align*}
\int_{F^{\times}\backslash \mathbb{A}_{F}^{\times}}
\omega'(z)\, |z|^{2s}\,
\pi(\mathrm{diag}(z,z,1))\,\varphi\, d^{\times}z.
\end{align*}

We insert the $\mathrm{GL}_{2}$-spectral decomposition on the left.  On the
right, we apply the $\mathrm{GL}_{1}$-spectral decomposition to the
$\beta\in F^{\times}$ part of the Whittaker expression.  Meromorphic
continuation gives Theorem~\ref{thmA}.

\subsubsection{Type \RNum{2}: partial additive Poisson summation}

Let $\delta \in G(F)$.  Automorphy and the Fourier expansion give
\begin{align*}
\sum_{\gamma \in F}
\int_{F\backslash \mathbb{A}_{F}}
\varphi(n_{12}(a))\overline{\psi(\gamma a)}\, da
=
\sum_{\gamma \in F}
\int_{F\backslash \mathbb{A}_{F}}
\varphi(n_{12}(a)\delta)\overline{\psi(\gamma a)}\, da.
\end{align*}

For $g'\in G'(\mathbb{A}_{F})$, automorphy gives
\begin{align*}
\varphi(\diag(g',1)w_{2})
=\varphi(w_{2}^{-1}\diag(g',1)w_{2}).
\end{align*}
Take $\delta=w_{2}$ in the Fourier expansion in $a$ and use the
Bruhat decomposition to obtain
\begin{multline}\label{1.16}
\sum_{\gamma \in F^{\times}}
\int_{F\backslash \mathbb{A}_{F}}
\varphi\bigl(n_{12}(a)\diag(\gamma,1,1)\bigr)\overline{\psi(a)}\, da
- \mathcal{I}_{\mathrm{gen}}(\varphi)
\\
=
\sum_{\gamma \in F^{\times}}
\int_{F\backslash \mathbb{A}_{F}}
\varphi\bigl(n_{12}(a)\diag(\gamma,1,1)w_{2}\bigr)
\overline{\psi(a)}\, da
- \mathcal{I}_{\mathrm{gen}}\bigl(\pi(w_{2})\varphi\bigr),
\end{multline}
where
\begin{align*}
\mathcal{I}_{\mathrm{gen}}(\varphi)
:= \sum_{\alpha \in F^{\times}}
   \sum_{\gamma \in F^{\times}}
      W_{\varphi}\left(
      \diag(\alpha \gamma,\gamma,1)
      \right).
\end{align*}

Suppose $\Re(s_{2})>\Re(s_{1})+5\geq 10$.  We apply \eqref{1.16} to the integral
\begin{align*}
\int_{F^{\times}\backslash \mathbb{A}_{F}^{\times}}
\int_{F^{\times}\backslash \mathbb{A}_{F}^{\times}}
\omega'(z)\, |z|^{2s_{1}}\,
\eta(y)\, |y|^{s_{1}+s_{2}}\,
\pi(\mathrm{diag}(yz, z, 1))\varphi\,
d^{\times}y\, d^{\times}z.
\end{align*}

We insert the $\mathrm{GL}_{2}$-spectral decomposition on both sides.
Meromorphic continuation then gives Theorem~\ref{thmii}.

Combining the Fourier expansions with multiplicative summation also
gives a type~\RNum{3} identity, relating two $\mathrm{GL}_{2}$-spectral
expansions with constant and degenerate contributions.
Appendix~\ref{app:further-reciprocity} gives the type~\RNum{3}
spectral formula in its region of absolute convergence and discusses
higher ranks.

\subsection{Structure of the Paper}
Parts I--III prove the reciprocity formulas, construct the local
weights and their dual transforms, and give the arithmetic applications,
respectively. Appendices~\ref{app:fourth-local}--\ref{app:further-reciprocity}
contain the local fourth-moment calculations, the polynomial
Archimedean bounds, and the further reciprocity formulas.
The List of Notation and Transforms before the references collects the
parameter changes, local weights, and normalizations used across these sections.

\subsection{Notation}
\subsubsection{Number fields and measures}
Let $F$ be a number field with ring of integers $\mathcal{O}_{F}$,
absolute norm $N_{F}$, different $\mathfrak{D}_{F}$, absolute discriminant
$D_{F}=N_{F}(\mathfrak{D}_{F})$, and adele ring $\mathbb{A}_{F}$.  Write $\Sigma_{F}$ for its set of places and
$F_{v}$ for the completion at $v$.  If $v<\infty$, let $\mathcal{O}_{v}$ be
the ring of integers, $\mathfrak{p}_{v}$ its maximal ideal, and $\varpi_{v}$ a
fixed uniformizer.  We normalize the valuation by
$e_{v}(\varpi_{v})=1$ and put
\begin{align*}
 q_{v}=|\mathcal{O}_{v}/\mathfrak{p}_{v}|,
 \qquad d_{v}=e_{v}(\mathfrak{D}_{F}).
\end{align*}
For an integral ideal
$\mathfrak{n}$, the notation $v\mid\mathfrak{n}$ means
$\mathfrak{n}\mathcal{O}_{v}\subseteq\mathfrak{p}_{v}$.

We use $v\mid\infty$ for an Archimedean place and $v<\infty$ for a
non-Archimedean place.  Let $|\cdot|_{v}$ be the normalized absolute value
on $F_{v}$; at a complex place, $|z|_{v}=z\overline{z}$.
On $\mathbb{A}_{F}^{\times}$ and its quotients we write
$|x|_{\mathbb{A}_{F}}=\prod_{v}|x_{v}|_{v}$, usually abbreviated to $|x|$.  We likewise abbreviate products
and tensor products over $v\in\Sigma_{F}$ to $\prod_{v}$ and $\otimes_{v}$.

Let $\psi=\prod_{v}\psi_{v}$ be the additive character $\psi=\psi_{\mathbb{Q}}\circ \Tr_{F}$ of $F\backslash\mathbb{A}_{F}$, where $\Tr_{F}$ is the trace map, and $\psi_{\mathbb{Q}}$ is the additive character of $\mathbb{Q}\backslash\mathbb{A}_{\mathbb{Q}}$ taking $x\mapsto e^{2\pi i x}$ on $\mathbb{R}$. At $v<\infty$, $\psi_{v}$ has conductor $\mathfrak{p}_{v}^{-d_{v}}$. Globally, we have the orthogonality relation 
\begin{equation}\label{e1.1}
\int_{F\backslash\mathbb{A}_{F}}\psi(c)dc=0.	
\end{equation}

For $v\in\Sigma_{F}$, let $dt_{v}$ be the additive Haar measure on $F_{v}$
self-dual relative to $\psi_{v}$.  Then $dt=\prod_{v}dt_{v}$ is the standard
Tamagawa measure on $\mathbb{A}_{F}$.  Put
$d^{\times} t_{v}=\zeta_{v}(1)dt_{v}/|t_{v}|_{v}$, where $\zeta_{v}(s)$ is
$(1-q_{v}^{-s})^{-1}$ at a finite place, $\Gamma_{\mathbb{R}}(s):=\pi^{-s/2}\Gamma(s/2)$
at a real place, and $\Gamma_{\mathbb{C}}(s):=2(2\pi)^{-s}\Gamma(s)$ at a complex place.
At $v<\infty$ we then have
\begin{align*}
 \Vol(\mathcal{O}_{v}^{\times},d^{\times} t_{v})
 =\Vol(\mathcal{O}_{v},dt_{v})=q_{v}^{-d_{v}/2}.
\end{align*}
Moreover,
\begin{align*}
 \Vol(F\backslash\mathbb{A}_{F};dt)=1,
 \qquad
 \rho_{F}:=\Vol(C_{F};d^{\times} t)
 =\Res_{s=1}\zeta_{F}(s)\prod_{v\mid\infty}\zeta_{v}(1),
\end{align*}
where the measure on the subgroup $\mathbb{A}_{F}^{(1)}$ of norm-one
ideles is induced using $dt/t$ in the norm direction, and
$\zeta_{F}(s)=\prod_{v<\infty}\zeta_{v}(s)$ is the finite Dedekind zeta
function. Thus this volume equals $\Res_{s=1}\zeta_{F}(s)$ when
$F$ is totally real. Denote by
$\widehat{C_{F}}$ the Pontryagin dual of
$C_{F}$. Fix a splitting of the idele
norm and extend these characters trivially on its positive real factor.
	
\subsubsection{\texorpdfstring{The Haar measure on $G(\mathbb{A}_{F})$}{The Haar measure on G over the adeles}}\label{sec:matrix-conventions}

Throughout the paper we set $G=\mathrm{GL}_{3}$ and $G'=\mathrm{GL}_{2}$.  Let
$B$ and $B'$ be their upper triangular subgroups, with unipotent
radicals $N$ and $N'$. We view $G'$ in $G$ by
$g'\mapsto\diag(g',1)$. Write $Z$ and $Z'$ for their centers, and
$\overline{G'}=Z'\backslash G'$.
Let $K=\prod_{v}K_{v}$ be a maximal compact subgroup of $G(\mathbb{A}_{F})$,
where $K_{v}=\mathrm{U}_{3}(\mathbb{C})$ if $v$ is complex,
$K_{v}=\mathrm{O}_{3}(\mathbb{R})$ if $v$ is real, and
$K_{v}=G(\mathcal{O}_{v})$ if $v<\infty$. Set
$K_{v}'=K_{v}\cap G'(F_{v})$ and $K'=\prod_{v}K_{v}'$. For
$v<\infty$ and $m\in\mathbb{Z}_{\geq0}$, define
\begin{align*}
K_{0,v}'[m]:=\Big\{\begin{pmatrix}
a&b\\
c&d
\end{pmatrix}\in G'(\mathcal{O}_{v}):\ c\in \varpi_{v}^{m}\mathcal{O}_{v}\Big\}.
\end{align*}

For a matrix $x$, $x_{ij}$ denotes its entry in row $i$ and column $j$.
Blank entries in matrices denote zero.
We write $a(y):=\diag(y,1)$ and $n(b):=\left(\begin{smallmatrix}1&b\\&1\end{smallmatrix}\right)$.
Define the Weyl elements
\begin{align*}
w'=\begin{pmatrix}
& 1\\
1& 
\end{pmatrix}, \ \
w_{1}=\begin{pmatrix}
& 1\\
1& \\
&& 1
\end{pmatrix} \ \ \text{and}\ \ w_{2}=\begin{pmatrix}
1& \\
&& 1 \\
&1& 
\end{pmatrix}.
\end{align*}
Write $\widetilde{w}:=w_{1}w_{2}w_{1}$ for the long Weyl element, and $g^{\iota}$ for the inverse transpose of $g$.

Unipotent measures are products of the self-dual additive measures
in the standard matrix coordinates. For $v\in\Sigma_{F}$, use the
Iwasawa coordinates
\begin{align*}
g_{v}=z_{v}\diag(a_{1,v},a_{2,v},1)u_{v}k_{v},
\qquad z_{v},a_{1,v},a_{2,v}\in F_{v}^{\times},
\quad u_{v}\in N(F_{v}),\quad k_{v}\in K_{v}.
\end{align*}
With $dk_{v}$ the probability Haar measure on $K_{v}$, we use
\begin{align*}
dg_{v}=d^{\times} z_{v}\,d^{\times} a_{1,v}\,d^{\times} a_{2,v}\,du_{v}\,dk_{v}.
\end{align*}

For $G'(F_{v})$, let $dk_{v}'$ be probability Haar measure on $K_{v}'$.
In the coordinate
$g_{v}'=z_{v}'\diag(a_{v}',1)
n(b_{v}')k_{v}'$,
use $dg_{v}'=d^{\times} z_{v}'\,d^{\times} a_{v}'\,db_{v}'\,dk_{v}'$.
The equivalent coordinates with $K_{v}'$ on the left give
\begin{equation}\label{e1.4}
 g_{v}'=z_{v}'k_{v}'\diag(a_{v}',1)
 n(b_{v}'),\qquad
 dg_{v}'=|a_{v}'|_{v}\,d^{\times} z_{v}'\,dk_{v}'\,d^{\times} a_{v}'\,db_{v}',
\end{equation}

\begin{equation}\label{e1.5}
 g_{v}'=z_{v}'k_{v}'\diag(a_{v}',1)
 \begin{pmatrix}1&\\b_{v}'&1\end{pmatrix},\qquad
 dg_{v}'=|a_{v}'|_{v}^{-1}\,d^{\times} z_{v}'\,dk_{v}'\,d^{\times} a_{v}'\,db_{v}'.
\end{equation}

\subsubsection{Automorphic data}\label{sec:automorphic-data-conventions}
 
Let $H$ be an algebraic group and $[H] := H(F) \backslash H(\mathbb{A}_{F})$ be the automorphic quotient.  If $v<\infty$ and $\pi_{v}$ is an
irreducible generic representation of $G'(F_{v})$ with central character
$\omega_{\pi_{v}}$, let $r_{\pi_{v}}$ denote its conductor exponent.  It is $0$ precisely when $\pi_{v}$ is spherical. Otherwise,
it is the least $r\geq1$ for which $\pi_{v}$ contains a nonzero vector
$\phi$ satisfying
\begin{align*}
 \phi(gk)=\omega_{\pi_{v}}(k_{22})\phi(g),\qquad k\in K_{0,v}'[r].
\end{align*}
For representations
on the other general linear groups, $r_{\pi_{v}}$ has its standard
conductor meaning.

For a global representation or Hecke character $\rho$, let
$\mathfrak{c}(\rho):=\prod_{v<\infty}\mathfrak{p}_{v}^{r_{\rho_{v}}}$
be its finite conductor ideal. Write
$C_{v}(\rho)=C(\rho_{v})$ for its local analytic conductor, with
$C_{v}(\rho)=q_{v}^{r_{\rho_{v}}}$ at finite places. Put
\begin{align*}
 C_{\fin}(\rho)=\prod_{v<\infty}C_{v}(\rho),\qquad
 C_{\infty}(\rho)=\prod_{v\mid\infty}C_{v}(\rho),\qquad
 C(\rho)=C_{\fin}(\rho)C_{\infty}(\rho).
\end{align*}
For representations or characters defined only locally, we use $C(\rho_{v})$.
We write $L(s,\rho)=\prod_{v<\infty}L_{v}(s,\rho_{v})$ and $\Lambda(s,\rho)=\prod_{v}L_{v}(s,\rho_{v})$, also for Rankin--Selberg $L$-functions.
For a nonzero integral ideal $\mathfrak{a}$, a superscript
$(\mathfrak{a})$ means that the Euler factors at $v\mid\mathfrak{a}$ are omitted.
For a standard $L$-function, $\lambda_{\rho}(\mathfrak{n})$ denotes its
Dirichlet coefficient.

Let $\omega$ be a unitary character of
$F^{\times}\backslash\mathbb{A}_{F}^{\times}$.  Denote by
$\mathcal{A}_{0}([H],\omega)$ the set of unitary cuspidal automorphic
representations of $H(\mathbb{A}_{F})$ with central character $\omega$,
for $H=G$ or $G'$.

For unitary Hecke characters $\xi_{1},\xi_{2}$,
let $\Ind(\xi_{1} \otimes \xi_{2})$ denote the normalized parabolic induction
from $B'(\mathbb{A}_{F})$ to $G'(\mathbb{A}_{F})$ associated with
$\xi_{1} \otimes \xi_{2}$.  We write
$\xi_{1} \boxplus \xi_{2} = \Ind(\xi_{1} \otimes \xi_{2})$.

\subsubsection{Other conventions} 
All circular contours are positively oriented.

Throughout, we follow the $\varepsilon$-convention: $\varepsilon$ is
always positive and may be taken as small as needed, but its value may
vary from one occurrence to another.

We allow repetitions in a displayed finite collection
$\{a_{1},\ldots,a_{n}\}$, but sums over this collection run over its
distinct elements, without multiplicity.

In Euler factorizations, $A\propto B$ means that $A$ is $B$ times
the omitted local factors.

\medskip

\textbf{Acknowledgements}
I wish to express my sincere gratitude to Henryk Iwaniec, Alex Kontorovich,
Chung-Hang Kwan, Wing Hong Leung, Stephen D. Miller, Peter Sarnak, and
Matthew Young for their helpful comments and valuable suggestions.

This work was initiated during my final year at Princeton University and
developed in substantial part with Princeton's continued hospitality
through August 2025.  The first version was completed during my time at
the California Institute of Technology. I gratefully acknowledge the support
of both institutions. I also thank Dinakar Ramakrishnan for his continued
interest in this work and for his encouragement throughout more than a year
of revisions.

\part{Global Fourier Identities and Spectral Reciprocity}

\section{Fourier Reciprocity Formulas}\label{sec2}
Partial Fourier expansions give the identities of types~\RNum{1} and~\RNum{2},
including the constant and degenerate terms for Eisenstein data.

\subsection{Automorphic Forms}
Let $B$ be the upper triangular subgroup of $G$, and let $P$ and $Q$
be the standard parabolic subgroups of types $(2,1)$ and $(1,2)$,
respectively.

Let $\omega$ be a unitary Hecke character of
$F^{\times}\backslash\mathbb{A}_{F}^{\times}$.  We consider a generic
automorphic representation $\pi$ of $[G]$ with central character $\omega$
in each of the following three cases.
\begin{itemize}
\item $\pi$ is \textit{cuspidal}, and $\varphi$ is a cusp form in $\pi$.
\item $\pi=\chi_{1}|\cdot|^{\nu_{1}}\boxplus
\chi_{2}|\cdot|^{\nu_{2}}\boxplus\chi_{3}|\cdot|^{\nu_{3}}$ is induced from
unitary Hecke characters $\chi_{i}$, where
$\boldsymbol{\nu}:=(\nu_{1},\nu_{2},\nu_{3})\in \mathbb{C}^{3}$ and
$\chi_{1}\chi_{2}\chi_{3}|\cdot|^{\nu_{1}+\nu_{2}+\nu_{3}}=\omega$.  We call this
the \textit{minimal-parabolic} case.

Let $f(\cdot;\boldsymbol{\chi},\boldsymbol{\nu})$ be a section in $\pi$.
Thus, for $t=\diag(t_{1},t_{2},t_{3})\in B(\mathbb{A}_{F})$,
$u\in N(\mathbb{A}_{F})$, and $g\in G(\mathbb{A}_{F})$, it satisfies
\begin{equation}\label{2.9}
f(tug;\boldsymbol{\chi},\boldsymbol{\nu})=|t_{1}|^{1+\nu_{1}}|t_{2}|^{\nu_{2}}|t_{3}|^{-1+\nu_{3}}\chi_{1}(t_{1})\chi_{2}(t_{2})\chi_{3}(t_{3})f(g;\boldsymbol{\chi},\boldsymbol{\nu}).
\end{equation}

Define the Eisenstein series
\begin{equation}\label{f2.5}
\varphi(g)=E(g,\boldsymbol{\chi},\boldsymbol{\nu}):=\sum_{\delta\in B(F)\backslash G(F)}f(\delta g;\boldsymbol{\chi},\boldsymbol{\nu}).
\end{equation}

When $\Re(\nu_{1}-\nu_{2})>1$ and $\Re(\nu_{2}-\nu_{3})>1$, the Eisenstein
series converges absolutely and uniformly on compact subsets of
$G(\mathbb{A}_{F})$.  It continues meromorphically to
$\boldsymbol{\nu}\in\mathbb{C}^{3}$.

\item $\pi=\chi|\cdot|^{\mu_{1}}\boxplus
\sigma|\cdot|^{\mu_{2}/2}$, where
$\boldsymbol{\mu}=(\mu_{1},\mu_{2})\in\mathbb{C}^{2}$, $\chi$ is a unitary
Hecke character, and $\sigma$ is a unitary cuspidal automorphic
representation of $\mathrm{GL}_{2}/F$ with central character
$\omega_{\sigma}$.  Here
$\chi\omega_{\sigma}|\cdot|^{\mu_{1}+\mu_{2}}=\omega$.  We call this the
\textit{maximal-parabolic} case.

Let $f(\cdot;\sigma,\chi,\boldsymbol{\mu})$ be a section in this
normalized induced representation. Thus
\begin{equation}\label{e2.3}
f\left(\begin{pmatrix}
a& \mathfrak{b}\\
& g'
\end{pmatrix}g
;\sigma,\chi,\boldsymbol{\mu}\right)=|a|^{1+\mu_{1}}|\det g'|^{\frac{-1+\mu_{2}}{2}}\chi(a)\sigma(g')f(g;\sigma,\chi,\boldsymbol{\mu}),
\end{equation}
where $a\in \mathbb{A}_{F}^{\times}$,
$\mathfrak{b}\in M_{1\times 2}(\mathbb{A}_{F})$,
$g'\in G'(\mathbb{A}_{F})$, and $g\in G(\mathbb{A}_{F})$. Here $\sigma(g')$ acts on the inducing representation, by right
translation in its $\mathrm{GL}_{2}$ variable.  We use the same notation
for the vector-valued section and its scalar realization, obtained
by evaluation at $I_{2}$ in that variable.

Define the Eisenstein series 
\begin{equation}\label{e2.4}
\varphi(g):=E(g,\sigma,\chi,\boldsymbol{\mu}):=\sum_{\delta\in Q(F)\backslash G(F)}f(\delta g;\sigma,\chi,\boldsymbol{\mu}).
\end{equation}

When $\Re(\mu_{1}-\mu_{2}/2)>3/2$, this Eisenstein series converges
absolutely and uniformly on compact subsets of $G(\mathbb{A}_{F})$.  It
continues meromorphically to $\boldsymbol{\mu}\in\mathbb{C}^{2}$.

\end{itemize}

\subsection{Fourier Expansions of Automorphic Forms}\label{sec2.2}
Let $B'$ be the upper triangular subgroup of $G'$ and $N'$ its
unipotent radical. Put
\begin{align*}
&B_{0}'(F):=\bigg\{\begin{pmatrix}
a & b\\
& 1
\end{pmatrix}:\ a\in F^{\times},\ b\in F\bigg\},\\
&B_{0}'^{\dag}(F):=\bigg\{\begin{pmatrix}
1 & b\\
& d
\end{pmatrix}:\ d\in F^{\times},\ b\in F\bigg\}.
\end{align*}

Let $\pi$ be a generic automorphic representation of $[G]$, and let
$\varphi$ be an automorphic form in $\pi$.  For $g\in G(\mathbb{A}_{F})$,
let
\begin{equation}\label{w2.5}
W_{\varphi}(g):=\int_{(F\backslash\mathbb{A}_{F})^{3}}\varphi\left(\begin{pmatrix}
1& a & b\\
& 1& c\\
&& 1
\end{pmatrix}g\right)\overline{\psi}(a+c)dadbdc
\end{equation}
be the Whittaker function of $\varphi$ relative to $\psi$.  Define also the
degenerate Whittaker functions
\begin{align*}
&W_{\varphi}^{\mathrm{degen}}(g):=\int_{(F\backslash\mathbb{A}_{F})^{3}}\varphi\left(\begin{pmatrix}
1& a& b\\
&1& c\\
&& 1
\end{pmatrix}g\right)\overline{\psi(c)}dadbdc,\\
&W_{\varphi}^{\mathrm{degen},\dag}(g):=\int_{(F\backslash\mathbb{A}_{F})^{3}}\varphi\left(\begin{pmatrix}
1& a& b\\
&1& c\\
&& 1
\end{pmatrix}g\right)\overline{\psi(a)}dadbdc.
\end{align*}

For the forms considered here, polynomial growth of the Fourier
coefficients and rapid decay of the Whittaker functions give absolute
and locally uniform convergence of the Fourier sums away from
Eisenstein poles, also after differentiation in $g$. We may therefore
integrate the sums termwise over the compact unipotent quotients.

For $g\in G(\mathbb{A}_{F})$, the Fourier expansion
\cite[\textsection 2.1]{Yuk93} is
\begin{equation}\label{2.1}
\varphi(g)=\mathcal{F}_{\mathrm{const}}\varphi(g)+\mathcal{F}_{\mathrm{degen}}\varphi(g)+\mathcal{F}_{\mathrm{gen}}\varphi(g),
\end{equation}
where
\begin{align*}
&\mathcal{F}_{\mathrm{const}}\varphi(g):=\int_{(F\backslash\mathbb{A}_{F})^{2}}\varphi\left(n_{13}(b)n_{23}(c)g\right)dbdc,\\
&\mathcal{F}_{\mathrm{degen}}\varphi(g):=\sum_{\delta\in B_{0}'(F)\backslash G'(F)}W_{\varphi}^{\mathrm{degen}}\left(\diag(\delta,1)g\right),\\
&\mathcal{F}_{\mathrm{gen}}\varphi(g):=\sum_{\delta\in N'(F)\backslash G'(F)}W_{\varphi}\left(\diag(\delta,1)g\right).
\end{align*}

Apply \eqref{2.1} to $\varphi^{\iota}(g):=\varphi(g^{\iota})$ at $\widetilde{w} g^{\iota}$. Automorphy under $\widetilde{w}$, with the character signs absorbed by reindexing the rational diagonal factors, gives the opposite expansion
\begin{equation}\label{2.2}
\varphi(g)=\mathcal{F}_{\mathrm{const}}^{\dag}\varphi(g)+\mathcal{F}_{\mathrm{degen}}^{\dag}\varphi(g)+\mathcal{F}_{\mathrm{gen}}^{\dag}\varphi(g),
\end{equation}
where 
\begin{align*}
&\mathcal{F}_{\mathrm{const}}^{\dag}\varphi(g):=\int_{(F\backslash\mathbb{A}_{F})^{2}}\varphi\left(n_{12}(a)n_{13}(b)g\right)dadb,\\
&\mathcal{F}_{\mathrm{degen}}^{\dag}\varphi(g):=\sum_{\delta\in B_{0}'^{\dag}(F)\backslash G'(F)}W_{\varphi}^{\mathrm{degen},\dag}\left(\diag(1,\delta)g\right),\\
&\mathcal{F}_{\mathrm{gen}}^{\dag}\varphi(g):=\sum_{\delta\in N'(F)\backslash G'(F)}W_{\varphi}\left(\diag(1,\delta)g\right).
\end{align*}

The two expansions give the following identity.
\begin{lemma}
Let $\pi$ be a generic automorphic representation of $[G]$ and $\varphi$ be an automorphic form in $\pi$.  Then for $g\in G(\mathbb{A}_{F})$ and $\delta\in G(F)$, we have  
\begin{equation}\label{2.3}
\mathcal{F}_{\mathrm{gen}}\varphi(g)=\mathcal{F}_{\mathrm{gen}}^{\dag}\varphi(\delta g)+\mathcal{F}_{\mathrm{const}}^{\dag}\varphi(\delta g)-\mathcal{F}_{\mathrm{const}}\varphi(g)+\mathcal{F}_{\mathrm{degen}}^{\dag}\varphi(\delta g)-\mathcal{F}_{\mathrm{degen}}\varphi(g).
\end{equation} 
\end{lemma} 

\begin{proof}
Compare \eqref{2.1} at $g$ with \eqref{2.2} at $\delta g$, using
$\varphi(\delta g)=\varphi(g)$.
\end{proof}

In particular, when $\varphi$ is cuspidal, all constant and degenerate
terms vanish. Taking $\delta=I_{3}$ in \eqref{2.3} gives
\begin{align*}
\mathcal{F}_{\mathrm{gen}}\varphi(g)=\mathcal{F}_{\mathrm{gen}}^{\dag}\varphi(g).
\end{align*} 

\subsection{Fourier Reciprocity Formula of Type \RNum{1}}
For a continuous function $h$ on $G(\mathbb{A}_{F})$ that is left invariant under $\diag(N'(F),1)$, we define    
\begin{equation}\label{f2.9}
\mathcal{W}(h):=\int_{F\backslash\mathbb{A}_{F}}h\left(n_{12}(a)\right)\overline{\psi(a)}da.
\end{equation}

\begin{prop}\label{prop2.2}
Let $\varphi$ be a generic automorphic form on $[G]$ with central character $\omega$. Then
\begin{equation}\label{e2.20}
\mathcal{W}(\mathcal{F}_{\mathrm{gen}}\varphi)=\mathcal{W}_{\mathrm{gen}}(\varphi)+\mathcal{W}_{\mathrm{degen}}^{\dag}(\varphi)-\mathcal{W}_{\mathrm{degen}}(\varphi),
\end{equation}
where 	
\begin{align*}
&\mathcal{W}_{\mathrm{gen}}(\varphi):=\sum_{\beta\in F}\sum_{\alpha\in F^{\times}}W_{\varphi}\left(n_{32}(\beta)\diag(\alpha,\alpha,1)\right),\\
&\mathcal{W}_{\mathrm{degen}}^{\dag}(\varphi):=\sum_{\alpha\in F^{\times}}W_{\varphi}^{\mathrm{degen},\dag}\left(w_{2}n_{23}(-1)\diag(\alpha,\alpha,1)\right),\\
&\mathcal{W}_{\mathrm{degen}}(\varphi):=\sum_{\alpha\in F^{\times}}\int_{\mathbb{A}_{F}}W_{\varphi}^{\mathrm{degen}}\left(w_{1}n_{12}(c')\diag(\alpha,\alpha,1)\right)\overline{\psi(c')}dc'.
\end{align*}
\end{prop} 

See \textsection\ref{sec2.3.2} for the proof.

\subsubsection{The constant and degenerate terms}
\begin{lemma}
Let $\varphi$ be a generic automorphic form on $[G]$ with central character $\omega$. Then 
\begin{multline}\label{2.7}
\mathcal{W}(\mathcal{F}_{\mathrm{const}}^{\dag}\varphi-\mathcal{F}_{\mathrm{const}}\varphi+\mathcal{F}_{\mathrm{degen}}^{\dag}\varphi-\mathcal{F}_{\mathrm{degen}}\varphi)=\mathcal{W}_{\mathrm{degen}}^{\dag}(\varphi)-\mathcal{W}_{\mathrm{degen}}(\varphi).
\end{multline}
\end{lemma}

\begin{proof}
The shift $a\mapsto a-c'$ in the constant term gives
\begin{equation}\label{e2.10}
 \mathcal{W}(\mathcal{F}_{\mathrm{const}}^{\dag}\varphi)=0.
\end{equation}
The other constant term is
\begin{equation}\label{e2.15}
 \mathcal{W}(\mathcal{F}_{\mathrm{const}}\varphi)
 =W_{\varphi}^{\mathrm{degen},\dag}(I_{3}).
\end{equation}

The Bruhat decomposition
\begin{align*}
 G'(F)=B'(F)\bigsqcup B'(F)w'N'(F)
\end{align*}
gives the representatives
\begin{equation}\label{2.5}
 B_{0}'^{\dag}(F)\backslash G'(F)
 =\{\alpha'I_{2}:\ \alpha'\in F^{\times}\}
 \bigsqcup\left\{w'\begin{pmatrix}\alpha&\beta\\&\alpha\end{pmatrix}:\ 
 \alpha\in F^{\times},\ \beta\in F\right\}.
\end{equation}
Write $\mathcal{F}_{\mathrm{degen},1}^{\dag}\varphi$ for the sum
over $\alpha'I_{2}$ and $\mathcal{F}_{\mathrm{degen},2}^{\dag}\varphi$
for the sum over $(\alpha,\beta)$ in \eqref{2.5}.
In the sum over $\alpha'I_{2}$, the shift $a\mapsto a-\alpha'^{-1}c'$
leaves $\psi((\alpha'^{-1}-1)c')$. Integration in $c'$ leaves only
$\alpha'=1$, giving \eqref{e2.15}.
For the sum over $(\alpha,\beta)$, use
\begin{align*}
 w_{2}\begin{pmatrix}1&&\\&\alpha&\beta\\&&\alpha\end{pmatrix}
 n_{12}(c')
 =n_{12}(-\alpha^{-2}\beta c')n_{13}(\alpha^{-1}c')
 w_{2}\begin{pmatrix}1&&\\&\alpha&\beta\\&&\alpha\end{pmatrix}.
\end{align*}
After $a\mapsto a+\alpha^{-2}\beta c'$, the $c'$-integral is
\begin{align*}
 \int_{F\backslash\mathbb{A}_{F}}\psi(-(1+\alpha^{-2}\beta)c')\,dc'
 =\mathbf{1}_{\beta=-\alpha^{2}}.
\end{align*} Invariance under $\diag(1,1,\alpha)\in G(F)$
and the rational center then gives
\begin{align*}
 \mathcal{W}(\mathcal{F}_{\mathrm{degen},2}^{\dag}\varphi)
 =\mathcal{W}_{\mathrm{degen}}^{\dag}(\varphi).
\end{align*}
Consequently,
\begin{equation}\label{eq2.16}
 \mathcal{W}(\mathcal{F}_{\mathrm{degen}}^{\dag}\varphi)
 -\mathcal{W}(\mathcal{F}_{\mathrm{const}}\varphi)
 =\mathcal{W}_{\mathrm{degen}}^{\dag}(\varphi).
\end{equation}

For $B_{0}'(F)\backslash G'(F)$ the two Bruhat families are
$\alpha'I_{2}$ and $\alpha w'n(\beta)$, with
$\alpha',\alpha\in F^{\times}$ and $\beta\in F$.
Denote their contributions by
$\mathcal{F}_{\mathrm{degen},1}\varphi$ and
$\mathcal{F}_{\mathrm{degen},2}\varphi$.
The first is invariant under right $n_{12}(c')$, so
\begin{align*}
 \mathcal{W}(\mathcal{F}_{\mathrm{degen},1}\varphi)=0.
\end{align*}
For the family $\alpha w'n(\beta)$, the unipotents combine as
$n(\beta+c')$. Put $c=\beta+c'$. Since $\psi(\beta)=1$,
the sum and the integral over $F\backslash\mathbb{A}_{F}$ give
one integral over $\mathbb{A}_{F}$ against $\overline{\psi}(c)$.
Commuting the scalar $\alpha I_{2}$ with $w'n(c)$ gives
\begin{equation}\label{2.19}
 \mathcal{W}(\mathcal{F}_{\mathrm{degen},2}\varphi)
 =\mathcal{W}_{\mathrm{degen}}(\varphi).
\end{equation}
Together with \eqref{e2.10} and \eqref{eq2.16}, these identities
prove \eqref{2.7}.
\end{proof}

\subsubsection{Proof of Proposition \ref{prop2.2}}\label{sec2.3.2}
The Bruhat decomposition gives
\begin{equation}\label{eq2.21}
N'(F)\backslash \mathrm{GL}_{2}(F)=\bigsqcup_{\alpha_{1}',\alpha_{2}'\in F^{\times}}\bigg\{\begin{pmatrix}
&\alpha_{1}'\\
\alpha_{2}'
\end{pmatrix} \bigg\}\bigsqcup \bigsqcup_{\substack{\alpha_{1}, \alpha_{2}\in F^{\times}\\
\beta\in F }}\bigg\{\begin{pmatrix}
\alpha_{1}\\
\beta & \alpha_{2}
\end{pmatrix}\bigg\}.
\end{equation}

Let $\mathcal{F}_{\mathrm{gen},1}^{\dag}\varphi$ and $\mathcal{F}_{\mathrm{gen},2}^{\dag}\varphi$ be the two sums in \eqref{eq2.21}. For the first family, conjugation sends $n_{12}(c')$ to $n_{13}(\alpha_{2}'^{-1}c')$, so its contribution vanishes. For the second, it sends $n_{12}(c')$ to $n_{12}(\alpha_{1}^{-1}c')$. Integration against $\overline{\psi(c')}$ forces $\alpha_{1}=1$, giving
\begin{equation}\label{2.23}
\mathcal{W}(\mathcal{F}_{\mathrm{gen},1}^{\dag}\varphi)\equiv 0,\qquad \mathcal{W}(\mathcal{F}_{\mathrm{gen},2}^{\dag}\varphi)=\sum_{\alpha_{2}\in F^{\times}}\sum_{\beta\in F}W_{\varphi}\left(n_{32}(\beta)\diag(1,1,\alpha_{2})\right).
\end{equation}

In the second sum in \eqref{2.23} we replace $\alpha_{2}$ by its inverse.
Rational central invariance gives $\mathcal{W}_{\mathrm{gen}}(\varphi)$.  Thus \eqref{e2.20} follows from
\eqref{2.3}, \eqref{2.7}, and \eqref{2.23}.

\subsection{Fourier Reciprocity Formula of Type \RNum{2}}\label{sec2.4}
For a continuous function $h$ on $G(\mathbb{A}_{F})$ that is left invariant under $\diag(N'(F),1)$, define, whenever the sum converges absolutely,
\begin{align*}
\mathcal{I}(h):=\sum_{\alpha\in F^{\times}}\int_{F\backslash\mathbb{A}_{F}}h\left(n_{12}(a)\diag(\alpha,1,1)\right)\overline{\psi(a)}da.
\end{align*}

\begin{prop}\label{prop2.5}
Let $\pi$ be generic automorphic on $[G]$, and let $\varphi\in\pi$. Then
\begin{multline}\label{ew2.26}
\mathcal{I}(\mathcal{F}_{\mathrm{gen}}\varphi)+\mathcal{I}_{\mathrm{degen}}(\varphi)-\mathcal{I}_{\mathrm{gen}}(\varphi)\\
=\mathcal{I}(\mathcal{F}_{\mathrm{gen}}\pi(w_{2})\varphi)+\mathcal{I}_{\mathrm{degen}}(\pi(w_{2})\varphi)-\mathcal{I}_{\mathrm{gen}}(\pi(w_{2})\varphi),
\end{multline}	
where 
\begin{align*}
&\mathcal{I}_{\mathrm{gen}}(\varphi):=\sum_{\alpha\in F^{\times}}\sum_{\gamma\in F^{\times}}W_{\varphi}\left(\diag(\alpha\gamma,\gamma,1)\right),\\
&\mathcal{I}_{\mathrm{degen}}(\varphi):=\sum_{\alpha\in F^{\times}}\sum_{\gamma\in F^{\times}}\int_{\mathbb{A}_{F}}W_{\varphi}^{\mathrm{degen}}\left(w_{1}n_{12}(a)\diag(\alpha\gamma,\gamma,1)\right)\overline{\psi(a)}da.
\end{align*}
\end{prop}
See \textsection\ref{sect2.4.2} for the proof.

\subsubsection{Fourier expansion and automorphy}
For continuous $h$ on $G(\mathbb{A}_{F})$ left invariant under
$\diag(N'(F),1)$, put
\begin{align*}
\mathcal{I}_{0}(h):=\int_{F\backslash\mathbb{A}_{F}}h\left(n_{12}(a)\right)da.
\end{align*}  

\begin{lemma}\label{lem2.4}
Let $\pi$ be a generic automorphic representation of $[G]$ with central character $\omega$, $\varphi\in\pi$, and $\delta\in G(F)$. Then
\begin{align*}
\mathcal{I}(\varphi)+\mathcal{I}_{0}(\varphi)=\mathcal{I}(\pi(\delta)\varphi)+\mathcal{I}_{0}(\pi(\delta)\varphi).
\end{align*}
\end{lemma}

\begin{proof}
Fourier expansion gives
\begin{equation}\label{e2.31}
\varphi(I_{3})=\mathcal{I}(\varphi)+\mathcal{I}_{0}(\varphi).
\end{equation}
For $\pi(\delta)\varphi$, the left-hand side of \eqref{e2.31} is
$\varphi(\delta)=\varphi(I_{3})$.
\end{proof}

\subsubsection{Proof of Proposition \ref{prop2.5}}\label{sect2.4.2}
Expand $\mathcal{I}(\varphi)$ by \eqref{2.1} and
$\mathcal{I}_{0}(\varphi)$ by \eqref{2.2}, summing the Bruhat terms in
$\mathcal{I}$ over each right $N'(F)$-orbit first.
The subgroup defining $\mathcal{F}_{\mathrm{const}}^{\dag}$ is
normalized by $w_{2}$, so automorphy gives
\begin{equation}\label{f2.36}
 \mathcal{I}_{0}(\mathcal{F}_{\mathrm{const}}^{\dag}\pi(w_{2})\varphi)
 =\mathcal{I}_{0}(\mathcal{F}_{\mathrm{const}}^{\dag}\varphi).
\end{equation}
Applying \eqref{2.19} to the right translates in $\mathcal{I}$ gives
\begin{equation}\label{fc2.37}
 \mathcal{I}(\mathcal{F}_{\mathrm{degen}}\varphi)
 =\mathcal{I}_{\mathrm{degen}}(\varphi).
\end{equation}

For the opposite generic sum, use the two Bruhat families in
\eqref{eq2.21}. The lower-triangular family contributes zero:
its $a$-integral has the factor
$\int_{F\backslash\mathbb{A}_{F}}\psi(\gamma^{-1}a)\,da=0$.
The other family gives
\begin{align*}
 \mathcal{I}_{0}(\mathcal{F}_{\mathrm{gen}}^{\dag}\pi(w_{2})\varphi)
 =\sum_{\alpha,\gamma\in F^{\times}}
 W_{\varphi}(\diag(\alpha,\gamma,1))
 =\mathcal{I}_{\mathrm{gen}}(\varphi).
\end{align*}
The last equality uses $(\alpha,\gamma)\mapsto(\alpha\gamma,\gamma)$.
In the opposite degenerate sum, the $a$-integral vanishes on the
scalar family in \eqref{2.5} and forces $\beta=0$ in the family indexed
by $(\alpha,\beta)$. Invariance under the rational center reduces the sum
to $\sum_{\alpha\in F^{\times}}W_{\varphi}^{\mathrm{degen},\dag}(\diag(\alpha,1,1))$.
Thus
\begin{equation}\label{e2.39}
 \mathcal{I}_{0}(\mathcal{F}_{\mathrm{degen}}^{\dag}\pi(w_{2})\varphi)
 =\mathcal{I}(\mathcal{F}_{\mathrm{const}}\varphi).
\end{equation}

In the equality from Lemma~\ref{lem2.4} with $\delta=w_{2}$,
\eqref{f2.36} cancels the two $\mathcal{F}_{\mathrm{const}}^{\dag}$ contributions.
By \eqref{e2.39} and $w_{2}^{2}=I_{3}$, the $\mathcal{I}$-term of
$\mathcal{F}_{\mathrm{const}}$ on each side cancels the
$\mathcal{I}_{0}$-term of $\mathcal{F}_{\mathrm{degen}}^{\dag}$ on the other.
This leaves
\begin{align*}
 \mathcal{I}(\mathcal{F}_{\mathrm{gen}}\varphi)
 +\mathcal{I}_{\mathrm{degen}}(\varphi)
 +\mathcal{I}_{\mathrm{gen}}(\pi(w_{2})\varphi),
\end{align*}
and the same expression with $\varphi$ replaced by $\pi(w_{2})\varphi$.
Their equality proves \eqref{ew2.26}.

\section{Spectral Decomposition in the Convergence Region}\label{sec3}
Use the spectral measure of \textsection\ref{sec:spectral-normalization}
and the local Rankin--Selberg normalizations of \textsection\ref{sec3.2}.
For $\pi$ as in \textsection\ref{sec2}, with central character $\omega$, put
\begin{equation}\label{f2.2}
\begin{cases}
\mathcal{R}_{\mathrm{max}}:=\big\{(\mu_{1},\mu_{2})\in\mathbb{C}^{2}:\ \Re(\mu_{1}-\mu_{2}/2)>3/2\big\},\\
\mathcal{R}_{\mathrm{min}}:=\big\{(\nu_{1},\nu_{2},\nu_{3})\in\mathbb{C}^{3}:\ \Re(\nu_{1}-\nu_{2})>1,\ \Re(\nu_{2}-\nu_{3})>2\big\}.
\end{cases}
\end{equation}

\begin{defn}\label{defn3.2}
For the three types in \textsection\ref{sec2}, call $\pi$
\textit{convergent} if it is cuspidal, or if its inducing parameters
lie in $\mathcal{R}_{\mathrm{min}}$ or $\mathcal{R}_{\mathrm{max}}$,
respectively.

Let $\varphi$ be a vector in $\pi$.  In the minimal- and
maximal-parabolic cases, $\varphi$ is the Eisenstein series
\eqref{f2.5} and \eqref{e2.4}, respectively. We use the convergence conditions
\begin{equation}\label{f3.2}
\begin{aligned}
&\min\{\Re(s_{1})-1/2,\Re(s_{2}-s_{1})-1\}\\
&\qquad >
\begin{cases}
0,&\pi\text{ cuspidal},\\
\max_{i}|\Re(\nu_{i})|,&\pi\text{ minimal-parabolic},\\
\max\{|\Re(\mu_{1})|,|\Re(\mu_{2})|/2\},&\pi\text{ maximal-parabolic}.
\end{cases}
\end{aligned}
\end{equation} 
\end{defn}

\subsection{\texorpdfstring{The $\mathrm{GL}_{2}$-spectrum}{The GL2 Spectrum}}
For unitary Hecke characters $\xi_{1},\xi_{2}$, let
$H(\xi_{1},\xi_{2})$ be the space of functions $h\in L^{2}(K')$ satisfying
\begin{align*}
h\left(\begin{pmatrix}
a& b\\
&d
\end{pmatrix}k\right)=\xi_{1}(a)\xi_{2}(d)h(k),\ \forall\ k\in K',\ 
 \begin{pmatrix}
a& b\\
&d
\end{pmatrix}\in B'(\mathbb{A}_{F})\cap K'.
\end{align*}

For $K'$-finite $h\in H(\xi_{1},\xi_{2})$ and $\lambda\in \mathbb{C}$, define
\begin{align*}
S(h)\left(\begin{pmatrix}
a& b\\
&d
\end{pmatrix}k,\lambda\right):=\xi_{1}(a)\xi_{2}(d)\Big|\frac{a}{d}\Big|^{\lambda+1/2}h(k),
\end{align*}
where $a,d\in \mathbb{A}_{F}^{\times}$, $b\in \mathbb{A}_{F}$, and $k\in K'$.

Define the associated Eisenstein series by
\begin{align*}
E(g,h,\lambda):=\sum_{\gamma\in B'(F)\backslash G'(F)}S(h)(\gamma g,\lambda).
\end{align*}

The series $E(g,h,\lambda)$ converges absolutely for
$\Re(\lambda)>1/2$ and extends meromorphically to $\mathbb{C}$. In
$\Re(\lambda)\geq0$, its possible poles are simple and satisfy
$\xi_{1}\xi_{2}^{-1}|\cdot|^{2\lambda}=|\cdot|$.

For $\lambda\in i\mathbb{R}$, use the invariant inner product
\begin{equation}\label{eq:inducing-global-norm}
\langle h_{1},h_{2}\rangle_{\mathrm{Ind}}
:=\rho_{F}\int_{K'}h_{1}(k)\overline{h_{2}(k)}\,dk.
\end{equation}

Fix an orthonormal basis
$\mathfrak{B}(\xi_{1},\xi_{2})$ of $K'$-finite functions in
$H(\xi_{1},\xi_{2})$. Let $\omega'$ be a unitary Hecke character. For
$\sigma\in \mathcal{A}_{0}([G'],\omega')$, let
$\mathfrak{B}(\sigma)$ be an orthonormal basis of $\sigma$.
For a function $\varphi'$ on $[G']$ with central character $\omega'$
that is square-integrable modulo the center, the $L^{2}$ spectral expansion is
\begin{equation}\label{3.2}
\varphi'(g)=\varphi_{\mathrm{res}}'(g)+\varphi_{\mathrm{cusp}}'(g)+\varphi_{\mathrm{Eis}}'(g),
\end{equation}
where 
\begin{align*}
&\varphi_{\mathrm{res}}'(g)=\frac{1}{\Vol([\overline{G'}])}\sum_{\substack{\xi\in \widehat{F^{\times}\backslash\mathbb{A}_{F}^{\times}},\ \xi^{2}=\omega'}}\langle\varphi',\xi\circ\det\rangle\xi(\det g),
\\
&\varphi_{\mathrm{cusp}}'(g)=\sum_{\sigma\in \mathcal{A}_{0}([G'],\omega')}\sum_{\phi\in\mathfrak{B}(\sigma)}\langle\varphi',\phi\rangle\phi(g),\\
&\varphi_{\mathrm{Eis}}'(g)=\sum_{\substack{\xi\in \widehat{C_{F}}}} \frac{1}{4\pi i}\int_{i\mathbb{R}}\sum_{h\in \mathfrak{B}(\xi,\overline{\xi}\omega')}\langle\varphi',E(\cdot,h,\lambda)\rangle E(g,h,\lambda)d\lambda.
\end{align*}

\subsubsection{Spectral measure and the adjoint denominator}\label{sec:spectral-normalization}
\phantomsection\label{def:spectral-measure}
Write $d\mu_{\sigma}$ for the unweighted measure in \eqref{3.2}. For an integrable function $A$ on the generic spectrum with central character $\omega'$, write
\begin{align*}
 \int A(\sigma)\,d\mu_{\sigma}
=\sum_{\sigma\in\mathcal{A}_{0}([G'],\omega')}A(\sigma)
+\frac{1}{4\pi i}\sum_{\xi\in\widehat{C_{F}}}
\int_{i\mathbb{R}}A(\sigma_{\xi,\lambda})\,d\lambda,
\end{align*}
where
$\sigma_{\xi,\lambda}
=\xi|\cdot|^{\lambda}\boxplus\overline{\xi}\omega'|\cdot|^{-\lambda}$.

If $H\in L^{2}(F^{\times}\backslash\mathbb{A}_{F}^{\times})$ and
\begin{align*}
\widehat{H}(\xi,\lambda)
:=\int_{F^{\times}\backslash\mathbb{A}_{F}^{\times}}
H(y)\xi(y)|y|^{\lambda}\,d^{\times} y,
\end{align*}
then the $\mathrm{GL}_{1}$-spectral expansion, understood in $L^{2}$, reads
\begin{equation}\label{eq:gl1-plancherel}
H(x)=\frac{1}{2\pi i\rho_{F}}
\sum_{\xi\in\widehat{C_{F}}}
\int_{i\mathbb{R}}\widehat{H}(\xi,\lambda)
\overline{\xi(x)}|x|^{-\lambda}\,d\lambda.
\end{equation}

Define
\begin{equation}\label{eq:regularized-adjoint}
L^{\mathrm{reg}}(1,\sigma,\Ad):=
\begin{cases}
L(1,\sigma,\Ad),&\sigma\text{ cuspidal},\\
\underset{z=1}{\Res}\zeta_{F}(z)\,
L(1+2\lambda,\xi^{2}\overline{\omega'})
L(1-2\lambda,\overline{\xi}^{\,2}\omega'),
&\sigma=\sigma_{\xi,\lambda}.
\end{cases}
\end{equation}
For Eisenstein data, \eqref{eq:regularized-adjoint} is the residue of
$L(s,\sigma,\Ad)$ at $s=1$ for generic $\lambda$; at exceptional
parameters use meromorphic continuation in $\lambda$.
On $\lambda=it$ it equals
$\underset{z=1}{\Res}\zeta_{F}(z)\,
|L(1+2it,\xi^{2}\overline{\omega'})|^{2}$.
In global spectral moments below, the shorthand $L(1,\sigma,\Ad)$
in a denominator always means \eqref{eq:regularized-adjoint}.
Local factors $L_{v}(1,\sigma_{v},\Ad)$ have their usual meaning.

\subsection{Rankin--Selberg Integrals}\label{sec3.2}
Let $\pi=\otimes_{v}'\pi_{v}$ and $\sigma=\otimes_{v}'\sigma_{v}$ be generic
automorphic representations of $[G]$ and $[G']$, respectively, and let
$\varphi\in\pi$ and $\phi\in\sigma$ be automorphic forms.  The Whittaker
function $W_{\varphi}$ is defined in \eqref{w2.5}.  In parallel, define
\begin{align*}
W_{\phi}(g'):=\int_{F\backslash\mathbb{A}_{F}}\phi\left(n(c)g'\right)\psi(c)dc,\ \ g'\in G'(\mathbb{A}_{F}). 
\end{align*}

Let $\chi$ be a unitary Hecke character. For $\Re(s)$ sufficiently large,
define the Rankin--Selberg integrals
\begin{align*}
&\Psi(1/2+s,W_{\varphi},W_{\phi})
:=\int_{N'(\mathbb{A}_{F})\backslash G'(\mathbb{A}_{F})}
W_{\varphi}(\diag(g',1))W_{\phi}(g')|\det g'|^{s}\,dg',\\
&\Psi(1/2+s,W_{\phi},\chi)
:=\int_{\mathbb{A}_{F}^{\times}}W_{\phi}(a(y))\chi(y)|y|^{s}\,d^{\times} y.
\end{align*}
These integrals converge absolutely in right half-planes and continue
meromorphically, representing $\Lambda(1/2+s,\pi\times\sigma)$ and
$\Lambda(1/2+s,\sigma\times\chi)$, respectively; see \cite{JPSS83,Jac09}.
Write $\Psi_{v}$ for the local integrals and $\widetilde{\Psi}$ for the
meromorphic continuations of the global integrals $\Psi$.

For unitary generic $\sigma_{v}$, local orthonormal bases use the Kirillov
inner product at Archimedean places and, for $v<\infty$, the inner product
\cite[(2.2)]{MV10}
\begin{align*}
\langle U,V\rangle_{v}
:=\frac{\zeta_{v}(2)}{\zeta_{v}(1)L_{v}(1,\sigma_{v},\Ad)}
\int_{F_{v}^{\times}}U(a(y))\overline{V(a(y))}\,d^{\times}y.
\end{align*}

\phantomsection\label{def:jacquet-euler-factor}
For minimal Eisenstein data $\pi=\boxplus_{i=1}^{3}\chi_{i}|\cdot|^{\nu_{i}}$, put
\begin{align*}
 \mathbf{L}(\pi_{v}):=\prod_{i<j}L_{v}(1+\nu_{i}-\nu_{j},\chi_{i,v}\overline{\chi}_{j,v}),
 \qquad \mathbf{L}_{\fin}(\pi):=\prod_{v<\infty}\mathbf{L}(\pi_{v}).
\end{align*}
At finite places where $\pi_{v}$ and $\psi_{v}$ are unramified, the flat spherical section normalized at $I_{3}$ gives $W_{v}(I_{3})=\mathbf{L}(\pi_{v})^{-1}$. At the other places use the Jacquet integrals, multiplied by the constants in \eqref{eq:jacquet-normalization-transport}.

\subsubsection{Local functional equations for $\mathrm{GL}_{3}\times\mathrm{GL}_{2}$}
For decomposable data, write the global Whittaker functions as restricted tensor products
$W_{\varphi}=\prod_{v}W_{v}$ and $W_{\phi}=\prod_{v}W_{v}'$.

For each $v\in\Sigma_{F}$, $g\in G(F_{v})$, and $g'\in G'(F_{v})$, let
\begin{align*}
\widetilde{W}_{v}(g):=W_{v}(\widetilde{w}g^{\iota}),\qquad \text{and}\ 
\widetilde{W}_{v}'(g'):=W_{v}'(w'g'^{\iota}).
\end{align*}

With the contragredient Whittaker coefficients taken relative to the
opposite character, we have
$W_{\varphi^{\iota}}=\prod_{v}\widetilde{W}_{v}$ and
$W_{\phi^{\iota}}=\prod_{v}\widetilde{W}_{v}'$.
The quotient $\Psi_{v}(1/2+s,W_{v},W_{v}')/L_{v}(1/2+s,\pi_{v}\times\sigma_{v})$
extends to an entire function of $s$. The local functional equation is
\begin{align*}
\frac{\Psi_{v}(1/2+s,\widetilde{W}_{v},\widetilde{W}_{v}')}
{L_{v}(1/2+s,\widetilde{\pi}_{v}\times\widetilde{\sigma}_{v})}
=\varepsilon(1/2-s,\pi_{v}\times\sigma_{v},\psi_{v})
\frac{\Psi_{v}(1/2-s,W_{v},W_{v}')}{L_{v}(1/2-s,\pi_{v}\times\sigma_{v})},
\end{align*}
where the $\varepsilon$-factor is taken relative to $\psi_{v}$.

\subsubsection{Local functional equations for $\mathrm{GL}_{2}\times\mathrm{GL}_{1}$}
Write $\chi=\otimes_{v}\chi_{v}$. Standard local theory shows that
$\Psi_{v}(1/2+s,W_{v}',\chi_{v})/L_{v}(1/2+s,\sigma_{v}\times\chi_{v})$
extends to an entire function of $s$. The local functional equation is
\begin{align*}
\frac{\Psi_{v}(1/2+s,\widetilde{W}_{v}',\overline{\chi}_{v})}
{L_{v}(1/2+s,\widetilde{\sigma}_{v}\times\overline{\chi}_{v})}
=\chi_{v}(-1)\varepsilon(1/2-s,\sigma_{v}\times\chi_{v},\overline{\psi}_{v})
\frac{\Psi_{v}(1/2-s,W_{v}',\chi_{v})}{L_{v}(1/2-s,\sigma_{v}\times\chi_{v})}.
\end{align*}
Here the $\varepsilon$-factor is taken relative to $\overline{\psi}_{v}$, since $W_{v}'$ is in the $\overline{\psi}_{v}$-Whittaker model. The factors $\chi_{v}(-1)$ have global product $1$.

\subsection{Spectral Reciprocity Formula of Type \RNum{1}}\label{sec3.3}
\subsubsection{An auxiliary function} 
Let $s_{0}\ggg1$ be real and suppose that $\Re(s)-s_{0}\ggg1$. Here and below, $\ggg1$ means sufficiently large.
For $x\in\mathbb{A}_{F}^{\times}$, define
\begin{equation}\label{cf3.9}
J_{s_{0},s}(x):=|x|^{s_{0}+1/2}\int_{\mathbb{A}_{F}^{\times}}
 \sum_{\beta\in F^{\times}}
 W_{\varphi^{\iota}}\bigl(\diag(z,1,1)n_{21}(\beta xz)\widetilde{w}\bigr)
 \omega\omega'(z)|z|^{2s}\,d^{\times}z.
\end{equation}

\begin{lemma}\label{lem3.3}
Suppose $s_{0}\ggg 1$ and $\Re(s)-s_{0}\ggg 1$. Then $J_{s_{0},s}(\cdot)\in L^{2}(F^{\times}\backslash\mathbb{A}_{F}^{\times})$. 
\end{lemma}

\begin{proof}
By linearity, take $W_{\varphi^{\iota}}$ factorizable.
\begin{itemize}
\item At an Archimedean place, Iwasawa decomposition and
\cite[Proposition~3.1]{Jac09} give, for sufficiently large $l_{1}$
and every positive integer $l_{2}$,
\begin{equation}\label{3.17}
W_{\varphi^{\iota},v}\left(\diag(z_{v},1,1)n_{21}(y_{v})\widetilde{w}\right)
 \ll_{\varphi,l_{1},l_{2}}|z_{v}|_{v}^{-l_{1}}
 (1+|z_{v}|_{v}+|y_{v}|_{v})^{-l_{2}}.
\end{equation}

\item At finite $v$, right invariance under a compact open subgroup
of $K_{v}$ and Iwasawa decomposition give compact support in some
$C_{v}\subset F_{v}^{2}$. If $\varphi$ is right-invariant by $K_{v}$ and $\psi_{v}$ is unramified, then
\begin{equation}\label{3.18}
W_{\varphi^{\iota},v}\left(\diag(z_{v},1,1)n_{21}(y_{v})\widetilde{w}\right)=\mathbf{1}_{\mathcal{O}_{v}}(y_{v})W_{\varphi^{\iota},v}\left(\diag(z_{v},1,1)\right). 
\end{equation}

At an arbitrary finite place, Whittaker asymptotics give, for sufficiently large $l_{1}$,
\begin{equation}\label{3.19}
W_{\varphi^{\iota},v}\left(\diag(z_{v},1,1)n_{21}(y_{v})\widetilde{w}\right)\ll |z_{v}|_{v}^{-l_{1}}\mathbf{1}_{C_{v}}(y_{v},z_{v}),
\end{equation}
where the implied constant depends on $\varphi$.
\end{itemize}

Unfolding one copy of the $\beta$-sum in $|J|^{2}$ and making the change $x\mapsto xz'^{-1}$ gives
\begin{multline*}
\int_{F^{\times}\backslash\mathbb{A}_{F}^{\times}}|J_{s_{0},s}(x)|^{2}d^{\times}x
=\int_{\mathbb{A}_{F}^{\times}}\int_{\mathbb{A}_{F}^{\times}}\int_{\mathbb{A}_{F}^{\times}}\sum_{\beta\in F^{\times}}\\
W_{\varphi^{\iota}}\bigl(\diag(z,1,1)n_{21}(xz/z')\widetilde{w}\bigr)
\overline{W_{\varphi^{\iota}}\left(\diag(z',1,1)n_{21}(\beta x)\widetilde{w}\right)}
\omega\omega'(z)\overline{\omega}\overline{\omega}'(z')\\
|z|^{2s}|z'|^{2\overline{s}-2s_{0}-1}|x|^{2s_{0}+1}
\,d^{\times}x\,d^{\times}z\,d^{\times}z'.
\end{multline*}
For absolute convergence, take $l_{1}$ large enough to absorb the
fixed character factors and the growth of the spherical coefficients.
Use \eqref{3.17}--\eqref{3.19} and bound the sum over
$\beta\neq0$ by Poisson summation in the variable $\beta x$,
which gives a factor $|x|^{-1}$. After $x\mapsto xz'/z$,
the powers of $|x|,|z|,|z'|$ are respectively $2s_{0}$,
$2\Re(s)-2s_{0}-l_{1}$, and $2\Re(s)-l_{1}-1$.
These exceed $1$ for $s_{0},\Re(s)-s_{0}\ggg1$.
Take $l_{2}$ large to ensure convergence at infinity. Hence $J_{s_{0},s}(\cdot)\in L^{2}(F^{\times}\backslash\mathbb{A}_{F}^{\times})$.
\end{proof}

\subsubsection{Spectral reciprocity formula of type \RNum{1}}
Let $\Re(s)\ggg 1$. Define 
\begin{equation}\label{s3.3}
J_{\mathrm{spec}}(s,\varphi,\omega'):=J_{\mathrm{cusp}}(s,\varphi,\omega')+J_{\mathrm{Eis}}(s,\varphi,\omega'),
\end{equation}
where 
\begin{align*}
&J_{\mathrm{cusp}}(s,\varphi,\omega'):=\sum_{\sigma\in \mathcal{A}_{0}([G'],\omega')}\sum_{\phi\in\mathfrak{B}(\sigma)}\Psi(1/2+s,W_{\varphi},W_{\phi})\overline{W_{\phi}}(I_{2}),\\
&J_{\mathrm{Eis}}(s,\varphi,\omega'):=\sum_{\xi\in \widehat{C_{F}}} \frac{1}{4\pi i}\int_{i\mathbb{R}}\sum_{h\in \mathfrak{B}(\xi,\overline{\xi}\omega')} \Psi\left(1/2+s,W_{\varphi},W_{E(\cdot,h,\lambda)}\right)\overline{W_{E(\cdot,h,\lambda)}}(I_{2})d\lambda.
\end{align*}

\begin{thm}\label{thm3.3.}
Suppose $\pi$ is convergent in the sense of Definition~\ref{defn3.2},
and let $\varphi\in\pi$.  If $s_{0}\ggg1$ and
$\Re(s)-s_{0}\ggg1$, then
\begin{equation}\label{e3.3}
 J_{\mathrm{spec}}(s,\varphi,\omega')
 =\bigl(J_{\mathrm{sing}}+J_{\mathrm{dual}}
 +J_{\mathrm{degen}}^{\dag}-J_{\mathrm{degen}}\bigr)(s,\varphi,\omega'),
\end{equation}
where
\begin{align*}
J_{\mathrm{sing}}(s,\varphi,\omega')
&:=\int_{\mathbb{A}_{F}^{\times}}
 W_{\varphi^{\iota}}\bigl(\diag(z,1,1)\widetilde{w}\bigr)
 \omega\omega'(z)|z|^{2s}\,d^{\times}z,\\
J_{\mathrm{dual}}(s,\varphi,\omega')
&:=\frac{1}{2\pi i\rho_{F}}\sum_{\xi\in\widehat{C_{F}}}
 \int_{s_{0}-i\infty}^{s_{0}+i\infty}\iint_{(\mathbb{A}_{F}^{\times})^{2}}
 W_{\varphi^{\iota}}\bigl(\diag(z,1,1)n_{21}(y)\widetilde{w}\bigr)\\
&\qquad \overline{\xi}\omega\omega'(z)|z|^{2s-\lambda-\frac{1}{2}}
 \xi(y)|y|^{\lambda+\frac{1}{2}}\,d^{\times}z\,d^{\times}y\,d\lambda,\\
J_{\mathrm{degen}}^{\dag}(s,\varphi,\omega')
&:=\int_{\mathbb{A}_{F}^{\times}}
 W_{\varphi}^{\mathrm{degen},\dag}\bigl(w_{2}n_{23}(-1)\diag(z,z,1)\bigr)
 \omega'(z)|z|^{2s}\,d^{\times}z,\\
J_{\mathrm{degen}}(s,\varphi,\omega')
&:=\int_{\mathbb{A}_{F}^{\times}}\int_{\mathbb{A}_{F}}
 W_{\varphi}^{\mathrm{degen}}\bigl(w_{1}n_{12}(c')\diag(z,z,1)\bigr)
 \overline{\psi(c')}\omega'(z)|z|^{2s}\,dc'\,d^{\times}z.
\end{align*}
\end{thm}

\begin{proof}
The rapid decay of Whittaker functions and the assumption $\Re(s)\ggg 1$
imply that the function
\begin{align*}
g'\mapsto \varphi'_{s}(g'):=|\det g'|^{s}
 \int_{F^{\times}\backslash\mathbb{A}_{F}^{\times}}
 \mathcal{F}_{\mathrm{gen}}\varphi\bigl(\diag(zg',1)\bigr)
 \omega'(z)|z|^{2s}\,d^{\times}z
\end{align*}
is square-integrable modulo the center and has central character
$\overline{\omega'}$. We apply \eqref{3.2} with that character and the basis
$\{\overline{\phi}:\ \phi\in\mathfrak{B}(\sigma)\}$, where
$\omega_{\sigma}=\omega'$.  For such a $\phi$, we have
$\langle \varphi'_{s},\overline{\phi}\rangle=0$ if
$\phi=\xi\circ\det$ for some Hecke character $\xi$; if $\phi$ is a cusp
form or an Eisenstein series, unfolding gives
\begin{align*}
\langle \varphi'_{s},\overline{\phi}\rangle=\int_{[G']}\mathcal{F}_{\mathrm{gen}}\varphi\left(\diag(g',1)\right)\phi(g')|\det g'|^{s}dg'=\Psi(1/2+s,W_{\varphi},W_{\phi}).
\end{align*} 

Consider the function
\begin{align*}
J(s,\varphi,\omega'):=\int_{F\backslash\mathbb{A}_{F}}
 \int_{F^{\times}\backslash\mathbb{A}_{F}^{\times}}
 \mathcal{F}_{\mathrm{gen}}\varphi\bigl(n_{12}(c)\diag(z,z,1)\bigr)
 \omega'(z)|z|^{2s}\overline{\psi(c)}\,d^{\times}z\,dc.
\end{align*}
The additive quotient is compact, and the generic Whittaker expansion
decreases rapidly in $z$. Thus the integral converges absolutely for $\Re(s)\ggg1$.

Rapid Whittaker decay also bounds the invariant derivatives of
$\varphi'_{s}$ in $L^{2}$. Sobolev estimates give convergence of the
spectral expansion uniformly on compact subsets.
Termwise integration over $F\backslash\mathbb{A}_{F}$ gives
\begin{equation}\label{3.3}
J(s,\varphi,\omega')=J_{\mathrm{spec}}(s,\varphi,\omega').
\end{equation}

Integrate Proposition~\ref{prop2.2} against
$\omega'(z)|z|^{2s}d^{\times} z$ over
$F^{\times}\backslash\mathbb{A}_{F}^{\times}$.
By \eqref{f2.9} and absolute convergence, the generic term unfolds
to a sum over $\beta\in F$. The degenerate contribution is
$J_{\mathrm{degen}}^{\dag}(s,\varphi,\omega')
-J_{\mathrm{degen}}(s,\varphi,\omega')$.
The contragredient relation gives
\begin{equation}\label{3.6}
 W_{\varphi}\bigl(n_{32}(\beta)\diag(z,z,1)\bigr)
 =\omega(z)W_{\varphi^{\iota}}\bigl(\diag(z,1,1)n_{21}(-\beta z)\widetilde{w}\bigr).
\end{equation}
After $\beta\mapsto-\beta$, splitting this sum into
$\beta=0$ and $\beta\neq0$ gives
$J_{\mathrm{sing}}(s,\varphi,\omega')$ and $J_{s_{0},s}(1)$,
respectively, with $J$ as in \eqref{cf3.9}.

The bounds in Lemma~\ref{lem3.3} remain valid after any fixed number
of derivatives in the Archimedean variables, when the convergence
parameters are sufficiently large. The function $J$ is smooth and invariant
under an open subgroup of the finite idele units. Integration by parts
therefore gives absolute convergence in \eqref{eq:gl1-plancherel}.
In the Mellin coefficient of $J$, insert \eqref{cf3.9} and unfold the
$\beta$-sum against the integral over $F^{\times}\backslash\mathbb{A}_{F}^{\times}$.
After $y\mapsto z^{-1}y$, shift the Mellin variable by writing the
old $\lambda$ as $\lambda-s_{0}$, so the new line is $\Re(\lambda)=s_{0}$.
Taking $x=1$ in \eqref{eq:gl1-plancherel} gives
$J_{s_{0},s}(1)=J_{\mathrm{dual}}(s,\varphi,\omega')$.

For $s_{0},\Re(s)-s_{0}\ggg1$, the preceding argument gives absolute
convergence of $J_{\mathrm{dual}}$, and Lemma~\ref{lem7.1} gives that
of $J_{\mathrm{sing}}$. For $J_{\mathrm{degen}}$, use
Propositions~\ref{prop5.1.}, \ref{prop5.6.}, \ref{prop5.8}, and~\ref{prop6.3}.
For $J_{\mathrm{degen}}^{\dag}$, use
Propositions~\ref{prop5.5}, \ref{prop5.8.}, \ref{prop5.12}, and~\ref{prop6.4}.
\end{proof}

\subsubsection{Unitary data}
For cuspidal $\pi$, both degenerate terms in \eqref{e3.3} vanish.
For minimal or maximal Eisenstein data, write $\pi_{\mathbf{u}}$ and $\varphi_{\mathbf{u}}$ for the representations and Eisenstein series of \textsection\ref{sec2}, with $(\mathbf{u},d)=(\boldsymbol{\nu},3)$ or $(\boldsymbol{\mu},2)$, respectively. Put $\pi=\pi_{\mathbf{0}}$ and $\varphi=\varphi_{\mathbf{0}}$.
For either family, let $F(\mathbf{u})$ be meromorphic near $\mathbf{u}=0$, with no poles on the contours below. Define its constant term by
\begin{equation}\label{eq:inducing-constant-term}
\operatorname{CT}_{\mathbf{u}=0}F(\mathbf{u})
:=\frac{1}{(2\pi i)^{d}}\int_{|u_{1}|=\varepsilon}\cdots
\int_{|u_{d}|=d\varepsilon}F(\mathbf{u})
\frac{du_{1}\cdots du_{d}}{u_{1}\cdots u_{d}},
\end{equation}
where $0<\varepsilon<10^{-3}$ is sufficiently small. For the Mellin integrals below, take this constant term for $\Re(s)\ggg1$ and then continue meromorphically in $s$.
Take $F(\mathbf{u})$ to be the meromorphic continuation of $J_{\mathrm{degen}}(s,\varphi_{\mathbf{u}},\omega')$ or $J_{\mathrm{degen}}^{\dag}(s,\varphi_{\mathbf{u}},\omega')$. Denote the corresponding constant terms by $J_{\mathrm{degen}}^{\heartsuit}(s,\varphi,\omega')$ and $J_{\mathrm{degen}}^{\dag,\heartsuit}(s,\varphi,\omega')$, respectively.
In the minimal case, use Propositions~\ref{prop5.1.},
\ref{prop5.6.}, \ref{prop5.8} for $J_{\mathrm{degen}}$ and
Propositions~\ref{prop5.5}, \ref{prop5.8.}, \ref{prop5.12} for
$J_{\mathrm{degen}}^{\dag}$; in the maximal case use
Propositions~\ref{prop6.3} and~\ref{prop6.4}.
For $s_{0}$ and $\Re(s)-s_{0}$ sufficiently large, the spectral,
singular and dual terms are holomorphic in $\mathbf{u}$ near zero,
by Proposition~\ref{prop4.1}, Lemma~\ref{lem7.1} and the locally
uniform bounds of Lemma~\ref{lem3.3}. Their constant terms are
their values at zero. Applying \eqref{eq:inducing-constant-term} to
Theorem~\ref{thm3.3.} gives \eqref{e3.3}, with
$J_{\mathrm{degen}}^{\heartsuit}$ and
$J_{\mathrm{degen}}^{\dag,\heartsuit}$ as the degenerate terms.

\subsection{Spectral Reciprocity Formula of Type \RNum{2}}
Let $\eta$ be a unitary Hecke character, and let
$\mathbf{s}=(s_{1},s_{2})$. For $\Re(s_{1}),\Re(s_{2})\ggg1$, define
\begin{equation}\label{3.25}
I_{\mathrm{spec}}(\mathbf{s},\varphi,\omega',\eta):=I_{\mathrm{cusp}}(\mathbf{s},\varphi,\omega',\eta)+I_{\mathrm{Eis}}(\mathbf{s},\varphi,\omega',\eta),
\end{equation}
where
\begin{align*}
I_{\mathrm{cusp}}(\mathbf{s},\varphi,\omega',\eta):=\sum_{\substack{\pi'\in \mathcal{A}_{0}([G'],\omega')\\
\phi\in\mathfrak{B}(\pi')}}\Psi(1/2+s_{1},W_{\varphi},W_{\phi})\overline{\Psi(1/2+\overline{s_{2}},W_{\phi},\overline{\eta})},
\end{align*}
and  
\begin{multline*}
I_{\mathrm{Eis}}(\mathbf{s},\varphi,\omega',\eta):=\sum_{\xi\in \widehat{C_{F}}} \frac{1}{4\pi i}\int_{i\mathbb{R}}\sum_{h\in \mathfrak{B}(\xi,\overline{\xi}\omega')} \Psi\left(1/2+s_{1},W_{\varphi},W_{E(\cdot,h,\lambda)}\right)\\
\overline{\Psi(1/2+\overline{s_{2}},W_{E(\cdot,h,-\overline{\lambda})},\overline{\eta})}d\lambda.
\end{multline*} 

We extend this notation to norm twists of unitary characters by
\begin{equation}\label{eq:moment-character-twists}
I_{\mathrm{spec}}((s_{1},s_{2}),\varphi,
\omega_{0}'|\cdot|^{u},\eta_{0}|\cdot|^{v})
:=I_{\mathrm{spec}}((s_{1}+u/2,s_{2}+v-u/2),
\varphi,\omega_{0}',\eta_{0}),
\end{equation}
where $\omega_{0}'$ and $\eta_{0}$ are unitary. The same rule holds for
$I$, $I_{\mathrm{gen}}$, and $I_{\mathrm{degen}}$ below, since the
twists are absorbed into $|z|^{2s_{1}}|y|^{s_{1}+s_{2}}$.

\begin{thm}\label{thme3.3}
Suppose $\pi$ is convergent in the sense of Definition~\ref{defn3.2},
and let $\varphi\in\pi$.  Put
$\mathbf{s}=(s_{1},s_{2})$ and
$\mathbf{s}^{\vee}:=(\frac{s_{2}-s_{1}}{2},\frac{3s_{1}+s_{2}}{2})$.
If $\Re(s_{1})\ggg1$ and $\Re(s_{2}-s_{1})\ggg1$, with constants allowed to
depend on $\varphi$, then
\begin{multline}\label{e3.8}
I_{\mathrm{spec}}(\mathbf{s},\varphi,\omega',\eta)
=I_{\mathrm{spec}}(\mathbf{s}^{\vee},\pi(w_{2})\varphi,\omega^{-1}\overline{\omega}'\eta,\eta)\\
+I_{\mathrm{degen}}(\mathbf{s}^{\vee},\pi(w_{2})\varphi,\omega^{-1}\overline{\omega}'\eta,\eta)
-I_{\mathrm{gen}}(\mathbf{s}^{\vee},\pi(w_{2})\varphi,\omega^{-1}\overline{\omega}'\eta,\eta)\\
-I_{\mathrm{degen}}(\mathbf{s},\varphi,\omega',\eta)+I_{\mathrm{gen}}(\mathbf{s},\varphi,\omega',\eta),
\end{multline}
where 
\begin{align*}
I_{\mathrm{gen}}(\mathbf{s},\varphi,\omega',\eta):=\int_{\mathbb{A}_{F}^{\times}}\int_{\mathbb{A}_{F}^{\times}}W_{\varphi}\left(\diag(yz,z,1)\right)\omega'(z)|z|^{2s_{1}}\eta(y)|y|^{s_{1}+s_{2}}d^{\times}zd^{\times}y,
\end{align*}
and 
\begin{multline*}
I_{\mathrm{degen}}(\mathbf{s},\varphi,\omega',\eta):=\int_{\mathbb{A}_{F}^{\times}}\int_{\mathbb{A}_{F}^{\times}}\int_{\mathbb{A}_{F}}W_{\varphi}^{\mathrm{degen}}\left(w_{1}n_{12}(a)\diag(yz,z,1)\right)\\
 \overline{\psi(a)}da\omega'(z)|z|^{2s_{1}}d^{\times}z\eta(y)|y|^{s_{1}+s_{2}}d^{\times}y.
\end{multline*}
\end{thm}

\begin{proof}
To apply Proposition~\ref{prop2.5}, we first identify
$I_{\mathrm{spec}}$ with the Fourier integral
\begin{multline}\label{e3.13}
I(\mathbf{s},\varphi,\omega',\eta):=\int_{\mathbb{A}_{F}^{\times}}\int_{F^{\times}\backslash\mathbb{A}_{F}^{\times}}\int_{F\backslash\mathbb{A}_{F}}\mathcal{F}_{\mathrm{gen}}\varphi\left(n_{12}(c)\diag(zy,z,1)\right)\\
\overline{\psi(c)}dc\omega'(z)|z|^{2s_{1}}d^{\times}z\eta(y)|y|^{s_{1}+s_{2}}d^{\times}y.
\end{multline}

In \eqref{e3.13}, take the additive Fourier coefficient first.
Use the central Mellin integral $\varphi'_{s_{1}}$ of
\textsection\ref{sec3.3} with character $\omega'$ to write
\begin{equation}\label{e3.14}
I(\mathbf{s},\varphi,\omega',\eta)=\int_{\mathbb{A}_{F}^{\times}}\int_{F\backslash\mathbb{A}_{F}}\varphi'_{s_{1}}\left(n(c)a(y)\right)\overline{\psi(c)}dc\eta(y)|y|^{s_{2}}d^{\times}y.
\end{equation}

Use \eqref{3.2} with central character $\overline{\omega'}$ for $\varphi'_{s_{1}}$. The coefficients of $\overline{\phi}$ are
$\Psi(1/2+s_{1},W_{\varphi},W_{\phi})$, which decay rapidly by
Lemma~\ref{lem4.3}.
The bounds for $W_{\phi}$ used there give polynomial bounds at fixed finite level
for the absolute integrals defining $\Psi(1/2+\overline{s_{2}},W_{\phi},\overline{\eta})$ when $\Re(s_{2})$
is sufficiently large. Thus \eqref{e3.14} converges absolutely and
termwise integration is justified, locally uniformly for
$\Re(s_{1}),\Re(s_{2})\ggg1$. Unfolding the two periods gives
\begin{equation}\label{3.9}
I(\mathbf{s},\varphi,\omega',\eta)=I_{\mathrm{spec}}(\mathbf{s},\varphi,\omega',\eta).
\end{equation}

Use the identity
\begin{align*}
\pi(w_{2})\pi(\diag(zy,z,1))\varphi=\pi(\diag(zy,1,z))\pi(w_{2})\varphi=\omega(z)\pi(\diag(y,z^{-1},1))\pi(w_{2})\varphi.
\end{align*}

To write the dual terms of Proposition~\ref{prop2.5} as periods of
$\pi(w_{2})\varphi$, change $(z,y)\mapsto(z^{-1},zy)$. The character and
modulus factors become
\begin{align*}
\omega^{-1}\overline{\omega'}\eta(z)
|z|^{s_{2}-s_{1}}\eta(y)|y|^{s_{1}+s_{2}}.
\end{align*}
Since $2s_{1}^{\vee}=s_{2}-s_{1}$ and $s_{1}^{\vee}+s_{2}^{\vee}=s_{1}+s_{2}$,
the transformed integrals have parameters $\mathbf{s}^{\vee}$ in \eqref{e3.8}.
Moreover, $s_{2}^{\vee}-s_{1}^{\vee}=2s_{1}$, so the dual variables satisfy the same convergence conditions.
If $\omega=\omega_{0}|\cdot|^{u}$ with $\omega_{0}$ unitary, \eqref{eq:moment-character-twists} rewrites the dual spectral term with Mellin parameters $(s_{1}^{\vee}-u/2,s_{2}^{\vee}+u/2)$ and unitary central character $\omega_{0}^{-1}\overline{\omega'}\eta$. The convergence parameters are chosen large enough to allow this fixed shift.

Absolute convergence of the degenerate terms follows from
Propositions~\ref{prop5.1}, \ref{prop5.4}, \ref{prop5.6} and
\ref{prop6.2}, and that of the generic terms from
Proposition~\ref{prop5.4.}, applied to $\mathbf{s}$ and
$\mathbf{s}^{\vee}$ with the character shift
\eqref{eq:moment-character-twists}.
\end{proof}

\subsubsection{Unitary data}\label{sec3.4.2}
For unitary cuspidal $\pi$, \eqref{e3.8} holds for
$\Re(s_{1}),\Re(s_{2}-s_{1})\ggg1$, with both degenerate terms zero.

For Eisenstein data, use the minimal or maximal family of \textsection\ref{sec3.3}, and let $\omega_{\mathbf{u}}$ be the central character of $\pi_{\mathbf{u}}$. Propositions~\ref{prop5.1}, \ref{prop5.4}, \ref{prop5.6} in the minimal case, and Proposition~\ref{prop6.2} in the maximal case, continue both degenerate terms meromorphically. Take their constant terms at $\mathbf{u}=0$ as in \eqref{eq:inducing-constant-term}. On the dual side, take the constant term of
\begin{align*}
I_{\mathrm{degen}}^{\heartsuit}(\mathbf{s}^{\vee},
\pi_{\mathbf{u}}(w_{2})\varphi_{\mathbf{u}},
\omega_{\mathbf{u}}^{-1}\overline{\omega}'\eta,\eta).
\end{align*}
Absorb the norm twist of $\omega_{\mathbf{u}}^{-1}$ into the Mellin parameters using \eqref{eq:moment-character-twists}.

For $\Re(s_{1}),\Re(s_{2}-s_{1})\ggg1$, the proofs of Propositions~\ref{prop4.4} and~\ref{prop5.4.} show that the spectral and generic terms are holomorphic in $\mathbf{u}$ near zero. Taking constant terms therefore gives \eqref{e3.8} for both unitary Eisenstein families. For these unitary data, the generic integrals also converge in
\begin{align*}
\Re(s_{2})>\Re(s_{1})>0.
\end{align*}

\section{Meromorphic Continuation of the Spectral Side}\label{sec4}
We continue the spectral integrals of \textsection\ref{sec3} to a
neighborhood of the central parameters. The residues from shifting
the contours appear in the identities of \textsection\ref{sec8}.

\subsection{\texorpdfstring{Meromorphic Continuation of $J_{\mathrm{spec}}(s,\varphi,\omega')$}{Meromorphic Continuation of J\_spec}}\label{sec4.1}

\begin{defn}
Let $\pi$ be a unitary generic automorphic representation of $G(\mathbb{A}_{F})$ with central character $\omega$. 
\begin{itemize}
\item If $\pi$ is cuspidal, let
$\mathfrak{X}_{\pi,\omega'}^{\mathrm{dual}}
=\mathfrak{X}_{\pi,\omega'}^{+}=\mathfrak{X}_{\pi,\omega'}^{-}$ be the empty set.
\item If $\pi=\sigma\boxplus \chi$, where $\chi$ is a unitary Hecke character, and $\sigma$ is a unitary cuspidal representation of $G'(\mathbb{A}_{F})$ with central character $\overline{\chi}\omega$, we let $\mathfrak{X}_{\pi,\omega'}^{\mathrm{dual}}=\{\overline{\chi}\omega\omega'\}$, $\mathfrak{X}_{\pi,\omega'}^{+}:=\{\overline{\chi}\}$ and $\mathfrak{X}_{\pi,\omega'}^{-}:=\{\chi\omega'\}$.
\item If $\pi=\chi_{1}\boxplus \chi_{2}\boxplus\chi_{3}$, where $\chi_{1}$ and
$\chi_{2}$ are unitary Hecke characters and
$\chi_{3}:=\overline{\chi}_{1}\overline{\chi}_{2}\omega$, define
$\mathfrak{X}_{\pi,\omega'}^{\mathrm{dual}}
:=\{\overline{\chi}_{1}\omega\omega',\overline{\chi}_{2}\omega\omega',
\overline{\chi}_{3}\omega\omega'\}$,
$\mathfrak{X}_{\pi,\omega'}^{+}
:=\{\overline{\chi}_{1},\overline{\chi}_{2},\overline{\chi}_{3}\}$, and
$\mathfrak{X}_{\pi,\omega'}^{-}
:=\{\chi_{1}\omega',\chi_{2}\omega',\chi_{3}\omega'\}$.
\end{itemize}
\end{defn}

For the continuous spectrum put
\begin{align*}
F_{\xi}(s,\lambda):=\sum_{h\in\mathfrak{B}(\xi,\overline{\xi}\omega')}
\widetilde{\Psi}(1/2+s,W_{\varphi},W_{E(\cdot,h,\lambda)})
\overline{W_{E(\cdot,h,-\overline{\lambda})}}(I_{2}).
\end{align*}
For $\xi=\xi_{0}|\cdot|^{it}$ with $\xi_{0}\in\widehat{C_{F}}$,
norm twists are absorbed into the parameter:
\begin{align*}
 F_{\xi_{0}|\cdot|^{it}}(s,\lambda)
 :=F_{\xi_{0}}(s,\lambda+it).
\end{align*}
We use this convention also for $\widetilde{\mathcal{G}}$ and
$\widetilde{\mathcal{J}}$ below. Thus the characters in the residue
sets need not belong to the chosen extension of $\widehat{C_{F}}$.

\begin{prop}\label{prop4.1}
Let $\pi$ be unitary generic automorphic on $[G]$ with central
character $\omega$, and let $\omega'$ be a unitary Hecke character.
The integral $J_{\mathrm{spec}}(s,\varphi,\omega')$ converges absolutely
for $\Re(s)>1/2$ and extends meromorphically to $\Re(s)>-1/2$.
For $|\Re(s)|<1/2$, its continuation satisfies
\begin{align*}
J_{\mathrm{spec}}^{\heartsuit}(s,\varphi,\omega')=J_{\mathrm{cusp}}^{\heartsuit}(s,\varphi,\omega')+J_{\mathrm{Eis}}^{\heartsuit}(s,\varphi,\omega')+R_{\RNum{1}}^{+}(s,\varphi,\omega')-R_{\RNum{1}}^{-}(s,\varphi,\omega'),
\end{align*}
where $J_{\mathrm{cusp}}^{\heartsuit}$ is the cuspidal sum in
\eqref{s3.3} with $\Psi$ replaced by $\widetilde{\Psi}$, and
\begin{align*}
&J_{\mathrm{Eis}}^{\heartsuit}(s,\varphi,\omega')
:=\frac{1}{4\pi i}\sum_{\xi}\int_{i\mathbb{R}}F_{\xi}(s,\lambda)\,d\lambda,\\
&R_{\RNum{1}}^{+}(s,\varphi,\omega')
:=\frac{1}{2}\sum_{\xi\in\mathfrak{X}_{\pi,\omega'}^{+}}
\underset{\lambda=1/2-s}{\Res}F_{\xi}(s,\lambda),\\
&R_{\RNum{1}}^{-}(s,\varphi,\omega')
:=\frac{1}{2}\sum_{\xi\in\mathfrak{X}_{\pi,\omega'}^{-}}
\underset{\lambda=s-1/2}{\Res}F_{\xi}(s,\lambda).
\end{align*}
\end{prop}

Proposition~\ref{prop4.1} follows from Corollary~\ref{cor4.4} and Lemma~\ref{lem4.2} below.

\subsubsection{Estimates of $\widetilde{\Psi}(1/2+s,W_{\varphi},W_{\phi})$}
Let $\sigma\in \mathcal{A}_{0}([G'],\omega')$ or $\sigma=\xi|\cdot|^{\lambda}\boxplus \overline{\xi}\omega'|\cdot|^{-\lambda}$ with $\lambda\in i\mathbb{R}$. In the cuspidal case take $\phi\in\mathfrak{B}(\sigma)$; in the Eisenstein case take $\phi=E(\cdot,h,\lambda)$ with $h\in\mathfrak{B}(\xi,\overline{\xi}\omega')$.

We choose the bases to consist of eigenvectors of a positive self-adjoint elliptic operator $\Delta_{\infty}$, formed from the group and maximal-compact Casimir operators with $1$ added, and invariant under $g\mapsto g^{\iota}$. Write
\begin{align*}
\Delta_{\infty}\phi=\lambda_{\phi,\infty}\phi,
\qquad \lambda_{\phi,\infty}\geq1.
\end{align*}
By the $K_{\fin}$-finiteness of $\varphi$, choose a compact open subgroup $K_{\fin}^{*}\subseteq\prod_{v<\infty}K_{v}'$ fixing $\varphi$ under the embedding $G'\hookrightarrow G$. Choose the bases to respect the $K_{\fin}^{*}$-fixed subspaces. Hence,
\begin{equation}\label{4.1.}
\widetilde{\Psi}(1/2+s,W_{\varphi},W_{\phi})=\widetilde{\Psi}(1/2+s,W_{\varphi},W_{\phi})\mathbf{1}_{\phi\in \sigma^{K_{\fin}^{*}}}.
\end{equation}

\begin{lemma}\label{lem4.3}
Let $\Re(s)\ggg1$ and $m\in\mathbb{Z}_{\geq1}$. For a cuspidal
or Eisenstein basis vector $\phi$ with central character $\omega'$,
\begin{equation}\label{f4.1}
\big|\widetilde{\Psi}(1/2+s,W_{\varphi},W_{\phi})\big|
\ll \lambda_{\phi,\infty}^{-m}\mathbf{1}_{\phi\in\sigma^{K_{\fin}^{*}}},
\end{equation}
where the implied constant depends on $m$, $\varphi$, $s$, and $\omega'$, locally uniformly in $s$.
\end{lemma}

\begin{proof}
For $\Re(s)\ggg1$, Whittaker decay at infinity and power bounds at zero justify transferring powers of $\Delta_{\infty}$ to $W_{\varphi}(\diag(g',1))|\det g'|^{s}$ by integration by parts. The bounds for $W_{\phi}$ in \cite[\textsection3.2, \textsection3.4]{Mag18}, with the finite level and $\omega'$ fixed, grow at most as a fixed power of $\lambda_{\phi,\infty}$. If the bound for the absolute integral loses a fixed power
$\lambda_{\phi,\infty}^{b}$, transfer at least $m+b$ powers of
$\Delta_{\infty}$ and apply \eqref{4.1.} to obtain \eqref{f4.1}. On each vertical line with sufficiently large real part, the implied constant grows at most polynomially in $|\Im(s)|$.
\end{proof}

\subsubsection{Absolute convergence}
Consider the set 
\begin{align*}
\boldsymbol{\xi}(\pi,\omega'):=
\Big\{\left.\xi\right|_{\mathbb{A}_{F}^{(1)}}:\ 
\xi\in\mathfrak{X}_{\pi,\omega'}^{+}\cup
\mathfrak{X}_{\pi,\omega'}^{-}\Big\},
\end{align*}
so $|\boldsymbol{\xi}(\pi,\omega')|\leq6$. Let
\begin{align*}
J_{\mathrm{spec}}^{\mathrm{reg},\heartsuit}(s,\varphi,\omega')
:=J_{\mathrm{cusp}}^{\heartsuit}(s,\varphi,\omega')
+\frac{1}{4\pi i}\sum_{\xi\notin\boldsymbol{\xi}(\pi,\omega')}
\int_{i\mathbb{R}}F_{\xi}(s,\lambda)\,d\lambda.
\end{align*}

\begin{cor}\label{cor4.4}
The integral $J_{\mathrm{spec}}^{\mathrm{reg},\heartsuit}(s,\varphi,\omega')$ converges absolutely away from the poles of its summands and defines a meromorphic function on $\mathbb{C}$.
\end{cor}

\begin{proof}
Following \cite[\textsection7]{Yan25}, use
\cite[\textsection3.2, \textsection3.4]{Mag18} to obtain
\begin{equation}\label{eq4.9}
|W_{\phi}(I_{2})|\ll\lambda_{\phi,\infty}^{O_{F}(1)}.
\end{equation}
Lemma~\ref{lem4.3}, \eqref{eq4.9}, and Weyl's law at fixed level give absolute convergence for $\Re(s)$ sufficiently large. By the global functional equation,
\begin{align*}
\big|\widetilde{\Psi}(1/2+s,W_{\varphi},W_{\phi})\big|
=\big|\widetilde{\Psi}(1/2-s,\widetilde{W}_{\varphi},\widetilde{W}_{\phi})\big|,
\end{align*}
the spectral expansion also converges absolutely for $\Re(s)$ sufficiently negative.

The excluded inducing characters account for all moving Eisenstein poles. The remaining poles involve only finitely many cuspidal representations. After removing these poles by a fixed polynomial, the maximum principle on vertical strips, using \cite[Theorem~2.1(ii)]{Jac09}, gives, for every $m\geq1$,
\begin{equation}\label{eq4.10}
\big|\widetilde{\Psi}(1/2+s,W_{\varphi},W_{\phi})\big|
\ll\lambda_{\phi,\infty}^{-m}\mathbf{1}_{\phi\in\sigma^{K_{\fin}^{*}}},
\end{equation}
locally uniformly away from the poles. Together with \eqref{eq4.9} and Weyl's law, this proves the assertion.
\end{proof}

\subsubsection{Choice of contours}
The function $F_{\xi}(s,\lambda)$ contains the reciprocal Hecke
$L$-factors
\begin{align*}
 \Lambda(1+2\lambda,\xi^{2}\overline{\omega'})^{-1},
 \qquad
 \Lambda(1-2\lambda,\overline{\xi}^{2}\omega')^{-1}.
\end{align*}
The classical zero-free region for Hecke $L$-functions and the finiteness of
$\boldsymbol{\xi}(\pi,\omega')$ give a common constant $c>0$ such that these
reciprocals have no poles in
\begin{equation}\label{eq3.12}
 |\Re(\lambda)|<\frac{c}{10\log(3+|\Im(\lambda)|^{2})},
 \qquad \xi\in\boldsymbol{\xi}(\pi,\omega').
\end{equation}
Choose $c$ sufficiently small that $|\Re(\lambda)|<10^{-1}$
throughout this region. Put
\begin{align*}
 \mathcal{C}^{\pm}=\left\{\lambda=\pm
 \frac{c}{10\log(3+t^{2})}+it:\ t\in\mathbb{R}\right\},
\end{align*}
with both contours oriented upward.

\subsubsection{Meromorphic continuation of $J_{\mathrm{Eis}}^{\mathrm{sing}}(s,\varphi,\omega')$}\label{sec4.1.4}

For $\Re(s)\ggg1$, define
\begin{align*}
J_{\mathrm{Eis}}^{\mathrm{sing}}(s,\varphi,\omega')
=(4\pi i)^{-1}\sum_{\xi\in\boldsymbol{\xi}(\pi,\omega')}
\int_{i\mathbb{R}}F_{\xi}(s,\lambda)\,d\lambda.
\end{align*}

\begin{lemma}\label{lem4.2}
The series $J_{\mathrm{Eis}}^{\mathrm{sing}}(s,\varphi,\omega')$
converges absolutely for $\Re(s)\ggg1$ and continues
meromorphically to $\Re(s)>-1/2$. For $|\Re(s)|<1/2$,
\begin{equation}\label{e4.3}
J_{\mathrm{Eis}}^{\mathrm{sing},\heartsuit}(s,\varphi,\omega')
=\frac{1}{4\pi i}\sum_{\xi\in\boldsymbol{\xi}(\pi,\omega')}
\int_{i\mathbb{R}}F_{\xi}(s,\lambda)\,d\lambda 
+R_{\RNum{1}}^{+}(s,\varphi,\omega')-R_{\RNum{1}}^{-}(s,\varphi,\omega').
\end{equation}
In this strip, the integral over $i\mathbb{R}$ converges absolutely
and locally uniformly away from the poles.
\end{lemma}

\begin{proof}
We shift contours as in \cite[\textsection 9]{Yan25}.
As $\Re(s)$ decreases through $1/2$, the poles crossing $i\mathbb{R}$ are at $\lambda=s-1/2$ for $\xi\in\mathfrak{X}_{\pi,\omega'}^{-}$ and at $\lambda=1/2-s$ for $\xi\in\mathfrak{X}_{\pi,\omega'}^{+}$.
By \eqref{eq3.12}, the reciprocal $L$-factors have no poles between $\mathcal{C}^{-}$ and $\mathcal{C}^{+}$. To estimate $W_{E(\cdot,h,\lambda)}$ and $W_{E(\cdot,h,-\overline{\lambda})}$, use the Eisenstein functional equation for the parameter with negative real part. Their ramified and Archimedean local integrals converge absolutely,
with bounds $\ll (1+|\Re(\lambda)|^{-1})^{A}$ for a fixed $A$.
On $\mathcal{C}^{\pm}$, $|\Re(\lambda)|^{-1}\ll\log(3+|\Im(\lambda)|^{2})$. Use the unramified local formulas and Stirling's formula on
$\mathcal{C}^{\pm}$, and integrate by parts as in the bound \eqref{eq4.10}
for the Rankin--Selberg pairing. This gives rapid decay of
$F_{\xi}(s,\lambda)$ in $\Im(\lambda)$ locally uniformly in $s$ away
from the poles.
Near any fixed point of $\Re(s)=1/2$, choose a sufficiently small
neighborhood $U$ so that all crossing poles lie between $\mathcal{C}^{-}$ and $\mathcal{C}^{+}$.
For $s\in U$ with $\Re(s)>1/2$, shifting the contour to $\mathcal{C}^{+}$, oriented upwards, gives
\begin{equation}\label{4.3}
J_{\mathrm{Eis}}^{\mathrm{sing}}(s,\varphi,\omega')
=\frac{1}{4\pi i}\sum_{\xi\in\boldsymbol{\xi}(\pi,\omega')}
\int_{\mathcal{C}^{+}}F_{\xi}(s,\lambda)\,d\lambda
-R_{\RNum{1}}^{-}(s,\varphi,\omega').
\end{equation}
The integral over $\mathcal{C}^{+}$ is holomorphic on $U$.
For $s\in U$ with $\Re(s)<1/2$, moving back to $i\mathbb{R}$ adds
$R_{\RNum{1}}^{+}$ and gives \eqref{e4.3}. The integral over $i\mathbb{R}$
for $\Re(s)>1/2$, \eqref{4.3} near $\Re(s)=1/2$, and \eqref{e4.3} for
$|\Re(s)|<1/2$ together give the continuation to $\Re(s)>-1/2$.
At a multiple pole, take the residue of $F_{\xi}$.
\end{proof}

\subsection{\texorpdfstring{Meromorphic Continuation of $I_{\mathrm{spec}}(\mathbf{s},\varphi,\omega',\eta)$}{Meromorphic Continuation of I\_spec}}\label{sec4.2}
Let $\pi$ be a unitary generic automorphic representation of $G(\mathbb{A}_{F})$ with central character $\omega$ and $\varphi\in \pi$. Let $\omega'$ and $\eta$ be unitary Hecke characters. We define   
\begin{multline*}
\widetilde{\mathcal{G}}(\mathbf{s},\lambda,\xi;\varphi,\omega',\eta):=\sum_{h\in \mathfrak{B}(\xi,\overline{\xi}\omega')}\widetilde{\Psi}\left(1/2+s_{1},W_{\varphi},W_{E(\cdot,h,\lambda)}\right)\\
\overline{\widetilde{\Psi}(1/2+\overline{s_{2}},W_{E(\cdot,h,-\overline{\lambda})},\overline{\eta})}.
\end{multline*} 

By the decay of Whittaker functions, $\widetilde{\mathcal{G}}(\mathbf{s},\lambda,\xi;\varphi,\omega',\eta)$ converges absolutely when $\Re(s_{1}\pm\lambda)\ggg 1$ and $\Re(s_{2}\pm\lambda)\ggg 1$. Moreover, $\widetilde{\mathcal{G}}(\mathbf{s},\lambda,\xi;\varphi,\omega',\eta)$ admits a meromorphic continuation to $(s_{1},s_{2},\lambda)\in \mathbb{C}^{3}$, satisfying
\begin{multline}\label{e4.6}
\widetilde{\mathcal{G}}(\mathbf{s},\lambda,\xi;\varphi,\omega',\eta)\propto \Lambda(1+2\lambda,\xi^{2}\omega'^{-1})^{-1}\Lambda(1/2+s_{1}+\lambda,\pi\times\xi)\\
\Lambda(1-2\lambda,\overline{\xi}^{2}\omega')^{-1}\Lambda(1/2+s_{1}-\lambda,\pi\times\overline{\xi}\omega')\Lambda(1/2+s_{2}-\lambda,\overline{\xi}\eta)\Lambda(1/2+s_{2}+\lambda,\xi\overline{\omega}'\eta).
\end{multline}
The omitted factor in \eqref{e4.6} is holomorphic near the contours below.

For both shifts, choose $\mathcal{C}^{\pm}$ by decreasing $c$ in
\eqref{eq3.12} if necessary, so that the reciprocal $L$-factors in
\eqref{e4.6} have no poles in that strip for any
\begin{align*}
 \xi\in \mathfrak{X}_{\pi,\omega'}^{+}\cup
 \mathfrak{X}_{\pi,\omega'}^{-}\cup
 \{\eta,\overline{\eta}\omega'\}.
\end{align*}
Impose the same condition on the restrictions of these characters to $C_{F}$.
For the second shift, restrict $s_{1}$ to a compact subset of
$|\Re(s_{1})|<1/2$ and narrow $\mathcal{C}^{\pm}$ so that
$|\Re(\lambda)|<1/2-|\Re(s_{1})|$ uniformly in $s_{1}$.
This avoids the poles involving $s_{1}$.

Define the following residue terms:
\begin{align*}
&R_{\RNum{2}}^{1,+}(\mathbf{s},\varphi,\omega',\eta):=\frac{1}{2}\sum_{\xi\in \mathfrak{X}_{\pi,\omega'}^{+}}\underset{\lambda=1/2-s_{1}}{\Res}\widetilde{\mathcal{G}}(\mathbf{s},\lambda,\xi;\varphi,\omega',\eta),\\
&R_{\RNum{2}}^{1,-}(\mathbf{s},\varphi,\omega',\eta):=\frac{1}{2}\sum_{\xi\in \mathfrak{X}_{\pi,\omega'}^{-}}\underset{\lambda=s_{1}-1/2}{\Res}\widetilde{\mathcal{G}}(\mathbf{s},\lambda,\xi;\varphi,\omega',\eta),\\
&R_{\RNum{2}}^{2,+}(\mathbf{s},\varphi,\omega',\eta):=\frac{1}{2}\underset{\lambda=1/2-s_{2}}{\Res}\widetilde{\mathcal{G}}(\mathbf{s},\lambda,\overline{\eta}\omega';\varphi,\omega',\eta),\\
&R_{\RNum{2}}^{2,-}(\mathbf{s},\varphi,\omega',\eta):=\frac{1}{2}\underset{\lambda=s_{2}-1/2}{\Res}\widetilde{\mathcal{G}}(\mathbf{s},\lambda,\eta;\varphi,\omega',\eta),
\end{align*}
Define the residue terms first for distinct poles and continue their signed sum meromorphically. In the central strips, coincident poles in the same integral have the same sign and are counted once as a multiple pole.

\begin{prop}\label{prop4.4}
Let $\pi$ be unitary generic automorphic on $[G]$ with central
character $\omega$. The integral
$I_{\mathrm{spec}}(\mathbf{s},\varphi,\omega',\eta)$ converges absolutely
for $\Re(s_{1}),\Re(s_{2})>1/2$ and extends meromorphically to
$\Re(s_{1}),\Re(s_{2})>-1/2$. For $|\Re(s_{1})|,|\Re(s_{2})|<1/2$,
its continuation satisfies
\begin{multline}\label{4.7}
I_{\mathrm{spec}}^{\heartsuit}(\mathbf{s},\varphi,\omega',\eta)=I_{\mathrm{cusp}}^{\heartsuit}(\mathbf{s},\varphi,\omega',\eta)+I_{\mathrm{Eis}}^{\heartsuit}(\mathbf{s},\varphi,\omega',\eta)+R_{\RNum{2}}^{1,+}(\mathbf{s},\varphi,\omega',\eta)\\
-R_{\RNum{2}}^{1,-}(\mathbf{s},\varphi,\omega',\eta)+R_{\RNum{2}}^{2,+}(\mathbf{s},\varphi,\omega',\eta)-R_{\RNum{2}}^{2,-}(\mathbf{s},\varphi,\omega',\eta),
\end{multline}
where 
\begin{align*}
&I_{\mathrm{cusp}}^{\heartsuit}(\mathbf{s},\varphi,\omega',\eta):=\sum_{\pi'\in \mathcal{A}_{0}([G'],\omega')}\sum_{\phi\in\mathfrak{B}(\pi')}\widetilde{\Psi}(1/2+s_{1},W_{\varphi},W_{\phi})\overline{\widetilde{\Psi}(1/2+\overline{s_{2}},W_{\phi},\overline{\eta})},\\
&I_{\mathrm{Eis}}^{\heartsuit}(\mathbf{s},\varphi,\omega',\eta):=\sum_{\xi\in \widehat{C_{F}}} \frac{1}{4\pi i}\int_{i\mathbb{R}} \widetilde{\mathcal{G}}(\mathbf{s},\lambda,\xi;\varphi,\omega',\eta)d\lambda.
\end{align*}
\end{prop}

\begin{proof}
For the contributing $K_{\fin}^{*}$-fixed vectors, the Whittaker bounds
in the proof of Lemma~\ref{lem4.3} give
$|\widetilde{\Psi}(1/2+\overline{s_{2}},W_{\phi},\overline{\eta})|
\ll(1+\lambda_{\phi,\infty})^{A}$ on a fixed line $\Re(s_{2})>10$. The
$\mathrm{GL}_{2}\times\mathrm{GL}_{1}$ functional equation gives a polynomial
bound on a left boundary line. The vertical-strip bounds
\cite[Theorem~2.1(ii)]{Jac09} allow the maximum principle after clearing
the poles, giving polynomial growth in $\lambda_{\phi,\infty}$ locally
uniformly in $s_{2}$ away from the poles. On $\mathcal{C}^{\pm}$ the same
argument uses the Eisenstein
Whittaker bounds in the proof of Lemma~\ref{lem4.2}. The rapid decay of
the first period in \eqref{eq4.10} and in that proof absorbs this growth.
Thus the sum and integral converge locally uniformly away from the poles,
and the horizontal integrals tend to zero.

For $\Re(s_{2})>10$, we move the $\lambda$-contour of
$I_{\mathrm{Eis}}^{\heartsuit}$ in \eqref{4.7} to $\mathcal{C}^{+}$ as in \eqref{4.3},
then back to
$i\mathbb{R}$ by \eqref{e4.3} at $s=s_{1}$. Neither of the two factors at $1/2+s_{2}\pm\lambda$
in \eqref{e4.6} has a pole in this shift. The moment continues to $\Re(s_{1})>-1/2$, with
$R_{\RNum{2}}^{1,+}-R_{\RNum{2}}^{1,-}$ for $|\Re(s_{1})|<1/2$.
Next keep $s_{1}$ in a compact subset of $|\Re(s_{1})|<1/2$ and continue in $s_{2}$. The cuspidal sum continues termwise. For the Eisenstein integral, shift to $\mathcal{C}^{+}$ as $\Re(s_{2})$ decreases through $1/2$, then return to $i\mathbb{R}$. Our choice of contours avoids the poles involving $s_{1}$.
The poles $(\xi,\lambda)=(\eta,s_{2}-1/2)$ and $(\overline{\eta}\omega',1/2-s_{2})$ contribute $-R_{\RNum{2}}^{2,-}$ and $R_{\RNum{2}}^{2,+}$, respectively, with the factor $1/2$ coming from the measure $d\lambda/(4\pi i)$.
The other possible poles of the last two factors in \eqref{e4.6} lie at $\lambda=s_{2}+1/2$ and $-s_{2}-1/2$. Their real parts have opposite signs for $\Re(s_{2})>-1/2$,
and neither pole crosses $i\mathbb{R}$.
The terms $R_{\RNum{2}}^{1,+}-R_{\RNum{2}}^{1,-}$ continue meromorphically in $s_{2}$, proving \eqref{4.7}.

Reversing the order of the shifts, with $s_{2}$ in a compact subset of $|\Re(s_{2})|<1/2$, covers the remaining region. Near any point of $\Re(s_{j})=1/2$, choose fixed contours avoiding the moving poles and include the residues crossed. These formulas agree on overlaps and give continuation to $\Re(s_{1}),\Re(s_{2})>-1/2$.
\end{proof}

\section{Minimal Eisenstein Data: Degenerate Terms}\label{sec5}
We compute the degenerate terms for minimal Eisenstein data using
Godement sections and local integral formulas.

For unitary Hecke characters $\chi_{1},\chi_{2},\chi_{3}$ and $\boldsymbol{\nu}=(\nu_{1},\nu_{2},\nu_{3})\in\mathcal{R}_{\mathrm{min}}$ as in \eqref{f2.2}, put $\pi=\chi_{1}\boxplus\chi_{2}\boxplus\chi_{3}$ and $\pi_{\boldsymbol{\nu}}=\chi_{1}|\cdot|^{\nu_{1}}\boxplus\chi_{2}|\cdot|^{\nu_{2}}\boxplus\chi_{3}|\cdot|^{\nu_{3}}$. Write $\varphi_{\boldsymbol{\nu}}=E(\cdot,\boldsymbol{\chi},\boldsymbol{\nu})\in\pi_{\boldsymbol{\nu}}$ for the Eisenstein series \eqref{f2.5}.

\subsection{Godement Sections}\label{sec5.1}
We work with Godement sections of the following form.  General sections
are treated by finite linear combinations and meromorphic scalar
multiples; see \cite[\textsection 7.1]{Yan23}:
\begin{multline}\label{e5.1}
f(g; \boldsymbol{\chi}, \boldsymbol{\nu})=\chi_{1}(\det g)|\det g|^{1+\nu_{1}}\int_{G'(\mathbb{A}_{F})}\Phi_{2\times 3}[(\mathbf{0},g')g]\\
h(g'^{-1})\chi_{1}(\det g')|\det g'|^{\frac{3}{2}+\nu_{1}}dg',
\end{multline}
where $\Phi_{2\times 3}$ is a Schwartz--Bruhat function on
$M_{2\times 3}(\mathbb{A}_{F})$, and $h(\cdot)$ is a section in
$\chi_{2}|\cdot|^{\nu_{2}}\boxplus\chi_{3}|\cdot|^{\nu_{3}}$. More precisely,
\begin{equation}\label{e5.2}
h(g'):=\chi_{2}(\det g')|\det g'|^{\frac{1}{2}+\nu_{2}}\int_{\mathbb{A}_{F}^{\times}}\Phi_{1\times 2}((0,t)g')\chi_{2}\overline{\chi}_{3}(t)|t|^{1+\nu_{2}-\nu_{3}}d^{\times}t,
\end{equation}
where $\Phi_{1\times 2}$ is a Schwartz--Bruhat function on
$M_{1\times 2}(\mathbb{A}_{F})$.

Since $\boldsymbol{\nu}\in\mathcal{R}_{\mathrm{min}}$, the integrals in \eqref{e5.1} and \eqref{e5.2} converge absolutely. Substituting \eqref{e5.2} into \eqref{e5.1} yields 
\begin{multline}\label{e5.3}
f(g; \boldsymbol{\chi}, \boldsymbol{\nu})=\chi_{1}(\det g)|\det g|^{1+\nu_{1}}\int_{G'(\mathbb{A}_{F})}\Phi_{2\times 3}[(\mathbf{0},g')g]\chi_{1}\overline{\chi}_{2}(\det g')\\
|\det g'|^{1+\nu_{1}-\nu_{2}}\int_{\mathbb{A}_{F}^{\times}}\Phi_{1\times 2}((0,t)g'^{-1})\chi_{2}\overline{\chi}_{3}(t)|t|^{1+\nu_{2}-\nu_{3}}d^{\times}tdg'.
\end{multline}

Tate's thesis and Iwasawa decomposition applied to \eqref{e5.3} give a
meromorphic continuation of $f(g;\boldsymbol{\chi},\boldsymbol{\nu})$ to
$\boldsymbol{\nu}\in\mathbb{C}^{3}$, with
\begin{align*}
f(g; \boldsymbol{\chi}, \boldsymbol{\nu})\propto \Lambda(1+\nu_{1}-\nu_{2},\chi_{1}\overline{\chi}_{2})\Lambda(1+\nu_{1}-\nu_{3},\chi_{1}\overline{\chi}_{3})\Lambda(1+\nu_{2}-\nu_{3},\chi_{2}\overline{\chi}_{3}).
\end{align*}

\begin{defn}\label{def:exceptional-places-minimal}
Let $S$ be the finite set such that $v\notin S$ precisely when
$v<\infty$,
$v\nmid\mathfrak{D}_{F}$, the section
$f_{v}(\cdot;\boldsymbol{\chi},\boldsymbol{\nu})$ is spherical,
$\Phi_{2\times 3,v}=\mathbf{1}_{(\mathcal{O}_{v})^{6}}$,
$\Phi_{1\times 2,v}=\mathbf{1}_{(\mathcal{O}_{v})^{2}}$, and the characters
$\omega_{v}'$ and $\eta_{v}$ are unramified.
\end{defn}

\subsection{Degenerate Whittaker Functions}
We use the degenerate Whittaker functions of \textsection\ref{sec2.2}. For simplicity, we abbreviate $f(g; \boldsymbol{\chi}, \boldsymbol{\nu})$ by $f_{\boldsymbol{\nu}}(g)$.

\subsubsection{The Bruhat decomposition}
Let $S_{3}:=\{I_{3}, w_{1}, w_{2}, w_{1}w_{2}, w_{2}w_{1}, w_{1}w_{2}w_{1}\}$ be the Weyl group of $G$ relative to $B$. For $w\in S_{3}$, let $N_{w}:=(N\cap w^{-1}Nw)\backslash N$. We have the Bruhat decomposition 
\begin{equation}\label{fc2.5}
G(F)=\bigsqcup_{w\in S_{3}}B(F)wN_{w}(F).
\end{equation}

The Bruhat decomposition \eqref{fc2.5} gives 
\begin{equation}\label{2.12}
W_{\varphi_{\boldsymbol{\nu}}}^{\mathrm{degen}}(g)=\sum_{w\in S_{3}}W_{\varphi_{\boldsymbol{\nu}}}^{\mathrm{degen}}(g;w),
\end{equation}
where 
\begin{align*}
W_{\varphi_{\boldsymbol{\nu}}}^{\mathrm{degen}}(g;w):=\int_{(F\backslash\mathbb{A}_{F})^{3}}\sum_{u\in N_{w}(F)}f_{\boldsymbol{\nu}}\left(wu\begin{pmatrix}
1& a& b\\
&1& c\\
&& 1
\end{pmatrix}g\right)\overline{\psi(c)}dadbdc.
\end{align*}

By \eqref{2.9} and the orthogonality relation \eqref{e1.1}, we have
\begin{align*}
W_{\varphi_{\boldsymbol{\nu}}}^{\mathrm{degen}}(g;w)\equiv 0,\ \ w\in \{I_{3}, w_{1}, w_{2}w_{1}\}.
\end{align*}
Thus \eqref{2.12} reduces to
\begin{equation}\label{2.13}
W_{\varphi_{\boldsymbol{\nu}}}^{\mathrm{degen}}(g)=W_{\varphi_{\boldsymbol{\nu}}}^{\mathrm{degen}}(g;w_{2})+W_{\varphi_{\boldsymbol{\nu}}}^{\mathrm{degen}}(g;w_{1}w_{2})+W_{\varphi_{\boldsymbol{\nu}}}^{\mathrm{degen}}(g;w_{1}w_{2}w_{1}).
\end{equation}

For $W_{\varphi_{\boldsymbol{\nu}}}^{\mathrm{degen},\dag}$,
\eqref{fc2.5} gives
\begin{equation}\label{f5.7}
W_{\varphi_{\boldsymbol{\nu}}}^{\mathrm{degen},\dag}(g)=\sum_{w\in S_{3}}W_{\varphi_{\boldsymbol{\nu}}}^{\mathrm{degen},\dag}(g;w),	
\end{equation}
where
\begin{align*}
W_{\varphi_{\boldsymbol{\nu}}}^{\mathrm{degen},\dag}(g;w):=\int_{(F\backslash\mathbb{A}_{F})^{3}}\sum_{u\in N_{w}(F)}f_{\boldsymbol{\nu}}\left(wu\begin{pmatrix}
1& a& b\\
&1& c\\
&& 1
\end{pmatrix}g\right)\overline{\psi(a)}dadbdc.
\end{align*}

Orthogonality likewise gives
$W_{\varphi_{\boldsymbol{\nu}}}^{\mathrm{degen},\dag}(g;w)\equiv0$ for
$w\in\{I_{3},w_{2},w_{1}w_{2}\}$. Hence \eqref{f5.7} becomes
\begin{equation}\label{f2.13}
W_{\varphi_{\boldsymbol{\nu}}}^{\mathrm{degen},\dag}(g)=W_{\varphi_{\boldsymbol{\nu}}}^{\mathrm{degen},\dag}(g;w_{1})+W_{\varphi_{\boldsymbol{\nu}}}^{\mathrm{degen},\dag}(g;w_{2}w_{1})+W_{\varphi_{\boldsymbol{\nu}}}^{\mathrm{degen},\dag}(g;w_{1}w_{2}w_{1}).
\end{equation}

Unfolding gives
\begin{align*}
&W_{\varphi_{\boldsymbol{\nu}}}^{\mathrm{degen}}(g;w_{2})=\int_{\mathbb{A}_{F}}f_{\boldsymbol{\nu}}\left(w_{2}n_{23}(c)g\right)\overline{\psi(c)}dc,\\
&W_{\varphi_{\boldsymbol{\nu}}}^{\mathrm{degen},\dag}(g;w_{1})=\int_{\mathbb{A}_{F}}f_{\boldsymbol{\nu}}\left(w_{1}n_{12}(a)g\right)\overline{\psi(a)}da,\\
&W_{\varphi_{\boldsymbol{\nu}}}^{\mathrm{degen}}(g;w_{1}w_{2})=\int_{\mathbb{A}_{F}}\int_{\mathbb{A}_{F}}f_{\boldsymbol{\nu}}\left(w_{1}w_{2}n_{13}(b)n_{23}(c)g\right)\overline{\psi(c)}dbdc,\\
&W_{\varphi_{\boldsymbol{\nu}}}^{\mathrm{degen},\dag}(g;w_{2}w_{1})=\int_{\mathbb{A}_{F}}\int_{\mathbb{A}_{F}}f_{\boldsymbol{\nu}}\left(w_{2}w_{1}n_{12}(a)n_{13}(b)g\right)\overline{\psi(a)}dbda,\\
&W_{\varphi_{\boldsymbol{\nu}}}^{\mathrm{degen}}(g;w_{1}w_{2}w_{1})=\int_{\mathbb{A}_{F}}\int_{\mathbb{A}_{F}}\int_{\mathbb{A}_{F}}f_{\boldsymbol{\nu}}\left(w_{1}w_{2}w_{1}\begin{pmatrix}
1& a& b\\
&1& c\\
&& 1
\end{pmatrix}g\right)\overline{\psi(c)}dadbdc,\\
&W_{\varphi_{\boldsymbol{\nu}}}^{\mathrm{degen},\dag}(g;w_{1}w_{2}w_{1})=\int_{\mathbb{A}_{F}}\int_{\mathbb{A}_{F}}\int_{\mathbb{A}_{F}}f_{\boldsymbol{\nu}}\left(w_{1}w_{2}w_{1}\begin{pmatrix}
1& a& b\\
&1& c\\
&& 1
\end{pmatrix}g\right)\overline{\psi(a)}dadbdc.
\end{align*}

\subsubsection{Decomposition of the degenerate terms}\label{sec5.2.2}
Let $w\in S_{3}$. Define 
\begin{align*}
J_{\mathrm{degen}}^{\dag}(s,\varphi_{\boldsymbol{\nu}},\omega';w):=\int_{\mathbb{A}_{F}^{\times}}W_{\varphi_{\boldsymbol{\nu}}}^{\mathrm{degen},\dag}\left(w_{2}n_{23}(-1)\diag(z,z,1);w\right)\omega'(z)|z|^{2s}d^{\times}z,
\end{align*}
\begin{multline*}
J_{\mathrm{degen}}(s,\varphi_{\boldsymbol{\nu}},\omega';w):=\int_{\mathbb{A}_{F}^{\times}}\int_{\mathbb{A}_{F}}W_{\varphi_{\boldsymbol{\nu}}}^{\mathrm{degen}}\left(w_{1}n_{12}(a)\diag(z,z,1);w\right)\\
\overline{\psi(a)}da\omega'(z)|z|^{2s}d^{\times}z,
\end{multline*}
\begin{multline*}
I_{\mathrm{degen}}(\mathbf{s},\varphi_{\boldsymbol{\nu}},\omega',\eta;w):=\int_{\mathbb{A}_{F}^{\times}}\int_{\mathbb{A}_{F}^{\times}}\int_{\mathbb{A}_{F}}W_{\varphi_{\boldsymbol{\nu}}}^{\mathrm{degen}}\left(w_{1}n_{12}(a)\diag(yz,z,1);w\right)\\
\overline{\psi(a)}da\omega'(z)|z|^{2s_{1}}d^{\times}z\eta(y)|y|^{s_{1}+s_{2}}d^{\times}y.
\end{multline*}

As a result of \eqref{2.13} and \eqref{f2.13}, we have
\begin{align*}
&J_{\mathrm{degen}}^{\dag}(s,\varphi_{\boldsymbol{\nu}},\omega')=\sum_{w\in\{w_{1},w_{2}w_{1},w_{1}w_{2}w_{1}\}}J_{\mathrm{degen}}^{\dag}(s,\varphi_{\boldsymbol{\nu}},\omega';w),\\
&J_{\mathrm{degen}}(s,\varphi_{\boldsymbol{\nu}},\omega')=\sum_{w\in\{w_{2},w_{1}w_{2},w_{1}w_{2}w_{1}\}}J_{\mathrm{degen}}(s,\varphi_{\boldsymbol{\nu}},\omega';w),\\
&I_{\mathrm{degen}}(\mathbf{s},\varphi_{\boldsymbol{\nu}},\omega',\eta)=\sum_{w\in\{w_{2},w_{1}w_{2},w_{1}w_{2}w_{1}\}}I_{\mathrm{degen}}(\mathbf{s},\varphi_{\boldsymbol{\nu}},\omega',\eta;w).
\end{align*}

\subsection{A Local Convolution of Kirillov Vectors}
Let $v\in\Sigma_{F}$, and let $\sigma_{1,v}$ and $\sigma_{2,v}$ be generic
representations of $G'(F_{v})$ with central characters $\omega_{1,v}$ and
$\omega_{2,v}$, respectively, and let $W_{1,v}$ and $W_{2,v}$ be
Whittaker functions. Let $s\in\mathbb{C}$ and define
\begin{equation}\label{f3.5}
\mathcal{R}_{v}(s,\sigma_{1,v},\sigma_{2,v}):=\int_{F_{v}^{\times}}W_{1,v}\left(a(y_{v})\right)W_{2,v}\left(a(y_{v})\right)|y_{v}|_{v}^{s}d^{\times}y_{v}.
\end{equation}

\begin{lemma}\label{lem3.2}
Suppose $v<\infty$, $v\nmid\mathfrak{D}_{F}$, and $W_{1,v}$ and
$W_{2,v}$ are spherical. Then
\begin{equation}\label{eq3.6}
\mathcal{R}_{v}(s,\sigma_{1,v},\sigma_{2,v})=\frac{W_{1,v}(I_{2})W_{2,v}(I_{2})L_{v}(1+s,\sigma_{1,v}\times\sigma_{2,v})}{L_{v}(2+2s,\omega_{1,v}\omega_{2,v})}.
\end{equation}

Suppose $\sigma_{1,v}=\chi_{1,v}|\cdot|^{\nu_{1}}\boxplus \chi_{2,v}|\cdot|^{\nu_{2}}$ and $\sigma_{2,v}=\chi_{3,v}|\cdot|^{\nu_{3}}\boxplus \chi_{4,v}|\cdot|^{\nu_{4}}$, where for $1\leq j\leq 4$, $\chi_{j,v}$ is a unitary character of $F_{v}^{\times}$ and $\nu_{j}\in \mathbb{C}$. The function $\mathcal{R}_{v}(s,\sigma_{1,v},\sigma_{2,v})$ converges absolutely in 
\begin{equation}\label{f5.12}
\Re(s)>-1+\max\{|\Re(\nu_{1})|,|\Re(\nu_{2})|\}+\max\{|\Re(\nu_{3})|,|\Re(\nu_{4})|\},
\end{equation}
and admits a meromorphic continuation to $(s,\nu_{1},\nu_{2},\nu_{3},\nu_{4})\in \mathbb{C}^{5}$. 

Suppose $\sigma_{1,v}=\chi_{1,v}|\cdot|^{\nu_{1}}\boxplus \chi_{2,v}|\cdot|^{\nu_{2}}$, and $\sigma_{2,v}$ is a unitary supercuspidal representation or Steinberg twist at a finite place, or a unitary discrete series at a real place, where for $1\leq j\leq 2$, $\chi_{j,v}$ is a unitary character of $F_{v}^{\times}$ and $\nu_{j}\in \mathbb{C}$. The function $\mathcal{R}_{v}(s,\sigma_{1,v},\sigma_{2,v})$ converges absolutely in 
\begin{equation}\label{equ5.13}
\Re(s)>-1+\max\{|\Re(\nu_{1})|,|\Re(\nu_{2})|\},
\end{equation}
and admits a meromorphic continuation to $(s,\nu_{1},\nu_{2})\in \mathbb{C}^{3}$.
\end{lemma}

\begin{proof}
The standard Whittaker bounds at zero and infinity give absolute
convergence in \eqref{f5.12} and \eqref{equ5.13}.
Casselman--Shalika and the Cauchy identity, with $\Vol(\mathcal{O}_{v}^{\times},d^{\times} y_{v})=1$, give \eqref{eq3.6}.

For meromorphic continuation, take $s_{0}$ and $\Re(s)-s_{0}$ sufficiently large. Parseval's identity on $F_{v}^{\times}$ gives
\begin{multline}\label{eq5.41}
\mathcal{R}_{v}(s,\sigma_{1,v},\sigma_{2,v})=c_{v}\sum_{\xi_{v}}\int_{(s_{0})}\Psi_{v}(1/2+\lambda,W_{1,v},\xi_{v})\\
\Psi_{v}(1/2+s-\lambda,W_{2,v},\overline{\xi}_{v})d\lambda,
\end{multline}
where $c_{v}$ depends only on the local measure normalization and $\xi_{v}$ ranges over all characters of $\mathcal{O}_{v}^{\times}$, $\{\pm1\}$, or $S^{1}$ at finite, real, or complex places, respectively. Extend them trivially on $\varpi_{v}^{\mathbb{Z}}$ at finite places and on $\mathbb{R}_{>0}$ at Archimedean places. The contour $(s_{0})$ is $s_{0}+i\mathbb{R}$ at Archimedean places and the segment
\begin{align*}
 [s_{0}-\pi i/\log q_{v},\ s_{0}+\pi i/\log q_{v}]
\end{align*}
at finite places.

The local periods in \eqref{eq5.41} continue meromorphically by
\textsection\ref{sec3.2}; the character sum is finite at finite places.
At infinity, absolute convergence and the local functional equations
bound the two edges of a vertical strip. The maximum principle and
Stirling's formula give polynomial bounds on the strip away from the
poles, also after fixed derivatives of the vectors.
Integration by parts then gives rapid decay in the character parameter
and $\Im(\lambda)$. Moving the contour and including the crossed
residues proves meromorphic continuation.
\end{proof}

\subsection{Integral Representation of Whittaker Functions}
Let $\Phi$ be a Schwartz--Bruhat function on
$\mathbb{A}_{F}\times\mathbb{A}_{F}$. Define the partial Fourier transform by
\begin{align*}
\mathcal{F}_{2}\Phi(t_{1},t_{2}):=\int_{\mathbb{A}_{F}}\Phi(t_{1},b')\psi(b't_{2})db',\ \ t_{1}, t_{2}\in \mathbb{A}_{F}.
\end{align*}

\begin{lemma}\label{lem5.1}
Let $\xi_{1}$ and $\xi_{2}$ be Hecke characters of $F^{\times}\backslash\mathbb{A}_{F}^{\times}$. For $y\in\mathbb{A}_{F}^{\times}$, the functions
\begin{align*}
&W_{1}\left(a(y)\right)=\xi_{1}(y)|y|^{\frac{1}{2}}\int_{\mathbb{A}_{F}}\int_{\mathbb{A}_{F}^{\times}}\Phi(yt,b't)\xi_{1}\xi_{2}^{-1}(t)|t|d^{\times}t\overline{\psi(b')}db',\\
&W_{2}\left(a(y)\right)=\xi_{2}(y)|y|^{\frac{1}{2}}\int_{\mathbb{A}_{F}}\int_{\mathbb{A}_{F}^{\times}}\Phi(t,b't)\xi_{1}\xi_{2}^{-1}(t)|t|d^{\times}t\overline{\psi(b'y)}db',\\
&W_{3}\left(a(y)\right)=\xi_{2}(y)|y|^{\frac{1}{2}}\int_{\mathbb{A}_{F}^{\times}}\mathcal{F}_{2}\Phi(t,-t^{-1}y)\xi_{1}\xi_{2}^{-1}(t)d^{\times}t
\end{align*}
are equal and define a vector in the Kirillov model of
$\xi_{1}\boxplus\xi_{2}$. The first two expressions converge absolutely
if $|\xi_{1}\xi_{2}^{-1}(t)|=|t|^{r}$ with $r>1$, and are then continued
in the inducing parameters; the third is absolutely convergent for
every fixed $y$.
\end{lemma}

\begin{proof}
The Godement--Whittaker formula \cite[\textsection 2.1.2]{Yan25b}, after a norm twist, gives $W_{1}$. The substitutions $t\mapsto y^{-1}t$ and $b'\mapsto yb'$ give $W_{2}$; the change $b'\mapsto t^{-1}u$ gives $W_{3}$, with Jacobian $|t|^{-1}$ and character $\psi(-uy/t)$.
For fixed $y\neq0$, the arguments $t$ and $-y/t$ give rapid decay at
both ends at Archimedean places and compact support in $F_{v}^{\times}$
at finite places. For $v\notin S$, the standard local factor
$\mathbf{1}_{\mathcal{O}_{v}^{2}}$ and $y\in\mathcal{O}_{v}^{\times}$ give
the integral over $\mathcal{O}_{v}^{\times}$ of the constant $1$.
This gives locally uniform convergence in the inducing parameters.
\end{proof}

For the Godement sections used in \textsection\ref{sec5}--\textsection\ref{sec6}, the integrals below converge locally uniformly in all parameters in the stated regions and are holomorphic there. We write $\mathcal{P}_{v}^{\heartsuit}$ for the meromorphic continuation of $\mathcal{P}_{v}^{\sharp}$, and $\mathcal{P}_{v}^{\dag,\heartsuit}$ for the holomorphic continuation of $\mathcal{P}_{v}^{\dag,\sharp}$, in all parameters. For $v\notin S$,
\begin{equation}\label{f5.15}
 \mathcal{P}_{v}^{\heartsuit}=1,\qquad
 \mathcal{P}_{v}^{\dag,\heartsuit}=1.
\end{equation}

\subsection{\texorpdfstring{The $w_{2}$ Term of $I_{\mathrm{degen}}$}{The w2 Term of I-degen}}\label{sec5.5}
Define 
\begin{multline*}
\mathcal{L}(\mathbf{s},\boldsymbol{\nu},\pi,\omega',\eta;w_{2}):=\Lambda(2+2\nu_{1}+2s_{1},\chi_{1}^{2}\omega')^{-1}\Lambda(1+\nu_{1}+\nu_{3}+2s_{1},\chi_{1}\chi_{3}\omega')\\
\Lambda(1+\nu_{1}+\nu_{2}+2s_{1},\chi_{1}\chi_{2}\omega')\Lambda(1+\nu_{1}-\nu_{2},\chi_{1}\overline{\chi}_{2})
\Lambda(1+\nu_{1}-\nu_{3},\chi_{1}\overline{\chi}_{3})\\
\Lambda(s_{2}-s_{1}-\nu_{1},\overline{\chi}_{1}\overline{\omega}'\eta)\Lambda(1+\nu_{1}+s_{1}+s_{2},\chi_{1}\eta).
\end{multline*}

For each $v\in\Sigma_{F}$, let
$\mathcal{L}_{v}(\mathbf{s},\boldsymbol{\nu},\pi,\omega',\eta;w_{2})$ be the
local factor of $\mathcal{L}(\mathbf{s},\boldsymbol{\nu},\pi,\omega',\eta;w_{2})$. Define the normalized local factor $\mathcal{P}_{v}^{\sharp}(\mathbf{s},\varphi_{\boldsymbol{\nu}},\omega',\eta;w_{2})$ by 
\begin{multline}\label{eq:local-w2-degenerate}
 \frac{1}{\mathcal{L}_{v}(\mathbf{s},\boldsymbol{\nu},\pi,\omega',\eta;w_{2})}\int_{(F_{v}^{\times})^{2}}\int_{(F_{v})^{2}}
 f_{\boldsymbol{\nu},v}\left(
 \begin{pmatrix}1&&\\&1&\\a_{v}&c_{v}&1\end{pmatrix}w_{2}w_{1}\right)\\
\overline{\psi_{v}(a_{v}y_{v}+c_{v}z_{v})}
\chi_{1,v}\chi_{3,v}\omega_{v}'(z_{v})|z_{v}|_{v}^{1+\nu_{1}+\nu_{3}+2s_{1}}\\
\overline{\chi}_{1,v}\overline{\omega}_{v}'\eta_{v}(y_{v})|y_{v}|_{v}^{s_{2}-s_{1}-\nu_{1}}dc_{v}da_{v}d^{\times}z_{v}d^{\times}y_{v}.
\end{multline}

In \textsection\ref{sec5.5}--\textsection\ref{sec5.13}, the additive
integrals are evaluated first as Fourier integrals, initially in the
convergence region of the Godement sections and then by continuation.
The convergence assertions refer to the subsequent integrations over $\mathbb{A}_{F}^{\times}$.
\begin{prop}\label{prop5.1}
$I_{\mathrm{degen}}(\mathbf{s},\varphi_{\boldsymbol{\nu}},\omega',\eta;w_{2})$
converges absolutely in the region 
\begin{equation}\label{eq5.13}
\begin{cases}
\Re(\nu_{1}-\nu_{3})>0,\ \ \Re(\nu_{1}-\nu_{2})>0,\\
\Re(2s_{1}+\nu_{1}+\nu_{2})>0,\ \ \Re(2s_{1}+\nu_{1}+\nu_{3})>0,\\
\Re(s_{2}+s_{1}+\nu_{1})>0,\ \ \Re(s_{2}-s_{1}-\nu_{1})>1.
\end{cases}
\end{equation}
\end{prop}

\begin{proof}
In the unfolded integral, make the simultaneous changes
\begin{align*}
 a\mapsto ay,\qquad c\mapsto cz,\qquad z\mapsto z/y.
\end{align*}
The measure contributes $|yz|$. Moving $\diag(z/y,1,z)$
to the left gives the factors in \eqref{2.9}.

Insert \eqref{e5.3} with
$g'=k\diag(t_{1},t_{2})\left(\begin{smallmatrix}1&\\b'&1\end{smallmatrix}\right)$
and Haar measure \eqref{e1.5}. Keeping $t_{1},t_{2}$ fixed, make the simultaneous changes of variables
\begin{align*}
 (a,c,b',y,z,t)
 \longmapsto\left(t_{2}^{-1}a,\ t_{2}^{-1}(c+t_{1}b'),\ -t_{1}t_{2}^{-1}b',\
 t_{2}y,\ t_{2}t_{1}^{-1}z,\ t_{2}t\right).
\end{align*}
The additive Jacobian is $|t_{1}t_{2}^{-3}|$, and the character becomes
$\overline{\psi}(ay+ct_{1}^{-1}z+b'z)$. 
For $t_{1},t_{1}'\in\mathbb{A}_{F}$ and $k\in K'$, define
\begin{multline*}
h(t_{1},t_{1}';\mathbf{s},\boldsymbol{\nu},
k):=\int_{\mathbb{A}_{F}^{\times}}\int_{\mathbb{A}_{F}^{\times}}\int_{\mathbb{A}_{F}}\int_{\mathbb{A}_{F}}\Phi_{2\times 3}\left(k\begin{pmatrix}
& t_{1} &  \\
a & c & t_{2}
\end{pmatrix}w_{2}w_{1}\right)\overline{\psi(ct_{1}')}dc\\
\overline{\psi(ay)}da\overline{\chi}_{1}\overline{\omega}'\eta(y)|y|^{s_{2}-s_{1}-\nu_{1}}d^{\times}y
\chi_{1}\eta(t_{2})|t_{2}|^{1+\nu_{1}+s_{1}+s_{2}}d^{\times}t_{2}.
\end{multline*}

Taking the Fourier transforms in $a$ and $c$ first leaves a
Schwartz--Bruhat function of $(t_{1},y,t_{1}',t_{2})$. Its integrals
in $y$ and $t_{2}$ converge absolutely when
\begin{align*}
 \Re(s_{2}-s_{1}-\nu_{1})>1,\quad \Re(s_{2}+s_{1}+\nu_{1})>0
\end{align*}
and define a Schwartz--Bruhat function of $(t_{1},t_{1}')$, locally
uniformly in these parameters and $k$.  Tate's thesis gives a
meromorphic continuation of
$h(t_{1},t_{1}';\mathbf{s},\boldsymbol{\nu},k)$ to
$(\mathbf{s},\boldsymbol{\nu})\in\mathbb{C}^{5}$, with
\begin{align*}
h(t_{1},t_{1}';\mathbf{s},\boldsymbol{\nu},
k)\propto \Lambda(s_{2}-s_{1}-\nu_{1},\overline{\chi}_{1}\overline{\omega}'\eta)\Lambda(1+\nu_{1}+s_{1}+s_{2},\chi_{1}\eta).
\end{align*}

Moreover, at a place $v\notin S$, we have
\begin{multline}\label{equa5.23}
h_{v}(t_{1,v},t_{1,v}';\mathbf{s},\boldsymbol{\nu},
k_{v})=\mathbf{1}_{\mathcal{O}_{v}}(t_{1,v})\mathbf{1}_{\mathcal{O}_{v}}(t_{1,v}')L_{v}(s_{2}-s_{1}-\nu_{1},\overline{\chi}_{1}\overline{\omega}'\eta)\\
L_{v}(1+\nu_{1}+s_{1}+s_{2},\chi_{1}\eta).
\end{multline}

These changes of variables give
\begin{align*}
I_{\mathrm{degen}}(\mathbf{s},\varphi_{\boldsymbol{\nu}},\omega',\eta;w_{2})
=\int_{K'}\chi_{1}\overline{\chi}_{2}(\det k) 
\int_{\mathbb{A}_{F}^{\times}}W_{2}(a(z);k)W_{3}(a(z);k)\,d^{\times} z\,dk,
\end{align*}
where the function $W_{2}\left(a(z);k\right)$ is defined by 
\begin{equation}\label{fc5.22}
\chi_{3}(z)|z|^{\frac{1}{2}+\nu_{3}}\int_{\mathbb{A}_{F}}\int_{\mathbb{A}_{F}^{\times}}\Phi_{1\times 2}((tb',t)k^{-1})\chi_{2}\overline{\chi}_{3}(t)|t|^{1+\nu_{2}-\nu_{3}}d^{\times}t\overline{\psi(b'z)}db',
\end{equation}
and the function $W_{3}\left(a(z);k\right)$ is defined by 
\begin{equation}\label{fc5.23}
\chi_{1}\omega'(z)|z|^{\frac{1}{2}+\nu_{1}+2s_{1}}\int_{\mathbb{A}_{F}^{\times}}h(t_{1},t_{1}^{-1}z;\mathbf{s},\boldsymbol{\nu},
k)\overline{\chi}_{2}\overline{\chi}_{3}\overline{\omega}'(t_{1})|t_{1}|^{-\nu_{2}-\nu_{3}-2s_{1}}d^{\times}t_{1}.
\end{equation}

By Lemma \ref{lem5.1}, $W_{2}\left(a(z);k\right)$ is a vector in the Kirillov model of $\pi_{1}:=\chi_{3}|\cdot|^{\nu_{3}}\boxplus \chi_{2}|\cdot|^{\nu_{2}}$, and $W_{3}\left(a(z);k\right)$ is a vector in the Kirillov model of
$\pi_{2}:=\chi_{1}\omega'|\cdot|^{\nu_{1}+2s_{1}}\boxplus
\chi_{1}\overline{\chi}_{2}\overline{\chi}_{3}|\cdot|^{\nu_{1}-\nu_{2}-\nu_{3}}$.
For $W_{3}$, use Lemma~\ref{lem5.1} with the Schwartz--Bruhat
function whose second partial Fourier transform is
$h(t,-u;\mathbf{s},\boldsymbol{\nu},k)$.
Moreover, at a place $v\notin S$, we have
\begin{align*}
W_{j,v}\left(a(z_{v});k_{v}\right)=W_{j,v}\left(a(z_{v});I_{2}\right),\ \ j=2,3,
\end{align*}
which are spherical Kirillov functions normalized by $W_{2,v}(I_{2};I_{2})=1$ and 
\begin{equation}\label{cf5.24}
W_{3,v}(I_{2};I_{2})=h_{v}(1,1;\mathbf{s},\boldsymbol{\nu},
I_{2}).
\end{equation}
Lemma \ref{lem3.2} and \eqref{cf5.24} give
\begin{equation}\label{cf5.25}
I_{\mathrm{degen},v}(\mathbf{s},\varphi_{\boldsymbol{\nu}},\omega',\eta;w_{2})=\frac{h_{v}(1,1;\mathbf{s},\boldsymbol{\nu},
I_{2})L_{v}(1,\pi_{1,v}\times\pi_{2,v})}{L_{v}(2+2\nu_{1}+2s_{1},\chi_{1,v}^{2}\omega_{v}')}.
\end{equation}

Equations~\eqref{equa5.23} and \eqref{cf5.25} give \eqref{f5.15}
by multiplicativity of $L_{v}(1,\pi_{1,v}\times\pi_{2,v})$.
 
For $v\in S$, consider the local convolution of $W_{2,v}$ and
$W_{3,v}$.  At zero, the four possible power exponents in their
product are $1+\alpha$, where
\begin{equation}\label{eq:minimal-w2-kirillov-exponents}
 \alpha\in\{\nu_{1}+\nu_{3}+2s_{1},\ \nu_{1}+\nu_{2}+2s_{1},\
             \nu_{1}-\nu_{2},\ \nu_{1}-\nu_{3}\}.
\end{equation}
Repeated exponents may introduce logarithms.  Since
$\Re(\alpha)>0$ in \eqref{eq5.13}, each term is locally integrable
against $d^{\times} z_{v}$.  At infinity the Whittaker functions have
rapid decay at Archimedean places and bounded support at finite
places.  These bounds are uniform for $k_{v}\in K_{v}'$ and for parameters
in compact subsets of \eqref{eq5.13}; they prove absolute convergence
of the integrals in $z_{v}$. These integrals continue meromorphically
by Lemma~\ref{lem3.2}.

In \eqref{f5.15}, the product over $v\notin S$ converges absolutely
in \eqref{eq5.13}, and the local quotients in \eqref{f5.15}
at $v\in S$ continue meromorphically.
\end{proof}

\subsection{\texorpdfstring{The $w_{2}$ Term of $J_{\mathrm{degen}}$}{The w2 Term of J-degen}}
Define 
\begin{multline*}
\mathcal{L}(s,\boldsymbol{\nu},\pi,\omega';w_{2}):=\Lambda(2+2\nu_{1}+2s,\chi_{1}^{2}\omega')^{-1}\Lambda(1+\nu_{1}+\nu_{3}+2s,\chi_{1}\chi_{3}\omega')\\
\Lambda(1+\nu_{1}+\nu_{2}+2s,\chi_{1}\chi_{2}\omega')\Lambda(1+\nu_{1}-\nu_{2},\chi_{1}\overline{\chi}_{2})
\Lambda(1+\nu_{1}-\nu_{3},\chi_{1}\overline{\chi}_{3}).
\end{multline*}

Define the normalized local factor $\mathcal{P}_{v}^{\sharp}(s,\varphi_{\boldsymbol{\nu}},\omega';w_{2})$ by
\begin{multline*}
\mathcal{P}_{v}^{\sharp}(s,\varphi_{\boldsymbol{\nu}},\omega';w_{2}):=\mathcal{L}_{v}(s,\boldsymbol{\nu},\pi,\omega';w_{2})^{-1}\int_{F_{v}^{\times}}\int_{F_{v}}\int_{F_{v}}\overline{\psi_{v}(a_{v}+c_{v}z_{v})}\\
f_{\boldsymbol{\nu},v}\left(
 \begin{pmatrix}1&&\\&1&\\a_{v}&c_{v}&1\end{pmatrix}w_{2}w_{1}\right)dc_{v}da_{v}
\chi_{1,v}\chi_{3,v}\omega_{v}'(z_{v})|z_{v}|_{v}^{1+\nu_{1}+\nu_{3}+2s}d^{\times}z_{v}.
\end{multline*}

\begin{prop}\label{prop5.1.}
$J_{\mathrm{degen}}(s,\varphi_{\boldsymbol{\nu}},\omega';w_{2})$
converges absolutely in the region 
\begin{equation}\label{eq5.18-w2}
\begin{cases}
\Re(\nu_{1}-\nu_{3})>0,\ \ \Re(\nu_{1}-\nu_{2})>0,\\
\Re(2s+\nu_{1}+\nu_{2})>0,\ \ \Re(2s+\nu_{1}+\nu_{3})>0.
\end{cases}
\end{equation}
\end{prop}

\begin{proof}
In the unfolded $w_{2}$ term, make the change $c\mapsto cz$, with
Jacobian $|z|$. Move $\diag(z,1,z)$ to the left and use \eqref{2.9}.

In \eqref{e5.3} we use the coordinates and the substitution in
$(a,c,b',t)$ from the proof of Proposition~\ref{prop5.1}. At fixed $t_{1},t_{2}$
we replace $z$ by $t_{2}z$. The Jacobian is $|t_{1}t_{2}^{-3}|$ and the phase
is $\overline{\psi}(at_{2}^{-1}+cz+b't_{1}z)$.

For $t_{1},t_{1}'\in\mathbb{A}_{F}$, define
\begin{multline*}
h(t_{1},t_{1}';s,\boldsymbol{\nu},
k):=\int_{\mathbb{A}_{F}}\int_{\mathbb{A}_{F}^{\times}}\int_{\mathbb{A}_{F}}\Phi_{2\times 3}\left(k\begin{pmatrix}
& t_{1} &  \\
a & c & t_{2}
\end{pmatrix}w_{2}w_{1}\right)\overline{\psi(at_{2}^{-1})}da\\
\chi_{1}^{2}\omega'(t_{2})|t_{2}|^{1+2\nu_{1}+2s}d^{\times}t_{2}
\overline{\psi(ct_{1}')}dc.
\end{multline*}
Fourier transformation in $a$ and $c$ gives a Schwartz--Bruhat function
of $(t_{1},t_{1}')$, entire in $(s,\boldsymbol{\nu})$; the argument
$t_{2}^{-1}$ gives rapid decay as $t_{2}\to0$.

The change $z\mapsto t_{1}^{-1}z$ turns the last Fourier character into
$\overline{\psi}(b'z)$ and gives
\begin{align*}
J_{\mathrm{degen}}(s,\varphi_{\boldsymbol{\nu}},\omega';w_{2})
=\int_{K'}\chi_{1}\overline{\chi}_{2}(\det k)
\int_{\mathbb{A}_{F}^{\times}}W_{2}(a(z);k)W_{3}(a(z);k)\,d^{\times} z\,dk.
\end{align*}
Here $W_{2}$ is given by \eqref{fc5.22}, and $W_{3}$ by
\eqref{fc5.23} with $s_{1}=s$ and $h$ defined above.
The exponents \eqref{eq:minimal-w2-kirillov-exponents} at $s_{1}=s$
give local integrability at zero in \eqref{eq5.18-w2}.
Lemma~\ref{lem3.2} gives meromorphic continuation of the local
Rankin--Selberg integrals.
For $v\notin S$, the support requires $t_{2},t_{2}^{-1}\in\mathcal{O}_{v}$,
so the remaining $t_{2}$-integral is over $\mathcal{O}_{v}^{\times}$, of
volume $1$. Hence
$h_{v}(t_{1},t_{1}';s,\boldsymbol{\nu},k_{v})
=\mathbf{1}_{\mathcal{O}_{v}}(t_{1})\mathbf{1}_{\mathcal{O}_{v}}(t_{1}')$.  Thus $W_{2,v}(I_{2};k_{v})=W_{3,v}(I_{2};k_{v})=1$.
The spherical convolution gives the four numerator factors and the
denominator in $\mathcal{L}$, proving \eqref{f5.15}.
\end{proof}

\subsection{\texorpdfstring{The $w_{1}$ Term of $J_{\mathrm{degen}}^{\dag}$}{The w1 Term of J-degen dagger}}\label{sec5.7}
Define 
\begin{multline*}
\mathcal{L}^{\dag}(s,\boldsymbol{\nu},\pi,\omega';w_{1}):=\Lambda(1+\nu_{2}-\nu_{3},\chi_{2}\overline{\chi}_{3})\\
\Lambda(1+\nu_{1}-\nu_{3},\chi_{1}\overline{\chi}_{3})\Lambda(1+\nu_{1}+\nu_{2}+2s,\chi_{1}\chi_{2}\omega').
\end{multline*}

For each $v\in\Sigma_{F}$, let
$\mathcal{L}_{v}^{\dag}(s,\boldsymbol{\nu},\pi,\omega';w_{1})$ be the local
factor of $\mathcal{L}^{\dag}(s,\boldsymbol{\nu},\pi,\omega';w_{1})$. Define the normalized local factor $\mathcal{P}_{v}^{\dag,\sharp}(s,\varphi_{\boldsymbol{\nu}},\omega';w_{1})$ by 
\begin{multline*}
\mathcal{P}_{v}^{\dag,\sharp}(s,\varphi_{\boldsymbol{\nu}},\omega';w_{1}):=\mathcal{L}_{v}^{\dag}(s,\boldsymbol{\nu},\pi,\omega';w_{1})^{-1}\int_{F_{v}^{\times}}\int_{F_{v}}
f_{\boldsymbol{\nu},v}\left(
 \begin{pmatrix}1&&\\a_{v}&1&\\-z_{v}^{-1}&&1\end{pmatrix}w_{1}w_{2}\right)\\
\overline{\psi_{v}(a_{v}z_{v})}da_{v}\chi_{2,v}\chi_{3,v}\omega_{v}'(z_{v})|z_{v}|_{v}^{\nu_{2}+\nu_{3}+2s}d^{\times}z_{v}.
\end{multline*}

\begin{prop}\label{prop5.5}
The integral $J_{\mathrm{degen}}^{\dag}(s,\varphi_{\boldsymbol{\nu}},\omega';w_{1})$
converges absolutely, with the additive Fourier integral evaluated first,
in the region
\begin{equation}\label{eq5.18-dag-w1}
 \Re(\nu_{2}-\nu_{3})>0,\quad \Re(\nu_{1}-\nu_{3})>0,\quad \Re(2s+\nu_{1}+\nu_{2})>0.
\end{equation}
\end{prop}

\begin{proof}
In the unfolded $w_{1}$ term, make the change $a\mapsto az$, with
Jacobian $|z|$. Move $\diag(1,z,z)$ to the left and use \eqref{2.9}.
Insert \eqref{e5.3}, using
$g'=k\diag(t_{1},t_{2})n(b')$
and the Haar measure \eqref{e1.4}, and make the simultaneous changes
\begin{align*}
a\mapsto t_{1}^{-1}(a+b'z^{-1}),\qquad
b'\mapsto t_{1}^{-1}b',\qquad t\mapsto t_{2}t.
\end{align*}
The additive Jacobian is $|t_{1}|^{-2}$, and the character becomes
$\overline{\psi}(at_{1}^{-1}z+b't_{1}^{-1})$.

Here we use the partial Fourier transform
\begin{align*}
 h(r,x,u;t_{1},t_{1}',k)
 :=\int_{\mathbb{A}_{F}^{2}}
 \Phi_{2\times3}\left(k\begin{pmatrix}
 a&t_{1}&b'\\-u&&r
 \end{pmatrix}w_{1}w_{2}\right)
 \overline{\psi}(ax+b't_{1}')\,da\,db'.
\end{align*}
It is Schwartz--Bruhat in all five additive variables. Make the simultaneous changes $z\mapsto t_{1}x$, $t_{2}\mapsto t_{1}xu$ and define
\begin{align*}
H(x,u;s,\boldsymbol{\nu},k)
:=\int_{\mathbb{A}_{F}^{\times}}
h(t_{1}xu,x,u;t_{1},t_{1}^{-1},k)
\chi_{1}^{2}\omega'(t_{1})|t_{1}|^{1+2\nu_{1}+2s}\,d^{\times} t_{1}.
\end{align*}
The decay of $h$ in $t_{1},t_{1}^{-1},x,u$ makes $H$ Schwartz--Bruhat
in $(x,u)$ and entire in $(s,\boldsymbol{\nu})$, locally uniformly
in the parameters and $k$. We obtain
\begin{multline*}
J_{\mathrm{degen}}^{\dag}(s,\varphi_{\boldsymbol{\nu}},\omega';w_{1})
=\int_{K'}\chi_{1}\overline{\chi}_{2}(\det k)\int_{\mathbb{A}_{F}^{\times}}
\Phi_{1\times2}((0,t)k^{-1})\chi_{2}\overline{\chi}_{3}(t)
|t|^{1+\nu_{2}-\nu_{3}}\,d^{\times} t\\
\int_{(\mathbb{A}_{F}^{\times})^{2}}
H(x,u;s,\boldsymbol{\nu},k)
\chi_{1}\chi_{2}\omega'(x)|x|^{1+\nu_{1}+\nu_{2}+2s}
\chi_{1}\overline{\chi}_{3}(u)|u|^{1+\nu_{1}-\nu_{3}}
\,d^{\times} x\,d^{\times} u\,dk.
\end{multline*}
The three Tate integrals converge absolutely in \eqref{eq5.18-dag-w1}.
Normalizing these local Tate integrals gives an entire quotient by
$\mathcal{L}_{v}^{\dag}$, also after integration over $K_{v}'$.

For $v\notin S$, Fourier self-duality and the pair $t_{1},t_{1}^{-1}$
restrict the $t_{1}$-integral to $\mathcal{O}_{v}^{\times}$, of volume $1$,
and give
\begin{align*}
H_{v}(x,u;s,\boldsymbol{\nu},k_{v})
=\mathbf{1}_{\mathcal{O}_{v}}(x)\mathbf{1}_{\mathcal{O}_{v}}(u).
\end{align*}
The three local Tate integrals therefore equal their standard
$L$-factors, proving \eqref{f5.15}.
\end{proof}

\subsection{\texorpdfstring{The $w_{1}w_{2}$ Term of $I_{\mathrm{degen}}$}{The w1w2 Term of I-degen}}\label{sec5.8}
Define 
\begin{multline*}
\mathcal{L}(\mathbf{s},\boldsymbol{\nu},\pi,\omega',\eta;w_{1}w_{2}):=\Lambda(2+2\nu_{2}+2s_{1},\chi_{2}^{2}\omega')^{-1}\\
\Lambda(1+\nu_{2}+\nu_{3}+2s_{1},\chi_{2}\chi_{3}
\omega')\Lambda(1+\nu_{2}-\nu_{3},\chi_{2}\overline{\chi}_{3})\Lambda(1+\nu_{2}+s_{1}+s_{2},\chi_{2}\eta)\\\Lambda(s_{2}-s_{1}-\nu_{2},\overline{\chi}_{2}\overline{\omega}'\eta)
\Lambda(1+\nu_{1}+\nu_{2}+2s_{1},\chi_{1}\chi_{2}\omega')\Lambda(\nu_{1}-\nu_{2},\chi_{1}\overline{\chi}_{2}).
\end{multline*}

Define the normalized local factor
$\mathcal{P}_{v}^{\sharp}(\mathbf{s},\varphi_{\boldsymbol{\nu}},\omega',\eta;w_{1}w_{2})$ by
\begin{multline}\label{eq:local-w1w2-degenerate}
\mathcal{P}_{v}^{\sharp}(\mathbf{s},\varphi_{\boldsymbol{\nu}},\omega',\eta;w_{1}w_{2}):=\mathcal{L}_{v}(\mathbf{s},\boldsymbol{\nu},\pi,\omega',\eta;w_{1}w_{2})^{-1}\int_{F_{v}^{\times}}\int_{F_{v}^{\times}}\int_{F_{v}}\int_{F_{v}}\int_{F_{v}}\\
f_{\boldsymbol{\nu},v}\left(\begin{pmatrix}
1& & \\
b_{v}&1& \\
c_{v}& a_{v}& 1
\end{pmatrix}w_{1}w_{2}w_{1}\right)\overline{\psi_{v}(a_{v}y_{v}+c_{v}z_{v})}db_{v}da_{v}dc_{v}\\
\chi_{2,v}\chi_{3,v}
\omega_{v}'(z_{v})|z_{v}|_{v}^{1+\nu_{2}+\nu_{3}+2s_{1}}\overline{\chi}_{2,v}\overline{\omega}_{v}'\eta_{v}(y_{v})|y_{v}|_{v}^{s_{2}-s_{1}-\nu_{2}}d^{\times}z_{v}d^{\times}y_{v}.
\end{multline}

\begin{prop}\label{prop5.4}
$I_{\mathrm{degen}}(\mathbf{s},\varphi_{\boldsymbol{\nu}},\omega',\eta;w_{1}w_{2})$
converges absolutely in the region 
\begin{equation}\label{f5.19}
\begin{cases}
\Re(\nu_{1}-\nu_{2})>1,\ \ \Re(\nu_{2}-\nu_{3})>0,\\
\Re(2s_{1}+\nu_{1}+\nu_{2})>0,\ \ \Re(2s_{1}+\nu_{2}+\nu_{3})>0,\\
\Re(s_{1}+s_{2}+\nu_{1})>0,\ \ \Re(s_{2}-s_{1}-\nu_{2})>1
\end{cases}
\end{equation}
\end{prop}

\begin{proof}
In the unfolded integral, make the simultaneous changes
\begin{align*}
 a\mapsto ya,\qquad b\mapsto(z/y)b,\qquad c\mapsto zc,\qquad z\mapsto z/y.
\end{align*}
The additive Jacobian is $|z|^{2}$. Move $\diag(1,z/y,z)$ to the left
and use \eqref{2.9}. Insert \eqref{e5.3} in the Iwasawa coordinates
$g'=k\diag(t_{1},t_{2})\left(\begin{smallmatrix}1&\\b'&1\end{smallmatrix}\right)$,
and, with $t_{1},t_{2}$ fixed, make the simultaneous substitution

\begin{align*}
 (a,b,c,b',y,z,t)
 \longmapsto\left(\frac{a+t_{1}b'}{t_{2}},\frac{b}{t_{1}},
 \frac{c+bb'}{t_{2}},-\frac{t_{1}b'}{t_{2}},t_{2}y,t_{2}z,t_{2}t\right).
\end{align*}
The additive Jacobian is $|t_{2}|^{-3}$ and the character becomes
$\overline{\psi}(ay+cz+b't_{1}y+bb'z)$.

Define the auxiliary function
\begin{multline}\label{5.19}
h((t_{1},y),(z,z');k,\nu_{1},\mathbf{s}):=\int_{\mathbb{A}_{F}^{\times}}
\int_{(\mathbb{A}_{F})^{3}}\Phi_{2\times 3}\left(k\begin{pmatrix}
b & t_{1} &  \\
c & a & t_{2}
\end{pmatrix}w_{1}w_{2}w_{1}\right)\\
\overline{\psi(ay+cz)}dadc\overline{\psi(bz')}db\chi_{1}\eta(t_{2})|t_{2}|^{1+\nu_{1}+s_{1}+s_{2}}d^{\times}t_{2}. 
\end{multline}

By Tate's thesis,
$\Lambda(1+\nu_{1}+s_{1}+s_{2},\chi_{1}\eta)^{-1}
h((t_{1},y),(z,z');k,\nu_{1},\mathbf{s})$ is holomorphic in
$(\nu_{1},s_{1},s_{2})\in\mathbb{C}^{3}$ and is Schwartz--Bruhat in
$(t_{1},y,z,z')\in\mathbb{A}_{F}^{4}$.

Collect the $t_{2},a,b,c$ integrals in $h$ by \eqref{5.19}. Then
change $y\mapsto t_{1}^{-1}y$, followed by
$b'\mapsto y^{-1}b'$, $z\mapsto yz$ and
$t\mapsto ty$. The change in $b'$ contributes $|y|^{-1}$;
the multiplicative measures are unchanged.
The character is now $\overline{\psi}(b')$, the $y$-weight is
$\chi_{2}\eta(y)|y|^{1+\nu_{2}+s_{1}+s_{2}}$, and the $t_{1}$-weight is
$\chi_{1}\omega'\overline{\eta}(t_{1})|t_{1}|^{\nu_{1}+s_{1}-s_{2}}$.
Thus
\begin{multline}\label{c5.25}
I_{\mathrm{degen}}(\mathbf{s},\varphi_{\boldsymbol{\nu}},\omega',\eta;w_{1}w_{2})=\int_{K'}\chi_{1}\overline{\chi}_{2}(\det k)\int_{\mathbb{A}_{F}^{\times}}\int_{\mathbb{A}_{F}}\\
\chi_{2}(y)|y|^{\frac{1}{2}+\nu_{2}}\int_{\mathbb{A}_{F}^{\times}}W(y,(zy,b'z))\chi_{2}\chi_{3}
\omega'(z)|z|^{1+\nu_{2}+\nu_{3}+2s_{1}}d^{\times}z\\
\int_{\mathbb{A}_{F}^{\times}}\Phi_{1\times 2}((ty,tb')w'k^{-1})\chi_{2}\overline{\chi}_{3}(t)|t|^{1+\nu_{2}-\nu_{3}}d^{\times}t\overline{\psi(b')}db'd^{\times}ydk.
\end{multline}
Here, for $y, y'\in \mathbb{A}_{F}^{\times}$, the function $W(y,(zy',b'z))$ is defined by  
\begin{align*}
\eta(y)|y|^{\frac{1}{2}+s_{1}+s_{2}}\int_{\mathbb{A}_{F}^{\times}}h((t_{1},t_{1}^{-1}y),(zy',b'z);k,\nu_{1},\mathbf{s})\chi_{1}\omega'\overline{\eta}(t_{1})|t_{1}|^{\nu_{1}+s_{1}-s_{2}}
d^{\times}t_{1}.
\end{align*}

By Lemma~\ref{lem5.1} applied to the Fourier transform in the second
additive variable of $h$, $y\mapsto W(y,(zy',b'z))$ is a Whittaker
vector in the principal series $\pi_{1}$ induced by
$\eta|\cdot|^{s_{1}+s_{2}}$ and $\chi_{1}\omega'|\cdot|^{\nu_{1}+2s_{1}}$.

Define
\begin{multline*}
\mathfrak{F}_{1}(g';y,k):=\chi_{2}(\det g')|\det g'|^{\frac{1}{2}+\nu_{2}}\int_{\mathbb{A}_{F}^{\times}}W(y,(0,z)g')\chi_{2}\chi_{3}\\
\omega'(z)|z|^{1+\nu_{2}+\nu_{3}+2s_{1}}d^{\times}z
\int_{\mathbb{A}_{F}^{\times}}\Phi_{1\times 2}((0,t)g'w'k^{-1})\chi_{2}\overline{\chi}_{3}(t)|t|^{1+\nu_{2}-\nu_{3}}d^{\times}t,
\end{multline*}
which is a continuous function in $k\in K'$ and $g'\in G'(\mathbb{A}_{F})$ in the region 
\begin{align*}
 \Re(\nu_{2}+\nu_{3}+2s_{1})>0,\quad \Re(\nu_{2}-\nu_{3})>0.
\end{align*} 

For $t_{1},t_{2}\in\mathbb{A}_{F}^{\times}$ and $u\in\mathbb{A}_{F}$, a direct
calculation gives
\begin{align*}
\mathfrak{F}_{1}\left(\begin{pmatrix}
t_{1}& u\\
& t_{2}
\end{pmatrix}g';y,k\right)=\chi_{2}(t_{1})|t_{1}|^{\frac{1}{2}+\nu_{2}}
\overline{\chi}_{2}\overline{\omega}'(t_{2})|t_{2}|^{-\frac{3}{2}-\nu_{2}-2s_{1}}\mathfrak{F}_{1}(g';y,k).
\end{align*}
Hence, $\mathfrak{F}_{1}(\cdot;y,k)$ is a section in $\chi_{2}|\cdot|^{\nu_{2}}\boxplus \overline{\chi}_{2}\overline{\omega}'|\cdot|^{-1-\nu_{2}-2s_{1}}$. Let 
\begin{align*}
W_{\pi_{1},\pi_{2}}\left(a(y),
g';k\right):=\int_{\mathbb{A}_{F}}\mathfrak{F}_{1}\left(w'n(b')g';y,k\right)\overline{\psi(b')}db'.
\end{align*}

As a function of $y'$, $W_{\pi_{1},\pi_{2}}(a(y),a(y');k)$
lies in the Kirillov model of
$\pi_{2}=\chi_{2}|\cdot|^{\nu_{2}}\boxplus
\overline{\chi}_{2}\overline{\omega}'|\cdot|^{-1-\nu_{2}-2s_{1}}$.
Evaluating its defining Jacquet integral at $y'=y$ gives the integrand
in \eqref{c5.25}, hence
\begin{align*}
I_{\mathrm{degen}}(\mathbf{s},\varphi_{\boldsymbol{\nu}},\omega',\eta;w_{1}w_{2})=\int_{K'}\chi_{1}\overline{\chi}_{2}(\det k)\int_{\mathbb{A}_{F}^{\times}}W_{\pi_{1},\pi_{2}}\left(a(y),
a(y);k\right)d^{\times}ydk.	
\end{align*}

At finite places, $W_{\pi_{1},\pi_{2}}$ is a finite sum of products of
Whittaker vectors for $\pi_{1}$ and $\pi_{2}$. At Archimedean places its
$K_{v}'$-type coefficients decrease rapidly, also after differentiation.
Thus Lemma~\ref{lem3.2} applies termwise, locally uniformly away from poles.

The four leading powers at zero in the diagonal product are
$|y_{v}|_{v}^{1+\alpha}$, up to logarithms, with
\begin{align*}
 \alpha\in\{s_{1}+s_{2}+\nu_{2},\ s_{2}-s_{1}-\nu_{2}-1,\
             \nu_{1}+\nu_{2}+2s_{1},\ \nu_{1}-\nu_{2}-1\}.
\end{align*}
Here $\Re(\alpha)>0$ in \eqref{f5.19}.
The four powers are integrable at zero. Whittaker decay at infinity
and the absolute convergence of the Tate integrals defining $h$ give
convergence in \eqref{f5.19}. Lemma~\ref{lem3.2} and Tate continuation
give meromorphic continuation. The unramified integrals give
\begin{align*}
I_{\mathrm{degen}}(\mathbf{s},\varphi_{\boldsymbol{\nu}},\omega',\eta;w_{1}w_{2})
\propto\mathcal{L}(\mathbf{s},\boldsymbol{\nu},\pi,\omega',\eta;w_{1}w_{2}).
\end{align*}
The proportionality factor is the product of the meromorphic local
quotients by the respective $\mathcal{L}$-factors.

Let $v\notin S$. Then the local section
$\mathfrak{F}_{1,v}(\cdot;1,k_{v})$ is spherical and satisfies
\begin{align*}
\frac{\mathfrak{F}_{1,v}(I_{2};1,k_{v})}{L_{v}(1+\nu_{1}+s_{1}+s_{2},\chi_{1}\eta)}=L_{v}(1+\nu_{2}+\nu_{3}+2s_{1},\chi_{2}\chi_{3}
\omega')L_{v}(1+\nu_{2}-\nu_{3},\chi_{2}\overline{\chi}_{3}).
\end{align*}
For $v\notin S$, $W_{\pi_{1},\pi_{2}}$ is spherical in both variables,
so its diagonal integral is given by \eqref{eq3.6} times
its value at $(I_{2},I_{2})$.
The spherical Jacquet integral for $\pi_{2,v}$ divides the section
value by $L_{v}(2+2\nu_{2}+2s_{1},\chi_{2}^{2}\omega')$, so
\begin{multline}\label{eq5.27}
W_{\pi_{1},\pi_{2},v}\left(I_{2},
I_{2};k_{v}\right)=L_{v}(2+2\nu_{2}+2s_{1},\chi_{2}^{2}\omega')^{-1}L_{v}(1+\nu_{1}+s_{1}+s_{2},\chi_{1}\eta)\\
L_{v}(1+\nu_{2}+\nu_{3}+2s_{1},\chi_{2}\chi_{3}
\omega')L_{v}(1+\nu_{2}-\nu_{3},\chi_{2}\overline{\chi}_{3}).
\end{multline}

The diagonal integral of $W_{\pi_{1},\pi_{2},v}$ has denominator
$L_{v}(2,\omega_{\pi_{1}}\omega_{\pi_{2}})$ by \eqref{eq3.6} at $s=0$, where $\omega_{\pi_{1}}\omega_{\pi_{2}}
=\chi_{1}\eta|\cdot|^{\nu_{1}+s_{1}+s_{2}-1}$. It therefore cancels $L_{v}(1+\nu_{1}+s_{1}+s_{2},\chi_{1}\eta)$ in
\eqref{eq5.27}. After this cancellation, expand $L_{v}(1,\pi_{1,v}\times\pi_{2,v})$
by multiplicativity to get $\mathcal{L}_{v}(\mathbf{s},\boldsymbol{\nu},\pi,\omega',\eta;w_{1}w_{2})$.
Taking the finite product over $S$ gives the global assertion.
\end{proof}

\subsection{\texorpdfstring{The $w_{1}w_{2}$ Term of $J_{\mathrm{degen}}$}{The w1w2 Term of J-degen}}
Define 
\begin{multline*}
\mathcal{L}(s,\boldsymbol{\nu},\pi,\omega';w_{1}w_{2}):= \Lambda(2+2\nu_{2}+2s,\chi_{2}^{2}\omega')^{-1}\Lambda(1+\nu_{2}+\nu_{3}+2s,\chi_{2}\chi_{3}\omega')\\
\Lambda(1+\nu_{2}-\nu_{3},\chi_{2}\overline{\chi}_{3})\Lambda(1+\nu_{1}+\nu_{2}+2s,\chi_{1}\chi_{2}
\omega')\Lambda(\nu_{1}-\nu_{2},\chi_{1}\overline{\chi}_{2}).
\end{multline*}

Define the normalized local factor $\mathcal{P}_{v}^{\sharp}(s,\varphi_{\boldsymbol{\nu}},\omega';w_{1}w_{2})$ by
\begin{multline*}
\mathcal{P}_{v}^{\sharp}(s,\varphi_{\boldsymbol{\nu}},\omega';w_{1}w_{2}):=\mathcal{L}_{v}(s,\boldsymbol{\nu},\pi,\omega';w_{1}w_{2})^{-1}\int_{F_{v}^{\times}}\int_{F_{v}}\int_{F_{v}}\int_{F_{v}}\\
f_{\boldsymbol{\nu},v}\left(\begin{pmatrix}
1& & \\
b_{v}&1& \\
c_{v}& a_{v}& 1
\end{pmatrix}w_{1}w_{2}w_{1}\right)\overline{\psi_{v}(a_{v}+c_{v}z_{v})}db_{v}da_{v}dc_{v}\\
\chi_{2,v}\chi_{3,v}\omega_{v}'(z_{v})|z_{v}|_{v}^{1+\nu_{2}+\nu_{3}+2s}d^{\times}z_{v}.
\end{multline*}

\begin{prop}\label{prop5.6.}
$J_{\mathrm{degen}}(s,\varphi_{\boldsymbol{\nu}},\omega';w_{1}w_{2})$
converges absolutely in the region 
\begin{equation}\label{eq5.33.}
 \Re(\nu_{2}+\nu_{3}+2s)>0,\quad \Re(\nu_{2}-\nu_{3})>0,\quad \Re(\nu_{1}-\nu_{2})>1,
\end{equation}
\end{prop}

\begin{proof}
Insert \eqref{e5.3} with $g'=k\diag(t_{1},t_{2})n(b')$ and measure
\eqref{e1.4}. With $t_{1},t_{2}$ fixed, substitute
\begin{align*}
 (a,b,c,b',z,t)
 \longmapsto\left(\frac{a+t_{1}b'}{t_{2}},\frac{b}{t_{1}},
 \frac{c+bb'}{t_{2}},-\frac{t_{1}b'}{t_{2}},t_{2}z,t_{2}t\right).
\end{align*}
The additive Jacobian is $|t_{2}|^{-3}$ and the character becomes
$\overline{\psi}(at_{2}^{-1}+cz+b't_{1}t_{2}^{-1}+bb'z)$. For fixed $k$, define the partial Fourier transform $\Phi_{1\times5}(t_{1},t_{2},t_{2}^{-1},(z,b'z);k)$ by
\begin{align*}
\int_{\mathbb{A}_{F}}\int_{\mathbb{A}_{F}}\int_{\mathbb{A}_{F}}\Phi_{2\times 3}\left(k\begin{pmatrix}
b & t_{1} & \\
c& a & t_{2}
\end{pmatrix}w_{1}w_{2}w_{1}\right)
\overline{\psi(at_{2}^{-1})}da\overline{\psi(cz)}dc\overline{\psi(bb'z)}db.
\end{align*}

Since $\Phi_{1\times5}$ is Schwartz--Bruhat, the additional
substitution $b'\mapsto t^{-1}b'$ gives absolute
convergence in \eqref{eq5.33.}.

Make the simultaneous changes $b'\mapsto t_{1}^{-1}t_{2}b'$,
$z\mapsto t_{1}t_{2}^{-1}z$ and
$t\mapsto t_{1}t_{2}^{-1}t$. The additive Jacobian is
$|t_{1}^{-1}t_{2}|$, so the last Fourier character becomes
$\overline{\psi}(b')$. The $\Phi_{1\times5}$-argument becomes
$(t_{1},t_{2},t_{2}^{-1},(zt_{1}t_{2}^{-1},b'z);k)$.

For $g'\in G'(\mathbb{A}_{F})$, set
\begin{multline*}
\mathfrak{F}_{2}(g';t_{1},t_{2},t_{2}^{-1},k):=\chi_{1}\chi_{2}
\omega'(\det g')|\det g'|^{1+\nu_{1}+\nu_{2}+2s}\\
\int_{\mathbb{A}_{F}^{\times}}
\Phi_{1\times5}(t_{1},t_{2},t_{2}^{-1},(0,z)g';k)
\chi_{2}\chi_{3}\omega'(z)|z|^{1+\nu_{2}+\nu_{3}+2s}d^{\times}z\\
\int_{\mathbb{A}_{F}^{\times}}\Phi_{1\times 2}((0,t)g'w'k^{-1})\chi_{2}\overline{\chi}_{3}(t)|t|^{1+\nu_{2}-\nu_{3}}d^{\times}t.
\end{multline*}

The change $t_{1}\mapsto t_{2}t_{1}$ and the Jacquet integral of
$\mathfrak{F}_{2}$ give
\begin{multline}\label{f5.34}
J_{\mathrm{degen}}(s,\varphi_{\boldsymbol{\nu}},\omega';w_{1}w_{2})=\int_{K'}\chi_{1}\overline{\chi}_{2}(\det k)\int_{\mathbb{A}_{F}^{\times}}\int_{\mathbb{A}_{F}^{\times}}\\
W_{2}\left(a(t_{1}),t_{1}t_{2},t_{2},t_{2}^{-1},k\right)d^{\times}t_{1}
\chi_{1}^{2}
\omega'(t_{2})|t_{2}|^{1+2\nu_{1}+2s}d^{\times}t_{2}dk,
\end{multline}
where for each $y\in \mathbb{A}_{F}^{\times}$, 
\begin{align*}
W_{2}\left(a(y),t_{1}t_{2},t_{2},t_{2}^{-1},k\right):=\int_{\mathbb{A}_{F}}\mathfrak{F}_{2}\left(w'n(u)a(y);t_{1}t_{2},t_{2},t_{2}^{-1},k\right)
\overline{\psi(u)}du.
\end{align*}

The Godement transformation law shows that $\mathfrak{F}_{2}(\cdot;t_{1}t_{2},t_{2},t_{2}^{-1},k)$ is a section in the induced representation 
$\pi_{2}':=\chi_{1}\chi_{2}
\omega'|\cdot|^{\frac{1}{2}+\nu_{1}+\nu_{2}+2s}\boxplus 
\chi_{1}\overline{\chi}_{2}|\cdot|^{-\frac{1}{2}+\nu_{1}-\nu_{2}}$. As a result, the function 
\begin{align*}
y\mapsto W_{2}\left(a(y),t_{1}t_{2},t_{2},t_{2}^{-1},k\right)
\end{align*}
is a vector in the Kirillov model of $\pi_{2}'$.

In \eqref{eq5.33.}, $\mathfrak{F}_{2}(\cdot;r,t_{2},t_{2}^{-1},k)$ is
Schwartz--Bruhat in $r$, also after differentiation in $g'$.
Fourier inversion identifies evaluation at $r=t_{1}t_{2}$ with an integral
of right translations by $n(\xi t_{2})$, since these multiply the
Kirillov function by $\psi(\xi t_{1}t_{2})$. Thus
\begin{align*}
t_{1}\mapsto W_{2}\left(a(t_{1}),t_{1}t_{2},t_{2},t_{2}^{-1},k\right)
\end{align*}
is also a vector in the Kirillov model of $\pi_{2}'$. Furthermore, for $v\notin S$, 
\begin{align*}
W_{2,v}\left(a(t_{1,v}),t_{1,v}t_{2,v},t_{2,v},t_{2,v}^{-1},k_{v}\right)=\mathbf{1}_{\mathcal{O}_{v}^{\times}}(t_{2,v})W_{2,v}^{\circ}\left(a(t_{1,v})\right),
\end{align*}
where $W_{2,v}^{\circ}$ is the Jacquet integral of the spherical
section $\mathfrak{F}_{2,v}(\cdot;1,1,1,k_{v})$; its value at $I_{2}$ is
computed below. On this support, $(t_{1,v},t_{2,v})\in
\mathcal{O}_{v}\times\mathcal{O}_{v}^{\times}$, so
$\mathbf{1}_{\mathcal{O}_{v}}(t_{1,v}t_{2,v})=1$.
For $v\notin S$, the $t_{2,v}$-integral in \eqref{f5.34}
is over $\mathcal{O}_{v}^{\times}$ of volume one. Evaluating
$\int_{F_{v}^{\times}}W_{2,v}^{\circ}(a(t_{1,v}))\,d^{\times}t_{1,v}$ gives
\begin{equation}\label{e5.36}
 J_{\mathrm{degen},v}(s,\varphi_{\boldsymbol{\nu}},\omega';w_{1}w_{2})
 =W_{2,v}^{\circ}(I_{2})L_{v}(1/2,\pi_{2,v}'),
\end{equation}
The two Tate integrals in $\mathfrak{F}_{2,v}(\cdot;1,1,1,k_{v})$
give $L_{v}(1+\nu_{2}+\nu_{3}+2s,\chi_{2}\chi_{3}\omega')$ and
$L_{v}(1+\nu_{2}-\nu_{3},\chi_{2}\overline{\chi}_{3})$.  Applying the spherical
Jacquet integral gives
\begin{equation}\label{e5.34}
W_{2,v}^{\circ}(I_{2})=\frac{L_{v}(1+\nu_{2}+\nu_{3}+2s,\chi_{2}\chi_{3}\omega')
L_{v}(1+\nu_{2}-\nu_{3},\chi_{2}\overline{\chi}_{3})}
{L_{v}(2+2\nu_{2}+2s,\chi_{2}^{2}\omega')}.
\end{equation}
At a place $v\in S$, the local Tate integrals, the Jacquet integral,
and the final $\mathrm{GL}_{2}\times\mathrm{GL}_{1}$ zeta integral give
meromorphic continuation. The dependence on both $t_{2}$ and $t_{2}^{-1}$
ensures rapid decay at zero and infinity, also after differentiating $\mathfrak{F}_{2}$ in $g'$,
so the remaining $t_{2}$-integration preserves this continuation.
Dividing by $\mathcal{L}_{v}$ therefore gives a meromorphic family
$\mathcal{P}_{v}^{\heartsuit}$.

Equations \eqref{e5.36} and \eqref{e5.34} give \eqref{f5.15}.
\end{proof}

\subsection{\texorpdfstring{The $w_{2}w_{1}$ Term of $J_{\mathrm{degen}}^{\dag}$}{The w2w1 Term of J-degen dagger}}
Define
\begin{multline*}
\mathcal{L}^{\dag}(s,\boldsymbol{\nu},\pi,\omega';w_{2}w_{1}):=\Lambda(\nu_{2}-\nu_{3},\chi_{2}\overline{\chi}_{3})\\
\Lambda(1+\nu_{1}-\nu_{2},\chi_{1}\overline{\chi}_{2})\Lambda(1+\nu_{1}+\nu_{3}+2s,\chi_{1}\chi_{3}\omega').
\end{multline*}

Define the normalized local factor $\mathcal{P}_{v}^{\dag,\sharp}(s,\varphi_{\boldsymbol{\nu}},\omega';w_{2}w_{1})$ by
\begin{multline*}
\mathcal{P}_{v}^{\dag,\sharp}(s,\varphi_{\boldsymbol{\nu}},\omega';w_{2}w_{1}):=\mathcal{L}_{v}^{\dag}(s,\boldsymbol{\nu},\pi,\omega';w_{2}w_{1})^{-1}\int_{F_{v}^{\times}}\int_{F_{v}}\int_{F_{v}}
\overline{\psi_{v}(a_{v}z_{v}+b_{v})}\\
f_{\boldsymbol{\nu},v}\left(\begin{pmatrix}
1& & \\
-z_{v}^{-1}&1& \\
a_{v}& b_{v}& 1
\end{pmatrix}w_{1}w_{2}w_{1}\right)da_{v}db_{v}\chi_{2,v}\chi_{3,v}\omega_{v}'(z_{v})|z_{v}|_{v}^{\nu_{2}+\nu_{3}+2s}d^{\times}z_{v}.
\end{multline*}

\begin{prop}\label{prop5.8.}
$J_{\mathrm{degen}}^{\dag}(s,\varphi_{\boldsymbol{\nu}},\omega';w_{2}w_{1})$ converges absolutely in the region 
\begin{equation}\label{eq5.44}
 \Re(\nu_{2}-\nu_{3})>1,\quad \Re(\nu_{1}-\nu_{2})>0,\quad \Re(2s+\nu_{1}+\nu_{3})>0.
\end{equation}
\end{prop}

\begin{proof}
Make the change $a\mapsto za+b$, with Jacobian $|z|$, in the
unipotent integral. Move $\diag(1,z,z)$ to the left and use \eqref{2.9}.
In \eqref{e5.3}, use
$g'=k\diag(t_{1},t_{2})\left(\begin{smallmatrix}1&\\b'&1\end{smallmatrix}\right)$
and Haar measure \eqref{e1.5}. Make the changes
$a\mapsto a+z^{-1}b'$, $b\mapsto b-b'$,
$a\mapsto t_{2}^{-1}a$, $b\mapsto t_{2}^{-1}b$, followed by the simultaneous changes
$t\mapsto t_{2}t$, $b'\mapsto-t^{-1}t_{1}t_{2}^{-1}b'$. Their total additive Jacobian is
$|t_{1}t^{-1}t_{2}^{-3}|$. Define the partial Fourier transform
\begin{align*}
 h(r,x,u;t_{2},t_{2}',k)
 :=\int_{\mathbb{A}_{F}^{2}}
 \Phi_{2\times3}\left(k\begin{pmatrix}
 -u&r&\\a&b&t_{2}
 \end{pmatrix}\widetilde{w}\right)
 \overline{\psi}(ax+bt_{2}')\,da\,db.
\end{align*}
This is Schwartz--Bruhat in the five additive variables. Make the simultaneous changes of variables
$z\mapsto t_{2}x$ and $t_{1}\mapsto t_{2}xu$; multiplicative Haar
measure is unchanged.  Set
\begin{align*}
 &\phi_{k}(t):=\int_{\mathbb{A}_{F}}
 \Phi_{1\times2}((b',t)k^{-1})\,db',\\
 &H(x,u;s,\boldsymbol{\nu},k):=\int_{\mathbb{A}_{F}^{\times}}
 h(t_{2}xu,x,u;t_{2},t_{2}^{-1},k)
 \chi_{1}^{2}\omega'(t_{2})|t_{2}|^{1+2\nu_{1}+2s}\,d^{\times} t_{2}.
\end{align*}
The decay of $h$ in $t_{2},t_{2}^{-1},x,u$ makes $H$ Schwartz--Bruhat in
$(x,u)$ and entire in $(s,\boldsymbol{\nu})$, locally uniformly in the
parameters and $k$.
\begin{multline*}
 J_{\mathrm{degen}}^{\dag}(s,\varphi_{\boldsymbol{\nu}},\omega';w_{2}w_{1})
 =\int_{K'}\chi_{1}\overline{\chi}_{2}(\det k)\left[\int_{\mathbb{A}_{F}^{\times}}
 \phi_{k}(t)\chi_{2}\overline{\chi}_{3}(t)|t|^{\nu_{2}-\nu_{3}}\,d^{\times} t\right]\\
 \int_{(\mathbb{A}_{F}^{\times})^{2}}
 H(x,u;s,\boldsymbol{\nu},k)
 \chi_{1}\chi_{3}\omega'(x)|x|^{1+\nu_{1}+\nu_{3}+2s}
 \chi_{1}\overline{\chi}_{2}(u)|u|^{1+\nu_{1}-\nu_{2}}
 \,d^{\times} x\,d^{\times} u\,dk.
\end{multline*}
These three Tate integrals give the factors in
$\mathcal{L}^{\dag}(s,\boldsymbol{\nu},\pi,\omega';w_{2}w_{1})$, and their
global region of absolute convergence is \eqref{eq5.44}.
Locally, \eqref{e5.3} also contributes the constant
$\chi_{1,v}(\det\widetilde{w})=\chi_{1,v}(-1)$; its global product is
$1$.  Division by the three local $L$-factors gives an entire
function by the normalized Tate integrals, and integration over
$K_{v}'$ preserves holomorphy.

For $v\notin S$, self-duality gives
$\phi_{k_{v}}=\mathbf{1}_{\mathcal{O}_{v}}$; the pair $t_{2},t_{2}^{-1}$
restricts $t_{2}$ to $\mathcal{O}_{v}^{\times}$, of volume $1$. Hence
$H_{v}(x,u;s,\boldsymbol{\nu},k_{v})
=\mathbf{1}_{\mathcal{O}_{v}}(x)\mathbf{1}_{\mathcal{O}_{v}}(u)$.
Since $\chi_{1,v}(-1)=1$, the normalized local factor is $1$.
\end{proof}

\subsection{\texorpdfstring{The $w_{1}w_{2}w_{1}$ Term of $I_{\mathrm{degen}}$}{The w1w2w1 Term of I-degen}}\label{sec5.11}
Define
\begin{multline*}
\mathcal{L}(\mathbf{s},\boldsymbol{\nu},\pi,\omega',\eta;w_{1}w_{2}w_{1}):=\Lambda(2+2\nu_{3}+2s_{1},\chi_{3}^{2}\omega')^{-1}\Lambda(s_{2}-s_{1}-\nu_{3},\overline{\chi}_{3}\overline{\omega}'\eta)\\\Lambda(1+\nu_{1}+\nu_{3}+2s_{1},\chi_{1}\chi_{3}\omega')
\Lambda(\nu_{1}-\nu_{3},\chi_{1}\overline{\chi}_{3})\Lambda(1+\nu_{3}+s_{1}+s_{2},\chi_{3}\eta)\\
\Lambda(1+\nu_{2}+\nu_{3}+2s_{1},\chi_{2}\chi_{3}\omega')\Lambda(\nu_{2}-\nu_{3},\chi_{2}\overline{\chi}_{3}).
\end{multline*}
 
Define the normalized local factor
$\mathcal{P}_{v}^{\sharp}(\mathbf{s},\varphi_{\boldsymbol{\nu}},\omega',\eta;w_{1}w_{2}w_{1})$ by
\begin{multline}\label{eq:local-long-weyl-degenerate}
\mathcal{P}_{v}^{\sharp}(\mathbf{s},\varphi_{\boldsymbol{\nu}},\omega',\eta;w_{1}w_{2}w_{1})
:=\mathcal{L}_{v}(\mathbf{s},\boldsymbol{\nu},\pi,\omega',\eta;w_{1}w_{2}w_{1})^{-1}
\int_{(F_{v}^{\times})^{2}}\int_{F_{v}^{4}}\\
f_{\boldsymbol{\nu},v}\left(\begin{pmatrix}
1& & \\
c_{v}&1& \\
b_{v}& a_{v}& 1
\end{pmatrix}n_{23}(c_{v}')w_{1}w_{2}\right)\overline{\psi_{v}(c_{v}z_{v}+c_{v}'y_{v})}dc_{v}dc_{v}'da_{v}db_{v}\\
\chi_{2,v}\chi_{3,v}\omega_{v}'(z_{v})|z_{v}|_{v}^{1+\nu_{2}+\nu_{3}+2s_{1}}\overline{\chi}_{3,v}\overline{\omega}_{v}'\eta_{v}(y_{v})|y_{v}|_{v}^{s_{2}-s_{1}-\nu_{3}}d^{\times}z_{v}d^{\times}y_{v}.
\end{multline}

\begin{prop}\label{prop5.6}
$I_{\mathrm{degen}}(\mathbf{s},\varphi_{\boldsymbol{\nu}},\omega',\eta;w_{1}w_{2}w_{1})$ converges absolutely in the following region 
\begin{equation}\label{e5.33}
\begin{cases}
\Re(\nu_{1}-\nu_{3})>1,\ \ \Re(\nu_{2}-\nu_{3})>1\\
\Re(\nu_{2}+\nu_{3}+2s_{1})>0,\ \ \Re(2s_{1}+2\nu_{3})>0,\\
\Re(s_{2}-s_{1}-\nu_{3})>1
\end{cases}	
\end{equation}
\end{prop}
See \textsection\ref{sec5.11.2} for the proof.
\subsubsection{An auxiliary integral}
For $(\mathbf{s},\boldsymbol{\nu})\in \mathbb{C}^{5}$, let  
\begin{multline*}
\mathbf{L}(\mathbf{s},\boldsymbol{\nu}):=\Lambda(2+2\nu_{3}+2s_{1},\chi_{3}^{2}\omega')^{-1}\Lambda(1+\nu_{3}+s_{1}+s_{2},\chi_{3}\eta)
\Lambda(\nu_{1}-\nu_{3},\chi_{1}\overline{\chi}_{3})\\
\Lambda(s_{2}-s_{1}-\nu_{3},\overline{\chi}_{3}\overline{\omega}'\eta)\Lambda(1+\nu_{1}+\nu_{3}+2s_{1},\chi_{1}\chi_{3}\omega').
\end{multline*}

Let $\Phi_{*}(\cdot;k)=\bigotimes_{v}\Phi_{*,v}(\cdot;k_{v})$ be a Schwartz--Bruhat function on
$M_{2\times2}(\mathbb{A}_{F})$, depending smoothly on $k\in K'$.
Define 
\begin{multline}\label{e5.21}
\mathcal{I}(\mathbf{s},\boldsymbol{\nu};k)=\int_{\mathbb{A}_{F}^{\times}}
\int_{\mathbb{A}_{F}^{\times}}\int_{\mathbb{A}_{F}^{\times}}h(t_{1},t_{2},y;k)\chi_{1}\overline{\chi}_{3}(t_{2})|t_{2}|^{\nu_{1}-\nu_{3}}d^{\times}t_{2}\\ 
\overline{\chi}_{3}\overline{\omega}'\eta(y)|y|^{s_{2}-s_{1}-\nu_{3}}d^{\times}y\chi_{1}\eta(t_{1})|t_{1}|^{1+\nu_{1}+s_{1}+s_{2}}d^{\times}t_{1},
\end{multline}
where 
\begin{equation}\label{e5.25}
h(t_{1},t_{2},y;k):=\int_{\mathbb{A}_{F}}\int_{\mathbb{A}_{F}}\Phi_{*}\left(\begin{pmatrix}
t_{1} & c' \\
at_{1} & ac'+t_{2}
\end{pmatrix};k
\right)\overline{\psi(c'y)}dc'da.
\end{equation}

For each $v\in\Sigma_{F}$, denote by $\mathbf{L}_{v}(\mathbf{s},\boldsymbol{\nu})$ and $\mathcal{I}_{v}(\mathbf{s},\boldsymbol{\nu};k)$ the local components of $\mathbf{L}(\mathbf{s},\boldsymbol{\nu})$ and $\mathcal{I}(\mathbf{s},\boldsymbol{\nu};k)$, respectively. Let $\mathcal{P}_{v}^{\sharp}(\mathbf{s},\boldsymbol{\nu};k):=\mathbf{L}_{v}(\mathbf{s},\boldsymbol{\nu})^{-1}\mathcal{I}_{v}(\mathbf{s},\boldsymbol{\nu};k)$.

\begin{lemma}\label{lem5.6}
$\mathcal{I}(\mathbf{s},\boldsymbol{\nu};k)$ converges absolutely for
\begin{equation}\label{5.25}
 \Re(\nu_{1}-\nu_{3})>1,\quad \Re(s_{2}-s_{1}-\nu_{3})>1,\quad \Re(2s_{1}+2\nu_{3})>0.
\end{equation}

For each $v\in\Sigma_{F}$, the function $\mathcal{P}_{v}^{\sharp}(\mathbf{s},\boldsymbol{\nu};k)$ is meromorphic in the variables $(\mathbf{s},\boldsymbol{\nu})\in \mathbb{C}^{5}$. Moreover, when $v<\infty$, $v\nmid\mathfrak{D}_{F}$, the characters $\chi_{i,v}$, $1\leq i\leq 3$, $\omega_{v}'$, and $\eta_{v}$ are all unramified, and $\Phi_{*,v}(\cdot;k_{v})=\mathbf{1}_{(\mathcal{O}_{v})^{4}}$, we have
\begin{equation}\label{e.5.23}
\mathcal{P}_{v}^{\sharp}(\mathbf{s},\boldsymbol{\nu};k)\equiv 1.
\end{equation}
In particular, the global function
$\mathcal{P}^{\sharp}(\mathbf{s},\boldsymbol{\nu};k)
:=\prod_{v}\mathcal{P}_{v}^{\sharp}(\mathbf{s},\boldsymbol{\nu};k)$ is a finite product and is meromorphic.

The meromorphic continuation of $\mathcal{I}(\mathbf{s},\boldsymbol{\nu};k)$ to $\mathbb{C}^{5}$ is
\begin{align*}
\mathcal{I}^{\heartsuit}(\mathbf{s},\boldsymbol{\nu};k)=\mathcal{P}^{\sharp}(\mathbf{s},\boldsymbol{\nu};k)\mathbf{L}(\mathbf{s},\boldsymbol{\nu}).
\end{align*}
\end{lemma}

\begin{proof}
Write the Mellin exponents of $t_{2}$, $y$ and $t_{1}$ in \eqref{e5.21} as
\begin{align*}
 \alpha=\nu_{1}-\nu_{3},\qquad
 \beta=s_{2}-s_{1}-\nu_{3},\qquad
 \kappa=1+\nu_{1}+s_{1}+s_{2}.
\end{align*}
For $t_{1}\neq0$, $h(t_{1},\cdot,\cdot;k)$ is Schwartz--Bruhat.
Split the $a$-integral at $|a|_{v}=1$.  In the second part make the change $c'\mapsto a^{-1}(u-t_{2})$;
the Jacobian is $|a|_{v}^{-1}$ and the character becomes
$\psi_{v}((t_{2}-u)y/a)$.  Schwartz decay and Fourier integration in
$u$ then give, at an Archimedean place and for
$A=\Re(\alpha)>0$, $B=\Re(\beta)>0$,
\begin{align*}
 \int_{(F_{v}^{\times})^{2}}|h_{v}(t_{1},t_{2},y;k_{v})|
 |t_{2}|_{v}^{A}|y|_{v}^{B}\,d^{\times} t_{2}\,d^{\times} y
 \ll_{l} (1+|t_{1}|_{v})^{-l}(1+|t_{1}|_{v}^{-A-B}).
\end{align*}
The bound is locally uniform in $A,B$ and $k_{v}$. At finite places,
$|t_{1}|_{v}\ll1$ on the support. For $|a|_{v}>1$, the
$t_{2},y$ integrals are $O(|a|_{v}^{A+B-1})$ after the change in $c'$.
Including $|a|_{v}\leq1$, sum the bounds using
$1+\sum_{1<q_{v}^{j}\ll|t_{1}|_{v}^{-1}}q_{v}^{j(A+B)}
\ll1+|t_{1}|_{v}^{-A-B}$.
For $v\notin S$ and $m=e_{v}(t_{1})$,
\begin{multline*}
 \int_{(F_{v}^{\times})^{2}}|h_{v}(t_{1},t_{2},y;k_{v})|
 |t_{2}|_{v}^{A}|y|_{v}^{B}\,d^{\times} t_{2}\,d^{\times} y\\
 \leq \zeta_{v}(A)\zeta_{v}(B)\mathbf{1}_{m\geq0}
 \left(1+(1-q_{v}^{-1})\sum_{i=1}^{m}q_{v}^{(A+B)i}\right).
\end{multline*}  Summing these bounds over $m$ proves absolute convergence of the
global integral when
$\Re(\alpha),\Re(\beta),\Re(\kappa-\alpha-\beta)>1$.
Since $\kappa-\alpha-\beta=1+2s_{1}+2\nu_{3}$, this proves
convergence and holomorphy in \eqref{5.25}.

Apply the Tate functional equation locally, first in $t_{2}$ and then
in $y$. The three local integrals have the common convergence region
\begin{align*}
 0<\Re(\alpha)<1,\qquad 0<\Re(\beta)<1,\quad \Re(\kappa)>3.
\end{align*}
The substitution $a\mapsto t_{1}^{-1}b$ gives convergence near $t_{1}=0$; rapid decay gives convergence at infinity.
Thus these local functional equations can be integrated in the remaining variables.
Their gamma factors are given below and have global product $1$.
Taking Euler products gives the meromorphic identity
\begin{multline}\label{5.30}
 \mathcal{I}(\mathbf{s},\boldsymbol{\nu};k)
 =\int_{(\mathbb{A}_{F}^{\times})^{3}}\mathcal{F}h(t_{1},t_{2},y;k)
 \overline{\chi}_{1}\chi_{3}(t_{2})|t_{2}|^{1-\alpha}\\
\chi_{3}\omega'\overline{\eta}(y)|y|^{1-\beta}\chi_{1}\eta(t_{1})|t_{1}|^{\kappa}
 \,d^{\times} t_{2}\,d^{\times} y\,d^{\times} t_{1},
\end{multline}
which is absolutely convergent in
\begin{equation}\label{e5.31}
 \Re(\alpha)<0,\qquad \Re(\beta)<0,\quad \Re(\kappa)>1.
\end{equation}
Here $\mathcal{F}h$ is the partial Fourier transform in the last two
variables. Fourier inversion in \eqref{e5.25} gives
\begin{equation}\label{e5.72}
 \mathcal{F}h(t_{1},t_{2},y;k)
 =\int_{\mathbb{A}_{F}^{2}}
 \Phi_{*}\left(\begin{pmatrix}t_{1}&y\\at_{1}&u\end{pmatrix};k\right)
 \overline{\psi}(at_{2}y)\psi(ut_{2})\,da\,du.
\end{equation}
After $a\mapsto t_{1}^{-1}b$, the right-hand side is $|t_{1}|^{-1}$ times a
Schwartz--Bruhat function evaluated at $(t_{1},y,-t_{2}y/t_{1},t_{2})$,
which gives local convergence in \eqref{e5.31}.

By the two local Tate functional equations, $\mathcal{I}_{v}$ is the
local integral in \eqref{5.30} multiplied by
\begin{equation}\label{eq:auxiliary-local-tate-factors}
 \gamma(\alpha,\chi_{1,v}\overline{\chi}_{3,v},\psi_{v})^{-1}
 \gamma(\beta,\overline{\chi}_{3,v}\overline{\omega_{v}'}\eta_{v},\psi_{v})^{-1}.
\end{equation}
Substitute \eqref{e5.72} into \eqref{5.30} and make the change of variable $y\mapsto t_{2}^{-1}y$.
For $y,y'\in\mathbb{A}_{F}^{\times}$, set
\begin{multline*}
W_{4}'\left(a(y)
,a(y')\right):=\overline{\eta}(y')|y'|^{\frac{1}{2}-s_{2}-s_{1}}\int_{\mathbb{A}_{F}}\int_{\mathbb{A}_{F}^{\times}}W_{3}'\left(a(y)
;t_{1},a\right)\\
\chi_{1}\eta(t_{1})|t_{1}|^{1+\nu_{1}+s_{1}+s_{2}}d^{\times}t_{1}\overline{\psi(ay')}da,
\end{multline*}
where $W_{3}'\left(a(y)
;t_{1},a\right)$ is defined by 
\begin{multline*}
\chi_{3}\omega'(y)|y|^{\frac{1}{2}+2s_{1}+\nu_{3}}\int_{\mathbb{A}_{F}^{\times}}\int_{\mathbb{A}_{F}}\Phi_{*}\left(\begin{pmatrix}
t_{1} & t_{2}^{-1}y \\
at_{1} & u
\end{pmatrix};k
\right)
\psi(ut_{2})du\\
\overline{\chi}_{1}\overline{\omega}'\eta(t_{2})|t_{2}|^{s_{2}-s_{1}-\nu_{1}}d^{\times}t_{2}.
\end{multline*}

Hence 
\begin{equation}\label{5.34}
\mathcal{I}(\mathbf{s},\boldsymbol{\nu};k)=\int_{\mathbb{A}_{F}^{\times}}
W_{4}'\left(a(y)
,a(y)\right)d^{\times}y.
\end{equation}

After $u\mapsto-u$ and $t_{2}\mapsto y/t_{2}$, apply Lemma~\ref{lem5.1}
with $t_{1},a$ fixed. Then $W_{3}'$ lies in the Kirillov model of
$\pi_{3}=\chi_{3}\omega'|\cdot|^{2s_{1}+\nu_{3}}\boxplus
\overline{\chi}_{1}\chi_{3}\eta|\cdot|^{-\nu_{1}+\nu_{3}+s_{1}+s_{2}}$.
As a function of $y'$, $W_{4}'(a(y),a(y'))$ lies in the Kirillov model of
$\pi_{4}=\overline{\eta}|\cdot|^{-s_{1}-s_{2}}\boxplus
\chi_{1}|\cdot|^{\nu_{1}}$.

For $\Phi_{*}$ factored by columns, \eqref{5.34} is the convolution
of these two Kirillov vectors. The general case follows by linearity at finite places
and continuity in the completed Schwartz tensor product at infinity. Lemma~\ref{lem3.2}, together with the factor in
\eqref{eq:auxiliary-local-tate-factors}, therefore continues the
original local integral meromorphically from \eqref{e5.31}.
For $v\notin S$, both Kirillov vectors in \eqref{5.34} are
spherical and equal $1$ at $I_{2}$. The local integral is then
given by \eqref{eq3.6} as
\begin{align*}
L_{v}(2,\omega_{\pi_{3,v}}\omega_{\pi_{4,v}})^{-1}L_{v}(1,\pi_{3,v}\times\pi_{4,v}).
\end{align*}
The two inverse gamma factors in \eqref{eq:auxiliary-local-tate-factors}
replace the numerator factors at $1-\alpha$ and $1-\beta$ by those at $\alpha$ and $\beta$, respectively. Thus
$\mathcal{I}_{v}=\mathbf{L}_{v}(\mathbf{s},\boldsymbol{\nu})$, proving
\eqref{e.5.23}. The normalized local factors are meromorphic and equal $1$ for $v\notin S$,
so their product gives the asserted global continuation.
\end{proof}

\subsubsection{Proof of Proposition \ref{prop5.6}}\label{sec5.11.2}
By linearity and continuity in the Schwartz functions, it suffices to take factorizable Godement data. Then
\begin{multline}\label{5.21}
I_{\mathrm{degen}}(\mathbf{s},\varphi_{\boldsymbol{\nu}},\omega',\eta;w_{1}w_{2}w_{1})=\int_{\mathbb{A}_{F}^{\times}}\int_{\mathbb{A}_{F}^{\times}}\int_{\mathbb{A}_{F}}\int_{\mathbb{A}_{F}}\int_{\mathbb{A}_{F}}\int_{\mathbb{A}_{F}}f_{\boldsymbol{\nu}}(g)\\
\overline{\psi(c+c')}dadbdcdc' \omega'(z)|z|^{2s_{1}}\eta(y)|y|^{s_{1}+s_{2}}d^{\times}zd^{\times}y,
\end{multline}
where
\begin{align*}
g=\begin{pmatrix}
1& & \\
c&1& \\
b&a& 1
\end{pmatrix}n_{23}(c')\diag(1,yz,z)w_{1}w_{2}.
\end{align*}

Use
$g'=k\diag(t_{1},t_{2})\left(\begin{smallmatrix}1&\\b'&1\end{smallmatrix}\right)$
in \eqref{e5.3} and substitute into \eqref{5.21}, with Haar measure
\eqref{e1.5}. The translations $(a,b)\mapsto(a-b',b-cb')$, followed
by $b'\mapsto-b'$, have Jacobian $1$ and leave the additive character unchanged. The compact integral contains the factor
$\chi_{1}\overline{\chi}_{2}(\det k)$.

Change $t\mapsto t_{2}t$ and then $b'\mapsto t_{1}t_{2}^{-1}t^{-1}b'$,
so that $\Phi_{1\times2}$ has argument $(b',t)k^{-1}$.
With the multiplicative variables fixed, change
\begin{align*}
 (c,b,c',a)\mapsto(t_{1}^{-1}c,t_{2}^{-1}b,yc',t_{2}^{-1}y^{-1}z^{-1}a).
\end{align*}
Together with the change in $b'$, these changes contribute
$|t^{-1}t_{2}^{-3}z^{-1}|$ to the additive measure. Then make
$t_{2}\mapsto z^{-1}t_{2}$, $z\mapsto y^{-1}z$, and $t_{1}\mapsto z^{-1}t_{1}$.
The factors involving $\Phi_{2\times3}$ and $\psi$ become
\begin{align*}
 \Phi_{2\times3}\left(k\begin{pmatrix}c&t_{1}&t_{1}c'\\b&a&t_{2}+ac'\end{pmatrix}w_{1}w_{2}\right),
 \qquad \overline{\psi}(ct_{1}^{-1}z+c'y).
\end{align*}

Let $(t_{1}',t_{2}',t_{3}',t_{4}')\in \mathbb{A}_{F}\times\mathbb{A}_{F}\times\mathbb{A}_{F}\times\mathbb{A}_{F}$. Define 
\begin{multline}\label{eq:long-weyl-Phi}
\Phi_{*}\left(\begin{pmatrix}
t_{1}'& t_{2}'\\
t_{3}'& t_{4}'
\end{pmatrix};k
\right):=\int_{\mathbb{A}_{F}^{\times}}\int_{\mathbb{A}_{F}}\int_{\mathbb{A}_{F}}\Phi_{2\times 3}\left(k\begin{pmatrix}
c& t_{1}'& t_{2}'\\
b& t_{3}'& t_{4}'
\end{pmatrix}
w_{1}w_{2}\right)\overline{\psi(cz)}dcdb\\
\chi_{2}\chi_{3}\omega'(z)|z|^{1+\nu_{2}+\nu_{3}+2s_{1}}d^{\times}z.
\end{multline}

The integral $\Phi_{*}\left(\begin{pmatrix}
t_{1}'& t_{2}'\\
t_{3}'& t_{4}'
\end{pmatrix};k
\right)$ converges absolutely in the region
$\Re(\nu_{2}+\nu_{3}+2s_{1})>0$.  The $z$-integration is a Tate zeta
integral whose unramified local factor is
$L_{v}(1+\nu_{2}+\nu_{3}+2s_{1},\chi_{2,v}\chi_{3,v}\omega_{v}')$.
In this region, $\Phi_{*}(\cdot;k)$ is Schwartz--Bruhat on
$M_{2\times 2}(\mathbb{A}_{F})$. At each place, dividing the local
$\Phi_{*,v}(\cdot;k_{v})$ by
$L_{v}(1+\nu_{2}+\nu_{3}+2s_{1},\chi_{2,v}\chi_{3,v}\omega_{v}')$
gives an entire family of Schwartz--Bruhat functions by Tate's theorem.

Next make the simultaneous changes of variables $c'\mapsto t_{1}^{-1}c'$,
$a\mapsto t_{1}a$, $y\mapsto t_{1}y$, and $z\mapsto t_{1}z$.
The additive Jacobians cancel, and the $t_{1}$-weight becomes
$\chi_{1}\eta(t_{1})|t_{1}|^{1+\nu_{1}+s_{1}+s_{2}}$.  Thus 
\begin{multline}\label{e5.75}
I_{\mathrm{degen}}(\mathbf{s},\varphi_{\boldsymbol{\nu}},\omega',\eta;w_{1}w_{2}w_{1})=\int_{K'}\chi_{1}\overline{\chi}_{2}(\det k)\int_{\mathbb{A}_{F}}\int_{\mathbb{A}_{F}^{\times}}\Phi_{1\times 2}((b',t)k^{-1})\\
\chi_{2}\overline{\chi}_{3}(t)|t|^{\nu_{2}-\nu_{3}}d^{\times}tdb'\mathcal{I}(\mathbf{s},\boldsymbol{\nu};k)
dk,
\end{multline}
where $\mathcal{I}$ is defined by \eqref{e5.21}--\eqref{e5.25}
with $\Phi_{*}$ given by \eqref{eq:long-weyl-Phi}.

By Lemma~\ref{lem5.6} and \eqref{e5.75}, the integral
$I_{\mathrm{degen}}(\mathbf{s},\varphi_{\boldsymbol{\nu}},\omega',\eta;w_{1}w_{2}w_{1})$
converges when $(\mathbf{s},\boldsymbol{\nu})$ lies in \eqref{e5.33}.
The two additional factors are
$\Lambda(\nu_{2}-\nu_{3},\chi_{2}\overline{\chi}_{3})$, from the Tate integral
in \eqref{e5.75}, and
$\Lambda(1+\nu_{2}+\nu_{3}+2s_{1},\chi_{2}\chi_{3}\omega')$, from the
$z$-integration in the definition of $\Phi_{*}$. Multiplying these by
$\mathbf{L}(\mathbf{s},\boldsymbol{\nu})$ gives
$\mathcal{L}(\mathbf{s},\boldsymbol{\nu},\pi,\omega',\eta;w_{1}w_{2}w_{1})$.
For $v\notin S$, the self-dual additive measure gives
\begin{align*}
 \Phi_{*,v}(X;k_{v})
 =L_{v}(1+\nu_{2}+\nu_{3}+2s_{1},\chi_{2,v}\chi_{3,v}\omega_{v}')
  \mathbf{1}_{M_{2}(\mathcal{O}_{v})}(X).
\end{align*}
The $\Phi_{1\times2,v}$ integral equals
$L_{v}(\nu_{2}-\nu_{3},\chi_{2,v}\overline{\chi}_{3,v})$.
Lemma~\ref{lem5.6} and $\operatorname{vol}(K_{v}')=1$ therefore give
$\mathcal{L}_{v}$, proving \eqref{f5.15}.

At the remaining places, Lemma~\ref{lem5.6} applied to the normalized
$\Phi_{*,v}$ gives meromorphic local quotients. Their product with
$\mathcal{L}$ agrees with the original integral in \eqref{e5.33},
giving its continuation.

\subsection{\texorpdfstring{The $w_{1}w_{2}w_{1}$ Term of $J_{\mathrm{degen}}$}{The w1w2w1 Term of J-degen}}
Define
\begin{multline*}
\mathcal{L}(s,\boldsymbol{\nu},\pi,\omega';w_{1}w_{2}w_{1}):= \Lambda(2+2\nu_{3}+2s,\chi_{3}^{2}\omega')^{-1}\Lambda(1+\nu_{1}+\nu_{3}+2s,\chi_{1}\chi_{3}\omega')
\\
\Lambda(\nu_{1}-\nu_{3},\chi_{1}\overline{\chi}_{3})\Lambda(1+\nu_{2}+\nu_{3}+2s,\chi_{2}\chi_{3}\omega')\Lambda(\nu_{2}-\nu_{3},\chi_{2}\overline{\chi}_{3}).
\end{multline*}

Define the normalized local factor
$\mathcal{P}_{v}^{\sharp}(s,\varphi_{\boldsymbol{\nu}},\omega';w_{1}w_{2}w_{1})$ by
\begin{multline*}
\mathcal{P}_{v}^{\sharp}(s,\varphi_{\boldsymbol{\nu}},\omega';w_{1}w_{2}w_{1}):=\mathcal{L}_{v}(s,\boldsymbol{\nu},\pi,\omega';w_{1}w_{2}w_{1})^{-1}\int_{F_{v}^{\times}}\int_{F_{v}}\int_{F_{v}}\int_{F_{v}}\int_{F_{v}}\\
f_{\boldsymbol{\nu},v}\left(\begin{pmatrix}
1& & \\
c_{v}&1& \\
b_{v}& a_{v}& 1
\end{pmatrix}n_{23}(c_{v}')w_{1}w_{2}\right)\overline{\psi_{v}(c_{v}z_{v}+c_{v}')}dc_{v}dc_{v}'da_{v}db_{v}\\
\chi_{2,v}\chi_{3,v}\omega_{v}'(z_{v})|z_{v}|_{v}^{1+\nu_{2}+\nu_{3}+2s}d^{\times}z_{v}.
\end{multline*}

\begin{prop}\label{prop5.8}
$J_{\mathrm{degen}}(s,\varphi_{\boldsymbol{\nu}},\omega';w_{1}w_{2}w_{1})$ converges absolutely in the following region 
\begin{equation}\label{e5.33.}
\begin{cases}
\Re(\nu_{2}-\nu_{3})>1,\ \ \Re(\nu_{1}-\nu_{3})>1,\\
\Re(\nu_{2}+\nu_{3}+2s)>0,\ \ \Re(2\nu_{3}+2s)>0.
\end{cases}	
\end{equation}
\end{prop}

\begin{proof}
Set $y=1$ and $s_{1}=s$ in \eqref{5.21} before making the changes
that lead to \eqref{e5.75}, and omit the $y$-integration. With
$\Phi_{*}$ from \eqref{eq:long-weyl-Phi} at $s_{1}=s$, the final change
$z\mapsto t_{1}z$ gives
\begin{multline}\label{e5.41}
J_{\mathrm{degen}}(s,\varphi_{\boldsymbol{\nu}},\omega';w_{1}w_{2}w_{1})=\int_{K'}\chi_{1}\overline{\chi}_{2}(\det k)\int_{\mathbb{A}_{F}}\int_{\mathbb{A}_{F}^{\times}}\Phi_{1\times 2}((b',t)k^{-1})\\
\chi_{2}\overline{\chi}_{3}(t)|t|^{\nu_{2}-\nu_{3}}d^{\times}tdb'\int_{\mathbb{A}_{F}}\int_{\mathbb{A}_{F}^{\times}}\int_{\mathbb{A}_{F}}\int_{\mathbb{A}_{F}^{\times}}
\Phi_{*}\left(\begin{pmatrix}
t_{1} & t_{1}c' \\
a & ac'+t_{2}
\end{pmatrix};k\right)\\
\chi_{1}\overline{\chi}_{3}(t_{2})|t_{2}|^{\nu_{1}-\nu_{3}}d^{\times}t_{2}da\chi_{1}\chi_{3}\omega'(t_{1})|t_{1}|^{1+\nu_{1}+\nu_{3}+2s}d^{\times}t_{1}\overline{\psi(c')}dc'
dk,
\end{multline}
Locally, division by
$L_{v}(1+\nu_{2}+\nu_{3}+2s,\chi_{2,v}\chi_{3,v}\omega_{v}')$
gives an entire Schwartz--Bruhat family. The initial region is
$\Re(\nu_{2}+\nu_{3}+2s)>0$.

Put $\alpha=\nu_{1}-\nu_{3}$ and
$\kappa=1+\nu_{1}+\nu_{3}+2s$. Take absolute values in the inner four integrals of
\eqref{e5.41} and make the changes $c'\mapsto t_{1}^{-1}c'$ and
$t_{2}\mapsto t_{2}-at_{1}^{-1}c'$. Up to a fixed measure constant, the local integral becomes
\begin{align*}
 \int_{F_{v}^{4}}\left|\Phi_{*,v}\left(
 \begin{pmatrix}t_{1}&c'\\a&t_{2}\end{pmatrix};k_{v}\right)\right|
 |t_{1}t_{2}-ac'|_{v}^{\Re(\alpha)-1}
 |t_{1}|_{v}^{\Re(\kappa-\alpha)-1}\,dt_{1}\,dc'\,da\,dt_{2}.
\end{align*}
For $\Re(\alpha)>1$ and $\Re(\kappa-\alpha)>1$, both powers are positive,
so compact support at finite places and rapid decay at Archimedean places give convergence.
For $\Phi_{*,v}=\mathbf{1}_{M_{2}(\mathcal{O}_{v})}$, summing over valuations bounds the local integral by
\begin{align*}
 \zeta_{v}(\Re(\alpha))\zeta_{v}(\Re(\kappa))
 \zeta_{v}(\Re(\kappa-\alpha)).
\end{align*}
Together with the two outer Tate integrals, this proves absolute convergence in \eqref{e5.33.}.

For fixed $c'$ and $k$, the $(t_{1},t_{2},a)$-integral in
\eqref{e5.41} converges for $\Re(\alpha)>1$, $\Re(\kappa)>1$ and
$\Re(\nu_{2}+\nu_{3}+2s)>0$. The Tate integral argument of
Lemma~\ref{lem5.6} continues the $(t_{1},t_{2},a)$-integral meromorphically
in $\alpha$ for $\Re(\kappa)>1$ and $\Re(\nu_{2}+\nu_{3}+2s)>0$.

Moreover, in the region
\begin{equation}\label{5.41}
\begin{cases}
\Re(\nu_{1}-\nu_{3})<0,\ \ \Re(\nu_{2}+\nu_{3}+2s)>0,\\
\Re(\nu_{1}+\nu_{3}+2s)>0,\ \ \Re(\nu_{2}-\nu_{3})>1,
\end{cases}
\end{equation}
the continuation of the $(t_{1},t_{2},a)$-integral in \eqref{e5.41} is given by the convergent integral
\begin{multline}\label{5.42}
\int_{\mathbb{A}_{F}^{\times}}\int_{\mathbb{A}_{F}^{\times}}\mathcal{F}\Phi_{*}\left(\begin{pmatrix}
t_{1} & t_{1}c' \\
t_{2} & t_{2}c'
\end{pmatrix};k\right)
\overline{\chi}_{1}\chi_{3}(t_{2})|t_{2}|^{1-\nu_{1}+\nu_{3}}d^{\times}t_{2}\\
\chi_{1}\chi_{3}\omega'(t_{1})|t_{1}|^{1+\nu_{1}+\nu_{3}+2s}d^{\times}t_{1},
\end{multline}
where 
\begin{align*}
\mathcal{F}\Phi_{*}\left(\begin{pmatrix}
t_{1} & t_{1}c' \\
t_{2} & t_{2}c'
\end{pmatrix};k\right):=\int_{\mathbb{A}_{F}}\int_{\mathbb{A}_{F}}\Phi_{*}\left(\begin{pmatrix}
t_{1} & t_{1}c' \\
a & u
\end{pmatrix};k\right)\psi(ut_{2})du\overline{\psi(at_{2}c')}da.
\end{align*}

The local Tate functional equation gives the factor $\gamma(\alpha,\chi_{1,v}\overline{\chi}_{3,v},\psi_{v})^{-1}$ on the right of \eqref{5.42}. This factor is independent of $c'$ and $k_{v}$ and has global product $1$.
For $0<\Re(\alpha)<1$ and $\Re(\kappa)>3$, absolute convergence permits integration against $\overline{\psi_{v}(c')}\,dc'$.
Substitute \eqref{5.42} for the $(t_{1},t_{2},a)$-integral in
\eqref{e5.41} to get, in the region \eqref{5.41},
\begin{multline}\label{5.44}
J_{\mathrm{degen}}^{\heartsuit}(s,\varphi_{\boldsymbol{\nu}},\omega';w_{1}w_{2}w_{1})=\int_{K'}\chi_{1}\overline{\chi}_{2}(\det k)\int_{\mathbb{A}_{F}}\int_{\mathbb{A}_{F}^{\times}}\Phi_{1\times 2}((b',t)k^{-1})\\
\chi_{2}\overline{\chi}_{3}(t)|t|^{\nu_{2}-\nu_{3}}d^{\times}tdb'\int_{\mathbb{A}_{F}}\mathfrak{F}_{3}\left(w'n(c');k\right)
\overline{\psi(c')}dc'
dk,
\end{multline}
where
\begin{multline*}
\mathfrak{F}_{3}(g';k):=|\det g'|^{\frac{1}{2}}\int_{\mathbb{A}_{F}^{\times}}\int_{\mathbb{A}_{F}^{\times}}\mathcal{F}\Phi_{*}\left(\begin{pmatrix}
 & t_{1} \\
 & t_{2}
\end{pmatrix}g';k\right)
\overline{\chi}_{1}\chi_{3}(t_{2})|t_{2}|^{1-\nu_{1}+\nu_{3}}d^{\times}t_{2}\\
\chi_{1}\chi_{3}\omega'(t_{1})|t_{1}|^{1+\nu_{1}+\nu_{3}+2s}d^{\times}t_{1}.
\end{multline*}

The integral converges and defines a section $\mathfrak{F}_{3}(\cdot;k)$ of
$\mathbf{1}\boxplus\overline{\chi}_{3}^{2}\overline{\omega}'|\cdot|^{-1-2\nu_{3}-2s}$, continuous in $(g',k)$.
For $v\notin S$,
\begin{multline*}
\mathfrak{F}_{3,v}(I_{2};k)=
L_{v}(1-\nu_{1}+\nu_{3},\overline{\chi}_{1}\chi_{3})
L_{v}(1+\nu_{1}+\nu_{3}+2s,\chi_{1}\chi_{3}\omega')\\
L_{v}(1+\nu_{2}+\nu_{3}+2s,\chi_{2}\chi_{3}\omega'),
\end{multline*}
and
\begin{align*}
 \int_{F_{v}}\mathfrak{F}_{3,v}\left(w'n(c');k_{v}\right)
 \overline{\psi}_{v}(c')\,dc'
 =\frac{\mathfrak{F}_{3,v}(I_{2};k_{v})}
 {L_{v}(2+2\nu_{3}+2s,\chi_{3,v}^{2}\omega_{v}')}.
\end{align*}
At the other places, the Tate integrals defining $\mathfrak{F}_{3,v}$ and its Jacquet integral are meromorphic, hence so are the normalized local factors.

To recover the original local integral, multiply by
$\gamma(\alpha,\chi_{1,v}\overline{\chi}_{3,v},\psi_{v})^{-1}$.
For $v\notin S$, this factor is
\begin{align*}
 \gamma(\alpha,\chi_{1,v}\overline{\chi}_{3,v},\psi_{v})^{-1}
 =\frac{L_{v}(\nu_{1}-\nu_{3},\chi_{1,v}\overline{\chi}_{3,v})}
 {L_{v}(1-\nu_{1}+\nu_{3},\overline{\chi}_{1,v}\chi_{3,v})}.
\end{align*}
The first line of \eqref{5.44} contributes
$L_{v}(\nu_{2}-\nu_{3},\chi_{2,v}\overline{\chi}_{3,v})$.
This gives \eqref{f5.15}. The Euler product agrees with the integral
in \eqref{e5.33.} and gives its continuation.
\end{proof}

\subsection{\texorpdfstring{The $w_{1}w_{2}w_{1}$ Term of $J_{\mathrm{degen}}^{\dag}$}{The w1w2w1 Term of J-degen dagger}}\label{sec5.13}
Define
\begin{multline*}
\mathcal{L}^{\dag}(s,\boldsymbol{\nu},\pi,\omega';w_{1}w_{2}w_{1}):=\Lambda(\nu_{1}-\nu_{3},\chi_{1}\overline{\chi}_{3})\Lambda(\nu_{1}-\nu_{2},\chi_{1}\overline{\chi}_{2})\\
\Lambda(1+\nu_{2}+\nu_{3}+2s,\chi_{2}\chi_{3}\omega').
\end{multline*}

Define the normalized local factor
$\mathcal{P}_{v}^{\dag,\sharp}(s,\varphi_{\boldsymbol{\nu}},\omega';w_{1}w_{2}w_{1})$ by
\begin{multline*}
\mathcal{P}_{v}^{\dag,\sharp}(s,\varphi_{\boldsymbol{\nu}},\omega';w_{1}w_{2}w_{1}):=\mathcal{L}_{v}^{\dag}(s,\boldsymbol{\nu},\pi,\omega';w_{1}w_{2}w_{1})^{-1}\int_{F_{v}^{\times}}\int_{F_{v}}\int_{F_{v}}\int_{F_{v}}\\
f_{\boldsymbol{\nu},v}\left(\begin{pmatrix}
1& & \\
c_{v}&1& \\
b_{v}& a_{v}& 1
\end{pmatrix}n_{12}(-z_{v}^{-1})w_{2}w_{1}\right)
\overline{\psi_{v}(a_{v}z_{v})}da_{v}db_{v}dc_{v}\\
\chi_{1,v}\chi_{3,v}\omega_{v}'(z_{v})|z_{v}|_{v}^{\nu_{1}+\nu_{3}+2s}d^{\times}z_{v}.
\end{multline*}

\begin{prop}\label{prop5.12}
$J_{\mathrm{degen}}^{\dag}(s,\varphi_{\boldsymbol{\nu}},\omega';w_{1}w_{2}w_{1})$ converges absolutely in the region 
\begin{equation}\label{5.71}
 \Re(\nu_{1}-\nu_{3})>1,\quad \Re(\nu_{1}-\nu_{2})>1,\quad \Re(2s+\nu_{2}+\nu_{3})>0.
\end{equation}
\end{prop}

\begin{proof}
In the defining unipotent integral make the simultaneous changes
$a\mapsto az$, $c\mapsto z^{-1}c$, with additive Jacobian $1$.
Move $\diag(z,1,z)$ to the left and use \eqref{2.9}. Insert \eqref{e5.3}
with $g'=k\diag(t_{1},t_{2})\left(\begin{smallmatrix}1&\\b'&1\end{smallmatrix}\right)$ and Haar measure \eqref{e1.5}, and make the successive changes
$b'\mapsto-b'$, $b\mapsto b+b'c$, $a\mapsto a+b'+bz^{-1}$,
$b'\mapsto t_{1}b'$, $c\mapsto t_{1}^{-1}c$,
$b\mapsto t_{2}^{-1}b$, $a\mapsto t_{2}^{-1}a$, and
$t_{1}\mapsto z^{-1}t_{1}$. Their total additive Jacobian is
$|t_{2}|^{-2}$. The factors involving $\Phi_{2\times3}$ and $\psi$ become
\begin{align*}
 \Phi_{2\times3}\left(k\begin{pmatrix}c&z^{-1}(t_{1}-c)&\\b&a&t_{2}\end{pmatrix}w_{2}w_{1}\right),
 \qquad \overline{\psi}(at_{2}^{-1}z+bt_{2}^{-1}+b't_{1}).
\end{align*}
To recover the integral over $u_{1}+zu_{2}=t_{1}$ at fixed $t_{1}$,
use the partial Fourier transform
\begin{multline*}
 \widehat{\Phi}_{k}(c,d,\xi,\eta,t_{2})
 :=\int_{\mathbb{A}_{F}^{4}}
 \Phi_{2\times3}\left(k\begin{pmatrix}u_{1}&u_{2}&\\b&a&t_{2}\end{pmatrix}w_{2}w_{1}\right)\\
 \psi(u_{1}c+u_{2}d-a\eta-b\xi)\,du_{1}\,du_{2}\,da\,db,
\end{multline*}
and put
\begin{align*}
&\phi_{k}(x,u):=\int_{\mathbb{A}_{F}}
 \Phi_{1\times2}((r,x)k^{-1})\overline{\psi}(ru)\,dr,\\
 &h(t_{1},z,\xi,\eta,t_{2};k)
:=\int_{\mathbb{A}_{F}}
 \widehat{\Phi}_{k}(c,cz,\xi,\eta,t_{2})\overline{\psi}(ct_{1})\,dc.
\end{align*}
At fixed $u_{1}$, the change from $u_{2}$ to $u_{1}+zu_{2}$ has
Jacobian $|z|$. Fourier inversion in $u_{1}+zu_{2}$ therefore gives
$|z|h(t_{1},z,\xi,\eta,t_{2};k)$ for the integral in $du_{1}$ on
$u_{1}+zu_{2}=t_{1}$.

In the $\Phi_{1\times2}$ integral we replace $b'$ by $t^{-1}r$
at fixed $t$, with Jacobian $|t|^{-1}$. At fixed $t_{2}$ we then
substitute $(t,t_{1})=(t_{2}x,t_{2}xu)$. With
\begin{multline*}
 H(x,u,z;s,\boldsymbol{\nu},k):=\phi_{k}(x,u)
 \int_{\mathbb{A}_{F}^{\times}}
 h(t_{2}xu,z,t_{2}^{-1},t_{2}^{-1}z,t_{2};k)\\
 \chi_{1}^{2}\overline{\chi}_{2}\overline{\chi}_{3}(t_{2})
 |t_{2}|^{2\nu_{1}-\nu_{2}-\nu_{3}}\,d^{\times} t_{2},
\end{multline*}
the integral is
\begin{multline}\label{5.96}
 J_{\mathrm{degen}}^{\dag}(s,\varphi_{\boldsymbol{\nu}},\omega';w_{1}w_{2}w_{1})
 =\int_{K'}\chi_{1}\overline{\chi}_{2}(\det k)
 \int_{(\mathbb{A}_{F}^{\times})^{3}}H(x,u,z;s,\boldsymbol{\nu},k)\\
 \chi_{1}\overline{\chi}_{3}(x)|x|^{\nu_{1}-\nu_{3}}
 \chi_{1}\overline{\chi}_{2}(u)|u|^{\nu_{1}-\nu_{2}}
 \chi_{2}\chi_{3}\omega'(z)|z|^{1+\nu_{2}+\nu_{3}+2s}
 \,d^{\times} x\,d^{\times} u\,d^{\times} z\,dk.
\end{multline}

In addition to the decay of $h$ in $t_{2},t_{2}^{-1}$, use
$(1+|t_{2}|_{v})(1+|t_{2}^{-1}z|_{v})\geq1+|z|_{v}$ to control $z$.
Thus $H$ is Schwartz--Bruhat in $(x,u,z)$ and entire in
$(s,\boldsymbol{\nu})$, locally uniformly in the parameters and $k$.

Thus \eqref{5.96} consists of three Tate integrals, with precisely
the factors in $\mathcal{L}^{\dag}$ and global convergence region
\eqref{5.71}. Division by the three local $L$-factors leaves an
entire family, also after compact integration.

For $v\notin S$, $\phi_{k_{v}}(x,u)=\mathbf{1}_{\mathcal{O}_{v}}(x)
\mathbf{1}_{\mathcal{O}_{v}}(u)$, and Fourier inversion gives
\begin{align*}
 h_{v}(t_{1},z,t_{2}^{-1},t_{2}^{-1}z,t_{2};k_{v})
 =\mathbf{1}_{\mathcal{O}_{v}^{\times}}(t_{2})
  \mathbf{1}_{\mathcal{O}_{v}}(t_{1})\mathbf{1}_{\mathcal{O}_{v}}(z).
\end{align*}
The $t_{2}$-integral is over $\mathcal{O}_{v}^{\times}$, of volume $1$, so
$H_{v}(x,u,z)=\mathbf{1}_{\mathcal{O}_{v}^{3}}(x,u,z)$. Since
$\chi_{1,v}\overline{\chi}_{2,v}(\det k_{v})=1$, this proves \eqref{f5.15}.
\end{proof}

\section{Maximal Eisenstein Data: Degenerate Terms}\label{sec6}
We compute the degenerate terms for Eisenstein series induced from a
cuspidal $\mathrm{GL}_{2}$ representation and a character.

Let $\boldsymbol{\mu}=(\mu_{1},\mu_{2})\in\mathcal{R}_{\mathrm{max}}$; see \eqref{f2.2}. Let $\pi=\chi\boxplus \sigma$ and $\pi_{\boldsymbol{\mu}}=\chi|\cdot|^{\mu_{1}}\boxplus \sigma|\cdot|^{\mu_{2}/2}$, where $\chi$ is a unitary Hecke character, and $\sigma$ is a unitary cuspidal automorphic representation of $\mathrm{GL}_{2}$ over $F$, with central character $\omega_{\sigma}$. %
Put $\omega:=\chi\omega_{\sigma}$, so that the central character of
$\pi_{\boldsymbol{\mu}}$ is $\omega|\cdot|^{\mu_{1}+\mu_{2}}$.

Let $\varphi_{\boldsymbol{\mu}}=E(\cdot,\sigma,\chi,\boldsymbol{\mu})
\in\pi_{\boldsymbol{\mu}}$ be the Eisenstein series \eqref{e2.4}.

\subsection{Godement Sections}\label{sec6.1}
Use the Godement section of \cite[\textsection 7.1]{Yan23}
\begin{multline}\label{6.1}
f(g;\sigma,\chi,\boldsymbol{\mu})=\chi(\det g)|\det g|^{1+\mu_{1}}\int_{G'(\mathbb{A}_{F})}\Phi_{2\times 3}[(\mathbf{0},g')g]\\
\phi_{\mu_{2}}(g'^{-1})\chi(\det g')|\det g'|^{3/2+\mu_{1}}dg',
\end{multline}
where $\Phi_{2\times 3}$ is a Schwartz--Bruhat function on
$M_{2\times 3}(\mathbb{A}_{F})$, and $\phi_{\mu_{2}}$ is an automorphic
vector in $\sigma|\cdot|^{\mu_{2}/2}$.

Changing variables in \eqref{6.1} gives \eqref{e2.3}.
The section \eqref{6.1} admits a meromorphic continuation to
$\boldsymbol{\mu}\in\mathbb{C}^{2}$, with
\begin{align*}
f(g;\sigma,\chi,\boldsymbol{\mu})\propto
\Lambda(1+\mu_{1}-\mu_{2}/2,\sigma\times\chi^{2}\overline{\omega}),
\end{align*}
since $\omega=\chi\omega_{\sigma}$ and
$\widetilde{\sigma}\simeq\sigma\otimes\overline{\omega}_{\sigma}$.

For the inducing form, use
\begin{align*}
 W_{\phi_{\mu_{2}}}(g'):=\int_{F\backslash\mathbb{A}_{F}}
 \phi_{\mu_{2}}(n(c)g')\overline{\psi}(c)\,dc.
\end{align*}
Let $W_{\phi_{\mu_{2}},v}$ be the associated local Whittaker functions. Define
\begin{multline}\label{e6.3}
W_{f_{\boldsymbol{\mu},v}}(g_{v}):=\chi_{v}(\det g_{v})|\det g_{v}|_{v}^{1+\mu_{1}}\int_{G'(F_{v})}\Phi_{2\times 3,v}[(\mathbf{0},g_{v}')g_{v}]\\
W_{\phi_{\mu_{2}},v}(g_{v}'^{-1})\chi_{v}(\det g_{v}')|\det g_{v}'|_{v}^{3/2+\mu_{1}}dg_{v}'.
\end{multline}
Normalize $W_{\phi_{\mu_{2}},v}$ so that
$W_{\phi_{\mu_{2}}}=\prod_{v}W_{\phi_{\mu_{2}},v}$.
Taking the Fourier coefficient in the inducing $\mathrm{GL}_{2}$ variable
of \eqref{6.1} then gives
\begin{align*}
 \int_{F\backslash\mathbb{A}_{F}}
 f_{\boldsymbol{\mu}}\left(n_{23}(c)g\right)\overline{\psi}(c)\,dc
 =\prod_{v}W_{f_{\boldsymbol{\mu},v}}(g_{v}).
\end{align*}

\begin{defn}\label{def:exceptional-places-maximal}
Let $S$ be the finite set such that $v\notin S$ precisely when
$v<\infty$,
$v\nmid\mathfrak{D}_{F}$, the section
$f_{v}(\cdot;\sigma,\chi,\boldsymbol{\mu})$ is spherical,
$\Phi_{2\times 3,v}=\mathbf{1}_{(\mathcal{O}_{v})^{6}}$,
$W_{\phi_{\mu_{2}},v}$ is spherical with
$W_{\phi_{\mu_{2}},v}(I_{2})=1$, and the characters
$\chi_{v}$, $\omega_{v}'$, and $\eta_{v}$ are unramified.
\end{defn}

\subsection{Degenerate Whittaker Functions}
We use the degenerate Whittaker functions of \textsection\ref{sec2.2}
and the Bruhat decomposition relative to $Q(F)$: 
\begin{align*}
G(F)=Q(F)\bigsqcup Q(F)w_{1}N_{w_{1}}(F)\bigsqcup Q(F)w_{1}w_{2}N_{w_{1}w_{2}}(F).
\end{align*}

For simplicity, we abbreviate $f(g;\sigma,\chi,\boldsymbol{\mu})$ by $f_{\boldsymbol{\mu}}(g)$.
For the character $\psi(c)$, the $w_{1}$-cell vanishes by
orthogonality along the unipotent radical of $Q$, and the
$w_{1}w_{2}$-cell contains the constant term of $\phi_{\mu_{2}}$, which
vanishes by cuspidality.  The identity cell therefore gives 
\begin{equation}\label{6.5}
W_{\varphi_{\boldsymbol{\mu}}}^{\mathrm{degen}}(g)=\int_{F\backslash\mathbb{A}_{F}}f_{\boldsymbol{\mu}}\left(n_{23}(c)g\right)\overline{\psi(c)}dc.
\end{equation}

For the opposite coefficient, with character $\psi(a)$, the identity
cell vanishes by orthogonality and the $w_{1}$-cell by cuspidality.
The $w_{1}w_{2}$-cell gives
\begin{equation}\label{e6.6}
W_{\varphi_{\boldsymbol{\mu}}}^{\mathrm{degen},\dag}(g)=\int_{F\backslash\mathbb{A}_{F}}\int_{\mathbb{A}_{F}}\int_{\mathbb{A}_{F}}f_{\boldsymbol{\mu}}\left(w_{1}w_{2}\begin{pmatrix}
1&a& b \\
&1& c\\
&& 1
\end{pmatrix}g\right)dbdc\overline{\psi(a)}da.	
\end{equation}

\subsection{\texorpdfstring{The $I_{\mathrm{degen}}$ Term}{The I-degen Term}}\label{sec6.3}
Define  
\begin{multline*}
\mathcal{L}(\mathbf{s},\boldsymbol{\mu},\pi,\omega',\eta):=\Lambda(2+2\mu_{1}+2s_{1},\chi^{2}\omega')^{-1}\Lambda(1+\mu_{1}+\mu_{2}/2+2s_{1},\sigma\times\chi\omega')\\
\Lambda(1+\mu_{1}-\mu_{2}/2,\sigma\times\chi^{2}\overline{\omega})\Lambda(1+\mu_{1}+s_{1}+s_{2},\chi\eta)\Lambda(s_{2}-\mu_{1}-s_{1},\eta\overline{\chi}\overline{\omega}').
\end{multline*}

Define the normalized local factor $\mathcal{P}_{v}^{\sharp}(\mathbf{s},\varphi_{\boldsymbol{\mu}},\omega',\eta)$ by
\begin{multline*}
\mathcal{P}_{v}^{\sharp}(\mathbf{s},\varphi_{\boldsymbol{\mu}},\omega',\eta):=\frac{1}{\mathcal{L}_{v}(\mathbf{s},\boldsymbol{\mu},\pi,\omega',\eta)}
\int_{(F_{v}^{\times})^{2}}\int_{F_{v}}
W_{f_{\boldsymbol{\mu},v}}\left(w_{1}n_{12}(a_{v})\diag(y_{v}z_{v},z_{v},1)\right)\\
\overline{\psi_{v}(a_{v})}da_{v}\omega_{v}'(z_{v})|z_{v}|_{v}^{2s_{1}}d^{\times}z_{v}\eta_{v}(y_{v})|y_{v}|_{v}^{s_{1}+s_{2}}d^{\times}y_{v}.
\end{multline*}
Here $W_{f_{\boldsymbol{\mu},v}}$ is defined by \eqref{e6.3}; the
additive Fourier integral is evaluated before integration in $y_{v},z_{v}$.
 
\begin{prop}\label{prop6.2}
$I_{\mathrm{degen}}(\mathbf{s},\varphi_{\boldsymbol{\mu}},\omega',\eta)$ converges absolutely in the region 
\begin{equation}\label{6.7}
\begin{cases}
\Re(\mu_{1}+s_{1}+s_{2})>0,\ \ \Re(\mu_{1}-\mu_{2}/2)>0,\\
\Re(s_{2}-s_{1}-\mu_{1})>1,\ \ \Re(\mu_{1}+\mu_{2}/2+2s_{1})>0.
\end{cases}
\end{equation}
\end{prop}

\begin{proof}
Insert \eqref{6.1} and \eqref{6.5} in $I_{\mathrm{degen}}$ and use
$g'=k\diag(t_{1},t_{2})n(b')$ with Haar measure \eqref{e1.4}. Make the simultaneous changes
\begin{equation}\label{eq:maximal-w2-variables}
 (a,b',y,z,t_{1},t_{2})\mapsto
 \left(ya,r^{-1}zb',ty,r^{-1}y^{-1}z,rz^{-1},r\right).
\end{equation}
For fixed $t,y,z,r$, the changes in $a,b'$ give the factor $|yz/r|$.
The map on the other four variables preserves multiplicative Haar measure.
The inducing Whittaker function contributes
$\overline{\omega}_{\sigma}(r)|r|^{-\mu_{2}}\overline{\psi}(b'z/r)$. The local determinant factor is
$\chi_{v}(-1)$, with global product $1$. Put
\begin{align*}
h(t,y,r,u;k):=\int_{\mathbb{A}_{F}^{2}}
\Phi_{2\times3}\left(k\begin{pmatrix}a&t&b'\\&&r\end{pmatrix}w_{1}\right)
\overline{\psi}(ay+b'u)\,da\,db'.
\end{align*}
This is Schwartz--Bruhat in $(t,y,r,u)$, uniformly for $k\in K'$.
Define
\begin{multline}\label{6.15}
H_{I}(r,u;k):=\int_{(\mathbb{A}_{F}^{\times})^{2}}h(t,y,r,u;k)
\chi\eta(t)|t|^{1+\mu_{1}+s_{1}+s_{2}}\\
\eta\overline{\chi}\overline{\omega}'(y)|y|^{s_{2}-\mu_{1}-s_{1}}
\,d^{\times} t\,d^{\times} y,
\end{multline}
and
\begin{equation}\label{6.13}
W(a(z);k):=\chi\omega'(z)|z|^{1/2+\mu_{1}+2s_{1}}
\int_{\mathbb{A}_{F}^{\times}}H_{I}(r,r^{-1}z;k)
\chi\overline{\omega}\overline{\omega}'(r)|r|^{-\mu_{2}-2s_{1}}\,d^{\times} r.
\end{equation}
These changes give
\begin{equation}\label{6.14}
I_{\mathrm{degen}}(\mathbf{s},\varphi_{\boldsymbol{\mu}},\omega',\eta)
=\int_{K'}\chi(\det k)\int_{\mathbb{A}_{F}^{\times}}
W(a(z);k)W_{\phi_{\mu_{2}}}(a(z)k^{-1})\,d^{\times} z\,dk.
\end{equation}

By Tate's thesis, dividing the local integral \eqref{6.15} by the two factors
$L_{v}(1+\mu_{1}+s_{1}+s_{2},\chi_{v}\eta_{v})$ and
$L_{v}(s_{2}-\mu_{1}-s_{1},\eta_{v}\overline{\chi}_{v}\overline{\omega}_{v}')$
gives an entire Schwartz--Bruhat family in $(r,u)$.

By Lemma~\ref{lem5.1} for $\mathcal{F}_{2}H_{I}$, the function $W$ in
\eqref{6.13} is a meromorphic family of Whittaker functions in
$\chi\omega'|\cdot|^{\mu_{1}+2s_{1}}\boxplus\chi^{2}\overline{\omega}|\cdot|^{\mu_{1}-\mu_{2}}$. For $v\notin S$, $h_{v}$ is the characteristic function of
$\mathcal{O}_{v}^{4}$, so $W_{v}$ is the product of the two Tate
factors and the spherical Whittaker function normalized
at $I_{2}$. Equation~\eqref{eq3.6} at $s=0$ evaluates the convolution
\eqref{6.14} with $\sigma|\cdot|^{\mu_{2}/2}$, giving the
$\mathrm{GL}_{2}\times\mathrm{GL}_{2}$ quotient of $L$-factors in $\mathcal{L}_{v}$.
In \eqref{6.7}, the two Tate and two $\mathrm{GL}_{2}$ arguments have
real part greater than one. Whittaker bounds and absolute convergence
of the Euler products give convergence in \eqref{6.7}. The quotients
in \eqref{f5.15} continue meromorphically by local Tate and
Rankin--Selberg theory.
\end{proof}

\subsection{\texorpdfstring{The $J_{\mathrm{degen}}$ Term}{The J-degen Term}}
Define 
\begin{multline*}
\mathcal{L}(s,\boldsymbol{\mu},\pi,\omega'):=\Lambda(2+2\mu_{1}+2s,\chi^{2}\omega')^{-1}\Lambda(1+\mu_{1}-\mu_{2}/2,\sigma\times\chi^{2}\overline{\omega})\\\Lambda(1+\mu_{1}+\mu_{2}/2+2s,\sigma\times \chi\omega').
\end{multline*}

Define the normalized local factor $\mathcal{P}_{v}^{\sharp}(s,\varphi_{\boldsymbol{\mu}},\omega')$ by
\begin{multline*}
\mathcal{P}_{v}^{\sharp}(s,\varphi_{\boldsymbol{\mu}},\omega'):=\mathcal{L}_{v}(s,\boldsymbol{\mu},\pi,\omega')^{-1}\int_{F_{v}^{\times}}\int_{F_{v}}
W_{f_{\boldsymbol{\mu},v}}\left(w_{1}n_{12}(a_{v})\diag(z_{v},z_{v},1)\right)\\
\overline{\psi_{v}(a_{v})}da_{v}\omega_{v}'(z_{v})|z_{v}|_{v}^{2s}d^{\times}z_{v}.
\end{multline*}

\begin{prop}\label{prop6.3}
$J_{\mathrm{degen}}(s,\varphi_{\boldsymbol{\mu}},\omega')$ converges absolutely in the region 
\begin{equation}\label{e6.17}
 \Re(\mu_{1}-\mu_{2}/2)>0,\quad \Re(\mu_{1}+\mu_{2}/2+2s)>0.
\end{equation}
\end{prop}

\begin{proof}
Let $h(t,y,r,u;k)$ be the Schwartz function in the proof
of Proposition~\ref{prop6.2}, and put
\begin{equation}\label{6.18}
H_{J}(r,u;k):=\int_{\mathbb{A}_{F}^{\times}}h(t,t^{-1},r,u;k)
\chi^{2}\omega'(t)|t|^{1+2\mu_{1}+2s}\,d^{\times} t.
\end{equation}
In the Godement expansion of $J_{\mathrm{degen}}$, use
$g'=k\diag(t_{1},t_{2})n(b')$ and set $y=t^{-1}$ on the right of
\eqref{eq:maximal-w2-variables}, so the original $y$ is $1$:
\begin{align*}
 (a,b',z,t_{1},t_{2})\mapsto
 \left(t^{-1}a,r^{-1}zb',tr^{-1}z,rz^{-1},r\right).
\end{align*}
The additive Jacobian is $|z/(tr)|$. The inducing Whittaker function
contributes $\overline{\omega}_{\sigma}(r)|r|^{-\mu_{2}}\overline{\psi}(b'z/r)$.
Thus \eqref{6.13}--\eqref{6.14}
apply to $J_{\mathrm{degen}}$ with $s_{1}=s$ and $H_{I}$ replaced by $H_{J}$.
The rapid decay of $h(t,t^{-1},r,u;k)$ at $t=0$ and infinity, also after derivatives in $(r,u)$, shows that \eqref{6.18} is entire in $(\mu_{1},s)$ and Schwartz--Bruhat in $(r,u)$, uniformly for $k\in K'$.

For $v\notin S$, $H_{J,v}(r,u;k_{v})=\mathbf{1}_{\mathcal{O}_{v}}(r)\mathbf{1}_{\mathcal{O}_{v}}(u)$, since $t,t^{-1}\in\mathcal{O}_{v}$ force $t\in\mathcal{O}_{v}^{\times}$.
Lemma~\ref{lem5.1} therefore gives a Kirillov vector in $\chi\omega'|\cdot|^{\mu_{1}+2s}
\boxplus\chi^{2}\overline{\omega}|\cdot|^{\mu_{1}-\mu_{2}}$, spherical with value $1$ at $I_{2}$ for $v\notin S$.
Applying \eqref{eq3.6} at $s=0$ to the convolution with $\sigma|\cdot|^{\mu_{2}/2}$ gives the three factors defining $\mathcal{L}(s,\boldsymbol{\mu},\pi,\omega')$.
Whittaker bounds and absolute Euler-product convergence give
\eqref{e6.17}. Local Rankin--Selberg continuation gives the
meromorphic identity \eqref{f5.15}.
\end{proof}

\subsection{\texorpdfstring{The $J_{\mathrm{degen}}^{\dag}$ Term}{The J-degen dagger Term}}
Put
\begin{align*}
\mathcal{L}^{\dag}(s,\boldsymbol{\mu},\pi,\omega'):=\Lambda(1+\mu_{2}+2s,\overline{\chi}\omega\omega')\Lambda(\mu_{1}-\mu_{2}/2,\sigma\times \chi^{2}\overline{\omega}).
\end{align*}

Define the normalized local factor by
\begin{multline}\label{e6.23}
\mathcal{P}_{v}^{\dag,\sharp}(s,\varphi_{\boldsymbol{\mu}},\omega'):=\mathcal{L}_{v}^{\dag}(s,\boldsymbol{\mu},\pi,\omega')^{-1}\int_{F_{v}^{\times}}\chi_{v}\omega_{v}'(z_{v})|z_{v}|_{v}^{1+\mu_{1}+2s}\int_{F_{v}}\int_{F_{v}}\\
W_{f_{\boldsymbol{\mu},v}}\left(\diag(1,z_{v},1)\begin{pmatrix}
1&&-z_{v}^{-1}\\
b_{v}&1&-z_{v}^{-1}b_{v}\\
z_{v}c_{v}&&1-c_{v}
\end{pmatrix}w_{1}\right)\,db_{v}\,dc_{v}\,d^{\times}z_{v}.
\end{multline}

\begin{prop}\label{prop6.4}
$J_{\mathrm{degen}}^{\dag}(s,\varphi_{\boldsymbol{\mu}},\omega')$ converges absolutely in the region 
\begin{equation}\label{6.22}
 \Re(\mu_{2}+2s)>0,\quad \Re(\mu_{1}-\mu_{2}/2)>1.
\end{equation}
\end{prop}

\begin{proof}
In \eqref{e6.6}, make the change $b\mapsto b+ac$ and move $\diag(z,1,1)$ to the left. Take the Fourier coefficient in $a$ and divide the local integral
by $\mathcal{L}_{v}^{\dag}$ as in \eqref{e6.23}.
Insert \eqref{e6.3}, use $g'=k\diag(t_{1},t_{2})n(b')$ and Haar measure \eqref{e1.4}, and make the successive changes $b\mapsto b-b'c$, $b'\mapsto b'+b$, followed by the simultaneous changes
\begin{align*}
 b\mapsto(t_{1}z)^{-1}b,\qquad b'\mapsto t_{1}^{-1}b',\qquad
 c\mapsto(t_{2}z)^{-1}c.
\end{align*}
The two shears have determinant one. The three rescalings give
$|t_{1}^{-2}t_{2}^{-1}z^{-2}|$ in $db\,db'\,dc$. Next change
$t_{1}\mapsto z^{-1}t_{1}$ in the multiplicative integral.
The central scalar $t_{2}^{-1}$ in the inducing Whittaker function gives $\overline{\omega}_{\sigma}(t_{2})|t_{2}|^{-\mu_{2}}$.
Fourier inversion along $u_{1}+zu_{2}=zt_{2}$ contributes $|z|$.
The change $t_{2}\mapsto z^{-1}t_{2}$ gives the following integral,
since $\prod_{v}\chi_{v}(-1)=1$:
\begin{multline*}
J_{\mathrm{degen}}^{\dag}(s,\varphi_{\boldsymbol{\mu}},\omega')
=\int_{K'}\chi(\det k)\int_{(\mathbb{A}_{F}^{\times})^{3}}
W_{\phi_{\mu_{2}}}(a(t_{2}t_{1}^{-1})k^{-1})h(t_{1},z,t_{2};k)\\
\chi^{2}\overline{\omega}(t_{2})|t_{2}|^{-1/2+\mu_{1}-\mu_{2}}
\overline{\chi}\omega\omega'(z)|z|^{1+\mu_{2}+2s}
\chi(t_{1})|t_{1}|^{1/2+\mu_{1}}
\,d^{\times} t_{2}\,d^{\times} z\,d^{\times} t_{1}\,dk,
\end{multline*}
where
\begin{align*}
h(t,z,u;k):=\int_{\mathbb{A}_{F}}
\widehat{\Phi}_{k}(t;t^{-1},t^{-1}z,c,cz)\overline{\psi}(cu)\,dc.
\end{align*}
Here $\widehat{\Phi}_{k}(t;\xi,\zeta,\eta,\eta')$ is the partial Fourier
transform of
\begin{align*}
(b,b',u_{1},u_{2})\longmapsto
\Phi_{2\times3}\left(k\begin{pmatrix}b&t&b'\\u_{1}&&u_{2}\end{pmatrix}w_{1}\right)
\end{align*}
with character $\overline{\psi}(b\xi+b'\zeta)\psi(u_{1}\eta+u_{2}\eta')$.
It is Schwartz--Bruhat in all five variables, uniformly for $k\in K'$.

At Archimedean places, integration by parts in $c$ and the inequality
$1+|z|_{v}\leq(1+|t|_{v})(1+|t^{-1}z|_{v})$ give, for every $l$,
\begin{align*}
 |h_{v}(t,z,u;k_{v})|\ll_{l}
 (1+|t|_{v}+|t|_{v}^{-1}+|z|_{v}+|u|_{v})^{-l}.
\end{align*}
The rapid decay of $\widehat{\Phi}_{k_{v}}$ absorbs the powers of
$t^{-1},z,c$ from differentiation, so all fixed $(t,z,u)$-derivatives
also decay rapidly, uniformly in $k_{v}\in K_{v}'$. Extend $h_{v}$ by zero at
$t=0$.

At finite places, $\widehat{\Phi}_{k_{v}}$ is Schwartz--Bruhat and
$|t|_{v}\gg1$ on the support, so $h_{v}$ is Schwartz--Bruhat in $(t,z,u)$.
For $v\notin S$,
\begin{equation}\label{6.29}
h_{v}(t,z,u;k_{v})=\mathbf{1}_{\mathcal{O}_{v}^{\times}}(t)
\mathbf{1}_{\mathcal{O}_{v}}(z)\mathbf{1}_{\mathcal{O}_{v}}(u).
\end{equation}
Set $t_{2}=t_{1}y$. Since
$1+|y|_{v}\leq(1+|t_{2}|_{v})(1+|t_{1}|_{v}^{-1})$, the rapid decay of $h_{v}$
and its derivatives gives a Schwartz--Bruhat function of $(z,y)$ after
integration in $t_{1}$, entire in
$(\mu_{1},\mu_{2})$, with local factor $\mathbf{1}_{\mathcal{O}_{v}^{2}}(z,y)$ for $v\notin S$. Its $z$-integration is a Tate integral with factor $L_{v}(1+\mu_{2}+2s,\overline{\chi}_{v}\omega_{v}\omega_{v}')$. The $y$-integration against $W_{\phi_{\mu_{2}}}$ is a
$\mathrm{GL}_{2}\times\mathrm{GL}_{1}$ integral with factor $L_{v}(\mu_{1}-\mu_{2}/2,\sigma_{v}\times\chi_{v}^{2}\overline{\omega}_{v})$. Tate and Rankin--Selberg theory give convergence in \eqref{6.22}.
Dividing the local $z$-integral by its $L$-factor gives a Schwartz
function of $y$. Multiply $W_{\phi_{\mu_{2}}}$ in the Kirillov model
by this function and apply the normalized Rankin--Selberg integral
to get an entire quotient by $\mathcal{L}_{v}^{\dag}$. For $v\notin S$, the spherical $z$- and
$y$-integrals just evaluated give $\mathcal{L}_{v}^{\dag}$ before
normalization, so $\mathcal{P}_{v}^{\dag,\heartsuit}=1$.
\end{proof}

\section{Meromorphic Continuation of the Dual Side}\label{sect5}
We establish the meromorphic continuation of $J_{\mathrm{sing}}$,
$J_{\mathrm{dual}}$, and $I_{\mathrm{gen}}$ needed in \textsection\ref{sec8}.

\subsection{\texorpdfstring{Meromorphic Continuation of $J_{\mathrm{sing}}(s,\varphi,\omega')$ }{Meromorphic Continuation of J\_sing}}
Local Rankin--Selberg theory gives:
\begin{lemma}\label{lem7.1}
Let $\pi$ be generic automorphic on $[G]$ with central character
$\omega$, and let $\varphi\in\pi$. The singular term
$J_{\mathrm{sing}}(s,\varphi,\omega')$ converges absolutely for
$\Re(s)\ggg1$ and extends meromorphically to $\mathbb{C}$, with
\begin{align*}
J_{\mathrm{sing}}^{\heartsuit}(s,\varphi,\omega')\propto \Lambda(1+2s,\widetilde{\pi}\times\omega\omega').
\end{align*}
More precisely, each local zeta integral divided by its local
$L$-factor is entire in $s$.
For $v<\infty$ with $v\nmid\mathfrak{D}_{F}$, if $\pi_{v}$ and $\omega_{v}'$
are unramified and $\varphi$ is right-$K_{v}$-invariant, then
\begin{align*}
J_{\mathrm{sing},v}^{\heartsuit}(s,\varphi,\omega')=W_{\varphi,v}(I_{3})L_{v}(1+2s,\widetilde{\pi}\times\omega\omega').	
\end{align*}
\end{lemma}

\subsection{\texorpdfstring{Meromorphic Continuation of $J_{\mathrm{dual}}(s,\varphi,\omega')$ }{Meromorphic Continuation of J\_dual}}\label{sec7.2}
First take $W_{\varphi^{\iota}}$ factorizable; the estimates below give the general case by linearity and continuity.
For $\Re(\lambda)=s_{0}\ggg1$ and $\Re(s)-s_{0}\ggg1$, the proof
of Lemma~\ref{lem3.3} gives absolute convergence and the factorization
\begin{align*}
\mathcal{J}(s,\lambda;\varphi,\xi,\omega')=\prod_{v}\mathcal{J}_{v}(s,\lambda;\varphi,\xi,\omega'),
\end{align*}
where  
\begin{multline}\label{f7.4}
\mathcal{J}_{v}(s,\lambda;\varphi,\xi,\omega'):=\int_{F_{v}^{\times}}\int_{F_{v}^{\times}}W_{\varphi^{\iota},v}\left(\diag(z_{v},1,1)n_{21}(y_{v})\widetilde{w}\right)\\
\overline{\xi}_{v}\omega_{v}\omega_{v}'(z_{v})|z_{v}|_{v}^{2s-\lambda-\frac{1}{2}}\xi_{v}(y_{v})|y_{v}|_{v}^{\lambda+\frac{1}{2}}d^{\times}z_{v}d^{\times}y_{v}.
\end{multline}

Write $\widetilde{\mathcal{J}}$ for the meromorphic continuation of $\mathcal{J}$ constructed in \textsection\ref{sec4.3.1} below.
\begin{prop}\label{prop4.5}
Let $\pi$ be one of the unitary data in \textsection\ref{sec4.1}, with central
character $\omega$, let $\omega'$ be unitary, and let $\varphi\in\pi$.
The integral $J_{\mathrm{dual}}(s,\varphi,\omega')$ extends
meromorphically to $\Re(s)>-1/4$. For $|\Re(s)|<1/4$,
\begin{equation}\label{e4.11}
\widetilde{J}_{\mathrm{dual}}(s,\varphi,\omega')=J_{\mathrm{dual}}^{\heartsuit}(s,\varphi,\omega')+R_{\mathrm{dual}}^{+}(s,\varphi,\omega')-R_{\mathrm{dual}}^{-}(s,\varphi,\omega'),
\end{equation}
where
\begin{align*}
&J_{\mathrm{dual}}^{\heartsuit}(s,\varphi,\omega'):=\frac{1}{2\pi i\rho_{F}}\sum_{\xi\in \widehat{C_{F}}}\int_{i\mathbb{R}}\widetilde{\mathcal{J}}(s,\lambda;\varphi,\xi,\omega')d\lambda,\\
&R_{\mathrm{dual}}^{+}(s,\varphi,\omega'):=\rho_{F}^{-1}\underset{\lambda=1/2}{\Res}\widetilde{\mathcal{J}}(s,\lambda;\varphi,\mathbf{1},\omega'),\\
&R_{\mathrm{dual}}^{-}(s,\varphi,\omega'):=\rho_{F}^{-1}\sum_{\xi\in \mathfrak{X}_{\pi,\omega'}^{\mathrm{dual}}}\underset{\lambda=2s-1/2}{\Res}\widetilde{\mathcal{J}}(s,\lambda;\varphi,\xi,\omega').
\end{align*}
\end{prop}
See \textsection\ref{sec:proof-dual-continuation} for the proof.

\subsubsection{Meromorphic continuation of $\mathcal{J}(s,\lambda;\varphi,\xi,\omega')$}\label{sec4.3.1}

\begin{lemma}\label{lem7.3}
Suppose $v<\infty$, $v\nmid \mathfrak{D}_{F}$, the characters
$\xi_{v}$ and $\omega_{v}'$ are unramified, and
$W_{\varphi^{\iota},v}$ is spherical. Then
\begin{equation}\label{4.10}
\mathcal{J}_{v}(s,\lambda;\varphi,\xi,\omega')=W_{\varphi^{\iota},v}(I_{3})L_{v}(1/2+2s-\lambda,\widetilde{\pi}\times \overline{\xi}\omega\omega')L_{v}(1/2+\lambda,\xi).
\end{equation}

At a general place $v$, there is an entire function
$\mathcal{Q}_{v}(s,\lambda;\varphi,\xi,\omega')$ such that
\begin{equation}\label{4.11}
\mathcal{J}_{v}(s,\lambda;\varphi,\xi,\omega')
=\mathcal{Q}_{v}(s,\lambda;\varphi,\xi,\omega')L_{v}(1/2+2s-\lambda,\widetilde{\pi}\times
\overline{\xi}\omega\omega')L_{v}(1/2+\lambda,\xi).
\end{equation}
\end{lemma}

\begin{proof}
The $z_{v}$-integral divided by its $L$-factor is entire in $(s,\lambda)$
and continuous in $W_{\varphi^{\iota},v}$. The identity
\begin{align*}
 \diag(z_{v},1,1)n_{21}(y_{v})\widetilde{w}n_{31}(t)
 =n_{13}(z_{v}t)n_{23}(y_{v}t)\diag(z_{v},1,1)n_{21}(y_{v})\widetilde{w}
\end{align*}
makes right translation by $n_{31}(t)$ act on the integrand by
$\overline{\psi}_{v}(y_{v}t)$. Differentiation in $t$ and moderate growth
at the Archimedean places, and smoothness at the finite places, give
Schwartz--Bruhat bounds in $y_{v}$ for the normalized $z_{v}$-integral,
locally uniformly in $(s,\lambda)$. Apply Tate's theorem for
\eqref{4.11}, and use \eqref{3.18}, Tate's formula and
Casselman--Shalika for \eqref{4.10}.
\end{proof}

\phantomsection\label{def:continued-dual-integral}
Lemma~\ref{lem7.3} gives a meromorphic continuation $\widetilde{\mathcal{J}}(s,\lambda;\varphi,\xi,\omega')$ of $\mathcal{J}(s,\lambda;\varphi,\xi,\omega')$ to $(s,\lambda)\in\mathbb{C}^{2}$, with
\begin{equation}\label{e4.10}
\widetilde{\mathcal{J}}(s,\lambda;\varphi,\xi,\omega')\propto \Lambda(1/2+2s-\lambda,\widetilde{\pi}\times \overline{\xi}\omega\omega')\Lambda(1/2+\lambda,\xi).
\end{equation}
Outside a fixed finite set of places, \eqref{3.18} makes the local integral zero if $\xi_{v}$ is ramified; otherwise \eqref{4.10} applies. Thus the product of the nontrivial $\mathcal{Q}_{v}$ is finite and entire, and all poles come from the two global $L$-factors in \eqref{e4.10}.

\subsubsection{Proof of Proposition \ref{prop4.5}}\label{sec:proof-dual-continuation}
For $\Re(\lambda)\ggg1$ and $\Re(s)-\Re(\lambda)\ggg1$, right translation by $n_{23}(t)$ shifts $y_{v}$ by $t$ in \eqref{f7.4}. Together with the $n_{31}(t)$ identity in Lemma~\ref{lem7.3}, integration by parts in $y_{v}$ gives, for $m\geq1$, finitely many $X_{m,j}$ in the enveloping algebra of $G(F_{\infty})$ such that
\begin{equation}\label{eq7.5}
\big|\widetilde{\mathcal{J}}(s,\lambda;\varphi,\xi,\omega')\big|
\ll_{m} C_{\infty}(\xi|\cdot|^{\lambda+1/2})^{-m}
\sum_{j}\big|\widetilde{\mathcal{J}}
(s,\lambda;\pi(X_{m,j})\varphi,\xi,\omega')\big|.
\end{equation}
The underlying integration-by-parts identities continue meromorphically.

At finite places, smoothness gives a compact open subgroup $D_{\fin}^{*}\subseteq\prod_{v<\infty}\mathcal{O}_{v}^{\times}$ such that
\begin{equation}\label{eq7.6}
\widetilde{\mathcal{J}}(s,\lambda;\varphi,\xi,\omega')=\widetilde{\mathcal{J}}(s,\lambda;\varphi,\xi,\omega')\mathbf{1}_{\text{$\xi$ is right-$D_{\fin}^{*}$-invariant}}.
\end{equation}

Choose $S$ to contain the Archimedean places, the places where the
fixed data are ramified, and the finite places where
$D_{\fin}^{*}$ does not contain $\mathcal{O}_{v}^{\times}$. By \eqref{eq7.6},
$\xi_{v}$ is unramified for $v\notin S$ in every nonzero term.
At such $v$, \eqref{3.18} replaces the Whittaker factor by its
value at $y_{v}=0$ times $\mathbf{1}_{\mathcal{O}_{v}}(y_{v})$. Integrating
in $y_{v}$ contributes $L_{v}(1/2+\lambda,\xi_{v})$.

The remaining adelic $z$-integral, with $y_{v}=0$ for $v\notin S$,
is a $\mathrm{GL}_{3}\times\mathrm{GL}_{1}$ integral at $1/2+2s-\lambda$
with character $\overline{\xi}\omega\omega'$. For $\Re(1/2+2s-\lambda)$ sufficiently large, bound the integral
by absolute convergence. For $\Re(1/2+2s-\lambda)$ sufficiently negative,
apply the global functional equation and bound the dual integral by
absolute convergence. Both absolute bounds are independent of the unitary character $\xi$. After removing the finitely many global poles by a fixed polynomial, apply the maximum principle on the strip. The resulting polynomial growth is uniform in $\xi$, with order independent of $X_{m,j}$ for fixed $\pi$ and a fixed strip.

Moderate growth under right translation and the $n_{31}(t)$ identity above then give rapid decay in the Archimedean $y_{v}$, with the same polynomial growth in the parameters. At finite places, smoothness gives bounded support in $y_{v}$. The remaining $y_{v}$-integrals, for $v\in S$, therefore converge absolutely when $\Re(\lambda)>-1/2$. The polynomial bounds for $L^{S}(1/2+\lambda,\xi)$, together with \eqref{eq7.5}--\eqref{eq7.6}, now give locally uniform absolute convergence away from poles and justify the contour shifts.

For $\Re(s)\ggg1$, shifting from $\Re(\lambda)=s_{0}$ to $0$ crosses only $(\xi,\lambda)=(\mathbf{1},1/2)$ and gives
$J_{\mathrm{dual}}=J_{\mathrm{dual}}^{\heartsuit}+R_{\mathrm{dual}}^{+}$,
with residue coefficient $\rho_{F}^{-1}$. As $\Re(s)$ decreases through $1/4$, the poles at $\lambda=2s-1/2$ cross in the opposite direction, giving $-R_{\mathrm{dual}}^{-}$. The inducing characters are those of \textsection\ref{sec4.1}. The other possible poles of the two global $L$-factors, at $\lambda=-1/2$ and $2s+1/2$, are not crossed for $\Re(s)>-1/4$, and \eqref{e4.10} excludes other poles. Near a crossing, keep the local $\lambda$-contour fixed on the side
that avoids the moving poles, then return to $i\mathbb{R}$ and add the
crossed residues. This continues the integral to $\Re(s)>-1/4$ and
gives \eqref{e4.11} for $|\Re(s)|<1/4$.

\subsection{\texorpdfstring{Meromorphic Continuation of $I_{\mathrm{gen}}(\mathbf{s},\varphi,\omega',\eta)$}{Meromorphic Continuation of I\_gen}}\label{sec7.3}
Let $\pi$ be a generic automorphic representation of $[G]$ with
central character $\omega$. We continue the generic integral defined
in Theorem~\ref{thme3.3}.

\begin{prop}\label{prop5.4.}
Suppose $v<\infty$, $v\nmid\mathfrak{D}_{F}$, $\omega_{v}'$ and
$\eta_{v}$ are unramified, and $\varphi$ is right-$K_{v}$-invariant. Then
\begin{multline}\label{e7.5}
I_{\mathrm{gen},v}(\mathbf{s},\varphi,\omega',\eta)=W_{\varphi,v}(I_{3})L_{v}(2+3s_{1}+s_{2},\omega_{v}\omega_{v}'\eta_{v})^{-1}\\
L_{v}(1+s_{1}+s_{2},\pi_{v}\times\eta_{v})L_{v}(1+2s_{1},\widetilde{\pi}_{v}\times\omega_{v}\omega_{v}').
\end{multline}

$I_{\mathrm{gen}}(\mathbf{s},\varphi,\omega',\eta)$ converges
absolutely for $\Re(s_{1}),\Re(s_{2})\ggg1$ and has a meromorphic continuation
$I_{\mathrm{gen}}^{\heartsuit}(\mathbf{s},\varphi,\omega',\eta)$ to $\mathbb{C}^{2}$.

\end{prop}

\begin{proof}
By linearity, take $W_{\varphi}=\prod_{v}W_{v}$. The local integral
$I_{\mathrm{gen},v}(\mathbf{s},\varphi,\omega',\eta)$ at $v\in\Sigma_{F}$ is
\begin{equation}\label{7.5}
\int_{F_{v}^{\times}}\int_{F_{v}^{\times}}W_{v}\left(\diag(y_{v}z_{v},z_{v},1)\right)\omega_{v}'(z_{v})|z_{v}|_{v}^{2s_{1}}\eta_{v}(y_{v})|y_{v}|_{v}^{s_{1}+s_{2}}d^{\times}z_{v}d^{\times}y_{v}.
\end{equation}

For the minimal Godement section in \textsection\ref{sec5.1}, the associated Whittaker function is the Jacquet integral
\begin{equation}\label{f7.11}
W_{v}(g)=\int_{F_{v}^{3}}f_{v}(\widetilde{w} n(a,b,c)g)
\overline{\psi}_{v}(a+c)\,da\,db\,dc.
\end{equation}
Here $n(a,b,c)$ is the upper unipotent matrix in \eqref{w2.5}.
We treat the following local types.
\begin{itemize}
\item Suppose $\pi_{v}=\chi_{1,v}|\cdot|^{\nu_{1}}\boxplus \chi_{2,v}|\cdot|^{\nu_{2}}\boxplus\chi_{3,v}|\cdot|^{\nu_{3}}$ is a principal series. Use the local Godement section \eqref{e5.3}.

Substitute the local version of \eqref{e5.3} and \eqref{f7.11} into \eqref{7.5} and make the simultaneous changes
\begin{align*}
 b_{v}'\mapsto  -b_{v}',\qquad
 a_{v}\mapsto a_{v}+b_{v}',\qquad
 b_{v}\mapsto t_{2,v}^{-1}b_{v}+b_{v}'c_{v}.
\end{align*}
The additive Jacobian is $|t_{2,v}|_{v}^{-1}$.
The determinant contributes $\chi_{1,v}(-1)$ because
$\det\widetilde{w}=-1$, and the compact integration contributes
$\chi_{1,v}\overline{\chi}_{2,v}(\det k_{v})$ from \eqref{e5.3}.
Now make the substitution
\begin{align*}
 (c_{v},a_{v},b_{v}',t_{v})\mapsto(c_{v}/t_{1,v},a_{v}/(t_{2,v}z_{v}),(t_{1,v}/t_{2,v})b_{v}',t_{2,v}t_{v}).
\end{align*}
The additional additive Jacobian is $|t_{2,v}^{-2}z_{v}^{-1}|_{v}$.
Use $t_{1,v}^{-1}$, $(t_{2,v}z_{v})^{-1}$, $t_{1,v}/t_{2,v}$, and $t_{2,v}y_{v}z_{v}$ as the
new variables in place of $(t_{1,v},t_{2,v},z_{v},y_{v})$, preserving the product
of multiplicative measures. The character in $b_{v}'$ is now
$\overline{\psi}_{v}(b_{v}'z_{v})$, so the $(b_{v}',t_{v})$-integral gives
$W_{2,v}$ from \eqref{fc5.22} times the torus factor
$\chi_{3,v}(z_{v})|z_{v}|_{v}^{1/2+\nu_{3}}$.
Thus
\begin{multline*}
I_{\mathrm{gen},v}(\mathbf{s},\varphi,\omega',\eta)
=\chi_{1,v}(-1)\int_{K_{v}'}
\chi_{1,v}\overline{\chi}_{2,v}(\det k_{v})
\int_{F_{v}^{\times}}W_{2,v}(a(z_{v});k_{v})\chi_{1,v}\omega_{v}'(z_{v})\\
|z_{v}|_{v}^{1/2+\nu_{1}+2s_{1}}
\int_{F_{v}^{\times}}h(t_{2,v},t_{2,v}^{-1}z_{v};k_{v})
\overline{\chi}_{1,v}\overline{\omega}_{v}'\eta_{v}(t_{2,v})|t_{2,v}|_{v}^{s_{2}-s_{1}-\nu_{1}}
\,d^{\times} t_{2,v}\,d^{\times} z_{v}\,dk_{v},
\end{multline*}
where $h(t_{2,v},z_{v};k_{v})$ is defined by 
\begin{multline*}
\int_{F_{v}}
\int_{F_{v}^{\times}}\int_{F_{v}^{\times}}\int_{F_{v}}\int_{F_{v}}
\Phi_{2\times 3,v}\left(k\begin{pmatrix}
c_{v} & z_{v} &  \\
b_{v} & a_{v} & y_{v}
\end{pmatrix}w_{1}w_{2}w_{1}\right)\overline{\psi_{v}(c_{v}t_{1,v})}dc_{v}db_{v}\\
\chi_{2,v}\chi_{3,v}\omega_{v}'(t_{1,v})|t_{1,v}|_{v}^{1+\nu_{2}+\nu_{3}+2s_{1}}d^{\times}t_{1,v}\\
\chi_{1,v}\eta_{v}(y_{v})|y_{v}|_{v}^{1+\nu_{1}+s_{1}+s_{2}}d^{\times}y_{v}
\overline{\psi_{v}(a_{v}t_{2,v})}da_{v}.
\end{multline*}

The $t_{1,v}$- and $y_{v}$-integrals contribute
$L_{v}(1+\nu_{2}+\nu_{3}+2s_{1},\chi_{2,v}\chi_{3,v}\omega_{v}')$
and $L_{v}(1+\nu_{1}+s_{1}+s_{2},\chi_{1,v}\eta_{v})$ in their common
region of absolute convergence. After division by these factors, $h$ is entire in the parameters and Schwartz--Bruhat in $(t_{2,v},z_{v})$, with locally uniform derivative and decay bounds at Archimedean places and fixed compact support at finite places.
By Lemma~\ref{lem5.1}, the integral in $(t_{2,v},z_{v})$ is
\eqref{f3.5} at $s=0$ for Kirillov vectors in
\begin{align*}
\chi_{1,v}\omega_{v}'|\cdot|_{v}^{\nu_{1}+2s_{1}}
 \boxplus\eta_{v}|\cdot|_{v}^{s_{1}+s_{2}},\qquad
\chi_{2,v}|\cdot|_{v}^{\nu_{2}}\boxplus\chi_{3,v}|\cdot|_{v}^{\nu_{3}}.
\end{align*}
For $\Re(s_{1}),\Re(s_{2})\ggg1$, Whittaker bounds at zero give absolute convergence;
Lemma~\ref{lem3.2} gives meromorphic continuation. For $v\notin S$, $h$ equals the product of its two
Tate factors and
$\mathbf{1}_{\mathcal{O}_{v}}(t_{2,v})\mathbf{1}_{\mathcal{O}_{v}}(z_{v})$.
Both determinant characters are $1$ there, so \eqref{eq3.6} gives
\eqref{e7.5}.

\item Suppose $\pi_{v}=\chi_{v}|\cdot|_{v}^{\mu_{1}}\boxplus \sigma_{v}|\cdot|_{v}^{\mu_{2}/2}$, where $\sigma_{v}$ is an irreducible generic unitary representation of $G'(F_{v})$. Use the local Godement section \eqref{6.1}, with inducing vector
$\phi_{\mu_{2},v}\in\sigma_{v}|\cdot|_{v}^{\mu_{2}/2}$. 
Unfolding the open $Q$-cell and using
$w_{1}w_{2}n_{12}(a)=n_{23}(a)w_{1}w_{2}$ to take the inducing Fourier
coefficient gives
\begin{align*}
 W_{v}(g)=\int_{F_{v}^{2}}W_{f_{\boldsymbol{\mu},v}}\left(w_{1}w_{2}
 n_{13}(b_{v})n_{23}(c_{v})g\right)
 \overline{\psi}_{v}(c_{v})\,db_{v}\,dc_{v}.
\end{align*}
Insert \eqref{e6.3}, use $g_{v}'=k_{v}\diag(t_{1,v},t_{2,v})n(b_{v}')$ and Haar measure \eqref{e1.4}, and make the change $b_{v}\mapsto b_{v}-b_{v}'c_{v}$.
Next make the simultaneous changes
\begin{multline*}
 (b_{v},c_{v},b_{v}',y_{v},z_{v},t_{1,v},t_{2,v})\\
 \mapsto\left(z_{v}y_{v}b_{v},z_{v}c_{v},t_{2,v}^{-1}y_{v}b_{v}',
 t_{1,v}t_{2,v}^{-1}y_{v},t_{2,v}z_{v},z_{v}^{-1}y_{v}^{-1},z_{v}^{-1}\right).
\end{multline*}
The additive Jacobian is $|z_{v}^{2}y_{v}^{2}/t_{2,v}|_{v}$, and the multiplicative measures are unchanged. The factor $\omega_{\sigma_{v}}(z_{v})|z_{v}|_{v}^{\mu_{2}}\overline{\psi}_{v}(b_{v}'y_{v}/t_{2,v})$ comes from the inducing Whittaker function.

Let $t_{2,v}, y_{v}'\in F_{v}$.  In the initial domain of convergence, define
the Schwartz--Bruhat function
\begin{multline*}
h_{v}(t_{2,v},y_{v}';k_{v}):=\int_{(F_{v}^{\times})^{2}}\int_{(F_{v})^{3}}\Phi_{2\times 3,v}\left(k_{v}\begin{pmatrix}
b_{v} &t_{1,v}& b_{v}'\\
c_{v}&& t_{2,v}
\end{pmatrix}w_{1}w_{2}\right)\overline{\psi_{v}(b_{v}'y_{v}')}\\
db_{v}'db_{v}\overline{\psi_{v}(c_{v}z_{v})}dc_{v}
\omega_{\sigma_{v}}\omega_{v}'(z_{v})|z_{v}|_{v}^{1+\mu_{2}+2s_{1}}d^{\times}z_{v}\\
\chi_{v}\eta_{v}(t_{1,v})|t_{1,v}|_{v}^{1+\mu_{1}+s_{1}+s_{2}}d^{\times}t_{1,v}.
\end{multline*}

By Tate's thesis, $h_{v}(t_{2,v},y_{v}';k_{v})$ converges absolutely in the region 
\begin{align*}
 \Re(\mu_{2}+2s_{1})>0,\quad \Re(\mu_{1}+s_{1}+s_{2})>0
\end{align*}
and the $z_{v}$- and $t_{1,v}$-integrals contribute the Tate factors
$L_{v}(1+\mu_{2}+2s_{1},\omega_{\sigma_{v}}\omega_{v}')$ and
$L_{v}(1+\mu_{1}+s_{1}+s_{2},\chi_{v}\eta_{v})$. After division by these factors, $h_{v}$ is entire in $(\mathbf{s},\boldsymbol{\mu})$ and Schwartz--Bruhat in $(t_{2,v},y_{v}')$, with the same locally uniform bounds.

These changes give
\begin{multline}\label{e7.18}
I_{\mathrm{gen},v}(\mathbf{s},\varphi,\omega',\eta)=\int_{K_{v}'}\chi_{v}(\det k_{v})\int_{F_{v}^{\times}}W_{\phi_{\mu_{2}},v}\left(a(y_{v})
k_{v}^{-1}\right)\\
W_{3,v}'\left(a(y_{v});
k_{v}\right)|y_{v}|_{v}^{s_{1}+s_{2}}d^{\times}y_{v}
dk_{v},
\end{multline}
where $W_{3,v}'\left(a(y_{v});
k_{v}\right)$ is defined by 
\begin{equation}\label{7.18}
\eta_{v}(y_{v})|y_{v}|_{v}^{\frac{1}{2}}\int_{F_{v}^{\times}}h_{v}(t_{2,v},t_{2,v}^{-1}y_{v};k_{v})\omega_{v}'\chi_{v}\overline{\eta}_{v}(t_{2,v})|t_{2,v}|_{v}^{\mu_{1}+s_{1}-s_{2}}d^{\times}t_{2,v}.
\end{equation}

By Lemma~\ref{lem5.1}, \eqref{7.18} gives a Kirillov vector in
$\pi_{3,v}'=\eta_{v}\boxplus\chi_{v}\omega_{v}'
|\cdot|_{v}^{\mu_{1}+s_{1}-s_{2}}$, depending meromorphically on $(\mathbf{s},\boldsymbol{\mu})$. Lemma~\ref{lem3.2} applied to
\eqref{e7.18} gives absolute convergence for
$\Re(s_{1}),\Re(s_{2})\ggg1$ and meromorphic continuation to
$(\mathbf{s},\boldsymbol{\mu})\in\mathbb{C}^{4}$.
For $v\notin S$, \eqref{eq3.6} gives
\begin{align*}
I_{\mathrm{gen},v}(\mathbf{s},\varphi,\omega',\eta)=\frac{W_{\phi_{\mu_{2}},v}(I_{2})W_{3,v}'(I_{2};I_{2})L_{v}(1+\mu_{2}/2+s_{1}+s_{2},\pi_{3,v}'\times \sigma_{v})}{L_{v}(2+3s_{1}+s_{2},\omega_{v}\omega_{v}'\eta_{v})},
\end{align*}
which amounts to \eqref{e7.5}. 

For $v\in S$, the normalized local factors are meromorphic.

\item Suppose $v<\infty$ and $\pi_{v}$ is supercuspidal. Its Whittaker
functions are compactly supported modulo $N(F_{v})Z(F_{v})$.
Thus $W_{v}(\diag(y_{v}z_{v},z_{v},1))$ is compactly supported as a function of
$(y_{v},z_{v})\in(F_{v}^{\times})^{2}$.

The integral \eqref{7.5} therefore converges absolutely for all
$(s_{1},s_{2})\in\mathbb{C}^{2}$ and is entire.
\end{itemize}

Steinberg twists are treated as subrepresentations of principal series.
Multiply the Euler product of \eqref{e7.5} over $v\notin S$ by
the meromorphic local integrals at $v\in S$ to continue the global
integral. Whittaker bounds and absolute convergence of the Euler product
give convergence for $\Re(s_{1}),\Re(s_{2})\ggg1$.
\end{proof}

\subsection{Integral Formulas for the Residues}
Let $\pi$ and $\varphi$ be as in Proposition~\ref{prop4.1}, and let
$\omega'$ and $\eta$ be unitary Hecke characters. We take the
Rankin--Selberg residue of each period and then apply Parseval to
the sum over the Eisenstein basis. The resulting integrals
involve only $W_{\varphi}$ and the degenerate Whittaker functions
of $\varphi$.

Throughout this subsection, integrals are continued meromorphically
from the stated regions of absolute convergence, with inner
integrals evaluated before outer ones. We form the residue sums
defined in \textsection\ref{sec4} and \textsection\ref{sec7.2} at
distinct simple poles before specializing to coincident poles.
For Eisenstein data, we take the constant term
in the inducing parameters as in \eqref{eq:inducing-constant-term}.
At the central parameters, we
use the constant-term contours of \textsection\ref{sec8}.

\subsubsection{The Rankin--Selberg period}
For a unitary Hecke character $\xi$ and $k\in K'$, put
\begin{multline}\label{eq:residue-initial-section}
 \mathfrak{F}(s,\lambda,\xi;k)
 :=\int_{(\mathbb{A}_{F}^{\times})^{2}}\int_{\mathbb{A}_{F}}
 W_{\varphi}(\diag(w'n(u)\diag(a,d)k,1))\\
 \xi(a)\overline{\xi}\omega'(d)
 |a|^{s+\lambda-1/2}|d|^{s-\lambda+1/2}
 \,du\,d^{\times}a\,d^{\times}d.
\end{multline}
A sufficient region of absolute convergence for unitary data is
\begin{align*}
\Re(\lambda)>1,\qquad \Re(s)>\Re(\lambda)+2.
\end{align*}
Write $\mathfrak{F}_{v}$ for \eqref{eq:residue-initial-section} over $F_{v}$
with $W_{\varphi}$ replaced by $W_{v}$.
A sufficient local region is $\Re(s)>\Re(\lambda)>0$.
Extend $\mathfrak{F}$ to
$G'(\mathbb{A}_{F})$ by
\begin{align*}
 \mathfrak{F}\left(s,\lambda,\xi;
 \begin{pmatrix}a&b\\&d\end{pmatrix}k\right)
 =\overline{\xi}(a)\xi\overline{\omega'}(d)
 |a/d|^{1/2-\lambda}\mathfrak{F}(s,\lambda,\xi;k),
\end{align*}
and likewise locally. The period identity below gives its
meromorphic continuation; away from the poles it is smooth in $k$.

Let $f_{\lambda,h}$ be the section of
$\xi|\cdot|^{\lambda}\boxplus\overline{\xi}\omega'|\cdot|^{-\lambda}$
with $f_{\lambda,h}|_{K'}=h$. Unfolding the Jacquet integral of
$E(\cdot,h,\lambda)$ and using Iwasawa decomposition gives
\begin{align*}
 &\Psi(1/2+s,W_{\varphi},W_{E(\cdot,h,\lambda)})
 =\int_{G'(\mathbb{A}_{F})}W_{\varphi}(\diag(w'g,1))
 f_{\lambda,h}(g)|\det g|^{s}\,dg,\\
 &\Psi(1/2+s,W_{\varphi},W_{E(\cdot,h,\lambda)})
 =\int_{K'}\mathfrak{F}(s,\lambda,\xi;k)h(k)\,dk.
\end{align*}

\subsubsection{Integrals in the original automorphic data}
Put
\begin{equation}\label{eq:residue-original-B}
 \Psi(g';s):=|\det g'|^{s}
 \int_{\mathbb{A}_{F}^{\times}}
 W_{\varphi}^{\mathrm{degen},\dag}
 (w_{2}\diag(ag',1))\omega'(a)|a|^{2s}d^{\times}a.
\end{equation}
For unitary data this converges absolutely for $\Re(s)>1$.
For Eisenstein data, replace $\varphi$ in
\eqref{eq:residue-original-B} by the Eisenstein series
$\varphi_{\boldsymbol{\mu}}$ or $\varphi_{\boldsymbol{\nu}}$ of
\textsection\ref{sec2}. Sufficient regions of absolute convergence are
\begin{align*}
 &\Re(2s+\mu_{1}+\mu_{2}/2)>2,\qquad\text{for maximal data},\\
 &\min_{j<k}\Re(2s+\nu_{j}+\nu_{k})>2,
 \qquad\text{for minimal data}.
\end{align*}
$\Psi(g';s)=0$ when $\varphi$ is cuspidal. For
$W_{\varphi^{\iota}}^{\mathrm{degen},\dag}$ below, use the definition
in \textsection\ref{sec2.2} with $\psi(a)$ in place of
$\overline{\psi}(a)$.

\begin{prop}\label{prop:residues-initial-data}
We have 
\begin{align*}
 &R_{\RNum{1}}^{+}(s,\varphi,\omega')
 =\frac{1}{2}\int_{\mathbb{A}_{F}}
 \Psi(w'n(u);s)\overline{\psi}(u)\,du,\\
 &R_{\RNum{2}}^{1,+}(\mathbf{s},\varphi,\omega',\eta)
 =\frac{1}{2}\int_{\mathbb{A}_{F}^{\times}}
 \left(\int_{\mathbb{A}_{F}}
 \Psi(w'n(u)a(y);s_{1})\overline{\psi}(u)\,du\right)
 \eta(y)|y|^{s_{2}}d^{\times}y.
\end{align*}
The $u$-integral defining $R_{\RNum{1}}^{+}$ converges absolutely for $\Re(s)<-1$.
For $R_{\RNum{2}}^{1,+}$, take $\Re(s_{1})<-1$ and $\Re(s_{1}+s_{2})>3$.
For cuspidal and Eisenstein data,
\begin{multline}\label{eq:residue-global-typeII-second}
 R_{\RNum{2}}^{2,+}(\mathbf{s},\varphi,\omega',\eta)
 =\frac{1}{2}\int_{(\mathbb{A}_{F}^{\times})^{2}}\int_{\mathbb{A}_{F}}
 W_{\varphi}(\diag(w'n(u)\diag(a,d)w',1))\\
 \overline{\eta}\omega'(a)\eta(d)
 |a|^{s_{1}-s_{2}}|d|^{s_{1}+s_{2}}\,du\,d^{\times}a\,d^{\times}d,
\end{multline}
initially for $\Re(s_{2})<-1$ and $\Re(s_{1}+s_{2})>4$, and
\begin{equation}\label{eq:residue-global-dual-positive}
 R_{\mathrm{dual}}^{+}(s,\varphi,\omega')
 =\int_{\mathbb{A}_{F}^{\times}}\int_{\mathbb{A}_{F}}
 W_{\varphi^{\iota}}(\diag(z,1,1)n_{21}(y)\widetilde{w})\omega\omega'(z)|z|^{2s-1}\,dy\,d^{\times}z,
\end{equation}
initially for $\Re(s)>2$. Finally,
\begin{multline*}
 R_{\mathrm{dual}}^{-}(s,\varphi,\omega')
 =-\int_{\mathbb{A}_{F}^{\times}}\int_{\mathbb{A}_{F}}
 W_{\varphi^{\iota}}^{\mathrm{degen},\dag}
 (w_{2}w_{1}n_{13}(c)\diag(y,1,1)n_{21}(y)\widetilde{w})\\
 \omega\omega'(y)|y|^{2s}\,dc\,d^{\times}y,
\end{multline*}
initially for $\Re(s)>2$. The other spectral
residues satisfy
\begin{equation}\label{eq:opposite-spectral-residues}
 R_{\RNum{1}}^{-}=-R_{\RNum{1}}^{+},\qquad
 R_{\RNum{2}}^{j,-}=-R_{\RNum{2}}^{j,+},\quad j=1,2.
\end{equation}
\end{prop}

\begin{proof}
For minimal data, write
$W_{\varphi_{\boldsymbol{\nu}}}^{\mathrm{degen},\dag}(g;i)$ for the
three terms of \eqref{f2.13} with Weyl elements
$w_{1}w_{2}w_{1},w_{2}w_{1},w_{1}$, respectively. Thus
\begin{align*}
 &W_{\varphi_{\boldsymbol{\nu}}}^{\mathrm{degen},\dag}(g)
 =\sum_{i=1}^{3}W_{\varphi_{\boldsymbol{\nu}}}^{\mathrm{degen},\dag}(g;i),\\
 &W_{\varphi_{\boldsymbol{\nu}}}^{\mathrm{degen},\dag}
 (\diag(1,1,t)g;i)
 =\chi_{i}(t)|t|^{-1+\nu_{i}}
 W_{\varphi_{\boldsymbol{\nu}}}^{\mathrm{degen},\dag}(g;i).
\end{align*}
For maximal data there is only $i=1$: use the coefficient
in \eqref{e6.6}, with $(\chi_{1},\nu_{1})=(\chi,\mu_{1})$.
Put
\begin{align*}
 \Psi_{i}(g';s):=|\det g'|^{s}
 \int_{\mathbb{A}_{F}^{\times}}
 W_{\varphi}^{\mathrm{degen},\dag}
 (w_{2}\diag(ag',1);i)\omega'(a)|a|^{2s}d^{\times}a.
\end{align*}
Then $\Psi=\sum_{i}\Psi_{i}$. The identity
$\diag(a,1,a)w_{2}=w_{2}\diag(a,a,1)$ and homogeneity give
\begin{align*}
 \Psi_{i}\left(\begin{pmatrix}a&b\\&d\end{pmatrix}g';s\right)
 =\overline{\chi}_{i}\,\overline{\omega'}(a)\chi_{i}(d)
 |a/d|^{1-s-\nu_{i}}\Psi_{i}(g';s).
\end{align*}

Unfolding the Eisenstein series with $\chi_{i}|\cdot|^{\nu_{i}}$ in
the last diagonal position gives
\begin{align*}
 W_{\varphi}(g)=\int_{\mathbb{A}_{F}^{2}}
 W_{\varphi}^{\mathrm{degen},\dag}
 (w_{2}w_{1}n_{12}(b)n_{13}(c)g;i)
 \overline{\psi}(b)\,db\,dc.
\end{align*}
Use the Eisenstein functional equation for each $i$ in the minimal
case and unfold \eqref{e2.4} in the maximal case.
The matrix identity
\begin{align*}
 w_{2}w_{1}n_{12}(b)n_{13}(c)\diag(t,1,1)
 =\diag(1,1,t)w_{2}w_{1}n_{12}(b/t)n_{13}(c/t)
\end{align*}
and the additive Jacobian $|t|^{2}$ give
\begin{multline*}
 \int_{\mathbb{A}_{F}^{\times}}
 W_{\varphi}(\diag(t,1,1)g)\overline{\chi}_{i}(t)|t|^{r-1}d^{\times}t
 =\int_{\mathbb{A}_{F}^{\times}}\int_{\mathbb{A}_{F}^{2}}\\
 W_{\varphi}^{\mathrm{degen},\dag}
 (w_{2}w_{1}n_{12}(b)n_{13}(c)g;i)
 \overline{\psi}(tb)|t|^{r+\nu_{i}}\,db\,dc\,d^{\times}t.
\end{multline*}
Take the Fourier transform in $b$ before the Tate residue at
$r+\nu_{i}=1$. The residue gives $\rho_{F}$ times the integral of the
transform against the additive measure $dt$, which by Fourier
inversion is its original function at $b=0$. Thus
\begin{multline}\label{eq:residue-standard-period}
 \underset{r=1-\nu_{i}}{\Res}
 \int_{\mathbb{A}_{F}^{\times}}
 W_{\varphi}(\diag(t,1,1)g)
 \overline{\chi}_{i}(t)|t|^{r-1}d^{\times}t\\
 =\rho_{F}\int_{\mathbb{A}_{F}}
 W_{\varphi}^{\mathrm{degen},\dag}(w_{2}w_{1}n_{13}(c)g;i)\,dc.
\end{multline}
The $t$-integral initially converges when
$\min_{j}\Re(r+\nu_{j})>2$ for minimal data, and when both
$\Re(r+\mu_{1})>2$ and $\Re(r+\mu_{2}/2)>2$ for maximal data.

The $c$-integral in \eqref{eq:residue-standard-period} converges
absolutely by the Whittaker bounds. For
$\min_{j\ne i}\Re(\nu_{j}-\nu_{i})\ggg1$ in the minimal case or
$\Re(\mu_{2}/2-\mu_{1})\ggg1$ in the maximal case, these bounds are
integrable in $(b,c)$, uniformly on a small circle around
$r=1-\nu_{i}$ after Tate continuation, so we may take the residue
under the integrals.

To calculate the residue of \eqref{eq:residue-initial-section},
put $v=ud/a$ and use
\begin{align*}
w'n(u)\diag(a,d) =\diag(d,a)w'n(v),\qquad
 du=|a/d|\,dv.
\end{align*}
At $\xi=\chi_{i}\omega'$, apply
\eqref{eq:residue-standard-period} with $r=1/2+s-\lambda$ to the
$d$-integral in
\eqref{eq:residue-initial-section}. The change from $r$ to
$\lambda$ contributes a minus sign. To evaluate the remaining
integrals, use the defining character of
$W_{\varphi}^{\mathrm{degen},\dag}$:
\begin{multline*}
 W_{\varphi}^{\mathrm{degen},\dag}
 \left(\begin{pmatrix}a&av&\\&&1\\&1&c\end{pmatrix}
 \diag(k,1);i\right)\\
 =\overline{\psi}(avc)
 W_{\varphi}^{\mathrm{degen},\dag}
 \left(\begin{pmatrix}a&&\\&&1\\&1&c\end{pmatrix}
 \diag(k,1);i\right).
\end{multline*}
Put $u=av$ in the outer integral, so $dv=|a|^{-1}du$, and
apply Fourier inversion at $c=0$:
\begin{multline*}
 \int_{\mathbb{A}_{F}}
 \left(\int_{\mathbb{A}_{F}}
 W_{\varphi}^{\mathrm{degen},\dag}
 \left(\begin{pmatrix}a&av&\\&&1\\&1&c\end{pmatrix}
 \diag(k,1);i\right)\,dc\right)dv\\
 =|a|^{-1}W_{\varphi}^{\mathrm{degen},\dag}
 (\diag(a,1,1)w_{2}\diag(k,1);i).
\end{multline*}
For fixed $a$, decay in $c$ and its derivatives justifies taking
the $c$-integral first and applying Fourier inversion. The
$a$-integral converges absolutely for $\Re(s)$ sufficiently large.
Using
$\diag(a,1,a)w_{2}=w_{2}\diag(a,a,1)$ and the transformation law of
$W_{\varphi}^{\mathrm{degen},\dag}(\cdot;i)$, the remaining
$a$-integral is $\Psi_{i}(k;s)$. Hence
\begin{align*}
 &\underset{\lambda=s-1/2+\nu_{i}}{\Res}
 \mathfrak{F}(s,\lambda,\chi_{i}\omega';k)
 =-\rho_{F}\Psi_{i}(k;s),\\
 &\underset{\lambda=s-1/2+\nu_{i}}{\Res}
 \widetilde{\Psi}(1/2+s,W_{\varphi},W_{E(\cdot,h,\lambda)})
 =-\rho_{F}\int_{K'}\Psi_{i}(k;s)h(k)\,dk.
\end{align*}
Parseval with the norm \eqref{eq:inducing-global-norm} gives
\begin{align*}
\sum_{h\in\mathfrak{B}(\xi,\overline{\xi}\omega')}
 \left(\rho_{F}\int_{K'}\Psi_{i}(k;s)h(k)\,dk\right)
 \overline{h(k')}
 =\Psi_{i}(k';s),\qquad \xi=\chi_{i}\omega'.
\end{align*}
Integration by parts on the Archimedean factors of $K'$ gives
rapid decay of the basis coefficients, uniformly on a small
circle around $\lambda=s-1/2+\nu_{i}$.
The series and its derivatives in $k'$ converge uniformly, so
we may take the residue term by term. Extend both sides by
their Borel transformation laws. Integration at $w'n(u)$
against $\overline{\psi}(u)\,du$ gives
\begin{align*}
 R_{\RNum{1}}^{-}
 =-\frac{1}{2}\sum_{i}\int_{\mathbb{A}_{F}}
 \Psi_{i}(w'n(u);s)\overline{\psi}(u)\,du.
\end{align*}
The Eisenstein functional equation exchanges $(\xi,\lambda)$ with
$(\overline{\xi}\omega',-\lambda)$. Since $\lambda\mapsto-\lambda$
reverses the residue sign, this proves
\eqref{eq:opposite-spectral-residues}. Together with
$\sum_{i}\Psi_{i}=\Psi$, it gives the formula for $R_{\RNum{1}}^{+}$.
For $R_{\RNum{2}}^{1,+}$, evaluate the Parseval sum at
$w'n(u)a(y)$ and take the $u$-Fourier coefficient before the
$y$-Mellin integral.

Parseval in the Eisenstein basis gives
\begin{align*}
 &F_{\xi}(s,\lambda)
 =\rho_{F}^{-1}\int_{\mathbb{A}_{F}}
 \mathfrak{F}(s,\lambda,\xi;w'n(u))\overline{\psi}(u)\,du,\\
 &\widetilde{\mathcal{G}}(\mathbf{s},\lambda,\xi)
 =\rho_{F}^{-1}\int_{\mathbb{A}_{F}^{\times}}
 \left(\int_{\mathbb{A}_{F}}
 \mathfrak{F}(s_{1},\lambda,\xi;w'n(u)a(y))
 \overline{\psi}(u)\,du\right)\eta(y)|y|^{s_{2}}d^{\times}y.
\end{align*}
The $u$-integral converges absolutely for $\Re(\lambda)<-1$, and
the $y$-integral for $\Re(s_{2})>2+|\Re(\lambda)|$. These regions and
the uniform bounds for the basis coefficients justify termwise
integration.

In the expression for $\widetilde{\mathcal{G}}$, put $u=yt$.
The $y$-integral becomes
\begin{align*}
 \int_{\mathbb{A}_{F}^{\times}}
 \left(\int_{\mathbb{A}_{F}}
 \mathfrak{F}(s_{1},\lambda,\xi;w'n(t))\overline{\psi}(yt)\,dt\right)
 \xi\overline{\omega'}\eta(y)
 |y|^{s_{2}+\lambda+1/2}d^{\times}y.
\end{align*}
At $\xi=\overline{\eta}\omega'$ and
$\lambda=1/2-s_{2}$, its Tate residue is
$\rho_{F}\mathfrak{F}(s_{1},1/2-s_{2},\overline{\eta}\omega';w')$.
For $\Re(\lambda)$ sufficiently negative,
$\mathfrak{F}(s_{1},\lambda,\xi;w'n(t))$ and sufficiently many of
its Archimedean derivatives in $t$ are integrable. This justifies
taking the Tate residue. Multiplying the residue
by $1/(2\rho_{F})$ gives \eqref{eq:residue-global-typeII-second}.

In the adelic integral $\mathcal{J}$ obtained from the local
integrals in \eqref{f7.4}, the Tate residue of the $y$-integral at
$\lambda=1/2$ gives \eqref{eq:residue-global-dual-positive}.
For $R_{\mathrm{dual}}^{-}$, apply
\eqref{eq:residue-standard-period} to the contragredient inducing
family, with $\chi_{i}$ replaced by $\overline{\chi}_{i}$ and $\psi$ by
$\overline{\psi}$. For minimal data, the matrix identity and
homogeneity give
\begin{align*}
 &w_{2}w_{1}n_{13}(c)\diag(y,1,1)
 =\diag(1,1,y)w_{2}w_{1}n_{13}(c/y),\\
 &W_{\varphi_{\boldsymbol{\nu}}^{\iota}}^{\mathrm{degen},\dag}
 (\diag(1,1,y)g;i)
 =\overline{\chi}_{i}(y)|y|^{-1-\nu_{i}}
 W_{\varphi_{\boldsymbol{\nu}}^{\iota}}^{\mathrm{degen},\dag}(g;i).
\end{align*}
Substituting $c\mapsto yc$ therefore gives
\begin{multline*}
 \int_{\mathbb{A}_{F}}
 W_{\varphi_{\boldsymbol{\nu}}^{\iota}}^{\mathrm{degen},\dag}
 (w_{2}w_{1}n_{13}(c)\diag(y,1,1)g;i)\,dc\\
 =\overline{\chi}_{i}(y)|y|^{-\nu_{i}}\int_{\mathbb{A}_{F}}
 W_{\varphi_{\boldsymbol{\nu}}^{\iota}}^{\mathrm{degen},\dag}
 (w_{2}w_{1}n_{13}(c)g;i)\,dc.
\end{multline*}
The central character is
$\omega_{\boldsymbol{\nu}}=\chi_{1}\chi_{2}\chi_{3}
|\cdot|^{\nu_{1}+\nu_{2}+\nu_{3}}$. For
$\xi=\overline{\chi}_{i}\chi_{1}\chi_{2}\chi_{3}\omega'$, the $z$-integral
has Mellin parameter
$r=1/2+2s-\lambda+\nu_{1}+\nu_{2}+\nu_{3}$.
Its $i$-th residue is at $r=1+\nu_{i}$, and $dr/d\lambda=-1$.
The resulting $y$-weight is
$\overline{\chi}_{i}\omega_{\boldsymbol{\nu}}\omega'(y)
|y|^{2s-\nu_{i}}$.
Take $g=n_{21}(y)\widetilde{w}$ in the $c$-integral and apply
\eqref{eq:inducing-constant-term} to the sum over $i$ to obtain
$R_{\mathrm{dual}}^{-}$.
For maximal data, replace $\nu_{i}$ by $\mu_{1}$ and
$\nu_{1}+\nu_{2}+\nu_{3}$ by $\mu_{1}+\mu_{2}$.
Lemma~\ref{lem7.3} gives the decay needed for the remaining
$y$-integral when $\Re(s)>2$.
\end{proof}

\subsubsection{Euler products}
Suppose $W_{\varphi}=\prod_{v}W_{v}$. Choose a finite set $S$
containing the Archimedean places, the places dividing
$\mathfrak{D}_{F}$, those where $\xi$, $\omega'$, or $\eta$ is
ramified, and those where $W_{v}$ is not spherical. Normalize
$W_{v}(I_{3})=1$ for $v\notin S$. The local factors are
\begin{multline*}
\mathcal{Q}_{\RNum{1},v}(s,\lambda,\xi_{v})
 :=\frac{L_{v}(1+2\lambda,\xi_{v}^{2}\overline{\omega_{v}'})
 L_{v}(1-2\lambda,\overline{\xi}_{v}^{\,2}\omega_{v}')}
 {L_{v}(1/2+s+\lambda,\pi_{v}\times\xi_{v})
 L_{v}(1/2+s-\lambda,\pi_{v}\times\overline{\xi}_{v}\omega_{v}')}\\
\int_{F_{v}}
 \mathfrak{F}_{v}(s,\lambda,\xi_{v};w'n(u))\overline{\psi}_{v}(u)\,du,
\end{multline*}
\begin{multline*}
\mathcal{Q}_{\RNum{2},v}(\mathbf{s},\lambda,\xi_{v})
 :=\frac{L_{v}(1+2\lambda,\xi_{v}^{2}\overline{\omega_{v}'})
 L_{v}(1-2\lambda,\overline{\xi}_{v}^{\,2}\omega_{v}')}
 {L_{v}(1/2+s_{1}+\lambda,\pi_{v}\times\xi_{v})
 L_{v}(1/2+s_{1}-\lambda,\pi_{v}\times\overline{\xi}_{v}\omega_{v}')}\\
\frac{1}
 {L_{v}(1/2+s_{2}-\lambda,\overline{\xi}_{v}\eta_{v})
 L_{v}(1/2+s_{2}+\lambda,\xi_{v}\overline{\omega_{v}'}\eta_{v})}\\
\int_{F_{v}^{\times}}\left(\int_{F_{v}}
 \mathfrak{F}_{v}(s_{1},\lambda,\xi_{v};w'n(u)a(y))
 \overline{\psi}_{v}(u)\,du\right)\eta_{v}(y)|y|_{v}^{s_{2}}d^{\times}y.
\end{multline*}

The local $u$-integral converges for $\Re(\lambda)<-1$ and the
$y$-integral for $\Re(s_{2})>2+|\Re(\lambda)|$.
The spherical Rankin--Selberg, Jacquet, and Mellin
formulas give $\mathcal{Q}_{\RNum{1},v}=\mathcal{Q}_{\RNum{2},v}=1$
for $v\notin S$. Thus
\begin{equation}\label{eq:initial-data-typeI-density}
 F_{\xi}(s,\lambda)
 =\frac{\Lambda(1/2+s+\lambda,\pi\times\xi)
 \Lambda(1/2+s-\lambda,\pi\times\overline{\xi}\omega')}
 {\rho_{F}\Lambda(1+2\lambda,\xi^{2}\overline{\omega'})
 \Lambda(1-2\lambda,\overline{\xi}^{\,2}\omega')}
 \prod_{v\in S}\mathcal{Q}_{\RNum{1},v}(s,\lambda,\xi_{v}).
\end{equation}
For type~\RNum{2},
\begin{multline*}
 \widetilde{\mathcal{G}}(\mathbf{s},\lambda,\xi)
 =\frac{\Lambda(1/2+s_{1}+\lambda,\pi\times\xi)
 \Lambda(1/2+s_{1}-\lambda,\pi\times\overline{\xi}\omega')}
 {\rho_{F}\Lambda(1+2\lambda,\xi^{2}\overline{\omega'})
 \Lambda(1-2\lambda,\overline{\xi}^{\,2}\omega')}\\
 \Lambda(1/2+s_{2}-\lambda,\overline{\xi}\eta)
 \Lambda(1/2+s_{2}+\lambda,\xi\overline{\omega'}\eta)
 \prod_{v\in S}\mathcal{Q}_{\RNum{2},v}(\mathbf{s},\lambda,\xi_{v}).
\end{multline*}

For maximal data $\pi=\sigma\boxplus\chi$, the integral formulas
therefore give
\begin{equation}\label{eq:initial-data-typeI-residue}
 R_{\RNum{1}}^{+}
 =\frac{\Lambda(1,\sigma\times\overline{\chi})
 \Lambda(2s,\sigma\times\chi\omega')}
 {2\Lambda(2-2s,\overline{\chi}^{\,2}\overline{\omega'})}
 \prod_{v\in S}\mathcal{Q}_{\RNum{1},v}(s,1/2-s,\overline{\chi}_{v}),
\end{equation}
and 
\begin{multline*}
 R_{\RNum{2}}^{1,+}
 =\frac{\Lambda(1,\sigma\times\overline{\chi})
 \Lambda(2s_{1},\sigma\times\chi\omega')}
 {2\Lambda(2-2s_{1},\overline{\chi}^{\,2}\overline{\omega'})}\Lambda(s_{1}+s_{2},\chi\eta)
 \Lambda(1+s_{2}-s_{1},\overline{\chi}\overline{\omega'}\eta)\\
 \prod_{v\in S}\mathcal{Q}_{\RNum{2},v}
 (\mathbf{s},1/2-s_{1},\overline{\chi}_{v}).
\end{multline*}
For minimal data with distinct poles, sum over $i=1,2,3$ with
$\chi=\chi_{i}$ and $\sigma=\boxplus_{j\ne i}\chi_{j}$.
For cuspidal and Eisenstein data, we have 
\begin{align*}
 R_{\RNum{2}}^{2,+}
 =\frac{\Lambda(1+s_{1}-s_{2},\pi\times\overline{\eta}\omega')
 \Lambda(s_{1}+s_{2},\pi\times\eta)}
 {2\Lambda(2-2s_{2},\overline{\eta}^{\,2}\omega')}\prod_{v\in S}\mathcal{Q}_{\RNum{2},v}
 (\mathbf{s},1/2-s_{2},\overline{\eta}_{v}\omega_{v}'),
\end{align*}
and, with $\mathcal{Q}_{v}$ as in \eqref{4.11},
\begin{equation}\label{eq:initial-data-positive-dual-residue}
 R_{\mathrm{dual}}^{+}
 =\Lambda(2s,\widetilde{\pi}\times\omega\omega')
 \prod_{v\in S}\mathcal{Q}_{v}(s,1/2;\varphi,\mathbf{1},\omega').
\end{equation}
For maximal data,
\begin{equation}\label{eq:initial-data-negative-dual-residue}
 R_{\mathrm{dual}}^{-}
 =-\Lambda(1,\widetilde{\sigma}\times\chi)
 \Lambda(2s,\overline{\chi}\omega\omega')\prod_{v\in S}\mathcal{Q}_{v}
 (s,2s-1/2;\varphi,\overline{\chi}\omega\omega',\omega').
\end{equation}
For minimal data with distinct poles, again sum over $i=1,2,3$
with $\chi=\chi_{i}$ and $\sigma=\boxplus_{j\ne i}\chi_{j}$.
For $R_{\mathrm{dual}}^{+}$, take the residue of
$\Lambda(1/2+\lambda,\xi)$ at $\xi=\mathbf{1}$, $\lambda=1/2$.
For $R_{\mathrm{dual}}^{-}$, the pole is in
$\Lambda(1/2+2s-\lambda,\widetilde{\pi}\times\overline{\xi}\omega\omega')$
in the Euler product from \eqref{4.11}. The Tate residue is $\rho_{F}$
in its own variable, so the coefficient at $\lambda=2s-1/2$ is
$-\rho_{F}$.

\section{The Exact Spectral Reciprocity Formulas}\label{sec8}

Let $\omega, \omega'$ and $\eta$ be unitary Hecke characters of $F^{\times}\backslash\mathbb{A}_{F}^{\times}$. Let $\pi$ be a \textit{unitary} generic automorphic representation of $[G]$ with central character $\omega$. For Eisenstein data, use the unitary isobaric sums in \textsection\ref{sec4.1}. Let $\varphi\in\pi$.

\subsection{The Exact Spectral Reciprocity Formula of Type \RNum{1}}\label{sec8.1}
For Eisenstein data, first take the constant terms in the
inducing parameters for $\Re(s)\ggg1$, as in
\textsection\ref{sec3.3}, and then continue in $s$.
The degenerate terms are the sums over the Weyl elements in
\textsection\ref{sec5.2.2}, continued by the propositions of
\textsection\ref{sec5}; for maximal Eisenstein data use
Propositions~\ref{prop6.3} and~\ref{prop6.4}. They vanish for
cuspidal data. The spectral, singular, and dual terms are continued by
Proposition~\ref{prop4.1}, Lemma~\ref{lem7.1}, and
Proposition~\ref{prop4.5}, respectively.

\typei*
\begin{proof}
In \eqref{e3.3}, use Proposition~\ref{prop4.1} for $J_{\mathrm{spec}}$
and \eqref{e4.11} for $J_{\mathrm{dual}}$. Move
$R_{\RNum{1}}^{+}-R_{\RNum{1}}^{-}$ to the right and use
\eqref{eq:opposite-spectral-residues} to get \eqref{f8.1}.
\end{proof}

For arithmetic applications, we specialize \eqref{f8.1} at $s=0$.
The two spectral terms are holomorphic there by Proposition~\ref{prop4.1} and the unitary assumptions.
Since individual terms on the right of \eqref{f8.1} may have poles there, we take
their constant terms on one common circle
$|s|=\varepsilon$, sufficiently small to enclose no other pole. Define
\begin{align*}
&J_{*}^{\heartsuit}(\varphi,\omega'):=J_{*}^{\heartsuit}(s,\varphi,\omega')|_{s=0},\ \ *\in\big\{\mathrm{cusp}, \mathrm{Eis}\big\},\\
&h(\varphi,\omega'):=\frac{1}{2\pi i}\int_{|s|=\varepsilon}
\frac{h(s,\varphi,\omega')}{s}\,ds,
\quad\text{$h$ any term on the right of \eqref{f8.1}}.
\end{align*}

Then \eqref{f8.1} holds with the $s$-arguments omitted.

\subsection{The Exact Spectral Reciprocity Formula of Type \RNum{2}}\label{sec8.2}
Use \eqref{e3.8} for unitary data as in \textsection\ref{sec3.4.2},
where the Eisenstein constant terms are taken in the inducing parameters.
Propositions~\ref{prop4.4} and~\ref{prop5.4.} continue the spectral
and generic terms in $\mathbf{s}$, respectively.

\typeii*
\begin{proof}
Continue \eqref{e3.8} and insert
\eqref{4.7} on both spectral sides. The residues from the transformed
side enter with a plus sign, and those from the original side with
a minus sign. Using \eqref{eq:opposite-spectral-residues} gives
\eqref{8.1} for $|\Re(s_{i})|,|\Re(s_{i}^{\vee})|<1/2$, $i=1,2$.
\end{proof}

To specialize \eqref{8.1} at $\mathbf{s}=(0,0)$, choose
$0<r_{2}\ll r_{1}$ sufficiently small so that no term in \eqref{8.1}
has a pole on the product of $|s_{1}|=r_{1}$ and $|s_{2}|=r_{2}$. Define
\begin{align*}
\operatorname{CT}_{0} H
:=-\frac{1}{4\pi^{2}}
\int_{|s_{1}|=r_{1}}
\int_{|s_{2}|=r_{2}}
\frac{H(s_{1},s_{2})}{s_{1}s_{2}}\,ds_{2}\,ds_{1}.
\end{align*}
For the original terms of \eqref{8.1}, define
\begin{align*}
&I_{*}^{\heartsuit}(\varphi,\omega',\eta):=I_{*}^{\heartsuit}(\mathbf{s},\varphi,\omega',\eta)|_{\mathbf{s}=(0,0)},\ \ *\in\big\{\mathrm{cusp}, \mathrm{Eis}\big\},\\
&H(\varphi,\omega',\eta):=
\operatorname{CT}_{0}H(\mathbf{s},\varphi,\omega',\eta),\quad
H\in\{R_{\RNum{2}}^{j,\pm},I_{\mathrm{gen}}^{\heartsuit},
I_{\mathrm{degen}}^{\heartsuit}\},\quad j=1,2.
\end{align*}

For the transformed terms in \eqref{8.1}, keep this same $(s_{1},s_{2})$-contour.
If $H$ denotes any of
$R_{\RNum{2}}^{j,\pm}$, $I_{\mathrm{degen}}^{\heartsuit}$, or
$I_{\mathrm{gen}}^{\heartsuit}$, we therefore set
\begin{align*}
H(\pi(w_{2})\varphi,\overline{\omega}\overline{\omega}'\eta,\eta)
:=\operatorname{CT}_{0}
H(\mathbf{s}^{\vee},\pi(w_{2})\varphi,
\overline{\omega}\overline{\omega}'\eta,\eta).
\end{align*}
The cuspidal and Eisenstein terms with $\pi(w_{2})\varphi$ are holomorphic at
the origin and are evaluated there.

Taking $\operatorname{CT}_{0}$ of \eqref{8.1} then gives the identity
at the central parameters.

\part{Local Weights and Dual Transforms}
We compute explicit spectral weights and their dual transforms for type~\RNum{1} and type~\RNum{2} reciprocity and prove the estimates needed for the arithmetic applications in Part 3.

\section{Type \RNum{1}: Archimedean Weights and Dual Transforms}\label{sec9}
We construct the Archimedean weights used in Theorem~\ref{thm1.3}.
Theorem~\ref{thm9.5} bounds their dual transforms in terms of
$C(\xi_{v}|\cdot|_{v}^{\lambda})$ at real and complex places.
The finite-place construction is in
\textsection\ref{sec10}.

In \textsection\ref{sec9}--\textsection\ref{sec10} we use
$\overline{\psi}_{v}$ on $\mathrm{GL}_{3}$ and $\psi_{v}$ on
$\mathrm{GL}_{2}$, reversing \textsection\ref{sec3.2}.
Accordingly, \textsection\ref{sec13} uses reciprocity with
$\overline{\psi}$ in place of $\psi$.

Let $v\mid\infty$, put $d=[F_{v}:\mathbb{R}]$, and suppose $|\Re(s)|\leq 1/4$. Let $\pi_{v}$ and $\sigma_{v}$ be unitary generic representations of $G(F_{v})$ and $G'(F_{v})$, respectively. Let $\mathfrak{B}_{\sigma_{v}}$ be an orthonormal basis of the Whittaker model of $\sigma_{v}$ consisting of $K_{v}'$-type vectors. Let $W_{v}$ be a vector in the Whittaker model of $\pi_{v}$. Define
\begin{equation}\label{e9.1}
\mathcal{H}_{W_{v}}(\sigma_{v},s):=\sum_{W_{v}'}\overline{W_{v}'\left(I_{2}\right)}\int\, W_{v}\left(\diag(x_{v},1)\right)W_{v}'(x_{v})|\det x_{v}|_{v}^{s}dx_{v},
\end{equation}
where $W_{v}'\in \mathfrak{B}_{\sigma_{v}}$ and $x_{v}=(x_{ij})$ ranges over $N'(F_{v})\backslash G'(F_{v})$.

Let $\xi_{v}$ and $\omega_{v}'$ be unitary characters of $F_{v}^{\times}$ and $\lambda\in \mathbb{C}$. Define 
\begin{multline}\label{eqq9.2}
\mathcal{H}_{W_{v}}(\xi_{v},\lambda;\omega_{v}',s)
:=\int_{F_{v}^{\times}}\int_{F_{v}^{\times}}
 W_{v}\bigl(\diag(z_{v},z_{v},1)n_{32}(-y_{v})\bigr)\\
 \overline{\xi}_{v}\omega_{v}'(z_{v})|z_{v}|_{v}^{2s-\lambda-\frac{1}{2}}
 \xi_{v}(y_{v})|y_{v}|_{v}^{\lambda+\frac{1}{2}}\,d^{\times}z_{v}\,d^{\times}y_{v}.
\end{multline}

The integrals are first defined in their regions of absolute convergence
and then continued meromorphically. They give the local weights for
$J_{\mathrm{spec}}$ and $J_{\mathrm{dual}}$ in \eqref{e3.3},
respectively. See \textsection\ref{sec:type-I-dual-expansion} for the
change of Whittaker argument in the dual weight.

\subsection{Positive Analytic-Newvector Weights}
\subsubsection{Analytic conductors}\label{sec:archimedean-analytic-conductors}
Let $v\mid\infty$ and let $\rho_{v}$ be an irreducible representation of $\mathrm{GL}_{n}(F_{v})$, $n\geq1$. Write
\begin{align*}
 L_{v}(s,\rho_{v})=\prod_{j=1}^{n}\zeta_{v}(s+\mu_{j}),
 \qquad C(\rho_{v}):=\prod_{j=1}^{n}(1+|\mu_{j}|)^{d}.
\end{align*}
At any place $v$, let $\vartheta_{\rho_{v}}$ be the largest absolute real part of
the normalized Langlands inducing exponents.

Suppose $\pi_{v}$ and $\sigma_{v}$ are local components of cuspidal automorphic
representations of $G/F$ and $G'/F$, respectively.  By
\cite[Theorem~1]{BB11}, we have
$\vartheta_{\pi_{v}}\leq 5/14$ and
$\vartheta_{\sigma_{v}}\leq 7/64$. In Theorem~\ref{thm1.3}, cuspidal
$\pi$ is an adjoint lift by \cite[Theorem A]{Ram14}, so
$\vartheta_{\pi_{v}}\leq 2(7/64)=7/32$ by the $\mathrm{GL}_{2}$ bound.
For $\pi=\tau\boxplus\chi$, the same $\mathrm{GL}_{2}$ bound gives
$\vartheta_{\pi_{v}}=\vartheta_{\tau_{v}}\leq7/64$ since $\chi_{v}$ is unitary.

\subsubsection{Construction of $W_{v}$}\label{sec9.1.2}
Let $0<\varepsilon<10^{-8}$ and $C_{v}\geq100$. Choose fixed nonnegative
bumps $h_{1,v}\in C_{c}^{\infty}(G'(F_{v}))$ and
$h_{2,v}\in C_{c}^{\infty}(F_{v})$, supported in sufficiently small
neighborhoods of $I_{2}$ and $1$, with
$h_{1,v}(I_{2})=h_{2,v}(1)=1$. For $g_{v}=(g_{ij})$, require
$h_{1,v}(g_{v})=0$ unless
\begin{align*}
 |g_{11}-1|_{v},|g_{22}-1|_{v}\leq\varepsilon,
 \qquad \max\{|g_{12}|_{v},|g_{21}|_{v}\}\leq1,
\end{align*}
and $h_{2,v}(t_{v})=0$ unless $1/2\leq|t_{v}|_{v}\leq2$.
Put $h_{v}=(h_{1,v},h_{2,v})$ and
\begin{align*}
 h_{2,v}(g_{v};C_{v}):=C_{v}h_{2,v}(g_{21}C_{v}^{1/d}).
\end{align*}
On its support, $|g_{21}|_{v}\in[C_{v}^{-1}/2,2C_{v}^{-1}]$.
For $s\in\mathbb{C}$, let $W_{v}\in\pi_{v}$ be the Whittaker vector with restriction
\begin{multline}\label{eq9.2}
W_{v}(\diag(x_{v},1)):=\int_{F_{v}}\int_{G'(F_{v})}h_{1,v}\left(n(b_{v})x_{v}g_{v}\right)h_{2,v}(x_{v}g_{v};C_{v})\\
\overline{h_{1,v}(g_{v})h_{2,v}(g_{v};C_{v})}|\det g_{v}|_{v}^{s+\overline{s}}dg_{v}\psi_{v}(b_{v})db_{v}.
\end{multline}

The right side is smooth in $x_{v}$, compactly supported modulo
$N'(F_{v})$, and transforms on the left by $\overline{\psi}_{v}$.
The Kirillov-model theorem therefore gives a unique $W_{v}$ with this restriction.

\begin{lemma}
Let $W_{v}$ be defined by \eqref{eq9.2}. Then 
\begin{equation}\label{e9.2}
\mathcal{H}_{W_{v}}(\sigma_{v},s)=\sum_{W_{v}'\in \mathfrak{B}_{\sigma_{v}}}\bigg|\int_{G'(F_{v})}h_{1,v}(x_{v})h_{2,v}(x_{v};C_{v})W_{v}'(x_{v})|\det x_{v}|_{v}^{s}dx_{v}\bigg|^{2}.
\end{equation}
\end{lemma}

\begin{proof}
The square in \eqref{e9.2} comes from convolving the test function
with its adjoint.
Put $f(x)=h_{1,v}(x)h_{2,v}(x;C_{v})|\det x|_{v}^{s}$,
$f^{*}(g)=\overline{f(g^{-1})}$, and $\ell_{v}(W')=W'(I_{2})$.
After inserting \eqref{eq9.2},
$y=n_{12}(b_{v})x_{v}$ combines the two integrals into an integral
over $G'(F_{v})$ with weight
$(f*f^{*})(y)=\int f(yg)\overline{f(g)}\,dg$.
Thus the integral in \eqref{e9.1} is
\begin{equation}\label{9.2}
\ell_{v}\bigl(\sigma_{v}(f*f^{*})W_{v}'\bigr)
=\ell_{v}\bigl(\sigma_{v}(f)\sigma_{v}(f)^{*}W_{v}'\bigr).
\end{equation}
Convolution by $f$ smooths the evaluation functional
$\ell_{v}$, so $\ell_{v}\circ\sigma_{v}(f)$ is represented by a smooth
vector. By \eqref{9.2}, the sum in \eqref{e9.1} is the squared norm of this
vector, which is the right-hand side of \eqref{e9.2} by Parseval.
\end{proof}

\begin{cor}\label{cor9.2}
Let $0<\varepsilon_{1}<10^{-3}$, and suppose $\sigma_{v}$ has trivial
central character. For the vector constructed with $s=0$,
\begin{equation}\label{cf9.8}
\mathcal{H}_{W_{v}}(\sigma_{v},0)
\geq c_{h_{v}}-C_{h_{v},\varepsilon_{1}}
C_{v}^{-1/2}C(\sigma_{v})^{1/2+\varepsilon_{1}},
\end{equation}
where $c_{h_{v}}>0$ and $C_{h_{v},\varepsilon_{1}}>0$ depend only on
$h_{v}$ and $\varepsilon_{1}$.
\end{cor}

\begin{proof}
Take $s=0$ and first suppose $C_{v}$ is sufficiently large in terms of $h_{v}$. Define
\begin{multline*}
u(a_{v}):=|a_{v}|_{v}^{-1}\int_{K_{v}'}\int_{F_{v}}\int_{F_{v}^{\times}}h_{1,v}\left(z_{v}n(b_{v})a(a_{v})k_{v}\right)
h_{2,v}\left(z_{v}k_{v};C_{v}\right)\\
\psi_{v}(b_{v})\,d^{\times}z_{v}\,db_{v}\,dk_{v}.
\end{multline*}

Rescaling $k_{21}$ by $C_{v}^{1/d}$ contributes volume $\asymp C_{v}^{-1}$, canceling the factor $C_{v}$ in
$h_{2,v}(\cdot;C_{v})$. Thus $u$ has fixed support and
$\|u\|_{C^{j}}\ll_{j}1$. Choose the bump supports small enough that
$\Re(\psi_{v}(b_{v}))\geq1/2$. Integration over smaller neighborhoods gives
\begin{align*}
 \|u\|_{2}^{2}=\int_{F_{v}^{\times}}|u(a_{v})|^{2}\,d^{\times}a_{v}\asymp 1.
\end{align*}

Let $W_{1,v}'$ be the unit Whittaker vector with $ W_{1,v}'(\diag(a_{v},1))=\|u\|_{2}^{-1}\overline{u(a_{v})}$. The change $k_{v}\mapsto\diag(\epsilon^{-1},1)k_{v}$ gives
$u(\epsilon a_{v})=u(a_{v})$ for $|\epsilon|_{v}=1$.
Thus $W_{1,v}'$ is fixed by the diagonal subgroup of $K_{v}'$.

Write $x_{v}=z_{v}n(b_{v})a(a_{v})k_{v}$.
By Iwasawa decomposition and the definition of $u$,
\begin{equation}\label{equ9.8}
 \int_{G'(F_{v})}h_{1,v}(x_{v})h_{2,v}(x_{v};C_{v})
 \psi_{v}(b_{v})W_{1,v}'(\diag(a_{v},1))\,dx_{v}=\|u\|_{2}\gg1.
\end{equation}
On the support, $k_{v}$ is within $O(C_{v}^{-1/d})$ of the diagonal
subgroup. The bounds for $u$ and \cite[\textsection3.2.5, (3.19)]{MV10}
therefore give
\begin{align*}
 \|\sigma_{v}(k_{v})W_{1,v}'-W_{1,v}'\|
 \ll\min\{1,(C(\sigma_{v})/C_{v})^{1/d}\}
 \ll C_{v}^{-1/2}C(\sigma_{v})^{1/2}.
\end{align*}

After integration in $z_{v},b_{v}$, the coefficient of the $a_{v}$-pairing
has $L^{2}(d^{\times}a_{v})$-norm $O(C_{v})$, and its $k_{v}$-support has
volume $O(C_{v}^{-1})$. Cauchy--Schwarz in $a_{v}$ followed by
integration in $k_{v}$ gives
\begin{multline}\label{equ9.9}
 \int_{G'(F_{v})}h_{1,v}(x_{v})h_{2,v}(x_{v};C_{v})\psi_{v}(b_{v})\\
 \big[W_{1,v}'(\diag(a_{v},1))-W_{1,v}'(\diag(a_{v},1)k_{v})\big]\,dx_{v}
 \ll C_{v}^{-1/2}C(\sigma_{v})^{1/2+\varepsilon_{1}}.
\end{multline}

Since $\|W_{1,v}'\|=1$, \eqref{e9.2} and
\eqref{equ9.8}--\eqref{equ9.9} imply \eqref{cf9.8}.
For the remaining bounded values of $C_{v}$, use $C(\sigma_{v})\geq1$
and choose $C_{h_{v},\varepsilon_{1}}$ so that the right side of
\eqref{cf9.8} is nonpositive; \eqref{e9.2} then gives the result.
\end{proof}

\subsection{\texorpdfstring{Integral Representation of $W_{v}$}{Integral Representation of W\_v}}
For $N'(F_{v})\backslash G'(F_{v})$ with $x_{22}\neq0$, write
\begin{align*}
 x_{v}=\begin{pmatrix}a_{v}z_{v}&\\c_{v}z_{v}&z_{v}\end{pmatrix},
 \qquad a_{v},z_{v}\in F_{v}^{\times},\quad c_{v}\in F_{v}.
\end{align*}

By \eqref{eq9.2}, on the support we have $x_{v}=n(-b_{v})g_{2}g_{1}^{-1}$ with $g_{j}$ in the support of the bumps. Their diagonal entries are close to one and their lower-left entries are $O(C_{v}^{-1/d})$.
Multiplication gives $|a_{v}z_{v}|_{v},|z_{v}|_{v}\asymp1$ and
$c_{v}z_{v}=O(C_{v}^{-1/d})$. Thus
\begin{equation}\label{equ9.11}
 W_{v}(\diag(x_{v},1))=h_{3,v}(a_{v}z_{v},z_{v},C_{v}^{1/d}c_{v}z_{v}),
 \qquad h_{3,v}\in C_{c}^{\infty}(F_{v}^{3}).
\end{equation}

The functions $h_{3,v}$ have common compact support as $C_{v}\geq100$ and
$|\Re(s)|\leq1/4$ vary. The two $h_{2,v}$-factors and the rescaling of the lower-left variable
give $\|h_{3,v}\|_{C^{j}}\ll_{j} C_{v}^{2} C_{v}^{-1}=C_{v}$.
Choose $h_{4,v}\in C_{c}^{\infty}((F_{v}^{\times})^{2})$ equal to $1$ on the
projection of this support, and $h_{6,v}\in C_{c}^{\infty}(F_{v})$ equal to
$1$ on $|c_{v}|\leq1$, which contains the support in $c_{v}$ of $W_{v}$. At a complex place take $h_{6,v}$ radial.

Let $\mathcal{F}h_{3,v}$ be the Fourier transform with respect to $\overline{\psi}_{v}$. Fourier inversion gives
\begin{multline}\label{eq9.10}
 W_{v}(\diag(x_{v},1))
 =h_{6,v}(c_{v})h_{4,v}(a_{v}z_{v},z_{v})
 \int_{F_{v}^{3}}\mathcal{F}h_{3,v}(b_{1},b_{2},b_{3})\\
 \psi_{v}(a_{v}z_{v}b_{1}+z_{v}b_{2}+C_{v}^{1/d}c_{v}z_{v}b_{3})
 \,db_{1}db_{2}db_{3}.
\end{multline}
For $\mathbf{b}=(b_{1},b_{2})$, define the auxiliary vector by
\begin{equation}\label{eq9.11-real}
 W_{v}(\diag(x_{v},1);\mathbf{b})
 =h_{6,v}(c_{v})h_{4,v}(a_{v}z_{v},z_{v})
 \psi_{v}(a_{v}z_{v}b_{1}+z_{v}b_{2}).
\end{equation}
As in \eqref{eq9.2}, this defines a Whittaker vector in $\pi_{v}$, taking its restriction to be zero when $x_{22}=0$.
Every fixed-order coordinate derivative of this restriction grows polynomially in $1+|b_{1}|+|b_{2}|$.

At a complex place, write
$a_{v}=r_{1}e^{i\theta_{1}}$ and $z_{v}=r_{2}e^{i\theta_{2}}$, and put
\begin{align*}
 h_{5,v}(a_{v},z_{v})
 :=h_{4,v}(a_{v}z_{v},z_{v})\psi_{v}(a_{v}z_{v}b_{1}+z_{v}b_{2}).
\end{align*}
Write its angular Fourier expansion as
\begin{align*}
 h_{5,v}(r_{1}e^{i\theta_{1}},r_{2}e^{i\theta_{2}})
 =\sum_{n_{1},n_{2}\in\mathbb{Z}}e^{i(n_{1}\theta_{1}+n_{2}\theta_{2})}
 \mathcal{F}_{n_{1},n_{2}}h_{5,v}(r_{1},r_{2}),
\end{align*}
with Fourier coefficients normalized by $(2\pi)^{-2}d\theta_{1}d\theta_{2}$
on $[0,2\pi]^{2}$.
The radial supports lie in a fixed compact subset of
$(\mathbb{R}_{+}^{\times})^{2}$. Integration by parts in the two angles
gives, for integers $l,j_{1},j_{2}\geq0$,
\begin{multline}\label{eq9.19a}
 \left|(r_{1}\partial_{r_{1}})^{j_{1}}(r_{2}\partial_{r_{2}})^{j_{2}}
 \mathcal{F}_{n_{1},n_{2}}h_{5,v}(r_{1},r_{2})\right|\\
 \ll_{l,j_{1},j_{2}}(1+|b_{1}|+|b_{2}|)^{c_{0}}
 (1+|n_{1}|+|n_{2}|)^{-l},
\end{multline}
for some $c_{0}=c_{0}(l,j_{1},j_{2})$. Define
\begin{align*}
 W_{v}^{(n_{1},n_{2})}(\diag(x_{v},1);\mathbf{b})
 :=h_{6,v}(c_{v})e^{i(n_{1}\theta_{1}+n_{2}\theta_{2})}
 \mathcal{F}_{n_{1},n_{2}}h_{5,v}(r_{1},r_{2}).
\end{align*}
These restrictions have a common compact support. The series
\begin{align*}
  W_{v}(\cdot;\mathbf{b})
 =\sum_{n_{1},n_{2}\in\mathbb{Z}}W_{v}^{(n_{1},n_{2})}(\cdot;\mathbf{b})
\end{align*}
converges in the smooth Kirillov model by \eqref{eq9.19a}.

\begin{lemma}
For $g_{v}\in G(F_{v})$,
\begin{equation}\label{fc9.13}
 W_{v}(g_{v})
 =\int_{F_{v}^{3}}\mathcal{F}h_{3,v}(b_{1},b_{2},b_{3})W_{v}\bigl(g_{v}n_{13}(-C_{v}^{1/d}b_{3});
 \mathbf{b}\bigr)\,db_{1}db_{2}db_{3}.
\end{equation}
\end{lemma}

\begin{proof}
Right translation by $n_{13}(t)$ multiplies the restriction by
$\overline{\psi}_{v}(tx_{21})$.
Thus \eqref{eq9.10} gives \eqref{fc9.13} for $g_{v}=\diag(x_{v},1)$.
The vector $W_{v}(g_{v}n_{13}(-C_{v}^{1/d}b_{3});\mathbf{b})$ and its fixed
derivatives have polynomial norm bounds in $(b_{1},b_{2},b_{3})$.
Since $\mathcal{F}h_{3,v}$ decreases rapidly, the integral converges
in $\pi_{v}$ with all derivatives for fixed $C_{v}$.
Restriction to $\diag(G'(F_{v}),1)$ is injective on the Whittaker
model of $\pi_{v}$, so the identity holds on $G(F_{v})$.
\end{proof}

\subsection{\texorpdfstring{$G'$-spectral Components of $W_{v}(\cdot;\mathbf{b})$}{Spectral Components of W\_v on G'}}
To expand $x_{v}\mapsto W_{v}(\diag(x_{v},1);\mathbf{b})$, use the
$\overline{\psi}_{v}$-Whittaker model of $\sigma_{v}'$ and the
corresponding Whittaker--Plancherel measure.
We take $\gamma(s,\pi_{v}\times\widetilde{\sigma}_{v}';\overline{\psi}_{v})$
in the Rankin--Selberg functional equation.

Fix a positive self-adjoint elliptic differential operator
$\mathcal{D}_{v}$ of order two with $1$ added,
chosen to commute with $K_{v}'$, and choose an orthonormal basis
$\mathfrak{B}_{\sigma_{v}'}$ of its eigenvectors. Write
$\mathcal{D}_{v}W_{v}'=\lambda_{W_{v}'}W_{v}'$, with
$\lambda_{W_{v}'}\geq1$.

\stepcounter{thm}
\addtocounter{equation}{2}

\subsection{An Expansion of the Dual Weight}\label{sec:type-I-dual-expansion}
The matrix in the dual Whittaker function of \eqref{f7.4} becomes,
under \eqref{e2.25},
\begin{align*}
 \widetilde{w}\bigl(\diag(z_{v},1,1)n_{21}(y_{v})\widetilde{w}\bigr)^{\iota}
 =z_{v}^{-1}\diag(z_{v},z_{v},1)n_{32}(-y_{v}).
\end{align*}
The central factor $\omega_{v}(z_{v})^{-1}$ cancels
$\omega_{v}(z_{v})$ in \eqref{f7.4}. Hence
\begin{align*}
\mathcal{J}_{v}(s,\lambda;\varphi,\xi,\omega')=\mathcal{H}_{W_{v}}(\xi_{v},\lambda;\omega_{v}',s),
\end{align*}
where $\mathcal{H}_{W_{v}}(\xi_{v},\lambda;\omega_{v}',s)$ is defined by \eqref{eqq9.2}.
\begin{lemma}\label{lemma9.6}
For $W_{v}$ of \eqref{eq9.2}, we have, initially under absolute
convergence and then meromorphically in $(s,\lambda)$:
\begin{multline}\label{eq9.35}
\mathcal{H}_{W_{v}}(\xi_{v},\lambda;\omega_{v}',s)
=\overline{\xi}_{v}(C_{v}^{1/d})C_{v}^{-\lambda-\frac{1}{2}}\\
\int_{F_{v}^{2}}\int_{(F_{v}^{\times})^{2}}
W_{v}\left(\diag(z_{v},z_{v},1)n_{32}(-y_{v}C_{v}^{-1/d});\mathbf{b}\right)
\mathcal{F}_{12}h_{3,v}(b_{1},b_{2},y_{v})\\
\overline{\xi}_{v}\omega_{v}'(z_{v})|z_{v}|_{v}^{2s-\lambda-1/2}
\xi_{v}(y_{v})|y_{v}|_{v}^{\lambda+1/2}\,d^{\times}z_{v}d^{\times}y_{v}db_{1}db_{2},
\end{multline}
where $\mathcal{F}_{12}h_{3,v}(b_{1},b_{2},y_{v})$ is defined by 
\begin{align*}
\mathcal{F}_{12}h_{3,v}(b_{1},b_{2},y_{v})=\int_{F_{v}}\int_{F_{v}}h_{3,v}(c_{1},c_{2},y_{v})\overline{\psi_{v}(c_{1}b_{1}+c_{2}b_{2})}dc_{1}dc_{2}.
\end{align*}
\end{lemma}

\begin{proof}
Put $g_{v}=\diag(z_{v},z_{v},1)n_{32}(-y_{v})$. Since
$g_{v}n_{13}(t)g_{v}^{-1}=n_{12}(y_{v}t)n_{13}(z_{v}t)$,
we have $W_{v}(g_{v}n_{13}(t);\mathbf{b})=\overline{\psi}_{v}(y_{v}t)W_{v}(g_{v};\mathbf{b})$.
Thus \eqref{fc9.13} and Fourier inversion in $b_{3}$ give
\begin{align*}
 W_{v}(g_{v})=\int_{F_{v}^{2}}
 \mathcal{F}_{12}h_{3,v}(b_{1},b_{2},y_{v}C_{v}^{1/d})
 W_{v}(g_{v};\mathbf{b})\,db_{1}db_{2}.
\end{align*}
Insert this expression for $W_{v}(g_{v})$ in \eqref{eqq9.2}.
The substitution $y_{v}\mapsto C_{v}^{-1/d}y_{v}$ multiplies the $y_{v}$-weight by
$\overline{\xi}_{v}(C_{v}^{1/d})C_{v}^{-\lambda-1/2}$, giving \eqref{eq9.35}.
\end{proof}

\subsection{The Whittaker--Plancherel Decomposition}
For unitary $\chi_{v}$, let $\widehat{G'(F_{v})}_{\chi_{v}}$ denote the
generic tempered dual with central character $\chi_{v}$.
Use the Whittaker--Plancherel measure $d\mu_{\sigma_{v}',\chi_{v}}$
for the quotient of $dx_{v}$ by $d^{\times}z_{v}$.

For $g_{v}\in G(F_{v})$, apply Whittaker--Plancherel
\cite[Chapter~15]{Wal92} to
\begin{align*}
 x_{v}\longmapsto W_{v}(\diag(x_{v},1)g_{v};\mathbf{b})|\det x_{v}|_{v}^{-s'}.
\end{align*}
The displayed function and its fixed derivatives belong to
$L^{2}(N'\backslash G')$ for $|\Re(s')|<1/2-\vartheta_{\pi_{v}}$,
by the standard Whittaker bounds.
In this strip, the local functional equation gives
\begin{multline}\label{9.11}
\int_{N'(F_{v})\backslash G'(F_{v})}W_{v}\left(\diag(1,x_{v})\widetilde{w}g_{v};\mathbf{b}\right)\overline{W_{v}'(x_{v}w')}|\det x_{v}|_{v}^{-s'}dx_{v}\\
=\gamma(1/2-s',\pi_{v}\times\widetilde{\sigma}_{v}')\int_{N'(F_{v})\backslash G'(F_{v})}W_{v}\left(\diag(x_{v},1)g_{v};\mathbf{b}\right)\overline{W_{v}'(x_{v})}|\det x_{v}|_{v}^{-s'}dx_{v}.
\end{multline}

Conjugate $n_{32}(-C_{v}^{-1/d}y_{v})$ by $\diag(1,x_{v})\widetilde{w}$
in \eqref{9.11}. The resulting unipotent has $n_{23}$-coordinate zero
and $n_{12}$-coordinate $-C_{v}^{-1/d}y_{v}(x_{v}^{-1})_{11}$, giving
the factor $\psi_{v}(C_{v}^{-1/d}y_{v}(x_{v}^{-1})_{11})$.
Initially suppose
$-(1/2-\vartheta_{\pi_{v}})<\Re(s-\lambda/2-1/4)<0$, and, for $y_{v}\in F_{v}$, define
\begin{multline*}
 \mathcal{V}_{v}(y_{v}):=\int_{F_{v}^{\times}}
 \big[W_{v}(\diag(z_{v},z_{v},1)n_{32}(-y_{v}/C_{v}^{1/d});\mathbf{b})
       -W_{v}(\diag(z_{v},z_{v},1);\mathbf{b})\big]\\
 \overline{\xi}_{v}\omega_{v}'(z_{v})
 |z_{v}|_{v}^{2s-\lambda-1/2}\,d^{\times} z_{v}.
\end{multline*}
Apply Whittaker--Plancherel, the functional equation and the
$\mathrm{GL}_{1}(F_{v})$-spectral expansion in $z_{v}$, followed by
$m\in\mathbb{Z}_{\geq0}$ integrations by parts with
$\mathcal{D}_{v}$. Fourier inversion in $z_{v}$ with measure $d^{\times}z_{v}$
selects $\omega_{\sigma_{v}'}=\xi_{v}\overline{\omega_{v}'}$.
Right translation of the basis by $w'$ is
unitary and commutes with $\mathcal{D}_{v}$; it gives the evaluation
$W_{v}'(w')$ and leaves the gamma factor indexed by
$\widetilde{\sigma}_{v}'$.

Put $t=s-\lambda/2-1/4$ and $\chi_{v}=\xi_{v}\overline{\omega_{v}'}$.
In the $G'$-integral we replace $x_{v}$ by $C_{v}(\pi)^{-1/d}u_{v}$, giving
\begin{align*}
 |\det x_{v}|_{v}^{t}=C_{v}(\pi)^{-2t}|\det u_{v}|_{v}^{t},
 \qquad
 \overline{W_{v}'(x_{v})}=\chi_{v}(C_{v}(\pi)^{1/d})\overline{W_{v}'(u_{v})}.
\end{align*}
It commutes with $\mathcal{D}_{v}$. Let $\partial^{\boldsymbol{n}}$
mean $d^{n}/dy_{v}^{n}$ over $\mathbb{R}$, and
$\partial_{r}^{n_{1}}\partial_{\theta}^{n_{2}}$ over $\mathbb{C}$, where
$y_{v}=re^{i\theta}$. Then
\begin{equation}\label{e9.46}
\partial^{\boldsymbol{n}}\mathcal{V}_{v}(y_{v})
 =\frac{\chi_{v}(C_{v}(\pi)^{1/d})}{C_{v}(\pi)^{2t}}
 \int_{\widehat{G'(F_{v})}_{\chi_{v}}}
 \sum_{W_{v}'}
 \frac{\lambda_{W_{v}'}^{-m}\mathcal{P}_{\boldsymbol{n}}(\overline{W_{v}'})
       W_{v}'(w')}
 {\gamma(1/2+t,\pi_{v}\times\widetilde{\sigma}_{v}')}
 \,d\mu_{\sigma_{v}',\chi_{v}},
\end{equation}
where
\begin{multline}\label{eq9.47}
 \mathcal{P}_{\boldsymbol{n}}(\overline{W_{v}'})
 :=\int_{N'(F_{v})\backslash G'(F_{v})}\mathcal{D}_{v}^{m}\Big[
 \partial^{\boldsymbol{n}}
\bigl(\psi_{v}(y_{v}(C_{v}(\pi)/C_{v})^{1/d}
 (x_{v}^{-1})_{11})-1\bigr)\\
 W_{v}(\diag(1,C_{v}(\pi)^{-1/d}x_{v})\widetilde{w};\mathbf{b})
 |\det x_{v}|_{v}^{t}\Big]\overline{W_{v}'(x_{v})}\,dx_{v}.
\end{multline}

\subsection{Estimates of Derivatives}
We estimate $W_{v}(\diag(1,x_{v})\widetilde{w};\mathbf{b})$ and its
right $G'$-derivatives for the pairings in \eqref{eq9.47}.

\begin{prop}
Let $\delta\geq0$, $\mathbf{b}\in F_{v}\times F_{v}$, and let $U$ be a
fixed right $G'$-differential operator. Let $x_{v}=z_{v}\diag(a_{v},1)k_{v}$, where $z_{v}, a_{v}\in F_{v}^{\times}$ and $k_{v}=\begin{pmatrix}
k_{11}& k_{12}\\
k_{21} & k_{22}
\end{pmatrix}\in K_{v}'$. Then 
\begin{equation}\label{9.15}
U W_{v}\left(\diag(1,x_{v})\widetilde{w};\mathbf{b}\right)
\ll C_{v}(\pi)^{2\delta}|a_{v}|_{v}^{\frac{1}{2}}
\min\{1,|a_{v}z_{v}|_{v}\}^{2\delta}(1+|b_{1}|+|b_{2}|)^{c_{0}},	
\end{equation}
where $c_{0}$ and the implied constant depend on $U$ and $\delta$.
\end{prop}

\begin{proof}
Write $f(x_{v})=W_{v}(\diag(x_{v},1);\mathbf{b})$ and
$F(x_{v})=W_{v}(\diag(1,x_{v})\widetilde{w};\mathbf{b})$.
The gamma factor in \eqref{9.11} has absolute value one at $s'=0$.
Right translation by $w'$ is unitary and commutes with $\mathcal{D}_{v}$,
so the Plancherel formula gives
\begin{equation}\label{eq:opposite-cell-unitary-sobolev}
 \|\mathcal{D}_{v}^{m}F\|_{L^{2}(N'\backslash G')}^{2}
 =\|\mathcal{D}_{v}^{m}f\|_{2}^{2}.
\end{equation}
The Sobolev inequality, with quotient measure proportional to
$|a_{v}|_{v}^{-1}d^{\times} a_{v}d^{\times} z_{v}dk_{v}$, gives
\begin{equation}\label{eq:quotient-sobolev-evaluation}
 |F(z_{v}\diag(a_{v},1)k_{v})|
 \ll |a_{v}|_{v}^{1/2}\|\mathcal{D}_{v}^{m}F\|_{2}
\end{equation}
for fixed sufficiently large $m$.
Together with \eqref{eq9.11-real}, this proves the case $\delta=0$.

Put $g=\diag(1,x_{v})\widetilde{w}$. The element
\begin{align*}
 Y= a_{v}z_{v}|a_{v}z_{v}|_{v}^{-1/d}(k_{12}E_{31}+k_{11}E_{32})
\end{align*}
has bounded coefficients and satisfies
$\operatorname{Ad}(g)Y=|a_{v}z_{v}|_{v}^{-1/d}E_{12}$.
Differentiating the Whittaker character gives
\begin{align*}
 |W_{v}(g;\mathbf{b})|
 \ll_{l} |a_{v}z_{v}|_{v}^{l/d}|d\pi_{v}(Y)^{l}W_{v}(g;\mathbf{b})|.
\end{align*}
The support of $Y^{l}f$ lies in that of $f$, and its derivative bounds
are those for $f$ multiplied by $C_{v}(\pi)^{l/d}$.
Applying \eqref{eq:opposite-cell-unitary-sobolev} and
\eqref{eq:quotient-sobolev-evaluation}, together with the case $l=0$, gives
\begin{align*}
 |F(x_{v})|\ll_{l} |a_{v}|_{v}^{1/2}
 \min\{1,(C_{v}(\pi)|a_{v}z_{v}|_{v})^{l/d}\}
 (1+|b_{1}|+|b_{2}|)^{c_{0}}.
\end{align*}
Taking $l/d\geq2\delta$ proves \eqref{9.15} for $U=1$. For fixed right
$G'$-derivatives, conjugate by $\widetilde{w}$ to use the upper-left
$G'$-derivative bounds in \eqref{eq9.11-real}.
\end{proof}

To bound $W_{v}(g;\mathbf{b})$ for $|a_{v}|_{v}\geq1$, keep
$g=\diag(1,x_{v})\widetilde{w}$ and put
\begin{align*}
 Y_{a}=|a_{v}|_{v}^{1/d}a_{v}^{-1}\widetilde{w}^{-1}\diag(1,k_{v}^{-1})E_{23}
       \diag(1,k_{v})\widetilde{w}.
\end{align*}
The element $Y_{a}$ lies in the upper-left $\mathfrak{gl}_{2}$,
its coefficients are bounded, and
$\operatorname{Ad}(g)Y_{a}=|a_{v}|_{v}^{1/d}E_{23}$.
Differentiating the Whittaker character gives
\begin{align*}
 |W_{v}(g;\mathbf{b})|\ll_{l} |a_{v}|_{v}^{-l/d}|d\pi_{v}(Y_{a})^{l}W_{v}(g;\mathbf{b})|.
\end{align*}
By \eqref{eq9.11-real}, the resulting right $G'$-derivatives of $f$
have bounds independent of $\pi_{v}$.
Applying \eqref{eq:opposite-cell-unitary-sobolev} and
\eqref{eq:quotient-sobolev-evaluation} gives arbitrary decay in
$|a_{v}|_{v}\geq1$, uniformly in $\pi_{v}$.
In \eqref{9.15}, central scaling replaces
$|a_{v}z_{v}|_{v}$ by $C_{v}(\pi)^{-1}|a_{v}z_{v}|_{v}$, canceling the factor
$C_{v}(\pi)^{2\delta}$ in the bound for small $|a_{v}z_{v}|_{v}$. Combining this with the
large-$a_{v}$ bound, using twice the desired exponents and taking the
geometric mean, gives
\begin{multline}\label{eq:scaled-opposite-cell-decay}
 |U W_{v}(\diag(1,C_{v}(\pi)^{-1/d}x_{v})\widetilde{w};\mathbf{b})|\\
 \ll |a_{v}|_{v}^{1/2}\min\{1,|a_{v}z_{v}|_{v}^{m}\}
 \min\{1,|a_{v}|_{v}^{-l}\}(1+|b_{1}|+|b_{2}|)^{c_{0}}
\end{multline}
for arbitrary $m,l\geq0$ and any fixed right $G'$-differential
operator $U$. The constants depend on $m,l,U$ and the fixed cutoffs.

Use the identities for $Y$ and $Y_{a}$ in the $L^{2}$ bound
\eqref{eq:opposite-cell-unitary-sobolev} to obtain
\begin{multline}\label{eq:scaled-opposite-cell-weighted-unitary}
 \big\||a_{v}z_{v}|_{v}^{-m/d}|a_{v}|_{v}^{l/d}
 U W_{v}(\diag(1,C_{v}(\pi)^{-1/d}x_{v})\widetilde{w};\mathbf{b})\big\|_{L^{2}(N'\backslash G')}\\
 \ll (1+|b_{1}|+|b_{2}|)^{c_{0}},
 \qquad m,l\in\mathbb{Z}_{\geq0}.
\end{multline}
The factor $C_{v}(\pi)^{m/d}$ from $Y^{m}$ cancels
$C_{v}(\pi)^{-m/d}$ from the central scaling. The right
$G'$-derivatives are bounded independently of the conductor.

\subsection{Bounds for the Local Periods}
\stepcounter{thm}
\begin{lemma}\label{lem:normalized-archimedean-sum}
Fix $0\leq\vartheta_{0}<1/2$ and $\delta,\eta>0$ with
$\vartheta_{0}+\delta+\eta<1/2$.
Suppose $\vartheta_{\pi_{v}}\leq\vartheta_{0}$ and
\begin{align*}
 -1/2+\vartheta_{0}+\delta\leq\Re(t)\leq-\eta,
 \qquad t=s-\lambda/2-1/4.
\end{align*}
Write $A_{v}=C_{v}(\pi)/C_{v}$ and $\chi_{v}=\xi_{v}\overline{\omega_{v}'}$,
and define $\mathcal{P}_{\boldsymbol{n}}$ by \eqref{eq9.47} with $m=0$.
There is an integer $r\geq0$, independent of the derivative order
$n$, such that, on the $y_{v}$-support of $h_{3,v}$, for derivatives of
total order $n$ and some integer $c_{0}\geq0$,
\begin{multline}\label{eq:normalized-archimedean-polynomial}
 \int_{\widehat{G'(F_{v})}_{\chi_{v}}}
 \sum_{W_{v}'\in\mathfrak{B}_{\sigma_{v}'}}|W_{v}'(w')|
 \left|C_{v}(\pi)^{-2t}
 \gamma(1/2+t,\pi_{v}\times\widetilde{\sigma}_{v}')^{-1}\right|\\
 \left|\mathcal{P}_{\boldsymbol{n}}(\overline{W_{v}'})\right|
 d\mu_{\sigma_{v}',\chi_{v}}
 \ll (1+|t|+|b_{1}|+|b_{2}|)^{c_{0}}
 \max\{A_{v}^{1/d},A_{v}^{\max\{1,n+2r\}/d}\}.
\end{multline}
The constants are uniform in the unitary character $\chi_{v}$.
In particular, if $A_{v}\leq1$, the conductor factor is
$A_{v}^{1/d}$, for every fixed derivative order.
\end{lemma}

\begin{proof}
Let $h(x_{v})$ be the expression inside $\mathcal{D}_{v}^{m}$ in
\eqref{eq9.47}. Put
\begin{align*}
 h_{\chi_{v}}(x_{v})=\int_{F_{v}^{\times}}
 h(z_{v}x_{v})\overline{\chi_{v}(z_{v})}\,d^{\times} z_{v}.
\end{align*}
The function $h_{\chi_{v}}$ has central character $\chi_{v}$.
Write $x_{v}=z_{v}\diag(a_{v},1)k_{v}$ and $p_{v}=a_{v}z_{v}$. The first column of $x_{v}^{-1}$ has the scalar $p_{v}^{-1}$.
Differentiating $x_{v}\mapsto x_{v}\exp(tU)$ multiplies $x_{v}^{-1}$ on
the left, so the required derivatives of $(x_{v}^{-1})_{11}$ are
$O(|p_{v}|_{v}^{-1/d})$. The Leibniz rule and
\eqref{eq:scaled-opposite-cell-decay} bound $\mathcal{D}_{v}^{j}h$ by
sums, for $1\leq l\leq\max\{1,n+2j\}$, of
\begin{align*}
(1+|t|)^{2j}(1+|b_{1}|+|b_{2}|)^{c_{0}} A_{v}^{\frac{l}{d}}|a_{v}|_{v}^{\frac{1}{2}-\Re(t)}\min\{1,|a_{v}|_{v}^{-l_{2}}\}
 |p_{v}|_{v}^{2\Re(t)-\frac{l}{d}}\min\{1,|p_{v}|_{v}^{l_{1}}\}.
\end{align*}

For $l_{1},l_{2}$ large, this bound and
$W_{v}'(z_{v}\diag(a_{v},1)k_{v})\ll_{W_{v}',\eta'}|a_{v}|_{v}^{1/2-\eta'}$
at $|a_{v}|_{v}\leq1$, with $0<\eta'<\eta$
\cite[Proposition~3.2.3]{MV10}, give absolute convergence in \eqref{eq9.47}.
In $h_{\chi_{v}}(\diag(a_{v},1)k_{v})$, substitute
$z_{v}\mapsto a_{v}^{-1}p_{v}$.
The $p_{v}$-integral is bounded, leaving the factor
$|a_{v}|_{v}^{1/2-\Re(t)}\min\{1,|a_{v}|_{v}^{-l_{2}}\}$.
Its square is integrable against $|a_{v}|_{v}^{-1}d^{\times} a_{v}dk_{v}$:
the resulting power at $a_{v}=0$ is $-2\Re(t)>0$, and $l_{2}$ gives
convergence at infinity. Thus
\begin{equation}\label{eq:central-projection-sobolev}
 \|\mathcal{D}_{v}^{j}h_{\chi_{v}}\|_{2}
 \ll (1+|t|)^{2j}(1+|b_{1}|+|b_{2}|)^{c_{0}}
 \sum_{l=1}^{\max\{1,n+2j\}}A_{v}^{l/d}.
\end{equation}
These bounds justify differentiation under the integral and integration
by parts, uniformly in $\chi_{v}$.
The Plancherel formula with fixed central character now gives
\begin{equation}\label{eq:central-projection-plancherel}
 \int_{\widehat{G'(F_{v})}_{\chi_{v}}}
 \sum_{W_{v}'}\lambda_{W_{v}'}^{2j}
 |\mathcal{P}_{\boldsymbol{n}}(\overline{W_{v}'})|^{2}
 d\mu_{\sigma_{v}',\chi_{v}}
 =\|\mathcal{D}_{v}^{j}h_{\chi_{v}}\|_{2}^{2}.
\end{equation}

Stirling's formula and the Archimedean parameter formulas give
\begin{multline*}
 C_{v}(\pi)^{-2\Re(t)}
 |\gamma(1/2+t,\pi_{v}\times\widetilde{\sigma}_{v}')^{-1}|
 \ll [C(\sigma_{v}')(1+|\Im(t)|)^{2d}]^{-3\Re(t)}\\
 \ll (1+|t|)^{c_{1}}\lambda_{W_{v}'}^{c_{1}}.
\end{multline*}

For a fixed sufficiently large integer $q=q(F_{v})$, local Sobolev
evaluation and Plancherel give, uniformly in $\chi_{v}$,
\begin{align*}
 \int_{\widehat{G'(F_{v})}_{\chi_{v}}}
 \sum_{W_{v}'}\lambda_{W_{v}'}^{-2q}|W_{v}'(w')|^{2}
 d\mu_{\sigma_{v}',\chi_{v}}\ll1.
\end{align*}
Take $r\geq q+c_{1}$ independently of $n$. Cauchy--Schwarz and
\eqref{eq:central-projection-sobolev}--\eqref{eq:central-projection-plancherel},
with $j=r$, give \eqref{eq:normalized-archimedean-polynomial}.
If $\pi_{v}$ is of orthogonal type with trivial central character,
then, for every tempered $\rho_{v}$ on $G'(F_{v})$,
$C(\pi_{v}\times\rho_{v})\gg C_{v}(\pi)^{2}$.
If all summands are quadratic characters, $C_{v}(\pi)\asymp1$.
For $\chi\oplus\mathbf{1}\oplus\chi^{-1}$, at least two of
$\|u-b\|,\|u\|,\|u+b\|$ are $\geq\|b\|/2$, where $u,b$ record the
imaginary norm exponents and the angular exponents at a complex place.
For a real two-dimensional summand of $\rho_{v}$, use its norm parameter
$u$ and count its Gamma factor twice.
If $\pi_{v}$ has a real two-dimensional summand induced from angular
exponent $k\geq1$, a character twist preserves $k$. Against another
such summand of exponent $l$, the tensor product has exponents
$k+l,|k-l|,l$, with $(1+k+l)(1+|k-l|)(1+l)\gg(1+k)^{2}$.
With $\rho_{v}=\widetilde{\sigma}_{v}'\otimes|\cdot|_{v}^{i\Im(t)}$
and $\Re(t)<0$, Stirling's formula gives
\begin{align*}
 \left|C_{v}(\pi)^{-2t}
 \gamma(1/2+t,\pi_{v}\times\widetilde{\sigma}_{v}')^{-1}\right|
 \ll_{F,\vartheta_{0},\delta}1.
\end{align*}
Thus $c_{1}=0$, $r=q$, and the $t$-factor is $(1+|t|)^{2q}$
for every derivative order in $y_{v}$.
\end{proof}

\begin{lemma}\label{lem:unitary-residue-archimedean-endpoint}
Suppose $s-\lambda/2-1/4=0$ and
$(\xi_{v}\overline{\omega_{v}'})^{2}=1$.
Put $A_{v}=C_{v}(\pi)/C_{v}$.
On the fixed compact set in the $y_{v}$-variable, for every
$\epsilon>0$,
\begin{equation}\label{eq:unitary-residue-archimedean-endpoint}
 |\mathcal{V}_{v}(y_{v})|\ll_{\epsilon}
 C_{v}(\pi)^{\epsilon}(1+|b_{1}|+|b_{2}|)^{c_{0}}
 \big[A_{v}^{1/d}|y_{v}|+A_{v}|y_{v}|^{d}\big],
\end{equation}
where $c_{0}$ is a fixed exponent and $|y_{v}|$ is the Euclidean absolute value.
\end{lemma}

\begin{proof}
Let $\chi_{v}=\xi_{v}\overline{\omega_{v}'}$ and put
$B_{y}=A_{v}|y_{v}|^{d}$; the assertion is immediate if $y_{v}=0$.
Use $h_{\chi_{v}}$ from Lemma~\ref{lem:normalized-archimedean-sum},
with $t=0$ and derivative order zero.
For $p_{v}=a_{v}z_{v}$, we have
\begin{align*}
 \big|\psi_{v}(y_{v}A_{v}^{1/d}(x_{v}^{-1})_{11})-1\big|
 \ll\min\{1,B_{y}^{1/d}|p_{v}|_{v}^{-1/d}\}.
\end{align*}

For $|p_{v}|_{v}\leq1$, Cauchy--Schwarz and
\eqref{eq:scaled-opposite-cell-weighted-unitary}, with a sufficiently
large $m$, bound the central integral in $L^{2}(Z'N'\backslash G')$.
For $|p_{v}|_{v}\geq1$, Cauchy--Schwarz separates the character factor
from $W_{v}(\diag(1,C_{v}(\pi)^{-1/d}x_{v})\widetilde{w};\mathbf{b})$.
The character factor satisfies
\begin{align*}
 \int_{|p_{v}|_{v}\geq1}
 \min\{1,B_{y}^{2/d}|p_{v}|_{v}^{-2/d}\}\,d^{\times}p_{v}
 \ll\begin{cases}
 B_{y}^{2/d},&B_{y}\leq1,\\
 1+\log B_{y},&B_{y}\geq1.
 \end{cases}
\end{align*}
Consequently,
\begin{align*}
 \|h_{\chi_{v}}\|_{2}\ll (1+|b_{1}|+|b_{2}|)^{c_{0}}
 \begin{cases}
 B_{y}^{1/d},&B_{y}\leq1,\\
 (1+\log B_{y})^{1/2},&B_{y}\geq1.
 \end{cases}
\end{align*}
For an integer $j\geq1$, applying $\mathcal{D}_{v}^{j}$ to $h_{\chi_{v}}$ and using
\eqref{eq:scaled-opposite-cell-weighted-unitary} gives
\begin{align*}
 \|\mathcal{D}_{v}^{j}h_{\chi_{v}}\|_{2}
 \ll (1+|b_{1}|+|b_{2}|)^{c_{0}}\sum_{l=1}^{2j}B_{y}^{l/d}.
\end{align*}
The same estimates give convergence of $\mathcal{D}_{v}^{j}h_{\chi_{v}}$
in $L^{2}$ as $t\to0^{-}$. Indeed, the extra factor
$|p_{v}|_{v}^{2t}|a_{v}|_{v}^{-t}$ is absorbed by the weights in
\eqref{eq:scaled-opposite-cell-weighted-unitary} and the decay of the
phase difference for $|p_{v}|_{v}\geq1$. The polynomial bound for the
inverse gamma factor is uniform near $t=0$ by Stirling, so
Cauchy--Schwarz as in Lemma~\ref{lem:normalized-archimedean-sum} gives
\eqref{e9.46} at $t=0$ by continuity.

There the gamma factor in \eqref{e9.46} has absolute value one.
Sobolev's inequality gives
$|\mathcal{V}_{v}(y_{v})|\ll_{q}\|\mathcal{D}_{v}^{q}h_{\chi_{v}}\|_{2}$
for every real $q>d/2$.
H\"older's inequality gives, for $j>q$,
\begin{align*}
 \|\mathcal{D}_{v}^{q}h_{\chi_{v}}\|_{2}
 \leq \|h_{\chi_{v}}\|_{2}^{1-q/j}
 \|\mathcal{D}_{v}^{j}h_{\chi_{v}}\|_{2}^{q/j}.
\end{align*}
Apart from $(1+|b_{1}|+|b_{2}|)^{c_{0}}$, the interpolated bound is
$O(B_{y}^{1/d})$ for $B_{y}\leq1$ and
$O(B_{y}^{2q/d}(1+\log B_{y})^{(1-q/j)/2})$ for $B_{y}\geq1$.
Choose $2q/d<1+\epsilon$ and absorb the logarithm to obtain
\begin{align*}
 |\mathcal{V}_{v}(y_{v})|\ll_{\epsilon}
 (1+|b_{1}|+|b_{2}|)^{c_{0}}\big[B_{y}^{1/d}+B_{y}^{1+\epsilon}\big].
\end{align*}
Since $B_{y}\ll C_{v}(\pi)$ on the fixed $y_{v}$-support,
this proves \eqref{eq:unitary-residue-archimedean-endpoint}.
\end{proof}

\subsection{The Dual Weight}
Use $W_{v}$ constructed with $s=0$; below, $s$ denotes the
moment parameter. For the residue estimates, define $\mathcal{H}_{W_{v}}^{\circ}$
by replacing the Whittaker factor in \eqref{eq9.35} with
$W_{v}(\diag(z_{v},z_{v},1);\mathbf{b})$.
Since $h_{4,v}(z_{v},z_{v})=1$ on the support of
$h_{3,v}(z_{v},z_{v},y_{v})$, \eqref{eq9.11-real} and Fourier inversion
in $b_{1},b_{2}$ give
\begin{multline*}
\mathcal{H}_{W_{v}}^{\circ}(\xi_{v},\lambda;\omega_{v}',s)
 =\overline{\xi}_{v}(C_{v}^{1/d})C_{v}^{-\lambda-1/2}
 \int_{(F_{v}^{\times})^{2}}h_{3,v}(z_{v},z_{v},y_{v})\\
\overline{\xi}_{v}\omega_{v}'(z_{v})|z_{v}|_{v}^{2s-\lambda-1/2}
 \xi_{v}(y_{v})|y_{v}|_{v}^{\lambda+1/2}\,d^{\times} z_{v}d^{\times} y_{v}.
\end{multline*}

Lemma~\ref{lemma9.6} then gives
\begin{multline}\label{9.69}
\mathcal{H}_{W_{v}}(\xi_{v},\lambda;\omega_{v}',s)-\mathcal{H}_{W_{v}}^{\circ}(\xi_{v},\lambda;\omega_{v}',s)=\overline{\xi}_{v}(C_{v}^{1/d})C_{v}^{-\lambda-\frac{1}{2}}\\
\int_{F_{v}}\int_{F_{v}}\int_{F_{v}^{\times}}\mathcal{F}_{12}h_{3,v}(b_{1},b_{2},y_{v})\mathcal{V}_{v}(y_{v})
\xi_{v}(y_{v})|y_{v}|_{v}^{\lambda+\frac{1}{2}}d^{\times}y_{v}db_{1}db_{2}.
\end{multline}

Integrate by parts in $c_{1},c_{2}$ in $\mathcal{F}_{12}h_{3,v}$.
The compact support and the $O(C_{v})$ derivative bounds of $h_{3,v}$
following \eqref{equ9.11} give, for every $l\geq0$,
\begin{align*}
 \mathcal{F}_{12}h_{3,v}(b_{1},b_{2},y_{v})
 \ll_{l} C_{v}(1+|b_{1}|)^{-l}(1+|b_{2}|)^{-l}.
\end{align*}
The decay in $b_{1},b_{2}$ is unchanged after fixed derivatives in
$y_{v}$ over $\mathbb{R}$, or in $r,\theta$ for $y_{v}=re^{i\theta}$
over $\mathbb{C}$, with constants depending on the derivatives
and the fixed bumps of \textsection\ref{sec9.1.2}.

For $|\Re(\lambda)|\leq1/128$ and $|s|$ sufficiently small,
integration by parts gives
\begin{equation}\label{9.73}
 |\mathcal{H}_{W_{v}}^{\circ}(\xi_{v},\lambda;\omega_{v}',s)|
 \ll_{n}\frac{C_{v}^{1/2+\max\{0,-\Re(\lambda)\}}}
 {C(\xi_{v}|\cdot|_{v}^{\lambda})^{n}}.
\end{equation}
Absorb $\Im(\lambda)$ into $\xi_{v}$ and integrate by parts
$n$ times with $y_{v}\partial_{y_{v}}$ over $\mathbb{R}$, or
$2n$ times with $r\partial_{r},\partial_{\theta}$ over $\mathbb{C}$.
There is no boundary term at zero since $\Re(\lambda)+1/2>0$.
The fixed $z_{v}$-support away from zero and the $O(C_{v})$ derivative
bounds for $h_{3,v}$ give $O_{n}(C_{v} C(\xi_{v}|\cdot|_{v}^{\lambda})^{-n})$
for the double integral in $(z_{v},y_{v})$. Multiplying by
$C_{v}^{-\lambda-1/2}$ proves \eqref{9.73}.

\begin{thm}\label{thm9.5}
Suppose $\pi$ is as in Theorem~\ref{thm1.3}, and let $\xi_{v}$ and
$\omega_{v}'$ be unitary characters. For every $\epsilon>0$ there
is $\delta>0$ such that,
for integers $n\geq1$, $|s|\leq\delta$ and $|\Re(\lambda)|\leq\delta$,
\begin{equation}\label{fc9.37}
|\mathcal{H}_{W_{v}}(\xi_{v},\lambda;\omega_{v}',s)|
\ll_{F,n,\epsilon}
 \frac{C_{v}^{1/2+\max\{0,-\Re(\lambda)\}+\epsilon} [1+(C_{v}^{-1}C_{v}(\pi))^{n}]}
 {C(\xi_{v}|\cdot|_{v}^{\lambda})^{n}[C_{v}(\pi)C_{v}(\omega')]^{-\epsilon}}.
\end{equation}
\end{thm}

\begin{proof}
The defining integral \eqref{eqq9.2} gives
\begin{align*}
 \mathcal{H}_{W_{v}}(\xi_{v},\lambda;\omega_{v}',s)
 =\mathcal{H}_{W_{v}}(\xi_{v}|\cdot|_{v}^{i\Im(\lambda)},\Re(\lambda);
                   \omega_{v}',s).
\end{align*}
Put $\rho=\xi_{v}|\cdot|_{v}^{i\Im(\lambda)}$.
The following calculation uses only $\omega_{v}=\mathbf{1}_{v}$. Put
\begin{align*}
 B=1+C_{v}(\pi)/C_{v},\qquad S=C_{v}B^{2},\qquad V=\pi_{v}(w_{2}\widetilde{w})W_{v}.
\end{align*}
Write $v(p,q,r)=V\bigl(\begin{smallmatrix}p&&\\q&r&\\&&1\end{smallmatrix}\bigr)$.
For integers $m,l\geq0$, let $\mathcal{D}_{q}^{m}$ denote $m$ successive applications
of $q\partial_{q}$ over $\mathbb{R}$, or any composition of $m$ operators
chosen from $|q|\partial_{|q|}$ and $\partial_{\arg q}$ over $\mathbb{C}$.

By \eqref{equ9.11},
$C_{v}^{-1}W_{v}\bigl(\begin{smallmatrix}p&&\\q&r&\\&&1\end{smallmatrix}\bigr)$
has fixed compact support and bounded derivatives in $(p,r,C_{v}^{1/d}q)$,
with $p,r$ away from zero. Hence $\|W_{v}\|^{2}\ll C_{v}$.
Since $\omega_{v}=\mathbf{1}_{v}$, the central identities in
\cite[\textsection2, (2.21)]{Mol06} show that the representation
parameters enter the $E_{31}$ action only through $p$-derivatives.
Each $E_{31}$-derivative
therefore increases these bounds by at most $(C_{v}(\pi)+C_{v})^{1/d}$,
with support unchanged, also after fixed derivatives in the rescaled coordinates.
Conjugation by $(w_{2}\widetilde{w})^{-1}$ and unitarity then give
norm $O(C_{v}^{1/2}B^{m/d}S^{l/d})$ after applying $\mathcal{D}_{q}^{m}$ to $v$
and multiplying by $|p|_{v}^{l/d}$.
Apply Sobolev's inequality in $r$ at $r=1$ to bound the
$L^{2}$ norm in $(p,q)$. To bound the $L^{2}$ norm in $p$ at $q=q_{0}$,
also apply Sobolev's inequality in $q$ on intervals of width
$C_{v}^{1/d}$ in each real coordinate.
The conjugated derivatives $E_{11}$ and $C_{v}^{1/d}E_{13}$ preserve
the bounds in the rescaled coordinates on $W_{v}$. The $q$-intervals
have volume $C_{v}$,
giving an additional factor $C_{v}^{-1}$ in the squared bound.
Thus
\begin{equation}\label{eq:reflected-joint-trace}
 \int_{F_{v}^{\times}}\int_{F_{v}}
 |\mathcal{D}_{q}^{m}v(p,q,1)|^{2}|p|_{v}^{2l/d-1}\,dq\,d^{\times} p
 \ll_{m,l}C_{v}B^{2m/d}S^{2l/d},
\end{equation}
and
\begin{equation}\label{eq:reflected-point-trace}
 \sup_{q_{0}\in F_{v}}\int_{F_{v}^{\times}}
 |\mathcal{D}_{q}^{m}v(p,q_{0},1)|^{2}|p|_{v}^{2l/d-1}\,d^{\times} p
 \ll_{m,l}B^{2m/d}S^{2l/d}.
\end{equation}
For fixed unitary $\rho$, consider
$\mathcal{H}_{W_{v}}(\rho,\lambda;\omega_{v}',s)$. Put
$a_{0}=1/2+2s-\lambda$ and $b_{0}=1/2+\lambda$.

The $\mathrm{GL}_{3}\times\mathrm{GL}_{1}$ and Tate local functional
equations give
\begin{multline}\label{eq:reflected-dual-double-mellin}
 \gamma(a_{0},\widetilde{\pi}_{v}\otimes\overline{\rho}\omega_{v}',\psi_{v})
 \gamma(b_{0},\rho,\psi_{v})
 \mathcal{H}_{W_{v}}(\rho,\lambda;\omega_{v}',s)\\
 =\zeta_{v}(1)\rho(-1)
 \int_{F_{v}^{\times}}\int_{F_{v}}
 v(p,q,1)
 \rho(p)\overline{\omega_{v}'}(p)\overline{\rho}(q)|p|_{v}^{-a_{0}}|q|_{v}^{-b_{0}}
 \,dq\,d^{\times} p.
\end{multline}
Apply the first functional equation to the $z_{v}$-integral in
\eqref{f7.4}, with $y_{v}$ fixed. It gives a double integral over
$(p,q)\in F_{v}^{\times}\times F_{v}$, with Whittaker function
$g\mapsto W_{v}(gw_{2}\widetilde{w} n_{32}(-y_{v}))$.
Since $(w_{2}\widetilde{w})n_{32}(-y_{v})(w_{2}\widetilde{w})^{-1}=n_{13}(-y_{v})$, its Whittaker factor is
$\psi_{v}(qy_{v})$.
The Tate identity is
\begin{align*}
 \int_{F_{v}^{\times}}\rho(y_{v})|y_{v}|_{v}^{b_{0}}\psi_{v}(qy_{v})d^{\times} y_{v}
 =\zeta_{v}(1)\rho(-1)\gamma(b_{0},\rho,\psi_{v})^{-1}
   \overline{\rho}(q)|q|_{v}^{-b_{0}}.
\end{align*}
For $0<\Re(b_{0})<1$, this identity holds after integration against a
Schwartz function of $q$. Commuting the right translation by
$n_{13}(t)$ to the left in $v(p,q,1)$ gives the Whittaker character
$\overline{\psi_{v}(qt)}$. Differentiation in $t$ and moderate growth
give rapid decay in $q$ uniformly for $p$ in compact subsets of $F_{v}^{\times}$.
Apply the identity to this function and use
\eqref{eq:reflected-joint-trace}--\eqref{eq:reflected-point-trace}
to integrate over all $p$. The first functional equation and
meromorphic continuation give \eqref{eq:reflected-dual-double-mellin}.

Take $0<\alpha<1/128$, $|s|\leq\alpha/16$, and
$\Re(\lambda)=\alpha$. Apply Cauchy--Schwarz in $p$, using the $L^{2}$ weight $|p|_{v}^{2l/d-1}$
with $l=0$ on $|p|_{v}\leq S$ and large $l$ on $|p|_{v}>S$.
The factor from the other norm is $O(S^{\alpha-2\Re(s)-l/d})$
on the respective range. For $|q|_{v}\leq C_{v}$, use
\eqref{eq:reflected-point-trace} and integrate $|q|_{v}^{-1/2-\alpha}$.
For $|q|_{v}>C_{v}$, use \eqref{eq:reflected-joint-trace} and the
$L^{2}$ norm of $|q|_{v}^{-1/2-\alpha}$, which is $O_{\alpha}(C_{v}^{-\alpha})$.
The joint bound contributes $C_{v}^{1/2}$, so both ranges have the
factor $C_{v}^{1/2-\alpha}$. With $S=C_{v}B^{2}$,
\begin{multline*}
 \int_{F_{v}^{\times}}\int_{F_{v}}|\mathcal{D}_{q}^{m}v(p,q,1)|
 |p|_{v}^{-1/2+\alpha-2\Re(s)}|q|_{v}^{-1/2-\alpha}\,dq\,d^{\times} p\\
 \ll_{\alpha,m}C_{v}^{1/2-2\Re(s)}B^{m/d+2\alpha-4\Re(s)}.
\end{multline*}
Integrate by parts in $q$ with smooth cutoffs. Use $\Re(b_{0})<1$
to remove the cutoff at zero and the weighted $L^{1}$ estimate for
$\mathcal{D}_{q}^{m} v(p,q,1)$ to remove it at infinity.
Integrating by parts $n$ times over $\mathbb{R}$, or $2n$ times over
$\mathbb{C}$, gives the additional factor $B^{n}C(\rho)^{-n}$.

The odd-dimensional orthogonal local parameter contains a
quadratic character. The elementary conductor inequality for
characters therefore gives
$C(\widetilde{\pi}_{v}\otimes\overline{\rho}\omega_{v}'
 {|\cdot|_{v}^{2i\Im(s)}})\gg_{F} C(\rho)/C_{v}(\omega')$ for $|s|\leq\alpha/16$.
Stirling's formula, on these pole-free strips, consequently gives
\begin{align*}
 |\gamma(a_{0},\widetilde{\pi}_{v}\otimes\overline{\rho}\omega_{v}',\psi_{v})^{-1}
   \gamma(b_{0},\rho,\psi_{v})^{-1}|
 \ll_{\alpha} C_{v}(\omega')^{2\alpha}C(\rho)^{\alpha/8}.
\end{align*}
We have proved, uniformly on the line $\Re(\lambda)=\alpha$,
\begin{align*}
 |\mathcal{H}_{W_{v}}(\rho,\lambda;\omega_{v}',s)|
 \ll_{\alpha,n}C_{v}(\omega')^{2\alpha}C_{v}^{1/2+\alpha}B^{n+3\alpha}
             C(\rho|\cdot|_{v}^{i\Im(\lambda)})^{-n+\alpha/8}.
\end{align*}
On $\Re(\lambda)=-1/128$, Lemma~\ref{lem:normalized-archimedean-sum},
\eqref{9.69}, and \eqref{9.73}, followed by integration by parts, give
\begin{align*}
 |\mathcal{H}_{W_{v}}(\rho,\lambda;\omega_{v}',s)|
 \ll_{n} C_{v}^{1/2+1/128}B^{n+O_{F}(1)}
          C(\rho|\cdot|_{v}^{i\Im(\lambda)})^{-n}.
\end{align*}
On $\Re(\lambda)=\alpha$, raise the bounds for $n$ and $n+1$ to the
powers $1-\alpha/8$ and $\alpha/8$ and multiply them.
Both boundary bounds then have the factor
$C(\rho|\cdot|_{v}^{i\Im(\lambda)})^{-n}$, with an additional factor
$B^{\alpha/8}$ on the positive line.

Lemmas~\ref{lem7.3} and~\ref{lem:normalized-archimedean-sum} give
holomorphy and polynomial growth in $|\Im(\lambda)|$ on
$-1/128\leq\Re(\lambda)\leq\alpha$.
Multiply by a polynomial in $\lambda$ whose absolute value on this strip
is comparable to $C(\rho|\cdot|_{v}^{i\Im(\lambda)})^{n}$.
The maximum principle then gives, for $|\ell|\leq\alpha/2$,
\begin{equation}\label{eq:orthogonal-interpolated-dual-weight}
 |\mathcal{H}_{W_{v}}(\rho,\ell;\omega_{v}',s)|
 \ll_{\alpha,n}C_{v}(\omega')^{2\alpha}
 C_{v}^{1/2+\max\{0,-\ell\}+2\alpha}B^{n+O_{F}(\alpha)}C(\rho)^{-n}.
\end{equation}
Choose $\alpha$ sufficiently small in terms of $\epsilon$,
independently of $n$, and take $\delta\leq\alpha/16$.
Since $B\ll C_{v}(\pi)$ and $B^{n}\ll_{n}1+(C_{v}(\pi)/C_{v})^{n}$,
\eqref{eq:orthogonal-interpolated-dual-weight}, with
$\ell=\Re(\lambda)$ and $\Im(\lambda)$ absorbed into $\xi_{v}$,
proves \eqref{fc9.37}.
\end{proof}
By \eqref{9.73}, \eqref{fc9.37} also bounds
$\mathcal{H}_{W_{v}}-\mathcal{H}_{W_{v}}^{\circ}$.

\section{Type \RNum{1}: Non-Archimedean Weights and Dual Transforms}\label{sec10}
We construct nonnegative spectral weights for representations of bounded conductor.
Fix $v<\infty$, and use \eqref{e9.1} and \eqref{eqq9.2} over $F_{v}$.
For $h\in C_{c}^{\infty}(G'(F_{v}))$, let
$h^{*}(g)=\overline{h(g^{-1})}$. Then
\begin{align*}
 (h\ast h^{*})(x)=\int_{G'(F_{v})}h(xg)\overline{h(g)}\,dg.
\end{align*}
Let $\sigma_{v}$ be a unitary generic representation of $G'(F_{v})$
with central character $\omega_{v}'$. Put
$r_{v}=\max\{e_{v}(\mathfrak{q}),r_{\pi_{v}}\}$ and assume
$r_{\omega_{v}'}\leq e_{v}(\mathfrak{q})$.
Define
\begin{equation}\label{eq:finite-projector-function}
 h_{v}(a(\varpi_{v}^{-d_{v}})k)
 =\Vol(K_{0,v}'[r_{v}])^{-1}
 \begin{cases}
  1,&r_{v}=0,\\
  \overline{\omega_{v}'(k_{22})},&r_{v}\geq1,
 \end{cases}
 \qquad k\in K_{0,v}'[r_{v}],
\end{equation}
and set $h_{v}=0$ outside $a(\varpi_{v}^{-d_{v}})K_{0,v}'[r_{v}]$.
For $r_{v}\geq1$, $k_{22}$ is a unit and
$k\mapsto\omega_{v}'(k_{22})$ is a character, since
$r_{\omega_{v}'}\leq r_{v}$. Thus $\|h_{v}\|_{1}=1$.
The translation
\begin{align*}
 V(g)\longmapsto V\bigl(a(\varpi_{v}^{d_{v}})g\bigr)
\end{align*}
from the $\psi_{v}^{0}$-model to the $\psi_{v}$-model takes a vector with
value $1$ at $I_{2}$ to one with value $1$ at $a(\varpi_{v}^{-d_{v}})$;
see \textsection\ref{sec:finite-character-change}.

Set $g=a(\varpi_{v}^{-d_{v}})k$ in the convolution integral to obtain
\begin{equation}\label{eq:finite-ordinary-correlation}
 (h_{v}\ast h_{v}^{*})(x)=h_{v}\bigl(xa(\varpi_{v}^{-d_{v}})\bigr).
\end{equation}

For $x\in G'(F_{v})$, prescribe the restriction
\begin{equation}\label{c11.1}
 W_{v}(\diag(x,1)):=
 \int_{F_{v}}(h_{v}\ast h_{v}^{*})\left(
 n(b)x\right)\psi_{v}(b)\,db.
\end{equation}
This restriction is supported on
$N'(F_{v})a(\varpi_{v}^{-d_{v}})K_{0,v}'[r_{v}]a(\varpi_{v}^{d_{v}})$, and
\begin{equation}\label{eq:finite-projector-restriction}
 W_{v}\left(\diag\bigl(n(b)a(\varpi_{v}^{-d_{v}})ka(\varpi_{v}^{d_{v}}),1\bigr)\right)
 =q_{v}^{d_{v}/2}\overline{\psi}_{v}(b)
 h_{v}\bigl(a(\varpi_{v}^{-d_{v}})k\bigr).
\end{equation}
Here $k\in K_{0,v}'[r_{v}]$. The intersection of
$a(\varpi_{v}^{-d_{v}})K_{0,v}'[r_{v}]a(\varpi_{v}^{d_{v}})$ with $N'(F_{v})$
is $n(\mathfrak{p}_{v}^{-d_{v}})$, on which $\psi_{v}$ is trivial and whose
volume gives $q_{v}^{d_{v}/2}$.
The right side of \eqref{c11.1} is smooth, compactly supported
modulo $N'(F_{v})$, and transforms on the left by
$\overline{\psi}_{v}$. The Kirillov-model theorem gives a unique vector
$W_{v}$ in the Whittaker model of $\pi_{v}$ with this restriction.
At $v\nmid\mathfrak{q}\mathfrak{m}\mathfrak{D}_{F}$ take the
spherical vector $W_{v}^{\circ}$ with $W_{v}^{\circ}(I_{3})=1$;
at every place dividing $\mathfrak{D}_{F}$
use \eqref{eq:finite-projector-function}--\eqref{c11.1}.

\subsection{The Spectral Weight}
For $W_{v}$ defined by \eqref{c11.1}, the proof of \eqref{e9.2},
applied over $F_{v}$, gives for $s=0$
\begin{equation}\label{eq11.13}
 \mathcal{H}_{W_{v}}(\sigma_{v},0)
 =\sum_{W_{v}'\in\mathfrak{B}_{\sigma_{v}}}
 \bigg|\int_{G'(F_{v})}h_{v}(g)W_{v}'(g)\,dg\bigg|^{2}. 
\end{equation}
By admissibility, the space of $K_{v}'$-fixed vectors if $r_{v}=0$,
or of vectors transforming by $\omega_{v}'(k_{22})$ under
$K_{0,v}'[r_{v}]$ otherwise, is finite-dimensional.
By newform theory, $\mathcal{H}_{W_{v}}(\sigma_{v},0)\geq0$ and
$\mathcal{H}_{W_{v}}(\sigma_{v},0)=0$ for $r_{\sigma_{v}}>r_{v}$.
For $r_{\sigma_{v}}\leq e_{v}(\mathfrak{q})$, normalize the newvector
$V_{\sigma_{v}}^{\rm new}$ to have value $1$ at $a(\varpi_{v}^{-d_{v}})$.
Then $\int h_{v}(g)V_{\sigma_{v}}^{\rm new}(g)dg=1$.
Equation~\eqref{eq11.13} and Cauchy--Schwarz give
\begin{equation}\label{fc11.12}
 \mathcal{H}_{W_{v}}(\sigma_{v},0)
 \geq \|V_{\sigma_{v}}^{\rm new}\|_{v}^{-2},
\end{equation}
with the local invariant form used in the basis sum.

At $v\nmid\mathfrak{q}\mathfrak{m}\mathfrak{D}_{F}$, the spherical
vector gives zero unless $\sigma_{v}$ is unramified. In that case,
let $W_{v}'^{\circ}$ be the unit spherical vector for the local form
used in the basis sum. The Casselman--Shalika formula gives
\begin{equation}\label{fc11.2}
 \mathcal{H}_{W_{v}}(\sigma_{v},s)
 =W_{v}^{\circ}(I_{3})|W_{v}'^{\circ}(I_{2})|^{2}
 L_{v}(1/2+s,\pi_{v}\times\sigma_{v}).
\end{equation}

\subsection{The Dual Weight}\label{sec10-dual-weight}

\begin{lemma}\label{lem11.4}
Suppose $v\nmid\mathfrak{q}\mathfrak{m}\mathfrak{D}_{F}$ and
$-1/2<\Re(\lambda)<1/2+2\Re(s)-\vartheta_{\pi_{v}}$.  Then
\begin{align*}
 \mathcal{H}_{W_{v}}(\xi_{v},\lambda;\omega_{v}',s)
 =W_{v}^{\circ}(I_{3})
 L_{v}(1/2+2s-\lambda,\widetilde{\pi}_{v}\times
 \overline{\xi}_{v}\omega_{v}\omega_{v}')L_{v}(1/2+\lambda,\xi_{v})\,
 \mathbf{1}_{r_{\xi_{v}}=0}.
\end{align*}
\end{lemma}

\begin{proof}
In \eqref{f7.4}, Lemma~\ref{lem7.3} evaluates the period of
$\widetilde{W}_{v}$ after \eqref{3.6} with $\beta=-y_{v}/z_{v}$.
Sphericity gives $\widetilde{W}_{v}(I_{3})=W_{v}(I_{3})$; integration over
$\mathcal{O}_{v}^{\times}$ makes the $y_{v}$-integral zero for ramified $\xi_{v}$.
\end{proof}

\begin{lemma}\label{lem:finite-projector-depth}
The vector $W_{v}$ defined by \eqref{c11.1} is fixed by
$n_{32}(\mathfrak{p}_{v}^{r_{v}+r_{\omega_{v}'}+3+d_{v}})$.
\end{lemma}

\begin{proof}
Take $d_{v}=0$, and let $V_{v}^{\mathrm{new}}$ be the Whittaker newvector
of $\pi_{v}$ for $\overline{\psi}_{v}$, normalized by
$V_{v}^{\mathrm{new}}(I_{3})=1$. Write
$L_{v}(s,\pi_{v})=\prod_{i=1}^{3}(1-\alpha_{i}q_{v}^{-s})^{-1}$,
allowing zero $\alpha_{i}$.
Let $\mathbb{T}_{v}$ be the spherical Hecke operator on the upper-left
$\mathrm{GL}_{2}$, supported on
$K_{v}'a(\varpi_{v}^{-1})K_{v}'$ and normalized by
$\Vol(K_{v}')^{-1}$.
By \cite[Theorem~4.1]{Miy14} and the Cauchy identity for Schur polynomials,
the vector
\begin{align*}
 V_{v}=\prod_{i=1}^{3}
 \bigl(1-\alpha_{i}q_{v}^{-1}\mathbb{T}_{v}
 +\alpha_{i}^{2}q_{v}^{-1}\pi_{v}(\diag(\varpi_{v}^{-1},\varpi_{v}^{-1},1))\bigr)V_{v}^{\mathrm{new}}
\end{align*}
satisfies
\begin{align*}
 V_{v}(\diag(x,1))
 =\int_{F_{v}}\mathbf{1}_{K_{v}'}(n(b)x)\psi_{v}(b)\,db.
\end{align*}

Fourier inversion in $n_{13}$ and the Gauss-sum identity in $n_{23}$
give $W_{v}$, up to scalar, as a sum of
$\pi_{v}(n_{13}(b)n_{23}(b')\diag(g',1))V_{v}^{\mathrm{new}}$, where
$e_{v}(b)\geq-r_{v}$, $e_{v}(b')\geq-r_{\omega_{v}'}$,
and $g'\in G'(F_{v})$ has entries in $\varpi_{v}^{-3}\mathcal{O}_{v}$
with $(g')^{-1}$ integral.
Direct conjugation, using $r_{\pi_{v}},r_{\omega_{v}'}\leq r_{v}$,
shows that each translate is fixed by
$n_{32}(\mathfrak{p}_{v}^{r_{v}+r_{\omega_{v}'}+3})$.

For $d_{v}>0$, apply the $d_{v}=0$ bound with
$\psi_{v}^{0}$ from \textsection\ref{def:unramified-additive-character}.
Conjugation by $\diag(\varpi_{v}^{2d_{v}},\varpi_{v}^{d_{v}},1)$
sends $n_{32}(t)$ to $n_{32}(\varpi_{v}^{-d_{v}}t)$.
\end{proof}

\begin{lemma}\label{lem:finite-dual-support}
Let $v\mid\mathfrak{q}\mathfrak{m}\mathfrak{D}_{F}$ and put
$g(z,y)=\diag(z,z,1)n_{32}(-y)$. If $e_{v}(y)<r_{v}+d_{v}$, then
$W_{v}(g(z,y))=0$. For every $u\in\mathcal{O}_{v}^{\times}$,
\begin{equation}\label{eq:finite-dual-unit-covariance}
 W_{v}(g(zu,y))=\overline{\omega_{v}'(u)}W_{v}(g(z,y/u)).
\end{equation}
Put $m=r_{v}+r_{\omega_{v}'}+3+d_{v}$. Then
\begin{equation}\label{eq:finite-dual-character-cutoff}
 \mathcal{H}_{W_{v}}(\xi_{v},\lambda;\omega_{v}',s)=0
 \quad\text{if}\quad r_{\xi_{v}}>r_{\omega_{v}'}+3.
\end{equation}
The contribution of $e_{v}(y)\geq m$ to \eqref{eqq9.2} is exactly
\begin{equation}\label{eq:finite-dual-tate-tail}
q_{v}^{-d_{v}/2}\Vol(K_{0,v}'[r_{v}])^{-1}
 \xi_{v}(\varpi_{v})^{m} q_{v}^{-m(1/2+\lambda)}
 L_{v}(1/2+\lambda,\xi_{v})\mathbf{1}_{r_{\xi_{v}}=0}.
\end{equation}
\end{lemma}

\begin{proof}
By \eqref{eq:finite-projector-restriction}, $W_{v}(\diag(x,1))=0$
unless $x_{21}\in\mathfrak{p}_{v}^{r_{v}+d_{v}}$.
The identity
\begin{align*}
 W_{v}(\diag(x,1)n_{13}(b))
 =\overline{\psi}_{v}(x_{21}b)W_{v}(\diag(x,1))
\end{align*}
shows that $W_{v}$ is fixed by
$n_{13}(\mathfrak{p}_{v}^{-r_{v}-2d_{v}})$. Since
\begin{align*}
 g(z,y)n_{13}(b)
 =n_{12}(yb)n_{13}(zb)g(z,y),
\end{align*}
nonvanishing forces $\psi_{v}(yb)=1$ for every such $b$, which is
equivalent to $e_{v}(y)\geq r_{v}+d_{v}$.

The right transformation laws in \eqref{eq:finite-projector-restriction} also give
$\pi_{v}(\diag(u,1,1))W_{v}=W_{v}$ and
$\pi_{v}(\diag(1,u,1))W_{v}=\overline{\omega_{v}'(u)}W_{v}$.
Now use
\begin{align*}
 g(zu,y)\diag(u^{-1},1,1)
 =g(z,y/u)\diag(1,u,1)
\end{align*}
to obtain \eqref{eq:finite-dual-unit-covariance}.
In the inner $y$-integral of \eqref{eqq9.2}, the change
$y\mapsto uy$ gives the factor $\xi_{v}(u)\overline{\omega_{v}'(u)}$
under $z\mapsto zu$.

For $e_{v}(y)=k<m$, right invariance under $n_{32}(\mathfrak{p}_{v}^{m})$
makes $u\mapsto W_{v}(g(z,\varpi_{v}^{k}u))$ invariant under
$1+\mathfrak{p}_{v}^{m-k}$. Its integral against
$\xi_{v}(u)d^{\times}u$ over $\mathcal{O}_{v}^{\times}$ therefore vanishes
if $r_{\xi_{v}}>m-k$. For $k\geq m$,
$W_{v}(g(z,y))=W_{v}(\diag(z,z,1))$.
Equation~\eqref{eq:finite-projector-restriction} gives
\begin{align*}
 W_{v}(\diag(z,z,1))
 =\Vol(K_{0,v}'[r_{v}])^{-1}q_{v}^{d_{v}/2}
 \overline{\omega_{v}'(z)}\mathbf{1}_{\mathcal{O}_{v}^{\times}}(z).
\end{align*}
Thus the $z$-integral is $\Vol(K_{0,v}'[r_{v}])^{-1}$ for unramified
$\xi_{v}$ and zero otherwise. The set $e_{v}(y)=k$ has multiplicative
measure $q_{v}^{-d_{v}/2}$, so summing the geometric series over
$k\geq m$ proves \eqref{eq:finite-dual-tate-tail}, initially for
$\Re(\lambda)>-1/2$ and then meromorphically.
Together with $k\geq r_{v}+d_{v}$, it proves
\eqref{eq:finite-dual-character-cutoff}.
\end{proof}

For $r_{v}+d_{v}\leq k<m$, set
\begin{align*}
 W_{v,k,\xi}(g):=
 \int_{\mathcal{O}_{v}^{\times}}
 W_{v}\bigl(g n_{32}(-\varpi_{v}^{k}u)\bigr)\xi_{v}(u)\,d^{\times} u,
\end{align*}
and let $\widetilde{W}_{v,k,\xi}(g)
=W_{v,k,\xi}(\widetilde{w} g^{\iota})$.
Then, meromorphically in $(s,\lambda)$,
\begin{multline}\label{eq:finite-dual-shell-expansion}
 \mathcal{H}_{W_{v}}(\xi_{v},\lambda;\omega_{v}',s)
 =\frac{q_{v}^{-d_{v}/2}\xi_{v}(\varpi_{v})^{m}}{q_{v}^{m(1/2+\lambda)}\Vol(K_{0,v}'[r_{v}])}
 L_{v}(1/2+\lambda,\xi_{v})\mathbf{1}_{r_{\xi_{v}}=0}\\
 +\sum_{k=r_{v}+d_{v}}^{m-1}
 \frac{\xi_{v}(\varpi_{v})^{k}}{q_{v}^{k(1/2+\lambda)}}
 \int_{F_{v}^{\times}}
 \widetilde{W}_{v,k,\xi}\bigl(\diag(z,1,1)\widetilde{w}\bigr)\overline{\xi}_{v}\omega_{v}\omega_{v}'(z)
 |z|_{v}^{-1/2+2s-\lambda}\,d^{\times} z.
\end{multline}
In \eqref{eqq9.2}, write $y=\varpi_{v}^{k}u$, $u\in\mathcal{O}_{v}^{\times}$,
and use \eqref{3.6} with $\beta=0$ in the $z$-integral to obtain
\eqref{eq:finite-dual-tate-tail}.

\phantomsection\label{def:normalized-finite-dual-weight}
For $v\mid\mathfrak{q}\mathfrak{m}\mathfrak{D}_{F}$, define the normalized dual weight by
\begin{align*}
 \mathcal{H}_{W_{v}}^{\sharp}(\xi_{v},\lambda;\omega_{v}',s)
 :=\frac{\mathcal{H}_{W_{v}}(\xi_{v},\lambda;\omega_{v}',s)}
 {L_{v}(1/2+2s-\lambda,\widetilde{\pi}_{v}\times
  \overline{\xi}_{v}\omega_{v}\omega_{v}')L_{v}(1/2+\lambda,\xi_{v})}.
\end{align*}

\begin{lemma}\label{lem:spherical-projector-dual-mass}
Suppose $d_{v}=0$, $\pi_{v}$ and $\omega_{v}'$ are unramified, and
$r_{v}\geq1$.
For ramified unitary $\xi_{v}$ and $\lambda\in i\mathbb{R}$,
\begin{equation}\label{eq:spherical-projector-ramified-bound}
 |\mathcal{H}_{W_{v}}^{\sharp}(\xi_{v},\lambda;\omega_{v}',0)|
 =\frac{(q_{v}+1)^{2}}{q_{v}^{2}(q_{v}-1)}\mathbf{1}_{r_{v}=r_{\xi_{v}}=1}
 \ll q_{v}^{-1}.
\end{equation}
For an unramified unitary character $\xi_{v}$,
\begin{equation}\label{eq:spherical-projector-unramified-bound}
 |\mathcal{H}_{W_{v}}^{\sharp}(\xi_{v},\lambda;\omega_{v}',0)|\ll q_{v}^{r_{v}/2},
 \qquad \lambda\in i\mathbb{R}.
\end{equation}
Consequently, extending each character $\nu$ of $\mathcal{O}_{v}^{\times}$
by $\nu(\varpi_{v})=1$, we have
\begin{equation}\label{eq:spherical-projector-dual-mass}
 \sum_{\nu\in\widehat{\mathcal{O}_{v}^{\times}}}
 q_{v}^{(1+\varepsilon)r_{\nu}}
 \sup_{\substack{\rho\ \mathrm{unramified\ unitary}\\
                 \lambda\in i\mathbb{R}}}
 |\mathcal{H}_{W_{v}}^{\sharp}(\nu\rho,\lambda;\omega_{v}',0)|
 \ll_{\varepsilon} q_{v}^{r_{v}/2}+\mathbf{1}_{r_{v}=1}q_{v}^{1+\varepsilon}.
\end{equation}
All constants are uniform in the unramified unitary data.
For unitary $\xi_{v}$ and $|\lambda-1/2|,|\Re(s)|\leq\delta<10^{-2}$,
we have
\begin{equation}\label{eq:spherical-projector-residue-contour}
 |\mathcal{H}_{W_{v}}^{\sharp}(\xi_{v},\lambda;\omega_{v}',s)|
 \ll_{\delta} q_{v}^{(1/2-\Re(\lambda))r_{v}+\vartheta_{\pi_{v}}+9\delta}.
\end{equation}
\end{lemma}

\begin{proof}
For $g(z,y)=\diag(z,z,1)n_{32}(-y)$, use
$V_{v}$ from Lemma~\ref{lem:finite-projector-depth} and
$W_{v}$ from \eqref{eq:finite-projector-function}:
\begin{align*}
 W_{v}(g(z,y))=\Vol(K_{0,v}'[r_{v}])^{-1}
 \mathbf{1}_{e_{v}(y)\geq r_{v}}V_{v}(g(z,y)).
\end{align*}
Let $\alpha=(\alpha_{1},\alpha_{2},\alpha_{3})$ be the Satake parameters,
and write $s_{i,j,k}(\alpha)$ for the Schur polynomial; see \eqref{12.6}. Put
\begin{align*}
 F_{a,b}(z,y)=V_{v}^{\mathrm{new}}
 \bigl(g(z,y)\diag(\varpi_{v}^{-a},\varpi_{v}^{-b},1)\bigr).
\end{align*}
The right $\mathbb{T}_{v}$-action gives
$q_{v}F_{a,b+1}+F_{a+1,b}$ for $a\leq b$; the second term is zero
when $a=b$. Expanding the product defining $V_{v}$ in
Lemma~\ref{lem:finite-projector-depth} therefore gives
\begin{equation}\label{eq:spherical-projector-hecke-expansion}
 V_{v}(g(z,y))=
 \sum_{0\leq a\leq b\leq3}(-1)^{a+b}q_{v}^{-a}
 s_{(2^{a},1^{b-a})}(\alpha)F_{a,b}(z,y).
\end{equation}
Here $(2^{a},1^{b-a})$ lists $a$ twos and $b-a$ ones;
trailing zero indices are omitted.
For $k=e_{v}(y)$ and $\ell=e_{v}(z)$, if $k\geq b$ then
$F_{a,b}(z,y)=0$ for $\ell<b$, and otherwise
\begin{align*}
 F_{a,b}(z,y)=
 q_{v}^{-(\ell-a)}s_{\ell-a,\ell-b,0}(\alpha).
\end{align*}
If $k<b$, then $F_{a,b}(z,y)=0$ unless $a\leq k$ and $\ell\geq2k-b$.
When these inequalities hold, the Iwasawa decomposition of the lower-right block gives
\begin{align*}
 F_{a,b}(z,y)=\overline{\psi}_{v}(-z/y)s_{1,1,1}(\alpha)^{k-b}
 q_{v}^{-(\ell+b-a-k)}
 s_{\ell+b-a-k,\ell+b-2k,0}(\alpha).
\end{align*}
For $j=r_{\xi_{v}}>0$, write $y=\varpi_{v}^{k}u$ and integrate
$F_{a,b}(z,\varpi_{v}^{k}u)$ against $\xi_{v}(u)d^{\times}u$.
For $k\geq b$ the integrand apart from $\xi_{v}(u)$ is constant in $u$.
For $k<b$ the Gauss integral from $\overline{\psi}_{v}(-z/y)$
vanishes unless $\ell=k-j$.
Hence $1\leq r_{v}\leq k<b\leq3$ and $j\leq b-k$.
The two possible cases other than $k=j=1$ cancel in
\eqref{eq:spherical-projector-hecke-expansion}:
for $(k,j)=(1,2)$ use $s_{2,1,1}(\alpha)=s_{1,1,1}(\alpha)s_{1,0,0}(\alpha)$;
for $(k,j)=(2,1)$ use
$s_{2,0,0}(\alpha)-s_{1,0,0}(\alpha)^{2}+s_{1,1,0}(\alpha)=0$.
For $k=j=1$, the terms with $b=2,3$ leave respectively
$q_{v}^{-1}$ and $q_{v}^{-2}$, using
$s_{1,0,0}(\alpha)s_{1,1,0}(\alpha)-s_{2,1,0}(\alpha)=s_{1,1,1}(\alpha)$ and
$s_{2,2,1}(\alpha)=s_{1,1,1}(\alpha)s_{1,1,0}(\alpha)$.
The surviving integral is
\begin{align*}
 \mathcal{H}_{W_{v}}(\xi_{v},\lambda;\omega_{v}',s)
 =\frac{\xi_{v}(\varpi_{v})q_{v}^{-(1/2+\lambda)}}{\Vol(K_{0,v}'[1])}
 (q_{v}^{-1}+q_{v}^{-2})
 \int_{\mathcal{O}_{v}^{\times}}
 \xi_{v}(u)\overline{\psi}_{v}(-\varpi_{v}^{-1}u^{-1})\,d^{\times} u.
\end{align*}
The integral has modulus $q_{v}^{1/2}/(q_{v}-1)$.
Both local $L$-factors are one for ramified $\xi_{v}$.
Since $\Vol(K_{0,v}'[1])=(q_{v}+1)^{-1}$, this proves
\eqref{eq:spherical-projector-ramified-bound} and the vanishing.

For unramified $\xi_{v}$, put $ t=\xi_{v}(\varpi_{v})q_{v}^{-(1/2+\lambda)}$, $ X=\overline{\xi}_{v}\omega_{v}'(\varpi_{v})
 q_{v}^{-(1/2+2s-\lambda)}$, and  
\begin{align*}
D(X)=1-s_{1,1,0}(\alpha)X+s_{1,0,0}(\alpha)s_{1,1,1}(\alpha)X^{2}-s_{1,1,1}(\alpha)^{2}X^{3}.
\end{align*}
Then $D(X)$ is the reciprocal of the degree-three Euler factor
in $\mathcal{H}_{W_{v}}^{\sharp}$. For $k\geq2$, the displayed Schur formulas give
$V_{v}(g(z,y))=\mathbf{1}_{\ell=0}$.
For $k=1$, integrating $\overline{\psi}_{v}(-z/y)$ over $e_{v}(y)=1$ gives $1$ for
$\ell\geq1$, $-1/(q_{v}-1)$ for $\ell=0$, and zero for $\ell<0$.
Summing the Schur formulas gives $A_{q_{v}}(X)/D(X)$ for the
$z$-integral, after integrating over the unit parts of $y$ and $z$, where
\begin{align*}
 A_{q_{v}}(X):=
 s_{1,1,1}(\alpha)^{2}X^{3}+q_{v}^{-1}(1-s_{1,1,0}(\alpha)X)+q_{v}^{-2}
 +\left(1-\frac{q_{v}+1}{q_{v}(q_{v}-1)}\right)D(X).
\end{align*}
Use
\begin{align*}
\sum_{n\geq0}s_{n,n,0}(\alpha)X^{n}=D(X)^{-1},\qquad 
 \sum_{n\geq0}s_{n+1,n,0}(\alpha)X^{n}=\frac{s_{1,0,0}(\alpha)-s_{1,1,1}(\alpha)X}{D(X)}.
\end{align*}
Hence, meromorphically in $(s,\lambda)$,
\begin{equation}\label{eq:spherical-projector-exact-transform}
 \mathcal{H}_{W_{v}}^{\sharp}(\xi_{v},\lambda;\omega_{v}',s)=\Vol(K_{0,v}'[r_{v}])^{-1}
 \begin{cases}
 t(1-t)A_{q_{v}}(X)+t^{2}D(X),&r_{v}=1,\\
 t^{r_{v}}D(X),&r_{v}\geq2.
 \end{cases}
\end{equation}
For unitary spherical $\pi_{v}$, $|s_{1,1,1}(\alpha)|=1$ and
$|s_{1,0,0}(\alpha)|,|s_{1,1,0}(\alpha)|\leq3q_{v}^{1/2}$. On $s=0$, $\Re(\lambda)=0$,
$|t|=|X|=q_{v}^{-1/2}$ and $|D(X)|+|A_{q_{v}}(X)|\ll1$.
Using $\Vol(K_{0,v}'[r_{v}])^{-1}=(q_{v}+1)q_{v}^{r_{v}-1}$ proves
\eqref{eq:spherical-projector-unramified-bound}.
There are $q_{v}-2$ characters of $\mathcal{O}_{v}^{\times}$ of conductor exponent $1$.
Their sum and the unramified term prove
\eqref{eq:spherical-projector-dual-mass}.

For $|\lambda-1/2|\leq\delta$, $|t|=q_{v}^{-1/2-\Re(\lambda)}<1$ and
$|X|\leq q_{v}^{3\delta}$. Unitarity gives $|s_{1,1,1}(\alpha)|=1$ and
$|s_{1,0,0}(\alpha)|,|s_{1,1,0}(\alpha)|\leq3q_{v}^{\vartheta_{\pi_{v}}}$, so
$|D(X)|+|A_{q_{v}}(X)|\ll q_{v}^{\vartheta_{\pi_{v}}+9\delta}$.
For unramified $\xi_{v}$, the factor
$\Vol(K_{0,v}'[r_{v}])^{-1}|t|^{r_{v}}\ll q_{v}^{(1/2-\Re(\lambda))r_{v}}$
in \eqref{eq:spherical-projector-exact-transform} gives
\eqref{eq:spherical-projector-residue-contour} with this polynomial bound.
For ramified $\xi_{v}$,
$|\mathcal{H}_{W_{v}}^{\sharp}|\ll q_{v}^{-\Re(\lambda)}
\mathbf{1}_{r_{\xi_{v}}=r_{v}=1}$ by the Gauss integral.
\end{proof}

\begin{prop}\label{lem11.6}
Let $\pi$ be as in Theorem~\ref{thm1.3}, let $\omega_{v}'=\mathbf{1}$,
and let $v\mid\mathfrak{q}\mathfrak{m}\mathfrak{D}_{F}$. Take $W_{v}$ from
\eqref{eq:finite-projector-function}--\eqref{c11.1}.
For unramified unitary $\xi_{v}$ and $\Re(\lambda)=0$,
\begin{equation}\label{11.37}
 |\mathcal{H}_{W_{v}}^{\sharp}(\xi_{v},\lambda;\mathbf{1},0)|
 \ll_{F} q_{v}^{r_{v}/2}.
\end{equation}
For $d_{v}=0$ the implied constant is absolute; the factors
with $d_{v}>0$ depend only on $F$.
For $r_{v}=0$ and any unitary $\xi_{v}$, we also have
\begin{equation}\label{eq:finite-zero-level-dual-bound}
 |\mathcal{H}_{W_{v}}^{\sharp}(\xi_{v},\lambda;\mathbf{1},s)|
 \ll_{F}1,\qquad |\Re(s)|\leq1/10,\quad |\Re(\lambda)|\leq1.
\end{equation}
For ramified $\xi_{v}$, the support
restriction and formula are
\eqref{eq:finite-dual-character-cutoff} and
\eqref{eq:finite-dual-shell-expansion}.
\end{prop}

\begin{proof}
First take $d_{v}=0$ and $r_{v}\geq1$. Write
$L_{v}(s,\pi_{v})=\prod_{i=1}^{3}(1-\alpha_{i}q_{v}^{-s})^{-1}$, allowing
$\alpha_{i}=0$.
If at most one $\alpha_{i}$ is nonzero, then $L_{v}(s,\pi_{v})$ is either
$1$ or has degree one with $r_{\pi_{v}}\geq2$ and $|\alpha_{i}|\leq1$.
Since $r_{v}\geq r_{\pi_{v}}$, Lemma~\ref{lem:finite-projector-full-Tate-tail}
gives \eqref{11.37} with bound
$\Vol(K_{0,v}'[r_{v}])^{-1}q_{v}^{-r_{v}/2}\ll q_{v}^{r_{v}/2}$.

If exactly two $\alpha_{i}$ are nonzero, trivial central character excludes the type consisting
of one ramified character and two unramified characters. The remaining
local type is an unramified unitary Steinberg twist plus an unramified
character. Lemma~\ref{lem:maximal-Steinberg-unitary-dual} gives
\eqref{11.37} for this type.
If all three $\alpha_{i}$ are nonzero, $\pi_{v}$ is unramified and
\eqref{eq:spherical-projector-unramified-bound} applies.

For $r_{v}=0$, $\pi_{v}$ is unramified and $v\mid\mathfrak{D}_{F}$.
Work first with $\psi_{v}^{0}$
and its measures, so $W_{v}^{0}=V_{v}$ from
Lemma~\ref{lem:finite-projector-depth}.
For ramified $\xi_{v}$, the Gauss integral in the proof of
Lemma~\ref{lem:spherical-projector-dual-mass} leaves only
$0\leq e_{v}(y)<b\leq3$ and
$e_{v}(z)=e_{v}(y)-r_{\xi_{v}}\geq2e_{v}(y)-b$.
These finitely many shells give \eqref{eq:finite-zero-level-dual-bound}
for $W_{v}^{0}$, since both local $L$-factors are one and $q_{v}$ is fixed
at $v\mid\mathfrak{D}_{F}$.

Take $t,X,D(X)$ as defined in the proof of
Lemma~\ref{lem:spherical-projector-dual-mass} for unramified $\xi_{v}$.
Initially take $|t|<1$ and $|X|<q_{v}^{-1}$, so the series below
converge absolutely.
The contribution of $e_{v}(y)\geq1$ to the normalized transform is
obtained by multiplying \eqref{eq:spherical-projector-exact-transform}
at $r_{v}=1$ by $\Vol(K_{0,v}'[1])$, giving
$t(1-t)A_{q_{v}}(X)+t^{2}D(X)$.
For $e_{v}(y)=0$, only the terms $F_{0,b}$, $0\leq b\leq3$, occur
in \eqref{eq:spherical-projector-hecke-expansion}.
For $b\geq1$, integration over the units of $y$ gives the factor
$1$ for $e_{v}(z)\geq0$, $-1/(q_{v}-1)$ for $e_{v}(z)=-1$, and zero
for $e_{v}(z)<-1$.
Thus the $z$-integral in \eqref{eqq9.2} for $F_{0,0}$ is $D(X)^{-1}$;
for $F_{0,b}$, $1\leq b\leq3$, it is
\begin{align*}
 \frac{q_{v}^{-b}}{s_{1,1,1}(\alpha)^{b}}
 \bigg[X^{-b}\Big(D(X)^{-1}
 -\sum_{j=0}^{b-1}s_{j,j,0}(\alpha)X^{j}\Big)
 -\frac{s_{b-1,b-1,0}(\alpha)}{q_{v}-1}X^{-1}\bigg],
\end{align*}
using $\sum_{j\geq0}s_{j,j,0}(\alpha)X^{j}=D(X)^{-1}$.
Thus normalization multiplies these series by $(1-t)D(X)$,
leaving a finite Laurent polynomial in $t,X$.
Its coefficients involve only bounded-degree Schur polynomials,
powers of $q_{v}$, and powers of $s_{1,1,1}(\alpha)^{-1}$.
Since $|\alpha_{i}|\leq q_{v}^{1/2}$ and
$|s_{1,1,1}(\alpha)|=1$, they are bounded at this fixed place.
Also $t,X$ and their inverses are bounded in the region in
\eqref{eq:finite-zero-level-dual-bound}.
The identity continues meromorphically, proving
\eqref{eq:finite-zero-level-dual-bound} for $W_{v}^{0}$.

Finally, \eqref{eq:maximal-positive-residue-different-transport}
contributes $q_{v}^{d_{v}/2}$ on $\Re(\lambda)=0$, and a bounded factor
throughout the region in \eqref{eq:finite-zero-level-dual-bound}.
\end{proof}

\subsection{Change of Additive Character}\label{sec:finite-character-change}
\phantomsection\label{def:unramified-additive-character}
For $v<\infty$, put $\psi_{v}^{0}(x):=\psi_{v}(\varpi_{v}^{-d_{v}}x)$ 
which has conductor $\mathcal{O}_{v}$. Let $dt_{v}^{0}$ be self-dual for
$\psi_{v}^{0}$; a superscript $0$ on the other Haar measures means that
$dt_{v}^{0}$ is used in their definitions. Then
\begin{equation}\label{eq:finite-character-measures}
 dt_{v}=q_{v}^{-d_{v}/2}dt_{v}^{0},\qquad
 d^{\times}t_{v}=q_{v}^{-d_{v}/2}d^{\times}t_{v}^{0},\qquad
 dg_{v}'=q_{v}^{-3d_{v}/2}dg_{v}^{\prime,0}.
\end{equation}
In particular,
$\Vol(\mathcal{O}_{v},dt_{v}^{0})=\Vol(\mathcal{O}_{v}^{\times},d^{\times}t_{v}^{0})=1$.
For a unitary representation $\rho_{v}$ of $\mathrm{GL}_{n}(F_{v})$,
the change of additive character gives
\begin{equation}\label{eq:different-epsilon-factor}
 |\epsilon(s,\rho_{v},\psi_{v})|
 =q_{v}^{nd_{v}(1/2-\Re(s))}|\epsilon(s,\rho_{v},\psi_{v}^{0})|.
\end{equation}

For a section $f$ of
$\chi_{v}|\cdot|_{v}^{\lambda}\boxplus\chi_{v}|\cdot|_{v}^{-\lambda}$,
normalized induction and $x=\varpi_{v}^{-d_{v}}x_{0}$ give
\begin{flalign}\label{eq:different-GL2-Jacquet}
 &&
 \int_{F_{v}}f(w'n(x)a(\varpi_{v}^{-d_{v}}))\overline{\psi_{v}(x)}\,dx
 =\chi_{v}(\varpi_{v}^{-d_{v}})q_{v}^{-d_{v}\lambda}
 \int_{F_{v}}f(w'n(x))\overline{\psi_{v}^{0}(x)}\,dx^{0},&
\end{flalign}
as a meromorphic identity in $\lambda$.

For the projector vector in \eqref{c11.1}, let
$W_{v}^{0}$ be given by \eqref{eq:finite-projector-function}--\eqref{c11.1}
with $\psi_{v}^{0}$ and the superscript-$0$ measures.
In the Whittaker projection make the change $b\mapsto\varpi_{v}^{-d_{v}}b^{0}$; then
$db=q_{v}^{d_{v}/2}db^{0}$.
The convolution contributes $q_{v}^{3d_{v}/2}$, so injectivity of the
restriction to the mirabolic gives
\begin{align*}
 W_{v}(g)=q_{v}^{2d_{v}}W_{v}^{0}\bigl(\diag(\varpi_{v}^{2d_{v}},\varpi_{v}^{d_{v}},1)
 g\diag(\varpi_{v}^{-2d_{v}},\varpi_{v}^{-d_{v}},1)\bigr).
\end{align*}
This conjugation sends $g(z,y)$ to $g(z,\varpi_{v}^{-d_{v}}y)$. Changing
$y\mapsto\varpi_{v}^{d_{v}}y^{0}$ in \eqref{eqq9.2} and using
\eqref{eq:finite-character-measures} gives
\begin{equation}\label{eq:maximal-positive-residue-different-transport}
 \mathcal{H}_{W_{v}}(\xi_{v},\lambda;\omega_{v}',s)
 =\xi_{v}(\varpi_{v})^{d_{v}} q_{v}^{d_{v}(1/2-\lambda)}
 \mathcal{H}_{W_{v}^{0}}(\xi_{v},\lambda;\omega_{v}',s).
\end{equation}
The same identity holds for $\mathcal{H}^{\sharp}$, since the Euler
factors do not depend on the additive character.
At $(\xi_{v},\lambda)=(\mathbf{1},1/2)$ the transforms agree.

For the type~\RNum{2} weights, let $V^{0}$ be a $\psi_{v}^{0}$-Whittaker
function on $G(F_{v})$, and put
\begin{align*}
V(g)=V^{0}(\diag(\varpi_{v}^{2d_{v}},\varpi_{v}^{d_{v}},1)g).
\end{align*}
\begin{equation}\label{eq:different-raw-weight}
\mathcal{H}_{V}^{\psi_{v}}(\sigma_{v},\mathbf{s},\mathbf{1})
=q_{v}^{d_{v}(3s_{1}+s_{2}-2)}
\mathcal{H}_{V^{0}}^{\psi_{v}^{0}}(\sigma_{v},\mathbf{s},\mathbf{1}).
\end{equation}
In the two local periods use
$x\mapsto\diag(\varpi_{v}^{-2d_{v}},\varpi_{v}^{-d_{v}})x$ and
$a\mapsto\varpi_{v}^{-d_{v}}a$ with the corresponding unit bases and
self-dual measures.
Write $\mathcal{J}_{\psi_{v}}$ and $\mathcal{J}_{\psi_{v}^{0}}$ for the
Jacquet integrals with characters $\psi_{v}$ and $\psi_{v}^{0}$,
respectively, and the corresponding self-dual measures.
For the spherical section $f_{v}^{\circ}$ in
$\pi_{v}=\boxplus_{j=1}^{3}|\cdot|_{v}^{\nu_{j}}$, normalized by
$f_{v}^{\circ}(I_{3})=1$, conjugating the unipotent variable
in the Jacquet integral gives
\begin{equation}\label{eq:different-jacquet-scalar}
(\mathcal{J}_{\psi_{v}}f_{v}^{\circ})(g)
=q_{v}^{d_{v}(1/2+\nu_{3}-\nu_{1})}
(\mathcal{J}_{\psi_{v}^{0}}f_{v}^{\circ})
(\diag(\varpi_{v}^{2d_{v}},\varpi_{v}^{d_{v}},1)g).
\end{equation}
Equations~\eqref{eq:different-raw-weight} and
\eqref{eq:different-jacquet-scalar} remain valid after applying the same
right translation to the two $\mathrm{GL}_{3}$ Whittaker functions.
For the dual weights, replace $\mathbf{s}$ by $\mathbf{s}^{\vee}$.

\section{Type \RNum{2}: Archimedean Weights and Dual Transforms}\label{sec11}
We choose spherical spectral weights for type~\RNum{2} reciprocity and compute
their dual transforms using the Hecke integrals of all $K_{v}'$-types.
The real-place bounds in \textsection\ref{sec:G-normalized-real-transform}
and Appendix~\ref{app:central-polynomial-transform} are used in
\textsection\ref{sec15} to prove the fifth-moment estimate of
Theorem~\ref{thm1.5}, with polynomial dependence on the Archimedean cutoff.

Let $v\mid\infty$, and let $\pi_{v}$ and $\sigma_{v}$ be unitary generic
representations of $G(F_{v})$ and $G'(F_{v})$, respectively, in the
Whittaker models of \textsection\ref{sec9}. Let $\eta_{v}$ be a unitary
character of $F_{v}^{\times}$, and let $\mathfrak{B}_{\sigma_{v}}$ be an
orthonormal basis of $K_{v}'$-type vectors for the Kirillov inner
product. For a Whittaker vector $W_{v}$ of $\pi_{v}$, use the local
Rankin--Selberg integrals of \textsection\ref{sec3.2} to define
\begin{equation}\label{fc10.1}
\mathcal{H}_{W_{v}}(\sigma_{v},\mathbf{s},\eta_{v})
:=\sum_{W_{v}'\in\mathfrak{B}_{\sigma_{v}}}
\Psi_{v}(1/2+s_{1},W_{v},W_{v}')\cdot \overline{\Psi_{v}(1/2+\overline{s_{2}},W_{v}',\overline{\eta_{v}})}.
\end{equation}
Initially take $\Re(s_{i})\ggg1$, $i=1,2$, and continue
meromorphically to $|\Re(s_{i})|\leq1/4$.

Put $\mathbf{s}^{\vee}=((s_{2}-s_{1})/2,(3s_{1}+s_{2})/2)$. Define the dual
weight by right translation:
\begin{equation}\label{fc10.2}
\mathcal{H}_{W_{v}}^{\vee}(\sigma_{v},\mathbf{s}^{\vee},\eta_{v})
:=\mathcal{H}_{\pi_{v}(w_{2})W_{v}}(\sigma_{v},\mathbf{s}^{\vee},\eta_{v}),
\qquad \pi_{v}(w_{2})W_{v}(g)=W_{v}(gw_{2}).
\end{equation}
For $\omega_{\sigma_{v}}=\eta_{v}=\mathbf{1}$, left multiplication by
$\diag(1,-1,1)$ and $\diag(1,-1)$ reverses the Whittaker characters.
The changes $g'\mapsto\diag(1,-1)g'$ and $y\mapsto-y$ leave the periods
unchanged, so the weights agree with those in \textsection\ref{sec3.2}.

\subsection{Spherical Weights}\label{sec11.1}
Suppose $F_{v}\simeq\mathbb{R}$ and
$\omega_{v}=\omega_{v}'=\eta_{v}=\mathbf{1}$. Write
$\pi_{v}=\boxplus_{j=1}^{3}|\cdot|_{v}^{\nu_{j}}$, with
$\nu_{1}+\nu_{2}+\nu_{3}=0$, and choose its spherical Whittaker vector
$W_{v}^{\circ}$.
For bi-$K_{v}'$-invariant $h_{v}\in C_{c}^{\infty}(G'(F_{v}))$, set
\begin{equation}\label{fc11.3}
W_{v}\left(\diag(x_{v},1)\right):=\int_{G'(F_{v})}W_{v}^{\circ}\left(\diag(x_{v}y_{v},1)\right)h_{v}(y_{v}^{-1})|\det y_{v}|_{v}^{s_{1}}dy_{v}.
\end{equation}

In \eqref{fc10.1}, substitute $x_{v}\mapsto x_{v}y_{v}^{-1}$ and then
$y_{v}\mapsto y_{v}^{-1}$. The determinant factors cancel, and right
invariance of the quotient measure replaces $W_{v}'$ in the first
period by $\sigma_{v}(h_{v})W_{v}'$. Bi-$K_{v}'$-invariance of $h_{v}$ restricts the basis sum in \eqref{fc10.1}
to $(\sigma_{v})^{K_{v}'}$.

For unitary spherical
$\sigma_{v}=|\cdot|_{v}^{\nu_{\sigma_{v}}}\boxplus|\cdot|_{v}^{-\nu_{\sigma_{v}}}$,
the spherical transform $\mathcal{S}h_{v}(\nu_{\sigma_{v}})$ is the eigenvalue of
$\sigma_{v}(h_{v})$ on its $K_{v}'$-invariant line:
\begin{equation}\label{eq10.4}
\mathcal{S}h_{v}(\nu)
=\int_{\mathbb{R}^{\times}}\int_{\mathbb{R}^{\times}}\int_{\mathbb{R}}
h_{v}\left(\diag(za,z)
n(b)\right)
|a|^{1/2+\nu}\,db\,d^{\times} a\,d^{\times} z.
\end{equation}

\begin{lemma}\label{lem:real-rs-evaluation}
Let $\beta_{v}$ be the normalization constant of $W_{v}^{\circ}$ in
\eqref{eq:real-jacquet-double-mellin}. Then, as a meromorphic identity
in $\mathbf{s}$,
\begin{equation}\label{eq:real-weight-evaluated}
\mathcal{H}_{W_{v}}(\sigma_{v},\mathbf{s},\mathbf{1})
=\frac{\beta_{v}}{\pi L_{v}(1,\sigma_{v},\Ad)}
L_{v}(1/2+s_{1},\pi_{v}\times\sigma_{v})L_{v}(1/2+s_{2},\sigma_{v})
\mathcal{S}h_{v}(\nu).
\end{equation}
\end{lemma}

\begin{proof}
Take $W_{v}'^{\circ}(\diag(a,1))=|a|^{1/2}K_{\nu}(2\pi|a|)$.
For $\Re(s_{1})$ and $\Re(s_{2})$ sufficiently large, the $K$-Bessel
integrals give $\Psi_{v}(1/2+s_{2},W_{v}'^{\circ},\mathbf{1})
=2^{-1}L_{v}(1/2+s_{2},\sigma_{v})$ and 
\begin{equation}\label{eq:real-bessel-norm}
\int_{\mathbb{R}^{\times}}|W_{v}'^{\circ}(\diag(a,1))|^{2}\,d^{\times}a=4^{-1}\pi L_{v}(1,\sigma_{v},\Ad).
\end{equation}
Stade's formula \cite[Theorem~3.4]{Sta01}, with our Haar measure, gives
\begin{equation}\label{eq:real-rs-evaluated}
\Psi_{v}(1/2+s_{1},W_{v}^{\circ},W_{v}'^{\circ})
=2^{-1}\beta_{v}L_{v}(1/2+s_{1},\pi_{v}\times\sigma_{v}).
\end{equation}
Dividing the product of these periods by the squared norm
\eqref{eq:real-bessel-norm} proves \eqref{eq:real-weight-evaluated}.
\end{proof}

By spherical Paley--Wiener, for compactly supported $h_{v}$,
$\mathcal{S}h_{v}$ is even, entire and rapidly decreasing with all
derivatives on closed vertical strips.

\subsection{\texorpdfstring{$K_{v}'$-isotypical Vectors and Hecke Integrals}{K\_v'-isotypical Vectors and Hecke Integrals}}
Assume $F_{v}\simeq\mathbb{R}$ and let $\sigma_{v}$ be generic unitary
with trivial central character. Write it as a subrepresentation of
$\sgn^{\epsilon}|\cdot|_{v}^{\nu}\boxplus
\sgn^{\epsilon}|\cdot|_{v}^{-\nu}$, with $\epsilon\in\{0,1\}$ and
$\nu\in i\mathbb{R}_{\geq0}$ or $\nu>0$.
For a discrete series, $2\nu$ is a positive odd integer and its
lowest weight is $2\nu+1$; take $\epsilon=0$ in this case.

\subsubsection{Construction of $\mathfrak{B}_{\sigma_{v}}$}\label{sec10.2.1}
Put
\begin{align*}
\mathcal{N}_{\sigma_{v}}:=\begin{cases}
2\mathbb{Z},\ \ & \text{if $\nu\in i\mathbb{R}$ or $0<\nu<1/2$,}\\
\{m\in 2\mathbb{Z}:\ |m|\geq 2\nu+1\},\ \ & \text{if $2\nu$ is a positive odd integer.}
\end{cases}
\end{align*}

Let $W_{\kappa,\nu}$ be the standard Whittaker function,
normalized by $W_{\kappa,\nu}(z)\sim e^{-z/2}z^{\kappa}$ at infinity.

By \cite[(2.16), Lemmas~4--5]{BM05}, take
\begin{equation} \label{10.8}
W_{n}\left(\diag(a,1)\right):=\frac{\sgn(a)^{\epsilon}i^{n\sgn(a)}W_{\sgn(a)n,\nu}(4\pi|a|)}{\sqrt{\Gamma(1/2-\nu+\sgn(a)n)\Gamma(1/2+\nu+\sgn(a)n)}}.
\end{equation}
Take the positive square root, and set the expression to zero if
one of the Gamma factors in the denominator has a pole.
The vector $W_{n}$ has rotation weight $2n$:
\begin{align*}
\sigma_{v}\left(\begin{pmatrix}
\cos\theta & \sin\theta\\
-\sin\theta & \cos\theta
\end{pmatrix}\right)W_{n}=e^{2in \theta}W_{n},\ \ 2n\in \mathcal{N}_{\sigma_{v}}.
\end{align*}

Thus we may take
\begin{align*}
 \mathfrak{B}_{\sigma_{v}}=\{W_{n}:\ 2n\in\mathcal{N}_{\sigma_{v}}\}.
\end{align*}

\subsubsection{The Hecke integral of $K_{v}'$-isotypical vectors}\label{sec:real-hecke-integral}
\phantomsection\label{def:real-hecke-integral}
For $\Re(\lambda)>-1/2+|\Re(\nu)|$, define the Hecke integral
\begin{align*}
\Psi(W_{n};\lambda):=\int_{\mathbb{R}^{\times}}
W_{n}(\diag(a,1))|a|^{\lambda}\,d^{\times} a.
\end{align*}
The Whittaker expansion at zero gives its meromorphic
continuation. Applying \cite[\textsection 6.10, (14)]{Bat53} and
Pfaff's identity to the two half-lines in \eqref{10.8} gives
\begin{multline*}
\Psi(W_{n};\lambda)=\frac{\Gamma(\lambda+1/2+\nu)
\Gamma(\lambda+1/2-\nu)}{(4\pi)^{\lambda}}\\
\sum_{\varepsilon\in\{\pm1\}}\frac{\varepsilon^{\epsilon}i^{\varepsilon n}{}_{2}F_{1}(\lambda+\nu+1/2,\lambda-\nu+1/2;\lambda-\varepsilon n+1;1/2)}{\Gamma(\lambda-\varepsilon n+1)\sqrt{\Gamma(1/2-\nu+\varepsilon n)\Gamma(1/2+\nu+\varepsilon n)}}.
\end{multline*}
Interpret ${}_{2}F_{1}(a,b;c;1/2)/\Gamma(c)$ by continuation at
nonpositive integral $c$. Terms excluded by the convention in
\eqref{10.8} are zero.

\begin{lemma}\label{lemm11.3}
Put $r=\lambda_{1}+s_{1}-s_{2}$. For fixed unitary $\sigma_{v}$, the sum
\begin{multline*}
\mathcal{S}_{\sigma_{v}}(\lambda_{1},\lambda_{2},\mathbf{s})
:=\sum_{2n\in\mathcal{N}_{\sigma_{v}}}
\frac{\Gamma(|n|+(1-r)/2)}{\Gamma(|n|+(1+r)/2)}\\
\Psi(W_{n};-\lambda_{2}+(s_{2}-s_{1})/2)
\overline{\Psi(W_{n};(3\overline{s_{1}}+\overline{s_{2}})/2)}
\end{multline*}
converges absolutely and locally uniformly in
\begin{equation}\label{eq:normalized-ktype-region}
\begin{cases}
 |\Re(-\lambda_{2}+(s_{2}-s_{1})/2)|,
 |\Re((3s_{1}+s_{2})/2)|<1/2-\vartheta_{\sigma_{v}},\   \Re(r)<1\\
 \delta:=\Re(r)-|\Re(-\lambda_{2}+(s_{2}-s_{1})/2)|
 -|\Re((3s_{1}+s_{2})/2)|>0.
\end{cases}
\end{equation}
\end{lemma}

\begin{proof}
For fixed $\sigma_{v}$, Fourier inversion and
\cite[(2.16), \textsection4]{BM05} in the principal-series case,
and the Laguerre generating function in the discrete-series case,
give, locally uniformly for real $u$ with $|u|<1/2-\vartheta_{\sigma_{v}}$,
\begin{equation}\label{eq:hecke-vertical}
 |\Psi(W_{n};u+it)|\ll
 (1+|t|)^{l}(1+|n|)^{-1/2+|u|}\log(2+|n|),
 \qquad t\in\mathbb{R},
\end{equation}
for some $l\geq0$. For principal series, use $x\mapsto\tan\theta$
and the powers at $0,\pm\pi/2$ in the resulting Fourier integral;
for discrete series, use Cauchy's coefficient formula and Stirling's
formula. Stirling's formula now bounds the summand locally uniformly
by $(1+|n|)^{-1-\delta}\log^{2}(2+|n|)$, which is summable.
\end{proof}

\subsection{A Fourier Transform of the Spherical Test Function}
For $x\in F_{v}\simeq\mathbb{R}$, define
\begin{equation}\label{eq10.17}
f_{v}(x):=\int_{0}^{\infty}\int_{0}^{\infty}h_{v}\left(n(x)\diag(y^{1/2}z,y^{-1/2}z)\right)d^{\times}zy^{s_{2}}d^{\times}y.
\end{equation} 

For compactly supported $h_{v}$, we have $f_{v}\in C_{c}^{\infty}(\mathbb{R})$,
with holomorphic dependence on $s_{2}$. Bi-$K_{v}'$-invariance gives
$f_{v}(x)=f_{v}(-x)$. Let
\begin{align*}
\mathcal{M}f_{v}(\lambda):=\int_{0}^{\infty} f_{v}(x)x^{\lambda}d^{\times}x,
\qquad \Re(\lambda)>0,
\end{align*}
be its Mellin transform, continued meromorphically using its expansion
at zero. Let 
\begin{equation}\label{c11.20}
\widehat{f}_{v}(y):=\int_{\mathbb{R}}f_{v}(x)e^{2\pi i yx}dx=\int_{\mathbb{R}}f_{v}(x)e^{-2\pi i yx}dx
\end{equation}
be the Fourier transform of $f_{v}$. By \eqref{eq10.4}, we have 
\begin{equation}\label{11.22}
\widehat{f}_{v}(0)=4^{-1}\mathcal{S}h_{v}(1/2+s_{2}). 	
\end{equation}  

In \eqref{eq10.4}, bi-$K'_{v}$-invariance gives four equal sign
contributions from $a,z$. In the positive quadrant put
$a=y$, $b=x/y$, and $z=z'/\sqrt{y}$. The multiplicative measures
are unchanged and $db=y^{-1}dx$. The matrix is now the argument
of $h_{v}$ in \eqref{eq10.17}, and the power is $y^{\nu-1/2}$.
Taking $\nu=1/2+s_{2}$ gives \eqref{11.22}.

\begin{lemma}\label{lemma11.4}
For $y>0$,
\begin{multline}\label{11.21}
\int_{0}^{\infty}h_{v}\left(n(x)\diag(y^{1/2}z,y^{-1/2}z)\right)d^{\times}z\\
=\frac{\sqrt{y}}{2\pi^{2}}\int_{\mathbb{R}}\int_{0}^{\infty}K_{it}(2\pi b)K_{it}(2\pi by)\cos(2\pi bx)db\mathcal{S}h_{v}(it)t\sinh(\pi t)dt.
\end{multline}
\end{lemma}

\begin{proof}
By the Cartan decomposition, the left-hand side of \eqref{11.21} is
$V(y+y^{-1}-2+y^{-1}x^{2})$, where
\begin{align*}
V(y+y^{-1}-2):=\int_{0}^{\infty}
h_{v}(\diag(y^{1/2}z,y^{-1/2}z))\,d^{\times} z.
\end{align*}

Using the positive-quadrant substitution from the proof of
\eqref{11.22} in \eqref{eq10.4} gives
\begin{align*}
 \tfrac{1}{4}\mathcal{S} h_{v}(\nu)
 =\int_{0}^{\infty}\int_{\mathbb{R}}
 V(y+y^{-1}-2+x^{2}/y)y^{\nu-1/2}\,dx\,d^{\times}y.
\end{align*}
At fixed $y$, the substitution $x=\sqrt{y}\,u$ gives the Selberg
transform of $V$. The inversion formula \cite[(2.24)]{Zag79} gives
\begin{equation}\label{10.36}
V(u)=\frac{1}{16\pi }\int_{\mathbb{R}} P_{-\frac{1}{2}+it}(1+u/2)\mathcal{S}h_{v}(it)t\tanh(\pi t)dt,
\end{equation}
where $P_{-1/2+\nu}(\cdot)$ is the Legendre function of the first kind.

By \cite[\textsection 1.12, (51)]{EMOT54},
\begin{multline}\label{eq10.40}
P_{-1/2+it}\left(\frac{y^{2}+1+x^{2}}{2y}\right)\\
=\frac{8\sqrt{y}\cosh(\pi t)}{\pi}\int_{0}^{\infty}K_{it}(2\pi b)K_{it}(2\pi by)\cos(2\pi bx)db.
\end{multline}

Multiplying \eqref{10.36} and \eqref{eq10.40} gives
\eqref{11.21}, since $\tanh(\pi t)\cosh(\pi t)=\sinh(\pi t)$.

For the integral exchanges below, the standard representation
\begin{align*}
K_{it}(x)=\frac{1}{2}\int_{\mathbb{R}}e^{-x\cosh u+itu}\,du
\end{align*}
implies
\begin{equation}\label{eq:real-k-bessel-majorant}
|K_{it}(x)|\ll e^{-\pi|t|/2}
K_{0}\left(x\sin\frac{1}{1+|t|}\right),\qquad x>0.
\end{equation}
For $t\neq0$, shift $u\mapsto u+i\sgn(t)(\pi/2-(1+|t|)^{-1})$.
The factor $e^{itu}$ then contributes
$e^{-\pi|t|/2+|t|/(1+|t|)}\leq e^{1-\pi|t|/2}$.
The case $t=0$ follows from the monotonicity of
$K_{0}$. Since $\sin((1+|t|)^{-1})\asymp(1+|t|)^{-1}$,
substituting $x=(1+|t|)y$ in the weighted $K_{0}$-integrals gives
bounds by powers of $1+|t|$. Two $K$-Bessel factors cancel the exponential growth of
$\sinh(\pi t)$, and $\mathcal{S}h_{v}(it)$ is rapidly decreasing.
\end{proof}

\begin{lemma}
Let $f_{v}$ and $\widehat{f}_{v}$ be defined by \eqref{eq10.17} and
\eqref{c11.20}. For $\Re(s_{2})>-1/2$ and $y\neq0$,
\begin{multline}\label{eq11.26}
\widehat{f}_{v}(y)=\frac{|y|^{-\frac{1}{2}-s_{2}}}{16\pi^{5/2+s_{2}}}\int_{\mathbb{R}}\, \mathcal{S}h_{v}(it)t\sinh(\pi t)K_{it}(2\pi |y|)\\
\Gamma\left(\frac{s_{2}+\frac{1}{2}+it}{2}\right)
\Gamma\left(\frac{s_{2}+\frac{1}{2}-it}{2}\right)dt.
\end{multline}
If also $\Re(\lambda)>1/2+\Re(s_{2})$, then
\begin{multline}\label{fc11.29}
\mathcal{M}\widehat{f}_{v}(\lambda):=\int_{0}^{\infty} \widehat{f}_{v}(y)y^{\lambda}d^{\times}y=\frac{\pi^{-\lambda}}{64\pi^{2}}\int_{\mathbb{R}}\, \mathcal{S}h_{v}(it)t\sinh(\pi t)\\
\Gamma\left(\frac{\lambda-\frac{1}{2}-s_{2}+it}{2}\right)
\Gamma\left(\frac{\lambda-\frac{1}{2}-s_{2}-it}{2}\right)
\Gamma\left(\frac{s_{2}+\frac{1}{2}+it}{2}\right)
\Gamma\left(\frac{s_{2}+\frac{1}{2}-it}{2}\right)dt.
\end{multline}
For $|\Re(s_{2})|<1/2$ and $0<\Re(\lambda)<1/2-\Re(s_{2})$,
\begin{multline}\label{eq10.37}
\mathcal{M}f_{v}(\lambda)=\frac{\pi^{-5/2}\Gamma(\lambda/2)}{64\Gamma((1-\lambda)/2)}\int_{\mathbb{R}}\, 
\Gamma(1/4+s_{2}/2+it/2)\Gamma(1/4+s_{2}/2-it/2)\\
\Gamma(1/4-s_{2}/2+it/2-\lambda/2)
\Gamma(1/4-s_{2}/2-it/2-\lambda/2)\mathcal{S}h_{v}(it)t\sinh(\pi t)dt.
\end{multline}
\end{lemma}

\begin{proof}
We use the standard Mellin transform \cite[\textsection 6.561, (16)]{GR14}:
\begin{equation}\label{10.42}
\int_{0}^{\infty} K_{it}(2\pi y)y^{\alpha}\,d^{\times} y
=\frac{\pi^{-\alpha}}{4}
\Gamma\left(\frac{\alpha+it}{2}\right)
\Gamma\left(\frac{\alpha-it}{2}\right),\qquad \Re(\alpha)>0.
\end{equation}
First take $|\Re(s_{2})|<1/2$. Substitute \eqref{11.21} in
\eqref{eq10.17} and apply \eqref{10.42} to the $y$-integral,
with $\alpha=1/2+s_{2}$. The remaining $b$-integral contains
$b^{-1/2-s_{2}}K_{it}(2\pi b)$ and converges at $b=0$ since
$\Re(s_{2})<1/2$. The bound \eqref{eq:real-k-bessel-majorant}
justifies changing the order of integration in $y,b,t$.
The $1/4$ in \eqref{10.42} and the $1/2$ from extending the
$b$-integral evenly to $\mathbb{R}$ give the scalar in \eqref{eq11.26}
after Fourier inversion.
For fixed $y\neq0$, both sides are holomorphic in $\Re(s_{2})>-1/2$,
which gives its extension to that half-plane.
Taking the Mellin transform of \eqref{eq11.26}, use \eqref{10.42}
at $\alpha=\lambda-1/2-s_{2}$ to get \eqref{fc11.29}.
The bound \eqref{eq:real-k-bessel-majorant}, Stirling's formula and
rapid decay of $\mathcal{S}h_{v}(it)$ justify interchanging the
$y$- and $t$-integrals for $\Re(s_{2})>-1/2$ and
$\Re(\lambda)>1/2+\Re(s_{2})$.

To compute $\mathcal{M} f_{v}(\lambda)$, use evenness and the cosine Mellin identity
\begin{align*}
\mathcal{M} f_{v}(\lambda)
=2(2\pi)^{-\lambda}\Gamma(\lambda)\cos(\pi\lambda/2)
\mathcal{M}\widehat{f}_{v}(1-\lambda),\qquad 0<\Re(\lambda)<1.
\end{align*}
The hypotheses give $\Re(1-\lambda)>1/2+\Re(s_{2})$.
Substituting \eqref{fc11.29} at $1-\lambda$ and using
$2^{1-\lambda}\Gamma(\lambda)\cos(\pi\lambda/2)
=\sqrt{\pi}\,\Gamma(\lambda/2)/\Gamma((1-\lambda)/2)$ gives \eqref{eq10.37}.
\end{proof}

\subsection{The Dual Weight}\label{sec11.4}
Take $\eta_{v}=\mathbf{1}$ in \eqref{fc10.2} and use the vector
\eqref{fc11.3}.

\subsubsection{Simplification of the dual weight}
Insert \eqref{fc11.3} into \eqref{fc10.2} and use Iwasawa decomposition,
the sphericity of $W_{v}^{\circ}$ and unitarity of $\sigma_{v}$ to write
\begin{multline}\label{10.3}
\mathcal{H}_{W_{v}}^{\vee}(\sigma_{v},\mathbf{s}^{\vee},\mathbf{1})=\sum_{W_{v}'\in \mathfrak{B}_{\sigma_{v}}}
\int_{F_{v}^{\times}}\overline{W_{v}'}\left(\diag(a_{2},1)\right)|a_{2}|^{\frac{3s_{1}+s_{2}}{2}}d^{\times}a_{2}\\
\int_{F_{v}}\int_{F_{v}^{\times}}\int_{F_{v}^{\times}}h_{v}\left(n(-b_{1})\diag(a_{1},z_{1})\right)|a_{1}|^{s_{2}}|z_{1}|^{-s_{2}}d^{\times}a_{1}d^{\times}z_{1}
\mathcal{I}_{W_{v}'}(b_{1})db_{1},
\end{multline}
where
\begin{align*}
\mathcal{I}_{W_{v}'}(b_{1})
:=\Psi_{v}((1+s_{2}-s_{1})/2,\pi_{v}(n_{13}(b_{1}))W_{v}^{\circ},W_{v}').
\end{align*}

Use $\beta_{v}$ from Lemma~\ref{lem:real-rs-evaluation} and
$\mathcal{S}_{\sigma_{v}}$ from Lemma~\ref{lemm11.3}. For the weight in
Theorem~\ref{thmC}, take $\beta_{v}=\mathbf{L}(\pi_{v})^{-1}$. Put
\begin{equation}\label{eq:real-double-Mellin-Gamma}
 G_{\boldsymbol{\nu}}(\lambda_{1},\lambda_{2})
 :=\frac{\prod_{j=1}^{3}\Gamma((\lambda_{1}-\nu_{j})/2)
 \Gamma((\lambda_{2}+\nu_{j})/2)}{\Gamma((\lambda_{1}+\lambda_{2})/2)}.
\end{equation}
With $r=\lambda_{1}+s_{1}-s_{2}$, define
\begin{multline*}
\mathcal{I}_{W_{v}}(\sigma_{v},\mathbf{s},t)
:=\frac{1}{(2\pi i)^{2}}\int_{(c_{1})}\int_{(c_{2})}
G_{\boldsymbol{\nu}}(\lambda_{1},\lambda_{2})
\frac{\Gamma(r/2)\Gamma(1/4+it/2-\lambda_{1}/2-s_{1}/2)}{\Gamma((1-r)/2)}\\
\Gamma(1/4-it/2-\lambda_{1}/2-s_{1}/2)
\mathcal{S}_{\sigma_{v}}(\lambda_{1},\lambda_{2},\mathbf{s})
\pi^{-\lambda_{2}}\,d\lambda_{2}\,d\lambda_{1}.
\end{multline*}
Initially choose $c_{1},c_{2}$ so that \eqref{eq:normalized-ktype-region}
holds on $\Re(\lambda_{j})=c_{j}$, $j=1,2$, and
\begin{equation}\label{eq:dual-common-contours}
\begin{cases}
 c_{1},c_{2}>\vartheta_{\pi_{v}},\quad -1/2<\Re(s_{2})<0,\\
 1/2<c_{1}+\Re(s_{1}-s_{2})<1/2-\Re(s_{2}).
\end{cases}
\end{equation}

\begin{thm}\label{prop10.2}
Under \eqref{eq:normalized-ktype-region} and
\eqref{eq:dual-common-contours},
\begin{multline}\label{10.25}
\mathcal{H}_{W_{v}}^{\vee}(\sigma_{v},\mathbf{s}^{\vee},\mathbf{1})=\frac{\beta_{v}\pi^{s_{1}-s_{2}}}{32\pi^{7/2}}\int_{\mathbb{R}}\, 
\Gamma(1/4+s_{2}/2+it/2)\Gamma(1/4+s_{2}/2-it/2)\\
\mathcal{I}_{W_{v}}(\sigma_{v},\mathbf{s},t)\cdot \mathcal{S}h_{v}(it)t\sinh(\pi t)dt.
\end{multline}
If $\vartheta_{\pi_{v}}<3/8$ and $\vartheta_{\sigma_{v}}<1/8$, the same
identity holds at $\mathbf{s}=\mathbf{0}$ with $c_{1}=15/32$ and $c_{2}=3/8$.
\end{thm}

\begin{proof}
For $W_{n}$ of \textsection\ref{sec10.2.1}, the rotation integral is
\begin{align*}
 \frac{1}{2\pi}\int_{0}^{2\pi}e^{2\pi izb_{1}\sin\theta+2in\theta}\,d\theta
 =J_{2n}(2\pi zb_{1}).
\end{align*}
The other component of $K_{v}'=\mathrm{O}(2)$ gives the same
$a$-integral after $a\mapsto-a$. Hence
\begin{multline}\label{10.4}
\mathcal{I}_{W_{n}}(b_{1})=\int_{F_{v}^{\times}}\int_{F_{v}^{\times}}
W_{v}^{\circ}(\diag(az,z,1))W_{n}(\diag(a,1))\\
J_{2n}(2\pi zb_{1})|a|^{(s_{2}-s_{1})/2-1}
|z|^{s_{2}-s_{1}}\,d^{\times} a\,d^{\times} z.
\end{multline}
In \eqref{10.3}, use $|a_{1}|=y^{1/2}z'$ and
$|z_{1}|=y^{-1/2}z'$, which preserve the product of the multiplicative
measures on each of the four sign ranges. The integral in $a_{1},z_{1}$
is then $4f_{v}(-b_{1})=4f_{v}(b_{1})$ by \eqref{eq10.17}.
Insert \eqref{10.4} and use evenness in $z$ to obtain
\begin{multline}\label{10.7}
\mathcal{H}_{W_{v}}^{\vee}(\sigma_{v},\mathbf{s}^{\vee},\mathbf{1})
=8\sum_{2n\in\mathcal{N}_{\sigma_{v}}}
\overline{\Psi(W_{n};(3\overline{s_{1}}+\overline{s_{2}})/2)}\\
\int_{F_{v}}f_{v}(b_{1})\int_{0}^{\infty}\int_{F_{v}^{\times}}
W_{v}^{\circ}(\diag(z|a|,z,1))W_{n}(\diag(a,1))\\
J_{2n}(2\pi zb_{1})|a|^{(s_{2}-s_{1})/2-1}z^{s_{2}-s_{1}}
\,d^{\times} a\,d^{\times} z\,db_{1}.
\end{multline}
For $a,z>0$, we recover the Whittaker function in \eqref{10.7}
by inverting \eqref{eq:real-jacquet-double-mellin}:
\begin{multline}\label{10.19}
W_{v}^{\circ}(\diag(az,z,1))
=\frac{\beta_{v}}{4(2\pi i)^{2}}\int_{(c_{1})}\int_{(c_{2})}
G_{\boldsymbol{\nu}}(\lambda_{1},\lambda_{2})\\
\pi^{-\lambda_{1}-\lambda_{2}}
 z^{1-\lambda_{1}}a^{1-\lambda_{2}}\,d\lambda_{2}\,d\lambda_{1},
\end{multline}
where $c_{1},c_{2}>\max_{j}|\Re(\nu_{j})|$.
By \cite[\textsection 6.561, (14)]{GR14} and evenness of $f_{v}$,
\begin{equation}\label{10.21}
\int_{\mathbb{R}}\int_{0}^{\infty}
f_{v}(b)J_{2n}(2\pi zb)z^{1-\lambda}\,d^{\times} z\,db
=\frac{\Gamma(|n|+1/2-\lambda/2)\mathcal{M} f_{v}(\lambda)}
{\pi^{1-\lambda}\Gamma(|n|+1/2+\lambda/2)}.
\end{equation}
The integrals converge absolutely for $1/2<\Re(\lambda)<1$.

Insert \eqref{10.19} into \eqref{10.7}. The $b_{1},z$-integral is
\eqref{10.21} at $\lambda=\lambda_{1}+s_{1}-s_{2}$, and the $a$-integral
is $\Psi(W_{n};(s_{2}-s_{1})/2-\lambda_{2})$.
For $\lambda_{j}=c_{j}+i\tau_{j}$, \eqref{eq:hecke-vertical} and Stirling's
formula bound the $n$-th term, including $G_{\boldsymbol{\nu}}$, by
\begin{align*}
 (1+|\tau_{1}|+|\tau_{2}|)^{l}
 e^{-\pi(|\tau_{1}|+|\tau_{2}|)/2}
 (1+|n|)^{-1-\delta}\log^{2}(2+|n|).
\end{align*}
Since $\delta>0$ and $\mathcal{M}f_{v}(r)$ is bounded on these lines,
the sum and the two integrals may be interchanged.
Insert \eqref{eq10.37} and exchange the $t$-integral with the
sum and Mellin integrals to obtain \eqref{10.25}. Stirling's formula
and rapid decay of $\mathcal{S} h_{v}(it)$ justify the exchange.

For the central value, take
\begin{align*}
s_{1}=s_{2}=-\eta,\qquad c_{1}=15/32+\eta,\qquad c_{2}=3/8,
\qquad 0\leq\eta\leq1/16.
\end{align*}
At $\eta=1/16$ the initial conditions hold. Continue the right-hand
side of \eqref{10.25} along this path, keeping
$0<\Re(r)<1/2-\Re(s_{2})$ as in \eqref{eq10.37}.
The lower bound $\Re(r)>1/2$ was needed only for the absolute
integral in \eqref{10.21}. Along the path,
\eqref{eq:normalized-ktype-region} remains strict, with
$\delta=3/32-\eta\geq1/32$. The bounds for the $n$-sum and the
$t,\lambda_{1},\lambda_{2}$ integrals are locally uniform, which gives
\eqref{10.25} at $\eta=0$.
\end{proof}

\subsection{Bounds for the Dual Weight}
\label{sec:G-normalized-real-transform}

Take $\pi_{v}=\mathbf{1}\boxplus\mathbf{1}\boxplus\mathbf{1}$ and
$\boldsymbol{\nu}=\mathbf{s}=\mathbf{0}$, so
$\beta_{v}=\mathbf{L}(\pi_{v})^{-1}=1$ in Theorem~\ref{thmC}.
\phantomsection\label{def:real-spherical-parameter}
Write $\sigma_{t}=|\cdot|^{it}\boxplus|\cdot|^{-it}$.

For $h_{v},W_{v}$ as in \eqref{fc11.3} at $s_{1}=0$,
Lemma~\ref{lem:real-rs-evaluation} gives
\begin{equation}\label{eq:G-as-L-ratio}
 \pi\mathcal{H}_{W_{v}}(\sigma_{t},\mathbf{0},\mathbf{1})
 =\frac{L_{v}(1/2,\sigma_{t})^{4}}
 {L_{v}(1,\sigma_{t},\Ad)}\mathcal{S}h_{v}(it).
\end{equation}

Fix a sufficiently large integer $m$. For $R\geq2$, let $H_{R}$ be even
and holomorphic in $|\Im(t)|\leq m+1$. Suppose it vanishes to order
at least four at $t=\pm i(j+1/2)$ for $0\leq j\leq m$, is nonnegative for
$t\in\mathbb{R}$ and for $t=i\eta$, $|\eta|<1/2$, and satisfies,
uniformly in that strip and for every $j\geq0$,
\begin{equation}\label{eq:R-admissible-test}
 H_{R}^{(j)}(t)\ll_{m,j}
 e^{-(\Re(t)/R)^{2}}
 \left(1+\frac{|t|^{2}}{R^{2}}\right)^{j}
 \left(\frac{1+|t|}{R}\right)^{m}(1+|t|)^{-j}.
\end{equation}
For example, take
\begin{align*}
 H_{R}(t)=e^{-2t^{2}/R^{2}}\prod_{j=0}^{m}
 \left(\frac{t^{2}+(j+1/2)^{2}}{R^{2}}\right)^{4}.
\end{align*}
Then $H_{R}(t)\geq e^{-8}$ for $R\leq|t|\leq2R$ on the real axis.
The orders of the zeros may be increased if needed.

\phantomsection\label{def:polynomial-spherical-test}
We choose $h_{v,R}$ with
\begin{align*}
 \mathcal{S}h_{v,R}(it)
 =H_{R}(t)\frac{L_{v}(1,\sigma_{t},\Ad)}{L_{v}(1/2,\sigma_{t})^{4}}.
\end{align*}
Let $V$ be given by \eqref{10.36} with this prescribed transform.
\phantomsection\label{def:central-test-cutoff}
Fix an even $\alpha\in C_{c}^{\infty}(\mathbb{R}^{\times})$ with
$\int_{\mathbb{R}^{\times}}\alpha(z)\,d^{\times}z=1$, and set
\begin{align*}
 h_{v,R}(z\diag(y^{1/2},y^{-1/2}))
 :=2\alpha(z)V(y+y^{-1}-2),\qquad y>0,\quad z\in\mathbb{R}^{\times}.
\end{align*}
Extend by bi-$K_{v}'$-invariance. Its integral over $z>0$ is
$V(y+y^{-1}-2)$, so it has the prescribed transform.
The estimates in Appendix~\ref{app:central-polynomial-transform}
justify using $h_{v}=h_{v,R}$ in \eqref{fc11.3} and the
dual-weight formula \eqref{10.7}.

\begin{defn}\label{def:real-dual-transform}
Take $h_{v}=h_{v,R}$ in \eqref{fc11.3} at $\mathbf{s}=\mathbf{0}$.
For generic unitary $\sigma$ with trivial central character, set
\begin{align*}
 (\mathscr{T}H_{R})(\sigma):=
 \pi\mathcal{H}_{W_{v}}^{\vee}(\sigma,\mathbf{0},\mathbf{1}).
\end{align*}
For $\sigma=\sgn^{\epsilon}|\cdot|^{iu}\boxplus\sgn^{\epsilon}|\cdot|^{-iu}$,
with $u\in\mathbb{R}$ and $\epsilon\in\{0,1\}$, write this value as
$(\mathscr{T}_{\epsilon}H_{R})(u)$, omitting the subscript when $\epsilon=0$.
For a complementary series
$\sigma=\sgn^{\epsilon}|\cdot|^{\eta}\boxplus\sgn^{\epsilon}|\cdot|^{-\eta}$,
with $\eta\in\mathbb{R}$ and $0<|\eta|<1/2$, use continuation to $u=i\eta$.
For a discrete series $\sigma$ of lowest weight $2l$,
write $(\mathscr{T}H_{R})^{\mathrm{disc}}(2l)$. The Fourier factors are
\begin{align*}
 C_{0}(z):=2(2\pi)^{-z}\Gamma(z)\cos(\pi z/2),\quad
 C_{1}(z):=-2i(2\pi)^{-z}\Gamma(z)\sin(\pi z/2).
\end{align*}
\end{defn}

Let $f_{v,R}$ be $f_{v}$ in \eqref{eq10.17} at $h_{v}=h_{v,R}$ and $s_{2}=0$.
The key step is to sum over all $K_{v}'$-types before estimating the
principal-series dual weight. The Jacobi--Anger identity in the proof
of Lemma~\ref{lem:real-principal-dual-integral} replaces this sum by
the cosine kernel \eqref{eq:principal-cosine-kernel}.
The subsequent integral evaluations, Gamma cancellations and uniform
hypergeometric bounds give the polynomial estimates of
Proposition~\ref{prop:direct-polynomial-central-transform}, used in
the fifth-moment bound in \textsection\ref{sec15.2}.

\begin{lemma}\label{lem:real-principal-dual-integral}
For $u\in\mathbb{R}$ and $\epsilon\in\{0,1\}$,
\begin{multline}\label{eq:direct-central-principal-transform}
 (\mathscr{T}_{\epsilon}H_{R})(u)
 =8\int_{\mathbb{R}}f_{v,R}(b)
 \int_{0}^{\infty}\int_{\mathbb{R}^{\times}}
 W_{v}^{\circ}(\diag(z|a|,z,1))\\
 \mathcal{K}_{iu,\epsilon}(a,2\pi zb)
 |a|^{-1}\,d^{\times} a\,d^{\times} z\,db,
\end{multline}
where, for $a\neq0$ and $b\in\mathbb{R}$,
\begin{multline}\label{eq:principal-cosine-kernel}
 \mathcal{K}_{iu,\epsilon}(a,b)
 :=\overline{C_{\epsilon}(1/2-iu)}\,|a|^{1/2-iu}(\sgn a)^{\epsilon}
 \int_{\mathbb{R}}\frac{e^{-2\pi iax}}{(1+x^{2})^{1/2+iu}}\\
 \int_{-\pi/2}^{\pi/2}
 (\sgn\theta)^{\epsilon}|\sin\theta\cos\theta|^{-1/2-iu}
 \cos(b\sin(\arctan x-\theta))\,d\theta\,dx.
\end{multline}
\end{lemma}

\begin{proof}
By \eqref{10.7} at $\mathbf{s}=\mathbf{0}$ and
Definition~\ref{def:real-dual-transform},
\eqref{eq:direct-central-principal-transform} holds with
\begin{align*}
 \mathcal{K}_{iu,\epsilon}(a,b)=
 \pi\sum_{n\in\mathbb{Z}}\overline{\Psi(W_{n};0)}
 W_{n}(\diag(a,1))J_{2n}(b).
\end{align*}
Here $W_{n}$ is given by \eqref{10.8} with $\nu=iu$ and parity
$\epsilon$, and $J_{2n}$ is the Bessel function of the first kind.

Multiplying $W_{n}$ by a scalar of absolute value one leaves
$\overline{\Psi(W_{n};0)}W_{n}$ unchanged. We therefore use the Fourier
integral in \cite[(2.16), \textsection4]{BM05}. Fourier inversion
followed by $x=\tan\theta$ expresses $\Psi(W_{n};0)$ as an integral
over $(-\pi/2,\pi/2)$. Substitute these expressions and use
the even part of the Jacobi--Anger expansion
\begin{align*}
 \sum_{n\in\mathbb{Z}}J_{2n}(b)e^{2in\theta}
 =\cos(b\sin\theta)
\end{align*}
to obtain \eqref{eq:principal-cosine-kernel}.
For fixed $a,b$, the rapid decay of $J_{2n}(b)$ permits termwise
integration on bounded $x$-intervals; integration by parts then
justifies passing to the integral over $\mathbb{R}$.
\end{proof}

Define
\begin{equation}\label{eq:central-input-Ft}
\begin{aligned}
 A(t)&:=\Gamma(1/4+it/2)\Gamma(1/4-it/2),\\
 F_{t}(b)&:={}_{2}F_{1}(1/4+it/2,1/4-it/2;1/2;-b^{2}).
\end{aligned}
\end{equation}

For $D>0$, define
\begin{equation}\label{eq:central-K-Fourier-Mellin-kernel}
 B_{u,\epsilon}(D,x):=\int_{\mathbb{R}}|a|^{-1/2-iu}(\sgn a)^{\epsilon}
 K_{0}(2\pi D|a|)e^{-2\pi iax}\,da,
\end{equation}
which converges absolutely. Define
\begin{multline}\label{eq:central-omega-definition}
 \Omega_{u,\epsilon}(t):=\int_{0}^{\infty}\int_{\mathbb{R}}
 \int_{-\pi/2}^{\pi/2}
 \frac{B_{u,\epsilon}(\sqrt{1+s^{-1}},x)(\sgn\theta)^{\epsilon}|\sin\theta\cos\theta|^{-1/2-iu}}
 {(1+x^{2})^{1/2+iu}}\\
 |\sin(\arctan x-\theta)|^{-1}F_{t}\left(\sqrt{\frac{1+s}{\sin^{2}(\arctan x-\theta)}-1}\right)
 \,d\theta\,dx\,\frac{ds}{s}.
\end{multline}

\begin{prop}\label{prop:central-integral-transform}
For $u\in\mathbb{R}$ and $\epsilon\in\{0,1\}$,
\begin{equation}\label{eq:direct-threefold-central-transform}
 (\mathscr{T}_{\epsilon}H_{R})(u)
 =\frac{\overline{C_{\epsilon}(1/2-iu)}}{4\pi^{3/2}}
 \int_{\mathbb{R}}H_{R}(t)\frac{t\tanh(\pi t)}{A(t)}
 \Omega_{u,\epsilon}(t)\,dt.
\end{equation}
The integral in \eqref{eq:central-omega-definition} converges absolutely
for real $t$. For $|t|\geq1$, its integral of absolute values is
\begin{align*}
 \ll (1+|u|)^{-1/2}\log^{2}(2+|u|).
\end{align*}
\end{prop}

\begin{proof}
For $h_{v}=h_{v,R}$ and $s_{2}=0$, we recover $f_{v,R}(b)$ from
$\widehat{f}_{v}$ in \eqref{eq11.26} by \cite[\textsection6.699, (4)]{GR14}:
\begin{equation}\label{eq:direct-central-input-hypergeometric}
 f_{v,R}(b)=\frac{1}{32\pi}\int_{\mathbb{R}}
 H_{R}(t)\frac{t\tanh(\pi t)}{A(t)^{2}}F_{t}(b)\,dt.
\end{equation}
For $a,z>0$, Mellin inversion of
\eqref{eq:real-jacquet-double-mellin} takes the form
\begin{equation}\label{eq:central-Whittaker-two-K-kernel}
 W_{v}^{\circ}(\diag(az,z,1))
 =4az\int_{0}^{\infty}
 K_{0}\left(2\pi z\sqrt{1+s}\right)
 K_{0}\left(2\pi a\sqrt{1+s^{-1}}\right)\frac{ds}{s}.
\end{equation}
We also have
\begin{multline*}
\int_{0}^{\infty} K_{0}\left(2\pi z\sqrt{1+s}\right)
 \cos(2\pi zb\sin(\arctan x-\theta))\,dz\\
 =\frac{1}{4\sqrt{1+s+b^{2}\sin^{2}(\arctan x-\theta)}}.
\end{multline*}
For real $t$ and $c>0$, we have
\begin{equation}\label{eq:central-hypergeometric-kernel-pairing}
 \int_{\mathbb{R}}\frac{F_{t}(b)}{\sqrt{c^{2}+b^{2}}}\,db
 =\pi^{-1/2}A(t)\cdot 
 {}_{2}F_{1}(1/4+it/2,
1/4-it/2;1/2;1-c^{2}).
\end{equation}
For $|b|\geq2$, \cite[\textsection2.1.4, (17)--(18)]{Bat53} gives
$\partial_{b}^{j}F_{t}(b)=O_{t,j}(|b|^{-1/2-j}(1+\log|b|))$
for each fixed $j\geq0$. Thus the integral converges absolutely and is
analytic for $c>0$. Integration by parts shows that both sides of
\eqref{eq:central-hypergeometric-kernel-pairing} satisfy the same
differential equation in $c$. At $c=1$, they are analytic and have
the same value $A(t)/\sqrt{\pi}$ by Pfaff's transformation and Gauss's
summation formula \cite[\textsection2.1.4, (22); \textsection2.1.3, (14)]{Bat53}.
Uniqueness of the analytic solution proves the identity.

To integrate against $H_{R}(t)$, put $g_{t}(r):=(\cosh r)^{1/2}F_{t}(\sinh r)$.
The hypergeometric
differential equation becomes
\begin{equation}\label{eq:direct-central-input-wave-equation}
 g_{t}''(r)+\bigl(t^{2}-\tfrac{1}{4}\operatorname{sech}^{2}r\bigr)g_{t}(r)=0,
 \qquad g_{t}(0)=1,\quad g_{t}'(0)=0.
\end{equation}
Its integral equation is
\begin{align*}
 g_{t}(r)=\cos(tr)+\frac{1}{4t}\int_{0}^{r}
 \sin(t(r-s))\operatorname{sech}^{2}s\,g_{t}(s)\,ds,\qquad t\neq0.
\end{align*}
Using $\int_{0}^{\infty}\operatorname{sech}^{2}s\,ds=1$, we get,
for real $r,t$ with $|t|\geq1$ and $j\geq0$,
\begin{equation}\label{eq:direct-central-input-wave-bound}
 |g_{t}(r)-\cos(tr)|\ll |t|^{-1},\qquad
 |g_{t}^{(j)}(r)|\ll_{j}(1+|t|)^{j}.
\end{equation}
For the derivative bounds, differentiate once and then use
\eqref{eq:direct-central-input-wave-equation}.
For $|t|\geq1$, this gives $F_{t}(b)\ll(1+b^{2})^{-1/4}$; for
$|t|\leq1$, an additional factor $\log(2+|b|)$ suffices.
Together with the Gaussian decay of $H_{R}$ for fixed $R$, these bounds
justify interchanging the $b$- and $t$-integrals in
\eqref{eq:direct-central-input-hypergeometric}. Thus, for $c>0$ and
$0<d\leq c$, applying \eqref{eq:central-hypergeometric-kernel-pairing}
with $c/d$ in place of $c$ gives
\begin{equation}\label{eq:central-input-kernel-pairing}
 \int_{\mathbb{R}}\frac{f_{v,R}(b)}{\sqrt{c^{2}+d^{2}b^{2}}}\,db
 =\frac{1}{32\pi^{3/2}d}\int_{\mathbb{R}}
 H_{R}(t)\frac{t\tanh(\pi t)}{A(t)}
 F_{t}\left(\sqrt{c^{2}/d^{2}-1}\right)\,dt.
\end{equation}

Together with Lemma~\ref{lem:real-principal-dual-integral},
\eqref{eq:central-input-kernel-pairing} and
\eqref{eq:central-K-Fourier-Mellin-kernel}, these identities yield
\eqref{eq:direct-threefold-central-transform}.
For real $u$, the Gamma reflection formula gives
\begin{equation}\label{eq:principal-unitary-Fourier-factor}
 |C_{\epsilon}(1/2-iu)|=1,\qquad
 \overline{C_{\epsilon}(1/2-iu)}
 =(-1)^{\epsilon}C_{\epsilon}(1/2+iu).
\end{equation}

To justify exchanging the integrals, insert $K_{0}(y)=\int_{0}^{\infty}e^{-y\cosh v}\,dv$ into
\eqref{eq:central-K-Fourier-Mellin-kernel} and first integrate in $a$.
For real $x$ and $D\geq1$, split the $v$-integral at
$\log(2+|x|/D)$. Using
$|\Gamma(1/2-iu)|^{2}=\pi/\cosh(\pi u)$, we obtain
\begin{equation}\label{eq:central-K-Fourier-Mellin-bound}
 |B_{u,\epsilon}(D,x)|\ll \frac{\log(2+|x|/D)}{D^{1/2}
 (1+|x|/D)^{1/2}}\exp\big(-|u|(\pi/2-\arctan(|x|/D))\big).
\end{equation}
For $x\neq0$,
\begin{align*}
\int_{-\pi/2}^{\pi/2}
 |\sin\theta\cos\theta\sin(\arctan x-\theta)|^{-1/2}\,d\theta\ll 1+\log(2+|x|)+\log(1+|x|^{-1}).
\end{align*}
Also, by \eqref{eq:direct-central-input-wave-bound},
\begin{align*}
\left|F_{t}\left(\sqrt{\frac{1+s}{\sin^{2}(\arctan x-\theta)}-1}\right)\right|\ll (1+s)^{-1/4}|\sin(\arctan x-\theta)|^{1/2},\quad |t|\geq1.
\end{align*}
Set $x=(1+s^{-1})^{1/2}y$. The exponential in
\eqref{eq:central-K-Fourier-Mellin-bound} is
$\leq e^{-\pi|u|/4}$ for $|y|\leq1$ and
$\leq e^{-c|u|/|y|}$ for $|y|\geq1$, with absolute $c>0$.
The integral of the absolute value in $x,\theta$ is therefore
\begin{align*}
 \ll (1+s)^{-1/4}(1+s^{-1})^{-1/4}(1+\tfrac{1}{2}\log(1+s^{-1}))^{2}
 (1+|u|)^{-1/2}\log^{2}(2+|u|).
\end{align*}
Its $s$-integral converges, since
$(1+s)^{-1/4}(1+s^{-1})^{-1/4}=s^{1/4}(1+s)^{-1/2}$.
For $|t|\leq1$, the integral equation gives $|g_{t}(r)|\ll1+|r|$;
this adds only integrable logarithmic factors. Convergence is therefore locally uniform in $t$,
and Gaussian decay permits the final $t$-integration for fixed $R$.
First carry out the calculation on bounded $x$-intervals.
Absolute convergence then justifies integration over all $x\in\mathbb{R}$,
proving \eqref{eq:direct-threefold-central-transform} for $f_{v}=f_{v,R}$.
\end{proof}

\begin{prop}\label{prop:direct-polynomial-central-transform}
Let $H_{R}$ be even and holomorphic in $|\Im(t)|\leq m+1$,
with zeros of order at least four at $t=\pm i(j+1/2)$ for
$0\leq j\leq m$, and assume \eqref{eq:R-admissible-test}.
For $\epsilon\in\{0,1\}$ and $u\in\mathbb{R}$,
\begin{equation}\label{eq:G-transform-tempered-bound}
 |(\mathscr{T}_{\epsilon}H_{R})(u)|\ll_{m} R^{2}(1+|u|/R)^{-3}.
\end{equation}
For $l\geq1$,
\begin{equation}\label{eq:G-transform-discrete-bound}
 |(\mathscr{T} H_{R})^{\mathrm{disc}}(2l)|
 \ll_{m}\begin{cases}1,&l=1,\\
 \min\{1,R^{4}[l(l-1)]^{-2}\},&l\geq2.
 \end{cases}
\end{equation}
Fix $\vartheta_{\infty}<1/2$. For $|\eta|\leq\vartheta_{\infty}$,
\begin{equation}\label{eq:G-transform-complementary-bound}
 |(\mathscr{T}_{\epsilon}H_{R})(i\eta)|\ll_{m,\vartheta_{\infty}}R^{2}.
\end{equation}
\end{prop}

\begin{proof}
In \eqref{eq:direct-central-cosine-period}, the $L^{1}$ bound of
Lemma~\ref{lem:principal-series-periods} gives $O_{m}(R^{2})$ for real $u$. Three integrations by parts give
$O_{m}(R^{5}|u|^{-3})$ when $u\neq0$, proving
\eqref{eq:G-transform-tempered-bound}. For $u=i\eta$,
$|\cos(i\eta\tau)|\leq e^{\vartheta_{\infty}\tau}$ gives
\eqref{eq:G-transform-complementary-bound}.

Lemma~\ref{lem:discrete-series-coefficient} gives \eqref{eq:G-transform-discrete-bound}.
\end{proof}

\section{Type \RNum{2}: Non-Archimedean Weights and Dual Transforms}\label{sec12}
We compute the finite-place spectral weights and dual transforms
for the arithmetic applications.

Let $\mathfrak{q},\mathfrak{n}\subseteq\mathcal{O}_{F}$ satisfy
$\mathfrak{n}+\mathfrak{q}=\mathcal{O}_{F}$. Throughout this section,
$v<\infty$ and $\pi_{v}$ is an irreducible unramified generic
representation of $G(F_{v})$ with trivial central character.
Normalize the spherical Whittaker vector by
$W_{v}^{\circ}(\diag(\varpi_{v}^{-2d_{v}},\varpi_{v}^{-d_{v}},1))=1$.
The change of additive character in
\textsection\ref{sec:finite-character-change} extends the formulas
at $d_{v}=0$ to general $d_{v}$.

\subsection{\texorpdfstring{Construction of $W_{v}$}{Construction of W\_v}}\label{sec11.2.1}
Assume $\omega_{v}'=\eta_{v}=\mathbf{1}$, and put
$r_{v}=e_{v}(\mathfrak{q})$ and $n_{v}=e_{v}(\mathfrak{n})$.

If $v\nmid\mathfrak{q}$, put
\begin{equation}\label{fc11.41}
W_{v}:=\sum_{\ell=0}^{n_{v}}q_{v}^{-n_{v}/2+(2\ell-n_{v})s_{2}}\sum_{\alpha\in \mathcal{O}_{v}/\mathfrak{p}_{v}^{\ell}}\pi_{v}(w_{2}n_{13}(\alpha\varpi_{v}^{-\ell}))W_{v}^{\circ}.
\end{equation}
In particular, $W_{v}=W_{v}^{\circ}$ when $v\nmid\mathfrak{q}\mathfrak{n}$.
For $m\geq0$, define
\begin{equation}\label{c12.2}
W_{v}^{[m]}:=\sum_{\alpha\in \mathcal{O}_{v}/\mathfrak{p}_{v}^{m}}
\pi_{v}(n_{13}(\alpha\varpi_{v}^{-m}))W_{v}^{\circ}.
\end{equation}
Then $W_{v}^{[0]}=W_{v}^{\circ}$. If $v\mid\mathfrak{q}$, take $W_{v}=W_{v}^{[r_{v}]}$.

\subsection{Hecke Eigenvalues and Old Forms}\label{sec:finite-newvector}
Let $\sigma_{v}$ be an irreducible generic unitary representation of
$G'(F_{v})$ with central character $\omega_{v}'=\mathbf{1}$, and let
$W_{v}'^{\circ}$ be its unit newvector for the local inner product of
\textsection\ref{sec3.2}.

\subsubsection{Hecke eigenvalues}
For $\ell\in\mathbb{Z}$, define
\begin{align*}
\lambda_{\sigma_{v}}(\mathfrak{p}_{v}^{\ell})
:=q_{v}^{\ell/2}
\frac{W_{v}'^{\circ}(a(\varpi_{v}^{\ell-d_{v}}))}
     {W_{v}'^{\circ}(a(\varpi_{v}^{-d_{v}}))}.
\end{align*}
In all cases $\lambda_{\sigma_{v}}(\mathfrak{p}_{v}^{0})=1$ and
$\lambda_{\sigma_{v}}(\mathfrak{p}_{v}^{\ell})=0$ for $\ell<0$.
Globally, for a generic automorphic representation $\sigma=\otimes_{v}'\sigma_{v}$, we define the Hecke eigenvalue by 
\begin{equation}\label{eq12.3}
\lambda_{\sigma}(\mathfrak{n}):=\prod_{v\mid\mathfrak{n}}\lambda_{\sigma_{v}}(\mathfrak{p}_{v}^{e_{v}(\mathfrak{n})}).	
\end{equation}

\subsubsection{Normalization factors}
For $d_{v}=0$, the unit newvector satisfies
\begin{equation}\label{eq:finite-newvector-value}
|W_{v}'^{\circ}(I_{2})|^{2}=
\begin{cases}
1,&r_{\sigma_{v}}=0,\\
(1+q_{v}^{-1})L_{v}(1,\sigma_{v},\Ad)/\zeta_{v}(2),&r_{\sigma_{v}}=1,\\
(1+q_{v}^{-1})L_{v}(1,\sigma_{v},\Ad),&r_{\sigma_{v}}\geq2.
\end{cases}
\end{equation}
The first case follows from \eqref{eq3.6} at $s=0$ and the inner
product in \textsection\ref{sec3.2}. For the other two, use
$|\lambda_{\sigma_{v}}(\mathfrak{p}_{v}^{n})|^{2}=q_{v}^{-n}$ when $r_{\sigma_{v}}=1$, and
$\lambda_{\sigma_{v}}(\mathfrak{p}_{v}^{n})=0$ for $n>0$ when $r_{\sigma_{v}}\geq2$.

For $d_{v}>0$, use the change of additive character in
\textsection\ref{sec:finite-character-change}.

\subsubsection{An orthonormal oldform basis}\label{sec12-oldform-basis}
We use the local form of the oldform orthogonalization
in \cite[Lemma~2]{BM15}. For $l\geq r_{\sigma_{v}}$, an orthonormal basis of
$\sigma_{v}^{K_{0,v}'[l]}$ is
$\mathfrak{B}_{\sigma_{v}}^{K_{0,v}'[l]}=\big\{W_{m,v}':\ 0\leq m\leq l-r_{\sigma_{v}}\big\}$, where 
\begin{equation}\label{c12.4}
W_{m,v}':=\sum_{i=0}^{m}\xi_{\sigma_{v}}(\mathfrak{p}_{v}^{i},\mathfrak{p}_{v}^{m})q_{v}^{(i-m)/2}\sigma_{v}\left(\diag(1,\varpi_{v}^{i})\right)W_{v}'^{\circ}.
\end{equation}

For $0\leq i\leq m$, the coefficients are
\begin{flalign}\label{eq:oldform-basis-coefficients}
 && \xi_{\sigma_{v}}(\mathfrak{p}_{v}^{i},\mathfrak{p}_{v}^{m})=
\begin{cases}
1,&i=m=0,\\
-q_{v}^{1/2}(1-\alpha_{\sigma_{v}}^{2})^{-1/2}\alpha_{\sigma_{v}},&i=0,\ m=1,\\
(1-\alpha_{\sigma_{v}}^{2})^{-1/2},&i=m=1,\\
[(1-\delta_{\sigma_{v}}q_{v}^{-2})(1-\alpha_{\sigma_{v}}^{2})]^{-1/2}\delta_{\sigma_{v}},&i=m-2,\ m\geq2,\\
-[(1-\delta_{\sigma_{v}}q_{v}^{-2})(1-\alpha_{\sigma_{v}}^{2})]^{-1/2}\lambda_{\sigma_{v}}(\mathfrak{p}_{v}),&i=m-1,\ m\geq2,\\
[(1-\delta_{\sigma_{v}}q_{v}^{-2})(1-\alpha_{\sigma_{v}}^{2})]^{-1/2},&i=m,\ m\geq2,\\
0,&0\leq i\leq m-3.
\end{cases}&
\end{flalign}
Here $\delta_{\sigma_{v}}:=\mathbf{1}_{r_{\sigma_{v}}=0}$ and
$\alpha_{\sigma_{v}}:=q_{v}^{-1/2}(1+\delta_{\sigma_{v}}q_{v}^{-1})^{-1}
\lambda_{\sigma_{v}}(\mathfrak{p}_{v})$.

\subsubsection{Hecke eigenvalues of $\pi_{v}$}
Write $L_{v}(s,\pi_{v})=\prod_{j=1}^{3}(1-\alpha_{j}q_{v}^{-s})^{-1}$. For $i,j\geq0$, the Casselman--Shalika formula is
\begin{equation}\label{12.5}
W_{v}^{\circ}\left(\diag(\varpi_{v}^{i+j-2d_{v}},\varpi_{v}^{j-d_{v}},1)\right)=q_{v}^{-i-j}\lambda_{\pi_{v}}(\mathfrak{p}_{v}^{i+j},\mathfrak{p}_{v}^{j}),
\end{equation}
where $\lambda_{\pi_{v}}(\mathfrak{p}_{v}^{i+j},\mathfrak{p}_{v}^{j})=s_{i+j,j,0}(\alpha_{1},\alpha_{2},\alpha_{3})$ is the Schur polynomial defined by 
\begin{equation}\label{12.6}
\frac{1}{
(\alpha_{1}-\alpha_{2})(\alpha_{1}-\alpha_{3})(\alpha_{2}-\alpha_{3})
}\det 
\begin{pmatrix}
\alpha_{1}^{i+j+2} & \alpha_{2}^{i+j+2} & \alpha_{3}^{i+j+2}\\
\alpha_{1}^{j+1}   & \alpha_{2}^{j+1}   & \alpha_{3}^{j+1}\\
1                  & 1                  & 1
\end{pmatrix}.
\end{equation}
At coincident Satake parameters, the quotient is evaluated by continuity.

\subsection{The Spectral Weight}
We evaluate \eqref{fc10.1} for the vectors of
\textsection\ref{sec11.2.1}. In this and the next subsection, assume
$d_{v}=0$; \eqref{eq:different-raw-weight} gives the formulas
for general $d_{v}$. The spectral calculations initially require
\begin{align*}
 \Re(s_{1})>-1/2+\vartheta_{\pi_{v}}+\vartheta_{\sigma_{v}},\quad \Re(s_{2})>-1/2+\vartheta_{\sigma_{v}}.
\end{align*}
All identities extend meromorphically in the parameters.

\subsubsection{Unramified calculation}
\begin{lemma}
If $v\nmid\mathfrak{q}$, then
\begin{multline}\label{fc11.42}
\mathcal{H}_{W_{v}}(\sigma_{v},\mathbf{s},\mathbf{1})=\lambda_{\sigma_{v}}(\mathfrak{p}_{v}^{n_{v}})|W_{v}'^{\circ}(I_{2})|^{2}\mathbf{1}_{r_{\sigma_{v}}=0}\\
L_{v}(1/2+s_{1},\pi_{v}\times\sigma_{v})L_{v}(1/2+s_{2},\sigma_{v}).
\end{multline} 
\end{lemma}

\begin{proof}
For $\ell\geq0$, put
\begin{align*}
W_{v}^{(\ell)}(g_{v}):=q_{v}^{-\ell/2}\sum_{\alpha\in \mathcal{O}_{v}/\mathfrak{p}_{v}^{\ell}}W_{v}^{\circ}(g_{v}n_{12}(\alpha\varpi_{v}^{-\ell})).
\end{align*}

Make the change $x_{v}\mapsto x_{v} n(-\alpha\varpi_{v}^{-\ell})$ in the Rankin--Selberg period.
Since $W_{v}^{\circ}$ is spherical, only $K_{v}'$-fixed vectors
contribute to the $\mathrm{GL}_{3}\times\mathrm{GL}_{2}$ integral. It vanishes if $r_{\sigma_{v}}>0$;
otherwise its value at $W_{v}'^{\circ}$ is
$W_{v}'^{\circ}(I_{2})L_{v}(1/2+s_{1},\pi_{v}\times\sigma_{v})$.
In the $\mathrm{GL}_{2}\times\mathrm{GL}_{1}$ integral, additive orthogonality gives
$q_{v}^{\ell}\mathbf{1}_{e_{v}(a_{v})\geq\ell}$.
The Hecke recurrence then gives
\begin{multline}\label{11.46}
\mathcal{H}_{W_{v}^{(\ell)}}(\sigma_{v},\mathbf{s},\mathbf{1})=
q_{v}^{-s_{2}\ell}\cdot \big[\lambda_{\sigma_{v}}(\mathfrak{p}_{v}^{\ell})-q_{v}^{-1/2-s_{2}}\lambda_{\sigma_{v}}(\mathfrak{p}_{v}^{\ell-1})\big]\\
|W_{v}'^{\circ}(I_{2})|^{2}L_{v}(1/2+s_{1},\pi_{v}\times\sigma_{v})L_{v}(1/2+s_{2},\sigma_{v})\cdot \mathbf{1}_{r_{\sigma_{v}}=0}.
\end{multline} 

Using $w_{2}n_{13}(b)w_{2}=n_{12}(b)$ and sphericity, substitute
\eqref{11.46} into \eqref{fc11.41}. The terms cancel except for
$\lambda_{\sigma_{v}}(\mathfrak{p}_{v}^{n_{v}})$, proving \eqref{fc11.42}.
\end{proof}

\subsubsection{Ramified calculation}
For $v\mid\mathfrak{q}$, the calculation in
\cite[(29), (43)]{Nun23}, applied to $W_{v}^{[r_{v}]}$, yields
\begin{multline}\label{equ12.11}
\mathcal{H}_{W_{v}}(\sigma_{v},\mathbf{s},\mathbf{1})=|W_{v}'^{\circ}(I_{2})|^{2}L_{v}(1/2+s_{2},\sigma_{v})\\
q_{v}^{r_{v}}\sum_{\substack{d_{1} + d_{2} = r_{v}\\
d_{1},\, d_{2}\geq 0}} 
\Vol(K_{0,v}'[d_{2}])
\sum_{l=0}^{d_{2} - r_{\sigma_{v}}} q_{v}^{-l}\sum_{k_{1}=0}^{l} \sum_{k_{2}=0}^{l}
\xi_{\sigma_{v}}(\mathfrak{p}_{v}^{k_{1}}, \mathfrak{p}_{v}^{l})\xi_{\sigma_{v}}(\mathfrak{p}_{v}^{k_{2}}, \mathfrak{p}_{v}^{l})
\\
q_{v}^{k_{1}(1/2-s_{2})}q_{v}^{k_{2}(1/2-s_{1})}q_{v}^{-2d_{1}(1/2+s_{1})}
\sum_{\substack{l_{1},l_{2}\geq0\\\min(l_{1},d_{2})=0}}
\frac{
\lambda_{\pi_{v}}(\mathfrak{p}_{v}^{l_{2}+k_{2}+l_{1}+d_{1}}, \mathfrak{p}_{v}^{l_{1}+d_{1}})
\overline{\lambda_{\sigma_{v}}(\mathfrak{p}_{v}^{l_{2}})}
}{
q_{v}^{(2l_{1} + l_{2})(1/2+s_{1})}
}.
\end{multline}

For $r_{\sigma_{v}}>r_{v}$, the $l$-sum in \eqref{equ12.11} is empty.
For $r_{\sigma_{v}}=r_{v}\geq1$, take $d_{2}=r_{v}$ and $d_{1}=l=k_{1}=k_{2}=0$.
The $(l_{1},l_{2})$-sum is $L_{v}(1/2+s_{1},\pi_{v}\times\sigma_{v})$, giving
\begin{multline}\label{fc12.16}
\mathcal{H}_{W_{v}}(\sigma_{v},\mathbf{s},\mathbf{1})
=q_{v}^{r_{v}}\Vol(K_{0,v}'[r_{v}])|W_{v}'^{\circ}(I_{2})|^{2}\\
L_{v}(1/2+s_{1},\pi_{v}\times\sigma_{v})L_{v}(1/2+s_{2},\sigma_{v}).
\end{multline}

\subsection{The Dual Weight}\label{sec12.4}
The initial convergence conditions are
\begin{align*}
 \Re((s_{2}-s_{1})/2)>-1/2+\vartheta_{\pi_{v}}+\vartheta_{\sigma_{v}},\quad \Re((3s_{1}+s_{2})/2)>-1/2+\vartheta_{\sigma_{v}}.
\end{align*}

Suppose $v\nmid\mathfrak{q}$.
Applying $\pi_{v}(w_{2})$ in \eqref{fc11.41} and using $w_{2}^{2}=I_{3}$ gives
\begin{equation}\label{eq:dual-level-superposition}
 \mathcal{H}_{W_{v}}^{\vee}(\sigma_{v},\mathbf{s}^{\vee},\mathbf{1})
 =q_{v}^{-n_{v}(1/2+s_{2})}\sum_{\ell=0}^{n_{v}}q_{v}^{2\ell s_{2}}
 \mathcal{H}_{W_{v}^{[\ell]}}(\sigma_{v},\mathbf{s}^{\vee},\mathbf{1}).
\end{equation}
For $\ell\geq1$, use \eqref{equ12.11} with level $\ell$;
for $\ell=0$, use \eqref{fc11.42} with $n_{v}=0$,
both at $\mathbf{s}^{\vee}$.

If $v\mid\mathfrak{q}$, the identity
$w_{2}n_{13}(b)w_{2}=n_{12}(b)$ and sphericity give
$\pi_{v}(w_{2})W_{v}=q_{v}^{r_{v}/2}W_{v}^{(r_{v})}$.
Equation~\eqref{11.46} at $\mathbf{s}^{\vee}$ therefore gives
\begin{multline}\label{eq12.14}
\mathcal{H}_{W_{v}}^{\vee}(\sigma_{v},\mathbf{s}^{\vee},\mathbf{1})=q_{v}^{(1-3s_{1}-s_{2})r_{v}/2}\big[\lambda_{\sigma_{v}}(\mathfrak{p}_{v}^{r_{v}})-q_{v}^{-(1+3s_{1}+s_{2})/2}\lambda_{\sigma_{v}}(\mathfrak{p}_{v}^{r_{v}-1})\big]\\
|W_{v}'^{\circ}(I_{2})|^{2}L_{v}((1+s_{2}-s_{1})/2,\pi_{v}\times\sigma_{v})L_{v}((1+3s_{1}+s_{2})/2,\sigma_{v})\mathbf{1}_{r_{\sigma_{v}}=0}.
\end{multline}

\part{Arithmetic Applications}
We combine the spectral reciprocity formulas of Part~1 with the
explicit local weights and dual-transform estimates of Part~2.
Type~\RNum{1} gives the first-moment bound for cuspidal and
maximal Eisenstein $\mathrm{GL}_{3}$ data. Type~\RNum{2} gives the twisted
fourth moment, the fifth-moment bound and the mixed-moment applications.

\section{\texorpdfstring{The First Moment of $\mathrm{GL}_{3}\times\mathrm{GL}_{2}$ $L$-functions}{The First Moment of GL3 x GL2 L-functions}}\label{sec13}
We prove Theorem~\ref{thm1.3} using the type~\RNum{1} formula and the
local estimates of \textsection\ref{sec9} and \textsection\ref{sec10}.
We first treat cuspidal $\pi$; maximal Eisenstein data are treated in
\textsection\ref{sec13.maximal}.

Suppose $\pi$ is cuspidal and take $\omega=\omega'=\mathbf{1}$.
Use $\mathfrak{m}$ and $Q$ from Theorem~\ref{thm1.3}, and fix a
sufficiently small $\varepsilon>0$.
At infinity, choose the vector \eqref{eq9.2} at $s=0$ with scale $C_{v}$.
For $C_{v}$ sufficiently large in terms of $F,\varepsilon$,
Corollary~\ref{cor9.2} gives a fixed positive lower bound on
$C_{v}(\sigma)\leq C_{v}^{1-\varepsilon}$.
Replacing each test scale by a sufficiently large multiple of
$C_{v}^{1/(1-\varepsilon)}$ covers $C_{v}(\sigma)\leq C_{v}$ in
\eqref{fc1.1}, after decreasing the initial $\varepsilon$.
At a finite place, use the vector of \textsection\ref{sec10} with
$r_{v}=\max\{e_{v}(\mathfrak{q}),r_{\pi_{v}}\}$, taking the spherical vector
outside $\mathfrak{q}\mathfrak{m}\mathfrak{D}_{F}$.

Let $\varphi$ denote the cusp form in $\pi$ whose
$\overline{\psi}$-Whittaker coefficient is
$W=\bigotimes_{v}' W_{v}$, with the convention of \textsection\ref{sec9}.
Theorem~\ref{thmA} gives
\begin{equation}\label{13.1}
J_{\mathrm{spec}}^{\heartsuit}(\varphi,\omega')=J_{\mathrm{sing}}^{\heartsuit}(\varphi,\omega')+J_{\mathrm{dual}}^{\heartsuit}(\varphi,\omega')+R_{\mathrm{dual}}^{+}(\varphi,\omega').
\end{equation}
The degenerate terms vanish by cuspidality. The character sets for
$R_{\RNum{1}}^{\pm}$ and $R_{\mathrm{dual}}^{-}$ are empty, since
every standard twist of $\pi$ is entire.
The singular term and residue are holomorphic at $s=0$.
Theorem~\ref{thm9.5} and the support conditions in \textsection\ref{sec10}
give locally uniform convergence of the dual integral there,
so \eqref{13.1} holds at $s=0$.

\subsection{\texorpdfstring{Bounds for $R_{\mathrm{dual}}^{+}(\varphi,\omega')$}{Bounds for R\_dual+}}
\begin{lemma}\label{lem13.1}
For cuspidal $\pi$ with $\omega=\omega'=\mathbf{1}$,
\begin{equation}\label{e13.2}
 |R_{\mathrm{dual}}^{+}(\varphi,\mathbf{1})|
 \ll_{F,\varepsilon} Q^{1/2+\varepsilon}\mathbf{C}_{\infty}^{1/2+\varepsilon}
 C(\pi)^{1/4+\varepsilon}
 \prod_{v\mid\infty}[1+(C_{v}^{-1}C_{v}(\pi))^{2}].
\end{equation}
\end{lemma}

\begin{proof}
By Lemma~\ref{lem11.4}, $\widetilde{\mathcal{J}}$ vanishes unless
$\xi_{v}$ is unramified at every $v\nmid\mathfrak{q}\mathfrak{m}\mathfrak{D}_{F}$.
For the remaining $\xi$, \eqref{eqq9.2} gives
\begin{multline}\label{13.7}
\widetilde{\mathcal{J}}(s,\lambda;\varphi,\xi,\omega')
=L(1/2+2s-\lambda,\widetilde{\pi}\times\overline{\xi}\omega\omega')
 L(1/2+\lambda,\xi)\\
\prod_{v\mid\infty}\mathcal{H}_{W_{v}}
 (\xi_{v},\lambda;\omega_{v}',s)
\prod_{v\mid\mathfrak{q}\mathfrak{m}\mathfrak{D}_{F}}
 \mathcal{H}_{W_{v}}^{\sharp}(\xi_{v},\lambda;\omega_{v}',s),
\end{multline}
with the normalized finite local weights $\mathcal{H}_{W_{v}}^{\sharp}$
from \textsection\ref{def:normalized-finite-dual-weight}.
Take the residue at $\lambda=1/2$ first. Since $\pi$ is
cuspidal, the only pole there is that of $\zeta_{F}(1/2+\lambda)$.
Consequently,
\begin{multline}\label{13.8}
 R_{\mathrm{dual}}^{+}(s)=\rho_{F}^{-1}L(2s,\widetilde{\pi})\Res_{u=1}\zeta_{F}(u)
 \prod_{v\mid\infty}\mathcal{H}_{W_{v}}(\mathbf{1}_{v},1/2;\mathbf{1}_{v},s)\\
\prod_{v\mid\mathfrak{q}\mathfrak{m}\mathfrak{D}_{F}}
 \mathcal{H}_{W_{v}}^{\sharp}(\mathbf{1}_{v},1/2;\mathbf{1}_{v},s).
\end{multline}
The functional equation for $L(2s,\widetilde{\pi})$ puts
$\epsilon(2s,\widetilde{\pi}_{v},\psi_{v})$ next to each finite weight in
\eqref{13.8}. For sufficiently small fixed $\delta>0$,
\begin{align*}
 \prod_{v\mid\mathfrak{q}\mathfrak{m}\mathfrak{D}_{F}}
 |\epsilon(2s,\widetilde{\pi}_{v},\psi_{v})
   \mathcal{H}_{W_{v}}^{\sharp}(\mathbf{1}_{v},1/2;\mathbf{1}_{v},s)|
 \ll_{F,\varepsilon}Q^{1/2+\varepsilon}
 C_{\mathrm{fin}}(\pi)^{1/4+\varepsilon},
 \qquad |s|\leq\delta.
\end{align*}
Away from the different, the local bound
$q_{v}^{r_{v}/2+r_{\pi_{v}}/4+O(\delta r_{v})}$ follows for degrees at most
two from Lemmas~\ref{lem:finite-projector-full-Tate-tail}
and~\ref{lem:maximal-Steinberg-positive-residue}.
In degree three, $r_{\pi_{v}}=0$ and
\eqref{eq:spherical-projector-residue-contour} at $\lambda=1/2$ gives the
same bound since $\vartheta_{\pi_{v}}\leq r_{v}/2$.
For $v\mid\mathfrak{D}_{F}$ and $r_{v}>0$, the transform is unchanged by
\eqref{eq:maximal-positive-residue-different-transport} at
$(\xi_{v},\lambda)=(\mathbf{1}_{v},1/2)$. By \eqref{eq:different-epsilon-factor}
with $n=3$, the epsilon factor changes in modulus by
$q_{v}^{3d_{v}(1/2-2\Re(s))}$, which is bounded in terms of $F$.
If $r_{v}=0$, use \eqref{eq:finite-zero-level-dual-bound} and the same
change of epsilon factor.

At infinity, the local $\mathrm{GL}_{3}\times\mathrm{GL}_{1}$ functional equation gives
\begin{align*}
 \gamma(2s,\widetilde{\pi}_{v},\psi_{v})
 \mathcal{H}_{W_{v}}(\mathbf{1}_{v},1/2;\mathbf{1}_{v},s)
 =\zeta_{v}(1)\int_{F_{v}^{\times}}
 W_{v}(\diag(a,1,1)w_{2}\widetilde{w})
 |a|_{v}^{-2s}\,d^{\times}a.
\end{align*}
Use \eqref{eq:reflected-point-trace} at $m=q_{0}=0$, with $p=a$.
Apply Cauchy--Schwarz on $|a|_{v}\leq S$ with $l=0$ and on
$|a|_{v}>S$ with $l$ large, where $S=C_{v}(1+C_{v}(\pi)/C_{v})^{2}$.
The norm of $|a|_{v}^{1/2-2\Re(s)-l/d}$ in $L^{2}(d^{\times} a)$
on these two ranges is $O(S^{1/2-2\Re(s)-l/d})$.
Multiplying by $O(S^{l/d})$ from \eqref{eq:reflected-point-trace}
gives $O(S^{1/2-2\Re(s)})$. For sufficiently small $\delta$,
this yields
\begin{equation}\label{eq:positive-residue-archimedean-bound}
 |\gamma(2s,\widetilde{\pi}_{v},\psi_{v})
 \mathcal{H}_{W_{v}}(\mathbf{1}_{v},1/2;\mathbf{1}_{v},s)| \ll_{F,\varepsilon}C_{v}^{1/2+\varepsilon}C_{v}(\pi)^{1/4+\varepsilon}
 [1+(C_{v}(\pi)/C_{v})^{2}].
\end{equation}
This local argument uses only $\omega_{v}=\omega_{v}'=\mathbf{1}_{v}$,
so it also applies to the maximal Eisenstein data below.
After the functional equation, $L(1-2s,\pi)\ll C(\pi)^{\varepsilon}$
on a sufficiently small circle about $s=0$. Multiply the finite bound
by \eqref{eq:positive-residue-archimedean-bound} over $v\mid\infty$.
Cauchy's formula for \eqref{13.8} then gives \eqref{e13.2}.
\end{proof}

\subsection{Proof of Theorem \ref{thm1.3}}\label{sec13.2}
For cuspidal $\sigma$, the central values are nonnegative by
\cite[Theorem~1.1]{Lap03}. For
$\sigma=\xi|\cdot|^{it}\boxplus\overline{\xi}|\cdot|^{-it}$,
the value is $|L(1/2+it,\pi\otimes\xi)|^{2}$.
The local spectral weights are nonnegative by
\eqref{e9.2} and \eqref{eq11.13}. By \cite[Lemma~2.2.3]{MV10},
the global basis sum factors into the local sums of
\textsection\ref{sec3.2}, divided by
$L^{\mathrm{reg}}(1,\sigma,\Ad)$, up to positive constants
depending only on $F$.

Corollary~\ref{cor9.2} gives a fixed positive lower bound at infinity,
and \eqref{fc11.12} and \eqref{eq:finite-newvector-value} give
a product lower bound $\gg_{F,\varepsilon}Q^{-\varepsilon}$
at $v\mid\mathfrak{q}\mathfrak{m}\mathfrak{D}_{F}$; the places
dividing the different contribute only an $F$-dependent constant.
Passing from the unramified product in \eqref{fc11.2} to
$L(1/2,\pi\times\sigma)$ introduces the positive factors $L_{v}(1/2,\pi_{v}\times\sigma_{v})^{-1}$
at $v\mid\mathfrak{q}\mathfrak{m}\mathfrak{D}_{F}$. Their product is
$\gg_{F,\varepsilon}Q^{-\varepsilon}$, since
$\vartheta_{\pi_{v}}+\vartheta_{\sigma_{v}}<1/2$.
Exclude
\begin{align*}
 \mathcal{E}_{0}(\mathfrak{q})
 :=\{\sigma_{\xi,it}\in\mathcal{F}(\mathfrak{q};\mathbf{1}):\ 
       \xi^{2}=\mathbf{1},\ |t|\leq1\}.
\end{align*}
Outside this set, the adjoint bound and \eqref{eq:regularized-adjoint} give
$L^{\mathrm{reg}}(1,\sigma,\Ad)\ll_{F,\varepsilon}
C(\sigma)^{\varepsilon}$. Positivity therefore gives
\begin{equation}\label{13.2}
J_{\mathrm{spec}}^{\heartsuit}(\varphi,\mathbf{1})\gg
(C(\pi)\mathbf{C}_{\infty}Q)^{-\varepsilon}
\int_{\substack{\sigma\in\mathcal{F}(\mathfrak{q};\mathbf{1})
                 \setminus\mathcal{E}_{0}(\mathfrak{q})\\
 C_{v}(\sigma)\leq C_{v}^{1-\varepsilon},\ v\mid\infty}}
L(1/2,\pi\times\sigma)d\mu_{\sigma}.
\end{equation}

For $\sigma_{\xi,it}\in\mathcal{E}_{0}(\mathfrak{q})$, we have
$2r_{\xi_{v}}\leq e_{v}(\mathfrak{q})\leq r_{v}$.
Unramified twists preserve the conductor exponent. If $\xi_{v}$ is
ramified and the Weil parameter of $\pi_{v}$ is nontrivial on every
ramification group where $\xi_{v}$ is nontrivial, its trivial determinant
makes the fixed space on each such group at most one-dimensional.
Thus the tame term and each wild codimension increase by at most one,
giving $r_{\pi_{v}\otimes\xi_{v}}\leq r_{\pi_{v}}+r_{\xi_{v}}$.
Otherwise the twist has no invariants on the groups where $\xi_{v}$
is nontrivial and is trivial on the higher groups, so its exponent
is $3r_{\xi_{v}}$. Hence $r_{\pi_{v}\otimes\xi_{v}}\leq 3r_{v}/2$ and
$C_{\mathrm{fin}}(\pi\otimes\xi)\leq Q^{3/2}$ in both cases of
Theorem~\ref{thm1.3}.

The Eisenstein central value is $|L(1/2+it,\pi\otimes\xi)|^{2}$ in both
cases. For $2+1$ data, use $L(s,\pi\otimes\xi)=L(s,\tau\otimes\xi)L(s,\chi\xi)$
and the product of the two conductors. Since $|t|\leq1$ and
$C_{\infty}(\pi\otimes\xi)\asymp_{F}C_{\infty}(\pi)$, convexity gives
\begin{equation}\label{eq:exceptional-quadratic-convexity}
 |L(1/2,\pi\times\sigma_{\xi,it})|
 \ll_{F,\varepsilon}
 \bigl(Q^{3/4}C_{\infty}(\pi)^{1/2}\bigr)^{1+\varepsilon}.
\end{equation}
The arithmetic--geometric mean inequality bounds this by the first two
terms of \eqref{fc1.1}, with the product in its first term omitted.
There are $O_{F,\varepsilon}(Q^{\varepsilon})$ quadratic characters,
and $|t|\leq1$.

For the singular term, \eqref{eq:finite-projector-restriction} gives,
at $v\mid\mathfrak{q}\mathfrak{m}\mathfrak{D}_{F}$ and for every $s$,
\begin{equation}\label{eq:singular-finite-exact}
 J_{\mathrm{sing},v}(s)=\Vol(K_{0,v}'[r_{v}])^{-1}
 \ll q_{v}^{r_{v}}.
\end{equation}
At infinity, the compact support in \eqref{eq9.2}
gives $J_{\mathrm{sing},v}(0)\ll_{F,\varepsilon}C_{v}^{1+\varepsilon}$.
By Lemma~\ref{lem7.1}, the remaining finite Euler product is
$L(1,\widetilde{\pi})$ with the factors at
$v\mid\mathfrak{q}\mathfrak{m}\mathfrak{D}_{F}$ omitted, and is
$O(C(\pi)^{\varepsilon})$. Hence
\begin{equation}\label{13.9}
\left|J_{\mathrm{sing}}^{\heartsuit}(\varphi,\mathbf{1})\right|\ll
C(\pi)^{\varepsilon}\mathbf{C}_{\infty}^{1+\varepsilon}
Q^{1+\varepsilon}.
\end{equation}

Combine \eqref{13.1}, Lemma~\ref{lem13.1}, \eqref{13.2},
Corollary~\ref{cor:trivial-central-character-dual-moment},
\eqref{13.9}, and the bound for the integral over $\mathcal{E}_{0}(\mathfrak{q})$. Replacing $C_{v}$ by a sufficiently large
multiple of $C_{v}^{1/(1-\varepsilon)}$ and taking a smaller initial
$\varepsilon$ proves \eqref{fc1.1} for cuspidal
$\pi$, with the product in its first term omitted.

\subsection{Maximal Eisenstein Data}\label{sec13.maximal}
Take $\pi=\tau\boxplus\chi$ as in Theorem~\ref{thm1.3}, and set
$\omega=\omega'=\mathbf{1}$. Thus $\tau$ is unitary cuspidal of
orthogonal type, and
\begin{align*}
 &\widetilde{\tau}\simeq\tau,\qquad
 \chi^{2}=\mathbf{1},\qquad \omega_{\tau}=\chi,\qquad
 \tau\otimes\chi\simeq\tau,\\
 &L(1/2,\pi\times\sigma)
 =L(1/2,\tau\times\sigma)L(1/2,\sigma\otimes\chi).
\end{align*}

\subsubsection{The test function}\label{sec:maximal-test-function}

Let $S$ consist of the Archimedean places and the finite places dividing
$\mathfrak{q}\mathfrak{m}\mathfrak{D}_{F}$, and let $L^{S}$ omit their factors.
For $v\in S$, let $W_{v}$ be the $\overline{\psi}_{v}$-Whittaker vector
of $\pi_{v}$ chosen in \textsection\ref{sec9}--\textsection\ref{sec10},
with restriction \eqref{eq9.2} at infinity and \eqref{c11.1} at finite
places, taking $s=0$, scales $C_{v}$, and
$r_{v}=\max\{e_{v}(\mathfrak{q}),r_{\pi_{v}}\}$.
The function $x\mapsto W_{v}(\diag(x,1))$ is smooth and compactly
supported modulo $N'(F_{v})$. Keep these restrictions fixed as
$\boldsymbol{\mu}$ and the moment parameter $s$ vary.

\begin{lemma}\label{lem:maximal-test-sections}
For every $v\in S$, there is a holomorphic family of sections
$f_{\boldsymbol{\mu},v}\in
\chi_{v}|\cdot|_{v}^{\mu_{1}}\boxplus\tau_{v}|\cdot|_{v}^{\mu_{2}/2}$,
$\boldsymbol{\mu}\in\mathbb{C}^{2}$, such that, for
$x\in G'(F_{v})$ and $g\in G(F_{v})$,
\begin{align*}
 \int_{F_{v}^{2}}f_{\boldsymbol{\mu},v}\left(w_{1}w_{2}
 n_{13}(b)n_{23}(c)\diag(x,1)\right)(I_{2})
 \psi_{v}(c)\,db\,dc=W_{v}(\diag(x,1)),
\end{align*}
and $f_{\boldsymbol{\mu},v}(g)=0$ if $(g^{-1})_{31}=0$. Here $\tau_{v}$ is realized in its $\overline{\psi}_{v}$-Whittaker model.
\end{lemma}

\begin{proof}
Fix $v\in S$. First suppose the restriction is supported in $x_{22}\neq0$.
For $c'\in F_{v}$ and $y_{2}\in F_{v}^{\times}$, let
$W'_{\mu_{2},v}(\cdot;c',y_{2})$ be the vector in $\tau_{v}$ with
Kirillov restriction
\begin{align*}
 W'_{\mu_{2},v}\left(\diag(y_{1},1);c',y_{2}\right)
 :=|y_{1}y_{2}|_{v}^{-\frac{1+\mu_{2}}{2}}
 W_{v}\left(\begin{pmatrix}y_{1}&&\\c'&y_{2}&\\&&1\end{pmatrix}\right),\quad y_{1}\in F_{v}^{\times}.
\end{align*}
The given restriction has fixed compact support in
$(y_{1},c',y_{2})\in F_{v}^{\times}\times F_{v}\times F_{v}^{\times}$;
in particular, $y_{1}$ and $y_{2}$ stay away from zero.
The Kirillov-model theorem therefore gives these vectors in $\tau_{v}$; see
\cite[Proposition~2.9 and Theorem~2.13]{JL70} at finite places and
\cite[\textsection3.2.5]{MV10} at infinity.
For $g'\in G'(F_{v})$, define
\begin{align*}
 &f_{\boldsymbol{\mu},v}\left(w_{1}w_{2}
 n_{13}(b)n_{23}(c)\right)(g')\\
 &:=\int_{F_{v}}\int_{F_{v}^{\times}}
 W'_{\mu_{2},v}\left(g'\begin{pmatrix}1&\\c'&y_{2}\end{pmatrix}^{-1};c',y_{2}\right)
 \overline{\psi}_{v}(c'b+y_{2}c)\,dy_{2}\,dc'.
\end{align*}

Every $g$ with $(g^{-1})_{31}\neq0$ can be written uniquely as
\begin{align*}
 g=\begin{pmatrix}a&\mathfrak{b}\\&g'\end{pmatrix}w_{1}w_{2}
 n_{13}(b)n_{23}(c),\qquad
 (g^{-1})_{31}=a^{-1},
\end{align*}
where $a\in F_{v}^{\times}$, $g'\in G'(F_{v})$,
$\mathfrak{b}\in M_{1\times2}(F_{v})$ and $b,c\in F_{v}$.
Define the section on this set by \eqref{e2.3}:
\begin{align*}
 f_{\boldsymbol{\mu},v}(g)
 :=\chi_{v}(a)|a|_{v}^{1+\mu_{1}}|\det g'|_{v}^{(-1+\mu_{2})/2}\tau_{v}(g')f_{\boldsymbol{\mu},v}\left(w_{1}w_{2}
 n_{13}(b)n_{23}(c)\right).
\end{align*}
On the remaining set, put
\begin{align*}
 f_{\boldsymbol{\mu},v}(g)=0\qquad\text{if }(g^{-1})_{31}=0.
\end{align*}

By continuity of the Kirillov-model inclusion, the integrand defining
$f_{\boldsymbol{\mu},v}$ is a smooth, compactly supported
$\tau_{v}$-valued function of $(c',y_{2})$, holomorphic in $\mu_{2}$.
Its Fourier transform is therefore compactly supported in $(b,c)$
at finite places and Schwartz at infinity.
These statements hold in every smooth seminorm of $\tau_{v}$,
locally uniformly in $\boldsymbol{\mu}$; in particular, the finite-place
support is independent of $g'$.

In the decomposition defining $f_{\boldsymbol{\mu},v}$,
\begin{align*}
 b=-\frac{(g^{-1})_{11}}{(g^{-1})_{31}},\qquad
 c=-\frac{(g^{-1})_{21}}{(g^{-1})_{31}}.
\end{align*}
At a point with $(g^{-1})_{31}=0$, the other two entries of the
first column of $g^{-1}$ cannot both be zero. Hence, near that point,
\begin{align*}
 \max\{|b|_{v},|c|_{v}\}
 =\frac{\max\{|(g^{-1})_{11}|_{v},|(g^{-1})_{21}|_{v}\}}
 {|(g^{-1})_{31}|_{v}}
 \asymp|(g^{-1})_{31}|_{v}^{-1}.
\end{align*}
At finite places, the bound on $(b,c)$ in the support therefore forces
the section to vanish near this point. At infinity, the inducing factors
in \eqref{e2.3} and the action
of $\tau_{v}(g')$ have polynomial growth in $|(g^{-1})_{31}|_{v}^{-1}$,
also after differentiation, by moderate growth of $\tau_{v}$.
The Schwartz estimates therefore imply, for every fixed differential
operator $\mathcal{D}$ and every $m\geq0$,
\begin{align*}
 \bigl\|\mathcal{D}f_{\boldsymbol{\mu},v}(g)\bigr\|
 \ll_{\mathcal{D},m,\|\cdot\|}|(g^{-1})_{31}|_{v}^{m},
\end{align*}
for every fixed smooth seminorm $\|\cdot\|$ of $\tau_{v}$,
locally uniformly in $\boldsymbol{\mu}$.
Thus all derivatives of the section extend by zero across the boundary.
The dependence on $\boldsymbol{\mu}$ remains holomorphic because the
defining integral and inducing factors are holomorphic, with the same
locally uniform bounds.

For $x=\begin{pmatrix}y_{1}&\\c'&y_{2}\end{pmatrix}$, we have $x\begin{pmatrix}1&\\c'&y_{2}\end{pmatrix}^{-1}=\diag(y_{1},1)$ and 
\begin{align*}
 w_{1}w_{2}\begin{pmatrix}1&&y_{1}b\\&1&c'b+y_{2}c\\&&1\end{pmatrix}\diag(x,1)
 =\diag(1,x)w_{1}w_{2}n_{13}(b)n_{23}(c).
\end{align*}

In the integral in the lemma, make the change
$(b,c)\mapsto(y_{1}b,c'b+y_{2}c)$, with Jacobian $|y_{1}y_{2}|_{v}$.
The matrix identity, \eqref{e2.3}, and Fourier inversion in $(b,c)$ give
\begin{multline*}
 \int_{F_{v}^{2}}f_{\boldsymbol{\mu},v}\left(w_{1}w_{2}
 n_{13}(b)n_{23}(c)\diag(x,1)\right)(I_{2})
 \psi_{v}(c)\,db\,dc\\
 =|y_{1}y_{2}|_{v}^{(1+\mu_{2})/2}\int_{F_{v}^{2}}
 f_{\boldsymbol{\mu},v}\left(w_{1}w_{2}
 n_{13}(b)n_{23}(c)\right)(x)
 \psi_{v}(c'b+y_{2}c)\,db\,dc\\
 =|y_{1}y_{2}|_{v}^{(1+\mu_{2})/2}W'_{\mu_{2},v}(\diag(y_{1},1);c',y_{2})
 =W_{v}(\diag(x,1)).
\end{multline*}
For any $x\in G'(F_{v})$ with $x_{22}\neq0$,
\begin{align*}
 x=n(x_{12}/x_{22})
 \begin{pmatrix}\det(x)/x_{22}&\\x_{21}&x_{22}\end{pmatrix}.
\end{align*}
Left multiplication by $n(x_{12}/x_{22})$ multiplies both the Jacquet
integral and $W_{v}$ by $\overline{\psi}_{v}(x_{12}/x_{22})$,
so the equality holds for these $x$ as well.
The compact support at finite places and the Schwartz bounds at infinity
give locally uniform absolute convergence of the Jacquet integral in $x$,
so it is continuous.
Since $x_{22}\neq0$ is dense in $G'(F_{v})$, the equality holds for every
$x\in G'(F_{v})$.

For a general restriction, use a partition of unity in the nonzero
row $(x_{21},x_{22})$ to write
\begin{align*}
 W_{v}(\diag(x,1))
 =W_{v}^{(1)}(\diag(x,1))+W_{v}^{(2)}(\diag(x,1)),
\end{align*}
where the two restrictions are smooth and compactly supported modulo
$N'(F_{v})$ in $x_{22}\neq0$ and $x_{21}\neq0$, respectively.
Both functions
\begin{align*}
 x\longmapsto W_{v}^{(1)}(\diag(x,1)),\qquad
 x\longmapsto W_{v}^{(2)}(\diag(xw',1))
\end{align*}
are supported in $x_{22}\neq0$, since $(xw')_{21}=x_{22}$.
Let $f_{\boldsymbol{\mu},v}^{(1)}$ and $f_{\boldsymbol{\mu},v}^{(2)}$
be the sections constructed above for these two functions, and set
\begin{align*}
 f_{\boldsymbol{\mu},v}(g)
 :=f_{\boldsymbol{\mu},v}^{(1)}(g)+f_{\boldsymbol{\mu},v}^{(2)}(gw_{1}).
\end{align*}
Since $\diag(x,1)w_{1}=\diag(xw',1)$ and $(w')^{2}=I_{2}$,
the integral in the lemma for this section is
\begin{align*}
 W_{v}^{(1)}(\diag(x,1))+W_{v}^{(2)}(\diag(xw'w',1))
 =W_{v}(\diag(x,1)).
\end{align*}
Finally, $((gw_{1})^{-1})_{31}=(g^{-1})_{31}$, so both summands
vanish when $(g^{-1})_{31}=0$.
\end{proof}

Use the sections of Lemma~\ref{lem:maximal-test-sections} at $v\in S$. At the remaining finite places take
the spherical sections in \eqref{6.1}, with
Jacquet integrals equal to $1$ at $I_{3}$.
Form the global section $f_{\boldsymbol{\mu}}$ in the automorphic
realization of $\tau$, and put, as in \eqref{e2.4},
\begin{align*}
 \varphi_{\boldsymbol{\mu}}(g)
 :=E(g,f_{\boldsymbol{\mu}})
 =\sum_{\delta\in Q(F)\backslash G(F)}f_{\boldsymbol{\mu}}(\delta g)(I_{2}),
 \qquad \varphi:=\varphi_{\mathbf{0}}.
\end{align*}
We use this family in \eqref{eq:inducing-constant-term}.
By Lemma~\ref{lem:maximal-test-sections}, injectivity of the Whittaker
restriction, and the spherical normalization, the local Whittaker
vectors of $\varphi$ are the test vectors $W_{v}$ of
\textsection\ref{sec9}--\textsection\ref{sec10}.

\subsubsection{The maximal Eisenstein case of Theorem~\ref{thm1.3}}\label{sec:proof-maximal-first-moment}
Apply \eqref{e3.3} to $\varphi_{\boldsymbol{\mu}}$, whose central
character is $|\cdot|^{\mu_{1}+\mu_{2}}$, in the convergence region.
Take the constant term at $\boldsymbol{\mu}=0$ as in
\eqref{eq:inducing-constant-term}. Continuing in $s$ by
Theorem~\ref{thmA}, remove the vanishing terms by
Lemma~\ref{lem:maximal-degenerate-vanishing}:
\begin{equation}\label{eq:maximal-first-moment-reduction}
 J_{\mathrm{cusp}}^{\heartsuit}+J_{\mathrm{Eis}}^{\heartsuit}
 =J_{\mathrm{sing}}^{\heartsuit}
  +J_{\mathrm{degen}}^{\dag,\heartsuit}
  +J_{\mathrm{dual}}^{\heartsuit}+R_{\mathrm{dual}}^{+}-R_{\mathrm{dual}}^{-}
  -R_{\RNum{1}}^{+}+R_{\RNum{1}}^{-}.
\end{equation}
This is a meromorphic identity in $s$; the arguments
$(s,\varphi,\mathbf{1})$ are suppressed.
The character sets for $R_{\mathrm{dual}}^{-}$, $R_{\RNum{1}}^{+}$,
and $R_{\RNum{1}}^{-}$ are all $\{\chi\}$.
By \eqref{eq:opposite-spectral-residues},
\begin{align*}
 R_{\RNum{1}}^{-}(s,\varphi,\mathbf{1})
 =-R_{\RNum{1}}^{+}(s,\varphi,\mathbf{1}).
\end{align*}

The $L$-factors in \eqref{eq:maximal-positive-spectral-residue},
\eqref{eq:maximal-negative-dual-residue} and \eqref{eq:maximal-positive-dual-residue}
have only $O_{F}(1)$ possible poles in $|s|\leq3\delta$, also after applying
\eqref{eq:maximal-residue-functional-equation}. Choose a common
radius $\rho\in[\delta,2\delta]$ at distance $\gg_{F}\delta$ from
all their moduli, and use the local estimates with $2\delta$ in place
of $\delta$.

\begin{prop}
Put
\begin{align*}
 \mathcal{M}_{2+1}(s):=
 J_{\mathrm{sing}}^{\heartsuit}(s)
 +J_{\mathrm{degen}}^{\dag,\heartsuit}(s)
 +R_{\mathrm{dual}}^{+}(s)-R_{\mathrm{dual}}^{-}(s)-2R_{\RNum{1}}^{+}(s).
\end{align*}
Then $\mathcal{M}_{2+1}$ is holomorphic near $s=0$ and satisfies
\begin{multline}\label{eq:maximal-correction-bound}
 |\mathcal{M}_{2+1}(0)|\ll_{F,\varepsilon}Q^{\varepsilon}C(\pi)^{\varepsilon}
\mathbf{C}_{\infty}^{\varepsilon} \Big[Q\mathbf{C}_{\infty}\prod_{v\mid\infty}(1+C_{v}^{-1}C_{v}(\pi))\\
 {}+Q^{1/2}\mathbf{C}_{\infty}^{1/2}
 C(\pi)^{1/4}
 \prod_{v\mid\infty}[1+(C_{v}^{-1}C_{v}(\pi))^{2}]
 {}+\mathbf{C}_{\infty}^{1/2}
 C(\pi)^{1/2}
 \prod_{v\mid\infty}(1+C_{v}^{-1}C_{v}(\pi))\Big].
\end{multline}
\end{prop}

\begin{proof}

Corollary~\ref{cor4.4} and \textsection\ref{sec4.1.4} give holomorphy
of the spectral terms near $s=0$. Theorem~\ref{thm9.5}, the finite
conductor restriction and convexity give locally uniform absolute
convergence of the dual integral there, whose two $L$-factors are
holomorphic near their central lines. Thus
\eqref{eq:maximal-first-moment-reduction} makes $\mathcal{M}_{2+1}$
holomorphic and gives, at $s=0$,
\begin{equation}\label{eq:maximal-first-moment-grouped}
 J_{\mathrm{cusp}}^{\heartsuit}+J_{\mathrm{Eis}}^{\heartsuit}
 =J_{\mathrm{dual}}^{\heartsuit}+\mathcal{M}_{2+1}(0).
\end{equation}
For $|s|=\rho$, apply \eqref{eq:maximal-singular-bound},
\eqref{eq:maximal-opposite-degenerate-bound},
Corollaries~\ref{cor:maximal-positive-residue-sharp-bound}
and~\ref{cor:maximal-negative-residue-covered-range}, and
\eqref{eq:maximal-spectral-residue-absorbed-bound} below.
Cauchy's formula gives \eqref{eq:maximal-correction-bound}.
The extra product in its first term comes from $R_{\mathrm{dual}}^{-}$,
which is absent for cuspidal $\pi$.
\end{proof}

\begin{proof}[Proof of the maximal Eisenstein assertion of Theorem~\ref{thm1.3}]
For $\varphi=E(\cdot,f_{\mathbf{0}})$,
combine \eqref{eq:maximal-correction-bound} and
Corollary~\ref{cor:trivial-central-character-dual-moment} in
\eqref{eq:maximal-first-moment-grouped} to bound the weighted first moment.
Cuspidal $\sigma$ is symplectic because its central character is trivial.
Since $\tau$ and $\chi$ are orthogonal, \cite[Theorem~1.1]{Lap03} gives
nonnegativity of both central values. On the continuous spectrum they are
absolute squares.
Apply \eqref{13.2} to this $\varphi$ and bound the integral over
$\mathcal{E}_{0}(\mathfrak{q})$ by \eqref{eq:exceptional-quadratic-convexity}
as in the cuspidal case.
Replace each $C_{v}$ by a sufficiently large multiple of
$C_{v}^{1/(1-\varepsilon)}$ to cover $C_{v}(\sigma)\leq C_{v}$ in
\eqref{fc1.1}, and take the initial $\varepsilon$ smaller.
\end{proof}

\subsubsection{The degenerate and singular terms}

\begin{lemma}\label{lem:maximal-degenerate-vanishing}
Let $f_{\boldsymbol{\mu}}$ have the local sections of
Lemma~\ref{lem:maximal-test-sections} at $v\in S$ and the normalized spherical
sections of \eqref{6.1} elsewhere. Then
\begin{align*}
 J_{\mathrm{degen}}(s,\varphi_{\boldsymbol{\mu}},\mathbf{1})\equiv0
\end{align*}
as a meromorphic identity in $(s,\boldsymbol{\mu})$. In particular, $J_{\mathrm{degen}}^{\heartsuit}(s,\varphi,\mathbf{1})=0$. 
\end{lemma}

\begin{proof}
Fix an Archimedean place $v$. For $z\in F_{v}^{\times}$ and
$c,c'\in F_{v}$, the matrix
\begin{align*}
 n_{23}(c)w_{1}
 \diag(z,z,1)n_{12}(c')
 =\begin{pmatrix}&z&\\z&zc'&c\\&&1\end{pmatrix}
\end{align*}
has inverse with $(3,1)$-entry zero, so $f_{\boldsymbol{\mu},v}$
vanishes there by definition. Hence \eqref{6.5}, with $\psi$
replaced by $\overline{\psi}$, gives
\begin{multline*}
 W_{\varphi_{\boldsymbol{\mu}}}^{\mathrm{degen}}\left(w_{1}
 \diag(z,z,1)n_{12}(c')\right)\\
 =\int_{F\backslash\mathbb{A}_{F}}f_{\boldsymbol{\mu}}\left(
 n_{23}(c)w_{1}
 \diag(z,z,1)n_{12}(c')\right)(I_{2})\psi(c)\,dc=0.
\end{multline*}
Thus $J_{\mathrm{degen}}(s,\varphi_{\boldsymbol{\mu}},\mathbf{1})=0$
in the convergence region, and then identically by meromorphic
continuation. Taking the constant term at $\boldsymbol{\mu}=0$
proves the last assertion.
\end{proof}

\begin{lemma}
For every $\varepsilon>0$, there is a sufficiently small
$\delta=\delta(\varepsilon)>0$ and a radius
$\rho\in[\delta,2\delta]$ such that
\begin{equation}\label{eq:maximal-opposite-degenerate-bound}
 \sup_{|s|=\rho}
 \big|J_{\mathrm{degen}}^{\dag,\heartsuit}(s,\varphi,\mathbf{1})\big|
 \ll_{F,\varepsilon}C(\pi)^{\varepsilon}
 Q^{1+\varepsilon}\mathbf{C}_{\infty}^{\varepsilon}.
\end{equation}
Moreover, $J_{\mathrm{degen}}^{\dag,\heartsuit}=0$ if
$r_{v}>\max\{1,r_{\chi_{v}}\}$ at some finite place dividing
$\mathfrak{q}\mathfrak{m}\mathfrak{D}_{F}$.
\end{lemma}

\begin{proof}
At a finite place, take $h_{v}$ to be the normalized characteristic
function of $a(\varpi_{v}^{-d_{v}})K_{0,v}'[r_{v}]$ in
\eqref{eq:finite-projector-function}, with $\omega_{v}'=\mathbf{1}_{v}$.

At infinity, take
$h_{v}(g)=C_{v}h_{1,v}(g)h_{2,v}(g_{21}C_{v}^{1/[F_{v}:\mathbb{R}]})$,
where $h_{1,v}$ and $h_{2,v}$ are the fixed bump functions of
\textsection\ref{sec9.1.2}. Use the sections of
Lemma~\ref{lem:maximal-test-sections}, with the fixed restrictions
$W_{v}(\diag(x,1))$. The unnormalized integral in \eqref{e6.23} is
$\mathcal{L}_{v}^{\dag}\mathcal{P}_{v}^{\dag,\sharp}$, where
$\mathcal{L}_{v}^{\dag}=L_{v}(1+\mu_{2}+2s,\chi_{v})L_{v}(\mu_{1}-\mu_{2}/2,\tau_{v})$.

The substitution $(z_{v},b_{v},c_{v})=(c,-cb_{1}/(1+cb_{2}),(1+cb_{2})^{-1})$
in \eqref{e6.23} and \eqref{e2.3} with
$\omega_{\tau_{v}}=\chi_{v}$ give
\begin{align*}
\mathcal{L}_{v}^{\dag}\mathcal{P}_{v}^{\dag,\sharp}
=\chi_{v}(-1)\int_{F_{v}^{\times}}\int_{F_{v}^{2}}
 \chi_{v}(c)|c|_{v}^{1+\mu_{2}+2s}f^{*}(b_{1},b_{2},c)\psi_{v}(cb_{1})\,db_{1}\,db_{2}\,d^{\times}c,
\end{align*}
where 
\begin{align*}
f^{*}(b_{1},b_{2},c):=f_{\boldsymbol{\mu},v}\left(w_{1}w_{2}
 n_{13}(b_{1})n_{23}(b_{2})\right)
 \left(\diag(1+cb_{2},1)\right)
 |1+cb_{2}|_{v}^{\mu_{1}-\mu_{2}/2-3/2}.
\end{align*}

The $\mathrm{GL}_{2}$ local functional equation of
\textsection\ref{sec3.2}, followed by Fourier inversion at
$1+cb_{2}\neq0$, gives 
\begin{multline*}
 f^{*}(b_{1},b_{2},c) =\frac{1}{\zeta_{v}(1)\gamma(\mu_{1}-\mu_{2}/2,\tau_{v},\overline{\psi}_{v})}
 \int_{F_{v}}\int_{F_{v}^{\times}}
\chi_{v}(y)|y|_{v}^{1/2-\mu_{1}+\mu_{2}/2}\\
 f_{\boldsymbol{\mu},v}\left(w_{1}w_{2}
 n_{13}(b_{1})n_{23}(b_{2})\right)
 \left(\begin{pmatrix}&y\\1&t\end{pmatrix}\right)\,d^{\times}y\psi_{v}(t(1+cb_{2}))\,dt.
\end{multline*}
The outer $t$-integral is interpreted by Fourier inversion.
Insert the formula for $f^{*}$ in the $(c,b_{1},b_{2})$-integral for
$\mathcal{L}_{v}^{\dag}\mathcal{P}_{v}^{\dag,\sharp}$ and put
$b'=-c^{2}y$ and $a'=ct$, with $d^{\times} y=d^{\times} b'$ and
$dt=|c|_{v}^{-1}da'$ at fixed $c$. Use
$\omega_{\tau_{v}}=\chi_{v}$ and $\chi_{v}^{2}=\mathbf{1}_{v}$ for
the character factors.
To evaluate the $b_{1},b_{2}$-integral, set
$x=\left(\begin{smallmatrix}&-b'/c\\c&a'\end{smallmatrix}\right)$,
so $\det x=b'$. In Lemma~\ref{lem:maximal-test-sections}, make the change
\begin{align*}
(u_{1},u_{2})=(-b'b_{2}/c,cb_{1}+a'b_{2}),\qquad
 du_{1}\,du_{2}=|b'|_{v}\,db_{1}\,db_{2}.
\end{align*}
We have
\begin{align*}
 w_{1}w_{2}n_{13}(u_{1})n_{23}(u_{2})\diag(x,1)
 =\diag(1,x)w_{1}w_{2}n_{13}(b_{1})n_{23}(b_{2}).
\end{align*}
By \eqref{e2.3}, the left factor $\diag(1,x)$ contributes
$|b'|_{v}^{(\mu_{2}-1)/2}$ and changes evaluation at $I_{2}$ to
evaluation at $x$. After including the Jacobian, divide by
$|b'|_{v}^{(1+\mu_{2})/2}$ to obtain
\begin{align*}
\int_{F_{v}^{2}}f_{\boldsymbol{\mu},v}\left(w_{1}w_{2}
 n_{13}(b_{1})n_{23}(b_{2})\right)(x)
 \psi_{v}(cb_{1}+a'b_{2})\,db_{1}db_{2}=|b'|_{v}^{-\frac{1+\mu_{2}}{2}}W_{v}(\diag(x,1)).
\end{align*}
In the unipotent integral defining $W_{v}$, put $a=bc$, with
$db=|c|_{v}^{-1}da$ and $\psi_{v}(b)=\psi_{v}(a/c)$:
\begin{align*}
 W_{v}\left(\begin{pmatrix}&-b'/c&\\c&a'&\\&&1\end{pmatrix}\right)
 =|c|_{v}^{-1}\int_{F_{v}}(h_{v}\ast h_{v}^{*})\left(
 \begin{pmatrix}a&(aa'-b')/c\\c&a'\end{pmatrix}\right)
 \psi_{v}(a/c)\,da.
\end{align*}
Substitution in the $(c,b_{1},b_{2})$-integral for
$\mathcal{L}_{v}^{\dag}\mathcal{P}_{v}^{\dag,\sharp}$ now yields
\begin{equation}\label{eq:maximal-opposite-local-integral}
 \mathcal{L}_{v}^{\dag}\mathcal{P}_{v}^{\dag,\sharp}
 =\frac{\mathcal{I}_{v}(s,\boldsymbol{\mu})}
 {\zeta_{v}(1)\gamma(\mu_{1}-\mu_{2}/2,\tau_{v},\overline{\psi}_{v})},
\end{equation}
where
\begin{multline}\label{eq:maximal-opposite-kernel}
 \mathcal{I}_{v}(s,\boldsymbol{\mu}):=
 \int_{(F_{v}^{\times})^{2}}\int_{F_{v}^{2}}
 (h_{v}\ast h_{v}^{*})\left(\begin{pmatrix}a&(aa'-b')/c\\c&a'\end{pmatrix}\right)
 \chi_{v}(b')|b'|_{v}^{-\mu_{1}}\\
 |c|_{v}^{2s+2\mu_{1}-2}
 \psi_{v}((a+a')/c)\,da\,da'\,d^{\times}b'\,d^{\times}c.
\end{multline}
The Schwartz bounds of the sections and the compact Kirillov support
of $W_{v}$ justify Fourier inversion and rearrangement, with the additive
integrations performed before the $c$-integration. We first work in the
common convergence strip of the local zeta integrals with $\Re(s)$ large,
then continue meromorphically in $(s,\boldsymbol{\mu})$.

At a finite place with $r_{v}\geq1$, where $r_{v}$ is defined in
\textsection\ref{sec:maximal-test-function}, change $b'=aa'-bc$ in
\eqref{eq:maximal-opposite-kernel}; then
$d^{\times}b'=\zeta_{v}(1)|c|_{v}|aa'-bc|_{v}^{-1}db$.
By \eqref{eq:finite-ordinary-correlation}, the support has
$a,a'\in\mathcal{O}_{v}^{\times}$,
$b\in\varpi_{v}^{-d_{v}}\mathcal{O}_{v}$ and
$e_{v}(c)\geq r_{v}+d_{v}$. Since $r_{\chi_{v}}\leq r_{v}$,
$\chi_{v}(aa'-bc)=\chi_{v}(aa')$, so the $a$- and $a'$-integrals
are both
\begin{align*}
 \int_{\mathcal{O}_{v}^{\times}}\chi_{v}(a)\psi_{v}(a/c)\,da.
\end{align*}
A nonzero Gauss integral requires $e_{v}(c)=d_{v}+r_{\chi_{v}}$ for
$r_{\chi_{v}}>0$ and $e_{v}(c)\leq d_{v}+1$ for $r_{\chi_{v}}=0$. Since
$e_{v}(c)\geq r_{v}+d_{v}$,
\begin{equation}\label{eq:maximal-opposite-finite-vanishing}
 r_{v}>\max\{1,r_{\chi_{v}}\}\quad\Longrightarrow\quad
 \mathcal{I}_{v}(s,\boldsymbol{\mu})=0.
\end{equation}
Otherwise $e_{v}(c)=r_{v}+d_{v}$ and
$\big|\int_{\mathcal{O}_{v}^{\times}}\chi_{v}(a)\psi_{v}(a/c)\,da\big|
\leq q_{v}^{-(r_{v}+d_{v})/2}$.
Using volume $q_{v}^{d_{v}/2}$ in $b$ and multiplicative volume
$q_{v}^{-d_{v}/2}$ in $c$, we obtain
\begin{align*}
 |\mathcal{I}_{v}(s,\mathbf{0})|\ll_{d_{v}}
 \Vol(K_{0,v}'[r_{v}])^{-1}q_{v}^{2|\Re(s)|(r_{v}+d_{v})}.
\end{align*}
If $r_{v}=0$, the characters are unramified and only
$e_{v}(c)=d_{v}$ or $d_{v}+1$ contributes; these places depend only on
$F$. Since only these finitely many valuations of $c$ contribute,
$\mathcal{I}_{v}(s,\boldsymbol{\mu})$ is holomorphic in $\boldsymbol{\mu}$
at every finite place.

For the Archimedean estimate, use
$x=\left(\begin{smallmatrix}a&(aa'-b')/c\\c&a'\end{smallmatrix}\right)$
in \eqref{eq:maximal-opposite-kernel}, so
$d^{\times} b'=\zeta_{v}(1)|c|_{v}|\det x|_{v}^{-1}\,dx_{12}$. On the support, $|x_{21}|_{v}\ll C_{v}^{-1}$, the entries
$x_{11},x_{12},x_{22}$ lie in fixed compact sets, and $\det x$ stays near $1$, where
$\chi_{v}$ is trivial. In $(h_{v}\ast h_{v}^{*})(x)=\int h_{v}(xg)\overline{h_{v}(g)}\,dg$,
$(xg)_{21}$ is independent of $x_{11}$, so derivatives in $x_{11}$
of $h_{v}(xg)$ fall on $h_{1,v}$. Since $\|h_{v}\|_{1}\ll1$ and $|\det x|_{v}\asymp1$,
\begin{align*}
 \left|\partial_{x_{11}}^{j}\big((h_{v}\ast h_{v}^{*})(x)|\det x|_{v}^{-1-\mu_{1}}\big)\right|
 \ll_{j} C_{v}
\end{align*}
for real derivatives in $x_{11}$ on the support, uniformly for
$\boldsymbol{\mu}$ in compact sets. Integrate by parts against
$\psi_{v}((a+a')/c)$ in $a=x_{11}$, $m$ times over $\mathbb{R}$ or $2m$
times over $\mathbb{C}$, to gain $|c|_{v}^{m}$. For $m\geq2$ and
$|\Re(s)|\leq2\delta$, we obtain
\begin{align*}
 |\mathcal{I}_{v}(s,\mathbf{0})|
 \ll_{m,\delta}C_{v}\int_{|c|_{v}\ll C_{v}^{-1}}
 |c|_{v}^{m+2\Re(s)-2}\,dc
 \ll_{m,\delta}C_{v}^{2-m-2\Re(s)}.
\end{align*}
For $\boldsymbol{\mu}$ near $\mathbf{0}$, the integral bound uses
$m+2\Re(s+\mu_{1})-2$ as the exponent of $|c|_{v}$. Its local uniform
convergence gives holomorphy in $\boldsymbol{\mu}$ for $|\Re(s)|\leq2\delta$.

For $v\notin S$, use the local factors in Proposition~\ref{prop6.4}.
Since $\widetilde{\tau}\simeq\tau$, the partial functional equation gives
\begin{align*}
 L^{S}(\mu_{1}-\mu_{2}/2,\tau)
 \prod_{v\in S}\gamma(\mu_{1}-\mu_{2}/2,\tau_{v},\overline{\psi}_{v})^{-1}
 =L^{S}(1-\mu_{1}+\mu_{2}/2,\tau).
\end{align*}
Multiply this identity by $L^{S}(1+\mu_{2}+2s,\chi)$ and
$\prod_{v\in S}\mathcal{I}_{v}(s,\boldsymbol{\mu})/\zeta_{v}(1)$
to assemble the local integrals in \eqref{eq:maximal-opposite-local-integral}.
For $\Re(s)$ sufficiently large the product is holomorphic near
$\boldsymbol{\mu}=\mathbf{0}$, so the constant term is its value there. Continuing in $s$ gives
\begin{equation}\label{eq:maximal-opposite-global-assembly}
 J_{\mathrm{degen}}^{\dag,\heartsuit}(s,\varphi,\mathbf{1})
 =L^{S}(1,\tau)L^{S}(1+2s,\chi)
 \prod_{v\in S}\frac{\mathcal{I}_{v}(s,\mathbf{0})}{\zeta_{v}(1)}.
\end{equation}

For $|s|=\rho\geq\delta>0$, convexity near $1$, together with
$C(\tau)C(\chi)=C(\pi)$, bounds the two $L$-factors.
For the local product, use $m=2$ at infinity,
$\Vol(K_{0,v}'[r_{v}])^{-1}\ll_{F} q_{v}^{r_{v}}(1+q_{v}^{-1})$ at finite
places, and $\prod_{v<\infty}q_{v}^{r_{v}}=Q$. Thus
\begin{align*}
|L^{S}(1,\tau)L^{S}(1+2s,\chi)|
 \ll_{F,\varepsilon,\delta}Q^{\varepsilon}C(\pi)^{\varepsilon},\quad \prod_{v\in S}\Big|\frac{\mathcal{I}_{v}(s,\mathbf{0})}{\zeta_{v}(1)}\Big|
 \ll_{F,\delta}Q^{1+O(\delta)}\mathbf{C}_{\infty}^{4\delta},
\end{align*}
Taking $\delta$ sufficiently small proves
\eqref{eq:maximal-opposite-degenerate-bound}. The vanishing assertion
follows from \eqref{eq:maximal-opposite-finite-vanishing} and
\eqref{eq:maximal-opposite-global-assembly}.
\end{proof}

On the circle $|s|=\rho\geq\delta>0$,
\begin{equation}\label{eq:maximal-singular-bound}
 \sup_{|s|=\rho}|J_{\mathrm{sing}}^{\heartsuit}(s,\varphi,\mathbf{1})|
 \ll_{F,\varepsilon}C(\pi)^{\varepsilon}Q^{1+\varepsilon}
 \mathbf{C}_{\infty}^{1+\varepsilon}.
\end{equation}
Use \eqref{eq:singular-finite-exact} for $v\mid\mathfrak{q}\mathfrak{m}\mathfrak{D}_{F}$
and $J_{\mathrm{sing},v}(s)\ll C_{v}$ at the Archimedean places. The standard
bounds near $1$ on this circle applied to
the two factors give
$|L^{S}(1+2s,\tau)L^{S}(1+2s,\chi)|\ll_{F,\varepsilon}
(C(\tau)C(\chi))^{\varepsilon}=C(\pi)^{\varepsilon}$.

\subsubsection{The dual integral}
For a nonzero integral ideal $\mathfrak{m}$ and $T_{v}\geq1$, put
$M=N_{F}(\mathfrak{m})$, $T=\prod_{v\mid\infty}T_{v}$, and
\begin{align*}
\mathfrak{X}
 :=\{(\xi,\lambda)\in\widehat{C_{F}}\times i\mathbb{R}:\ 
    \mathfrak{c}(\xi)=\mathfrak{m},\quad
    T_{v}\leq C_{v}(\xi|\cdot|^{\lambda})<2T_{v}
    \text{ for all }v\mid\infty\}.
\end{align*}
We use counting measure in $\xi$ and Lebesgue measure on $i\mathbb{R}$.

\begin{lemma}\label{lem:dual-character-average}
Let $\Pi$ be unitary cuspidal on $\mathrm{GL}_{3}/F$, or
$\Pi=\tau\boxplus\chi$ with $\tau$ unitary cuspidal on
$\mathrm{GL}_{2}/F$ and $\chi$ a unitary Hecke character. Then
\begin{equation}\label{eq:dual-character-average}
 \int_{\mathfrak{X}}
 |L(1/2-\lambda,\Pi\otimes\overline{\xi})
   L(1/2+\lambda,\xi)|\,d\xi\,d\lambda 
 \ll_{F,\varepsilon}
 C(\Pi)^{1/4+\varepsilon}M^{5/4+\varepsilon}
 T^{5/4+\varepsilon}.
\end{equation}
\end{lemma}

\begin{proof}
The fixed-modulus hybrid large sieve
(cf.\ \cite[Proposition C.3]{Qi24}) gives
\begin{align*}
 \int_{\mathfrak{X}}
 \left|\sum_{N_{F}(\mathfrak{n})\leq Y}
       c(\mathfrak{n})\xi(\mathfrak{n})
       N_{F}(\mathfrak{n})^{\lambda}\right|^{2}\,d\xi\,d\lambda 
 \ll_{F,\varepsilon}(MT+Y)(MTY)^{\varepsilon}
 \sum_{N_{F}(\mathfrak{n})\leq Y}|c(\mathfrak{n})|^{2},
\end{align*}
for coefficients $c(\mathfrak{n})$ supported on ideals prime to $\mathfrak{m}$.

First suppose that $\Pi$ is cuspidal. By the standard Rankin--Selberg bound,
\begin{align*}
 \sum_{N_{F}(\mathfrak{n})\leq Y}
 \frac{|\lambda_{\Pi}(\mathfrak{n})|^{2}}{N_{F}(\mathfrak{n})}
 \ll_{F,\varepsilon}(C(\Pi)Y)^{\varepsilon}.
\end{align*}
The approximate functional equation
\cite[Theorem~5.3 and Proposition~5.4]{IK04} expresses the two
$L$-values as smooth sums of lengths
$\ll C(\Pi)^{1/2+\varepsilon}(MT)^{3/2+\varepsilon}$ and
$\ll(MT)^{1/2+\varepsilon}$, respectively, up to a negligible error.
After truncation, write each weight as a Mellin integral on
$\Re(w)=\delta>0$. By Stirling with a Gaussian auxiliary factor,
the Gamma and conductor factors are bounded by $(C(\Pi)MT)^{O(\delta)}$
times a rapidly decreasing function of $\Im(w)$, uniformly in $(\xi,\lambda)$.
Split each summation ideal as $\mathfrak{d}\mathfrak{b}$, with
$\mathfrak{d}$ supported on primes dividing $\mathfrak{m}$ and
$\mathfrak{b}$ coprime to $\mathfrak{m}$. The bounds towards Ramanujan
give $O_{\varepsilon}(M^{\varepsilon})$ for the sum over $\mathfrak{d}$
of the suprema of the absolute coefficients weighted by
$N_{F}(\mathfrak{d})^{-1/2-\delta}$.
The trivial character's polar terms have bounded $L^{2}$ norm by the Gaussian.
Apply Minkowski in $w,\mathfrak{d}$ and the large sieve to the polynomial
in $\mathfrak{b}$, whose coefficients other than the character are fixed.
With the Rankin--Selberg coefficient bound and small enough $\delta$, this gives
\begin{align*}
 &\int_{\mathfrak{X}}
  |L(1/2-\lambda,\Pi\otimes\overline{\xi})|^{2}\,d\xi\,d\lambda \ll_{F,\varepsilon}
 C(\Pi)^{1/2+\varepsilon}(MT)^{3/2+\varepsilon},\\
 &\int_{\mathfrak{X}}
 |L(1/2+\lambda,\xi)|^{2}\,d\xi\,d\lambda
 \ll_{F,\varepsilon}(MT)^{1+\varepsilon}.
\end{align*}
Cauchy--Schwarz proves
\eqref{eq:dual-character-average} in the cuspidal case.

Now let $\Pi=\tau\boxplus\chi$. The approximate functional equation
for $L(1/2-\lambda,\tau\otimes\overline{\xi})$ gives sums of length
$\ll C(\tau)^{1/2+\varepsilon}(MT)^{1+\varepsilon}$.
The same separation of the Mellin weights and the large sieve give
\begin{align*}
 \int_{\mathfrak{X}}
 |L(1/2-\lambda,\tau\otimes\overline{\xi})|^{2}\,d\xi\,d\lambda
 \ll_{F,\varepsilon}C(\tau)^{1/2+\varepsilon}(MT)^{1+\varepsilon}.
\end{align*}
For the fourth power of $L(1/2-\lambda,\chi\overline{\xi})
L(1/2+\lambda,\xi)$, \cite[Theorem~F]{Yan25b} bounds the integral
with the adjoint denominator at level $\mathfrak{c}(\chi)\mathfrak{m}^{2}$
and $C_{v}(\sigma_{v})\ll_{F} C_{v}(\chi)T_{v}^{2}$.
The factor $L(1+2\lambda,\xi^{2}\overline{\chi})$ in the denominator
can have a pole only when $\xi^{2}=\chi$ on $C_{F}$. In that case, write
$\chi=\xi^{2}|\cdot|^{2\lambda_{0}}$, $\lambda_{0}\in i\mathbb{R}$, so that
\begin{align*}
 L(1+2\lambda,\xi^{2}\overline{\chi})
 =\zeta_{F}(1+2\lambda-2\lambda_{0}).
\end{align*}
There are $O_{F,\varepsilon}(M^{\varepsilon})$ such characters.
By convexity, their total contribution to the integral of
$|L(1/2-\lambda,\chi\overline{\xi})L(1/2+\lambda,\xi)|^{4}$
over $|\lambda-\lambda_{0}|\leq1$ is
$O_{F,\varepsilon}((MT)^{2+\varepsilon})$.
On the rest of $\mathfrak{X}$ the bound near $1$ gives
$|L(1+2\lambda,\xi^{2}\overline{\chi})|^{2}
\ll_{F,\varepsilon}(C(\chi)MT)^{\varepsilon}$. Thus
\begin{align*}
 \int_{\mathfrak{X}}
 |L(1/2-\lambda,\chi\overline{\xi})
   L(1/2+\lambda,\xi)|^{4}\,d\xi\,d\lambda 
 \ll_{F,\varepsilon}C(\chi)^{1+\varepsilon}(MT)^{2+\varepsilon}.
\end{align*}
The family has measure $O_{F}(MT)$. H\"older's inequality with
exponents $2,4,4$ and $C(\Pi)=C(\tau)C(\chi)$ now prove
\eqref{eq:dual-character-average}.
\end{proof}

\begin{cor}\label{cor:trivial-central-character-dual-moment}
For $\pi$ as in Theorem~\ref{thm1.3} and any integral ideal
$\mathfrak{q}$,
\begin{align*}
 |J_{\mathrm{dual}}^{\heartsuit}(\varphi,\mathbf{1})|
 \ll_{F,\varepsilon}C(\pi)^{\varepsilon}
 Q^{1/2+\varepsilon}\mathbf{C}_{\infty}^{1/2+\varepsilon}
 C(\pi)^{1/4+\varepsilon}
 \prod_{v\mid\infty}[1+(C_{v}^{-1}C_{v}(\pi))^{2}].
\end{align*}
\end{cor}

\begin{proof}
Apply Lemma~\ref{lem:dual-character-average} to \eqref{13.7}
separately for each finite conductor and each range
$T_{v}\leq C_{v}(\xi|\cdot|^{\lambda})<2T_{v}$, where $T_{v}=1,2,4,\ldots$
at every Archimedean place. The decay in
Theorem~\ref{thm9.5}, with $n=2$, gives the convergent sums
$\sum_{T_{v}\geq1\ \mathrm{dyadic}}T_{v}^{-3/4+\varepsilon}$
at infinity. At each finite place, we must bound
$q_{v}^{(5/4+\varepsilon)r_{\xi_{v}}}
|\mathcal{H}_{W_{v}}^{\sharp}|$ for the allowed conductor exponents
$r_{\xi_{v}}$.

For $d_{v}=0$, Proposition~\ref{lem11.6} gives
$|\mathcal{H}_{W_{v}}^{\sharp}|\ll q_{v}^{r_{v}/2}$ for unramified
$\xi_{v}$. If $L_{v}(s,\pi_{v})$ has degree zero or one, the weight vanishes
for ramified $\xi_{v}$ by Lemma~\ref{lem:finite-projector-full-Tate-tail} and
\eqref{eq:maximal-degree-one-compact-lift}.
If its degree is two or three,
Lemma~\ref{lem:maximal-Steinberg-unitary-dual} and
\eqref{eq:spherical-projector-ramified-bound} allow a ramified
contribution only when $r_{v}=r_{\xi_{v}}=1$, with size $O(q_{v}^{-1})$.
Thus the ramified factor is
$O(q_{v}^{5/4+\varepsilon}q_{v}^{-1})
\ll q_{v}^{r_{v}/2+\varepsilon}$.
At $v\mid\mathfrak{D}_{F}$ with $r_{v}\geq1$ the same estimates apply
with the factor $q_{v}^{d_{v}/2}$ from
\eqref{eq:maximal-positive-residue-different-transport}.
For $r_{v}=0$, \eqref{eq:finite-zero-level-dual-bound} and
\eqref{eq:finite-dual-character-cutoff} bound the weighted local factor
by $O_{F}(1)$, since $r_{\xi_{v}}\leq3$ and $v\mid\mathfrak{D}_{F}$.
Summing the products of the local bounds over the
$O_{F,\varepsilon}(Q^{\varepsilon})$ supported finite conductors
gives $O_{F,\varepsilon}(Q^{1/2+\varepsilon})$.
Use $\widetilde{\pi}\simeq\pi$ in Lemma~\ref{lem:dual-character-average}.
\end{proof}

\subsubsection{The spectral residues}
Let $F_{\chi}(s,\lambda)$ be the basis sum defined in
\textsection\ref{sec4.1}. Define 
\begin{multline}\label{eq:maximal-spectral-density-factors}
\mathscr{P}_{\chi}(s,\lambda):=F_{\chi}(s,\lambda)
\frac{\Lambda(1+2\lambda,\mathbf{1})\Lambda(1-2\lambda,\mathbf{1})}
{\Lambda(1/2+s+\lambda,\mathbf{1})\Lambda(1/2+s-\lambda,\mathbf{1})}\\
\frac{1}{\Lambda(1/2+s+\lambda,\tau)\Lambda(1/2+s-\lambda,\tau)}.
\end{multline}

\begin{lemma}
We have $\mathscr{P}_{\chi}(s,\lambda)=\rho_{F}^{-1}\prod_{v\in S}\mathscr{P}_{\chi,v}(s,\lambda)$,
where $\mathscr{P}_{\chi,v}(s,\lambda)$ is defined by the local version of
\eqref{eq:maximal-spectral-density-factors} and equals $1$ at $v\notin S$.
As a meromorphic identity in $s$,
\begin{equation}\label{eq:maximal-positive-spectral-residue}
R_{\RNum{1}}^{+}(s)
 =\frac{\Res_{z=1}\Lambda(z,\mathbf{1})}{2}
\cdot  \frac{\Lambda(1,\tau)\Lambda(2s,\tau)}
        {\Lambda(2-2s,\mathbf{1})}
\cdot \mathscr{P}_{\chi}(s,1/2-s).
\end{equation}
The finite global factor satisfies, for $|s|\leq\delta$ and $\delta$ sufficiently small
in terms of $\varepsilon$,
\begin{equation}\label{eq:maximal-spectral-residue-finite-bound}
 \left|\frac{L(1,\tau)L(2s,\tau)}
 {L(2-2s,\mathbf{1})}\right|
 \ll_{F,\varepsilon}C(\tau)^{1/2+\varepsilon}.
\end{equation}
\end{lemma}

\begin{proof}
The product formula follows from \eqref{eq:initial-data-typeI-density},
with the factor $\rho_{F}^{-1}$ given by \eqref{eq:inducing-global-norm}.
Since $\pi\times\chi\simeq\tau\boxplus\mathbf{1}$,
\eqref{eq:initial-data-typeI-residue} gives
\eqref{eq:maximal-positive-spectral-residue}.
Finally, $\Re(2-2s)>1$ gives a uniform bound for
$L(2-2s,\mathbf{1})^{-1}$. The standard upper bound for $L(1,\tau)$
and convexity near zero for the cuspidal degree-two factor $L(2s,\tau)$
give \eqref{eq:maximal-spectral-residue-finite-bound}.
\end{proof}

At any place $v$, for a smooth vector $h$ in the compact model of
$\chi_{v}|\cdot|_{v}^{\lambda}\boxplus\chi_{v}|\cdot|_{v}^{-\lambda}$,
let $f_{\lambda,h}$ be the section with $f_{\lambda,h}|_{K_{v}'}=h$.
We use
\begin{equation}\label{eq:maximal-GL2-Jacquet}
 W_{\lambda,h}(g):=\int_{F_{v}}
 f_{\lambda,h}(w'n(x)g)\overline{\psi_{v}(x)}\,dx,
\end{equation}
initially for $\Re(\lambda)>0$, and then by continuation in $\lambda$.

\begin{lemma}
Let $v\mid\mathfrak{q}\mathfrak{m}\mathfrak{D}_{F}$ with $r_{v}\geq1$.
Let $\mathfrak{B}_{v}$ be an orthonormal basis in the compact model
for the norm $\int_{K_{v}'}|h(k)|^{2}\,dk$, and use
$W_{\lambda,h}$ from \eqref{eq:maximal-GL2-Jacquet}. Then
\begin{equation}\label{eq:maximal-finite-eisenstein-weight}
\sum_{h\in\mathfrak{B}_{v}}
 \int_{G'(F_{v})}h_{v}(g)W_{\lambda,h}(g)\,dg
\overline{\int_{G'(F_{v})}
 h_{v}(g)W_{-\overline{\lambda},h}(g)\,dg}
 \ll \mathbf{1}_{2r_{\chi_{v}}\leq r_{v}}.
\end{equation}
\end{lemma}

\begin{proof}
First take $d_{v}=0$. Let $P_{r_{v}}$ be the
$K_{0,v}'[r_{v}]$-projector. Its action
in the compact model is independent of $\lambda$, and
\begin{align*}
 \int_{G'(F_{v})}h_{v}(g)W_{\lambda,h}(g)\,dg=W_{\lambda,P_{r_{v}}h}(I_{2}).
\end{align*}
It is an orthogonal projector for the norm in the compact model. If its
image is nonzero, the conductor formula for a principal series gives
$2r_{\chi_{v}}\leq r_{v}$.

Decompose the Jacquet integral into $x\in\mathcal{O}_{v}$ and
$e_{v}(x)=-j$, $j\geq1$, using
\begin{align*}
 \begin{pmatrix}&1\\1&x\end{pmatrix}
 =\begin{pmatrix}-x^{-1}&1\\&x\end{pmatrix}
  \begin{pmatrix}1&\\x^{-1}&1\end{pmatrix}.
\end{align*}
The integral over $\mathcal{O}_{v}^{\times}$ is a Gauss integral for $\chi_{v}$.
If $r_{\chi_{v}}>0$, only $j=r_{\chi_{v}}\leq r_{v}$ contributes,
with absolute value $q_{v}^{r_{\chi_{v}}/2}$; the $x\in\mathcal{O}_{v}$ term is zero.
If $r_{\chi_{v}}=0$, only $x\in\mathcal{O}_{v}$ and $j=1$ contribute, with Gauss
integrals $1$ and $-1$. The inducing factor on the $j$-th set is
$q_{v}^{-j(1+2\lambda)}$.

The quotient measure is $dx/[(1+q_{v}^{-1})\max\{1,|x|_{v}\}^{2}]$.
Choose an orthonormal basis of the $K_{0,v}'[r_{v}]$-fixed space with each
vector supported on $\mathcal{O}_{v}$, on
$\{x:e_{v}(x)\leq-r_{v}\}\cup\{\infty\}$, or on a single
$\varpi_{v}^{-j}\mathcal{O}_{v}^{\times}$ with $1\leq j<r_{v}$.
The measures of these supports are
\begin{align*}
\mu_{0}=\frac{q_{v}}{q_{v}+1},\qquad  
\mu_{r_{v}}=\frac{q_{v}}{q_{v}+1}\cdot q_{v}^{-r_{v}},\qquad
\mu_{j}=\frac{q_{v}-1}{q_{v}+1}\cdot q_{v}^{-j},\quad 1\leq j<r_{v}.
\end{align*}

Each basis vector has constant absolute value on its support, equal to
$\mu_{j}^{-1/2}$ when the support has measure $\mu_{j}$;
disjoint supports give orthogonality. For $r_{\chi_{v}}>0$, only the vector
whose support meets $e_{v}(x)=-r_{\chi_{v}}$ contributes to
\eqref{eq:maximal-finite-eisenstein-weight}. The two inducing factors
have product $q_{v}^{-2r_{\chi_{v}}}$, so the sum equals
$q_{v}^{-r_{\chi_{v}}}/\mu_{r_{\chi_{v}}}\leq3$.
For $r_{\chi_{v}}=0$, \eqref{eq:maximal-finite-eisenstein-weight}
reduces to $\mu_{0}^{-1}+q_{v}^{-2}\mu_{1}^{-1}\leq3$.

For $d_{v}>0$, \eqref{eq:different-GL2-Jacquet} gives the factor
$\chi_{v}(\varpi_{v}^{-d_{v}})q_{v}^{-d_{v}\lambda}$. It cancels with the conjugate of
the factor at $-\overline{\lambda}$, since $|\chi_{v}(\varpi_{v}^{-d_{v}})|=1$.
\end{proof}
At $v\mid\mathfrak{q}\mathfrak{m}\mathfrak{D}_{F}$ with $r_{v}\geq1$,
the local factor in \eqref{eq:maximal-positive-spectral-residue}
is given in the basis $\mathfrak{B}_{v}$ by
\begin{multline}\label{eq:maximal-finite-residue-normalization}
\mathscr{P}_{\chi,v}(s,1/2-s)
 =\frac{L_{v}(2-2s,\mathbf{1})}{\zeta_{v}(1)L_{v}(1,\tau_{v})L_{v}(2s,\tau_{v})}\\
 \sum_{h\in\mathfrak{B}_{v}}
\int_{G'(F_{v})}h_{v}(g)W_{1/2-s,h}(g)\,dg\overline{\int_{G'(F_{v})}
 h_{v}(g)W_{\overline{s}-1/2,h}(g)\,dg}.
\end{multline}
By \eqref{eq:finite-projector-restriction}, $W_{v}(\diag(x,1))=0$
unless $|\det x|_{v}=1$, so the finite weight is unchanged by the $s$-twist.
Substitute $\lambda=1/2-s$ in \eqref{eq:maximal-spectral-density-factors}
and cancel the local $L$-factors to obtain
\eqref{eq:maximal-finite-residue-normalization}.

Let $\vartheta<1/2$ be an admissible finite-place Ramanujan exponent
for the cuspidal $\mathrm{GL}_{2}$ data. Write
$L_{v}(s,\tau_{v})=\prod_{j}(1-\alpha_{v,j}q_{v}^{-s})^{-1}$ and put
\begin{equation}\label{eq:maximal-local-exponent}
 \theta_{v}=\max\bigl(0,\max_{j}\log_{q_{v}}|\alpha_{v,j}|\bigr),
\end{equation}
with $\theta_{v}=0$ if the local factor is one. Using
$|\alpha_{v,1}\alpha_{v,2}|=1$ for unramified $\tau_{v}$ and
$|\alpha_{v,j}|\leq1$ for ramified $\tau_{v}$ gives
$|L_{v}(2s,\tau_{v})^{-1}|\ll q_{v}^{\theta_{v}+4\delta}$ for $|s|\leq\delta$.
Also $|L_{v}(2-2s,\mathbf{1})/(\zeta_{v}(1)L_{v}(1,\tau_{v}))|\ll1$, and
\eqref{eq:maximal-finite-eisenstein-weight} bounds the basis sum in
\eqref{eq:maximal-finite-residue-normalization} by $O(1)$. Hence
\begin{multline}\label{eq:maximal-finite-residue-bound}
 \prod_{v\in S,\ v<\infty}
 |\mathscr{P}_{\chi,v}(s,1/2-s)|
 \ll_{F,\varepsilon}Q^{\varepsilon}
 \prod_{v\in S,\ v<\infty}q_{v}^{\theta_{v}}\\
 \ll_{F,\varepsilon}Q^{\varepsilon}
 \min\{Q^{\vartheta},Q^{1/2}C_{\mathrm{fin}}(\tau)^{-1/4},
                 Q C_{\mathrm{fin}}(\tau)^{-1/2}\}.
\end{multline}
Here $\delta$ is sufficiently small in terms of $\varepsilon$.
Since $r_{\tau_{v}}\leq r_{\pi_{v}}\leq r_{v}$, and $\theta_{v}=0$ when
$\tau_{v}$ is ramified, we have, for $r_{v}\geq1$,
\begin{align*}
 \theta_{v}\leq\vartheta r_{v},\qquad
 r_{\tau_{v}}/2+\theta_{v}\leq r_{v},\qquad
 r_{\tau_{v}}+4\theta_{v}\leq2r_{v}.
\end{align*}
Multiplication gives all three estimates in
\eqref{eq:maximal-finite-residue-bound}, including
\begin{equation}\label{eq:maximal-spectral-residue-conductor-comparison}
 C_{\mathrm{fin}}(\tau)^{1/2}
 \prod_{v\in S,\ v<\infty}q_{v}^{\theta_{v}}
 \ll_{F} Q^{1/2}C_{\mathrm{fin}}(\tau)^{1/4}.
\end{equation}
The places with $r_{v}=0$ divide the fixed different and contribute
only to the constant.

The vanishing assertion following \eqref{eq:maximal-finite-eisenstein-weight}
implies $R_{\RNum{1}}^{\pm}=0$ if $2r_{\chi_{v}}>r_{v}$ at any
$v\mid\mathfrak{q}\mathfrak{m}\mathfrak{D}_{F}$.

\begin{lemma}
Let $v\mid\infty$.
For an orthonormal basis in the compact model, use $W_{\lambda,h}$
from \eqref{eq:maximal-GL2-Jacquet} to define
\begin{multline*}
\mathcal{K}_{v}(s,\lambda):=\sum_{h}
\int_{G'(F_{v})}h_{1,v}(g)h_{2,v}(g;C_{v})
 W_{\lambda,h}(g)|\det g|_{v}^{s}\,dg\\
\overline{\int_{G'(F_{v})}h_{1,v}(g)h_{2,v}(g;C_{v})W_{-\overline{\lambda},h}(g)|\det g|_{v}^{-\overline{s}}\,dg}.
\end{multline*}
For $|s|\leq\delta$ and $|\lambda-1/2|\leq\delta$, with
$0<\delta<10^{-2}$, $\mathcal{K}_{v}$ is holomorphic and satisfies
\begin{equation}\label{eq:maximal-arch-eisenstein-weight}
 |\mathcal{K}_{v}(s,\lambda)|\ll_{\delta,h_{v}}1.
\end{equation}
\end{lemma}

\begin{proof}
Put $d=[F_{v}:\mathbb{R}]$. Write $g=\diag(u,u^{-1})g_{0}$ and use quotient
measure for left multiplication by $\diag(u,u^{-1})$. On the common support of $h_{1,v}$ and $h_{2,v}$, $u$ stays in a
fixed neighborhood of $1$, while the quotient support has measure
$O(C_{v}^{-1})$. Every fixed-order $u$-derivative of
their product is $O(C_{v})$: differentiating $u^{-1}(g_{0})_{21}$
compensates for the scaling in $h_{2,v}$.

By normalized induction, the change $x\mapsto u^{2}y$ in
\eqref{eq:maximal-GL2-Jacquet} gives
\begin{align*}
 W_{\lambda,h}(\diag(u,u^{-1})g_{0})
 =|u|_{v}^{1-2\lambda}
 \int_{F_{v}}f_{\lambda,h}(w'n(y)g_{0})\psi_{v}(-u^{2}y)\,dy.
\end{align*}
The $u$-integration therefore gives the factor
\begin{align*}
 F_{g_{0},\lambda}(y)=\int_{F_{v}^{\times}}
 h_{1,v}(\diag(u,u^{-1})g_{0})h_{2,v}(\diag(u,u^{-1})g_{0};C_{v})|u|_{v}^{1-2\lambda}
 \psi_{v}(-u^{2}y)\,d^{\times}u.
\end{align*}
The map $u\mapsto u^{2}$ is a diffeomorphism on this neighborhood.
Hence $F_{g_{0},\lambda}$ is the Fourier transform of a compactly
supported smooth function. Plancherel and the derivative bounds give
\begin{align*}
 \left\|(1+|y|_{\mathbb{R}}^{2})^{l/2}
 F_{g_{0},\lambda}(y)\right\|_{L^{2}(dy)}
 \ll_{l,\delta}C_{v},\quad l\in\mathbb{Z}_{\geq 0},
\end{align*}
uniformly in $g_{0}$ for $|\lambda-1/2|\leq\delta$ or
$|\lambda+1/2|\leq\delta$.
On the support, $|\det g_{0}|_{v}=|\det g|_{v}\asymp1$, so
$\big||\det g_{0}|_{v}^{s}\big|\ll_{\delta}1$ for $|s|\leq\delta$.

Iwasawa decomposition and normalized induction give
\begin{align*}
f_{\lambda,h}(w'n(y))=(1+|y|_{\mathbb{R}}^{2})^{-d\lambda}f_{0,h}(w'n(y))
\end{align*}
and
\begin{align*}
\int_{K_{v}'}|h(k)|^{2}\,dk
 =\frac{1}{d\pi}\int_{F_{v}}|f_{0,h}(w'n(y))|^{2}\,dy.
\end{align*}
Here $|\cdot|_{\mathbb{R}}$ denotes Euclidean modulus.
Thus $(d\pi)^{-1/2}f_{0,h}(w'n(\cdot))$ form an orthonormal basis
of $L^{2}(F_{v},dy)$, and Parseval's identity gives
\begin{align*}
 \sum_{h}\left|\int_{F_{v}}
 f_{\lambda,h}(w'n(y))F_{g_{0},\lambda}(y)\,dy\right|^{2}=d\pi\int_{F_{v}}|F_{g_{0},\lambda}(y)|^{2}
 (1+|y|_{\mathbb{R}}^{2})^{-2d\Re(\lambda)}\,dy.
\end{align*}
For $|\lambda-1/2|\leq\delta$ or $|\lambda+1/2|\leq\delta$,
take an integer $l\geq d(1+2\delta)$. Since $-2d\Re(\lambda)\leq l$,
the weighted $L^{2}$ bound gives
\begin{align*}
\sum_{h}\left|\int_{F_{v}}
 f_{\lambda,h}(w'n(y))F_{g_{0},\lambda}(y)\,dy\right|^{2}\leq d\pi\left\|(1+|y|_{\mathbb{R}}^{2})^{l/2}
 F_{g_{0},\lambda}(y)\right\|_{L^{2}(dy)}^{2}
 \ll_{l,\delta}C_{v}^{2}.
\end{align*}
Right translation by $g_{0}$ has bounded operator norm in the compact model,
independently of the unitary characters, since $g_{0}$ ranges over a
fixed compact subset of $G'(F_{v})$. Integrating over the quotient
support therefore gives
\begin{align*}
 &\sum_{h}\left|\int_{G'(F_{v})}h_{1,v}(g)h_{2,v}(g;C_{v})
 W_{\lambda,h}(g)|\det g|_{v}^{s}\,dg\right|^{2}\ll1,\\
 &\sum_{h}\left|\int_{G'(F_{v})}h_{1,v}(g)h_{2,v}(g;C_{v})
 W_{-\overline{\lambda},h}(g)|\det g|_{v}^{-\overline{s}}\,dg\right|^{2}
 \ll1.
\end{align*}
Cauchy--Schwarz proves \eqref{eq:maximal-arch-eisenstein-weight}.
Both square sums are locally $O(1)$ after parameter differentiation,
so the series defining $\mathcal{K}_{v}$ is holomorphic.
\end{proof}

\begin{cor}
For $|s|=\rho$,
\begin{equation}\label{eq:maximal-spectral-residue-geometric-bound}
 |R_{\RNum{1}}^{\pm}(s)|\ll_{F,\varepsilon}
 C(\pi)^{\varepsilon} Q^{1/2+\varepsilon}
 \mathbf{C}_{\infty}^{3/4+\varepsilon}
 C(\pi)^{1/4+\varepsilon}
 \prod_{v\mid\infty}(1+C_{v}^{-1}C_{v}(\pi)).
\end{equation}
Consequently,
\begin{equation}\label{eq:maximal-spectral-residue-absorbed-bound}
|R_{\RNum{1}}^{\pm}(s)|\ll_{F,\varepsilon}C(\pi)^{\varepsilon}Q^{\varepsilon}\mathbf{C}_{\infty}^{\varepsilon}
\left[Q\mathbf{C}_{\infty}+\mathbf{C}_{\infty}^{1/2}C(\pi)^{1/2}\right]\prod_{v\mid\infty}(1+C_{v}^{-1}C_{v}(\pi)).
\end{equation}
\end{cor}
\begin{proof}
The Archimedean version of
\eqref{eq:maximal-finite-residue-normalization} is
\begin{align*}
 \mathscr{P}_{\chi,v}(s,1/2-s)
 =\mathcal{K}_{v}(s,1/2-s)\cdot 
 \frac{L_{v}(2-2s,\mathbf{1})}
{\zeta_{v}(1)L_{v}(1,\tau_{v})L_{v}(2s,\tau_{v})}.
\end{align*}
Substitute this expression for $\mathscr{P}_{\chi,v}(s,1/2-s)$
in \eqref{eq:maximal-positive-spectral-residue} before taking
absolute values. The $L_{v}$-factors cancel, leaving
$\zeta_{v}(1)^{-1}\mathcal{K}_{v}$. 

The bounds \eqref{eq:maximal-spectral-residue-finite-bound},
\eqref{eq:maximal-spectral-residue-conductor-comparison}, and
\eqref{eq:maximal-arch-eisenstein-weight} give
\begin{align*}
 |R_{\RNum{1}}^{\pm}(s)|
 \ll_{F,\varepsilon}C(\pi)^{\varepsilon}Q^{1/2+\varepsilon}
 C_{\mathrm{fin}}(\tau)^{1/4}C_{\infty}(\tau)^{1/2}.
\end{align*}
Since $C_{\mathrm{fin}}(\tau)\leq C_{\mathrm{fin}}(\pi)$ and
\begin{align*}
 C_{v}(\tau)^{1/2}
 \leq C_{v}^{3/4}C_{v}(\pi)^{1/4}(1+C_{v}(\pi)/C_{v}),
\end{align*}
we obtain \eqref{eq:maximal-spectral-residue-geometric-bound}.
Apply $2\sqrt{ab}\leq a+b$ to $a=Q\mathbf{C}_{\infty}$ and
$b=\mathbf{C}_{\infty}^{1/2}C(\pi)^{1/2}$ to obtain
\eqref{eq:maximal-spectral-residue-absorbed-bound}.
\end{proof}

\subsubsection{The dual residues and their local normalization}
In \eqref{eq:initial-data-negative-dual-residue} and
\eqref{eq:initial-data-positive-dual-residue}, use the finite
$L$-factor normalization of \eqref{13.7}. Then
\begin{multline}\label{eq:maximal-negative-dual-residue}
 R_{\mathrm{dual}}^{-}(s,\varphi,\mathbf{1})
 =-\rho_{F}^{-1}\underset{z=1}{\Res}\zeta_{F}(z)\,
 L(1,\tau)L(2s,\chi)\\
 \prod_{v\mid\infty}
 \mathcal{H}_{W_{v}}(\chi_{v},2s-1/2;\mathbf{1}_{v},s)
 \prod_{\substack{v\in S\\v<\infty}}
 \mathcal{H}_{W_{v}}^{\sharp}(\chi_{v},2s-1/2;\mathbf{1}_{v},s),
\end{multline}
and
\begin{multline}\label{eq:maximal-positive-dual-residue}
 R_{\mathrm{dual}}^{+}(s,\varphi,\mathbf{1})
 =\frac{\underset{z=1}{\Res}\zeta_{F}(z)}{\rho_{F}}L(2s,\pi)
 \prod_{v\mid\infty}
 \mathcal{H}_{W_{v}}(\mathbf{1}_{v},1/2;\mathbf{1}_{v},s)\\
 \prod_{\substack{v\in S\\v<\infty}}
 \mathcal{H}_{W_{v}}^{\sharp}(\mathbf{1}_{v},1/2;\mathbf{1}_{v},s).
\end{multline}

Apply the functional equation to $L(2s,\pi)$ before estimating
the local transforms:
\begin{equation}\label{eq:maximal-residue-functional-equation}
 L(2s,\pi)=L(1-2s,\widetilde{\pi})
 \prod_{\substack{v\in S\\v<\infty}}
 \epsilon(2s,\pi_{v},\psi_{v})
 \prod_{v\mid\infty}\gamma(2s,\pi_{v},\psi_{v}).
\end{equation}
Use \eqref{eq:maximal-residue-functional-equation} with $\pi=\chi$
for $L(2s,\chi)$. For $|s|=\rho$, with the radius
$\delta\leq\rho\leq2\delta$ chosen above, the standard bounds near $1$ give
\begin{align*}
 L(1-2s,\widetilde{\pi})\ll_{F,\varepsilon}C(\pi)^{\varepsilon},
 \qquad L(1-2s,\chi)\ll_{F,\varepsilon}C(\chi)^{\varepsilon}
\end{align*}
when $\delta$ is sufficiently small.
For $v\mid\mathfrak{q}\mathfrak{m}\mathfrak{D}_{F}$, if
$r_{\chi_{v}}>3$, then
\eqref{eq:finite-dual-character-cutoff} makes the local transform
identically zero, and hence $R_{\mathrm{dual}}^{-}$ vanishes.

At places dividing the different, use
\eqref{eq:maximal-positive-residue-different-transport} and
\eqref{eq:different-epsilon-factor}. The powers of $q_{v}$ in these
formulas are bounded in terms of $F$ for $|s|\leq2\delta$.

For $d_{v}=0$ and $r_{v}\geq r_{\pi_{v}}+3$,
Lemma~\ref{lem:finite-projector-full-Tate-tail} gives
\begin{align*}
 &\mathcal{H}_{W_{v}}^{\sharp}(\chi_{v},2s-1/2;\mathbf{1}_{v},s)
 =\frac{\chi_{v}(\varpi_{v})^{r_{v}}q_{v}^{-2sr_{v}}
            \mathbf{1}_{r_{\chi_{v}}=0}}
 {\Vol(K_{0,v}'[r_{v}])L_{v}(1,\tau_{v})\zeta_{v}(1)},\\
 &\mathcal{H}_{W_{v}}^{\sharp}(\mathbf{1}_{v},1/2;\mathbf{1}_{v},s)
 =\frac{q_{v}^{-r_{v}}}
 {\Vol(K_{0,v}'[r_{v}])L_{v}(2s,\pi_{v})}.
\end{align*}

\subsubsection{Finite projector identities}
We allow arbitrary unitary local data in the following lemmas.
We keep the coordinates $g(z,y)$ of Lemma~\ref{lem:finite-dual-support}.

\begin{lemma}\label{lem:maximal-degree-two-compact-lift}
Suppose $d_{v}=0$, $\pi_{v}=\tau_{v}\boxplus\chi_{v}$,
$\deg L_{v}(s,\pi_{v})=2$, and $r_{\pi_{v}}\geq1$.
For $V_{v}$ from Lemma~\ref{lem:finite-projector-depth},
\begin{equation}\label{eq:maximal-degree-two-compact-lift}
 V_{v}(g(z,y))=\mathbf{1}_{\mathcal{O}_{v}^{\times}}(z),
 \qquad e_{v}(y)\geq r_{\pi_{v}}+1.
\end{equation}
\end{lemma}
\begin{proof}
Write $L_{v}(s,\pi_{v})^{-1}=1-\beta_{1}q_{v}^{-s}+\beta_{2}q_{v}^{-2s}$.
With $F_{i,j}$ from the proof of Lemma~\ref{lem:spherical-projector-dual-mass},
the same expansion with two Satake parameters reads
\begin{multline*}
 V_{v}(g(z,y))=F_{0,0}-\beta_{1}F_{0,1}+\beta_{2}F_{0,2}
 +q_{v}^{-1}(\beta_{1}^{2}-\beta_{2})F_{1,1}\\
 -q_{v}^{-1}\beta_{1}\beta_{2}F_{1,2}+q_{v}^{-2}\beta_{2}^{2}F_{2,2}.
\end{multline*}
This expansion uses only the upper-left $K_{v}'$-invariance of $V_{v}^{\mathrm{new}}$.

Put $K_{1,v}(r_{\pi_{v}}):=\{k\in G(\mathcal{O}_{v}):\ k_{31},k_{32},k_{33}-1\in\mathfrak{p}_{v}^{r_{\pi_{v}}}\}$. Let $\mathcal{T}_{1,v}$, $\mathcal{T}_{2,v}$ be the $K_{1,v}(r_{\pi_{v}})$ Hecke
operators at $\diag(\varpi_{v},1,1)$ and
$\diag(\varpi_{v},\varpi_{v},1)$, normalized by $\Vol(K_{1,v}(r_{\pi_{v}}))^{-1}$.
By \cite[Lemma~3.1]{Miy14}, the first two families below represent the
right cosets for $\mathcal{T}_{1,v}$ and the last is for $\mathcal{T}_{2,v}$.
\begin{align*}
 \begin{pmatrix}\varpi_{v}&u&w\\&1&\\&&1\end{pmatrix},\qquad
 \begin{pmatrix}1&&\\&\varpi_{v}&w\\&&1\end{pmatrix},
 \qquad
 \begin{pmatrix}\varpi_{v}&&u\\&\varpi_{v}&w\\&&1\end{pmatrix},
 \qquad u,w\bmod\mathfrak{p}_{v}.
\end{align*}
Evaluation at identity, using the diagonal newvector formula,
gives 
\begin{align*}
\mathcal{T}_{1,v}V_{v}^{\mathrm{new}}=q_{v}\beta_{1}V_{v}^{\mathrm{new}},\qquad \mathcal{T}_{2,v}V_{v}^{\mathrm{new}}=q_{v}\beta_{2}V_{v}^{\mathrm{new}}.
\end{align*}

For $e_{v}(y)\geq r_{\pi_{v}}+1\geq2$, use these relations
at $g(z,y)\diag(\varpi_{v}^{-1},\varpi_{v}^{-2},1)$ and
$g(z,y)\diag(\varpi_{v}^{-2},\varpi_{v}^{-2},1)$, respectively.
The first $q_{v}^{2}$ cosets of $\mathcal{T}_{1,v}$ contribute $q_{v}^{2}F_{0,2}$;
the remaining cosets agree with the $\mathcal{T}_{2,v}$ cosets after
summing their $n_{13}$-coordinate. Their Whittaker factors are one,
since $y/\varpi_{v}^{2}\in\mathcal{O}_{v}$ and the Whittaker character is
trivial on $n_{13}$. Eliminating this common sum gives
\begin{align*}
 q_{v}\beta_{1}F_{1,2}=q_{v}^{2}F_{0,2}+\beta_{2}F_{2,2}.
\end{align*}
Using $\mathcal{T}_{1,v}$ at $g(z,y)\diag(\varpi_{v}^{-1},\varpi_{v}^{-1},1)$ also gives
\begin{align*}
 F_{0,1}=q_{v}^{-1}\beta_{1}F_{1,1},\qquad
 V_{v}(g(z,y))=F_{0,0}-q_{v}^{-1}\beta_{2}F_{1,1}.
\end{align*}

For $j=0,1$, we have
\begin{align*}
 g(z,y)\diag(\varpi_{v}^{-j},\varpi_{v}^{-j},1)
 =\diag(z\varpi_{v}^{-j},z\varpi_{v}^{-j},1)n_{32}(-y\varpi_{v}^{-j}).
\end{align*}
Since $e_{v}(y)-j\geq r_{\pi_{v}}$, the last factor fixes $V_{v}^{\mathrm{new}}$.
The diagonal newvector formula \cite[Theorem~4.1]{Miy14} gives
\begin{align*}
F_{0,0}=q_{v}^{-e_{v}(z)}\beta_{2}^{e_{v}(z)}\mathbf{1}_{e_{v}(z)\geq0},\qquad
 F_{1,1}=q_{v}^{1-e_{v}(z)}\beta_{2}^{e_{v}(z)-1}\mathbf{1}_{e_{v}(z)\geq1}.
\end{align*}
Thus $V_{v}(g(z,y))=F_{0,0}-q_{v}^{-1}\beta_{2}F_{1,1}
=\mathbf{1}_{\mathcal{O}_{v}^{\times}}(z)$, proving
\eqref{eq:maximal-degree-two-compact-lift}.
\end{proof}

\begin{lemma}\label{lem:finite-projector-full-Tate-tail}
Suppose $d_{v}=0$ and put $\ell=\deg L_{v}(s,\pi_{v})$. If
\begin{align*}
 r_{v}\geq\max\{r_{\pi_{v}}+\ell,\,
 r_{\omega_{v}'}+1,\,
 r_{\omega_{v}'}+r_{\omega_{v}\omega_{v}'}\},
\end{align*}
then, meromorphically in $(s,\lambda)$,
\begin{equation}\label{eq:finite-projector-full-Tate-tail}
 \mathcal{H}_{W_{v}}^{\sharp}(\xi_{v},\lambda;\omega_{v}',s)
 =\frac{\xi_{v}(\varpi_{v})^{r_{v}} q_{v}^{-r_{v}(1/2+\lambda)}}
 {\Vol(K_{0,v}'[r_{v}])L_{v}(1/2+2s-\lambda,
       \widetilde{\pi}_{v}\overline{\xi}_{v}\omega_{v}\omega_{v}')}
 \mathbf{1}_{r_{\xi_{v}}=0}.
\end{equation}
For unramified $\pi_{v}$, the term $r_{\pi_{v}}+\ell$ in this maximum can
be replaced by $2$. For $\ell=1$ and $r_{\pi_{v}}\geq2$, it can be replaced by
$r_{\pi_{v}}$; for $2+1$ data with $\ell=2$, by $r_{\pi_{v}}+1$.
If $r_{\omega_{v}'}=0$ and $\ell=0$, the formula holds at every
allowed level $r_{v}\geq r_{\pi_{v}}$.
\end{lemma}

\begin{proof}
For $V_{v}$ from Lemma~\ref{lem:finite-projector-depth}, conjugate
$n_{32}(\mathfrak{p}_{v}^{r_{\pi_{v}}+\ell})$ by the Hecke translates,
whose upper-left entries are in $\mathfrak{p}_{v}^{-\ell}$.
Thus $V_{v}$ is fixed by this subgroup and the upper-left $K_{v}'$.
Hence, for $e_{v}(y)\geq r_{v}$,
\begin{align*}
 V_{v}(g(z,y))=\mathbf{1}_{\mathcal{O}_{v}^{\times}}(z).
\end{align*}
The $n_{13}$-integration contributes $\Vol(K_{0,v}'[r_{v}])^{-1}$
and restricts $y$ to $e_{v}(y)\geq r_{v}$.
If $r_{\omega_{v}'}=0$, the integral over $z\in\mathcal{O}_{v}^{\times}$ is zero
when $\xi_{v}$ is ramified.
For unramified $\xi_{v}$, the $y$-integral is a geometric series starting at $r_{v}$.

Suppose $r_{\omega_{v}'}>0$.
For $e_{v}(y)=k\geq r_{v}>r_{\omega_{v}'}$ and $u\in\mathcal{O}_{v}^{\times}$,
$1-y\varpi_{v}^{-r_{\omega_{v}'}}u$ is a unit. Use
\begin{align*}
 g(z,y)n_{23}(b')
 =n_{23}\left(\frac{zb'}{1-yb'}\right)
 g(z,y(1-yb'))\diag(1,(1-yb')^{-1},1-yb').
\end{align*}
For $b'=\varpi_{v}^{-r_{\omega_{v}'}}u$, the diagonal matrix on the right
is $(1-yb')I_{3}$ times an element of the upper-left $K_{v}'$. Using the
Whittaker and central characters gives
\begin{align*}
 V_{v}(g(z,y)n_{23}(\varpi_{v}^{-r_{\omega_{v}'}}u))
 =\overline{\psi}_{v}\left(\frac{z\varpi_{v}^{-r_{\omega_{v}'}}u}
 {1-y\varpi_{v}^{-r_{\omega_{v}'}}u}\right)
 \omega_{v}(1-y\varpi_{v}^{-r_{\omega_{v}'}}u)
   \mathbf{1}_{\mathcal{O}_{v}^{\times}}(z).
\end{align*}
The integral over $z\in\mathcal{O}_{v}^{\times}$ is a Gauss integral for
$\omega_{v}'\overline{\xi}_{v}$ at additive conductor $r_{\omega_{v}'}$. It is zero
if $r_{\xi_{v}}>r_{\omega_{v}'}$. For $r_{\xi_{v}}\leq r_{\omega_{v}'}$, combining
the finite $u$-sum with the unit part of $y=\varpi_{v}^{k}y_{0}$ gives a
scalar multiple of
\begin{align*}
 \int_{\mathcal{O}_{v}^{\times}}
 \xi_{v}(t)\omega_{v}\omega_{v}'\overline{\xi}_{v}(1-\varpi_{v}^{k-r_{\omega_{v}'}}t)\,d^{\times} t.
\end{align*}
For unramified $\xi_{v}$, the two Gauss sums coincide, so the scalar is $1$.
Since $k-r_{\omega_{v}'}\geq r_{\omega_{v}\omega_{v}'}$, the factor $\omega_{v}\omega_{v}'(1-\varpi_{v}^{k-r_{\omega_{v}'}}t)$ is one.
The map
\begin{align*}
 t\longmapsto\frac{t}{1-\varpi_{v}^{k-r_{\omega_{v}'}}t}
\end{align*}
is a measure-preserving bijection of $\mathcal{O}_{v}^{\times}$,
with inverse $t\mapsto t/(1+\varpi_{v}^{k-r_{\omega_{v}'}}t)$.
The integral over $\mathcal{O}_{v}^{\times}$ is therefore $\mathbf{1}_{r_{\xi_{v}}=0}$.

Thus every $k\geq r_{v}$ contributes exactly
$\Vol(K_{0,v}'[r_{v}])^{-1}\xi_{v}(\varpi_{v})^{k}q_{v}^{-k(1/2+\lambda)}$ for
unramified $\xi_{v}$, and zero otherwise. For unramified $\xi_{v}$,
\begin{align*}
 \sum_{k=r_{v}}^{\infty}\xi_{v}(\varpi_{v})^{k}q_{v}^{-k(1/2+\lambda)}
 =\xi_{v}(\varpi_{v})^{r_{v}}q_{v}^{-r_{v}(1/2+\lambda)}L_{v}(1/2+\lambda,\xi_{v}).
\end{align*}
In $\mathcal{H}_{W_{v}}^{\sharp}$ the geometric series cancels
$L_{v}(1/2+\lambda,\xi_{v})$ in the denominator and leaves the denominator
$L_{v}(1/2+2s-\lambda,\widetilde{\pi}_{v}\overline{\xi}_{v}\omega_{v}\omega_{v}')$
in \eqref{eq:finite-projector-full-Tate-tail}.

For unramified $\pi_{v}$, the Schur calculation in
Lemma~\ref{lem:spherical-projector-dual-mass} gives
$V_{v}(g(z,y))=\mathbf{1}_{\mathcal{O}_{v}^{\times}}(z)$ already for
$e_{v}(y)\geq2$. For $\ell=1$, $r_{\pi_{v}}\geq2$, the same pointwise identity
holds for $e_{v}(y)\geq r_{\pi_{v}}$ by
\eqref{eq:maximal-degree-one-compact-lift}. For $2+1$ data
with $\ell=2$, it holds for $e_{v}(y)\geq r_{\pi_{v}}+1$ by
\eqref{eq:maximal-degree-two-compact-lift}.
The proof uses only this identity and the upper-left
$K_{v}'$-invariance, so it applies at these levels.
Finally $r_{\omega_{v}\omega_{v}'}=r_{\omega_{v}}\leq r_{\pi_{v}}$ when $r_{\omega_{v}'}=0$, proving the last assertion.
\end{proof}

\begin{lemma}\label{lem:maximal-higher-level-positive-residue}
Suppose $d_{v}=0$ and put $\ell=\deg L_{v}(s,\pi_{v})$.
If
\begin{align*}
 r_{v}\geq\max\{r_{\pi_{v}}+\ell,r_{\omega_{v}'}+1\},
\end{align*}
then
\begin{equation}\label{eq:maximal-higher-level-positive-residue}
 \mathcal{H}_{W_{v}}^{\sharp}(\mathbf{1},1/2;\omega_{v}',s)
 =\frac{q_{v}^{-r_{v}}\mathbf{1}_{r_{v}\geq r_{\omega_{v}'}+r_{\omega_{v}\omega_{v}'}}}
 {\Vol(K_{0,v}'[r_{v}])L_{v}(2s,\widetilde{\pi}_{v}\otimes\omega_{v}\omega_{v}')}.
\end{equation}
For $|\Re(s)|\leq\delta$,
\begin{equation}\label{eq:maximal-higher-level-positive-residue-bound}
\left|\epsilon(2s,\widetilde{\pi}_{v}\otimes\omega_{v}\omega_{v}',\psi_{v})
\mathcal{H}_{W_{v}}^{\sharp}(\mathbf{1},1/2;\omega_{v}',s)\right|\ll\frac{q_{v}^{(1/2+2\delta)r_{\widetilde{\pi}_{v}\otimes\omega_{v}\omega_{v}'}}}{|L_{v}(2s,\widetilde{\pi}_{v}\otimes\omega_{v}\omega_{v}')|}.
\end{equation}
\end{lemma}

\begin{proof}
For $r_{\omega_{v}'}=0$, the formula follows from
Lemma~\ref{lem:finite-projector-full-Tate-tail}, since
$r_{\omega_{v}\omega_{v}'}=r_{\omega_{v}}\leq r_{v}$.
For $r_{\omega_{v}'}>0$, the unit integral on $e_{v}(y)=k\geq r_{v}$ is
$A_{k-r_{\omega_{v}'}}$, where
\begin{align*}
 A_{j}=\int_{\mathcal{O}_{v}^{\times}}\omega_{v}\omega_{v}'(1-\varpi_{v}^{j}u)\,d^{\times} u
 =\begin{cases}
 0,&1\leq j<r_{\omega_{v}\omega_{v}'}-1,\\
 -(q_{v}-1)^{-1},&j=r_{\omega_{v}\omega_{v}'}-1\geq1,\\
 1,&j\geq r_{\omega_{v}\omega_{v}'},
 \end{cases}
\end{align*}
by subtraction of the character integrals over $1+\mathfrak{p}_{v}^{j}$
and $1+\mathfrak{p}_{v}^{j+1}$.
Summing over $k\geq r_{v}$ gives
\begin{align*}
 \zeta_{v}(1)^{-1}\sum_{k=r_{v}}^{\infty}q_{v}^{-k}A_{k-r_{\omega_{v}'}}
 =q_{v}^{-r_{v}}\mathbf{1}_{r_{v}\geq r_{\omega_{v}'}+r_{\omega_{v}\omega_{v}'}}:
\end{align*}
when $r_{v}<r_{\omega_{v}'}+r_{\omega_{v}\omega_{v}'}$, the term at $k=r_{\omega_{v}'}+r_{\omega_{v}\omega_{v}'}-1$ cancels the tail beginning
at $r_{\omega_{v}'}+r_{\omega_{v}\omega_{v}'}$.
Finally, $q_{v}^{-r_{v}}/\Vol(K_{0,v}'[r_{v}])=1+q_{v}^{-1}$ and the absolute value of the epsilon factor give \eqref{eq:maximal-higher-level-positive-residue-bound}.
\end{proof}

\begin{lemma}
Suppose $d_{v}=0$, $\deg L_{v}(s,\pi_{v})=1$, and $r_{\pi_{v}}\geq2$.
For $V_{v}$ from Lemma~\ref{lem:finite-projector-depth},
\begin{equation}\label{eq:maximal-degree-one-compact-lift}
 V_{v}(g(z,y))=\mathbf{1}_{\mathcal{O}_{v}^{\times}}(z),
 \qquad e_{v}(y)\geq r_{\pi_{v}}.
\end{equation}
Consequently, \eqref{eq:maximal-higher-level-positive-residue} and
its bound hold for $r_{v}\geq r_{\pi_{v}}$ with $r_{v}>r_{\omega_{v}'}$.
\end{lemma}

\begin{proof}
Write $L_{v}(s,\pi_{v})^{-1}=1-\alpha q_{v}^{-s}$.
With $F_{i,j}$ from the proof of Lemma~\ref{lem:spherical-projector-dual-mass},
the expansion with one Satake parameter is
\begin{align*}
 V_{v}(g(z,y))=F_{0,0}-\alpha F_{0,1}+\alpha^{2}q_{v}^{-1}F_{1,1}.
\end{align*}
Evaluate $\mathcal{T}_{1,v}V_{v}^{\mathrm{new}}=q_{v}\alpha V_{v}^{\mathrm{new}}$
at $g(z,y)\diag(\varpi_{v}^{-1},\varpi_{v}^{-1},1)$. The first $q_{v}^{2}$
cosets contribute $q_{v}^{2}F_{0,1}$. The others have lower-unipotent
parameter in $\mathfrak{p}_{v}^{r_{\pi_{v}}}$ and diagonal value
$V_{v}^{\mathrm{new}}(\diag(z/\varpi_{v},z,1))=0$. Thus
$F_{0,1}=q_{v}^{-1}\alpha F_{1,1}$ and the diagonal newvector formula
gives $V_{v}(g(z,y))=F_{0,0}=\mathbf{1}_{\mathcal{O}_{v}^{\times}}(z)$.
The calculation of \eqref{eq:maximal-higher-level-positive-residue}
now applies with this restriction for $r_{v}\geq r_{\pi_{v}}$ and
$r_{v}>r_{\omega_{v}'}$.
\end{proof}

\begin{lemma}\label{lem:maximal-Steinberg-positive-residue}
Suppose $d_{v}=0$, $\tau_{v}=\mathrm{St}\otimes\eta_{v}$ for an
unramified unitary character $\eta_{v}$ of $F_{v}^{\times}$, and
$\chi_{v}$, $\omega_{v}'$ are unramified.
For every $r_{v}\geq1$ and $|\Re(s)|\leq\delta<10^{-2}$,
\begin{equation}\label{eq:maximal-Steinberg-positive-residue}
\big|\epsilon(2s,\widetilde{\pi}_{v}\otimes\omega_{v}\omega_{v}',\psi_{v})
 \mathcal{H}_{W_{v}}^{\sharp}(\mathbf{1}_{v},1/2;\omega_{v}',s)\big|
 \ll q_{v}^{1/2+8\delta}.
\end{equation}
The implied constant is absolute; here $r_{\pi_{v}}=1$.
\end{lemma}

\begin{proof}
We have the embedding
\begin{align*}
 \pi_{v}\hookrightarrow\Ind_{B}^{\mathrm{GL}_{3}(F_{v})}
 (\eta_{v}|\cdot|_{v}^{1/2},\chi_{v},\eta_{v}|\cdot|_{v}^{-1/2}).
\end{align*}
Applying \eqref{eq:spherical-projector-hecke-expansion} to the principal
series' spherical Whittaker function with value $1$ at $I_{3}$ and parameters
$(\eta_{v}(\varpi_{v})q_{v}^{-1/2},\chi_{v}(\varpi_{v}),\eta_{v}(\varpi_{v})q_{v}^{1/2})$
gives $V_{v}$ because its restriction agrees with Lemma~\ref{lem:finite-projector-depth}
and determines the Whittaker function by
\cite[Proposition~3.2]{JS83}.
The Schur formulas in Lemma~\ref{lem:spherical-projector-dual-mass} give 
\begin{align*}
 V_{v}(g(z,y))=\mathbf{1}_{\mathcal{O}_{v}^{\times}}(z),
 \qquad e_{v}(y)\geq2.
\end{align*}

For unramified $\xi_{v}$, take $t$ and $X$ as defined in the proof of
Lemma~\ref{lem:spherical-projector-dual-mass}.
Put $a=(\eta_{v}\chi_{v})(\varpi_{v})q_{v}^{-1/2}$, $b=\eta_{v}^{2}(\varpi_{v})$,
and $D(X)=(1-bX)(1-aX)$.
With the Satake parameters above, the polynomial $A_{q_{v}}$ in
\eqref{eq:spherical-projector-exact-transform} satisfies
\begin{align*}
 \frac{A_{q_{v}}(X)}{1-q_{v}aX}
 =\left(1-\frac{q_{v}+1}{q_{v}(q_{v}-1)}\right)D(X)+q_{v}^{-1}+q_{v}^{-2}-q_{v}^{-1}bX-baX^{2}.
\end{align*}
In \eqref{eq:spherical-projector-exact-transform}, the integral is divided
by $D(X)^{-1}(1-q_{v}aX)^{-1}$, with the extra factor from the
contragredient inducing character $\eta_{v}^{-1}|\cdot|_{v}^{-1/2}$,
so dividing its right-hand side by $1-q_{v}aX$ normalizes the same
integral by $D(X)^{-1}$ and gives
\begin{equation}\label{eq:maximal-Steinberg-projector-transform}
 \mathcal{H}_{W_{v}}^{\sharp}(\xi_{v},\lambda;\omega_{v}',s)
 =\Vol(K_{0,v}'[r_{v}])^{-1}
 \begin{cases}
 t(1-t)\frac{A_{q_{v}}(X)}{1-q_{v}aX}+t^{2}D(X),&r_{v}=1,\\
 t^{r_{v}}D(X),&r_{v}\geq2.
 \end{cases}
\end{equation}
At $\lambda=1/2$, $\xi_{v}=\mathbf{1}_{v}$, we have $t=q_{v}^{-1}$,
$|a|=q_{v}^{-1/2}$, $|b|=1$ and $|X|\leq q_{v}^{2\delta}$. Hence
$|D(X)|+|A_{q_{v}}(X)/(1-q_{v}aX)|\ll q_{v}^{4\delta}$.
Also $\Vol(K_{0,v}'[r_{v}])^{-1}=(1+q_{v}^{-1})q_{v}^{r_{v}}$.
Thus the transform is $\ll q_{v}^{4\delta}$, uniformly in $r_{v}$.
Finally $|\epsilon(2s,\widetilde{\pi}_{v}\otimes\omega_{v}\omega_{v}',\psi_{v})|
=q_{v}^{1/2-2\Re(s)}\leq q_{v}^{1/2+2\delta}$, proving
\eqref{eq:maximal-Steinberg-positive-residue}.
\end{proof}

\begin{lemma}\label{lem:maximal-Steinberg-unitary-dual}
Suppose $d_{v}=0$, $\tau_{v}$ is an unramified unitary twist of
Steinberg, and $\chi_{v}$, $\omega_{v}'$ are unramified.
Let $r_{v}\geq1$ and $\Re(\lambda)=0$. For unramified unitary $\xi_{v}$,
$|\mathcal{H}_{W_{v}}^{\sharp}(\xi_{v},\lambda;\omega_{v}',0)|\ll q_{v}^{r_{v}/2}$.
For ramified unitary $\xi_{v}$,
$|\mathcal{H}_{W_{v}}^{\sharp}(\xi_{v},\lambda;\omega_{v}',0)|
\ll q_{v}^{-1}\mathbf{1}_{r_{v}=r_{\xi_{v}}=1}$.
\end{lemma}

\begin{proof}
For unramified $\xi_{v}$, \eqref{eq:maximal-Steinberg-projector-transform}
at $s=0$ gives the bound with $|t|=|X|=q_{v}^{-1/2}$ and
$|D(X)|+|A_{q_{v}}(X)/(1-q_{v}aX)|\ll1$.
For ramified $\xi_{v}$, use the Gauss integral as in
\eqref{eq:spherical-projector-ramified-bound}.
\end{proof}

\subsubsection{Finite local estimates for $R_{\mathrm{dual}}^{+}$}
\begin{lemma}\label{lem:maximal-positive-residue-local-bound}
Suppose $d_{v}=0$, $\pi_{v}=\tau_{v}\boxplus\chi_{v}$, and
$\omega_{v}=\omega_{v}'=\mathbf{1}$. For $r_{v}\geq\max\{1,r_{\pi_{v}}\}$ and
$|\Re(s)|\leq\delta<10^{-2}$,
\begin{align*}
\left|\epsilon(2s,\widetilde{\pi}_{v},\psi_{v})
 \mathcal{H}_{W_{v}}^{\sharp}(\mathbf{1}_{v},1/2;\mathbf{1}_{v},s)\right|
\ll q_{v}^{r_{\pi_{v}}/2+\theta_{v}+16\delta r_{v}} \leq q_{v}^{r_{v}/2+r_{\pi_{v}}/4+16\delta r_{v}},
\end{align*}
where $\theta_{v}$ is as in
\eqref{eq:maximal-local-exponent}.
\end{lemma}

\begin{proof}
If $L_{v}(s,\pi_{v})$ has degree zero,
Lemma~\ref{lem:finite-projector-full-Tate-tail} applies for every $r_{v}\geq r_{\pi_{v}}$.
If it is one, $r_{\pi_{v}}\geq2$; use
\eqref{eq:maximal-degree-one-compact-lift} at every allowed level.
In both cases,
\begin{align*}
 \mathcal{H}_{W_{v}}^{\sharp}(\mathbf{1}_{v},1/2;\mathbf{1}_{v},s)
 =\frac{q_{v}^{-r_{v}}}{\Vol(K_{0,v}'[r_{v}])L_{v}(2s,\widetilde{\pi}_{v})}.
\end{align*}
Here $|L_{v}(2s,\widetilde{\pi}_{v})^{-1}|\ll q_{v}^{2\delta}$ and
$q_{v}^{-r_{v}}/\Vol(K_{0,v}'[r_{v}])=1+q_{v}^{-1}$.
The epsilon factor has modulus $q_{v}^{r_{\pi_{v}}(1/2-2\Re(s))}$,
so their product is $O(q_{v}^{r_{\pi_{v}}/2+2\delta(r_{\pi_{v}}+1)})$.

If the degree is two, trivial central character excludes one ramified
inducing character and two unramified characters. The remaining
type is an unramified unitary Steinberg twist plus an unramified
character. Lemma~\ref{lem:maximal-Steinberg-positive-residue}
gives the estimate at every $r_{v}\geq1$, including the minimal level.

If the degree is three, $r_{\pi_{v}}=0$ and the epsilon factor is one.
Here $\vartheta_{\pi_{v}}=\theta_{v}$, so
\eqref{eq:spherical-projector-residue-contour} at
$(\xi_{v},\lambda)=(\mathbf{1}_{v},1/2)$ gives the bound.

Finally, $\theta_{v}=0$ unless $\tau_{v}\otimes\chi_{v}$ is unramified and
non-tempered. If $\tau_{v}\otimes\chi_{v}$ is unramified, trivial central character
forces $\chi_{v}$ and $\tau_{v}$ to be unramified as well.
Thus $r_{\pi_{v}}+4\theta_{v}\leq2r_{v}$: in the ramified case
use $r_{\pi_{v}}\leq r_{v}$, and in the unramified case use
$4\theta_{v}<2\leq2r_{v}$.
\end{proof}

\subsubsection{Finite local estimates for $R_{\mathrm{dual}}^{-}$}
\begin{cor}\label{cor:maximal-negative-residue-compact-range}
Suppose $d_{v}=0$, $\pi_{v}=\tau_{v}\boxplus\chi_{v}$, and
$\omega_{v}=\omega_{v}'=\mathbf{1}$. For every
$r_{v}\geq\max\{1,r_{\pi_{v}}\}$ and $|\Re(s)|\leq\delta<10^{-2}$,
\begin{align*}
 \left|\epsilon(2s,\overline{\chi}_{v},\psi_{v})
 \mathcal{H}_{W_{v}}^{\sharp}(\overline{\chi}_{v},2s-1/2;\mathbf{1}_{v},s)\right|
 \ll q_{v}^{r_{v}+8\delta r_{v}}.
\end{align*}
The local factor vanishes if $\overline{\chi}_{v}$ is ramified.
\end{cor}

\begin{proof}
If $L_{v}(s,\pi_{v})$ has degree zero or one,
Lemma~\ref{lem:finite-projector-full-Tate-tail} applies at every
allowed level; in degree one, use its refinement and
$r_{\pi_{v}}\geq2$. At
$(\xi_{v},\lambda)=(\overline{\chi}_{v},2s-1/2)$, it gives
\begin{align*}
 \mathcal{H}_{W_{v}}^{\sharp}(\overline{\chi}_{v},2s-1/2;\mathbf{1}_{v},s)
 =\frac{\overline{\chi}_{v}(\varpi_{v})^{r_{v}} q_{v}^{-2sr_{v}}
 \mathbf{1}_{r_{\chi_{v}}=0}}
 {\Vol(K_{0,v}'[r_{v}])L_{v}(1,\widetilde{\tau}_{v}\otimes\chi_{v})\zeta_{v}(1)}.
\end{align*}
Here we used
$\widetilde{\pi}_{v}\otimes\chi_{v}
=(\widetilde{\tau}_{v}\otimes\chi_{v})\boxplus\mathbf{1}_{v}$.
Since $|L_{v}(1,\widetilde{\tau}_{v}\otimes\chi_{v})^{-1}|\ll1$,
the weight is $O(q_{v}^{(1+2\delta)r_{v}}\mathbf{1}_{r_{\chi_{v}}=0})$.

If the degree is two, trivial central character forces $\tau_{v}$ to
be an unramified unitary twist of Steinberg and $\chi_{v}$ to be
unramified. For degree three, $\pi_{v}$ is unramified.
In \eqref{eq:maximal-Steinberg-projector-transform} for degree two and
\eqref{eq:spherical-projector-exact-transform} for degree three,
$(\xi_{v},\lambda)=(\overline{\chi}_{v},2s-1/2)$ gives
\begin{align*}
t=\overline{\chi}_{v}(\varpi_{v})q_{v}^{-2s},\qquad
 X=\chi_{v}(\varpi_{v})q_{v}^{-1}.
\end{align*}
For $r_{v}\geq2$ both formulas give $t^{r_{v}}\Vol(K_{0,v}'[r_{v}])^{-1}
L_{v}(1,\widetilde{\pi}_{v}\otimes\chi_{v})^{-1}$. For $r_{v}=1$ the polynomial factors are $O(1)$ at $|X|=q_{v}^{-1}$.
Together with $|t|\leq q_{v}^{2\delta}$ and
$\Vol(K_{0,v}'[1])^{-1}=q_{v}+1$, this gives
$\mathcal{H}_{W_{v}}^{\sharp}=O(q_{v}^{1+4\delta})$.
Whenever $\chi_{v}$ is unramified,
$\epsilon(2s,\overline{\chi}_{v},\psi_{v})=1$, proving both assertions.

At places dividing the different,
\eqref{eq:maximal-positive-residue-different-transport} at
$(\xi_{v},\lambda)=(\overline{\chi}_{v},2s-1/2)$ and
\eqref{eq:different-epsilon-factor} for $n=1$ give
\begin{equation}\label{eq:maximal-negative-residue-different-transport}
 |\epsilon(2s,\overline{\chi}_{v},\psi_{v})\mathcal{H}_{W_{v}}^{\sharp}|
 =q_{v}^{d_{v}(3/2-4\Re(s))}
 |\epsilon(2s,\overline{\chi}_{v},\psi_{v}^{0})\mathcal{H}_{W_{v}^{0}}^{\sharp}|,
\end{equation}
at $(\overline{\chi}_{v},2s-1/2;\mathbf{1}_{v},s)$. Here
$q_{v}^{d_{v}(3/2-4\Re(s))}\ll_{F}1$ for $|\Re(s)|\leq\delta$.
\end{proof}

\subsubsection{Archimedean estimates and global bounds}
\begin{lemma}\label{lem:maximal-negative-quadratic-archimedean}
For $|s|\leq\delta$ sufficiently small, at an Archimedean place,
\begin{align*}
 |\gamma(2s,\chi_{v},\psi_{v})
 \mathcal{H}_{W_{v}}(\chi_{v},2s-1/2;1,s)|
 \ll_{F,\delta,\epsilon}C_{v}(\pi)^{\epsilon}
 C_{v}^{1+2\delta}[1+C_{v}^{-1}C_{v}(\pi)].
\end{align*}
The product on the left is holomorphic near $s=0$.
\end{lemma}

\begin{proof}
To write $\mathcal{H}_{W_{v}}^{\circ}$ as a Tate integral, put
$d=[F_{v}:\mathbb{R}]$ and define
\begin{align*}
 \Phi_{v}^{\circ}(y)=\int_{F_{v}^{\times}}
 h_{3,v}(z,z,C_{v}^{1/d}y)
 \chi_{v}(z)\,d^{\times} z.
\end{align*}
Fourier inversion, as in \eqref{9.69}, gives
\begin{align*}
 \mathcal{H}_{W_{v}}^{\circ}(\chi_{v},2s-1/2;1,s)
 =\int_{F_{v}^{\times}}\Phi_{v}^{\circ}(y)\chi_{v}(y)|y|_{v}^{2s}\,d^{\times}y.
\end{align*}

By \eqref{equ9.11},
$t\mapsto C_{v}^{-1}\Phi_{v}^{\circ}(C_{v}^{-1/d}t)$
has fixed compact support and uniform derivative bounds. Hence the
Tate functional equation gives
\begin{multline*}
|\gamma(2s,\chi_{v},\psi_{v})
   \mathcal{H}_{W_{v}}^{\circ}(\chi_{v},2s-1/2;1,s)|=\left|\int_{F_{v}^{\times}}\widehat{\Phi}_{v}^{\circ}(x)
 \chi_{v}(x)|x|_{v}^{1-2s}\,d^{\times}x\right|\\
 \ll_{F,\delta}C_{v}^{1-2\Re(s)}.
\end{multline*}
The integral is holomorphic for $2\Re(s)<1$ by Fourier decay.

In \eqref{9.69}, use the rapid decay of $\mathcal{F}_{12}h_{3,v}$
in $(b_{1},b_{2})$ and apply
\eqref{eq:unitary-residue-archimedean-endpoint} to $\mathcal{V}_{v}$:
\begin{multline*}
\left|\mathcal{H}_{W_{v}}(\chi_{v},2s-1/2;1,s)
 -\mathcal{H}_{W_{v}}^{\circ}(\chi_{v},2s-1/2;1,s)\right|\\
\ll_{F,\delta,\epsilon}
 C_{v}(\pi)^{\epsilon}C_{v}^{1-2\Re(s)}
 [1+C_{v}^{-1}C_{v}(\pi)],\qquad |s|\leq\delta. 
\end{multline*}
On $\lambda=2s-1/2$, $\mathcal{V}_{v}$ is independent of $s$ and
satisfies $\mathcal{V}_{v}(y_{v})=O(|y_{v}|)$ as $y_{v}\to0$ by
\eqref{eq:unitary-residue-archimedean-endpoint}.
Thus the $y_{v}$-integral in \eqref{9.69} is holomorphic for
$\Re(s)>-1/(2d)$.
Since $\chi_{v}^{2}=1$, $\gamma(2s,\chi_{v},\psi_{v})$ is holomorphic
and bounded near $s=0$. Multiplying the bound for the difference by
this gamma factor and adding the bound for
$\gamma(2s,\chi_{v},\psi_{v})\mathcal{H}_{W_{v}}^{\circ}$ proves the lemma,
since $C_{v}^{1-2\Re(s)}\leq C_{v}^{1+2\delta}$ for $|s|\leq\delta$.
\end{proof}

\begin{cor}\label{cor:maximal-negative-residue-covered-range}
We have 
\begin{align*}
 \sup_{|s|=\rho}|R_{\mathrm{dual}}^{-}(s)|
 \ll_{F,\varepsilon}C(\pi)^{\varepsilon}
 Q^{1+\varepsilon}\mathbf{C}_{\infty}^{1+\varepsilon}\prod_{v\mid\infty}[1+C_{v}^{-1}C_{v}(\pi)].
\end{align*}
\end{cor}

\begin{proof}
In \eqref{eq:maximal-negative-dual-residue}, the functional
equation \eqref{eq:maximal-residue-functional-equation} for $\chi$
replaces $L(2s,\chi)$ by $L(1-2s,\chi)$ times local factors.
On $|s|=\rho\leq2\delta$, Corollary~\ref{cor:maximal-negative-residue-compact-range}
bounds the finite factors and Lemma~\ref{lem:maximal-negative-quadratic-archimedean}
bounds the Archimedean factors. Thus
\begin{align*}
 |R_{\mathrm{dual}}^{-}(s)|
 \ll_{F,\varepsilon}|L(1,\tau)L(1-2s,\chi)|C(\pi)^{\varepsilon}
 Q^{1+16\delta+\varepsilon}\mathbf{C}_{\infty}^{1+4\delta}
 \prod_{v\mid\infty}[1+C_{v}^{-1}C_{v}(\pi)].
\end{align*}
At $v\mid\mathfrak{D}_{F}$, the factor
$q_{v}^{d_{v}(3/2-4\Re(s))}$ in
\eqref{eq:maximal-negative-residue-different-transport} is $O_{F}(1)$
on $|s|=\rho$, since $q_{v}$ and $d_{v}$ depend only on $F$.
For $r_{v}=0$, \eqref{eq:finite-zero-level-dual-bound} at
$(\xi_{v},\lambda)=(\chi_{v},2s-1/2)$ and
\eqref{eq:different-epsilon-factor} with $n=1$ bound the local
factor by $O_{F}(1)$.
The standard bound
$|L(1,\tau)L(1-2s,\chi)|\ll_{F,\varepsilon}C(\pi)^{\varepsilon}$
on this circle proves the result, after taking $\delta$ and the
initial $\varepsilon$ sufficiently small.
\end{proof}

\begin{cor}\label{cor:maximal-positive-residue-sharp-bound}
For $|s|=\rho$ and sufficiently small $\delta$,
\begin{align*}
\sup_{|s|=\rho}|R_{\mathrm{dual}}^{+}(s)|
\ll_{F,\varepsilon}C(\pi)^{\varepsilon}
 \mathbf{C}_{\infty}^{1/2+\varepsilon}Q^{1/2+\varepsilon}
C(\pi)^{1/4+\varepsilon}\prod_{v\mid\infty}
 [1+(C_{v}^{-1}C_{v}(\pi))^{2}].
\end{align*}
\end{cor}

\begin{proof}
After the functional equation \eqref{eq:maximal-residue-functional-equation},
the residue \eqref{eq:maximal-positive-dual-residue} is bounded by
Lemma~\ref{lem:maximal-positive-residue-local-bound} at finite places
and \eqref{eq:positive-residue-archimedean-bound} at infinity. Thus
\begin{multline*}
 |R_{\mathrm{dual}}^{+}(s)|
 \ll_{F,\varepsilon}|L(1-2s,\widetilde{\pi})|
 Q^{1/2+32\delta+\varepsilon}
 \mathbf{C}_{\infty}^{1/2+\varepsilon}\\
 C_{\mathrm{fin}}(\pi)^{1/4}
 C_{\infty}(\pi)^{1/4+\varepsilon}
 \prod_{v\mid\infty}[1+(C_{v}^{-1}C_{v}(\pi))^{2}].
\end{multline*}
For $r_{v}=0$, use
\eqref{eq:finite-zero-level-dual-bound} at
$(\xi_{v},\lambda)=(\mathbf{1}_{v},1/2)$.
Use $C_{\mathrm{fin}}(\pi)C_{\infty}(\pi)=C(\pi)$ and
$L(1-2s,\widetilde{\pi})\ll_{F,\varepsilon}C(\pi)^{\varepsilon}$
on $|s|=\rho$, and take $\delta$ and the initial $\varepsilon$
sufficiently small.
\end{proof}

\section{\texorpdfstring{The Twisted Fourth Moment of $\mathrm{GL}_{2}$ $L$-functions}{The Twisted Fourth Moment of GL2 L-functions}}\label{sec14}
We specialize type~\RNum{2} to minimal Eisenstein data, using the
local evaluations in Appendix~\ref{app:fourth-local}.

Let $F$ be totally real. Initially take
$\boldsymbol{\nu}=(\nu_{1},\nu_{2},\nu_{3})\in\mathcal{R}_{\mathrm{min}}$
as in \eqref{f2.2}, with $\nu_{1}+\nu_{2}+\nu_{3}=0$, and put
$\pi=|\cdot|^{\nu_{1}}\boxplus|\cdot|^{\nu_{2}}\boxplus|\cdot|^{\nu_{3}}$.
Let $\mathfrak{q},\mathfrak{n}$ be nonzero integral ideals with
$\mathfrak{q}+\mathfrak{n}=\mathcal{O}_{F}$, and write
$r_{v}=e_{v}(\mathfrak{q})$, $n_{v}=e_{v}(\mathfrak{n})$.
At each $v\mid\infty$, choose
$h_{v}\in C_{c}^{\infty}(K_{v}'\backslash G'(F_{v})/K_{v}')$ and set
$\mathcal{T}_{v}=\mathcal{S} h_{v}$ from \eqref{eq10.4}.
Write $\mathcal{T}=(\mathcal{T}_{v})_{v\mid\infty}$.

Let $\varphi_{\boldsymbol{\nu}}=E(f_{\boldsymbol{\nu}})\in\pi$ be
the Eisenstein series of \eqref{f2.5}, where
$f_{\boldsymbol{\nu}}=\bigotimes_{v} f_{\boldsymbol{\nu},v}$ uses the flat
sections and operations in \eqref{fc14.3}, \eqref{c14.5}, and
\eqref{c14.6}. Its local Whittaker vectors are
$\mathcal{J}_{\psi_{v}}f_{\boldsymbol{\nu},v}$.
We normalize the local vectors as
$W_{v}=c_{v}\mathcal{J}_{\psi_{v}}f_{\boldsymbol{\nu},v}$, with $c_{v}$ as in
\eqref{eq:jacquet-normalization-transport}. At the real places the change
of Whittaker characters in \textsection\ref{sec11} identifies the weights
with those of the vector \eqref{fc11.3}.

For each place $v$, let $f_{\boldsymbol{\nu},v}^{\circ}$ denote the spherical flat section in $\pi_{v}$. On Iwasawa components 
\begin{align*}
g_{v}=\diag(a_{1,v},a_{2,v},a_{3,v})u_{v}k_{v},\ \ a_{i,v}\in F_{v}^{\times},\ 1\leq i\leq 3,\ u_{v}\in N(F_{v}),\ k_{v}\in K_{v},
\end{align*}  
the section $f_{\boldsymbol{\nu},v}^{\circ}$ is given by 
\begin{equation}\label{c14.1}
f_{\boldsymbol{\nu},v}^{\circ}(g_{v})=|a_{1,v}|_{v}^{1+\nu_{1}}|a_{2,v}|_{v}^{\nu_{2}}|a_{3,v}|_{v}^{-1+\nu_{3}}.
\end{equation} 

The local generic and degenerate integrals below use
$f_{\boldsymbol{\nu},v}$ with the original $\psi_{v}$ and measures;
the factors at places dividing the different are included in
Appendix~\ref{app:fourth-local}.

\subsection{The Global Fourth-Moment Formula}\label{sec14.8}
Put
\begin{align*}
\mathbf{s}=(s_{1},s_{2}),\qquad
 \mathbf{s}^{\vee}=(s_{1}',s_{2}')
 =\left(\frac{s_{2}-s_{1}}{2},\frac{3s_{1}+s_{2}}{2}\right).
\end{align*}
All identities below are continued meromorphically to
$|\Re(s_{1})|+|\Re(s_{2})|\leq10^{-1}$ and
$|\Re(\nu_{i})|\leq10^{-2}$, $1\leq i\leq3$.

\subsubsection{The spectral and dual sides}\label{sec14.8.1}
Use $\mathcal{L}_{\fin}$ from \textsection\ref{def:finite-fourth-moment-factor}.
\phantomsection\label{def:spectral-multiplicity}
Write $\epsilon_{\sigma}=1$ on the
cuspidal spectrum and $\epsilon_{\sigma}=2$ on the continuous spectrum.
For $v<\infty$, a superscript $0$ uses
$\psi_{v}^{0}$ from \textsection\ref{def:unramified-additive-character}
and the measures in \eqref{eq:finite-character-measures}. Take $W_{v}^{0}$ from
\eqref{fc11.41} or \eqref{c12.2}, according as $v\nmid\mathfrak{q}$ or
$v\mid\mathfrak{q}$, constructed from
the spherical $\psi_{v}^{0}$-Whittaker vector with value $1$ at $I_{3}$.
Use the local weight $\mathcal{H}_{W_{v}^{0}}$ of \eqref{fc10.1}.

\phantomsection\label{def:fourth-spectral-weights}
Define
\begin{multline*}
\mathcal{H}_{\sigma,\mathcal{T}}(\mathbf{s},\boldsymbol{\nu};\mathfrak{q},\mathfrak{n})
:=\epsilon_{\sigma}D_{F}^{3s_{1}+s_{2}-3/2+\nu_{3}-\nu_{1}}
\prod_{v\mid\infty}\pi \mathcal{H}_{\mathcal{J}_{\psi_{v}}f_{\boldsymbol{\nu},v}}(\sigma_{v},\mathbf{s},\mathbf{1})\\
\prod_{v\mid\mathfrak{q}}\frac{\mathcal{H}_{W_{v}^{0}}(\sigma_{v},\mathbf{s},\mathbf{1})}
{L_{v}(1/2+s_{1},\pi_{v}\times\sigma_{v})L_{v}(1/2+s_{2},\sigma_{v})},
\end{multline*}
\begin{multline*}
\mathcal{H}_{\sigma,\mathcal{T}}^{\vee}(\mathbf{s}^{\vee},\boldsymbol{\nu};\mathfrak{q},\mathfrak{n})
:=\epsilon_{\sigma}D_{F}^{3s_{1}'+s_{2}'-3/2+\nu_{3}-\nu_{1}}
\prod_{v\mid\infty}\pi \mathcal{H}_{\mathcal{J}_{\psi_{v}}f_{\boldsymbol{\nu},v}}^{\vee}(\sigma_{v},\mathbf{s}^{\vee},\mathbf{1})\\
\prod_{v\mid\mathfrak{n}}\frac{q_{v}^{(1/2+s_{2})n_{v}}\mathcal{H}_{W_{v}^{0}}^{\vee}(\sigma_{v},\mathbf{s}^{\vee},\mathbf{1})}
{L_{v}(1/2+s_{1}',\pi_{v}\times\sigma_{v})L_{v}(1/2+s_{2}',\sigma_{v})}.
\end{multline*}

Use \eqref{equ12.11} at $v\mid\mathfrak{q}$ and
\eqref{eq:dual-level-superposition} at $v\mid\mathfrak{n}$, both with
$\psi_{v}^{0}$ and its measures. Theorem~\ref{prop10.2} gives the real dual
weight with $\beta_{v}=\mathbf{L}(\pi_{v})^{-1}$.
On the original side, take $\beta_{v}=\mathbf{L}(\pi_{v})^{-1}$ in
\eqref{eq:real-weight-evaluated}:
\begin{align*}
\mathcal{H}_{\mathcal{J}_{\psi_{v}}f_{\boldsymbol{\nu},v}}(\sigma_{v},\mathbf{s},\mathbf{1})
=\frac{L_{v}(1/2+s_{1},\pi_{v}\times\sigma_{v})L_{v}(1/2+s_{2},\sigma_{v})
\mathcal{T}_{v}(\nu_{\sigma_{v}})}{\pi\,\mathbf{L}(\pi_{v}) L_{v}(1,\sigma_{v},\Ad)}.
\end{align*}
We also set
\begin{equation}\label{eq14.123}
\widetilde{\lambda}_{\sigma,\mathbf{s}}(\mathfrak{q}):=\prod_{v\mid\mathfrak{q}}\Big[\lambda_{\sigma_{v}}(\mathfrak{p}_{v}^{e_{v}(\mathfrak{q})})-q_{v}^{-\frac{1+3s_{1}+s_{2}}{2}}\lambda_{\sigma_{v}}(\mathfrak{p}_{v}^{e_{v}(\mathfrak{q})-1})\Big].
\end{equation}
Every $\sigma\in\mathcal{F}(\mathfrak{n};\mathbf{1})$ is unramified at
$v\mid\mathfrak{q}$, since $\mathfrak{q}+\mathfrak{n}=\mathcal{O}_{F}$.

The two moments are
\begin{align*}
&\widetilde{\mathcal{I}}_{\mathcal{T}}(\mathbf{s},\boldsymbol{\nu};\mathfrak{q},\mathfrak{n}):=\int_{\sigma\in\mathcal{F}^{\circ}(\mathfrak{q};\mathbf{1})}\lambda_{\sigma}(\mathfrak{n})\mathcal{L}_{\fin}(\mathbf{s},\pi,\sigma)\mathcal{H}_{\sigma,\mathcal{T}}(\mathbf{s},\boldsymbol{\nu};\mathfrak{q},\mathfrak{n})d\mu_{\sigma},\\ 
&\widetilde{\mathcal{I}}_{\mathcal{T}}^{\vee}(\mathbf{s}^{\vee},\boldsymbol{\nu};\mathfrak{q},\mathfrak{n}):=\frac{N_{F}(\mathfrak{q})^{\frac{1}{2}-s_{2}'}}{N_{F}(\mathfrak{n})^{\frac{1}{2}+s_{2}}}\int_{\sigma\in\mathcal{F}(\mathfrak{n};\mathbf{1})}\widetilde{\lambda}_{\sigma,\mathbf{s}}(\mathfrak{q})\mathcal{L}_{\fin}(\mathbf{s}^{\vee},\pi,\sigma)\mathcal{H}_{\sigma,\mathcal{T}}^{\vee}(\mathbf{s}^{\vee},\boldsymbol{\nu};\mathfrak{q},\mathfrak{n})d\mu_{\sigma}.
\end{align*}

\subsubsection{The degenerate terms}\label{sec14.8.2}
Let $\mathbb{W}=\{w_{2},w_{1}w_{2},w_{1}w_{2}w_{1}\}$.
For $w\in\mathbb{W}$, use $I_{\mathrm{degen}}^{\heartsuit}$ from
\textsection\ref{sec8.2}.

For $w=w_{2},w_{1}w_{2},w_{1}w_{2}w_{1}$, let
$\mathcal{L}_{\fin}(\mathbf{s},\boldsymbol{\nu};w)$ be the finite Euler
product in \textsection\ref{sec5.5}, \textsection\ref{sec5.8}, and
\textsection\ref{sec5.11}, respectively, specialized to
$\chi_{1}=\chi_{2}=\chi_{3}=\omega'=\eta=\mathbf{1}$.
Write $\mathcal{L}_{v}(\mathbf{s},\boldsymbol{\nu};w)$ for its local factor.

At finite $v\nmid\mathfrak{q}\mathfrak{n}\mathfrak{D}_{F}$,
\eqref{f5.15} with trivial characters and the normalization factor
$\mathbf{L}(\pi_{v})^{-1}$ give
\begin{align*}
\mathbf{L}(\pi_{v})I_{\mathrm{degen},v}
 (\mathbf{s},\varphi_{\boldsymbol{\nu}},\mathbf{1},\mathbf{1};w)
 =\mathcal{L}_{v}(\mathbf{s},\boldsymbol{\nu};w),
 \qquad w\in\mathbb{W}.
\end{align*}
For the dual integral, replace $\mathbf{s}$ by $\mathbf{s}^{\vee}$,
since $w_{2}$ fixes the spherical section.

Define
\begin{multline*}
\mathcal{H}_{\mathrm{degen},\mathcal{T}}(\mathbf{s},\boldsymbol{\nu};\mathfrak{q},\mathfrak{n};w)
:=\prod_{v\mid\mathfrak{q}\mathfrak{n}\mathfrak{D}_{F}}
\frac{\mathbf{L}(\pi_{v})I_{\mathrm{degen},v}(\mathbf{s},\varphi_{\boldsymbol{\nu}},\mathbf{1},\mathbf{1};w)}{\mathcal{L}_{v}(\mathbf{s},\boldsymbol{\nu};w)}\\
\prod_{v\mid\infty}I_{\mathrm{degen},v}(\mathbf{s},\varphi_{\boldsymbol{\nu}},\mathbf{1},\mathbf{1};w).
\end{multline*}

On the dual side, define
\begin{multline*}
\mathcal{H}_{\mathrm{degen},\mathcal{T}}^{\vee}
(\mathbf{s}^{\vee},\boldsymbol{\nu};\mathfrak{q},\mathfrak{n};w)
:=\prod_{v\mid\mathfrak{q}\mathfrak{n}\mathfrak{D}_{F}}
\frac{\mathbf{L}(\pi_{v})I_{\mathrm{degen},v}(\mathbf{s}^{\vee},\pi(w_{2})\varphi_{\boldsymbol{\nu}},\mathbf{1},\mathbf{1};w)}
{\mathcal{L}_{v}(\mathbf{s}^{\vee},\boldsymbol{\nu};w)}\\
\prod_{v\mid\infty}I_{\mathrm{degen},v}
(\mathbf{s}^{\vee},\pi(w_{2})\varphi_{\boldsymbol{\nu}},\mathbf{1},\mathbf{1};w).
\end{multline*}

Use \eqref{eq:finite-dual-degenerate-transport} and the specialization
for $v\mid\mathfrak{q}$ following it.
Multiplying the local integrals, with
$\mathcal{L}_{v}(\mathbf{s},\boldsymbol{\nu};w)/\mathbf{L}(\pi_{v})$
at $v\nmid\mathfrak{q}\mathfrak{n}\mathfrak{D}_{F}$, gives
\begin{equation}\label{e14.123}
I_{\mathrm{degen}}^{\heartsuit}(\mathbf{s},\varphi_{\boldsymbol{\nu}},\mathbf{1},\mathbf{1};w)=\mathbf{L}_{\fin}(\pi)^{-1}\mathcal{L}_{\fin}(\mathbf{s},\boldsymbol{\nu};w)\mathcal{H}_{\mathrm{degen},\mathcal{T}}(\mathbf{s},\boldsymbol{\nu};\mathfrak{q},\mathfrak{n};w).
\end{equation}

On the dual side,
\begin{equation}\label{e14.124}
I_{\mathrm{degen}}^{\heartsuit}(\mathbf{s}^{\vee},\pi(w_{2})\varphi_{\boldsymbol{\nu}},\mathbf{1},\mathbf{1};w)=\frac{\mathcal{L}_{\fin}(\mathbf{s}^{\vee},\boldsymbol{\nu};w)\mathcal{H}_{\mathrm{degen},\mathcal{T}}^{\vee}(\mathbf{s}^{\vee},\boldsymbol{\nu};\mathfrak{q},\mathfrak{n};w)}{\mathbf{L}_{\fin}(\pi)}.
\end{equation}

\subsubsection{The generic terms}\label{sec14.8.3}
Put
\begin{align*}
\mathcal{L}_{\fin}(\mathbf{s},\pi)
:=\frac{L(1+s_{1}+s_{2},\pi)L(1+2s_{1},\widetilde{\pi})}
{\zeta_{F}(2+3s_{1}+s_{2})},
\end{align*}
and write $\mathcal{L}_{v}(\mathbf{s},\pi_{v})$ for its local factor.
\phantomsection\label{def:fourth-generic-weights}
Define
\begin{align*}
\mathcal{H}_{\mathrm{gen},\mathcal{T}}(\mathbf{s},\boldsymbol{\nu};\mathfrak{q},\mathfrak{n})
:=\prod_{v\mid\infty}I_{\mathrm{gen},v}(\mathbf{s},\varphi_{\boldsymbol{\nu}},\mathbf{1},\mathbf{1})\prod_{v\mid\mathfrak{q}\mathfrak{n}\mathfrak{D}_{F}}\frac{\mathbf{L}(\pi_{v})I_{\mathrm{gen},v}(\mathbf{s},\varphi_{\boldsymbol{\nu}},\mathbf{1},\mathbf{1})}
{q_{v}^{r_{v}-(1/2+s_{2})n_{v}}\mathcal{L}_{v}(\mathbf{s},\pi_{v})},
\end{align*}
\begin{multline*}
\mathcal{H}_{\mathrm{gen},\mathcal{T}}^{\vee}(\mathbf{s}^{\vee},\boldsymbol{\nu};\mathfrak{q},\mathfrak{n})
:=\prod_{v\mid\infty}I_{\mathrm{gen},v}(\mathbf{s}^{\vee},\pi(w_{2})\varphi_{\boldsymbol{\nu}},\mathbf{1},\mathbf{1})\\
\prod_{v\mid\mathfrak{q}\mathfrak{n}\mathfrak{D}_{F}}\frac{\mathbf{L}(\pi_{v})I_{\mathrm{gen},v}(\mathbf{s}^{\vee},\pi(w_{2})\varphi_{\boldsymbol{\nu}},\mathbf{1},\mathbf{1})}
{q_{v}^{-(s_{1}+s_{2})r_{v}+(1/2+s_{2})n_{v}}\mathcal{L}_{v}(\mathbf{s}^{\vee},\pi_{v})}.
\end{multline*}

At finite $v\nmid\mathfrak{q}\mathfrak{n}\mathfrak{D}_{F}$, the section is spherical and fixed by $w_{2}$, and the
local generic integral is
$\mathbf{L}(\pi_{v})^{-1}\mathcal{L}_{v}(\mathbf{s},\pi_{v})$.
Multiplication gives
\begin{equation}\label{e14.125}
I_{\mathrm{gen}}^{\heartsuit}(\mathbf{s},\varphi_{\boldsymbol{\nu}},\mathbf{1},\mathbf{1})
=\mathbf{L}_{\fin}(\pi)^{-1}\mathcal{L}_{\fin}(\mathbf{s},\pi)
N_{F}(\mathfrak{q})N_{F}(\mathfrak{n})^{-1/2-s_{2}}
\mathcal{H}_{\mathrm{gen},\mathcal{T}}(\mathbf{s},\boldsymbol{\nu};\mathfrak{q},\mathfrak{n}),
\end{equation}
and 
\begin{equation}\label{e14.126}
I_{\mathrm{gen}}^{\heartsuit}(\mathbf{s}^{\vee},\pi(w_{2})\varphi_{\boldsymbol{\nu}},\mathbf{1},\mathbf{1})
=\frac{\mathcal{L}_{\fin}(\mathbf{s}^{\vee},\pi)N_{F}(\mathfrak{n})^{1/2+s_{2}}
\mathcal{H}_{\mathrm{gen},\mathcal{T}}^{\vee}(\mathbf{s}^{\vee},\boldsymbol{\nu};\mathfrak{q},\mathfrak{n})}{\mathbf{L}_{\fin}(\pi)
N_{F}(\mathfrak{q})^{s_{1}+s_{2}}}.
\end{equation}
Using $3s_{1}'+s_{2}'=2s_{2}$, we have
\begin{align*}
\mathcal{L}_{\fin}(\mathbf{s}^{\vee},\pi)
=\frac{L(1+s_{1}+s_{2},\pi)L(1+s_{2}-s_{1},\widetilde{\pi})}{\zeta_{F}(2+2s_{2})}.
\end{align*}

\subsubsection{The residual terms}\label{sec14.8.4}
Let $\sigma_{\lambda}=|\cdot|^{\lambda}\boxplus|\cdot|^{-\lambda}$,
with local components $\sigma_{\lambda,v}$.
\phantomsection\label{def:fourth-residue-densities}
Write the densities of the original and dual moments as
\begin{align*}
&D_{\mathcal{T}}(\mathbf{s},\lambda)
:=\lambda_{\sigma_{\lambda}}(\mathfrak{n})\mathcal{L}_{\fin}(\mathbf{s},\pi,\sigma_{\lambda})
\mathcal{H}_{\sigma_{\lambda},\mathcal{T}}(\mathbf{s},\boldsymbol{\nu};\mathfrak{q},\mathfrak{n}),\\
&D_{\mathcal{T}}^{\vee}(\mathbf{s}^{\vee},\lambda)
:=\frac{N_{F}(\mathfrak{q})^{1/2-s_{2}'}}{N_{F}(\mathfrak{n})^{1/2+s_{2}}}\cdot 
\widetilde{\lambda}_{\sigma_{\lambda},\mathbf{s}}(\mathfrak{q})
\mathcal{L}_{\fin}(\mathbf{s}^{\vee},\pi,\sigma_{\lambda})
\mathcal{H}_{\sigma_{\lambda},\mathcal{T}}^{\vee}(\mathbf{s}^{\vee},\boldsymbol{\nu};\mathfrak{q},\mathfrak{n}).
\end{align*}
For a meromorphic density $D(\lambda)$ with distinct moving poles, put
\begin{equation}\label{eq:residue-operator}
\mathscr{R}_{\mathbf{s},\boldsymbol{\nu}}[D]
:=\frac{1}{2}\sum_{\epsilon\in\{\pm1\}}\epsilon
\left(\sum_{i=1}^{3}\underset{\lambda=\epsilon(1/2-s_{1}-\nu_{i})}{\Res}D(\lambda)
+\underset{\lambda=\epsilon(1/2-s_{2})}{\Res}D(\lambda)\right).
\end{equation}
Here the moving poles come from
$\zeta_{F}(1/2+s_{1}+\nu_{i}\pm\lambda)$ and
$\zeta_{F}(1/2+s_{2}\pm\lambda)$, respectively.
If poles of sign $\epsilon$ coincide, use their combined contribution
\begin{equation}\label{eq:residue-cluster}
 \frac{\epsilon}{4\pi i}\int_{\mathcal{C}}D(\lambda)\,d\lambda,
\end{equation}
where $\mathcal{C}$ is a fixed positively oriented contour enclosing
these poles and no others for nearby parameters.
This continues \eqref{eq:residue-operator}, counting each multiple pole once.

For distinct poles, the residue terms of \textsection\ref{sec4.2} satisfy
\begin{equation}\label{14.125}
\begin{cases}
2D_{F}^{1/2}\mathbf{L}_{\fin}(\pi)R_{\RNum{2}}^{1,\epsilon}(\mathbf{s},\varphi_{\boldsymbol{\nu}},\mathbf{1},\mathbf{1})
=\displaystyle\frac{1}{2}\sum_{i=1}^{3}
\underset{\lambda=\epsilon(1/2-s_{1}-\nu_{i})}{\Res}D_{\mathcal{T}}(\mathbf{s},\lambda),\\[1em]
2D_{F}^{1/2}\mathbf{L}_{\fin}(\pi)R_{\RNum{2}}^{2,\epsilon}(\mathbf{s},\varphi_{\boldsymbol{\nu}},\mathbf{1},\mathbf{1})
=\displaystyle\frac{1}{2}\underset{\lambda=\epsilon(1/2-s_{2})}{\Res}D_{\mathcal{T}}(\mathbf{s},\lambda).
\end{cases}
\end{equation}
The dual families are obtained by replacing $\mathbf{s}$ with
$\mathbf{s}^{\vee}$ and $D_{\mathcal{T}}$ with $D_{\mathcal{T}}^{\vee}$.

\subsubsection{The explicit twisted fourth moment}\label{sec:proof-explicit-fourth}
\typeC*

\begin{proof}
For the spectral side of Theorem~\ref{thmii} at $\varphi_{\boldsymbol{\nu}}$,
the two Rankin--Selberg periods in Proposition~\ref{prop4.4} and the
inner products of \textsection\ref{sec3.2} give the density
\cite[(21)--(23)]{Nun23}
\begin{align*}
\frac{\epsilon_{\sigma}}
{2D_{F}^{1/2}L^{\mathrm{reg}}(1,\sigma,\Ad)}
\prod_{v\mid\infty}\pi \mathcal{H}_{\mathcal{J}_{\psi_{v}}f_{\boldsymbol{\nu},v}}(\sigma_{v},\mathbf{s},\mathbf{1})
\prod_{v<\infty}\mathcal{H}_{\mathcal{J}_{\psi_{v}}f_{\boldsymbol{\nu},v}}(\sigma_{v},\mathbf{s},\mathbf{1})\,d\mu_{\sigma}.
\end{align*}

At finite $v$, multiply the Jacquet scalar in
\eqref{eq:different-jacquet-scalar} by the weight factor in
\eqref{eq:different-raw-weight}:
\begin{align*}
 q_{v}^{d_{v}(1/2+\nu_{3}-\nu_{1})}q_{v}^{d_{v}(3s_{1}+s_{2}-2)}
 =q_{v}^{d_{v}(3s_{1}+s_{2}-3/2+\nu_{3}-\nu_{1})}.
\end{align*}
The spherical Jacquet vector with $\psi_{v}^{0}$ has value
$\mathbf{L}(\pi_{v})^{-1}$ at $I_{3}$. Multiplying these local factors gives
\begin{align*}
\prod_{v<\infty}\mathcal{H}_{\mathcal{J}_{\psi_{v}}f_{\boldsymbol{\nu},v}}(\sigma_{v},\mathbf{s},\mathbf{1})
=\mathbf{L}_{\fin}(\pi)^{-1}D_{F}^{3s_{1}+s_{2}-3/2+\nu_{3}-\nu_{1}}
\prod_{v<\infty}\mathcal{H}_{W_{v}^{0}}(\sigma_{v},\mathbf{s},\mathbf{1}).
\end{align*}
Multiplying by $2D_{F}^{1/2}\mathbf{L}_{\fin}(\pi)$ and using
\eqref{fc11.42}, \eqref{equ12.11}, and
Lemma~\ref{lem:real-rs-evaluation} gives
\begin{align*}
2D_{F}^{1/2}\mathbf{L}_{\fin}(\pi)
(I_{\mathrm{cusp}}^{\heartsuit}+I_{\mathrm{Eis}}^{\heartsuit})
(\mathbf{s},\varphi_{\boldsymbol{\nu}},\mathbf{1},\mathbf{1})
=\widetilde{\mathcal{I}}_{\mathcal{T}}(\mathbf{s},\boldsymbol{\nu};\mathfrak{q},\mathfrak{n}).
\end{align*}
For the dual term, the additive character change uses $\mathbf{s}^{\vee}$,
with $3s_{1}'+s_{2}'=2s_{2}$. Factor out $q_{v}^{r_{v}(1/2-s_{2}')}$ from
\eqref{eq12.14} at $v\mid\mathfrak{q}$ and $q_{v}^{-n_{v}(1/2+s_{2})}$
from \eqref{eq:dual-level-superposition} at $v\mid\mathfrak{n}$.

Proposition~\ref{prop4.4} writes the original spectral term as
\begin{align*}
I_{\mathrm{spec}}^{\heartsuit}
=I_{\mathrm{cusp}}^{\heartsuit}+I_{\mathrm{Eis}}^{\heartsuit}
+\sum_{j=1}^{2}\sum_{\epsilon\in\{\pm1\}}\epsilon R_{\RNum{2}}^{j,\epsilon}.
\end{align*}
The measure $d\lambda/(4\pi i)$ gives the factor $1/2$ in
\eqref{14.125}. Thus
\begin{equation}\label{14.126}
2D_{F}^{1/2}\mathbf{L}_{\fin}(\pi)I_{\mathrm{spec}}^{\heartsuit}
(\mathbf{s},\varphi_{\boldsymbol{\nu}},\mathbf{1},\mathbf{1})
=\widetilde{\mathcal{I}}_{\mathcal{T}}(\mathbf{s},\boldsymbol{\nu};\mathfrak{q},\mathfrak{n})
+\mathscr{R}_{\mathbf{s},\boldsymbol{\nu}}[D_{\mathcal{T}}(\mathbf{s},\cdot)],
\end{equation}
and 
\begin{equation}\label{14.127}
2D_{F}^{1/2}\mathbf{L}_{\fin}(\pi)I_{\mathrm{spec}}^{\heartsuit}
(\mathbf{s}^{\vee},\pi(w_{2})\varphi_{\boldsymbol{\nu}},\mathbf{1},\mathbf{1})
=\widetilde{\mathcal{I}}_{\mathcal{T}}^{\vee}(\mathbf{s}^{\vee},\boldsymbol{\nu};\mathfrak{q},\mathfrak{n})
+\mathscr{R}_{\mathbf{s}^{\vee},\boldsymbol{\nu}}[D_{\mathcal{T}}^{\vee}(\mathbf{s}^{\vee},\cdot)].
\end{equation}
Theorem~\ref{thmii} now gives \eqref{fc1.3} by
\eqref{e14.123}--\eqref{e14.126}, \eqref{14.126}, and \eqref{14.127}.
\end{proof}

\section{The Newform Sieve and Fifth-Moment Bounds}\label{sec15}
Let $F$ be totally real and $\mathfrak{q},\mathfrak{n}$ coprime
nonzero integral ideals. Use
$\mathbf{s}^{\vee}=((s_{2}-s_{1})/2,(3s_{1}+s_{2})/2)$ and
$\nu_{3}=-\nu_{1}-\nu_{2}$. Fix $0\leq\vartheta_{\infty}<1/2$, depending only
on $F$, so that every real complementary series component of an
automorphic $\sigma$ has real inducing exponents $\pm\eta_{v}$ with
$|\eta_{v}|\leq\vartheta_{\infty}$. Take
$0<\varepsilon<\min\{10^{-3},(1/2-\vartheta_{\infty})/10\}$.
Let $\mathcal{T}=(\mathcal{S}h_{v})_{v\mid\infty}$,
with $h_{v}$ as in \eqref{fc11.3}. For
$*\in\{\mathrm{gen}, \mathrm{degen}, \mathrm{res}\}$, define 
\begin{align*}
\widetilde{\mathcal{I}}_{*,\mathcal{T}}^{\heartsuit}(\mathfrak{q},\mathfrak{n}):=\frac{1}{(2\pi i)^{4}}\oint_{\substack{|s_{j}|=2^{-j}\varepsilon\\|\nu_{j}|=10^{-j}\varepsilon,\ j=1,2}}\frac{\widetilde{\mathcal{I}}_{*,\mathcal{T}}(\mathbf{s},\boldsymbol{\nu};\mathfrak{q},\mathfrak{n})}{s_{1}s_{2}\nu_{1}\nu_{2}}ds_{1}ds_{2}d\nu_{1}d\nu_{2}.
\end{align*}

By local uniform convergence of the spectral terms
(Proposition~\ref{prop4.4}), Cauchy's formula in Theorem~\ref{thmC} gives
\begin{equation}\label{e15.1}
\widetilde{\mathcal{I}}_{\mathcal{T}}^{\heartsuit}(\mathfrak{q},\mathfrak{n})=\widetilde{\mathcal{I}}_{\mathrm{gen},\mathcal{T}}^{\heartsuit}(\mathfrak{q},\mathfrak{n})+\widetilde{\mathcal{I}}_{\mathcal{T}}^{\vee,\heartsuit}(\mathfrak{q},\mathfrak{n})
+\widetilde{\mathcal{I}}_{\mathrm{degen},\mathcal{T}}^{\heartsuit}(\mathfrak{q},\mathfrak{n})
+\widetilde{\mathcal{I}}_{\mathrm{res},\mathcal{T}}^{\heartsuit}(\mathfrak{q},\mathfrak{n}),
\end{equation}
where 
\begin{align*}
&\widetilde{\mathcal{I}}_{\mathcal{T}}^{\heartsuit}(\mathfrak{q},\mathfrak{n}):=\int_{\sigma\in\mathcal{F}^{\circ}(\mathfrak{q};\mathbf{1})}\frac{\lambda_{\sigma}(\mathfrak{n})L(1/2,\sigma)^{4}}{L(1,\sigma,\Ad)}\cdot \mathcal{H}_{\sigma,\mathcal{T}}(\mathbf{0},\mathbf{0};\mathfrak{q},\mathfrak{n})d\mu_{\sigma},\\
&\widetilde{\mathcal{I}}_{\mathcal{T}}^{\vee,\heartsuit}(\mathfrak{q},\mathfrak{n}):=\frac{N_{F}(\mathfrak{q})^{\frac{1}{2}}}{N_{F}(\mathfrak{n})^{\frac{1}{2}}}\int_{\sigma\in\mathcal{F}(\mathfrak{n};\mathbf{1})}\frac{\widetilde{\lambda}_{\sigma,\mathbf{0}}(\mathfrak{q})L(1/2,\sigma)^{4}}{L(1,\sigma,\Ad)}\cdot \mathcal{H}_{\sigma,\mathcal{T}}^{\vee}(\mathbf{0},\mathbf{0};\mathfrak{q},\mathfrak{n})d\mu_{\sigma}.
\end{align*}

\subsection{Local and Global Newform Projectors}
Here $W_{v}^{[m]}$ uses $\psi_{v}^{0}$.

\begin{lemma}\label{lem:depth-kernel}
Let $\pi_{v}$ be irreducible unramified generic and $\sigma_{v}$ irreducible
unitary generic, both with trivial central character.
Assume $r_{\sigma_{v}}\geq2$ and write the Satake parameters of $\pi_{v}$ as
$\alpha_{1},\alpha_{2},\alpha_{3}$, with $\alpha_{1}\alpha_{2}\alpha_{3}=1$.
Define $b_{t}$ by
\begin{equation}\label{eq:depth-inverse}
\sum_{t\geq0}b_{t}z^{t}
=(1+z)^{-1}\prod_{j=1}^{3}(1-\alpha_{j}z)(1-\alpha_{j}^{-1}z).
\end{equation}
Then $|b_{t}|\leq\alpha_{\pi}(\mathfrak{p}_{v})$ and
$|b_{t}|\leq64q_{v}^{\vartheta_{\pi_{v}}\min(t,6)}$.
For every integer $r\geq2$,
\begin{equation}\label{eq:depth-projector}
\sum_{m=2}^{r} b_{r-m}\mathcal{H}_{W_{v}^{[m]}}(\sigma_{v},\mathbf{0},\mathbf{1})
=L_{v}(1,\sigma_{v},\Ad)\mathbf{1}_{r_{\sigma_{v}}=r}.
\end{equation}
\end{lemma}

\begin{proof}
Here
$\lambda_{\sigma_{v}}(\mathfrak{p}_{v}^{n})=\mathbf{1}_{n=0}$ and
$\xi_{\sigma_{v}}(\mathfrak{p}_{v}^{i},\mathfrak{p}_{v}^{l})=\mathbf{1}_{i=l}$.
In \eqref{equ12.11} at $r_{v}=m$, only $k_{1}=k_{2}=l$ and $l_{1}=l_{2}=0$ contribute.
Since $d_{2}\geq r_{\sigma_{v}}\geq2$ and $d_{1}+d_{2}=m$, we have
$q_{v}^{m}\Vol(K_{0,v}'[d_{2}])q_{v}^{-d_{1}}=q_{v}/(q_{v}+1)$.
For $m<r_{\sigma_{v}}$ the $l$-sum is empty. For $m\geq r_{\sigma_{v}}$, multiply
$q_{v}/(q_{v}+1)$ by $|W_{v}'^{\circ}(I_{2})|^{2}$ from
\eqref{eq:finite-newvector-value} to get
\begin{align*}
 \mathcal{H}_{W_{v}^{[m]}}(\sigma_{v},\mathbf{0},\mathbf{1})=L_{v}(1,\sigma_{v},\Ad)
 \sum_{\substack{l,d_{1}\geq0\\l+d_{1}\leq m-r_{\sigma_{v}}}}
 s_{l+d_{1},d_{1},0}(\alpha_{1},\alpha_{2},\alpha_{3}).
\end{align*}
The sum over $m\geq r_{\sigma_{v}}+l+d_{1}$ contributes
$z^{l+d_{1}}/(1-z)$, so \eqref{eq:two-row-schur-generating} at $A=B=z$ gives
\begin{align*}
 \sum_{m\geq r_{\sigma_{v}}}\mathcal{H}_{W_{v}^{[m]}}(\sigma_{v},\mathbf{0},\mathbf{1})z^{m-r_{\sigma_{v}}}
 =
 \frac{(1+z)L_{v}(1,\sigma_{v},\Ad)}{\prod_{j=1}^{3}(1-\alpha_{j}z)(1-\alpha_{j}^{-1}z)}.
\end{align*}
Multiplying by \eqref{eq:depth-inverse} and comparing coefficients proves
\eqref{eq:depth-projector}. The six linear factors in
\eqref{eq:depth-inverse} give $|b_{t}|\leq\alpha_{\pi}(\mathfrak{p}_{v})$ and
$|b_{t}|\leq\sum_{j\leq\min(t,6)}\binom{6}{j}q_{v}^{j\vartheta_{\pi_{v}}}
\leq64q_{v}^{\vartheta_{\pi_{v}}\min(t,6)}$.
\end{proof}

\begin{prop}\label{prop:depth-newform-projection}
Take $\pi=\mathbf{1}\boxplus\mathbf{1}\boxplus\mathbf{1}$ and $\boldsymbol{\nu}=\mathbf{0}$.
Suppose $r\geq2$ and $\mathfrak{n}+\mathfrak{p}=\mathcal{O}_{F}$. Then
\begin{multline}\label{eq:global-depth-projector}
\sum_{m=2}^{r} b_{r-m}\widetilde{\mathcal{I}}_{\mathcal{T}}^{\heartsuit}(\mathfrak{p}^{m},\mathfrak{n})=\int_{\substack{\sigma\in\mathcal{F}^{\circ}(\mathfrak{p}^{r};\mathbf{1})\\
r_{\sigma_{\mathfrak{p}}}=r}}
 \frac{\lambda_{\sigma}(\mathfrak{n})L(1/2,\sigma)^{4}}
 {L^{(\mathfrak{p})}(1,\sigma,\Ad)}
 \mathcal{H}_{\sigma,\mathcal{T}}(\mathbf{0},\mathbf{0};\mathcal{O}_{F},\mathcal{O}_{F})\,d\mu_{\sigma}\\
 {}+O_{F,\mathcal{T},\varepsilon}
 \bigl(N_{F}(\mathfrak{p})^{1+(r+1)\varepsilon}
 N_{F}(\mathfrak{n})^{\vartheta+\varepsilon}\bigr).
\end{multline}
\end{prop}

\begin{proof}
For $r_{\sigma_{\mathfrak{p}}}\geq2$, \eqref{eq:depth-projector}
restricts the integral to $r_{\sigma_{\mathfrak{p}}}=r$ and cancels
the local adjoint factor. For $r_{\sigma_{\mathfrak{p}}}\leq1$,
Lemma~\ref{lem:depth-low} at $\vartheta_{\pi_{\mathfrak{p}}}=0$ bounds
the finite weight. The terms with $r_{\sigma_{\mathfrak{p}}}\leq1$ are in
$\mathcal{F}(\mathfrak{p};\mathbf{1})$. Bound their spectral integral by
\eqref{15.14} at $\mathfrak{n}=\mathfrak{p}$, and use the Hecke bound for
$\lambda_{\sigma}(\mathfrak{n})$.
\end{proof}

\subsection{Proof of Theorem \ref{thm1.5}}\label{sec15.2}

We first bound the three correction terms in \eqref{e15.1}.
Uniformly for nonzero integral ideals $\mathfrak{q},\mathfrak{n}$ with
$\mathfrak{n}+\mathfrak{q}=\mathcal{O}_{F}$, we have
\begin{equation}\label{15.10}
\big|\widetilde{\mathcal{I}}_{\mathrm{gen},\mathcal{T}}^{\heartsuit}(\mathfrak{q},\mathfrak{n})\big|
\ll N_{F}(\mathfrak{q}\mathfrak{n})^{\varepsilon}
 \bigl[N_{F}(\mathfrak{q})N_{F}(\mathfrak{n})^{-1/2}+N_{F}(\mathfrak{n})^{1/2}\bigr],
\end{equation}
and 
\begin{equation}\label{15.11}
\big|\widetilde{\mathcal{I}}_{\mathrm{degen},\mathcal{T}}^{\heartsuit}(\mathfrak{q},\mathfrak{n})\big|
+\big|\widetilde{\mathcal{I}}_{\mathrm{res},\mathcal{T}}^{\heartsuit}(\mathfrak{q},\mathfrak{n})\big|\ll N_{F}(\mathfrak{q})^{\varepsilon}N_{F}(\mathfrak{n})^{\frac{1}{2}+\varepsilon}+N_{F}(\mathfrak{q})^{1+\varepsilon}N_{F}(\mathfrak{n})^{-\frac{1}{2}+\varepsilon}.
\end{equation}

Appendix~\ref{app:generic-degenerate-estimates} proves the generic
and degenerate estimates. Lemma~\ref{lem:fourth-moment-residue-bounds}
bounds the residues.

For spherical $\sigma_{v}=|\cdot|_{v}^{it_{v}}\boxplus|\cdot|_{v}^{-it_{v}}$,
take $H_{R_{v}}$ from the example after \eqref{eq:R-admissible-test},
with the fourth powers replaced by fifth powers.
For real $|t_{v}|>2$ take $R_{v}\leq |t_{v}|<2R_{v}$, $R_{v}=2,4,\ldots$.
For real $|t_{v}|\leq2$ and complementary parameters
$t_{v}=i\eta_{v}$, $|\eta_{v}|\leq\vartheta_{\infty}$, take $R_{v}=2$.
Put $\mathbf{R}=(R_{v})_{v\mid\infty}$ and $X=\prod_{v\mid\infty}R_{v}^{2}$.
By \eqref{eq:G-as-L-ratio}, the Archimedean weight is
\begin{align*}
\mathcal{W}_{\infty,\mathbf{R}}(\sigma)
 =\epsilon_{\sigma} D_{F}^{-3/2}\prod_{v\mid\infty}H_{R_{v}}(t_{v}),
\end{align*}
with $\epsilon_{\sigma}$ from \textsection\ref{def:spectral-multiplicity}.
By \eqref{eq:R-admissible-test}, \cite[Theorem~E]{Yan25b} and the adjoint
lower bound give the cuspidal estimate, and
\cite[Theorem~F]{Yan25b} gives the continuous one. Thus
\begin{equation}\label{eq:original-weighted-fourth}
 \int_{\sigma\in\mathcal{F}^{\circ}(\mathfrak{q};\mathbf{1})}
 \frac{L(1/2,\sigma)^{4}}{L(1,\sigma,\Ad)}
 |\mathcal{W}_{\infty,\mathbf{R}}(\sigma)|\,d\mu_{\sigma}
 \ll_{F,\varepsilon}(N_{F}(\mathfrak{q})X)^{1+\varepsilon}.
\end{equation}

\begin{lemma}\label{lem:hybrid-arch-bookkeeping}
Uniformly for $\mathfrak{q}+\mathfrak{n}=\mathcal{O}_{F}$, we have
\begin{equation}\label{eq:hybrid-dual-spectrum}
 |\widetilde{\mathcal{I}}_{\mathbf{R}}^{\vee,\heartsuit}|
 \ll_{F,\varepsilon}X^{2}N_{F}(\mathfrak{q})^{1/2+\vartheta}
 N_{F}(\mathfrak{n})^{1/2}(N_{F}(\mathfrak{q}\mathfrak{n})X)^{\varepsilon}.
\end{equation}
In particular,
\begin{multline}\label{eq:hybrid-twisted-fourth}
 |\widetilde{\mathcal{I}}_{\mathbf{R}}^{\heartsuit}
 (\mathfrak{q},\mathfrak{n})|
 \ll_{F,\varepsilon}(N_{F}(\mathfrak{q}\mathfrak{n})X)^{\varepsilon}
 N_{F}(\mathfrak{n})^{-1/2}\\
 \bigl[XN_{F}(\mathfrak{q})
 +X^{3/2}\bigl(N_{F}(\mathfrak{n})+N_{F}(\mathfrak{q})\bigr)
 +X^{2}N_{F}(\mathfrak{q})^{1/2+\vartheta}N_{F}(\mathfrak{n})\bigr].
\end{multline}
For $\mathfrak{q}=\mathfrak{p}^{r}$, $r\geq2$, the contribution from
$r_{\sigma_{\mathfrak{p}}}\leq1$ in \eqref{eq:global-depth-projector} is
\begin{equation}\label{eq:hybrid-depth-low-error}
 O_{F,\varepsilon}\bigl(N_{F}(\mathfrak{p})X
 (N_{F}(\mathfrak{q})X)^{\varepsilon}N_{F}(\mathfrak{n})^{\vartheta+\varepsilon}\bigr).
\end{equation}
\end{lemma}

\begin{proof}
After the replacement \eqref{eq:residue-weight-interpolation} at $v\mid\mathfrak{n}$,
the combined dual residues \eqref{eq:residue-cluster} on the fixed
circles \eqref{eq:dual-residue-circles} are holomorphic for
$|\nu_j|\leq10^{-j}\varepsilon$ ($j=1,2$), so their constant term in
$\boldsymbol{\nu}$ is their value at $\boldsymbol{\nu}=\mathbf{0}$.
Combine Lemma~\ref{lem:real-residue-contour-weights} with the original finite
residue estimate in the proof of Lemma~\ref{lem:fourth-moment-residue-bounds},
the dual estimate \eqref{eq:finite-residual-estimates}, and
the degenerate estimate \eqref{eq:finite-degenerate-estimates}.
The same real lemma gives zero for the dual generic term by taking
the $s_{2}$ constant term first on the $\boldsymbol{\nu}$ circles.

In \eqref{15.12} with $h_{v}=h_{v,R}$, the local weight is
$|(\mathscr{T} H_{R})(\sigma_{v})|$. On $U\leq C(\sigma_{v})^{1/2}<2U$,
\eqref{eq:G-transform-tempered-bound} and \eqref{eq:G-transform-discrete-bound}
give the respective bounds for tempered principal and discrete series.
With \cite[Theorems~E and F]{Yan25b}, the local sums are
\begin{align*}
 &\sum_{U=1,2,4,\ldots}U^{2+\varepsilon}
 R^{2}(1+U/R)^{-3}\ll R^{4+\varepsilon},\\
 &\sum_{U=1,2,4,\ldots}U^{2+\varepsilon}
 \min\{1,R^{4}U^{-4}\}\ll R^{2+\varepsilon}
\end{align*}
For complementary series, $C(\sigma_{v})^{1/2}\asymp1$ and
$|(\mathscr{T} H_{R})(\sigma_{v})|\ll R^{2}$ by \eqref{eq:G-transform-complementary-bound}.
Thus the integral is $\ll N_{F}(\mathfrak{n})^{1+\varepsilon}X^{2+\varepsilon}$.
Equation~\eqref{e15.1} gives \eqref{eq:hybrid-twisted-fourth}.

For the error in \eqref{eq:global-depth-projector}, replace \eqref{15.14}
in its proof by \eqref{eq:original-weighted-fourth} at level $\mathfrak{p}$.
\end{proof}

\begin{lemma}\label{lem:hybrid-afe-multiplier}
Put
\begin{align*}
 m_{z}(t):=\pi^{-z}
 \frac{\Gamma(1/4+it/2)\Gamma(1/4-it/2)}
 {\Gamma(1/4-z/2+it/2)\Gamma(1/4-z/2-it/2)}.
\end{align*}
In the definition of $h_{v,R}$ in \textsection\ref{def:polynomial-spherical-test}, replace each
$H_{R_{v}}$ by $H_{R_{v}}m_{z}$, with the fifth-order zeros chosen above,
and multiply the global terms by $e^{z^{2}}$.
For fixed $\Re(z)=c>0$, multiply the bounds in \eqref{eq:original-weighted-fourth}
and Lemma~\ref{lem:hybrid-arch-bookkeeping} by
\begin{equation}\label{eq:afe-multiplier-vertical-bound}
 X^{c/2}(1+|\Im(z)|)^{O(1)}
 e^{-(\Im(z))^{2}+[F:\mathbb{Q}]\pi|\Im(z)|/2},
\end{equation}
where the exponent in $(1+|\Im(z)|)^{O(1)}$ is independent of $X$ and the finite level.
\end{lemma}

\begin{proof}
See \textsection\ref{app:real-weights-moment-estimates}.
\end{proof}

For fixed $\mathcal{T}$, the dual term in \eqref{e15.1} satisfies
the following bound by \eqref{eq14.123} at $v\mid\mathfrak{q}$ and
\eqref{eq:dual-level-superposition} at $v\mid\mathfrak{n}$:
\begin{flalign}\label{15.12}
 && |\widetilde{\mathcal{I}}_{\mathcal{T}}^{\vee,\heartsuit}(\mathfrak{q},\mathfrak{n})|
\ll \frac{N_{F}(\mathfrak{q})^{\frac{1}{2}+\vartheta+\varepsilon}}{N_{F}(\mathfrak{n})^{\frac{1}{2}-\varepsilon}}\int_{\sigma\in\mathcal{F}(\mathfrak{n};\mathbf{1})}
\frac{L(1/2,\sigma)^{4}}{L(1,\sigma,\Ad)}
\prod_{v\mid\infty}
\pi|\mathcal{H}_{W_{v}}^{\vee}(\sigma_{v},\mathbf{0},\mathbf{1})|\,d\mu_{\sigma}.&
\end{flalign}
Here \eqref{eq14.123} gives $|\widetilde{\lambda}_{\sigma,\mathbf{0}}(\mathfrak{q})|
\ll_{\varepsilon} N_{F}(\mathfrak{q})^{\vartheta+\varepsilon}$.
At $\mathbf{s}=\mathbf{0}$, the locally normalized sum in
\eqref{eq:dual-level-superposition} is $O(q_{v}^{-n_{v}/2}(1+n_{v})^{6})$
by \eqref{fc11.42} at $n_{v}=0$ and
Lemmas~\ref{lem:depth-kernel} and~\ref{lem:depth-low} at $m=\ell\geq1$
for $\pi_{v}=\mathbf{1}\boxplus\mathbf{1}\boxplus\mathbf{1}$.
The product of $(1+n_{v})^{6}$ over $v\mid\mathfrak{n}$ is
$O_{\varepsilon}(N_{F}(\mathfrak{n})^{\varepsilon})$.

For $\ell>3$, dyadic summation in $C(\sigma_{v})^{1/2}$ using
\cite[Theorems~E and F]{Yan25b} gives
\begin{equation}\label{15.14}
\int_{\sigma\in\mathcal{F}(\mathfrak{n};\mathbf{1})}
\frac{L(1/2,\sigma)^{4}}{L(1,\sigma,\Ad)}
\prod_{v\mid\infty}C(\sigma_{v})^{-\ell/2}d\mu_{\sigma}
\ll N_{F}(\mathfrak{n})^{1+\varepsilon}.
\end{equation}

By \eqref{eq:uniform-real-dual-decay} at $\mathbf{s}=\mathbf{0}$,
each $\pi|\mathcal{H}_{W_{v}}^{\vee}(\sigma_{v},\mathbf{0},\mathbf{1})|$
in \eqref{15.12} is $O(C(\sigma_{v})^{-\ell/2})$.
Equations \eqref{15.14}, \eqref{15.10}, and \eqref{15.11} in \eqref{e15.1} give
\begin{align*}
 \widetilde{\mathcal{I}}_{\mathcal{T}}^{\heartsuit}(\mathfrak{q},\mathfrak{n})\ll N_{F}(\mathfrak{q})^{1+\varepsilon}N_{F}(\mathfrak{n})^{-\frac{1}{2}+\varepsilon}+N_{F}(\mathfrak{q})^{\frac{1}{2}+\vartheta+\varepsilon}N_{F}(\mathfrak{n})^{\frac{1}{2}+\varepsilon}.
\end{align*}

\begin{lemma}\label{lem:polynomial-newform-sieve}
Let $\mathfrak{a}$ be an integral ideal with $\mathfrak{a}+\mathfrak{n}=\mathcal{O}_{F}$, and let $S$ be a subset of prime
divisors of $\mathfrak{a}$. 
\begin{itemize}
\item Suppose $\mathfrak{a}$ is squarefree. Then 
\begin{multline}\label{eq:marked-exact-fourth-bound}
\Big|
\int_{\substack{\sigma\in\mathcal{F}^{\circ}(\mathfrak{a};\mathbf{1})\\
\mathfrak{c}(\sigma)=\mathfrak{a}}}
 \frac{\lambda_{\sigma}(\mathfrak{n})L(1/2,\sigma)^{4}}
 {L^{(\mathfrak{a})}(1,\sigma,\Ad)}
 \mathcal{W}_{\infty,\mathbf{R}}(\sigma) \prod_{v\in S}q_{v}^{-1/2}\lambda_{\sigma_{v}}(\mathfrak{p}_{v})
\,d\mu_{\sigma}\Big|\\
 \ll_{F,\varepsilon}\frac{(N_{F}(\mathfrak{a}\mathfrak{n})X)^{\varepsilon}}{
 N_{F}(\mathfrak{n})^{1/2}}\cdot \bigl[XN_{F}(\mathfrak{a})
 +X^{3/2}\bigl(N_{F}(\mathfrak{n})+N_{F}(\mathfrak{a})\bigr)
 +X^{2}N_{F}(\mathfrak{a})^{1/2+\vartheta}N_{F}(\mathfrak{n})\bigr].
\end{multline}
\item For $\mathfrak{a}=\mathfrak{p}^{r}$, $r\geq2$, use
\eqref{eq:marked-exact-fourth-bound} with $S=\varnothing$ and include
the error \eqref{eq:hybrid-depth-low-error} at $\mathfrak{q}=\mathfrak{a}$.
\item With $e^{z^{2}}\mathcal{W}_{\infty,\mathbf{R}}(\sigma)
\prod_{v\mid\infty}m_{z}(t_{v})$ in \eqref{eq:marked-exact-fourth-bound},
multiply both bounds by the factor in \eqref{eq:afe-multiplier-vertical-bound}.
\end{itemize}
\end{lemma}

\begin{proof}
The section at $v\mid\mathfrak{a}$ is \eqref{c14.6} with $r_{v}=1$.
At $v\in S$, replacing it by $q_{v}^{-s_{1}}f_{v}^{\mathrm{AL}}$ inserts
$q_{v}^{-1/2}\lambda_{\sigma_{v}}(\mathfrak{p}_{v})$ for $r_{\sigma_{v}}=1$
by \eqref{eq:signed-newvector-factor} at $\mathbf{s}=\mathbf{0}$.
For $r_{\sigma_{v}}=0$, \eqref{eq:signed-old-row} multiplies the
coefficients in \eqref{eq:prime-old-row-series} by $-q_{v}^{-1}$.
For the residues, apply \eqref{eq:residue-weight-interpolation} to
\eqref{eq:signed-old-row-parameters} before taking the constant term.
On the dual side use $\mathbf{s}^{\vee}$ and keep the translate
difference in \eqref{eq:signed-section-dual-identity} unchanged.
This and the generic and degenerate estimates in
\textsection\ref{app:squarefree-signed-section} bound the modified
local factors by $O(q_{v}^{-1+O(\varepsilon)})$ at $v\in S$.

In the dual integral, let $\mathfrak{m}$ be the product of the primes
in $S$ with $r_{\sigma_{v}}=1$. The other primes in $S$ have
$r_{\sigma_{v}}=0$, and the conductor divides $\mathfrak{n}\mathfrak{m}$
by \eqref{eq:signed-section-dual-identity}, \eqref{eq12.14} and
\eqref{eq:dual-level-superposition}. At $\mathbf{s}=\mathbf{0}$,
\eqref{eq:signed-newvector-factor} and Lemma~\ref{lem:prime-old-row-expansion} give
\begin{equation}\label{eq:signed-dual-level-factor}
 \prod_{v\mid\mathfrak{m}}q_{v}^{-1/2}|\lambda_{\sigma_{v}}(\mathfrak{p}_{v})|
 =N_{F}(\mathfrak{m})^{-1}.
\end{equation}
The dyadic estimate for $\prod_{v\mid\infty}|(\mathscr{T} H_{R_{v}})(\sigma_{v})|$
in \eqref{eq:hybrid-dual-spectrum} applies at level $\mathfrak{n}\mathfrak{m}$.
By \eqref{eq:signed-old-row} at $v\in S$ outside $\mathfrak{m}$
and \eqref{eq12.14} at unmarked primes, the dual integral is
\begin{align*}
 \ll X^{2}N_{F}(\mathfrak{n})^{-1/2}(N_{F}(\mathfrak{a}\mathfrak{n})X)^{\varepsilon}
 \prod_{\substack{v\mid\mathfrak{a}\;\,v\notin S}}q_{v}^{1/2+\vartheta}
 \sum_{\mathfrak{m}\mid\prod_{v\in S}\mathfrak{p}_{v}}
 \frac{N_{F}(\mathfrak{n}\mathfrak{m})^{1+\varepsilon}}{N_{F}(\mathfrak{m})}.
\end{align*}
This gives \eqref{eq:hybrid-dual-spectrum} and hence
\eqref{eq:hybrid-twisted-fourth} for the modified sections.

Write $\widetilde{\mathcal{I}}_{\mathbf{R},S}^{\heartsuit}$
for the central moment $\widetilde{\mathcal{I}}_{\mathbf{R}}^{\heartsuit}$
with $q_{v}^{-s_{1}}f_{v}^{\mathrm{AL}}$ at $v\in S$ dividing the level,
where $f_{v}^{\mathrm{AL}}$ is defined in \textsection\ref{app:squarefree-signed-section}.
For $\mathfrak{b}\mid\mathfrak{a}$ and $\mathfrak{l}$ whose prime divisors
all divide $\mathfrak{b}$, put
\begin{equation}\label{eq:squarefree-sieve-coefficients}
 c_{\mathfrak{b},S}(\mathfrak{l})
 :=\prod_{v\mid\mathfrak{b}}
 \begin{cases}
 -c_{v,e_{v}(\mathfrak{l})}(\pi_{v}),&v\notin S,\\
 q_{v}^{-1}c_{v,e_{v}(\mathfrak{l})}(\pi_{v}),&v\in S,
 \end{cases}
\end{equation}
with $c_{v,j}(\pi_{v})$ from \eqref{eq:mixed-prime-old-row-series} at
$\pi_{v}=\mathbf{1}\boxplus\mathbf{1}\boxplus\mathbf{1}$.
At each removed prime, subtract the coefficient of
$\lambda_{\sigma_{v}}(\mathfrak{p}_{v}^{j})$ in \eqref{eq:prime-old-row-series},
with the factor $-q_{v}^{-1}$ from \eqref{eq:signed-old-row} at $v\in S$.
Lemma~\ref{lem:prime-old-row-expansion} gives the factors on the exact-level spectrum, hence
\begin{multline}\label{eq:squarefree-newform-expansion}
 \prod_{v\mid\mathfrak{a}}(1-q_{v}^{-2})
 \int_{\substack{\sigma\in\mathcal{F}^{\circ}(\mathfrak{a};\mathbf{1})\\
 \mathfrak{c}(\sigma)=\mathfrak{a}}}
 \frac{\lambda_{\sigma}(\mathfrak{n})L(1/2,\sigma)^{4}}
 {L^{(\mathfrak{a})}(1,\sigma,\Ad)}
 \mathcal{W}_{\infty,\mathbf{R}}(\sigma)
 \prod_{v\in S}q_{v}^{-1/2}\lambda_{\sigma_{v}}(\mathfrak{p}_{v})\,d\mu_{\sigma}\\
 =\sum_{\mathfrak{b}\mid\mathfrak{a}}
 \sum_{\mathfrak{p}\mid\mathfrak{l}\Rightarrow\mathfrak{p}\mid\mathfrak{b}}
 c_{\mathfrak{b},S}(\mathfrak{l})
 \widetilde{\mathcal{I}}_{\mathbf{R},S}^{\heartsuit}
 (\mathfrak{a}/\mathfrak{b},\mathfrak{n}\mathfrak{l}).
\end{multline}

For $u,\eta>0$ and $m\in\mathbb{Z}_{\geq0}$,
\begin{equation}\label{eq:sieve-supported-Euler-sum}
 \sum_{\mathfrak{b}\mid\mathfrak{a}}\tau(\mathfrak{b})^{3}
 \sum_{\mathfrak{p}\mid\mathfrak{l}\Rightarrow\mathfrak{p}\mid\mathfrak{b}}
\frac{ \tau(\mathfrak{l})^{m}}{N_{F}(\mathfrak{l})^{u}}
 =\prod_{v\mid\mathfrak{a}}\Big(1+8\sum_{j\geq0}\frac{(j+1)^{m}}{q_{v}^{ju}}\Big)
 \ll_{F,u,m,\eta}N_{F}(\mathfrak{a})^{\eta}.
\end{equation}
By \eqref{eq:squarefree-sieve-coefficients} and \eqref{eq:prime-old-row-series},
$|c_{\mathfrak{b},S}(\mathfrak{l})|
\leq\tau(\mathfrak{b})^{3}N_{F}(\mathfrak{l})^{-1/2}$. In
\eqref{eq:sieve-supported-Euler-sum} take $u=\eta$ and $u=1-\eta$,
with $0<\eta<1/2$ and $L\geq1$:
\begin{align*}
 &\sum_{\mathfrak{b}\mid\mathfrak{a}}
 \sum_{\substack{\mathfrak{p}\mid\mathfrak{l}\Rightarrow\mathfrak{p}\mid\mathfrak{b}\\
 N_{F}(\mathfrak{l})\leq L}}
 |c_{\mathfrak{b},S}(\mathfrak{l})|N_{F}(\mathfrak{l})^{1/2+\eta}
 \ll_{F,\eta}N_{F}(\mathfrak{a})^{\eta}L^{2\eta},\\
 &\sum_{\mathfrak{b}\mid\mathfrak{a}}
 \sum_{\mathfrak{p}\mid\mathfrak{l}\Rightarrow\mathfrak{p}\mid\mathfrak{b}}
 |c_{\mathfrak{b},S}(\mathfrak{l})|N_{F}(\mathfrak{l})^{-1/2+\eta}
 \ll_{F,\eta}N_{F}(\mathfrak{a})^{\eta}.
\end{align*}
In \eqref{eq:squarefree-newform-expansion}, \eqref{eq:hybrid-twisted-fourth}
bounds each integral at level $\mathfrak{a}/\mathfrak{b}$ and twist
$\mathfrak{n}\mathfrak{l}$. The two coefficient sums above bound the part
with $N_{F}(\mathfrak{l})\leq L$ by $L^{2\eta}$ times
\eqref{eq:hybrid-twisted-fourth} at $\mathfrak{q}=\mathfrak{a}$, $\varepsilon=2\eta$.

For $N_{F}(\mathfrak{l})>L$, \eqref{eq:original-weighted-fourth} bounds
the integral at level $\mathfrak{a}/\mathfrak{b}$ after taking out the
Hecke bound and the $O(1)$ local weights from
Lemma~\ref{lem:prime-old-row-expansion}, \eqref{eq:signed-newvector-factor}
and \eqref{eq:signed-old-row}. The tail then follows from
\eqref{eq:sieve-supported-Euler-sum} at $(m,u)=(1,(1/2-\vartheta)/2)$:
\begin{align*}
 \sum_{\mathfrak{b}\mid\mathfrak{a}}
 \sum_{\substack{\mathfrak{p}\mid\mathfrak{l}\Rightarrow\mathfrak{p}\mid\mathfrak{b}\\
 N_{F}(\mathfrak{l})>L}}
 \big|c_{\mathfrak{b},S}(\mathfrak{l})\widetilde{\mathcal{I}}_{\mathbf{R},S}^{\heartsuit}
 (\mathfrak{a}/\mathfrak{b},\mathfrak{n}\mathfrak{l})\big|
 \ll_{F,\eta}(N_{F}(\mathfrak{a})X)^{1+\eta}
 N_{F}(\mathfrak{n})^{\vartheta+\eta}
 L^{-(1/2-\vartheta)/2}.
\end{align*}
Choose $\ell(1/2-\vartheta)>2(\vartheta+1)$, $L=(N_{F}(\mathfrak{a}\mathfrak{n})X)^{\ell}$ and $0<\eta<\min\big\{\frac{1}{2},\frac{\varepsilon}{2+2\ell}\big\}$. Then the tail is $O(XN_{F}(\mathfrak{a})N_{F}(\mathfrak{n})^{-1/2})$,
and the part of \eqref{eq:squarefree-newform-expansion} with
$N_{F}(\mathfrak{l})\leq L$ satisfies \eqref{eq:marked-exact-fourth-bound}.
For $\mathfrak{a}=\mathfrak{p}^{r}$, \eqref{eq:hybrid-twisted-fourth}
bounds each moment at level $\mathfrak{p}^{m}$ in \eqref{eq:global-depth-projector},
and $\sum_{m=2}^{r}|b_{r-m}|\leq64r$. The error is
\eqref{eq:hybrid-depth-low-error}. Lemma~\ref{lem:hybrid-afe-multiplier}
gives the additional factor for both estimates.
\end{proof}

Fix $\mathfrak{a}\mid\mathfrak{q}$ with $\mathfrak{a}\ne\mathcal{O}_{F}$
and consider the representations with $\mathfrak{c}(\sigma)=\mathfrak{a}$.
An Eisenstein representation of this conductor has ramified inducing
characters, so its completed standard $L$-function is entire.
Use the even factor $e^{z^{2}}\prod_{v\mid\infty}m_{z}(t_{v})m_{-z}(t_{v})$
in the symmetric approximate functional equation.
On $\Re(z)=c$ with $c>1/2+\vartheta$, the Gamma ratios cancel
$m_{-z}(t_{v})$, giving
\begin{multline}\label{eq:fifth-moment-root-number-removal}
 L(1/2,\sigma)^{5}
 =\frac{2L(1/2,\sigma)^{4}}{2\pi i}
 \int_{\Re(z)=c}e^{z^{2}}(D_{F}^{2}N_{F}(\mathfrak{a}))^{z/2}\\
 \prod_{v\mid\infty}m_{z}(t_{v})
 \sum_{\mathfrak{n}}\frac{\lambda_{\sigma}(\mathfrak{n})}
 {N_{F}(\mathfrak{n})^{1/2+z}}\frac{dz}{z}.
\end{multline}
We integrate \eqref{eq:fifth-moment-root-number-removal} against
$\mathcal{W}_{\infty,\mathbf{R}}(\sigma)d\mu_{\sigma}/L^{(\mathfrak{a})}(1,\sigma,\Ad)$
over $\mathfrak{c}(\sigma)=\mathfrak{a}$ in $\mathcal{F}^{\circ}(\mathfrak{a};\mathbf{1})$.
Absolute convergence follows from the Hecke bound for $c>1/2+\vartheta$ and
\eqref{eq:original-weighted-fourth} with \eqref{eq:afe-multiplier-vertical-bound}.
The omitted adjoint factors satisfy
$\prod_{v\mid\mathfrak{a}}L_{v}(1,\sigma_{v},\Ad)\ll_{\varepsilon} N_{F}(\mathfrak{a})^{\varepsilon}$.

For $v\mid\mathfrak{a}$, $L_{v}(1/2+z,\sigma_{v})=1$ if $r_{\sigma_{v}}\geq2$.
If $r_{\sigma_{v}}=1$, then for $\Re(z)>0$,
\begin{equation}\label{eq:bad-prime-AFE-sign-expansion}
 L_{v}(1/2+z,\sigma_{v})
 =\frac{1+q_{v}^{-1/2-z}\lambda_{\sigma_{v}}(\mathfrak{p}_{v})}{1-q_{v}^{-2-2z}}.
\end{equation}
For squarefree $\mathfrak{a}$, expanding the product of
\eqref{eq:bad-prime-AFE-sign-expansion} gives the integrals in
\eqref{eq:marked-exact-fourth-bound} with real weight
$e^{z^{2}}\mathcal{W}_{\infty,\mathbf{R}}(\sigma)\prod_{v\mid\infty}m_{z}(t_{v})$.
For fixed $\Re(z)>0$, their coefficients satisfy
\begin{equation}\label{eq:bad-prime-AFE-subset-bound}
 \sum_{S\subseteq\{v:v\mid\mathfrak{a}\}}
 \left|\frac{\prod_{v\in S}q_{v}^{-z}}{\prod_{v\mid\mathfrak{a}}(1-q_{v}^{-2-2z})}\right|
=\prod_{v\mid\mathfrak{a}}\frac{1+q_{v}^{-\Re(z)}}{|1-q_{v}^{-2-2z}|}
 \ll_{F,\varepsilon}N_{F}(\mathfrak{a})^{\varepsilon}.
\end{equation}
For $\mathfrak{a}=\mathfrak{p}^{r}$, $r\geq2$, the Euler factor is $1$ and
Lemma~\ref{lem:polynomial-newform-sieve} gives the integral bound with $S=\varnothing$.

For the terms with $N_{F}(\mathfrak{n})\asymp l$ and
$\mathfrak{n}+\mathfrak{a}=\mathcal{O}_{F}$, we move the $z$-integral in
\eqref{eq:fifth-moment-root-number-removal} to $\Re(z)=c>0$.
The conductor and $\mathfrak{n}$ powers give
$N_{F}(\mathfrak{n})^{-1/2-c}N_{F}(\mathfrak{a})^{c/2}$ up to a factor depending on $F,c$. The integral in $\sigma$ is bounded by
\eqref{eq:marked-exact-fourth-bound} times \eqref{eq:afe-multiplier-vertical-bound}.
After integrating in $\Im(z)$ and summing over the $O_{F}(l)$ ideals,
the contribution of this range is
\begin{multline}\label{eq:fifth-moment-dyadic-bound}
(N_{F}(\mathfrak{a})Xl)^{\varepsilon}
 \left(\frac{(N_{F}(\mathfrak{a})X)^{1/2}}{l}\right)^{c}
 \\
\bigl[XN_{F}(\mathfrak{a})
 +X^{3/2}(N_{F}(\mathfrak{a})+l)
 +X^{2}N_{F}(\mathfrak{a})^{1/2+\vartheta}l\bigr].
\end{multline}
Summing \eqref{eq:fifth-moment-dyadic-bound} with $c=\varepsilon$ up to
$l=(N_{F}(\mathfrak{a})X)^{1/2}$ and $c=3$ above gives
$O_{F,\varepsilon}(N_{F}(\mathfrak{a})^{1+\vartheta+\varepsilon}
\allowbreak X^{5/2+\varepsilon})$. At a tempered principal-series scale
$R_{v}$, \eqref{eq:G-transform-tempered-bound} bounds the dual transform by
$O(R_{v}^{2})$ without using oscillation in the cosine integral. The
Archimedean factor at $v$ in the fourth-moment bound is
$O_{\varepsilon}(R_{v}^{2+\varepsilon})$.
Their product over $v\mid\infty$ gives $X^{2}$ in \eqref{eq:hybrid-dual-spectrum},
up to $X^{\varepsilon}$. The factor $l$ in the last term of
\eqref{eq:fifth-moment-dyadic-bound} gives the additional $X^{1/2}$ at
$l\asymp(N_{F}(\mathfrak{a})X)^{1/2}$.

For $\mathfrak{a}=\mathfrak{p}^{r}$, $r\geq2$, the error
\eqref{eq:hybrid-depth-low-error} contributes
$N_{F}(\mathfrak{p})Xl^{1/2+\vartheta}$ inside the bracket of
\eqref{eq:fifth-moment-dyadic-bound}. Summing over $l$ gives
$O_{F,\varepsilon}(N_{F}(\mathfrak{p})X(N_{F}(\mathfrak{a})X)^{1/4+\vartheta/2+\varepsilon})$,
which is smaller than the preceding bound since
$N_{F}(\mathfrak{p})\leq N_{F}(\mathfrak{a})^{1/2}$.

Since $\mathcal{W}_{\infty,\mathbf{R}}(\sigma)L(1/2,\sigma)^{5}\geq0$,
keep only the cuspidal terms. On the chosen ranges and for nonzero
central values, $\mathcal{W}_{\infty,\mathbf{R}}\gg_{F}1$ and
$L^{(\mathfrak{a})}(1,\sigma,\Ad)\ll_{F,\varepsilon}(N_{F}(\mathfrak{a})X)^{\varepsilon}$
by the standard adjoint bound. Hence
\begin{align*}
 \sum_{\substack{\sigma\in\mathcal{F}_{\mathrm{cusp}}^{\circ}
 (\mathfrak{a};\mathbf{1})\\ \mathfrak{c}(\sigma)=\mathfrak{a}\\
 t_{v}\ \mathrm{in\ its\ chosen\ range}\ (v\mid\infty)}}
 L(1/2,\sigma)^{5}\ll_{F,\varepsilon}
 N_{F}(\mathfrak{a})^{1+\vartheta+\varepsilon}X^{5/2+\varepsilon}.
\end{align*}
At level one, $C_{\infty}(\sigma)\asymp X$ gives
$L(1/2,\sigma)\ll_{F,\varepsilon}X^{1/4+\varepsilon}$ by convexity.
The fourth moment is bounded by \cite[Theorem~E]{Yan25b} at level one
and $C(\sigma_{v})\ll R_{v}^{2}$.

Since $C_{\infty}(\sigma)\asymp X$, it suffices to sum over $X\ll_{F} T$. Sum over
$\mathfrak{a}\mid\mathfrak{q}$ and the dyadic choices of $\mathbf{R}$,
absorbing both counts into $(N_{F}(\mathfrak{q})T)^{\varepsilon}$.

\begin{lemma}\label{lem:fixed-level-fifth-moment}
Fix an integral ideal $\mathfrak{f}$ and let $T>10$.
If $\mathfrak{q}$ is squarefree and $\mathfrak{f}+\mathfrak{q}=\mathcal{O}_{F}$, then
\begin{align*}
 \sum_{\substack{\sigma\in\mathcal{F}_{\mathrm{cusp}}^{\circ}
 (\mathfrak{f}\mathfrak{q};\mathbf{1})\\ C_{\infty}(\sigma)\leq T}}
 L(1/2,\sigma)^{5}
 \ll_{F,\mathfrak{f},\varepsilon}
 N_{F}(\mathfrak{q})^{1+\vartheta+\varepsilon}T^{5/2+\varepsilon}.
\end{align*}
\end{lemma}

\begin{proof}
Fix $\mathfrak{a}=\mathfrak{f}'\mathfrak{q}'\ne\mathcal{O}_{F}$ with
$\mathfrak{f}'\mid\mathfrak{f}$ and $\mathfrak{q}'\mid\mathfrak{q}$, and
consider the representations with $\mathfrak{c}(\sigma)=\mathfrak{a}$.
At $v\mid\mathfrak{f}'$ with $e_{v}(\mathfrak{f}')\geq2$,
\eqref{eq:depth-projector} isolates $r_{\sigma_{v}}=e_{v}(\mathfrak{f}')$
among $r_{\sigma_{v}}\geq2$. The contributions with $r_{\sigma_{v}}=0,1$
are subtracted as follows.
For $v\mid\mathfrak{f}'$ and level exponent $m$,
let $\sigma_{\beta}$ have Satake parameters $\beta,\beta^{-1}$.
Equation~\eqref{equ12.11} at $r_{v}=m$ and the coefficients in
\eqref{eq:oldform-basis-coefficients} show that the local weight
$\mathcal{H}_{W_{v}^{[m]}}(\sigma_{\beta},\mathbf{0},\mathbf{1})/L_{v}(1/2,\sigma_{\beta})^{4}$ is symmetric in
$\beta,\beta^{-1}$ and holomorphic on
$N_{F}(\mathfrak{p}_{v})^{-1/2}<|\beta|<N_{F}(\mathfrak{p}_{v})^{1/2}$.
For $0<\delta<1/2-\vartheta$, Cauchy's estimate on
$|\beta|=N_{F}(\mathfrak{p}_{v})^{\pm(1/2-\delta)}$ and the Hecke
recurrence give, for
$N_{F}(\mathfrak{p}_{v})^{-\vartheta}\leq|\beta|\leq N_{F}(\mathfrak{p}_{v})^{\vartheta}$,
\begin{align*}
\frac{\mathcal{H}_{W_{v}^{[m]}}(\sigma_{\beta},\mathbf{0},\mathbf{1})}{L_{v}(1/2,\sigma_{\beta})^{4}}
 =\sum_{j\geq0}c_{\mathfrak{p}_{v},j}\lambda_{\sigma_{\beta}}(\mathfrak{p}_{v}^{j}),
 \qquad |c_{\mathfrak{p}_{v},j}|\ll_{\mathfrak{f},\delta}
 N_{F}(\mathfrak{p}_{v})^{-(1/2-\delta)j}.
\end{align*}
For $r_{\sigma_{v}}=0$, the term indexed by $j$ is a moment with
$\mathfrak{p}_{v}$ removed from the level and twist $\mathfrak{n}\mathfrak{p}_{v}^{j}$.
For $r_{\sigma_{v}}=1$, expand the weight in $1$ and
$\lambda_{\sigma_{v}}(\mathfrak{p}_{v})$. The constant term comes from
$f_{\boldsymbol{\nu},v}$ at $r_{v}=1$ in \eqref{c14.6}, and
$q_{v}^{1/2-s_{1}}f_{v}^{\mathrm{AL}}$ inserts $\lambda_{\sigma_{v}}(\mathfrak{p}_{v})$ by
\eqref{eq:signed-newvector-factor}.
By \eqref{eq:signed-old-row} at $\mathbf{s}=\mathbf{0}$, the normalized
$r_{\sigma_{v}}=0$ weights of $f_{\boldsymbol{\nu},v}$ and $q_{v}^{1/2-s_{1}}f_{v}^{\mathrm{AL}}$
are \eqref{eq:prime-old-row-series} and $-q_{v}^{-1/2}$ times that series.
We subtract each term at the level with $\mathfrak{p}_{v}$ removed and
twist $\mathfrak{n}\mathfrak{p}_{v}^{j}$. The two weights on
$r_{\sigma_{v}}=1$ are $(1-q_{v}^{-2})L_{v}(1,\sigma_{v},\Ad)$
times $1$ and $\lambda_{\sigma_{v}}(\mathfrak{p}_{v})$, respectively.
Dividing by $1-q_{v}^{-2}$ gives the exact exponent $1$ weights.
These also subtract the $r_{\sigma_{v}}=1$ contributions when $e_{v}(\mathfrak{f}')\geq2$.

Multiplying these local identities by the factors in
\eqref{eq:squarefree-newform-expansion} at $\mathfrak{q}'$ gives moments at
level $\mathfrak{d}\mathfrak{q}'/\mathfrak{b}$, with
$\mathfrak{d}\mid\mathfrak{f}'$ and $\mathfrak{b}\mid\mathfrak{q}'$.
The extra twist $\mathfrak{l}$ is supported on the primes removed from the
level, so $\mathfrak{n}\mathfrak{l}$ is coprime to
$\mathfrak{d}\mathfrak{q}'/\mathfrak{b}$. The bounds for the signed sections
in the proof of Lemma~\ref{lem:polynomial-newform-sieve} apply here as well,
since the marked primes still have level exponent $1$ and the local factors
at the other primes are unchanged.

By the coefficient bound at $\mathfrak{f}'$ and
\eqref{eq:sieve-supported-Euler-sum} with $u=\eta$ at $\mathfrak{q}'$,
the sum of absolute coefficients weighted by $N_{F}(\mathfrak{l})^{1/2+\eta}$
over $N_{F}(\mathfrak{l})\leq L$ is
$O_{\mathfrak{f},\delta,\eta}(N_{F}(\mathfrak{q}')^{\eta}L^{\delta+4\eta})$.
For the tail, put $\kappa=(1/2-\vartheta-\delta)/2$. The Hecke bound
and \eqref{eq:original-weighted-fourth} bound the tail by
\begin{align*}
 O_{F,\mathfrak{f},\delta,\eta}\big((N_{F}(\mathfrak{q}')X)^{1+2\eta}
 N_{F}(\mathfrak{n})^{\vartheta+\eta}L^{-\kappa}\big).
\end{align*}
For $N_{F}(\mathfrak{l})>L$, the inequality
$1\leq(N_{F}(\mathfrak{l})/L)^{\kappa}$ gives $L^{-\kappa}$ times the
coefficient sum with weight $\tau(\mathfrak{l})N_{F}(\mathfrak{l})^{\vartheta+\kappa}$.
It converges at the fixed primes and is
$O_{F,\delta,\eta}(N_{F}(\mathfrak{q}')^{\eta})$ at the other primes by
\eqref{eq:sieve-supported-Euler-sum} with $(u,m)=(1/2-\vartheta-\kappa,1)$.
With $L=(N_{F}(\mathfrak{q}'\mathfrak{n})X)^{\ell}$, choose
\begin{align*}
\ell\kappa>\vartheta+1,
 \qquad \ell(\delta+4\eta)+2\eta<\varepsilon.
\end{align*}
For $N_{F}(\mathfrak{l})\leq L$, apply \eqref{eq:hybrid-twisted-fourth}
at $(\mathfrak{d}\mathfrak{q}'/\mathfrak{b},\mathfrak{n}\mathfrak{l})$.
Together with the tail bound, this gives \eqref{eq:marked-exact-fourth-bound}
with $\mathfrak{a}$ replaced by $\mathfrak{q}'$ and an implied constant
depending also on $\mathfrak{f}$. Lemma~\ref{lem:hybrid-afe-multiplier}
supplies \eqref{eq:afe-multiplier-vertical-bound} for both bounds.

In \eqref{eq:fifth-moment-root-number-removal}, the conductor
$N_{F}(\mathfrak{f}')N_{F}(\mathfrak{q}')$ contributes only the additional
fixed factor $N_{F}(\mathfrak{f}')^{z/2}$.
For $v\mid\mathfrak{f}$ with $r_{\sigma_{v}}=1$, the term
$q_{v}^{-1/2-z}\lambda_{\sigma_{v}}(\mathfrak{p}_{v})$ in \eqref{eq:bad-prime-AFE-sign-expansion} uses
\eqref{eq:signed-newvector-factor} with the multiplier $q_{v}^{-s_{1}-z}$.
The dyadic contribution is therefore bounded by
\eqref{eq:fifth-moment-dyadic-bound} with $\mathfrak{q}'$ in place of
$\mathfrak{a}$, with constants depending on $\mathfrak{f}$.
Sum over $\mathfrak{f}'\mid\mathfrak{f}$ and $\mathfrak{q}'\mid\mathfrak{q}$,
including $O_{F,\varepsilon}(X^{5/4+\varepsilon})$ for each $\mathbf{R}$
at $\mathfrak{a}=\mathcal{O}_{F}$.
\end{proof}

For cuspidal $\sigma$ with spherical Archimedean components, trivial
central character and conductor dividing $\mathfrak{f}\mathfrak{q}$, where $\mathfrak{f}$ is fixed,
$\mathfrak{q}$ is squarefree, and $\mathfrak{f}+\mathfrak{q}=\mathcal{O}_{F}$,
Lemma~\ref{lem:fixed-level-fifth-moment} with
$T=10(1+C_{\infty}(\sigma))$ gives
\begin{equation}\label{eq:individual-hybrid-fifth}
 L(1/2,\sigma)\ll_{F,\mathfrak{f},\varepsilon}
 N_{F}(\mathfrak{q})^{(1+\vartheta)/5+\varepsilon}
 C_{\infty}(\sigma)^{1/2+\varepsilon}.
\end{equation}

\subsection{Mixed Moments for Cuspidal \texorpdfstring{$\mathrm{GL}_{3}$}{GL3} Data}\label{sec15.mixed}
Fix a self-dual cuspidal representation $\pi$ of $\mathrm{PGL}_{3}/F$,
unramified at all places. At real $v$, normalize the spherical Whittaker
vector $W_{\pi_{v}}^{\circ}$ by taking $\beta_{v}>0$ in
\eqref{eq:real-jacquet-double-mellin}.

Put $\mathcal{V}_{v}^{\circ}(\sigma_{v})
=\pi\mathcal{H}_{W_{\pi_{v}}^{\circ}}(\sigma_{v},\mathbf{0},\mathbf{1})$.

\begin{lemma}\label{lem:mixed-spherical-arch}
At every real place,
$\pi\mathcal{H}_{W_{\pi_{v}}^{\circ}}^{\vee}(\sigma_{v},\mathbf{0},\mathbf{1})
=\mathcal{V}_{v}^{\circ}(\sigma_{v})$.
For automorphic $\sigma$, this weight vanishes unless $\sigma_{v}$ is
spherical. If $\sigma_{v}$ is spherical, then
$0<\mathcal{V}_{v}^{\circ}(\sigma_{v})\ll_{\pi_{v},A}C(\sigma_{v})^{-A}$
for every $A>0$.
\end{lemma}

\begin{proof}
Since $w_{2}$ fixes $W_{\pi_{v}}^{\circ}$, \eqref{fc10.2} gives equality.
In \eqref{fc10.1}, only spherical vectors contribute.
In \eqref{eq:real-weight-evaluated} at $\mathbf{s}=\mathbf{0}$ and
$\mathcal{S} h_{v}=1$, the real Gamma arguments are positive since
$\vartheta_{\pi_{v}}+\vartheta_{\sigma_{v}}<1/2$ by $\pi=\Ad\tau$
\cite[Theorem A]{Ram14} and \cite{KS02}. The other Gamma factors pair
by complex conjugation. Stirling's formula gives the decay.
\end{proof}

\phantomsection\label{def:mixed-real-weight}
The mixed moment has real weight
\begin{align*}
 \mathcal{W}_{\infty}(\sigma)
 :=\epsilon_{\sigma}D_{F}^{-3/2}
 \prod_{v\mid\infty}\mathcal{V}_{v}^{\circ}(\sigma_{v}).
\end{align*}
\phantomsection\label{def:mixed-finite-weight}
For $m\geq0$, put
\begin{align*}
 \alpha_{v}(\sigma_{v};m):=
 \frac{\mathcal{H}_{W_{v}^{[m]}}(\sigma_{v},\mathbf{0},\mathbf{1})}
 {L_{v}(1/2,\pi_{v}\times\sigma_{v})L_{v}(1/2,\sigma_{v})}.
\end{align*}
In particular, $\alpha_{v}(\sigma_{v};0)=\mathbf{1}_{r_{\sigma_{v}}=0}$.
For $\mathfrak{a}+\mathfrak{n}=\mathcal{O}_{F}$, define
\begin{align*}
 \mathcal{M}_{\pi}(\mathfrak{a},\mathfrak{n})
 :=\int_{\sigma\in\mathcal{F}^{\circ}(\mathfrak{a};\mathbf{1})}
 \frac{\lambda_{\sigma}(\mathfrak{n})
 L(1/2,\pi\times\sigma)L(1/2,\sigma)\mathcal{W}_{\infty}(\sigma)}{L(1,\sigma,\Ad)}
 \prod_{v\mid\mathfrak{a}}\alpha_{v}(\sigma_{v};e_{v}(\mathfrak{a}))\,d\mu_{\sigma}.
\end{align*}
With $\mathfrak{q}=\mathfrak{a}$ in \eqref{8.1}, write
$\mathcal{M}_{\pi}^{\vee}(\mathfrak{a},\mathfrak{n})$ for the dual spectral
integral at $\mathbf{s}=\mathbf{0}$ in the normalization of $\mathcal{M}_{\pi}$.
Cuspidality leaves only the standard residues $R_{\RNum{2}}^{2,\pm}$.
By \eqref{eq:opposite-spectral-residues} their signed sum on each side is
twice $R_{\RNum{2}}^{2,+}$, taken at $\lambda=1/2-s_{2}$ or
$1/2-s_{2}'$ before setting $\mathbf{s}=\mathbf{0}$.

\phantomsection\label{def:mixed-residue-coefficient}
Define $g_{\pi}$ multiplicatively by
\begin{align*}
 g_{\pi}(\mathfrak{p}_{v}^{r})=\frac{1}{1-q_{v}^{-2}}
 \sum_{m=0}^{r}\mathcal{N}(\mathbf{0},\pi_{v};m),
\end{align*}
with $\mathcal{N}$ from \eqref{14.36}. Summing \eqref{14.36} at
$A=B=q_{v}^{-1}$ and using self-duality and $h_{m}=\lambda_{\pi}(\mathfrak{p}_{v}^{m})$
by \eqref{12.6} at $j=0$ gives
\begin{equation}\label{eq:mixed-generic-coefficient}
 \sum_{r\geq0}g_{\pi}(\mathfrak{p}_{v}^{r})u^{r}
 =\frac{1-(1+q_{v})^{-1}(1+\lambda_{\pi_{v}}(\mathfrak{p}_{v}))u
 +q_{v}^{-1}u^{2}}
 {1-u}\sum_{r\geq0}\lambda_{\pi}(\mathfrak{p}_{v}^{r})u^{r}.
\end{equation}
For a multiplicative majorant, put
\begin{align*}
 \sum_{r\geq0}g_{\pi}^{+}(\mathfrak{p}_{v}^{r})u^{r}
 =\frac{1+(1+q_{v})^{-1}|1+\lambda_{\pi_{v}}(\mathfrak{p}_{v})|u
 +q_{v}^{-1}u^{2}}{1-u}
 \sum_{r\geq0}|\lambda_{\pi}(\mathfrak{p}_{v}^{r})|u^{r}.
\end{align*}
By \eqref{eq:mixed-generic-coefficient}, $g_{\pi}^{+}\geq\max(1,|g_{\pi}|)$.
For every $\eta>0$,
\begin{equation}\label{eq:mixed-generic-average}
 \sum_{\mathfrak{n}}\frac{g_{\pi}^{+}(\mathfrak{n})}
 {N_{F}(\mathfrak{n})^{1+\eta}}\ll_{\pi,F,\eta}1.
\end{equation}
Cauchy--Schwarz and the Rankin--Selberg series for
$\pi\times\widetilde{\pi}$ give absolute convergence of
$\sum|\lambda_{\pi}(\mathfrak{n})|N_{F}(\mathfrak{n})^{-1-\eta}$.
At $u=q_{v}^{-1-\eta}$, the fraction multiplying
$\sum_{r\geq0}|\lambda_{\pi}(\mathfrak{p}_{v}^{r})|u^{r}$ is
$1+O(q_{v}^{-1-\eta})$, since
$|1+\lambda_{\pi_{v}}(\mathfrak{p}_{v})|/(1+q_{v})\ll q_{v}^{\vartheta_{\pi_{v}}-1}$.
Its product over $v$ converges absolutely.

\begin{lemma}\label{lem:mixed-full-spectrum-mass}
Let $\tau$ be a fixed unitary cuspidal representation of
$\mathrm{GL}_{2}/F$ and put $\pi=\Ad\tau$. Assume that $\pi$ is
unramified at all places; $\tau$ may be dihedral.
There is $\ell_{0}=\ell_{0}(\tau,F)>0$ such that, for $\ell>\ell_{0}$
and every nonzero integral ideal $\mathfrak{b}$,
\begin{equation}\label{eq:mixed-positive-full-spectrum-required}
 \int_{\sigma\in\mathcal{F}(\mathfrak{b};\mathbf{1})}
 \frac{L(1/2,\pi\times\sigma)L(1/2,\sigma)}{L(1,\sigma,\Ad)}
 \prod_{v\mid\infty}C(\sigma_{v})^{-\ell/2}\,d\mu_{\sigma}
 \ll_{\tau,F,\ell,\varepsilon}N_{F}(\mathfrak{b})^{1+\varepsilon}.
\end{equation}
\end{lemma}

\begin{proof}
In \eqref{8.1}, take trivial characters and $\mathfrak{n}=\mathcal{O}_{F}$.
For exact conductor $\mathfrak{a}$ and $C(\sigma_{v})^{1/2}\leq R_{v}$,
$R_{v}\geq2$, take the data of \textsection\ref{app:mixed-mass-nonspectral}.
By positivity as in \textsection\ref{sec13.2}, the local lower bounds
give $N_{F}(\mathfrak{a})^{-\varepsilon}$ times the integral in
\eqref{eq:mixed-positive-box-bound} as a lower bound for the original
spectral integral, up to an implied constant.

At $(r_{v},\mathbf{s})=(m,\mathbf{0})$, the first line of \eqref{eq12.14}
is $O(q_{v}^{m})$ and the local weight vanishes unless $r_{\sigma_{v}}=0$.
At real $v$, \eqref{eq:mixed-full-K-sobolev-bound} and
\eqref{eq:mixed-fixed-Whittaker-bound}, with $\Delta_{v}$ transferred to
$\pi_{v}(w_{2})W_{v}$, give
$|\mathcal{H}_{W_{v}}^{\vee}(\sigma_{v},\mathbf{0},\mathbf{1})|\allowbreak
\ll_{A} R_{v}^{O_{A}(1)}C(\sigma_{v})^{-A}$.
Convexity and the level-one Weyl bound give
$O(N_{F}(\mathfrak{a})^{1+\varepsilon}\allowbreak\prod_{v\mid\infty}R_{v}^{O(1)})$ for the dual spectral integral.
Section~\ref{app:mixed-mass-nonspectral} and
Lemma~\ref{lem:mixed-level-one-residue-bound} give this bound for the remaining terms of \eqref{8.1}. Thus
\begin{equation}\label{eq:mixed-positive-box-bound}
 \int_{\substack{\sigma\in\mathcal{F}(\mathfrak{a};\mathbf{1}),\ \mathfrak{c}(\sigma)=\mathfrak{a}\\
 C(\sigma_{v})^{1/2}\leq R_{v},\ v\mid\infty}}
 \frac{L(1/2,\pi\times\sigma)L(1/2,\sigma)}{L(1,\sigma,\Ad)}
 \,d\mu_{\sigma}
 \ll_{\pi,F,\varepsilon} N_{F}(\mathfrak{a})^{1+\varepsilon}
 \prod_{v\mid\infty}R_{v}^{O(1)}.
\end{equation}
The exponents of $R_{v}$ are independent of $\varepsilon$. Sum over
$\mathfrak{a}\mid\mathfrak{b}$ and dyadic ranges of $C(\sigma_{v})^{1/2}$ at each real place.
\end{proof}

Under the hypotheses of Lemma~\ref{lem:mixed-full-spectrum-mass},
for sufficiently large $\ell$, we have
\begin{multline}\label{eq:mixed-full-spectrum-finite-weights}
 \int_{\mathcal{F}(\mathfrak{b};\mathbf{1})}
 \frac{L(1/2,\pi\times\sigma)L(1/2,\sigma)}{L(1,\sigma,\Ad)}
 \prod_{v\mid\mathfrak{b}}|\alpha_{v}(\sigma_{v};e_{v}(\mathfrak{b}))|\\
 \prod_{v\mid\infty}C(\sigma_{v})^{-\ell/2}\,d\mu_{\sigma}
 \ll N_{F}(\mathfrak{b})^{1+\varepsilon}.
\end{multline}
For $m\geq1$ and $\eta>0$,
\begin{align*}
 \sum_{k=0}^{m} q_{v}^{k(1+\eta)}
 \sup_{r_{\sigma_{v}}=k}|\alpha_{v}(\sigma_{v};m)|
 \ll_{\pi,\eta} q_{v}^{m(1+\eta)}.
\end{align*}
For $k\geq2$, the Schur sum in the proof of Lemma~\ref{lem:depth-kernel}
and \eqref{eq:newform-Schur-bound} bound the summand divided by
$q_{v}^{m(1+\eta)}$ by
$O((1+m-k)^{5} q_{v}^{-(1+\eta-\vartheta_{\pi_{v}})(m-k)})$.
For $k\leq1$ and $m\geq2$, Lemma~\ref{lem:depth-low} gives
$O((1+m)^{5}q_{v}^{1+\eta+m\vartheta_{\pi_{v}}})$ before division,
using bounded local $L$-factors and \eqref{eq:finite-newvector-value}.
This is $O(q_{v}^{m(1+\eta)})$ since $\vartheta_{\pi_{v}}<1/4$.
For $m=1$, use Lemma~\ref{lem:prime-old-row-expansion}.
Sum \eqref{eq:mixed-positive-box-bound} over dyadic ranges of $C(\sigma_{v})^{1/2}$
and exact conductors $\mathfrak{a}\mid\mathfrak{b}$, applying the local
sum at $m=e_{v}(\mathfrak{b})$ and $\eta=\varepsilon/2$.

\begin{lemma}\label{lem:mixed-fixed-residues}
For the mixed moment $\mathcal{M}_{\pi}(\mathfrak{a},\mathfrak{n})$,
the original and dual standard residues satisfy
\begin{equation}\label{eq:mixed-residual-bound}
\begin{gathered}
\big|R_{\RNum{2}}^{2,+}(\mathbf{0},\varphi,\mathbf{1},\mathbf{1})\big|
 \ll_{\pi,F,\varepsilon}
 N_{F}(\mathfrak{a})^{1/2+\vartheta+\varepsilon}N_{F}(\mathfrak{n})^{1/2+\varepsilon},\\
\big|R_{\RNum{2}}^{2,+}(\mathbf{0},\pi(w_{2})\varphi,\mathbf{1},\mathbf{1})\big|
 \ll_{\pi,F,\varepsilon}
 N_{F}(\mathfrak{a})^{1+\varepsilon}N_{F}(\mathfrak{n})^{-1/2+\varepsilon}
 g_{\pi}^{+}(\mathfrak{n}).
\end{gathered}
\end{equation}
\end{lemma}

\begin{proof}
See \textsection\ref{app:mixed-residual-estimates}.
\end{proof}

The estimate $\mathcal{V}_{v}^{\circ}(\sigma_{v})\ll_{A} C(\sigma_{v})^{-A}$
in Lemma~\ref{lem:mixed-spherical-arch} gives the following
specialization of \eqref{eq:mixed-full-spectrum-finite-weights}
\begin{multline}\label{eq:mixed-absolute-row}
 \int_{\mathcal{F}(\mathfrak{a};\mathbf{1})}
 \frac{L(1/2,\pi\times\sigma)L(1/2,\sigma)}{L(1,\sigma,\Ad)}
 \prod_{v\mid\infty}\mathcal{V}_{v}^{\circ}(\sigma_{v})\\
 \prod_{v\mid\mathfrak{a}}|\alpha_{v}(\sigma_{v};e_{v}(\mathfrak{a}))|
 \,d\mu_{\sigma}
 \ll_{\pi,F,\varepsilon}N_{F}(\mathfrak{a})^{1+\varepsilon}.
\end{multline}

\begin{lemma}\label{lem:mixed-conductor-row}
Suppose $\mathfrak{a}+\mathfrak{n}=\mathcal{O}_{F}$. Then
\begin{multline}\label{eq:mixed-conductor-row}
|\mathcal{M}_{\pi}(\mathfrak{a},\mathfrak{n})|
 \ll_{\pi,F,\varepsilon}
 N_{F}(\mathfrak{a})^{1+\varepsilon}N_{F}(\mathfrak{n})^{-1/2+\varepsilon}
 g_{\pi}^{+}(\mathfrak{n})
 \\
+N_{F}(\mathfrak{a})^{1/2+\vartheta+\varepsilon}
  N_{F}(\mathfrak{n})^{1/2+\varepsilon}.
\end{multline}
\end{lemma}

\begin{proof}
The finite calculations in \eqref{14.37} and Lemma~\ref{lem14.6}
apply to the $\psi_{v}^{0}$-spherical vector with value $1$ at $I_{3}$, without the
Jacquet scalar $\mathbf{L}(\pi_{v})^{-1}$. Dividing by
$\mathcal{L}_{v}(\mathbf{0},\pi_{v})$ gives $q_{v}^{e_{v}(\mathfrak{a})}$
at $v\mid\mathfrak{a}$ and $q_{v}^{-n_{v}/2}g_{\pi}(\mathfrak{p}_{v}^{n_{v}})$
at $v\mid\mathfrak{n}$. In \eqref{e14.126}, use \eqref{14.45} and \eqref{14.36}
at $v\mid\mathfrak{a}$, and \eqref{14.46} at $v\mid\mathfrak{n}$.
With the norm factors of \eqref{e14.126}, this gives
$O_{\varepsilon}(N_{F}(\mathfrak{n})^{1/2}\prod_{v\mid\mathfrak{a}}
q_{v}^{e_{v}(\mathfrak{a})(\vartheta_{\pi_{v}}+\varepsilon)})$.
The two standard residues in \eqref{8.1} satisfy \eqref{eq:mixed-residual-bound}.
The Hecke bound used in \eqref{15.12}, at $\mathfrak{q}=\mathfrak{a}$,
gives $|\widetilde{\lambda}_{\sigma,\mathbf{0}}(\mathfrak{a})|
\ll N_{F}(\mathfrak{a})^{\vartheta+\varepsilon}$.
At $\mathbf{s}=\mathbf{0}$, factor out $q_{v}^{-n_{v}/2}$ from
\eqref{eq:dual-level-superposition}. For $\ell=e_{v}(\mathfrak{n}')$,
the weight is $\alpha_{v}(\sigma_{v};\ell)$ and vanishes unless
$r_{\sigma_{v}}\leq\ell$.
Summing \eqref{eq:mixed-absolute-row} at
$\mathfrak{a}=\mathfrak{n}'$ over $\mathfrak{n}'\mid\mathfrak{n}$, with the
factor $N_{F}(\mathfrak{a})^{1/2}N_{F}(\mathfrak{n})^{-1/2}$, gives
\begin{equation}\label{eq:mixed-dual-bound}
 |\mathcal{M}_{\pi}^{\vee}(\mathfrak{a},\mathfrak{n})|
 \ll N_{F}(\mathfrak{a})^{1/2+\vartheta+\varepsilon}
 N_{F}(\mathfrak{n})^{1/2+\varepsilon}.
\end{equation}
Combining \eqref{eq:mixed-dual-bound} and \eqref{eq:mixed-residual-bound}
with the two generic estimates in \eqref{8.1} proves
\eqref{eq:mixed-conductor-row}.
\end{proof}

\begin{prop}\label{prop:mixed-newform-sieve}
Suppose $\mathfrak{a}$ is squarefree, or $\mathfrak{a}=\mathfrak{p}^{r}$
with $\pi_{\mathfrak{p}}$ tempered, and $\mathfrak{a}+\mathfrak{n}=\mathcal{O}_{F}$.
Then
\begin{multline}\label{eq:general-mixed-newform-twisted}
\Big|
 \int_{\substack{\sigma\in\mathcal{F}^{\circ}(\mathfrak{a};\mathbf{1})\\
\mathfrak{c}(\sigma)=\mathfrak{a}}}
 \frac{\lambda_{\sigma}(\mathfrak{n})
 L(1/2,\pi\times\sigma)L(1/2,\sigma)}
 {L^{(\mathfrak{a})}(1,\sigma,\Ad)}
 \mathcal{W}_{\infty}(\sigma)d\mu_{\sigma}\Big|\\
 \ll_{\pi,F,\varepsilon}
 \alpha_{\pi}(\mathfrak{a})\left(
 N_{F}(\mathfrak{a})^{1+\varepsilon}
 N_{F}(\mathfrak{n})^{-1/2+\varepsilon}g_{\pi}^{+}(\mathfrak{n})
 +N_{F}(\mathfrak{a})^{1/2+\vartheta+\varepsilon}
 N_{F}(\mathfrak{n})^{1/2+\varepsilon}\right).
\end{multline}
For squarefree $\mathfrak{a}$, the bound also holds with
$\prod_{v\in S}q_{v}^{-1/2}\lambda_{\sigma_{v}}(\mathfrak{p}_{v})$ inserted in the integrand, for any subset
$S$ of the primes dividing $\mathfrak{a}$.
\end{prop}

\begin{proof}
For squarefree $\mathfrak{a}$, the expansion
\eqref{eq:squarefree-newform-expansion} holds for the mixed moment
with $c_{v,j}(\pi_{v})$ from \eqref{eq:mixed-prime-old-row-series}.
The translate $-q_{v}^{-1}\pi_{v}(\diag(w_{\mathfrak{p}_{v}},1))W_{v}^{[1]}$
at $v\in S$ contributes $q_{v}^{-1/2}\lambda_{\sigma_{v}}(\mathfrak{p}_{v})$
at $r_{\sigma_{v}}=1$ by \eqref{eq:signed-newvector-factor}.
At $r_{\sigma_{v}}=0$, \eqref{eq:signed-old-row} multiplies the series
coefficients by $-q_{v}^{-1}$, so subtraction gives $q_{v}^{-1}c_{v,j}(\pi_{v})$ at the level
with $\mathfrak{p}_{v}$ removed and the twist multiplied by $\mathfrak{p}_{v}^{j}$.

In the translated dual moment, group by $\ell=e_{v}(\mathfrak{n}')$
in \eqref{eq:dual-level-superposition} at $v\mid\mathfrak{n}$ and
$r_{\sigma_{v}}=\mathbf{1}_{v\mid\mathfrak{m}}$ at $v\in S$
by \eqref{eq:signed-section-dual-identity}. The conductor divides
$\mathfrak{n}'\mathfrak{m}$ and \eqref{eq:signed-dual-level-factor} gives
$N_{F}(\mathfrak{m})^{-1}$. Sum \eqref{eq:mixed-absolute-row} at level
$\mathfrak{n}'\mathfrak{m}$ over $\mathfrak{n}'\mid\mathfrak{n}$ and
$\mathfrak{m}\mid\prod_{v\in S}\mathfrak{p}_{v}$, with the prefactor
$N_{F}(\mathfrak{n})^{-1/2}$ and the factors in \eqref{eq:signed-old-row}
and \eqref{eq12.14}. The $N_{F}(\mathfrak{m})$ in the level bound
cancels the factor in \eqref{eq:signed-dual-level-factor}, giving the second term of
\eqref{eq:general-mixed-newform-twisted}. At $v\in S$, the generic
bounds follow from \eqref{eq:signed-generic-series} and
\eqref{eq:signed-dual-generic-transport} at $\mathbf{s}=\mathbf{0}$, and
the residue bounds from the moving limit of
\eqref{eq:signed-old-row-parameters} in \textsection\ref{app:mixed-residual-estimates}.
Together they give \eqref{eq:mixed-conductor-row} for the translated
moments with the factor $\alpha_{\pi}(\mathfrak{a})$.

By \eqref{eq:squarefree-sieve-coefficients},
$\alpha_{\pi}(\mathfrak{a}/\mathfrak{b})|c_{\mathfrak{b},S}(\mathfrak{l})|
\leq\alpha_{\pi}(\mathfrak{a})N_{F}(\mathfrak{l})^{-1/2}$. For small $\eta>0$
and $L\geq1$, this gives
\begin{equation}\label{eq:mixed-sieve-coefficient-bounds}
\begin{split}
 &\sum_{\mathfrak{b}\mid\mathfrak{a}}
 \sum_{\substack{\mathfrak{p}\mid\mathfrak{l}\Rightarrow\mathfrak{p}\mid\mathfrak{b}\\
 N_{F}(\mathfrak{l})\leq L}}
 \alpha_{\pi}(\mathfrak{a}/\mathfrak{b})|c_{\mathfrak{b},S}(\mathfrak{l})|
 N_{F}(\mathfrak{l})^{1/2+\eta}
 \ll_{\pi,F,\eta}\alpha_{\pi}(\mathfrak{a})N_{F}(\mathfrak{a})^{\eta}L^{2\eta},\\
 &\sum_{\mathfrak{b}\mid\mathfrak{a}}
 \sum_{\mathfrak{p}\mid\mathfrak{l}\Rightarrow\mathfrak{p}\mid\mathfrak{b}}
 \alpha_{\pi}(\mathfrak{a}/\mathfrak{b})|c_{\mathfrak{b},S}(\mathfrak{l})|
 g_{\pi}^{+}(\mathfrak{l})N_{F}(\mathfrak{l})^{-1/2+\eta}
 \ll_{\pi,F,\eta}\alpha_{\pi}(\mathfrak{a})N_{F}(\mathfrak{a})^{\eta}.
\end{split}
\end{equation}
The first estimate uses \eqref{eq:sieve-supported-Euler-sum} at
$(m,u)=(0,\eta)$ and $N_{F}(\mathfrak{l})^{2\eta}\leq L^{2\eta}$.
For the second, the local sum is
$1+\sum_{j\geq0}g_{\pi}^{+}(\mathfrak{p}_{v}^{j})q_{v}^{-j(1-\eta)}
=2+O_{\pi,\eta}(q_{v}^{-1+\eta+\vartheta_{\pi_{v}}})$.
Apply \eqref{eq:mixed-conductor-row} at level $\mathfrak{a}/\mathfrak{b}$
and twist $\mathfrak{n}\mathfrak{l}$. In \eqref{eq:mixed-sieve-coefficient-bounds},
use the second sum for the term containing $g_{\pi}^{+}(\mathfrak{n}\mathfrak{l})$
and the first for the other term.

At $\mathbf{s}=\mathbf{0}$, \eqref{eq:signed-newvector-factor} and
\eqref{eq:signed-old-row} give the factor $q_{v}^{-1}$ in the absolute
local weight at $v\in S$. For the tail, \eqref{eq:mixed-absolute-row}
and the Hecke bound give the $\mathfrak{l}$-dependent factor
$|c_{\mathfrak{b},S}(\mathfrak{l})|\tau(\mathfrak{l})N_{F}(\mathfrak{l})^{\vartheta}
\ll_{\eta}\alpha_{\pi}(\mathfrak{b})N_{F}(\mathfrak{l})^{-1/2+\vartheta+\eta}$.
Since $N_{F}(\mathfrak{l})>L$, bound one factor
$N_{F}(\mathfrak{l})^{-(1/2-\vartheta)/2}$ by $L^{-(1/2-\vartheta)/2}$
and apply \eqref{eq:sieve-supported-Euler-sum} with
$(m,u)=(0,(1/2-\vartheta)/2-\eta)$ to get
\begin{equation}\label{eq:mixed-sieve-tail}
 \sum_{\mathfrak{b}\mid\mathfrak{a}}
 \sum_{\substack{\mathfrak{p}\mid\mathfrak{l}\Rightarrow\mathfrak{p}\mid\mathfrak{b}\\
 N_{F}(\mathfrak{l})>L}}
 |c_{\mathfrak{b},S}(\mathfrak{l})|
 |\mathcal{M}_{\pi}(\mathfrak{a}/\mathfrak{b},\mathfrak{n}\mathfrak{l})|
 \ll_{\pi,F,\eta}\frac{\alpha_{\pi}(\mathfrak{a})
 N_{F}(\mathfrak{a})^{1+2\eta}N_{F}(\mathfrak{n})^{\vartheta+\eta}}{
 L^{(1/2-\vartheta)/2}}.
\end{equation}

With $L=N_{F}(\mathfrak{a}\mathfrak{n})^{\ell}$,
$\ell(1/2-\vartheta)/2>\vartheta+1$ and $2\eta(1+\ell)<\varepsilon$,
the bound \eqref{eq:mixed-sieve-tail}, also valid with the translations, is $O(\alpha_{\pi}(\mathfrak{a})N_{F}(\mathfrak{a})
N_{F}(\mathfrak{n})^{-1/2})$, and the truncated sum satisfies
\eqref{eq:general-mixed-newform-twisted}.

For $\mathfrak{a}=\mathfrak{p}^{r}$, $r\geq2$, the terms with
$2\leq r_{\sigma_{\mathfrak{p}}}<r$ in
$\sum_{m=2}^{r} b_{r-m}\mathcal{M}_{\pi}(\mathfrak{p}^{m},\mathfrak{n})$
vanish by \eqref{eq:depth-projector}. At $r_{\sigma_{\mathfrak{p}}}=r$
the local factor is $L_{\mathfrak{p}}(1,\sigma_{\mathfrak{p}},\Ad)$. Summing
\eqref{eq:mixed-conductor-row} with $|b_{r-m}|\leq\alpha_{\pi}(\mathfrak{p})$
gives \eqref{eq:general-mixed-newform-twisted} for this combination.
For $r_{\sigma_{\mathfrak{p}}}\leq1$, Lemma~\ref{lem:depth-low} at
$\vartheta_{\pi_{\mathfrak{p}}}=0$, the Hecke bound, and
Lemma~\ref{lem:mixed-full-spectrum-mass} at level $\mathfrak{p}$ give
\begin{multline*}
 \left|\int_{\mathcal{F}^{\circ}(\mathfrak{p};\mathbf{1})}
 \frac{\lambda_{\sigma}(\mathfrak{n})L(1/2,\pi\times\sigma)L(1/2,\sigma)}
 {L(1,\sigma,\Ad)}
 \mathcal{W}_{\infty}(\sigma)
 \sum_{m=2}^{r} b_{r-m}\alpha_{\mathfrak{p}}(\sigma_{\mathfrak{p}};m)\,
 d\mu_{\sigma}\right|\\
 \ll_{\pi,F,\varepsilon}
 N_{F}(\mathfrak{p}^{r})^{\varepsilon} N_{F}(\mathfrak{n})^{\vartheta+\varepsilon}
 N_{F}(\mathfrak{p})^{1+\varepsilon}
 \ll N_{F}(\mathfrak{a})^{1/2+2\varepsilon}N_{F}(\mathfrak{n})^{1/2+\varepsilon}.
\end{multline*}
Subtract this integral to obtain \eqref{eq:general-mixed-newform-twisted}.
\end{proof}

With $\beta_{v}$ from Lemma~\ref{lem:real-rs-evaluation}, put
\begin{align*}
G_{\pi_{v}}(t)=\beta_{v}
 \frac{L_{v}(1/2,\pi_{v}\times\sigma_{t})L_{v}(1/2,\sigma_{t})}
 {L_{v}(1,\sigma_{t},\Ad)},\qquad
 \sigma_{t}=|\cdot|^{it}\boxplus|\cdot|^{-it}.
\end{align*}
Writing $\nu_{1},\nu_{2},\nu_{3}$ for the Langlands parameters of $\pi_{v}$,
the spherical local factors give
\begin{align*}
 \frac{L_{v}(1/2,\pi_{v}\times\sigma_{t})}{L_{v}(1/2,\sigma_{t})^{3}}
 =\prod_{j=1}^{3}\prod_{\delta=\pm1}
 \frac{\Gamma_{\mathbb{R}}(1/2+\nu_{j}+\delta it)}
 {\Gamma_{\mathbb{R}}(1/2+\delta it)}.
\end{align*}
Since $\sum_{j}\nu_{j}=0$, Stirling's formula shows that this ratio
and its reciprocal have $j$-th derivatives
$O_{\pi_{v},j}((1+|t|)^{-j})$ on $\mathbb{R}$ and tend to $1$ as
$|t|\to\infty$. For complementary series,
$L_{v}(1/2,\pi_{v}\times\sigma_{v})/L_{v}(1/2,\sigma_{v})^{3}\asymp_{\pi_{v}}1$
by the uniform bound on the Gamma arguments in the proof of
Lemma~\ref{lem:mixed-spherical-arch}.

The test functions on $G'(\mathbb{R})$ are the lifts of the functions
on $\mathrm{PGL}_{2}(\mathbb{R})$ multiplied by $\alpha(|\det g|^{1/2})$,
with the cutoff $\alpha$ from \textsection\ref{def:central-test-cutoff}.

\begin{lemma}\label{lem:mixed-afe-realization}
Fix $l>0$ and $R_{0}>0$, and put
\begin{align*}
 S_{0}(t)=e^{-t^{2}/R_{0}^{2}}
 \prod_{\substack{j\geq0\\1/2+2j\leq l+1}}
 \bigl(t^{2}+(1/2+2j)^{2}\bigr)^{2}.
\end{align*}
Then $S_{0}(t)>0$ for real $t$ and for $t=i\eta$, $|\eta|<1/2$.
For fixed $\Re(z)=c>0$, the product
$S_{z}(t):=S_{0}(t)m_{z}(t)$ is even and holomorphic on
$|\Im(t)|<l+1$. For every $m,r\geq0$,
\begin{equation}\label{eq:mixed-afe-test-seminorms}
 \sup_{|\Im(t)|\leq l}(1+|\Re(t)|)^{m}
 |\partial_{t}^{r}S_{z}(t)|
 \ll_{l,R_{0},c,m,r}(1+|\Im(z)|)^{O(1)} e^{\pi|\Im(z)|/2}.
\end{equation}
Both $S_{0}$ and $S_{z}$ are limits of spherical transforms
of smooth compactly supported bi-$K'$-invariant functions, converging
in every bound in \eqref{eq:mixed-afe-test-seminorms}.
\end{lemma}

\begin{proof}
The double zeros of $S_{0}$ cancel the poles of $m_{z}$. The Stirling
and Cauchy bounds in Lemma~\ref{lem:hybrid-afe-multiplier} at $R=R_{0}$
give \eqref{eq:mixed-afe-test-seminorms}. Smooth truncation of the inverse
Fourier transforms followed by Abel inversion \cite[Chapter~IV, \textsection7]{Hel83}
gives the compact approximations. The strip $|\Im(t)|<l+1$ gives
convergence in \eqref{eq:mixed-afe-test-seminorms}.
\end{proof}

Let $h_{v,z}$ be the inverse spherical transform of $S_{z}$.
Let $W_{h_{v,z}}$ be the limit of the Whittaker vectors obtained from
the compactly supported test functions of Lemma~\ref{lem:mixed-afe-realization},
with restrictions given by \eqref{fc11.3} at $\mathbf{s}=\mathbf{0}$.

\phantomsection\label{sec:proof-mixed-moment}
\begin{proof}[Proof of Corollary~\ref{thm1.7}]
In \eqref{eq:mixed-full-spectrum-finite-weights} at
$\mathfrak{b}=\mathcal{O}_{F}$ and fixed $\ell>\ell_{0}$,
$\prod_{v\mid\infty}C(\sigma_{v})^{-\ell/2}\geq T^{-\ell/2}$ on
$C_{\infty}(\sigma)\leq T$, and
$L(1/2,\sigma)L(1,\sigma,\Ad)\ll_{F,\varepsilon}T^{1/4+\varepsilon}$.
Thus the level-one sum in \eqref{eq1.3} is $O_{\pi,F}(T^{O(1)})$.

At each real place, use the ranges specified before
\eqref{eq:original-weighted-fourth}. Fix $W_{v}=W_{h_{v,0}}$ for the
bounded range and $W_{v}=W_{h_{v,R_{v},0}}$ from
Corollary~\ref{cor:mixed-polynomial-input-decay} for a dyadic range.
For the entire spectral integral, set
\begin{align*}
 \mathcal{W}_{\infty,\mathbf{R}}(\sigma)
 :=\epsilon_{\sigma} D_{F}^{-3/2}
 \prod_{v\mid\infty}\pi\mathcal{H}_{W_{v}}(\sigma_{v},\mathbf{0},\mathbf{1}).
\end{align*}
On the spherical spectrum, \eqref{eq:real-weight-evaluated} gives
$G_{\pi_{v}}(t_{v})S_{0}(t_{v})$ for $W_{h_{v,0}}$ and $H_{R_{v}}(t_{v})$ for
$W_{h_{v,R_{v},0}}$, at every $t_{v}$. The factor vanishes on the
nonspherical spectrum by \eqref{fc10.1}. The vectors $W_{h_{v,z}}$
and $W_{h_{v,R_{v},z}}$ multiply these factors by $m_{z}(t_{v})$.

Take $\mathfrak{a}\mid\mathfrak{q}$, $\mathfrak{a}\ne\mathcal{O}_{F}$,
and $\mathfrak{a}+\mathfrak{n}=\mathcal{O}_{F}$. Put
$X=\prod_{v\mid\infty}R_{v}^{2}$ and $0<\Re(z)=c\leq2$.
For squarefree $\mathfrak{a}$, take any $S\subseteq\{v:v\mid\mathfrak{a}\}$.
For $\mathfrak{a}=\mathfrak{p}^{r}$, $r\geq2$, take $S=\varnothing$. Then
\begin{multline}\label{eq:mixed-polynomial-AFE-sieve}
 \Big|\int_{\substack{\sigma\in\mathcal{F}^{\circ}(\mathfrak{a};\mathbf{1})\\
\mathfrak{c}(\sigma)=\mathfrak{a}}}
 \frac{\lambda_{\sigma}(\mathfrak{n})
 L(1/2,\pi\times\sigma)L(1/2,\sigma)\mathcal{W}_{\infty,\mathbf{R}}(\sigma)}
 {L^{(\mathfrak{a})}(1,\sigma,\Ad)} 
 \prod_{v\mid\infty}m_{z}(t_{v})
 \prod_{v\in S}\frac{\lambda_{\sigma_{v}}(\mathfrak{p}_{v})}{q_{v}^{1/2}}d\mu_{\sigma}\Big|\\
 \ll_{\pi,F,\varepsilon}\alpha_{\pi}(\mathfrak{a})
 X^{c/2+O(1)}(1+|\Im(z)|)^{O(1)}
 e^{[F:\mathbb{Q}]\pi|\Im(z)|/2}\\
 \big(N_{F}(\mathfrak{a})^{1+\varepsilon}N_{F}(\mathfrak{n})^{-1/2+\varepsilon}g_{\pi}^{+}(\mathfrak{n})
       +N_{F}(\mathfrak{a})^{1/2+\vartheta+\varepsilon}N_{F}(\mathfrak{n})^{1/2+\varepsilon}\big).
\end{multline}
Corollary~\ref{cor:mixed-polynomial-input-decay} bounds the real terms
for $h_{v,R_{v},z}$. For $h_{v,z}$, the bound for $\mathscr{L}_{b}(h_{v,z})$
in Lemma~\ref{lem:mixed-fixed-real-dual-decay} controls the two generic
integrals \eqref{7.5} at $\mathbf{s}=\mathbf{0}$ by
\eqref{eq:mixed-fixed-Whittaker-bound}, the dual weight by
\eqref{eq:mixed-afe-local-decay}, and the dual residue quotient by
\eqref{eq:mixed-real-standard-pole-bound}. The original residues
vanish by \textsection\ref{app:mixed-residual-estimates}.
On the original spectrum, the strip estimate in the proof of
Lemma~\ref{lem:hybrid-afe-multiplier} for $H_{R_{v}}m_{z}$ and
\eqref{eq:mixed-afe-test-seminorms} for $S_{z}$ give the $X,z$ factor
in \eqref{eq:mixed-polynomial-AFE-sieve}, since $G_{\pi_{v}}$ is bounded
on the spherical spectrum here. These bounds include the decay
$\prod_{v\mid\infty}C(\sigma_{v})^{-\ell/2}$ for fixed $\ell>\ell_{0}$.

Apply \eqref{eq:mixed-conductor-row} and \eqref{eq:mixed-absolute-row}
at level $\mathfrak{a}/\mathfrak{b}$ and twist $\mathfrak{n}\mathfrak{l}$,
with the $X,z$ factor in \eqref{eq:mixed-polynomial-AFE-sieve}.
The sums \eqref{eq:mixed-sieve-coefficient-bounds} and
\eqref{eq:mixed-sieve-tail} use $L=N_{F}(\mathfrak{a}\mathfrak{n})^{\ell}$
with $\ell(1/2-\vartheta)/2>\vartheta+1$. At $\mathfrak{a}=\mathfrak{p}^{r}$,
\eqref{eq:depth-projector} keeps conductor $\mathfrak{p}^{r}$ among
$r_{\sigma_{\mathfrak{p}}}\geq2$. Lemma~\ref{lem:depth-low} at
$\vartheta_{\pi_{\mathfrak{p}}}=0$ bounds the terms with $r_{\sigma_{\mathfrak{p}}}\leq1$.
The polynomial exponents in \eqref{eq:mixed-polynomial-AFE-sieve}
are uniform for $0<c\leq2$ and independent of $\varepsilon$.

In $L(1/2,\pi\times\sigma)L(1/2,\sigma)^{2}$, expand one $L(1/2,\sigma)$
by \eqref{eq:fifth-moment-root-number-removal} on $\Re(z)=2$. Separate
the factors at $v\mid\mathfrak{a}$ by \eqref{eq:bad-prime-AFE-sign-expansion}.
For $l\leq N_{F}(\mathfrak{n})<2l$ and $(\mathfrak{n},\mathfrak{a})=1$,
we move the $z$-integral in \eqref{eq:fifth-moment-root-number-removal}
to $\Re(z)=c$ using holomorphy for $\Re(z)>0$ and Gaussian decay with
\eqref{eq:mixed-polynomial-AFE-sieve}. On this line
$(D_{F}^{2}N_{F}(\mathfrak{a})X)^{c/2}N_{F}(\mathfrak{n})^{-1/2-c}$ is
$\ll_{F} l^{-1/2}((N_{F}(\mathfrak{a})X)^{1/2}/l)^{c}$.
Integrating against $|e^{z^{2}}dz/z|$ and summing the coefficients by
\eqref{eq:bad-prime-AFE-subset-bound} for squarefree $\mathfrak{a}$
bounds the contribution of $l\leq N_{F}(\mathfrak{n})<2l$ by
\begin{align*}
 \ll\alpha_{\pi}(\mathfrak{a})X^{O_{\pi,F}(1)}
 \left(\frac{(N_{F}(\mathfrak{a})X)^{1/2}}{l}\right)^{c}
 \left(N_{F}(\mathfrak{a})^{1+\varepsilon}l^{\varepsilon}
 +N_{F}(\mathfrak{a})^{1/2+\vartheta+\varepsilon}l^{1+\varepsilon}\right).
\end{align*}
Here $\sum_{N_{F}(\mathfrak{n})\asymp l}g_{\pi}^{+}(\mathfrak{n})
\ll_{\varepsilon} l^{1+\varepsilon}$ by \eqref{eq:mixed-generic-average},
and there are $O_{F}(l)$ ideals in the range.
Summing over dyadic $l$ with $c=\varepsilon$ for
$l\leq(N_{F}(\mathfrak{a})X)^{1/2}$ and $c=2$ otherwise gives
$O(\alpha_{\pi}(\mathfrak{a})N_{F}(\mathfrak{a})^{1+\vartheta+O(\varepsilon)}
X^{O_{\pi,F}(1)})$.

Since $L(1/2,\pi\times\sigma)L(1/2,\sigma)^{2}
\mathcal{W}_{\infty,\mathbf{R}}(\sigma)\geq0$, we keep the cuspidal terms.
In the bounded range, $G_{\pi_{v}}(t_{v})S_{0}(t_{v})\gg_{\pi_{v}}1$ by
Lemmas~\ref{lem:mixed-spherical-arch} and~\ref{lem:mixed-afe-realization}.
In the dyadic range, $H_{R_{v}}(t_{v})\geq e^{-8}$ by the example after
\eqref{eq:R-admissible-test}. Thus $\mathcal{W}_{\infty,\mathbf{R}}(\sigma)
\gg_{\pi,F}1$ on the selected ranges, where $C_{\infty}(\sigma)\asymp X$ and
$L^{(\mathfrak{a})}(1,\sigma,\Ad)\ll_{F,\varepsilon}
(N_{F}(\mathfrak{a})X)^{\varepsilon}$.
Thus $C_{\infty}(\sigma)\leq T$ requires only the
$O_{F}((1+\log T)^{[F:\mathbb{Q}]})$ choices with $X\ll_{F} T$.
Summing over these choices and $\mathfrak{a}\mid\mathfrak{q}$ proves
\eqref{eq1.3}.
\end{proof}

\section{Amplification and Nonvanishing}\label{sec16}
\subsection{\texorpdfstring{Subconvexity of Triple Product $L$-functions}{Subconvexity of Triple Product L-functions}}\label{sec:proof-triple-product}
Let $\pi_{1}$ be as in Theorem~\ref{thm1.6}, let
$\mathfrak{q}$ be a nontrivial integral ideal, and let
$\pi_{2}\in\mathcal{F}_{\mathrm{cusp}}^{\circ}(\mathfrak{q};\mathbf{1})$
have conductor exactly $\mathfrak{q}$.

\subsubsection{Amplification}
For a prime ideal $\mathfrak{p}\nmid\mathfrak{q}$, let
$\ell_{\mathfrak{p}}$ be the smallest positive integer for which
$|\lambda_{\pi_{2}}(\mathfrak{p}^{\ell_{\mathfrak{p}}})|\geq10^{-1}$.
Since $\lambda_{\pi_{2}}(\mathfrak{p})^{2}=\lambda_{\pi_{2}}(\mathfrak{p}^{2})+1$,
we have $\ell_{\mathfrak{p}}\in\{1,2\}$.

Let $\mathfrak{L}$ be the set of prime ideals $\mathfrak{p}\nmid\mathfrak{q}$
with $L<N_{F}(\mathfrak{p})\leq2L$. For $L$ a fixed positive power of
$N_{F}(\mathfrak{q})$ the prime ideal theorem gives
$\#\mathfrak{L}\asymp_{F} L/\log L$ for large $L$. Put $\alpha_{\mathfrak{p}}
:=\overline{\lambda_{\pi_{2}}(\mathfrak{p}^{\ell_{\mathfrak{p}}})}/
|\lambda_{\pi_{2}}(\mathfrak{p}^{\ell_{\mathfrak{p}}})|$.
For unitary $\sigma$ with trivial central character, unramified
at the primes in $\mathfrak{L}$, the Hecke relations give
\begin{equation}\label{16.1}
\Big|\sum_{\mathfrak{p}\in \mathfrak{L}}\alpha_{\mathfrak{p}}\lambda_{\sigma}(\mathfrak{p}^{\ell_{\mathfrak{p}}})\Big|^{2}=\sum_{\substack{\mathfrak{p}\in \mathfrak{L}\\ 0\leq j\leq \ell_{\mathfrak{p}}}}\lambda_{\sigma}(\mathfrak{p}^{2j})
+\sum_{\substack{\mathfrak{p}_{1}\in \mathfrak{L}\\ \mathfrak{p}_{2}\in \mathfrak{L}\\ 
\mathfrak{p}_{1}\neq\mathfrak{p}_{2}}}\alpha_{\mathfrak{p}_{1}}\overline{\alpha_{\mathfrak{p}_{2}}}\lambda_{\sigma}(\mathfrak{p}_{1}^{\ell_{\mathfrak{p}_{1}}}\mathfrak{p}_{2}^{\ell_{\mathfrak{p}_{2}}}).
\end{equation}

Let $\mathcal{J}_{\mathfrak{L}}(f)$ denote
the right-hand side of \eqref{16.1} with each
$\lambda_{\sigma}(\mathfrak{n})$ replaced by $f(\mathfrak{n})$.

Take $\pi=\Ad\pi_{1}$ and the weight $\mathcal{W}_{\infty}$ from
\textsection\ref{def:mixed-real-weight}.
For $\mathfrak{n}+\mathfrak{q}=\mathcal{O}_{F}$, put
\begin{align*}
\mathcal{M}_{\mathfrak{q}}^{\rm new}(\mathfrak{n})
 :=\prod_{\mathfrak{p}_{v}\Vert\mathfrak{q}}(1-q_{v}^{-2})
 \int_{\substack{\sigma\in\mathcal{F}^{\circ}(\mathfrak{q};\mathbf{1})\\
\mathfrak{c}(\sigma)=\mathfrak{q}}}
 \frac{\lambda_{\sigma}(\mathfrak{n})
 L(1/2,\Ad\pi_{1}\times\sigma)L(1/2,\sigma)}
 {L^{(\mathfrak{q})}(1,\sigma,\Ad)}\mathcal{W}_{\infty}(\sigma)\,d\mu_{\sigma}.
\end{align*}

Both central values are nonnegative for cuspidal $\sigma$ by
\textsection\ref{sec13.2} and \cite[Theorem~1]{LR03}, and are absolute
squares on the continuous spectrum. With Lemma~\ref{lem:mixed-spherical-arch},
\eqref{16.1} thus gives a nonnegative integrand.
At $\sigma=\pi_{2}$, the amplifier has modulus
$\sum_{\mathfrak{p}\in\mathfrak{L}}|\lambda_{\pi_{2}}(\mathfrak{p}^{\ell_{\mathfrak{p}}})|
\gg L/\log L$, and $\mathcal{W}_{\infty}(\pi_{2})>0$ is fixed by its real components.
Since $\prod_{\mathfrak{p}_{v}\Vert\mathfrak{q}}(1-q_{v}^{-2})\geq\zeta_{F}(2)^{-1}$,
the adjoint upper bound and
$L(s,\pi_{1}\times\widetilde{\pi}_{1}\times\pi_{2})
=L(s,\pi_{2})L(s,\Ad\pi_{1}\times\pi_{2})$ give
\begin{equation}\label{16.9}
 \mathcal{J}_{\mathfrak{L}}
 (\mathcal{M}_{\mathfrak{q}}^{\rm new})
 \gg L^{2-\varepsilon}N_{F}(\mathfrak{q})^{-\varepsilon}
 L(1/2,\pi_{1}\times\widetilde{\pi}_{1}\times\pi_{2}).
\end{equation}

Here $\alpha_{\Ad\pi_{1}}(\mathfrak{q})\ll_{F,\varepsilon}N_{F}(\mathfrak{q})^{\varepsilon}$
by temperedness at $v\mid\mathfrak{q}$.
By \eqref{eq:mixed-generic-coefficient} and self-duality,
$g_{\pi}^{+}(\mathfrak{p}^{j})\ll(1+|\lambda_{\pi}(\mathfrak{p})|)^{j}$
for $j=1,2$. Rankin--Selberg bounds
$\sum_{\mathfrak{p}\in\mathfrak{L}}(1+|\lambda_{\pi}(\mathfrak{p})|)^{2}$
by $O(L^{1+\varepsilon})$, so the prime sum of
$g_{\pi}^{+}(\mathfrak{p}^{\ell_{\mathfrak{p}}})/N_{F}(\mathfrak{p})^{\ell_{\mathfrak{p}}/2}$
is $O(L^{1/2+\varepsilon})$.
For the first term of \eqref{eq:general-mixed-newform-twisted}, square
the prime bound for the off-diagonal twists and use
$g_{\pi}^{+}(\mathfrak{p}^{2j})\ll N_{F}(\mathfrak{p})^{j}$ for the
diagonal twists by the Satake bound. For the second term,
\eqref{16.1} has $O(L^{2})$ twists of norm $O(L^{4})$, so the sum of
$N_{F}(\mathfrak{n})^{1/2}$ is $O(L^{4})$. Taking $\mathfrak{a}=\mathfrak{q}$ gives
\begin{equation}\label{16.13}
 \left|\mathcal{J}_{\mathfrak{L}}
 (\mathcal{M}_{\mathfrak{q}}^{\rm new})\right|
 \ll L^{1+\varepsilon}N_{F}(\mathfrak{q})^{1+\varepsilon}
 +L^{4+\varepsilon}N_{F}(\mathfrak{q})^{1/2+\vartheta+\varepsilon}.
\end{equation}
After dividing by the factor in \eqref{16.9}, balance
$N_{F}(\mathfrak{q})/L$ and $N_{F}(\mathfrak{q})^{1/2+\vartheta}L^{2}$
at $L=N_{F}(\mathfrak{q})^{(1-2\vartheta)/6}$ to prove Theorem~\ref{thm1.6}.

\pagebreak
\subsection{Simultaneous Nonvanishing}\label{sec:simultaneous-nonvanishing}

\begin{lemma}
Let $\pi$ be as in Corollary~\ref{thm1.7}. For squarefree $\mathfrak{q}$
and the weight $\mathcal{W}_{\infty}$ in
\textsection\ref{def:mixed-real-weight},
\begin{multline}\label{e16.18}
 \int_{\sigma\in\mathcal{F}_{\mathrm{cusp}}^{\circ}(\mathfrak{q};\mathbf{1})}
 \frac{L(1/2,\pi\times\sigma)L(1/2,\sigma)}{L(1,\sigma,\Ad)}
 \mathcal{W}_{\infty}(\sigma)
 \prod_{v\mid\mathfrak{q}}\alpha_{v}(\sigma_{v};1)\,d\mu_{\sigma}\\
 =c_{\pi} N_{F}(\mathfrak{q})+O_{F,\pi,\varepsilon}
 (N_{F}(\mathfrak{q})^{1/2+\vartheta+\varepsilon}),
\end{multline}
with a constant $c_{\pi}>0$ independent of $\mathfrak{q}$.
\end{lemma}

\begin{proof}
In \eqref{8.1} take $\mathbf{s}=\mathbf{0}$ and $\mathfrak{n}=\mathcal{O}_{F}$.
Both degenerate terms vanish since $\pi$ is cuspidal.
At finite places use $\psi_{v}^{0}$ and $W_{v}^{\circ}(I_{3})=1$.
The change back to $\psi_{v}$ contributes factors independent of
$\mathfrak{q}$. The finite calculation in \eqref{14.37} at $r_{v}=1$
and \eqref{eq:real-jacquet-double-mellin} at infinity give the original
generic term as a positive multiple of
$N_{F}(\mathfrak{q})L(1,\pi)^{2}/\zeta_{F}(2)>0$.

By evenness, the signed half-residues at $\lambda=\pm(1/2-s_{2}')$
sum to the residue at $\lambda=1/2-s_{2}'$. As $\pi$ is unramified at
all places with trivial central character, its standard root number
is $+1$. The functional equation gives
$\Lambda(0,\pi)=D_{F}^{3/2}\Lambda(1,\pi)$, so
\eqref{eq:mixed-standard-completed-residue} is positive at
$\mathbf{s}^{\vee}=\mathbf{0}$. At $\lambda=1/2$ and $r_{v}=1$,
the finite factor in \eqref{eq12.14} after division by the two local
$L$-factors is $q_{v}$. The real quotient
\eqref{eq:mixed-real-standard-pole-value} is $1$ since $w_{2}$ fixes
the spherical vector. Thus the dual residue gives a second positive
multiple of $N_{F}(\mathfrak{q})$.

The original residue and dual generic term are
$O(N_{F}(\mathfrak{q})^{5/14+\varepsilon})$ by \eqref{14.36} at $r_{v}=1$
in \eqref{eq:mixed-finite-endpoint} and \eqref{14.45}, respectively. In the original Eisenstein integral, $2r_{\xi_{v}}=r_{\sigma_{v}}\leq1$
at $v\mid\mathfrak{q}$, so \eqref{eq:old-row-simplified} at $r_{\sigma_{v}}=0$
bounds the finite product by $O(N_{F}(\mathfrak{q})^{5/14+\varepsilon})$.
The real weight $\prod_{v\mid\infty}\mathcal{V}_{v}^{\circ}(\sigma_{v})$ makes the
remaining integral convergent. The dual spectrum contributes
$O(N_{F}(\mathfrak{q})^{1/2+\vartheta+\varepsilon})$ by
\eqref{eq:mixed-dual-bound} at $(\mathfrak{a},\mathfrak{n})=(\mathfrak{q},\mathcal{O}_{F})$.
These bounds give the error in \eqref{e16.18}.
\end{proof}

\phantomsection\label{sec:proof-simultaneous-nonvanishing}
\begin{proof}[Proof of Theorem~\ref{thm1.9}]
The cuspidal integral in \eqref{e16.18} keeps all conductors dividing
$\mathfrak{q}$.
For $v\mid\mathfrak{q}$,
$|\alpha_{v}(\sigma_{v};1)|\ll q_{v}^{(1-r_{\sigma_{v}})5/14}$ by
Lemma~\ref{lem:prime-old-row-expansion} and \eqref{eq:old-row-simplified},
using $L_{v}(1,\sigma_{v},\Ad)=\zeta_{v}(2)$ at $r_{\sigma_{v}}=1$.
Since $\prod_{v\mid\mathfrak{q}}q_{v}^{1-r_{\sigma_{v}}}
=N_{F}(\mathfrak{q})/C_{\fin}(\sigma)$, the adjoint lower bound gives
\begin{equation}\label{eq:nonvanishing-normalization-bound}
\frac{\prod_{v\mid\mathfrak{q}}|\alpha_{v}(\sigma_{v};1)|}
 {L(1,\sigma,\Ad)}\ll_{F,\pi,\varepsilon}
 (N_{F}(\mathfrak{q})/C_{\fin}(\sigma))^{5/14}
 \bigl(N_{F}(\mathfrak{q})C_{\infty}(\sigma)\bigr)^{\varepsilon}.
\end{equation}

Equation~\eqref{eq:individual-hybrid-fifth} at
$(\mathfrak{f},\mathfrak{q})=(\mathcal{O}_{F},\mathfrak{c}(\sigma))$
and convexity for $\pi\times\sigma$, of finite conductor
$C_{\fin}(\sigma)^{3}$, give
$|L(1/2,\pi\times\sigma)L(1/2,\sigma)|\ll
C_{\fin}(\sigma)^{19/20+\vartheta/5+\varepsilon}
C_{\infty}(\sigma)^{l_{0}}$ for some fixed $l_{0}$.
Let $0\leq j_{v}\leq1$ be smooth, equal to $1$ at
$x=\nu_{\sigma_{v}}^{2}$ for automorphic spherical $\sigma_{v}$ and supported in
$x<(1/2-\vartheta_{\pi_{v}})^{2}$. The uniform gap in the proof of
Lemma~\ref{lem:mixed-spherical-arch} allows this choice and keeps the
Gamma factors positive on the support. Set
\begin{align*}
 \alpha_{\infty}(\sigma_{\infty})
 =c\prod_{v\mid\infty}\mathcal{V}_{v}^{\circ}(\sigma_{v})
 j_{v}(\nu_{\sigma_{v}}^{2})(1+\nu_{\sigma_{v}}^{4})^{m}.
\end{align*}
The decay in Lemma~\ref{lem:mixed-spherical-arch} gives rapid decay
after the polynomial factor. Choose $c>0$ so that $\alpha_{\infty}\leq1$.
On the automorphic spherical spectrum, $j_{v}=1$ and
$C(\sigma_{v})\asymp(1+|\nu_{\sigma_{v}}|)^{2}$. Thus $2m\geq l_{0}+1$ gives
$\alpha_{\infty}(\sigma_{\infty})\gg C_{\infty}(\sigma)^{l_{0}+1}
\allowbreak\prod_{v\mid\infty}\mathcal{V}_{v}^{\circ}(\sigma_{v})$.

After \eqref{eq:nonvanishing-normalization-bound}, the exponent of
$C_{\fin}(\sigma)$ is $19/20+\vartheta/5-5/14>0$.
Replacing $C_{\fin}(\sigma)$ by $N_{F}(\mathfrak{q})$ bounds the
absolute contribution of each $\sigma$ to \eqref{e16.18} by
$O(N_{F}(\mathfrak{q})^{19/20+\vartheta/5+\varepsilon}
\alpha_{\infty}(\sigma_{\infty}))$.
Summing over $\sigma$ with $L(1/2,\pi\times\sigma)L(1/2,\sigma)\neq0$
and comparing with the main term in \eqref{e16.18} gives \eqref{1.5},
since $1-(19/20+\vartheta/5)=(1-4\vartheta)/20$.
\end{proof}

\clearpage
\appendix
\section{Local Evaluations for the Twisted Fourth Moment}\label{app:fourth-local}
Define $I_{v,w}(\mathbf{s};f_{v})$ by
\eqref{eq:local-w2-degenerate}, \eqref{eq:local-w1w2-degenerate}, and
\eqref{eq:local-long-weyl-degenerate} for $w=w_{2},w_{1}w_{2},w_{1}w_{2}w_{1}$,
respectively, with section $f_{v}$, trivial Hecke characters and
the factor $\mathcal{L}_{v}^{-1}$ omitted.

\subsection{Local Sections and Normalizations}
\subsubsection{Archimedean places}
We define the section by
\begin{equation}\label{fc14.3}
f_{\boldsymbol{\nu},v}\left(g_{v}\right):=\int_{G'(F_{v})}f_{\boldsymbol{\nu},v}^{\circ}\left(g_{v}\diag(y_{v},1)\right)h_{v}(y_{v}^{-1})|\det y_{v}|_{v}^{s_{1}}dy_{v}.
\end{equation}

\subsubsection{Non-Archimedean places: integral representation}
Write $n_{v}=e_{v}(\mathfrak{n})$ and $r_{v}=e_{v}(\mathfrak{q})$.
In \eqref{e5.3} with trivial characters and the indicators of
$M_{2\times3}(\mathcal{O}_{v})$ and $M_{1\times2}(\mathcal{O}_{v})$,
the multiplier for $f_{\boldsymbol{\nu},v}^{\circ}(I_{3})=1$ is
$q_{v}^{2d_{v}}/\mathbf{L}(\pi_{v})$ by \eqref{eq:finite-character-measures}.
For $v\nmid\mathfrak{q}$, set
\begin{equation}\label{c14.5}
f_{\boldsymbol{\nu},v}:=\sum_{\ell=0}^{n_{v}}q_{v}^{-\frac{n_{v}}{2}+(2\ell-n_{v})s_{2}}\sum_{\alpha\in \mathcal{O}_{v}/\mathfrak{p}_{v}^{\ell}}\pi_{v}(w_{2}n_{13}(\alpha\varpi_{v}^{-\ell}))f_{\boldsymbol{\nu},v}^{\circ}.
\end{equation}
For $v\mid\mathfrak{q}$, set
\begin{equation}\label{c14.6}
f_{\boldsymbol{\nu},v}:=\sum_{\alpha\in \mathcal{O}_{v}/\mathfrak{p}_{v}^{r_{v}}}\pi_{v}(n_{13}(\alpha\varpi_{v}^{-r_{v}}))f_{\boldsymbol{\nu},v}^{\circ}.
\end{equation}

\subsubsection{Jacquet normalization}
For the Jacquet integral \eqref{f7.11}, initially in
$\Re(\nu_{1})>\Re(\nu_{2})>\Re(\nu_{3})$, write
$W_{v}^{\circ}=c_{v}\mathcal{J}_{\psi_{v}}f_{\boldsymbol{\nu},v}^{\circ}$.
Commuting the Jacquet integral with the right translations in \eqref{fc14.3},
\eqref{c14.5}, and \eqref{c14.6} gives
\begin{equation}\label{eq:jacquet-normalization-transport}
W_{v}=c_{v}\mathcal{J}_{\psi_{v}} f_{\boldsymbol{\nu},v}.
\end{equation}
At a finite place, \textsection\ref{sec12} and
\eqref{eq:different-jacquet-scalar} give
$c_{v}=\mathbf{L}(\pi_{v})q_{v}^{-d_{v}(1/2+\nu_{3}-\nu_{1})}$.

\subsection{\texorpdfstring{The $w_{2}$ Degenerate Term}{The w2 Degenerate Term}}
Start in the convergence region \eqref{eq5.13}.

\subsubsection{Archimedean integrals}
\begin{lemma}\label{lem14.1}
Let $v$ be real. Then
\begin{multline}\label{14.4}
I_{v,w_{2}}(\mathbf{s};f_{\boldsymbol{\nu},v})
=\frac{\zeta_{v}(1+\nu_{1}+\nu_{2}+2s_{1})\zeta_{v}(1+\nu_{1}+\nu_{3}+2s_{1})}
{\zeta_{v}(1+\nu_{2}-\nu_{3})\zeta_{v}(2+2\nu_{1}+2s_{1})}\\
\zeta_{v}(s_{2}-s_{1}-\nu_{1})\zeta_{v}(1+\nu_{1}+s_{1}+s_{2})
\mathcal{S}h_{v}(1/2+\nu_{1}+s_{1}).
\end{multline}
\end{lemma}

\begin{proof}
Write the unipotent matrix in \eqref{eq:local-w2-degenerate} as
$n_{31}(a_{v})n_{32}(c_{v})$. Insert \eqref{fc14.3}, replace $g$ by
$w'g$, and write
$g=z'\left(\begin{smallmatrix}1&\\b'&1\end{smallmatrix}\right)\diag(a',1)k'$.
For fixed $g$ with $a',z'>0$, we substitute
\begin{align*}
 (a,c,y,z)\longmapsto(a/a'-b',\ z'c,\ a'y,\ z/z').
\end{align*}
Moving $\diag(a'z',1,z')$ to the left gives phase
$\overline{\psi}_{v}(ay+cz-a'b'y)$ and weight $|a'|^{1+s_{2}}$ after
the Jacobian $|z'/a'|$ and the section factor \eqref{c14.1}.
In the $g$-integral, at fixed $a'$ we substitute
$(b',z')\mapsto(b'/a',1/(z'\sqrt{a'}))$. The factor $|a'|^{-1}$
gives the weight $|a'|^{s_{2}}$ in \eqref{eq10.17}, hence $\widehat{f}_{v}(y)$.

Take
\begin{align*}
\begin{cases}
 \Re(\nu_{1}-\nu_{2})>-1,\quad \Re(\nu_{2}-\nu_{3})>0,\quad \Re(1+\nu_{1}+\nu_{3}+2s_{1})>0,\\
 \Re(\nu_{1}+s_{1})>-1/2,\quad
 \Re(s_{2}+1/2)>|\Re(1/2+\nu_{1}+s_{1})|.
\end{cases}
\end{align*}
Define
\begin{equation}\label{eq14.15}
F_{1}(y,z):=\int_{\mathbb{R}^{2}}
(1+b^{2})^{-(1+\nu_{1}-\nu_{2})/2}
(1+b^{2}+c^{2})^{-(1+\nu_{2}-\nu_{3})/2}e^{-2\pi i(by+cz)}\,db\,dc.
\end{equation}
In \eqref{eq14.15} take the limit with cutoff $e^{-\varepsilon(b^{2}+c^{2})}$
as $\varepsilon\to 0^{+}$. Evaluating \eqref{c14.1} at $n_{31}(b)n_{32}(c)$
gives this kernel, and summing the signs of $a',z',y,z$ gives
\begin{equation}\label{14.8}
I_{v,w_{2}}(\mathbf{s};f_{\boldsymbol{\nu},v})
=16\int_{0}^{\infty}\int_{0}^{\infty}
F_{1}(y,z)\widehat{f}_{v}(y)z^{1+\nu_{1}+\nu_{3}+2s_{1}}y^{s_{2}-s_{1}-\nu_{1}}
\,d^{\times} z\,d^{\times} y.
\end{equation}

For $\Re(P)>0$ and $\Re(2B+P)>1$, Gaussian integration gives
\begin{equation}\label{14.9}
\int_{0}^{\infty}\int_{\mathbb{R}}
(1+b^{2}+c^{2})^{-B}e^{-2\pi icz}\,dc\ z^{P}\,d^{\times} z
=\frac{\Gamma(P/2)\Gamma(B+(P-1)/2)}{2\pi^{P-1/2}\Gamma(B)(1+b^{2})^{B+(P-1)/2}}.
\end{equation}
Kontorovich--Lebedev inversion in the normalization \eqref{eq11.26}
gives, by holomorphy for $|\Re(d)|<\Re(s_{2}+1/2)$,
\begin{equation}\label{14.11}
\int_{0}^{\infty}\widehat{f}_{v}(y)y^{s_{2}+1/2}K_{d}(2\pi y)\,d^{\times} y 
=\frac{\mathcal{S}h_{v}(d)\Gamma \big(\frac{1/2+s_{2}+d}{2}\big)
\Gamma\big(\frac{1/2+s_{2}-d}{2}\big)}{16\pi^{1/2+s_{2}}}.
\end{equation}
In \eqref{14.8}, \eqref{14.9} evaluates the $z$-integral at
$(B,P)=((1+\nu_{2}-\nu_{3})/2,1+\nu_{1}+\nu_{3}+2s_{1})$. Fourier transforming
$(1+b^{2})^{-1-\nu_{1}-s_{1}}$ gives the $y$-integral \eqref{14.11} at
$d=1/2+\nu_{1}+s_{1}$.

\end{proof}

\subsubsection{Non-Archimedean integrals}\label{sec14.2.2}

\begin{lemma}\label{lem14.2}
Suppose $v\mid\mathfrak{q}$. Then 
\begin{multline}\label{e14.21}
I_{v,w_{2}}(\mathbf{s};f_{\boldsymbol{\nu},v})=\frac{q_{v}^{r_{v}+(s_{2}-s_{1}-\nu_{1}-2)d_{v}}\zeta_{v}(s_{2}-s_{1}-\nu_{1})}{\mathbf{L}(\pi_{v})}
\sum_{i'=0}^{\infty}\sum_{i=0}^{\infty}\sum_{j=r_{v}-d_{v}}^{\infty}
\\
\sum_{j'=0}^{j+d_{v}}\mathbf{1}_{i'+j\geq i+j'-d_{v}}q_{v}^{-i'(1+\nu_{1}-\nu_{2})}q_{v}^{j'(\nu_{3}+s_{1}-s_{2})}
q_{v}^{-i(\nu_{2}-\nu_{3})}q_{v}^{-j(1+\nu_{1}+\nu_{3}+2s_{1})}.
\end{multline}
\end{lemma}

\begin{proof}
In \eqref{eq:local-w2-degenerate} for the section \eqref{c14.6},
conjugation by $w_{2}w_{1}$ sends $n_{13}$ to $n_{32}$. The
sum of translates multiplies the $c$-Fourier coefficient by
$q_{v}^{r_{v}}\mathbf{1}_{e_{v}(z)\geq r_{v}-d_{v}}$
by orthogonality on $e_{v}(z)\geq-d_{v}$.

At fixed $z'$, substitute in the local integral \eqref{e5.3}
\begin{align*}
 (a',t,b',a,c,y,z)\mapsto
 \bigl((z')^{-1}a',z't,-a'b'/(z't), 
 (z')^{-1}a,(z')^{-1}(c+a'b'/t),z'y,z'z\bigr).
\end{align*}
With \eqref{e1.5} the substitution has measure factor $|(z')^{-2}/t|_{v}$.
The two Schwartz functions are the indicators of
$(a',t,z',a,c,b')\in\mathcal{O}_{v}^{6}$. Write $i',i,j',r,j$ for the
valuations of $a',t,z',y,z$, respectively. The phase
$\overline{\psi}_{v}(ay+cz+b'a'z/t)$ imposes $r,j\geq-d_{v}$ and
$i'+j\geq i-d_{v}$. The measure factor $q_{v}^{-4d_{v}}$ from \eqref{eq:finite-character-measures}
and the section multiplier give $q_{v}^{-2d_{v}}/\mathbf{L}(\pi_{v})$.
The five Mellin weights in these variables are
\begin{align*}
 q_{v}^{-i'(1+\nu_{1}-\nu_{2})-i(\nu_{2}-\nu_{3})-j'(1+\nu_{1}+s_{1}+s_{2})
       -r(s_{2}-s_{1}-\nu_{1})-j(1+\nu_{1}+\nu_{3}+2s_{1})}.
\end{align*}
Sum over $r\geq-d_{v}$, then replace $j$ by $j-j'$ in the remaining
sum. The conditions become $j\geq r_{v}-d_{v}$ and $j'\leq j+d_{v}$,
and the $j'$-weight is $q_{v}^{j'(\nu_{3}+s_{1}-s_{2})}$.
\end{proof}

\begin{lemma}\label{lem14.4}
Suppose $v<\infty$ and $v\nmid\mathfrak{q}$. Then 
\begin{multline*}
I_{v,w_{2}}(\mathbf{s};f_{\boldsymbol{\nu},v})=\mathbf{L}(\pi_{v})^{-1}\sum_{\ell=0}^{n_{v}}q_{v}^{-\frac{n_{v}}{2}+(2\ell-n_{v})s_{2}+\ell-2d_{v}}\sum_{r=-d_{v}}^{\infty}\sum_{i'=0}^{\infty}\sum_{i=0}^{\infty}\sum_{j'=0}^{\infty}\sum_{j=-d_{v}}^{\infty}\\
\mathbf{1}_{\substack{i'+j\geq i-d_{v}\\ r+j'\geq \ell-d_{v}}}q_{v}^{-i'(1+\nu_{1}-\nu_{2})}q_{v}^{-j'(1+\nu_{1}+s_{1}+s_{2})}q_{v}^{-i(\nu_{2}-\nu_{3})}
q_{v}^{-j(1+\nu_{1}+\nu_{3}+2s_{1})}q_{v}^{-r(s_{2}-s_{1}-\nu_{1})}.
\end{multline*}	
\end{lemma}

\begin{proof}
In the $\ell$th term of \eqref{c14.5}, conjugation by $w_{2}w_{1}w_{2}$
sends $n_{13}$ to $n_{31}$. Thus the level condition in the proof of
\eqref{e14.21} is $r+j'\geq\ell-d_{v}$. Keep both $r$ and $j$ in the
sum and multiply by $q_{v}^{-n_{v}/2+(2\ell-n_{v})s_{2}+\ell}$.
\end{proof}

\subsection{\texorpdfstring{The $w_{1}w_{2}$ Degenerate Term}{The w1w2 Degenerate Term}}
Start in the convergence region \eqref{f5.19}.

\subsubsection{Non-Archimedean integrals}\label{sec14.3.2}

\begin{lemma}\label{lemm14.5}
Suppose $v\mid\mathfrak{q}$. Then 
\begin{multline}\label{fc14.39}
I_{v,w_{1}w_{2}}(\mathbf{s};f_{\boldsymbol{\nu},v})=q_{v}^{-5d_{v}/2}\mathbf{L}(\pi_{v})^{-1}\zeta_{v}(1+\nu_{1}+s_{1}+s_{2})\\
\sum_{l=-d_{v}}^{\infty}\sum_{k=0}^{l+d_{v}}q_{v}^{-(1+\nu_{2}+s_{1}+s_{2})l-k(\nu_{1}+s_{1}-s_{2})}\sum_{i=-l-d_{v}}^{\infty}\sum_{j=-l}^{\infty}\mathbf{1}_{\min\{i+d_{v},j\}\leq d_{v}}\\
q_{v}^{r_{v}-\max\{0,r_{v}-d_{v}-i-l\}(1+\nu_{1}+s_{1}+s_{2})-(1+\nu_{2}+\nu_{3}+2s_{1})i-(1+\nu_{2}-\nu_{3})j+\min\{i+d_{v},j\}}.
\end{multline}
\end{lemma}

\begin{proof}
In the local integral of \eqref{5.19}, take trivial characters,
$k=I_{2}$, and replace the right factor $w_{1}w_{2}w_{1}$ by
$w_{1}w_{2}w_{1} n_{13}(\alpha\varpi_{v}^{-r_{v}})$, summed over
$\alpha\bmod\mathfrak{p}_{v}^{r_{v}}$. Integrate in $a,b,c$ and sum over
$e_{v}(t_{2})\geq\max(0,r_{v}-d_{v}-e_{v}(z))$:
\begin{multline}\label{eq14.41}
h_{v}((t_{1},y),(z,z');\nu_{1},s_{2})=q_{v}^{r_{v}-2d_{v}-\max\{0,r_{v}-d_{v}-e_{v}(z)\}(1+\nu_{1}+s_{1}+s_{2})}\\\mathbf{1}_{e_{v}(t_{1})\geq 0}\mathbf{1}_{e_{v}(y)\geq -d_{v}}\mathbf{1}_{e_{v}(z)\geq -d_{v}}\mathbf{1}_{e_{v}(z')\geq -d_{v}}\zeta_{v}(1+\nu_{1}+s_{1}+s_{2}).
\end{multline}

In \eqref{c5.25}, \eqref{eq14.41} at $(y/t_{1},zy,b'z)$ gives the
$b'$-integral $\mathbf{1}_{m\leq d_{v}}q_{v}^{m-d_{v}/2}$ with
$m=\min(e_{v}(z)+d_{v},e_{v}(t))$. Multiplying by the measure factor in
\eqref{eq:finite-character-measures} and the section normalization yields
\eqref{fc14.39}.
\end{proof}

\begin{remark}
We define $I_{v,w}^{0}(\mathbf{s};f_{v})$ by evaluating $I_{v,w}(\mathbf{s};f_{v})$ with
$\psi_{v}^{0}$ and the corresponding measures in \eqref{eq:finite-character-measures}. Scaling both outer Fourier variables $y,z$ by
$\varpi_{v}^{-d_{v}}$ gives the Mellin multiplier $q_{v}^{(1+\nu_{3}+s_{1}+s_{2})d_{v}}$.
The measure factors from \eqref{eq:finite-character-measures} are
$q_{v}^{-2d_{v}}$ for $w_{2}$ and $q_{v}^{-5d_{v}/2}$ for $w_{1}w_{2}$, so
\begin{equation}\label{eq:degen-different-check}
\begin{aligned}
 &I_{v,w_{2}}(\mathbf{s};f_{\boldsymbol{\nu},v})
 =q_{v}^{(\nu_{3}+s_{1}+s_{2}-1)d_{v}}I_{v,w_{2}}^{0}(\mathbf{s};f_{\boldsymbol{\nu},v}),\\
 &I_{v,w_{1}w_{2}}(\mathbf{s};f_{\boldsymbol{\nu},v})
 =q_{v}^{(\nu_{3}+s_{1}+s_{2}-3/2)d_{v}}I_{v,w_{1}w_{2}}^{0}(\mathbf{s};f_{\boldsymbol{\nu},v}).
\end{aligned}
\end{equation}
\end{remark}

\begin{lemma}\label{lemma14.7}
Suppose $v\nmid\mathfrak{q}$. Then 
\begin{flalign}\label{fc14.48}
 && I_{v,w_{1}w_{2}}(\mathbf{s};f_{\boldsymbol{\nu},v})=\frac{\zeta_{v}(1+\nu_{1}+s_{1}+s_{2})\zeta_{v}(1+\nu_{2}+\nu_{3}+2s_{1})
\zeta_{v}(1+\nu_{2}-\nu_{3})}{q_{v}^{(5/2+\nu_{2}-\nu_{3})d_{v}}\mathbf{L}(\pi_{v})\zeta_{v}(2+2\nu_{2}+2s_{1})}&\\
 && \sum_{\ell=0}^{n_{v}}q_{v}^{-n_{v}/2+(2\ell-n_{v})s_{2}}
\sum_{j=-d_{v}}^{\infty}q_{v}^{-(1+\nu_{2}+s_{1}+s_{2})j}
\sum_{k=0}^{j+d_{v}}q_{v}^{k(1+2\nu_{2}+2s_{1})}&\notag\\
 && \sum_{i=0}^{j+d_{v}}q_{v}^{\ell-\max\{0,\ell+i-j-d_{v}\}(1+\nu_{1}+s_{1}+s_{2})}q_{v}^{-i(\nu_{1}+s_{1}-s_{2})}.&\notag
\end{flalign}
\end{lemma}

\begin{proof}
For $\alpha\bmod\mathfrak{p}_{v}^{\ell}$, the $\ell$th term of
\eqref{c14.5} replaces the right factor in \eqref{5.19} by
$w_{1}w_{2}w_{1}w_{2}\allowbreak n_{13}(\alpha\varpi_{v}^{-\ell})$.
In the level factor of \eqref{fc14.39}, this replaces $r_{v}$ by
$\ell$ and $e_{v}(zy)$ by $e_{v}(y/t_{1})$, by \eqref{eq14.41}
and \eqref{c5.25}. For $-e_{v}(y)\leq m=\min(d_{v}+e_{v}(z),e_{v}(t))\leq d_{v}$,
sum over $d_{v}+e_{v}(z)=m\leq e_{v}(t)$ and $e_{v}(t)=m<d_{v}+e_{v}(z)$.
The two geometric sums give the zeta quotient in \eqref{fc14.48} for
$\Re(1+\nu_{2}+\nu_{3}+2s_{1}),\Re(1+\nu_{2}-\nu_{3})>0$.
In these sums we index by $e_{v}(z)+d_{v}$ and $k=d_{v}-m$,
giving the factor $q_{v}^{(\nu_{3}-\nu_{2})d_{v}}$. The remaining indices
are $j=e_{v}(y)$ and $i=e_{v}(t_{1})$, with the coefficient in \eqref{c14.5}.
\end{proof}

\subsection{The \texorpdfstring{$w_{1}w_{2}w_{1}$}{w1w2w1} Degenerate Term}\label{app:long-weyl}
Start in the convergence region \eqref{e5.33}.

 \subsubsection{Non-Archimedean integrals}\label{sec14.4.2}

\begin{lemma}\label{lem14.11}
Suppose $v\mid\mathfrak{q}$. Then
\begin{multline}\label{eq14.67}
 I_{v,w_{1}w_{2}w_{1}}(\mathbf{s};f_{\boldsymbol{\nu},v})
 =\frac{q_{v}^{r_{v}+(\nu_{2}+s_{1}+s_{2}-2)d_{v}}}{\mathbf{L}(\pi_{v})}\\
 \frac{\zeta_{v}(\nu_{2}-\nu_{3})\zeta_{v}(\nu_{1}-\nu_{3})
 \zeta_{v}(1+\nu_{2}+\nu_{3}+2s_{1})\zeta_{v}(s_{2}-s_{1}-\nu_{3})}
 {\zeta_{v}(1-\nu_{1}+\nu_{3})\zeta_{v}(1+s_{1}-s_{2}+\nu_{3})}\\
 \sum_{j,k\geq0}\sum_{i=0}^{j+k}
 q_{v}^{-(1-\nu_{1}+\nu_{3})j-(1+s_{1}-s_{2}+\nu_{3})k-(\nu_{1}+s_{1}+s_{2})i
 -(1+\nu_{2}+\nu_{3}+2s_{1})\max\{0,r_{v}-i\}}.
\end{multline}
Initially $\Re(\nu_{2}-\nu_{3})>0$, $0<\Re(\nu_{1}-\nu_{3})<1$,
$0<\Re(s_{2}-s_{1}-\nu_{3})<1$, $\Re(1+\nu_{2}+\nu_{3}+2s_{1})>0$, and
$\Re(\nu_{1}+s_{1}+s_{2})>0$. The identity extends meromorphically.
\end{lemma}

\begin{proof}
For $w_{1}w_{2}w_{1}$, scaling $y,z$ by $\varpi_{v}^{-d_{v}}$ gives
$q_{v}^{(1+\nu_{2}+s_{1}+s_{2})d_{v}}$ from the Mellin powers and $q_{v}^{-3d_{v}}$
from \eqref{eq:finite-character-measures}. Thus
\begin{equation}\label{eq:longest-different-scaling}
 I_{v,w_{1}w_{2}w_{1}}(\mathbf{s};f_{\boldsymbol{\nu},v})
 =q_{v}^{(\nu_{2}+s_{1}+s_{2}-2)d_{v}}I_{v,w_{1}w_{2}w_{1}}^{0}(\mathbf{s};f_{\boldsymbol{\nu},v}).
\end{equation}
For $I_{v,w_{1}w_{2}w_{1}}^{0}$ take the local factor of \eqref{e5.75} with
$\psi_{v}^{0}$, trivial Hecke characters, and the section \eqref{c14.6}.
The Tate integral in \eqref{e5.75} contributes $\zeta_{v}(\nu_{2}-\nu_{3})$.
In \eqref{e5.25} we subtract
$\alpha\varpi_{v}^{-r_{v}}(t_{1},at_{1})^{\mathsf{T}}$ from the second column
at fixed $a,t_{1}$. The $c$-Fourier
integral and character sum impose $e_{v}(z)\geq0$ and
$e_{v}(zt_{1})\geq r_{v}$. Thus the $z$-factor is
$q_{v}^{r_{v}-(1+\nu_{2}+\nu_{3}+2s_{1})\max(0,r_{v}-e_{v}(t_{1}))}\zeta_{v}(1+\nu_{2}+\nu_{3}+2s_{1})$.
For $0<\Re(\beta)<1$, Tate's identity is
\begin{equation}\label{14.73}
 \int_{F_{v}^{\times}}\mathbf{1}_{\mathcal{O}_{v}}(x+t)|t|_{v}^{\beta}\,d^{\times} t
 =\frac{\zeta_{v}(\beta)}{\zeta_{v}(1-\beta)}
 \int_{\mathcal{O}_{v}\setminus\{0\}}|u|_{v}^{1-\beta}
 \psi_{v}^{0}(xu)\,d^{\times} u.
\end{equation}
Equation \eqref{14.73} evaluates the $t_{2}$- and $y$-integrals at
$(x,\beta)=(ac',\nu_{1}-\nu_{3})$ and $(-au,s_{2}-s_{1}-\nu_{3})$, with the integral over
$c'\in\mathcal{O}_{v}$ between them. The $a$-integral contributes
$|t_{1}|_{v}^{-1}\mathbf{1}_{e_{v}(uv')\geq e_{v}(t_{1})}$, giving
$0\leq i\leq j+k$ for $(j,k,i)=(e_{v}(u),e_{v}(v'),e_{v}(t_{1}))$.
Together with \eqref{eq:longest-different-scaling} this gives \eqref{eq14.67}.
\end{proof}

\begin{lemma}\label{lem14.10}
Suppose $v\nmid\mathfrak{q}$. Then
\begin{multline}\label{e14.72}
 I_{v,w_{1}w_{2}w_{1}}(\mathbf{s};f_{\boldsymbol{\nu},v})
 =\frac{q_{v}^{(\nu_{2}+s_{1}+s_{2}-2)d_{v}}}{\mathbf{L}(\pi_{v})}\\
 \frac{\zeta_{v}(\nu_{2}-\nu_{3})\zeta_{v}(\nu_{1}-\nu_{3})
 \zeta_{v}(1+\nu_{2}+\nu_{3}+2s_{1})\zeta_{v}(s_{2}-s_{1}-\nu_{3})}{\zeta_{v}(1-\nu_{1}+\nu_{3})}\\
 \sum_{\ell=0}^{n_{v}}q_{v}^{-n_{v}/2+(2\ell-n_{v})s_{2}+\ell}
 \sum_{i,j\geq0}q_{v}^{-(1+\nu_{1}+s_{1}+s_{2})i-(1-\nu_{1}+\nu_{3})j+\min(i,j)}\\
 q_{v}^{-(s_{2}-s_{1}-\nu_{3})\max\{\ell-i,\min(0,j-i)\}}.
\end{multline}
The series converges in the region of Lemma~\ref{lem14.11} and extends
meromorphically.
\end{lemma}

\begin{proof}
With $\psi_{v}^{0}$, the $\ell$th translate in \eqref{c14.5} shifts $c'$
in \eqref{e5.75} by $\alpha\varpi_{v}^{-\ell}t_{1}$ after permuting the columns.
The character sum gives $q_{v}^{\ell}\mathbf{1}_{e_{v}(y)\geq\ell-e_{v}(t_{1})}$.
For the $z$ and $\Phi_{1\times2,v}$ integrals, the proof of
\eqref{eq14.67} applies with $r_{v}=0$.
Applying \eqref{14.73} to $t_{2}$ at $(x,\beta)=(ac',\nu_{1}-\nu_{3})$
and integrating $c'\in\mathcal{O}_{v}$ give
$\mathbf{1}_{\mathcal{O}_{v}}(y-au)$. For $i=e_{v}(t_{1})$ and $j=e_{v}(u)$,
the change $a\mapsto au$ sends $\mathfrak{p}_{v}^{-i}$ to $\mathfrak{p}_{v}^{j-i}$.
The $a$-integral is $q_{v}^{\min(i,j)}\mathbf{1}_{e_{v}(y)\geq\min(0,j-i)}$.
After division by $\zeta_{v}(s_{2}-s_{1}-\nu_{3})$, the $y$-sum leaves its term
at $e_{v}(y)=\max(\ell-i,\min(0,j-i))$.
Equation \eqref{e14.72} follows with the coefficients in \eqref{c14.5}
and the factor in \eqref{eq:longest-different-scaling}.
For continuation, split $j\geq i$ into $i<\ell$ and $i\geq\ell$,
and $j<i$ into $j<\ell$ and $j\geq\ell$.
\end{proof}

\subsection{The Dual Degenerate Terms}
\subsubsection{Non-Archimedean integrals}
For $v\nmid\mathfrak{q}$, the two copies of $w_{2}$ cancel when
$\pi_{v}(w_{2})$ acts on each summand of \eqref{c14.5}.
The dual integral is
\begin{multline}\label{eq:finite-dual-degenerate-transport}
 I_{v,w}(\mathbf{s}^{\vee};\pi_{v}(w_{2})f_{\boldsymbol{\nu},v})
 =\sum_{\ell=0}^{n_{v}}q_{v}^{-n_{v}/2+(2\ell-n_{v})s_{2}}\\
 I_{v,w}\left(\mathbf{s}^{\vee};
 \sum_{\alpha\in\mathcal{O}_{v}/\mathfrak{p}_{v}^{\ell}}
 \pi_{v}(n_{13}(\alpha\varpi_{v}^{-\ell}))f_{\boldsymbol{\nu},v}^{\circ}\right).
\end{multline}
Equations \eqref{e14.21}, \eqref{fc14.39}, and \eqref{eq14.67} give the
inner integrals at $(\mathbf{s},r_{v})=(\mathbf{s}^{\vee},\ell)$, including
$\ell=0$ when the inner section is spherical.
For $v\mid\mathfrak{q}$, $\pi_{v}(w_{2})f_{\boldsymbol{\nu},v}$ is the inner sum
of \eqref{c14.5} at $\ell=r_{v}$. In Lemma~\ref{lem14.4}, \eqref{fc14.48},
and \eqref{e14.72}, evaluate the $\ell=n_{v}=r_{v}$ term at $\mathbf{s}^{\vee}$
and divide by the coefficient $q_{v}^{-r_{v}/2+r_{v}s_{2}'}$ from
\eqref{c14.5}. The sum over translates contributes $q_{v}^{r_{v}}$.

\subsubsection{Archimedean integrals}

\begin{lemma}\label{lem:real-dual-w2}
Let $v$ be real. Set
\begin{align*}
\alpha_{0}=\tfrac{1/2+s_{2}}{2},\quad
 \alpha_{1}=\tfrac{1/2+\nu_{1}+s_{1}}{2},\quad
 \alpha_{2}=\tfrac{1/2-\nu_{2}-s_{1}}{2},\quad
 \alpha_{3}=\tfrac{1/2-\nu_{3}-s_{1}}{2}.
\end{align*}
Then
\begin{multline}\label{c14.4}
I_{v,w_{2}}(\mathbf{s}^{\vee};\pi_{v}(w_{2})f_{\boldsymbol{\nu},v}) =\frac{\Gamma(s_{1}-\nu_{1}/2)}{8\pi^{2+s_{1}+s_{2}+\nu_{3}}}
\bigg[\Gamma\left(\frac{1+\nu_{1}-\nu_{2}}{2}\right)
\Gamma\left(\frac{1+\nu_{2}-\nu_{3}}{2}\right)\\
\Gamma\left(\frac{1+\nu_{1}-\nu_{3}}{2}\right)
\Gamma\left(\frac{1+\nu_{1}-2s_{1}}{2}\right)\bigg]^{-1}\\
\int_{\mathbb{R}}\mathcal{S}h_{v}(it)t\sinh(\pi t)
\prod_{j=0}^{3}\prod_{\epsilon\in\{\pm1\}}
\Gamma\left(\alpha_{j}+\frac{\epsilon it}{2}\right)\,dt.
\end{multline}
The integral over $\mathbb{R}$ converges for $\Re(\alpha_{j})>0$,
$0\leq j\leq3$, and the identity extends meromorphically.
\end{lemma}

\begin{proof}
Insert \eqref{fc14.3} in \eqref{eq:local-w2-degenerate} at
$(\mathbf{s}^{\vee},\pi_{v}(w_{2})f_{\boldsymbol{\nu},v})$ and replace $g$ by $w'gw'$.
At fixed $g$, in the coordinates for \eqref{14.4}, substitute
\begin{align*}
 (a,c,y,z)\longmapsto(z'a,\ a'^{-1}c-b',\ z'^{-1}y,\ a'z).
\end{align*}
The phase is $\overline{\psi}_{v}(ay+cz-a'b'z)$ and the $a'$-weight is
$|a'|^{1+s_{2}}$ as in \eqref{14.8}. In the $g$-integral, replace
$(b',z')$ by $(a'^{-1}b',z'^{-1}|a'|^{-1/2})$ at fixed $a'$.
The Jacobian $|a'|^{-1}$ gives the weight $|a'|^{s_{2}}$ in \eqref{eq10.17}.
Thus \eqref{14.8} holds with $\widehat{f}_{v}(z)$ in place of $\widehat{f}_{v}(y)$
and weight $z^{1+\nu_{1}+\nu_{3}+s_{2}-s_{1}}y^{2s_{1}-\nu_{1}}$.

Take $\Re(\nu_{1}-\nu_{2})>-1$, $\Re(\nu_{2}-\nu_{3})>0$,
$\Re(2s_{1}-\nu_{1})>0$ and $\Re(s_{2})>-1/2$.
In the dual integral, \eqref{fc11.29} gives
\begin{align*}
 z^{1+\nu_{1}+\nu_{3}+s_{2}-s_{1}}\widehat{f}_{v}(z)
 =\frac{1}{2\pi i}\int_{(c)} \mathcal{M}\widehat{f}_{v}(\lambda)
 z^{1+\nu_{1}+\nu_{3}+s_{2}-s_{1}-\lambda}\,d\lambda .
\end{align*}

We separate the Barnes poles by taking
\begin{align*}
 1/2+\Re(s_{2})<c<\min_{j=2,3}\bigl\{
 \Re(1-\nu_{j}+s_{2}-s_{1}),\,\Re(1+\nu_{1}+s_{1}+s_{2})\bigr\}.
\end{align*}
Equation \eqref{14.9} evaluates the $z$-integral at $B=(1+\nu_{2}-\nu_{3})/2$
and $P=1+\nu_{1}+\nu_{3}+s_{2}-s_{1}-\lambda$.
The Fourier transform in $b$ leaves the factor
$2\pi^{1+\nu_{1}+(s_{2}-s_{1}-\lambda)/2}/\allowbreak\Gamma(1+\nu_{1}+(s_{2}-s_{1}-\lambda)/2)$
times the $y$-integral in \eqref{10.42}, with integrand
$y^{(1+3s_{1}+s_{2}-\lambda)/2}\allowbreak K_{(1+2\nu_{1}+s_{2}-s_{1}-\lambda)/2}(2\pi y)$
against $d^{\times} y$.
The Laplace bound in \eqref{14.9} and Stirling with \eqref{fc11.29} give
$O(e^{-\pi|\tau|/2}(1+|t|+|\tau|)^{m}|\mathcal{S} h_{v}(it)|)$
on $\lambda=c+i\tau$, so the integrals can be interchanged.

With $x=-\lambda/2$, the $x$-dependent factors form the quotient
\begin{align*}
 \frac{
 \Gamma(x+\tfrac{1+\nu_{1}+s_{1}+s_{2}}{2})
 \prod_{j=2}^{3}\Gamma(x+\tfrac{1-\nu_{j}+s_{2}-s_{1}}{2})
 \prod_{\epsilon=\pm1}\Gamma(-x+\tfrac{-1/2-s_{2}+\epsilon it}{2})
 }{\Gamma(x+1+\nu_{1}+\tfrac{s_{2}-s_{1}}{2})}.
\end{align*}
Barnes' second lemma gives \eqref{c14.4}. The factors in
\eqref{14.8}, \eqref{14.9}, \eqref{10.42}, and \eqref{fc11.29}, with
the Jacobian factor $2$ on the upward $x$-contour, give
$1/(8\pi^{2+s_{1}+s_{2}+\nu_{3}})$.
For continuation, we separate the poles $it=\pm(2\alpha_{j}+2m)$,
$m\geq0$, and sum the crossed residues before any collision.
\end{proof}

\begin{lemma}
At a real place take $f=f_{\boldsymbol{\nu},v}$ from \eqref{fc14.3}.
Write $f_{(\nu_{2},\nu_{1},\nu_{3})}$ and $f_{(\nu_{3},\nu_{1},\nu_{2})}$ for the
sections with permuted inducing parameters. Then
\begin{equation}\label{eq:real-weyl-reductions}
\begin{gathered}
 I_{v,w_{1}w_{2}}(\mathbf{s};f)
 =\frac{\zeta_{v}(\nu_{1}-\nu_{2})}{\zeta_{v}(1+\nu_{1}-\nu_{2})}
 I_{v,w_{2}}(\mathbf{s};f_{(\nu_{2},\nu_{1},\nu_{3})}),\\
 I_{v,w_{1}w_{2}w_{1}}(\mathbf{s};f)
 =\frac{\zeta_{v}(\nu_{2}-\nu_{3})\zeta_{v}(\nu_{1}-\nu_{3})}
 {\zeta_{v}(1+\nu_{2}-\nu_{3})\zeta_{v}(1+\nu_{1}-\nu_{3})}
 I_{v,w_{2}}(\mathbf{s};f_{(\nu_{3},\nu_{1},\nu_{2})}).
\end{gathered}
\end{equation}
For the dual identities, evaluate each integral at $\mathbf{s}^{\vee}$
and right-translate its section by $w_{2}$. The exponent of
$|\det y_{v}|_{v}$ in \eqref{fc14.3} is still $s_{1}$, not $s_{1}'$.
\end{lemma}

\begin{proof}
In \eqref{eq:local-w1w2-degenerate}, the lower unipotent matrix is
$n_{21}(b)n_{31}(c)n_{32}(a)$. Moving the first $w_{1}$ to the left
rewrites the argument of $f$ as $w_{1}n_{12}(b)n_{31}(a)n_{32}(c)w_{2}w_{1}$,
so the $b$-integral is the $w_{1}$ intertwiner. This gives
\eqref{eq:local-w2-degenerate} with $(a_{v},c_{v})=(a,c)$ and inducing
parameters $(\nu_{2},\nu_{1},\nu_{3})$, since $w_{1}$ exchanges the first two.

In \eqref{eq:local-long-weyl-degenerate}, the lower unipotent matrix is
$n_{21}(c)n_{31}(b)n_{32}(a)$. The $w_{2},w_{1}$ intertwiners evaluated
at $n_{31}(c')n_{32}(c)w_{2}w_{1}$ give
\begin{align*}
w_{2}n_{23}(a)w_{1}n_{12}(b)n_{31}(c')n_{32}(c)w_{2}w_{1}
=n_{21}(c)n_{31}(b+ac)n_{32}(a)n_{23}(c')w_{1}w_{2}.
\end{align*}
The $ac$ is the $(3,1)$ entry of $n_{32}(a)n_{21}(c)$.
Since the weight in \eqref{eq:local-long-weyl-degenerate} is independent
of $b$, translating $b$ gives \eqref{eq:local-w2-degenerate} with
$(a_{v},c_{v})=(c',c)$ and inducing parameters $(\nu_{3},\nu_{1},\nu_{2})$,
since $w_{2}$ swaps positions $2,3$ and then $w_{1}$ swaps positions $1,2$.
The displayed scalars are the spherical factors in \eqref{c14.1}.
With the cutoff of \eqref{eq14.15} in $(a,c)$ in
\eqref{eq:local-w1w2-degenerate} and $(c',c)$ in
\eqref{eq:local-long-weyl-degenerate}, we take the intertwining integrals
first by Fubini at each cutoff and then pass to the Fourier limit.
Since the $w_{1}$ and $w_{2}$ intertwiners integrate left translates,
they act on $f_{\boldsymbol{\nu},v}^{\circ}$ inside the $y_{v}$-integral
in \eqref{fc14.3}, also with $g_{v}$ replaced by $g_{v}w_{2}$, so
\eqref{14.4} and \eqref{c14.4} apply with the permuted inducing parameters.
\end{proof}

\subsection{The Generic Term}\label{sec14.6}
The generic integrals are \eqref{7.5} with trivial Hecke characters.
For Theorem~\ref{thmC}, their real normalization on both sides is
$\beta_{v}=\mathbf{L}(\pi_{v})^{-1}$.

\subsubsection{Archimedean integrals}
\begin{lemma}
Let $v$ be real and let $W_{v}$ be given by \eqref{fc11.3}, with
$W_{v}^{\circ}$ and $\beta_{v}$ as in Lemma~\ref{lem:real-rs-evaluation}
and $\sigma_{t}$ as in \textsection\ref{def:real-spherical-parameter}.
For $\Re(s_{2})>-1/2$ and
$\Re(1/2+s_{1})>\max_{j}|\Re(\nu_{j})|$, we have
\begin{multline}\label{e14.30}
I_{\mathrm{gen},v}(\mathbf{s};W_{v})
=\frac{\beta_{v}}{8\pi^{2}}\int_{\mathbb{R}}
L_{v}(1/2+s_{2},\sigma_{t})L_{v}(1/2+s_{1},\pi_{v}\times\sigma_{t})\\
\mathcal{S}h_{v}(it)t\sinh(\pi t)\,dt.
\end{multline}
The identity extends meromorphically.
\end{lemma}

\begin{proof}
In \eqref{7.5} with \eqref{fc11.3} inserted, we substitute
$(y,z)\mapsto(y/a',z/z')$ at fixed $g=z'\diag(a',1)n(b')k'$.
This gives the phase $e^{2\pi iyb'}$ and weight $|a'|^{-s_{2}}$.
Inverting $a',z'$ in the $g$-integral gives \eqref{eq10.17}, with
factor $4$ from their signs.
\begin{align*}
I_{\mathrm{gen},v}(\mathbf{s};W_{v})
=4\int_{(\mathbb{R}^{\times})^{2}}W_{v}^{\circ}(\diag(yz,z,1))
\widehat{f}_{v}(y)|z|^{2s_{1}}|y|^{s_{1}+s_{2}}
\,d^{\times} z\,d^{\times} y.
\end{align*}
With \eqref{eq11.26}, the $y$-weight against $K_{it}(2\pi|y|)$ is
$|y|^{s_{1}-1/2}$. The $y,z$-integral is
$\beta_{v} L_{v}(1/2+s_{1},\pi_{v}\times\sigma_{t})/2$ by
\eqref{eq:real-rs-evaluated} at $\nu=it$, and the Gamma product in
\eqref{eq11.26} is $\pi^{1/2+s_{2}}L_{v}(1/2+s_{2},\sigma_{t})$.
To continue \eqref{e14.30}, we move its $t$-contour between the
pole sets for $+it$ and $-it$ in
$\Gamma_{\mathbb{R}}(1/2+s_{1}+\nu_{j}\pm it)$ and
$\Gamma_{\mathbb{R}}(1/2+s_{2}\pm it)$ and add the crossed residues.
\end{proof}

\subsubsection{A non-Archimedean double series}
\begin{lemma}\label{lemma14.5}
Write $\alpha_{k}=q_{v}^{-\nu_{k}}$ for the Satake parameters of $\pi_{v}$. Let $\ell\in\mathbb{Z}_{\geq0}$, $A=q_{v}^{-1-2s_{1}}$ and $B=q_{v}^{-1-s_{1}-s_{2}}$.
Assume
\begin{align*}
 |B|\max_{k}|\alpha_{k}|<1,
 \qquad |A|\max_{i<j}|\alpha_{i}\alpha_{j}|<1.
\end{align*}
Then the double series converges absolutely, and
\begin{align*}
\sum_{i=\ell}^{\infty}\sum_{j=0}^{\infty}
\lambda_{\pi_{v}}(\mathfrak{p}_{v}^{\,i+j},\mathfrak{p}_{v}^{\,j})B^{i}A^{j}
=
\frac{B^{\ell}\mathcal{N}(\mathbf{s},\pi_{v};\ell)}
     {\prod_{k=1}^{3}(1-\alpha_{k}B)\prod_{1\leq i<j\leq 3}(1-\alpha_{i}\alpha_{j} A)},
\end{align*}
where $h_{m}$ is the complete homogeneous polynomial in
$(\alpha_{1},\alpha_{2},\alpha_{3})$, with $h_{0}=1$ and $h_{m}=0$
for $m<0$.
\begin{multline}\label{14.36}
 \mathcal{N}(\mathbf{s},\pi_{v};\ell)=
 B h_{\ell+1}+(1-\lambda_{\pi_{v}}(\mathfrak{p}_{v})B-AB)h_{\ell}\\
 +(B^{2}-A+A\lambda_{\pi_{v}}(\mathfrak{p}_{v})B)h_{\ell-1}-AB^{2}h_{\ell-2}.
\end{multline}
\end{lemma}

\begin{proof}
For formal variables $A,B$, the Pieri rule with $\alpha_{1}\alpha_{2}\alpha_{3}=1$ gives
 \begin{equation}\label{eq:two-row-schur-generating}
 \sum_{i,j\geq0}s_{i+j,j,0}(\alpha_{1},\alpha_{2},\alpha_{3})B^{i}A^{j}
 =\frac{1-AB}
 {\prod_{k=1}^{3}(1-\alpha_{k}B)
  \prod_{1\leq k<m\leq3}(1-\alpha_{k}\alpha_{m}A)}.
\end{equation}
For distinct $\alpha_{k}$, sum the determinant \eqref{12.6} over
$i\geq\ell$ and $j\geq0$ and clear the denominator in
\eqref{eq:two-row-schur-generating}. After division by $B^{\ell}$,
sum the two orders of the remaining indices for each first index $a$ to get
\begin{align*}
 \sum_{a=1}^{3}\frac{\alpha_{a}^{\ell+2}}{\prod_{b\ne a}(\alpha_{a}-\alpha_{b})}
 (1-(\lambda_{\pi_{v}}(\mathfrak{p}_{v})-\alpha_{a})B+B^{2}/\alpha_{a})(1-A/\alpha_{a}).
\end{align*}
Expanding the two factors gives \eqref{14.36} by
$h_{m}=\sum_{a} \alpha_{a}^{m+2}/\prod_{b\ne a}(\alpha_{a}-\alpha_{b})$ for $m\geq-2$.
Repeated parameters follow by continuity, since the series converges
locally uniformly under the two strict bounds in the lemma.
\end{proof}

\subsubsection{Non-Archimedean integrals}
For $d_{v}=0$, $(\mathcal{J}_{\psi_{v}}f_{\boldsymbol{\nu},v}^{\circ})(I_{3})
=\mathbf{L}(\pi_{v})^{-1}$ by \eqref{eq:jacquet-normalization-transport}.
At $v\mid\mathfrak{q}$ with $d_{v}=0$, each translate in \eqref{c14.6} contributes
\eqref{e7.5}, so
\begin{equation}\label{14.37}
I_{\mathrm{gen},v}(\mathbf{s},\varphi_{\boldsymbol{\nu}},\mathbf{1},\mathbf{1})
=\frac{q_{v}^{r_{v}}\mathcal{L}_{v}(\mathbf{s},\pi_{v})}{\mathbf{L}(\pi_{v})}.
\end{equation}

\begin{lemma}\label{lem14.6}
Let $v\nmid\mathfrak{q}$ and $d_{v}=0$. Then
\begin{align*}
I_{\mathrm{gen},v}(\mathbf{s},\varphi_{\boldsymbol{\nu}},\mathbf{1},\mathbf{1})
=\frac{L_{v}(1+s_{1}+s_{2},\pi_{v})L_{v}(1+2s_{1},\widetilde{\pi}_{v})}{\mathbf{L}(\pi_{v})}
 \sum_{\ell=0}^{n_{v}}\frac{\mathcal{N}(\mathbf{s},\pi_{v};\ell)}
 {q_{v}^{(1/2+s_{2})n_{v}+\ell(s_{1}-s_{2})}} ,
\end{align*}	
where $\mathcal{N}(\mathbf{s},\pi_{v};\ell)$ is defined by \eqref{14.36}. 
\end{lemma}

\begin{proof}
In \eqref{7.5}, the character sum in \eqref{c14.5} contributes
$q_{v}^{\ell}\mathbf{1}_{i\geq\ell}$ for $i=e_{v}(y_{v})$. Lemma~\ref{lemma14.5}
sums over $i\geq\ell$ and $j=e_{v}(z_{v})\geq0$. Its factor $q_{v}^{-\ell(1+s_{1}+s_{2})}$
combines with \eqref{c14.5} and $q_{v}^{\ell}$ to give
$q_{v}^{-(1/2+s_{2})n_{v}-\ell(s_{1}-s_{2})}$.
\end{proof}

For $d_{v}>0$, scaling $y,z$ by $\varpi_{v}^{-d_{v}}$ in \eqref{7.5}
multiplies the $d_{v}=0$ formulas by $q_{v}^{d_{v}(3s_{1}+s_{2}-1/2+\nu_{3}-\nu_{1})}$
by \eqref{eq:finite-character-measures} and \eqref{eq:different-jacquet-scalar}.

\subsection{The Dual Generic Term}\label{sec14.7}

\begin{lemma}
Let $v$ be real and let $W_{v}$ be given by \eqref{fc11.3}, with
$W_{v}^{\circ}$ and $\beta_{v}$ as in Lemma~\ref{lem:real-rs-evaluation}.
For
\begin{align*}
 \min\{\Re(1+s_{2}-s_{1}),\Re(1+s_{1}+s_{2})\}
 >\max_{j}|\Re(\nu_{j})|,
\end{align*}
we have
\begin{multline*}
I_{\mathrm{gen},v}(\mathbf{s}^{\vee};\pi(w_{2})W_{v})
=\frac{\beta_{v}\mathcal{S}h_{v}(1/2+s_{2})}{\Gamma_{\mathbb{R}}(2+2s_{2})}\\
\prod_{j=1}^{3}
\Gamma_{\mathbb{R}}(1+s_{2}-s_{1}-\nu_{j})
\Gamma_{\mathbb{R}}(1+s_{1}+s_{2}+\nu_{j}).
\end{multline*}
\end{lemma}

\begin{proof}
In \eqref{7.5} at $(\mathbf{s}^{\vee},\pi(w_{2})W_{v})$ we insert
\eqref{fc11.3}. At fixed $g$ in the coordinates for \eqref{e14.30},
we substitute $(y,z)\mapsto(y/(a'z'),z'z)$.
The $g$-integral against $h_{v}(g^{-1})|a'|^{-s_{2}}$ gives
$4\widehat{f}_{v}(0)$ by \eqref{eq10.17}. Summing the signs of $(y,z)$
gives $16\widehat{f}_{v}(0)$ times the positive Mellin integral
\begin{equation}\label{eq:real-jacquet-double-mellin}
\int_{0}^{\infty}\int_{0}^{\infty}
W_{v}^{\circ}(\diag(yz,z,1))
z^{\lambda_{1}-1}y^{\lambda_{2}-1}\,d^{\times} z\,d^{\times} y=\frac{\beta_{v}G_{\boldsymbol{\nu}}(\lambda_{1},\lambda_{2})}{4\pi^{\lambda_{1}+\lambda_{2}}}.
\end{equation}
The Vinogradov--Takhtadzhyan formula gives this Mellin integral
with \eqref{eq:real-double-Mellin-Gamma} for
$\Re(\lambda_{i})>\max_{j}|\Re(\nu_{j})|$.
This gives $c_{v}=\beta_{v}\mathbf{L}(\pi_{v})$ in
\eqref{eq:jacquet-normalization-transport} at a real place.
The dual Mellin exponents give
$(\lambda_{1},\lambda_{2})=(1+s_{2}-s_{1},1+s_{1}+s_{2})$, and
\eqref{11.22} evaluates $\widehat{f}_{v}(0)$.
\end{proof}

\begin{lemma}
For $v<\infty$ and $d_{v}=0$, we have
\begin{itemize}
\item If $v\nmid\mathfrak{q}$, then
\begin{equation}\label{14.46}
I_{\mathrm{gen},v}(\mathbf{s}^{\vee},\pi(w_{2})\varphi_{\boldsymbol{\nu}},\mathbf{1},\mathbf{1})
=\frac{q_{v}^{(1/2+s_{2})n_{v}}\mathcal{L}_{v}(\mathbf{s}^{\vee},\pi_{v})}{\mathbf{L}(\pi_{v})}
\sum_{\ell=0}^{n_{v}}q_{v}^{-(1+2s_{2})\ell}.
\end{equation}
\item If $v\mid\mathfrak{q}$, then
\begin{multline}\label{14.45}
I_{\mathrm{gen},v}(\mathbf{s}^{\vee},\pi(w_{2})\varphi_{\boldsymbol{\nu}},\mathbf{1},\mathbf{1})
=\frac{\mathcal{N}(\mathbf{s}^{\vee},\pi_{v};r_{v})}
 {q_{v}^{(s_{1}+s_{2})r_{v}}\mathbf{L}(\pi_{v})}\\
 L_{v}(1+s_{1}+s_{2},\pi_{v})L_{v}(1+s_{2}-s_{1},\widetilde{\pi}_{v}).
\end{multline}
\end{itemize}
For $d_{v}>0$, multiply \eqref{14.46} or \eqref{14.45}
by $q_{v}^{d_{v}(2s_{2}-1/2+\nu_{3}-\nu_{1})}$.
\end{lemma}

\begin{proof}
For $v\nmid\mathfrak{q}$, the $q_{v}^{\ell}$ in \eqref{14.37} at
$(\mathbf{s},r_{v})=(\mathbf{s}^{\vee},\ell)$ multiplies the coefficient
in \eqref{c14.5}. The product is $q_{v}^{(1+2s_{2})\ell-(1/2+s_{2})n_{v}}$.
Reindexing by $n_{v}-\ell$ gives \eqref{14.46}.
For $v\mid\mathfrak{q}$, $w_{2}n_{13}(b)w_{2}=n_{12}(b)$ gives
$q_{v}^{r_{v}}\mathbf{1}_{e_{v}(y_{v})\geq r_{v}}$ in \eqref{7.5}, so
\eqref{14.45} follows from Lemma~\ref{lemma14.5} at
$(\mathbf{s},\ell)=(\mathbf{s}^{\vee},r_{v})$ with factor
$q_{v}^{r_{v}}/\mathbf{L}(\pi_{v})$.
For $d_{v}>0$, the change of character in \eqref{14.37} replaces
the Mellin factor by $q_{v}^{d_{v}(3s_{1}'+s_{2}')}=q_{v}^{2d_{v}s_{2}}$.
\end{proof}

\subsection{Prime-Level Oldform Contributions}
\label{app:prime-old-row}
\begin{lemma}\label{lem:prime-old-row-expansion}
Let $\pi_{v}$ be unramified generic on $\mathrm{PGL}_{3}(F_{v})$ and
$\sigma_{v}$ unitary generic on $\mathrm{PGL}_{2}(F_{v})$ with $r_{\sigma_{v}}\leq1$.
If $r_{\sigma_{v}}=0$, the local weight from
\textsection\ref{def:mixed-finite-weight} is
\begin{align*}
 \alpha_{v}(\sigma_{v};1)=\sum_{j\geq0}c_{v,j}(\pi_{v})
 \lambda_{\sigma_{v}}(\mathfrak{p}_{v}^{j}),
\end{align*}
where the series converges absolutely and
\begin{equation}\label{eq:mixed-prime-old-row-series}
\begin{split}
 &c_{v,0}(\pi_{v})=1+\lambda_{\pi_{v}}(\mathfrak{p}_{v})
 +\lambda_{\widetilde{\pi}_{v}}(\mathfrak{p}_{v}),\\
 &c_{v,j}(\pi_{v})=\bigl(2+\lambda_{\pi_{v}}(\mathfrak{p}_{v})
 +\lambda_{\widetilde{\pi}_{v}}(\mathfrak{p}_{v})\bigr)(-q_{v}^{-1/2})^{j},
 \qquad j\geq1,
\end{split}
\end{equation}
and $|c_{v,j}(\pi_{v})|\ll\alpha_{\pi}(\mathfrak{p}_{v})q_{v}^{-j/2}$.
For $\pi_{v}=\mathbf{1}\boxplus\mathbf{1}\boxplus\mathbf{1}$, this is
\begin{equation}\label{eq:prime-old-row-series}
 \sum_{j\geq0}c_{v,j}(\pi_{v})\lambda_{\sigma_{v}}(\mathfrak{p}_{v}^{j})
 =7+8\sum_{j\geq1}(-1)^{j}q_{v}^{-j/2}
 \lambda_{\sigma_{v}}(\mathfrak{p}_{v}^{j}).
\end{equation}
For $r_{\sigma_{v}}=1$,
$\alpha_{v}(\sigma_{v};1)=(1-q_{v}^{-2})L_{v}(1,\sigma_{v},\Ad)$.
\end{lemma}

\begin{proof}
For $r_{\sigma_{v}}=0$, write $t=q_{v}^{-1/2}$ and
$x=\lambda_{\sigma_{v}}(\mathfrak{p}_{v})$. Put
\begin{align*}
 &S_{i}^{\pi}(t):=\sum_{j\geq0}s_{j+i,0,0}(\alpha_{1},\alpha_{2},\alpha_{3})
 \lambda_{\sigma_{v}}(\mathfrak{p}_{v}^{j})t^{j},\qquad i=0,1,\\
 &A_{\pi}(t):=\sum_{m,j\geq0}
 s_{j+m+1,m+1,0}(\alpha_{1},\alpha_{2},\alpha_{3})
 \lambda_{\sigma_{v}}(\mathfrak{p}_{v}^{j})t^{2m+j}.
\end{align*}
In \eqref{equ12.11} at $(r_{v},r_{\sigma_{v}})=(1,0)$, the two terms
are $A_{\pi}$ for $(d_{1},d_{2})=(1,0)$ and
$(S_{0}^{\pi}+S_{1}^{\pi})/(1+t^{2}+tx)$ for $(0,1)$.

By the Cauchy identity,
\begin{equation}\label{eq:old-row-P0}
 \frac{S_{0}^{\pi}(t)}{L_{v}(1/2,\pi_{v}\times\sigma_{v})}
 =1-\lambda_{\widetilde{\pi}_{v}}(\mathfrak{p}_{v})t^{2}+xt^{3},
\end{equation}
\begin{equation}\label{eq:old-row-P1}
 \frac{S_{1}^{\pi}(t)}{L_{v}(1/2,\pi_{v}\times\sigma_{v})}
 =\lambda_{\pi_{v}}(\mathfrak{p}_{v})-\lambda_{\widetilde{\pi}_{v}}(\mathfrak{p}_{v})xt+(x^{2}-1)t^{2}.
\end{equation}
\begin{equation}\label{eq:old-row-Api}
 \frac{A_{\pi}(t)}{L_{v}(1/2,\pi_{v}\times\sigma_{v})}
 =\lambda_{\widetilde{\pi}_{v}}(\mathfrak{p}_{v})-xt.
\end{equation}
Substituting the normalized series \eqref{eq:old-row-P0}--\eqref{eq:old-row-Api} gives
\begin{equation}\label{eq:old-row-simplified}
 \alpha_{v}(\sigma_{v};1)=
 \frac{\lambda_{\pi_{v}}(\mathfrak{p}_{v})+\lambda_{\widetilde{\pi}_{v}}(\mathfrak{p}_{v})+2}
 {1+q_{v}^{-1}+q_{v}^{-1/2}\lambda_{\sigma_{v}}(\mathfrak{p}_{v})}-1.
\end{equation}
The expansion $(1+t^{2}+tx)^{-1}=\sum_{j\geq0}(-t)^{j}\lambda_{\sigma_{v}}(\mathfrak{p}_{v}^{j})$
gives \eqref{eq:mixed-prime-old-row-series}, with
\eqref{eq:prime-old-row-series} at $\lambda_{\pi_{v}}(\mathfrak{p}_{v})
=\lambda_{\widetilde{\pi}_{v}}(\mathfrak{p}_{v})=3$.

At $r_{\sigma_{v}}=r_{v}=1$, \eqref{fc12.16} gives $\alpha_{v}(\sigma_{v};1)
=q_{v}|W_{v}'^{\circ}(I_{2})|^{2}/(q_{v}+1)=(1-q_{v}^{-2})L_{v}(1,\sigma_{v},\Ad)$
by \eqref{eq:finite-newvector-value}.
\end{proof}

\subsubsection{Conductor exponent at most one}

\begin{lemma}\label{lem:depth-low}
With the notation of Lemma~\ref{lem:depth-kernel}, assume $r_{\sigma_{v}}\leq1$.
Suppose $|\lambda_{\sigma_{v}}(\mathfrak{p}_{v}^{n})|\leq(n+1)q_{v}^{\theta n}$ with
$0\leq\theta<1/2-\vartheta_{\pi_{v}}$.
For $m\geq1$ and $r\geq2$,
\begin{align*}
 |\mathcal{H}_{W_{v}^{[m]}}(\sigma_{v},\mathbf{0},\mathbf{1})|&\ll (1+m)^{5}q_{v}^{m\vartheta_{\pi_{v}}}|W_{v}'^{\circ}(I_{2})|^{2},\\
 \bigg|\sum_{m=2}^{r} b_{r-m}\mathcal{H}_{W_{v}^{[m]}}(\sigma_{v},\mathbf{0},\mathbf{1})\bigg|
 &\ll_{\varepsilon} q_{v}^{r(\vartheta_{\pi_{v}}+\varepsilon)}|W_{v}'^{\circ}(I_{2})|^{2} .
\end{align*}
The constants depend on $\vartheta_{\pi_{v}},\theta$.
\end{lemma}

\begin{proof}
For $r_{\sigma_{v}}\leq1$, the Satake bound and
\eqref{eq:oldform-basis-coefficients} give
$1-|\alpha_{\sigma_{v}}|^{2}\gg_{\theta}1$ and
$\sum_{i=0}^{l} q_{v}^{(i-l)/2}|\xi_{\sigma_{v}}(\mathfrak{p}_{v}^{i},\mathfrak{p}_{v}^{l})|
\ll_{\theta}1$. Also $q_{v}^{m}\Vol(K_{0,v}'[d_{2}])q_{v}^{-d_{1}}\leq1$ for
$d_{1}+d_{2}=m$.
The tableau and dimension formulas, with $\alpha_{1}\alpha_{2}\alpha_{3}=1$, give
\begin{equation}\label{eq:newform-Schur-bound}
 |s_{a,b,0}(\alpha_{1},\alpha_{2},\alpha_{3})|
 \ll(1+a)^{3}\bigl(\max_{i,\pm}|\alpha_{i}^{\pm1}|\bigr)^{a},
 \qquad a\geq b\geq0.
\end{equation}

Put $\eta=1/2-\vartheta_{\pi_{v}}-\theta>0$. At $(\mathbf{s},r_{v})=(\mathbf{0},m)$
in \eqref{equ12.11}, $k_{2}+d_{1}\leq m$ and $L_{v}(1/2,\sigma_{v})\ll_{\theta}1$.
Thus \eqref{eq:newform-Schur-bound} yields
\begin{multline*}
\frac{|\mathcal{H}_{W_{v}^{[m]}}(\sigma_{v},\mathbf{0},\mathbf{1})|}{|W_{v}'^{\circ}(I_{2})|^{2}}
 \ll_{\theta} (1+m)^{2}q_{v}^{\vartheta_{\pi_{v}}m}
 \sum_{l_{1},l_{2}\geq0}\frac{(1+m+l_{1}+l_{2})^{3}(1+l_{2})}{q_{v}^{(1-\vartheta_{\pi_{v}})l_{1}+\eta l_{2}}}\\
 \ll_{\eta,\theta}(1+m)^{5}q_{v}^{\vartheta_{\pi_{v}}m}.
\end{multline*}
Summing this over $2\leq m\leq r$ with
$|b_{r-m}|\leq64q_{v}^{\vartheta_{\pi_{v}}(r-m)}$ gives the second estimate.
\end{proof}

\subsubsection{Local weights at the residue points}

\begin{lemma}\label{lem:original-moving-residue-weight}
Suppose $r_{v}\geq1$. At the positive poles in \eqref{14.125} with
$1-\alpha_{\sigma_{\lambda,v}}^{2}\ne0$, the quotient
\begin{align*}
 \frac{\mathcal{H}_{W_{v}^{[r_{v}]}}(\sigma_{\lambda,v},\mathbf{s},\mathbf{1})}
 {L_{v}(1/2+s_{1},\pi_{v}\times\sigma_{\lambda,v})L_{v}(1/2+s_{2},\sigma_{\lambda,v})}
\end{align*}
agrees with a holomorphic function on
$|\Re(s_{j})|,|\Re(\nu_{i})|<\delta$, $|\Re(\lambda)-1/2|<3\delta$,
with bound $O_{\delta}(q_{v}^{O(\delta r_{v})})$ for small $\delta>0$.
Reflect $\lambda$ for the negative points.
\end{lemma}

\begin{proof}
Summing $l_{2}$ in \eqref{equ12.11} by \eqref{12.6} and dividing
by $L_{v}(1/2+s_{1},\pi_{v}\times\sigma_{\lambda,v})$ leaves, in determinant
term $i$, the factors $1-q_{v}^{-1/2-s_{1}-\nu_{j}\pm\lambda}$ for $j\ne i$.
The polynomial in $q_{v}^{-s_{1}}$ has degree at most four, so its
coefficients depend only on the terms with $l_{2}\leq4$, also for
repeated Satake parameters.

For $l\geq1$, divide $\Psi_{v}(1/2+s_{2},W_{l,v}',\mathbf{1})$ by
$\Psi_{v}(1/2+s_{2},W_{v}'^{\circ},\mathbf{1})$ and clear the square-root
denominator in \eqref{eq:oldform-basis-coefficients} to get
\begin{align*}
 \begin{cases}
 q_{v}^{-s_{2}}-\alpha_{\sigma_{\lambda,v}}, & l=1,\\
 q_{v}^{-(l-2)s_{2}}(q_{v}^{-s_{2}}-q_{v}^{\lambda-1/2})
 (q_{v}^{-s_{2}}-q_{v}^{-\lambda-1/2}), & l\geq2 .
 \end{cases}
\end{align*}
In determinant term $i$ the oldform factor is the same expression
with $s_{2}$ replaced by $s_{1}+\nu_{i}$. At fixed $d_{1},d_{2},l,l_{1}$ with $l\geq1$,
let $N(\lambda)$ be the $q_{v}^{-l}\sum_{k_{1},k_{2},l_{2}}$ contribution in
\eqref{equ12.11} after the division by $L_{v}(1/2+s_{1},\pi_{v}\times\sigma_{\lambda,v})$
and clearing the common denominator
$(1-\alpha_{\sigma_{\lambda,v}}^{2})(1-q_{v}^{-2})^{\mathbf{1}_{l\geq2}}$ from
\eqref{eq:oldform-basis-coefficients}, with $q_{v}^{-2d_{1}(1/2+s_{1})}$ kept outside.

Put $M(\lambda)=(q_{v}^{\lambda-1/2}-q_{v}^{-s_{2}})
\prod_{i=1}^{3}(q_{v}^{\lambda-1/2}-q_{v}^{-s_{1}-\nu_{i}})$.
At $\lambda_{0}=1/2$ or $1/2+\pi i/\log q_{v}$, we have
$\alpha_{\sigma_{\lambda_{0},v}}=q_{v}^{\lambda_{0}-1/2}=\pm1$.
The standard numerator is divisible by $q_{v}^{\lambda_{0}-1/2}-q_{v}^{-s_{2}}$.
If $q_{v}^{-s_{1}-\nu_{j}}=q_{v}^{\lambda_{0}-1/2}$, the oldform factor
vanishes in determinant term $j$ and $1-q_{v}^{-1/2-s_{1}-\nu_{j}+\lambda_{0}}$
vanishes in every other term. Hence $M(\lambda_{0})$ divides
$N(\lambda_{0})$, also for repeated Satake parameters by polynomial
division in $q_{v}^{-s_{1}}$ with leading coefficient $-1$ in the product
over $i$. Replace the numerator by
\begin{equation}\label{eq:residue-weight-interpolation}
 N(\lambda)-M(\lambda)\sum_{\lambda_{0}\in\{1/2,\,1/2+\pi i/\log q_{v}\}}
 \frac{1+q_{v}^{\lambda+\lambda_{0}-1}}{2}
 \frac{N(\lambda_{0})}{M(\lambda_{0})}.
\end{equation}
At each denominator zero in the strip, $q_{v}^{\lambda-1/2}=\pm1$.
In \eqref{eq:residue-weight-interpolation} the sum selects the
$\lambda_{0}$ with the same sign, so the displayed difference is zero.
Away from the denominator zeros, the term subtracted from $N(\lambda)$
in \eqref{eq:residue-weight-interpolation} vanishes wherever $M(\lambda)=0$.
The residues agree for generic parameters, and their sum extends through
collisions by \eqref{eq:residue-cluster}.

By \eqref{eq:newform-Schur-bound} with $d_{1}+k_{2}\leq r_{v}$ and
$h_{l_{2}}(q_{v}^{\lambda-1/2},q_{v}^{-\lambda-1/2})\ll q_{v}^{O(\delta)}$
for $l_{2}\leq4$, the coefficient bound is
$O((1+r_{v}+l_{1})^{3}\allowbreak q_{v}^{\delta(r_{v}+l_{1})+O(\delta)})$.
For $l\geq1$ we have $l_{1}=0$. On $|\Re(s_{j})|=2\delta$ ($j=1,2$) with
$|\Re(\nu_{i})|\leq\delta$, we have $|M(\lambda_{0})|\gg_{\delta}1$.
On $|\Re(\lambda)-1/2|=3\delta$, we have
$|1-\alpha_{\sigma_{\lambda,v}}^{2}|\gg_{\delta}1$. Thus \eqref{eq:residue-weight-interpolation}
divided by $(1-\alpha_{\sigma_{\lambda,v}}^{2})(1-q_{v}^{-2})^{\mathbf{1}_{l\geq2}}$
is $O_{\delta}((1+r_{v})^{3}q_{v}^{O(\delta r_{v})})$ on these boundaries
simultaneously and inside by the maximum modulus principle.
For $l=0$, the factor $q_{v}^{-(1+2s_{1})l_{1}}$ makes the $l_{1}$ sum
geometrically convergent.
In \eqref{equ12.11} the measure factor satisfies
$q_{v}^{r_{v}-d_{1}}\Vol(K_{0,v}'[d_{2}])\leq1$,
$|q_{v}^{-2d_{1}s_{1}}|\leq q_{v}^{2\delta r_{v}}$, and
$|W_{v}'^{\circ}(I_{2})|^{2}=1$ by \eqref{eq:finite-newvector-value}.
Summing over $0\leq l\leq d_{2}\leq r_{v}$ with $d_{1}=r_{v}-d_{2}$ proves the bound.
\end{proof}

\subsection{Steinberg Signs in the Squarefree Sieve}
\label{app:squarefree-signed-section}

Suppose $\psi_{v}$ has conductor $\mathcal{O}_{v}$.
Let $f_{\boldsymbol{\nu},v}$ be the section
\eqref{c14.6} with $r_{v}=1$, and set
\begin{align*}
 w_{\mathfrak{p}_{v}}=\begin{pmatrix}&1\\-\varpi_{v}&\end{pmatrix},
 \qquad f_{v}^{\mathrm{AL}}=-q_{v}^{-1}\pi_{v}(\diag(w_{\mathfrak{p}_{v}},1))f_{\boldsymbol{\nu},v}.
\end{align*}
Applying the normalized spherical Hecke operator to $f_{\boldsymbol{\nu},v}^{\circ}$ gives
\begin{multline}\label{eq:signed-section-Hecke-identity}
 f_{v}^{\mathrm{AL}}=q_{v}^{-1}\sum_{a,b\bmod\mathfrak{p}_{v}}
 \pi_{v}\begin{pmatrix}\varpi_{v}&b&a\\&1&\\&&1\end{pmatrix}f_{\boldsymbol{\nu},v}^{\circ}\\
 {}-\lambda_{\pi_{v}}(\mathfrak{p}_{v})f_{\boldsymbol{\nu},v}^{\circ}
 +q_{v}^{-1}\pi_{v}(\diag(1,1,\varpi_{v}))f_{\boldsymbol{\nu},v}^{\circ}.
\end{multline}
Put $f_{v}^{\mathrm{AL},\vee}:=\pi_{v}(w_{2})f_{v}^{\mathrm{AL}}$. The double sum in
\eqref{eq:signed-section-Hecke-identity} is fixed by $\pi_{v}(w_{2})$, so
\begin{equation}\label{eq:signed-section-dual-identity}
 f_{v}^{\mathrm{AL},\vee}-f_{v}^{\mathrm{AL}}
 =q_{v}^{-1}\bigl(\pi_{v}(\diag(1,\varpi_{v},1))
 -\pi_{v}(\diag(1,1,\varpi_{v}))\bigr)f_{\boldsymbol{\nu},v}^{\circ}.
\end{equation}

In the first period of \eqref{fc10.1} at $r_{\sigma_{v}}=1$,
the change $g'\mapsto x_{v}w_{\mathfrak{p}_{v}}^{-1}$ gives the
Atkin--Lehner factor
\begin{equation}\label{eq:signed-newvector-factor}
 \mathcal{H}_{\mathcal{J}_{\psi_{v}}f_{v}^{\mathrm{AL}}}(\sigma_{v},\mathbf{s},\mathbf{1})
 =q_{v}^{s_{1}-1/2}\lambda_{\sigma_{v}}(\mathfrak{p}_{v})
 \mathcal{H}_{\mathcal{J}_{\psi_{v}}f_{\boldsymbol{\nu},v}}(\sigma_{v},\mathbf{s},\mathbf{1}).
\end{equation}
The translates in \eqref{eq:signed-section-dual-identity} have zero
weight for ramified $\sigma_{v}$. Thus \eqref{eq:signed-newvector-factor}
gives $q_{v}^{-1/2}\lambda_{\sigma_{v}}(\mathfrak{p}_{v})$ for $q_{v}^{-s_{1}}f_{v}^{\mathrm{AL}}$,
with an extra $q_{v}^{s_{1}'-s_{1}}$ for $q_{v}^{-s_{1}}f_{v}^{\mathrm{AL},\vee}$ at $\mathbf{s}^{\vee}$.
Both weights vanish for $r_{\sigma_{v}}\geq2$.

For unramified $\sigma_{v}$, put
\begin{align*}
f_{i}:=\frac{\Psi_{v}(1/2+s_{1},\mathcal{J}_{\psi_{v}}f_{\boldsymbol{\nu},v},
\sigma_{v}(\diag(1,\varpi_{v}^{i}))W_{v}'^{\circ})}
{\Psi_{v}(1/2+s_{1},\mathcal{J}_{\psi_{v}}f_{\boldsymbol{\nu},v}^{\circ},W_{v}'^{\circ})},
 \qquad i=0,1.
\end{align*}
At $r_{v}=1$ and $t=q_{v}^{-1/2-s_{1}}$, the $(d_{1},d_{2})=(1,0),(0,1)$
contributions before the oldform basis sum in \eqref{equ12.11} give
\begin{align*}
 f_{i} L_{v}(1/2+s_{1},\pi_{v}\times\sigma_{v})
 =\alpha_{\sigma_{v}}^{i} q_{v}^{-2s_{1}}A_{\pi}(t)
   +\frac{q_{v}^{-is_{1}}}{1+q_{v}^{-1}}S_{i}^{\pi}(t).
\end{align*}
Here $A_{\pi}$ sums over $m,j\geq0$ and $S_{i}^{\pi}$ over $j\geq0$.
Dividing these series by $L_{v}(1/2+s_{1},\pi_{v}\times\sigma_{v})$ gives the
right sides of \eqref{eq:old-row-Api} for $A_{\pi}$,
\eqref{eq:old-row-P0} for $S_{0}^{\pi}$, and \eqref{eq:old-row-P1} for $S_{1}^{\pi}$,
all at this $t$.
For the second period in \eqref{fc10.1}, the same $\mathrm{GL}_{2}$
translate multiplies $\Psi_{v}(1/2+s_{2},W_{v}'^{\circ},\mathbf{1})$ by $q_{v}^{-is_{2}}$.
Since $\sigma_{v}(w_{\mathfrak{p}_{v}})$ interchanges the two $\mathrm{GL}_{2}$
vectors in the numerators, \eqref{eq:oldform-basis-coefficients} at $m=0,1$ gives
\begin{equation}\label{eq:signed-old-row-parameters}
 \frac{\mathcal{H}_{\mathcal{J}_{\psi_{v}}f_{v}^{\mathrm{AL}}}}
 {\mathcal{H}_{\mathcal{J}_{\psi_{v}}f_{\boldsymbol{\nu},v}^{\circ}}}
 =-q_{v}^{s_{1}-1}\frac{(1-\alpha_{\sigma_{v}} q_{v}^{-s_{2}})f_{1}+(q_{v}^{-s_{2}}-\alpha_{\sigma_{v}})f_{0}}{1-\alpha_{\sigma_{v}}^{2}}.
\end{equation}
At $\mathbf{s}=\mathbf{0}$, both translate factors in
\eqref{eq:signed-section-dual-identity} equal $1$. For $r_{\sigma_{v}}=0$,
\eqref{fc11.42} at $n_{v}=0$ cancels the common Jacquet scalar and the two
local $L$-factors in the ratios, so
\begin{equation}\label{eq:signed-old-row}
\frac{\mathcal{H}_{\mathcal{J}_{\psi_{v}}f_{v}^{\mathrm{AL}}}(\sigma_{v},\mathbf{0},\mathbf{1})}
 {\mathcal{H}_{\mathcal{J}_{\psi_{v}}f_{\boldsymbol{\nu},v}^{\circ}}(\sigma_{v},\mathbf{0},\mathbf{1})}
 =\frac{\mathcal{H}_{\mathcal{J}_{\psi_{v}}f_{v}^{\mathrm{AL},\vee}}(\sigma_{v},\mathbf{0},\mathbf{1})}
 {\mathcal{H}_{\mathcal{J}_{\psi_{v}}f_{\boldsymbol{\nu},v}^{\circ}}(\sigma_{v},\mathbf{0},\mathbf{1})}
 =-q_{v}^{-1}\alpha_{v}(\sigma_{v};1).
\end{equation}

At either $\lambda_{0}$, the fraction in \eqref{eq:signed-old-row-parameters}
has numerator
$-q_{v}^{\lambda_{0}-1/2}M(\lambda_{0})/(1+q_{v}^{-1})$.
Thus \eqref{eq:residue-weight-interpolation} adds
$q_{v}^{\lambda-1/2}M(\lambda)/(1+q_{v}^{-1})$ to that numerator.
With this numerator, $q_{v}^{-s_{1}}$ times \eqref{eq:signed-old-row-parameters} is
$O_{\delta}(q_{v}^{-1+O(\delta)})$ on $|\Re(\lambda)-1/2|=3\delta$
and throughout the strip by the maximum modulus principle.

By \eqref{12.5}, the Jacquet integral of \eqref{eq:signed-section-Hecke-identity}
at $\diag(\varpi_{v}^{i+j},\varpi_{v}^{j},1)$ is $0$ for $i=0$ and
$-q_{v}^{-i-j}s_{i+j,j+1,0}/\mathbf{L}(\pi_{v})$ for $i\geq1$.
Thus \eqref{7.5} gives
\begin{equation}\label{eq:signed-generic-series}
 I_{\mathrm{gen},v}(\mathbf{s};\mathcal{J}_{\psi_{v}}f_{v}^{\mathrm{AL}})
 =-\mathbf{L}(\pi_{v})^{-1}\sum_{i=1}^{\infty}\sum_{j=0}^{\infty}
 \frac{s_{i+j,j+1,0}(q_{v}^{-\nu_{1}},q_{v}^{-\nu_{2}},q_{v}^{-\nu_{3}})}
 {q_{v}^{i(1+s_{1}+s_{2})+j(1+2s_{1})}}.
\end{equation}
By \eqref{eq:newform-Schur-bound} at $(a,b)=(i+j,j+1)$, with
$j=n-i$ and $1\leq i\leq n$,
\begin{align*}
|I_{\mathrm{gen},v}(\mathbf{s};\mathcal{J}_{\psi_{v}}f_{v}^{\mathrm{AL}})|
\ll\sum_{n\geq1}n(1+n)^{3}q_{v}^{-(1-\vartheta_{\pi_{v}}-O(\varepsilon))n}
\ll q_{v}^{-1+\vartheta_{\pi_{v}}+O(\varepsilon)}. 
\end{align*}
For the fixed cuspidal data at $\mathbf{s}=\mathbf{0}$, this is
$O(q_{v}^{-1}\alpha_{\pi}(\mathfrak{p}_{v}))$.
In \eqref{7.5}, substitute $(y,z)\mapsto(\varpi_{v} y,z/\varpi_{v})$ for
the first translate in \eqref{eq:signed-section-dual-identity} and
$z\mapsto\varpi_{v} z$ for the second. Since the central character is trivial,
\begin{equation}\label{eq:signed-dual-generic-transport}
 I_{\mathrm{gen},v}\bigl(\mathbf{s};
 \mathcal{J}_{\psi_{v}}(f_{v}^{\mathrm{AL},\vee}-f_{v}^{\mathrm{AL}})\bigr)
 =q_{v}^{-1}(q_{v}^{s_{1}-s_{2}}-q_{v}^{-2s_{1}})
 I_{\mathrm{gen},v}(\mathbf{s};\mathcal{J}_{\psi_{v}}f_{\boldsymbol{\nu},v}^{\circ}).
\end{equation}
On the minimal Eisenstein contours, \eqref{eq:signed-dual-generic-transport}
at $\mathbf{s}^{\vee}$ with $\vartheta_{\pi_{v}}=O(\varepsilon)$ contributes
$O(q_{v}^{-1+O(\varepsilon)})$ after multiplication by $q_{v}^{-s_{1}}$.

For the $w_{2}$ degenerate integral, Iwasawa decomposition gives
\begin{align*}
 &f_{\boldsymbol{\nu},v}^{\circ}(n_{31}(a)n_{32}(c)w_{2}w_{1})=\phi(a,c)
 =\max(1,|a|_{v})^{-1-\nu_{1}+\nu_{2}}
 \max(1,|a|_{v},|c|_{v})^{-1-\nu_{2}+\nu_{3}},\\
 &f_{v}^{\mathrm{AL}}(n_{31}(a)n_{32}(c)w_{2}w_{1})=-q_{v}^{-2-\nu_{1}}\sum_{\alpha\bmod\mathfrak{p}_{v}}
 \phi(\varpi_{v}a,c-\alpha a).
\end{align*}
If $n=e_{v}(z)\geq0$, $m=\max(0,-e_{v}(a))\leq n$, Fourier integration on
balls gives
\begin{align*}
 K_{z}(a):=\int\phi(a,c)\psi_{v}(-cz)\,dc
 =(1-q_{v}^{-1-\nu_{2}+\nu_{3}})q_{v}^{(-1-\nu_{1}+\nu_{3})m}
 \sum_{j=0}^{n-m}q_{v}^{(\nu_{3}-\nu_{2})j},
\end{align*}
and is zero otherwise. Changing $a'=\varpi_{v}a$, $c'=c-\alpha a$
in $f_{v}^{\mathrm{AL}}$ gives its two-variable Fourier coefficient
\begin{align*}
 \widehat{f_{v}^{\mathrm{AL}}}(y,z)=-q_{v}^{-\nu_{1}}\int K_{z}(a')\psi_{v}(-a'y/\varpi_{v})
 \mathbf{1}_{a'z\in\mathfrak{p}_{v}}\,da'.
\end{align*}
Since $K_{z}(a')$ is supported on $a'z\in\mathcal{O}_{v}$, impose the
indicator by subtracting the integral over $a'z\in\mathcal{O}_{v}^{\times}$.
Integrating in $y,z$ with the Mellin powers in \eqref{eq:local-w2-degenerate} gives
\begin{multline}\label{eq:signed-degenerate-ratio}
 \frac{I_{v,w_{2}}(\mathbf{s};f_{v}^{\mathrm{AL}})}
 {I_{v,w_{2}}(\mathbf{s};f_{\boldsymbol{\nu},v}^{\circ})}
 =-q_{v}^{s_{1}-s_{2}}\bigg[1-
 \frac{(1-q_{v}^{-1+s_{2}-s_{1}-\nu_{1}})
 (1-q_{v}^{-1-\nu_{1}-\nu_{3}-2s_{1}})}
 {1-q_{v}^{-2-2\nu_{1}-2s_{1}}}\\
 (1-q_{v}^{-1-\nu_{1}-\nu_{2}-2s_{1}})\bigg].
\end{multline}
In \eqref{eq:local-w2-degenerate} we rescale $(a,y)$ for the first
translate in \eqref{eq:signed-section-dual-identity} and $(c,z)$ for the second,
both by $(\varpi_{v}^{-1},\varpi_{v})$. The additive measure gives $q_{v}$, so
\begin{equation}\label{eq:signed-dual-degenerate-ratio}
 \frac{I_{v,w_{2}}(\mathbf{s};f_{v}^{\mathrm{AL},\vee})}
 {I_{v,w_{2}}(\mathbf{s};f_{\boldsymbol{\nu},v}^{\circ})}
 =\frac{I_{v,w_{2}}(\mathbf{s};f_{v}^{\mathrm{AL}})}
 {I_{v,w_{2}}(\mathbf{s};f_{\boldsymbol{\nu},v}^{\circ})}
 +q_{v}^{-1}(q_{v}^{s_{1}-s_{2}}-q_{v}^{-2s_{1}}).
\end{equation}

The spherical intertwining scalars cancel for $w_{1}w_{2}$ and $w_{1}w_{2}w_{1}$.
Their quotients are \eqref{eq:signed-degenerate-ratio} at $\mathbf{s}$
and \eqref{eq:signed-dual-degenerate-ratio} at $\mathbf{s}^{\vee}$,
with $\boldsymbol{\nu}$ replaced by $(\nu_{2},\nu_{1},\nu_{3})$ and
$(\nu_{3},\nu_{1},\nu_{2})$ for the two cells, respectively.
All six quotients are holomorphic and $O(q_{v}^{-1+O(\delta)})$ for
$|\Re(s_{j})|,|\Re(\nu_{i})|<\delta$ with $\delta$ small, since the
denominator $1-q_{v}^{-2-2\nu_{1}-2s_{1}}$ is bounded away from zero.
The factors in \eqref{eq:degen-different-check} and
\eqref{eq:longest-different-scaling} also cancel.

For $\pi_{v}=\mathbf{1}\boxplus\mathbf{1}\boxplus\mathbf{1}$, the generic
integrals, degenerate quotients and residue weights after
\eqref{eq:residue-weight-interpolation} are $O(q_{v}^{-1+O(\varepsilon)})$
on the circles of \textsection\ref{sec15} with $q_{v}^{-s_{1}}f_{v}^{\mathrm{AL}}$ on
the original side. On the dual side,
$q_{v}^{-s_{1}}$ stays outside the evaluation at $\mathbf{s}^{\vee}$.
The sum of the residues is unchanged by \eqref{eq:residue-weight-interpolation}.
By \eqref{eq:signed-section-dual-identity}, the dual quotient
is $q_{v}^{-s_{1}}$ times \eqref{eq:signed-old-row-parameters} at $\mathbf{s}^{\vee}$,
plus $q_{v}^{-1-s_{1}}(q_{v}^{s_{1}'-s_{2}'}-q_{v}^{-2s_{1}'})$.

\subsubsection{Central values of the mixed local weights}

\begin{lemma}
\label{lem:mixed-standard-pole-finite}
Suppose $\pi_{v}$ is self-dual, unramified and unitary on $\mathrm{PGL}_{3}(F_{v})$.
The weights $\alpha_{v}(\sigma_{\lambda,v};m)$ and the quotient
\eqref{eq:signed-old-row} are holomorphic near $\lambda=\pm1/2$.
At these points, with $g_{\pi}$ from \textsection\ref{def:mixed-residue-coefficient},
\begin{equation}\label{eq:mixed-finite-endpoint}
 \alpha_{v}(\sigma_{\pm1/2,v};m)
 =\frac{\mathcal{N}(\mathbf{0},\pi_{v};m)}{1-q_{v}^{-2}},\qquad
 \sum_{m=0}^{r} \alpha_{v}(\sigma_{\pm1/2,v};m)
 =g_{\pi}(\mathfrak{p}_{v}^{r}).
\end{equation}
\end{lemma}

\begin{proof}
At $\mathbf{s}=\mathbf{0}$, the standard numerator in the proof of
Lemma~\ref{lem:original-moving-residue-weight} cancels $1-q_{v}^{2\lambda-1}$
for $l\geq1$. The same $l=0$ coefficient estimate gives convergence
since $\max_{i}|\alpha_{i}^{\pm1}|<q_{v}^{1/2}$.
For \eqref{eq:signed-old-row}, the denominator in
\eqref{eq:old-row-simplified} is $2(1+q_{v}^{-1})$ at $\lambda=\pm1/2$.

For distinct Satake parameters $(\alpha,\alpha^{-1},1)$, only the
determinant term for the Satake parameter $1$ can contribute in
Lemma~\ref{lem:original-moving-residue-weight} at $(s_{1},\lambda)=(0,1/2)$.
Its oldform factor vanishes for $l\geq1$. In \eqref{equ12.11} at $r_{v}=m$, the
$d_{2}=0$ tail in $l_{1}$ and the terms $1\leq d_{2}\leq m$ with $l_{1}=0$
give the first two terms and the finite sum, respectively, in
\begin{multline*}
\alpha_{v}(\sigma_{1/2,v};m)
=\ h_{m}(\alpha,\alpha^{-1})-q_{v}^{-1}h_{m-1}(\alpha,\alpha^{-1})\\
+\frac{q_{v}}{q_{v}+1}(1-\alpha/q_{v})(1-\alpha^{-1}/q_{v})
 \sum_{d=0}^{m-1}h_{d}(\alpha,\alpha^{-1}),\qquad m\geq1.
\end{multline*}

The identity $h_{m}(\alpha,\alpha^{-1},1)=\sum_{d=0}^{m}h_{d}(\alpha,\alpha^{-1})$
reduces this sum to \eqref{14.36} at $A=B=q_{v}^{-1}$ divided by $1-q_{v}^{-2}$
(with $m=0$ from \eqref{fc11.42} at $n_{v}=0$).
Continuity handles repeated parameters and $\sigma_{-\lambda,v}\simeq\sigma_{\lambda,v}$
gives the value at $-1/2$.
At $\lambda=1/2-s_{2}$, the $l\geq2$ terms vanish. The quotients in
\eqref{eq:signed-old-row-parameters} give the $l=1$ term via
\begin{align*}
 \frac{(q_{v}^{-s_{2}}-\alpha_{\sigma_{\lambda,v}},\,
          1-\alpha_{\sigma_{\lambda,v}}q_{v}^{-s_{2}})}
 {1-\alpha_{\sigma_{\lambda,v}}^{2}}
 =\frac{1+q_{v}^{-1}}{1-q_{v}^{-2+2s_{2}}}(-q_{v}^{-1+s_{2}},1).
\end{align*}
At $(s_{1},\lambda)=(0,1/2)$ the $l=1$ Rankin--Selberg numerator vanishes
by the determinant calculation, giving \eqref{eq:mixed-finite-endpoint}
also along the moving pole. At depth one this gives $f_{1}=f_{0}$
and leaves only the $f_{0}$ term in the unsigned weight, so
\eqref{eq:signed-old-row-parameters} tends to
$-q_{v}^{-1}f_{0}=-q_{v}^{-1}\alpha_{v}(\sigma_{1/2,v};1)$.
The dual limit is the same since both translate factors from
\eqref{eq:signed-section-dual-identity} tend to $1$.
\end{proof}

\subsection{Bounds for the Generic and Degenerate Terms}
\label{app:generic-degenerate-estimates}

For small $\delta>0$, take $|\Re(s_{j})|\leq\delta$ and
$|\Re(\nu_{i})|\leq\delta$. Then $|h_{m}|\leq\binom{m+2}{2}q_{v}^{m\delta}$
in \eqref{14.36}.
Lemma~\ref{lem14.6} bounds the $v\mid\mathfrak{n}$ factor in \eqref{e14.125}
by $O((1+n_{v})^{3}\allowbreak q_{v}^{O(\delta n_{v})})$, and \eqref{14.45} bounds
the $v\mid\mathfrak{q}$ factor in \eqref{e14.126} by
$O((1+r_{v})^{2}q_{v}^{O(\delta r_{v})})$. The geometric sum in \eqref{14.46}
is $O(1)$. Multiplication gives \eqref{15.10}.

For the local quotients in \eqref{e14.123}, divide
$I_{v,w}(\mathbf{s};f_{v})$ by $\mathcal{L}_{v}(\mathbf{s},\boldsymbol{\nu};w)/\mathbf{L}(\pi_{v})$.
In \eqref{e14.21}, at fixed $n=i'+j+d_{v}$, sum $r_{v}\leq j+d_{v}\leq n$,
$0\leq j'\leq j+d_{v}$, $0\leq i\leq n-j'$. In \eqref{eq14.67}, at $n=j+k$,
sum $0\leq i,j\leq n$. In \eqref{fc14.39}, fix $i+l+d_{v}$ and $j+l$
with minimum $m$ and maximum $n$.
The $l+d_{v}\geq\max(k,m)$ and $k\geq0$ sums, normalized by
$\zeta_{v}(s_{2}-s_{1}-\nu_{2})\zeta_{v}(\nu_{1}-\nu_{2})$, give
\begin{align*}
 q_{v}^{-m(\nu_{1}-\nu_{2})}
 +\zeta_{v}(\nu_{1}-\nu_{2})^{-1}
 \sum_{k=0}^{m-1}q_{v}^{-m(s_{2}-s_{1}-\nu_{2})-k(\nu_{1}+s_{1}-s_{2})}.
\end{align*}
The normalized sums in \eqref{e14.21}, \eqref{eq14.67}, and \eqref{fc14.39}
are bounded by
\begin{align*}
\sum_{n\geq0}(1+n)^{3}q_{v}^{-\max(0,n-r_{v})+O(\delta(n+r_{v}))}
\allowbreak\ll_{F}(1+r_{v})^{4}q_{v}^{O(\delta r_{v})}. 
\end{align*}

In Lemma~\ref{lem14.4}, dividing the $r$-sum by $\zeta_{v}(s_{2}-s_{1}-\nu_{1})$
leaves the term at $r=\max(0,\ell-j')-d_{v}$. In \eqref{fc14.48}, put
$h=j+d_{v}-i$ after summing the geometric series in $k$. After division
by $q_{v}^{\ell}$ and normalization by the zeta factors, the term from
the upper limit of the $k$-sum contributes
$q_{v}^{-\ell(s_{2}-s_{1}-\nu_{2})}$ on $h\geq\ell$. For $0\leq h<\ell$,
the remaining factor $\zeta_{v}(s_{2}-s_{1}-\nu_{2})^{-1}$ is bounded and the
terms decay as $q_{v}^{-(1-O(\delta))(\ell-h)+O(\delta\ell)}$.
The term from the lower limit decays in $j+d_{v}$. The $(i,j)$-sum in \eqref{e14.72}
converges with decay $q_{v}^{-\max(i,j)+O(\delta(i+j+\ell))}$.
Dividing out the coefficient from \eqref{c14.5} leaves the level factor
$q_{v}^{\ell}$ and bound $O_{F,\delta}(q_{v}^{\ell+O(\delta\ell)})$.

The dual quotient in \eqref{e14.124} divides the local integral by
$\mathcal{L}_{v}(\mathbf{s}^{\vee},\boldsymbol{\nu};w)/\mathbf{L}(\pi_{v})$.
With \eqref{eq:finite-dual-degenerate-transport}, the bounds up to
$O_{F,\delta}(q_{v}^{O(\delta(r_{v}+n_{v}))})$ are
\begin{equation}\label{eq:finite-degenerate-estimates}
 \begin{array}{c|cc}
 &\text{Original}&\text{Dual}\\ \hline
 v\mid\mathfrak{q}&1&q_{v}^{r_{v}}\\
 v\nmid\mathfrak{q}&q_{v}^{n_{v}/2}&q_{v}^{-n_{v}/2}
 \end{array}
\end{equation}
With the real factors \eqref{14.4}, \eqref{c14.4}, and
\eqref{eq:real-weyl-reductions}, these bounds in
\eqref{e14.123}--\eqref{e14.124} give the degenerate estimate in \eqref{15.11}.

\subsubsection{Mixed-moment test data and bounds}\label{app:mixed-mass-nonspectral}

At real $v$, take a real unit Kirillov vector $u_{v}\in C_{c}^{\infty}(\mathbb{R}_{>0})$
and unit-mass bumps $h_{0,v}$ (fixed, central) and $h_{v}\geq0$ on $K_{v}'$
supported within $R_{v}^{-M}$ of the identity. For $y_{v},z_{v}>0$, set
\begin{align*}
&V_{v}(\diag(z_{v}a(y_{v})k_{v},1))
=y_{v}u_{v}(y_{v})h_{0,v}(z_{v})(h_{v}^{*}\ast h_{v})(k_{v}),\\
&W_{v}=\int_{\mathbb{R}} f_{v}(b)\pi_{v}(n_{12}(-b))V_{v}\,db,
 \qquad \mathcal{F}f_{v}=u_{v} .
\end{align*}
Fourier inversion in \eqref{fc10.1}, with \eqref{eq9.10}, gives
\begin{align*}
\mathcal{H}_{W_{v}}(\sigma_{v},\mathbf{0},\mathbf{1})=\|\sigma_{v}(h_{v})u_{v}\|_{2}^{2}
 \gg 1\qquad (C(\sigma_{v})^{1/2}\leq R_{v})
\end{align*}
by \eqref{equ9.9} for large $M$. The seminorms are polynomial in $R_{v}$,
also for the Eisenstein lifts by Lemma~\ref{lem:maximal-test-sections}
with $\overline{\psi}_{v}$ and induction in stages.

At finite $v$, put $r_{v}=e_{v}(\mathfrak{a})$. Take $W_{v}^{[0]}$ at $r_{v}=0$ and
$W_{v}^{[1]}$ plus a large positive multiple of $q_{v}^{\vartheta_{\pi_{v}}}W_{v}^{[0]}$
at $r_{v}=1$.
For $r_{v}\geq2$, add a large positive multiple of $q_{v}^{r_{v}(\vartheta_{\pi_{v}}+\delta)}$
times the vector for $r_{v}=1$ to $\sum_{m=2}^{r_{v}}b_{r_{v}-m}W_{v}^{[m]}$.
Write $W_{v}=\sum_{m=0}^{r_{v}}a_{v,m}W_{v}^{[m]}$ for this choice of vector.
By \eqref{eq:depth-projector} and Lemmas~\ref{lem:prime-old-row-expansion}
and~\ref{lem:depth-low}, $\sum_{m} a_{v,m}\alpha_{v}(\sigma_{v};m)\geq0$ and is
$\gg L_{v}(1,\sigma_{v},\Ad)$ if $r_{\sigma_{v}}=r_{v}\geq1$. For small $\delta>0$,
\eqref{eq:depth-inverse} gives
\begin{align*}
 \sum_{m=0}^{r_{v}}|a_{v,m}|q_{v}^{(1+\delta)m}\ll_{\delta} q_{v}^{(1+\delta)r_{v}}.
\end{align*}
The factor from \eqref{eq:different-raw-weight} at $\mathbf{s}=\mathbf{0}$ is $q_{v}^{-2d_{v}}$.

At real $v$, \eqref{7.5} is $O(R_{v}^{O(1)})$, and its dual vanishes
by $w_{2}n_{12}(b)w_{2}=n_{13}(b)$ and $\mathcal{F} f_{v}(0)=0$.
At finite depth $m$, the original integral is \eqref{14.37} at $r_{v}=m$.

At real $v$, divide the product of \eqref{fc5.22} and \eqref{fc5.23}
by the two local $L$-factors in \eqref{equa5.23}. By
\eqref{eq:minimal-w2-kirillov-exponents} it is $O(R_{v}^{O(1)}|z|^{1-O(\delta)})$
at zero and rapidly decaying at infinity. The $w_{2}$ integral is then
$O(R_{v}^{O(1)})$, as are the other two Weyl integrals by \eqref{eq:real-weyl-reductions}.
Equations \eqref{6.15}--\eqref{6.14} give the same bound for \eqref{6.1}.
On the dual side the matrix argument in \eqref{e5.3} and \eqref{6.1} is $(\mathbf{0},g')gw_{2}$.

At finite depth $m$, \eqref{eq:finite-degenerate-estimates} at
$(r_{v},n_{v})=(m,0)$ bounds the minimal integrals, with the characters
included in $q_{v}^{-\nu_{i}}$. The maximal integrals are $O_{F}(q_{v}^{(1+O(\delta))m})$.
The translates in $W_{v}^{[m]}$ contribute $q_{v}^{m}$. In the coordinates
of \eqref{eq:maximal-w2-variables} the original integral has
$e_{v}(tu)\geq m-d_{v}$, with $u=r^{-1}z$.
In the dual integral, the sum over $e_{v}(y)\geq\max(0,m-e_{v}(t))-d_{v}$
is $O_{F}(q_{v}^{O(\delta m)})$ after division by the $y$-Tate factor in
\eqref{6.15}. The $t,r,u$ sums have decay exponents
$\geq1-\vartheta_{\pi_{2,v}}-O(\delta)>0$.

\clearpage
\section{Polynomial Bounds for the Dual Weight}
\label{app:central-polynomial-transform}

With $m$ from \eqref{eq:R-admissible-test} and $A(t)$ from
\eqref{eq:central-input-Ft}, the zeros of $H_{R}$ cancel
the poles of
$\mathcal{S} h_{v,R}(it)=\pi^{2}H_{R}(t)/(A(t)^{4}\cosh(\pi t))$
at $t=\pm i(2j+3/2)$ within $|\Im(t)|<m+1$. On closed substrips,
Stirling's formula and the Gaussian bound for $H_{R}$
give rapid decay. Fourier inversion on $\Im(t)=a$
and Abel inversion \eqref{10.36} give, for $0<a<m+1$,
\begin{equation}\label{eq:real-inverse-Abel-decay}
 \left|h_{v,R}\left(z\diag(e^{r/2},e^{-r/2})\right)\right|
 \ll_{R,a}e^{-(a+1/2)r},\qquad r\geq0,\quad z\in\mathbb{R}^{\times}.
\end{equation}
This bound is integrable against the radial Haar density $O(e^{r})$ for $a>1/2$.

In \eqref{eq:direct-central-input-hypergeometric}, the two terms of
\cite[(15.8.2)]{DLMF} agree after $t\mapsto-t$. For the exponent
$-1/2-it$ we shift to $\Im(t)=-\gamma$, $5/2<\gamma<m+1$.
The zeros of $H_{R}$ cancel the crossed poles, so
$f_{v,R}(b)\ll_{R,\gamma}(1+|b|)^{-1/2-\gamma}$.
Lemma~\ref{lemm11.3} and Theorem~\ref{prop10.2} give convergence of the
principal-series basis sum. The weighted $L^{1}$ bound in Lemma~\ref{lem:principal-series-periods}
makes \eqref{eq:direct-central-cosine-period} holomorphic for $|\Im(u)|<1/2$.

\subsection{The Principal-Series Periods}

For $w>0$, define
\begin{align*}
V_{u,\epsilon}(w):=\int_{0}^{\pi}\int_{\mathbb{R}}
\frac{(\sgn Q_{\alpha,w}(y))^{\epsilon}}{|Q_{\alpha,w}(y)|^{1/2+iu}}\cdot  B_{u,\epsilon}(1,y)\,dy\,\frac{d\alpha}{\sin\alpha},
\end{align*}
with $B_{u,\epsilon}$ from \eqref{eq:central-K-Fourier-Mellin-kernel} and
\begin{align*}
Q_{\alpha,w}(y):=(\cot\alpha\,y-1)
 ((1+w\sin^{2}\alpha)y+w\sin\alpha\cos\alpha). 
\end{align*}

\begin{lemma}\label{lem:cosine-formula-V}
Take $C_{\epsilon}$ from Definition~\ref{def:real-dual-transform}.
For $w>0$, $u\in\mathbb{R}$ and $\epsilon\in\{0,1\}$,
\begin{multline*}
(1+w)^{iu}V_{u,\epsilon}(w)=\frac{8}{C_{\epsilon}(1/2+iu)}
 \int_{0}^{\infty}\cos(u\tau)\sum_{\kappa=\pm1}\kappa^{\epsilon}\\
 \int_{0}^{\infty} K_{\kappa}(\sqrt{w};\cosh v,2\cosh\tau)\,dv\,d\tau,
\end{multline*}
where $K_{\kappa}(z;c,X):=\int_{r_{1}}^{\infty} P_{\kappa}(x,z;c,X)^{-1/2}\,dx$
for $z\ne0$, $c>1$, $X>2$, and $r_{1}$ is the largest root of
\begin{align*}
 P_{\kappa}(x,z;c,X):=4(cx-\kappa)(x^{2}-1-z^{2})
 +(1+z^{2})(X+2\kappa)(x^{2}-1).
\end{align*}
\end{lemma}

\begin{proof}
Insert the Fourier identity
\begin{align*}
\frac{(\sgn Q_{\alpha,w}(y))^{\epsilon}}{|Q_{\alpha,w}(y)|^{1/2+iu}}=\frac{1}{C_{\epsilon}(1/2+iu)}
 \int_{\mathbb{R}}(\sgn\xi)^{\epsilon}|\xi|^{-1/2+iu}
 e^{-2\pi i\xi Q_{\alpha,w}(y)}\,d\xi
\end{align*}
in $V_{u,\epsilon}(w)$. Restrict to $\eta<|\xi|<\eta^{-1}$ and take
$\rho\to0^{+}$ after the $y$ integral with $e^{-\rho y^{2}}$.
On $\xi>0$, with $a$ from \eqref{eq:central-K-Fourier-Mellin-kernel},
the change $(a,\xi)\mapsto(\kappa r,re^{-\tau}/(1+w))$ with
$\kappa=\pm1$ and $r>0$ has $|da\,d\xi|=\xi\,dr\,d\tau$.
In $B_{u,\epsilon}(1,y)$, write
$K_{0}(2\pi r)=\int_{0}^{\infty} e^{-2\pi r\cosh v}\,dv$.
For $0<\alpha<\pi/2$, the Gaussian $y$ integral and the full Laplace
integral on $r>0$ give the integrand against $dv\,d\tau\,d\alpha/\sin\alpha$
\begin{align*}
 \frac{\kappa^{\epsilon}}{C_{\epsilon}(1/2+iu)}
 \frac{(1+w)^{-iu}e^{-iu\tau-i\pi/4}}{(k_{0}+wk_{1})^{1/2}},
\end{align*}
where $c=\cosh v$, $X=2\cosh\tau$ and
\begin{align*}
k_{0}=4c\cot\alpha-i(X-2\kappa),\qquad k_{1}=4c\sin\alpha\cos\alpha-i(X+2\kappa\cos2\alpha).
\end{align*}
For fixed $\kappa$, the two signs of $\xi$ give conjugate square-root
factors. With the principal branch,
\begin{align*}
 \Re\left(\int_{0}^{\pi/2}e^{-i\pi/4}(k_{0}+wk_{1})^{-1/2}
 \frac{d\alpha}{\sin\alpha}\right)=K_{\kappa}(\sqrt{w};c,X).
\end{align*}
With $\cot\alpha\mapsto-ix$, the differential is $-P_{\kappa}^{-1/2}\,dx$.
The signs of $P_{\kappa}$ at $-1,0,1$ are $+,-,-$, so its roots are
real and simple. Rotate the contour to $x:+\infty\to0$ with the branch
positive for $x>r_{1}$.

The truncated Fourier integral is $O_{u}(|Q_{\alpha,w}(y)|^{-1/2})$
uniformly in $\eta$. By \eqref{eq:central-K-Fourier-Mellin-bound},
splitting at the roots of $Q_{\alpha,w}$ gives
\begin{align*}
 \int_{\mathbb{R}}\frac{|B_{u,\epsilon}(1,y)|}{|Q_{\alpha,w}(y)|^{1/2}}\,dy
 \ll_{u,w}
 \begin{cases}
 \sqrt{\alpha}\log(2/\alpha),&\alpha\to0^{+},\\
 \log^{2}(2/|\alpha-\pi/2|),&\alpha\to(\pi/2)^{-}.
 \end{cases}
\end{align*}
Both endpoint bounds are integrable against $d\alpha/\sin\alpha$ on
$(0,\pi/2)$, and reflection gives the other half.
After the Gaussian integration the cutoff restricts $r$ to
$\eta(1+w)e^{\tau}<r<\eta^{-1}(1+w)e^{\tau}$. Integration by parts in $r$
bounds the truncated kernel by
$O(|k_{0}+wk_{1}|^{-1/2})\ll w^{-1/3}|k_{0}|^{-1/6}|k_{1}|^{-1/3}$.
Split the $\alpha$ integral at $\pi/4$ and, for $0<\tau<1$, also at $\tau^{2}/c$:
\begin{align*}
 \int_{0}^{\pi/2}|k_{0}|^{-1/6}|k_{1}|^{-1/3}\frac{d\alpha}{\sin\alpha}
 \ll
 \begin{cases}
 c^{-1/8}X^{-3/8},&\tau\ge1,\\
 c^{-1/6}+c^{-1/3}\tau^{-1/3},&0<\tau<1 .
 \end{cases}
\end{align*}
With $c=\cosh v$ and $X=2\cosh\tau$, these bounds allow dominated
convergence as $\eta\to0^{+}$.
The reflection $(\alpha,y)\mapsto(\pi-\alpha,-y)$ and the sum over
the signs of $\xi$ contribute $4$. Summing $e^{\pm iu\tau}$ on
$\tau>0$ gives $2\cos(u\tau)$.
\end{proof}

\begin{lemma}\label{lem:principal-series-periods}
Let $H_{R}$ satisfy the hypotheses of
Proposition~\ref{prop:direct-polynomial-central-transform}, and let
$\mathscr{T}_{\epsilon} H_{R}$ be as in Definition~\ref{def:real-dual-transform}.
There are functions $g_{\pm}$ on $[0,\infty)$ with $g_{\pm}'(0)=0$
such that, for $u\in\mathbb{R}$ and $\epsilon\in\{0,1\}$,
\begin{equation}\label{eq:direct-central-cosine-period}
 (\mathscr{T}_{\epsilon}H_{R})(u)=\frac{2}{\pi^{3/2}}
 \int_{0}^{\infty}\cos(u\tau)
 [g_{-}(\tau)+(-1)^{\epsilon}g_{+}(\tau)]\,d\tau.
\end{equation}
For $0\leq\vartheta<1/2$ they satisfy
\begin{align*}
 \int_{0}^{\infty}(|g_{-}(\tau)|+|g_{+}(\tau)|)e^{\vartheta\tau}\,d\tau
 \ll_{m,\vartheta}R^{2},\qquad 
 \int_{0}^{\infty}(|g_{-}'''(\tau)|+|g_{+}'''(\tau)|)\,d\tau
 \ll_{m} R^{5}.
\end{align*}
\end{lemma}

\begin{proof}
In \eqref{eq:central-omega-definition}, the substitutions
$\theta\mapsto\arctan x-\theta'$, $x\mapsto\sqrt{1+s^{-1}}\,y$,
$s\mapsto z^{2}$ and $\cot\theta'\mapsto v/\sqrt{1+z^{2}}$, followed by
$(z,v)\mapsto(b\sin\alpha,b\cos\alpha)$, give
\begin{equation}\label{eq:direct-radial-central-kernel}
 \Omega_{u,\epsilon}(t)=2\int_{0}^{\infty}
 F_{t}(b)(1+b^{2})^{iu}V_{u,\epsilon}(b^{2})\,db.
\end{equation}
With $A(t),F_{t}$ from \eqref{eq:central-input-Ft} and $K_{\pm1}$
from Lemma~\ref{lem:cosine-formula-V}, set
\begin{align*}
 m_{\pm}(t,X)&:=\frac{1}{A(t)}\int_{\mathbb{R}}F_{t}(z)
 \int_{0}^{\infty}K_{\pm1}(z;\cosh v,X)\,dv\,dz,\\
 g_{\pm}(\tau)&:=\int_{\mathbb{R}}H_{R}(t)t\tanh(\pi t)
 m_{\pm}(t,2\cosh\tau)\,dt.
\end{align*}

Insert Lemma~\ref{lem:cosine-formula-V} at $w=b^{2}$ into
\eqref{eq:direct-radial-central-kernel} and extend the even $b$ integral
to $\mathbb{R}$. In the cutoff estimate of Lemma~\ref{lem:cosine-formula-V},
$w^{-1/3}=b^{-2/3}$. The bounds after \eqref{eq:direct-central-input-wave-bound}
make $b^{-2/3}|F_{t}(b)|$ integrable in $b$, and $H_{R}$ supplies Gaussian
decay in $t$. Thus
\begin{align*}
 \frac{\Omega_{u,\epsilon}(t)}{A(t)}
 =\frac{8}{C_{\epsilon}(1/2+iu)}\int_{0}^{\infty}\cos(u\tau)
 [(-1)^{\epsilon} m_{-}(t,2\cosh\tau)+m_{+}(t,2\cosh\tau)]\,d\tau.
\end{align*}
In \eqref{eq:direct-threefold-central-transform}, the $t$ integral now
gives $(-1)^{\epsilon} g_{-}+g_{+}$. By \eqref{eq:principal-unitary-Fourier-factor},
\begin{align*}
\frac{\overline{C_{\epsilon}(1/2-iu)}}{4\pi^{3/2}}\cdot 
\frac{8}{C_{\epsilon}(1/2+iu)}=\frac{2(-1)^{\epsilon}}{\pi^{3/2}},
\end{align*}
which gives \eqref{eq:direct-central-cosine-period}.

\smallskip
\noindent\textit{Evaluation of $m_{\pm}$.} For $m_{-}$, write $\delta=X-2$
and $\theta=\delta\partial_{\delta}$. With $P=P_{-1}(x,z;c,X)$ from
Lemma~\ref{lem:cosine-formula-V}, put
$\mathcal{L}_{z}=\partial_{z}((1+z^{2})\partial_{z})$.
We use the identity
\begin{multline}\label{eq:central-period-certificate}
 \left[\theta^{3}-\frac{\theta}{4}+
 \frac{\delta}{4}(\theta+1/2)^{3}\right]P^{-1/2}
 -\mathcal{L}_{z}(\theta P^{-1/2})\\
 =\frac{\delta}{16}\left\{
 \partial_{x}((x^{2}-1)T_{0}P^{-5/2})+
 \sqrt{c^{2}-1}\partial_{c}(\sqrt{c^{2}-1}U_{0}P^{-5/2})\right\}.
\end{multline}
With $Z=1+z^{2}$, the polynomials on the right are
\begin{multline*}
 \frac{5}{8}T_{0}
 =2xc^{2}(x^{2}-Z)(7Z-5-5x^{2})\\
 +c\bigl[(30-70Z+5\delta Z)x^{4}
 +(20+2Z+(10-7\delta)Z^{2})x^{2}+Z(7Z-5)(\delta+4)\bigr]\\
 +x\bigl[(70-112Z)x^{2}+52Z^{2}-10Z+5\delta Z(Z-1)(x^{2}-1)\bigr],
\end{multline*}
and 
\begin{multline*}
 \frac{5}{24}U_{0}
 =2c\bigl[Z(5-7Z)+(-15+21Z-4Z^{2})x^{2}\bigr]\\
 +x\bigl[(14Z-10)x^{2}-22Z^{2}+38Z-20+\delta Z(4Z-5)(1-x^{2})\bigr].
\end{multline*}

Multiplication by $P^{7/2}$ verifies \eqref{eq:central-period-certificate}
as a polynomial identity. The bidegrees in $(x,c)$ are $(5,2)$ and $(3,1)$.
In \eqref{eq:central-period-certificate}, $P^{-1/2}$ and $P^{-5/2}$
are the values at $\lambda=-1/2,-5/2$ of the distribution
$\mathbf{1}_{P>0}P^{\lambda}$, extended meromorphically from $\Re(\lambda)>0$.
The right side of \eqref{eq:central-period-certificate} integrates to
zero over $x,c>1$ with measure $dx\,dc/\sqrt{c^{2}-1}$. The boundary terms
vanish at $x=1$ by $P<0$, at $c=1$ by $\sqrt{c^{2}-1}$, and at infinity
by the degree bounds.

For the roots $r_{3}<r_{2}<0<1<r_{1}$ of $P_{-1}$,
\cite[(19.16.1), (19.25.1)]{DLMF} give
\begin{align*}
 K_{-1}(z;c,2+\delta)
 =\frac{K\big(\sqrt{1-(r_{1}-r_{2})/(r_{1}-r_{3})}\big)}
 {\sqrt{c(r_{1}-r_{3})}},
\end{align*}
where $K(k):=\int_{0}^{\pi/2}(1-k^{2}\sin^{2}\varphi)^{-1/2}\,d\varphi$
is the complete elliptic integral of the first kind.
For $|z|\ge1$, the gaps $r_{1}-r_{2},r_{1}-r_{3}$ are $\asymp c,z^{2}/c$ for
$c\le|z|$ and $\asymp |z|,|z|$ otherwise. By \cite[(19.12.1)]{DLMF},
\begin{align*}
\int_{0}^{\infty} K_{-1}(z;\cosh v,2+\delta)\,dv\ll |z|^{-1}\log^{2}(2+|z|). 
\end{align*}
On disks of radius $\asymp |z|^{2}$ in $z^{2}$ and neighborhoods of compact
positive $\delta$ intervals, Cauchy's estimates give the same bound for
$(z\partial_{z})^{j}\theta^{l}\int_{0}^{\infty}K_{-1}(z;\cosh v,2+\delta)\,dv$. Integration by parts in $z$ and
$\mathcal{L}_{z}F_{t}=-(t^{2}+1/4)F_{t}$ from
\eqref{eq:direct-central-input-wave-equation} with $z\mapsto\sinh r$ give
\begin{equation}\label{eq:central-minus-period-ode}
 [\theta(\theta^{2}+t^{2})+\delta(\theta+1/2)^{3}/4]
 m_{-}(t,2+\delta)=0.
\end{equation}

For $t\ne0$, put $a=1/4+it/2$, $b=1/4-it/2$ and
\begin{align*}
f_{t}(\delta):={}_{2}F_{1}(b,b;1-it;-\delta/4). 
\end{align*}
Gauss's equation gives the solutions $\delta^{\mp it/2}f_{\pm t}$ of
\begin{align*}
[\theta^{2}+t^{2}/4+\delta(\theta+1/4)^{2}/4]y=0. 
\end{align*}
The product rule gives \eqref{eq:central-minus-period-ode} for their three
products $\delta^{-it}f_{t}^{2},\delta^{it}f_{-t}^{2},f_{t}f_{-t}$.
Their exponents $-it,it,0$ at $\delta=0$ are distinct, so they form a basis.

For $F_{t}$ in \eqref{eq:central-input-Ft}, the substitution $z\mapsto\sqrt{w}$
in the Mellin integral gives, by \cite[(15.14.1), $(a,b;1/2)$]{DLMF},
\begin{equation}\label{eq:central-input-mellin}
 2\int_{0}^{\infty} z^{2r-1}F_{t}(z)\,dz
=\frac{\pi^{1/2}\Gamma(r)\Gamma(a-r)\Gamma(b-r)}
 {A(t)\Gamma(1/2-r)},\qquad 0<\Re(r)<1/4.
\end{equation}

In \eqref{eq:central-period-certificate}, $c\mapsto-c$ and $\delta\mapsto-p$
with $p=X+2$ give
\begin{align*}
[\theta(\theta^{2}+t^{2})-p(\theta+1/2)^{3}/4]m_{+}(t,p-2)=0,
\end{align*}
where $\theta=p\partial_{p}$. At $c=1$ the root $x=1$ is simple,
with derivative $4+2(1+z^{2})(p-2)$, so the $c=1$ boundary term still vanishes.
As $x\to1^{+}$, rescale the $c$ interval where $P_{+1}>0$ to $[0,1]$.
The constant and linear terms of the smooth factor give
$B(1/2,-3/2)=B(3/2,-3/2)=0$, and the remainder gives an
$O(x-1)$ boundary term after multiplication by $x^{2}-1$.
For $X=2\cosh\tau$, the substitution $p\mapsto4\cosh^{2}r$ with $r=\tau/2$ gives
\begin{align*}
 E_{t}(r):=(\cosh r)^{it}{}_{2}F_{1}(a,a;1/2;-\sinh^{2}r)
 =(p/4)^{it/2}{}_{2}F_{1}(a,a;1/2;1-p/4).
\end{align*}
Gauss's equation gives
\begin{align*}
[\theta^{2}+t^{2}/4-p(\theta+1/4)^{2}/4]E_{t}(r)=0. 
\end{align*}
The product rule gives the $m_{+}$ equation for $E_{t}(r)^{2}$.

For $\kappa\in[-1,1]$, define $K_{\kappa}(z;c,X)$ by the integral in
Lemma~\ref{lem:cosine-formula-V}, extended to $X>0$ for $\kappa=0$.
For $X>0$, put
\begin{align*}
 m_{0}(t,X):=\frac{1}{A(t)}\int_{\mathbb{R}}F_{t}(z)
 \int_{0}^{\infty}\left.K_{\kappa}(z;\cosh v,X)\right|_{\kappa=0}\,dv\,dz.
\end{align*}
For $X\ge10$, the shift $x\mapsto r_{1}+y$ factors $P_{\kappa}$ as
$4cy(y+r_{1}-r_{2})(y+r_{1}-r_{3})$. The logarithmic $\kappa$-derivatives of
$2r_{1}-r_{2}-r_{3}$ and $(r_{1}-r_{2})(r_{1}-r_{3})$ at fixed $z,c,X$
are $O(X^{-1})$, so $|\partial_{\kappa}K_{\kappa}|\ll K_{\kappa}/X$.
With $|F_{t}(z)/A(t)|\ll_{t}\log(2+|z|)/(1+z^{2})^{1/4}$ from the bound
after \eqref{eq:direct-central-input-wave-bound} and the elliptic bounds
\cite[(19.16.1), (19.25.1), (19.12.1)]{DLMF},
\begin{align*}
K_{\kappa}(z;\cosh v,X)\ll
\begin{cases}
 \dfrac{\log(2+X(1+z^{2}))}{\sqrt{X(1+z^{2})}},&
 \cosh v\le X(1+z^{2}),\\
 (1+z^{2})^{-1/4}e^{-v/2},&\cosh v>X(1+z^{2}).
\end{cases}
\end{align*}
Therefore, 
\begin{equation}\label{B.7}
\int_{0}^{\infty}K_{\kappa}(z;\cosh v,X)\,dv
 \ll\frac{\log^{2}(2+X(1+z^{2}))}{\sqrt{X(1+z^{2})}}
 \ll\frac{\log^{2}X\,\log^{2}(2+|z|)}{\sqrt{X(1+z^{2})}}. 
\end{equation}

By $|\partial_{\kappa}K_{\kappa}|\ll K_{\kappa}/X$ and \eqref{B.7},
\begin{align*}
|m_{\pm}(t,X)-m_{0}(t,X)|
& =\left|\int_{\mathbb{R}}\frac{F_{t}(z)}{A(t)}\int_{0}^{\infty}
 \int_{0}^{\pm1}\partial_{\kappa}K_{\kappa}(z;\cosh v,X)\,d\kappa\,dv\,dz\right|\\
&\ll X^{-1}\int_{\mathbb{R}}\left|\frac{F_{t}(z)}{A(t)}\right|
 \sup_{|\kappa|\leq1}\int_{0}^{\infty}K_{\kappa}(z;\cosh v,X)\,dv\,dz\\
&\ll_{t}\frac{\log^{2}X}{X^{3/2}}\int_{\mathbb{R}}
 \frac{\log^{3}(2+|z|)}{(1+z^{2})^{3/4}}\,dz
 \ll_{t}\frac{\log^{2}X}{X^{3/2}}.
\end{align*}

For $0<\Re(q)<1/8$, we have 
\begin{equation}\label{eq:central-m0-mellin-expansion}
\int_{0}^{\infty}\frac{m_{0}(t,X)}{X^{1/2+q}}\,dX=\frac{1}{\sqrt{\pi}\,q^{3}}
 +\frac{\ell(t)}{\sqrt{\pi}\,q^{2}}
 +\frac{\ell(t)^{2}}{2\sqrt{\pi}\,q}
 +\text{holomorphic terms at }q=0,
\end{equation}
where $\ell(t)=-\log4-2\gamma-\psi(a)-\psi(b)$, $\psi=\Gamma'/\Gamma$
and $\gamma$ is Euler's constant.
With $X\mapsto4cx|x^{2}-1-z^{2}|Y/((1+z^{2})(x^{2}-1))$, the
$Y$ integral gives $B(q,1/2)$ on $Y>1$ for $1<x<\sqrt{1+z^{2}}$ and
$B(q,1/2-q)$ on $Y>0$ for $x>\sqrt{1+z^{2}}$
\cite[(5.12.1), (5.12.3)]{DLMF}. Their $O(q)$ difference contributes
a holomorphic term after the $v$ integral $\tfrac12 B(q/2,1/2)$. After $x\mapsto\cosh y$, the integral in
\eqref{eq:central-m0-mellin-expansion} equals, modulo a term holomorphic at $q=0$,
\begin{multline*}
 \frac{4^{-q}B(q,1/2-q)B(q/2,1/2)}{2A(t)}\int_{\mathbb{R}}F_{t}(z)(1+z^{2})^{q-1/2}\\
 \left(\tfrac12 B(q/2,1/2)-q(\operatorname{arsinh}z)^{2}\right)\,dz .
\end{multline*}
The term $-q(\operatorname{arsinh}z)^{2}$ comes from translating
$\log|\sinh y|$ by $\pm\operatorname{arsinh}z$ in the $x$ integral.

With $z\mapsto\sqrt{s/(1-s)}$ on $z>0$,
\cite[(15.8.1), (15.6.8), (15.4.20)]{DLMF} give
\begin{align*}
 \int_{\mathbb{R}}F_{t}(z)(1+z^{2})^{q-1/2}\,dz
 =\frac{\sqrt{\pi}\,\Gamma(a-q)\Gamma(b-q)}{\Gamma(1/2-q)^{2}}.
\end{align*}
At $q=0$, $z\mapsto\sinh u$ gives $\int_{\mathbb{R}}F_{t}(\sinh u)\,du=A(t)/\sqrt{\pi}$.
By \cite[(14.12.7), $m=0$, $\nu=-2b$ at $\pm i\sinh u$, (5.12.7), (5.12.1)]{DLMF},
\begin{align*}
 \left.\frac{\sqrt{\pi}}{A(t)}\partial_{\sigma}^{2}
 \int_{\mathbb{R}}e^{\sigma u}F_{t}(\sinh u)\,du\right|_{\sigma=0}
 =\frac{\psi'(a)+\psi'(b)}{2}-\frac{\pi^{2}}{4}.
\end{align*}
Multiply this value by $-1/\sqrt{\pi}$ for its contribution to the
$q^{-1}$ coefficient in \eqref{eq:central-m0-mellin-expansion}.

Both period equations are regular singular at infinity with exponent
$-1/2$ of multiplicity three. Since $m_{\pm}-m_{0}=O_{t}(X^{-3/2}\log^{2}X)$,
the principal part in \eqref{eq:central-m0-mellin-expansion} determines
the quadratic in $\log X$ in the leading term:
\begin{align*}
m_{\pm}(t,X)=\frac{X^{-1/2}(\log X+\ell(t))^{2}}{2\sqrt{\pi}}
 +o_{t}(X^{-1/2}),\qquad X\to\infty.
\end{align*}
For $m_{-}$, \cite[(15.8.8), $m=0$]{DLMF} at $(b,b;1-it)$ and $(a,a;1+it)$
with argument $-(X-2)/4$ gives the same expansion for the right side of
\eqref{eq:direct-central-minus-period}. Thus
\begin{equation}\label{eq:direct-central-minus-period}
 m_{-}(t,2+\delta)
 =D(t)\delta^{-it}f_{t}(\delta)^{2}
 +D(-t)\delta^{it}f_{-t}(\delta)^{2}
 +\frac{\sqrt{\pi}}{4t\tanh(\pi t)}f_{t}(\delta)f_{-t}(\delta),
\end{equation}
where
$D(t):=\sin^{2}(\pi a)4^{2a-1}\Gamma(b)^{2}\Gamma(it)^{2}/(2\sqrt{\pi}\,\Gamma(a)^{2})$.
At $t=0$, \eqref{eq:direct-central-minus-period} is read by continuity.

As $X\to\infty$, $E_{t}(\tau/2)$ has the expansion
$\sqrt{2\pi}\,X^{-1/4}(\log X+\ell(t))/A(t)$ with error
$O_{t}(X^{-5/4}\log X)$ by \cite[(15.8.8)]{DLMF} at $(a,a;1/2)$ and $-(X-2)/4$.
Squaring and matching the expansion of $m_{+}$ gives
\begin{equation}\label{eq:direct-central-plus-period}
 m_{+}(t,2\cosh\tau)
 =\frac{A(t)^{2}}{4\pi^{3/2}}E_{t}(\tau/2)^{2}.
\end{equation}

\noindent\textit{Bounds for $g_{-}$.} In \eqref{eq:direct-central-minus-period},
$\delta=2\cosh\tau-2$ and $|t\tanh(\pi t)D(t)|=\sqrt{\pi}/8$.
Euler's integral \cite[(15.6.1)]{DLMF} bounds the derivatives of $f_{t}$
for $j\leq3$. For $0<\tau\leq1$, the $D(\pm t)$ terms in
\eqref{eq:direct-central-minus-period} are integrated on $\Im(t)=\pm9/4$.
\begin{align*}
 &|(\delta\partial_{\delta})^{j}f_{t}(\delta)|
 \ll_{j}(1+|t|)^{j+1/2}(1+\tau)e^{-\tau/4},
 \quad \tau>1,\\
 &\left|\partial_{\tau}^{j}\int_{\mathbb{R}\pm9i/4}
 H_{R}(t)t\tanh(\pi t)D(\pm t)\delta^{\mp it}f_{\pm t}(\delta)^{2}\,dt\right|
 \ll_{m}R^{j+2}\tau^{9/2-j},\\
 &|\partial_{\tau}^{j}(f_{t}(\delta)f_{-t}(\delta))|
 \ll_{j}(1+|t|)^{j+1}, \quad 0<\tau\leq1.
\end{align*}
On the shifted lines, $|\delta^{\mp it}|=\delta^{9/4}$.
The zeros of $H_{R}$ leave only even residues at $\pm i,\pm2i$, with
$\partial_{\tau}^{j}$ bounds $O_{m}(R^{-m})$.
Thus $g_{-}'(0)=0$ and $|g_{-}^{(j)}(\tau)|\ll_{m}R^{j+2}(1+\tau)^{2}e^{-\tau/2}$.

\noindent\textit{Bounds for $g_{+}$.} For $h_{t}=\sqrt{\cosh r}E_{t}(r)$,
$h_{t}''=(t^{2}+1/4)\operatorname{sech}^{2}r\,h_{t}$ with data $(1,0)$.
With $r_{0}=\log(2+|t|)$, comparison on $[0,r_{0}]$, the integral equation
from $r_{0}$, and Stirling for $A(t)$ give
\begin{align*}
 &h_{t}(r)\leq\sqrt{\cosh r}\exp(\sqrt{t^{2}+1/4}\arctan(\sinh r)),
 \quad 0\leq r\leq r_{0},\\
 &(t^{2}+1/4)\operatorname{sech}^{2}(r_{0}+s)\ll e^{-2s},
 \quad s\geq0,\\
 &|A(t)E_{t}^{(j)}(r)|\ll_{j}(1+|t|)^{j-1/2}
 \begin{cases}
  1,&0\leq r\leq r_{0},\\
  (1+r-r_{0})e^{-(r-r_{0})/2},&r>r_{0}.
 \end{cases}
\end{align*}
From \eqref{eq:direct-central-plus-period} at $r=\tau/2$,
$g_{+}'(0)=0$, $\|g_{+}'''\|_{L^{1}}\ll_{m}R^{4}\log(2+R)$ and
$\|e^{\vartheta\tau}g_{+}\|_{L^{1}}\ll_{m,\vartheta}R^{2}$ for $0\leq\vartheta<1/2$.
\end{proof}

\subsection{The Discrete-Series Coefficient}\label{app:discrete-series-coefficient}

For $0<\eta<1/2$, the vectors $W_{n}$ in \eqref{10.8}
and their Hecke integrals $\Psi(W_{n};\lambda)$ from
\textsection\ref{def:real-hecke-integral} satisfy
\begin{equation}\label{eq:uniform-kirillov-periods}
 |\Psi(W_{n};0)|\ll1,\qquad
 \int_{\mathbb{R}^{\times}}|W_{n}(\diag(a,1))|^{2}
 |a|^{-1+2\eta}\,d^{\times}a\ll_{\eta}1,
\end{equation}
for tempered principal and discrete series. The Hecke integral of
\cite[(2.16)]{BM05} gives the bound for $\Psi(W_{n};0)$.
For complementary series with $|\eta_{v}|\leq\theta_{0}<\eta$,
the normalization in \cite[(4.16)]{BM05} gives
$|\Psi(W_{n};0)|\ll_{\theta_{0}}(1+|n|)^{\theta_{0}}$, and the integral satisfies
\begin{align*}
 \int_{\mathbb{R}^{\times}}|W_{n}(\diag(a,1))|^{2}
 |a|^{-1+2\eta}\,d^{\times}a\ll_{\eta,\theta_{0}}(1+|n|)^{2\theta_{0}}.
\end{align*}
For principal and complementary series, put
$v_{n,\nu}(x)=(1+x^{2})^{-1/2-\nu}e^{2in\arctan x}$ with
$0\leq\Re(\nu)\leq\theta_{0}$. The Fourier formula
\cite[(2.16), (4.16)]{BM05} and Stirling's formula give
\begin{align*}
 |W_{n}(\diag(a,1))|\ll_{\theta_{0}}
 (1+|n|)^{\Re(\nu)}|a|^{1/2-\Re(\nu)}|\widehat v_{n,\nu}(a)|.
\end{align*}
For $q>1/(\eta-\theta_{0})$, Hausdorff--Young bounds
$\|\widehat v_{n,\nu}\|_{q}$ uniformly since
$|v_{n,\nu}(x)|\leq(1+x^{2})^{-1/2}$.
On $|a|\leq1$, H\"older's inequality with
$|a|^{-1+2(\eta-\theta_{0})}\in L^{q/(q-2)}(-1,1)$ gives the claimed bounds,
with $\theta_{0}=0$ for the tempered principal series.
For discrete series of lowest weight $2l$, the Laguerre form of
\eqref{10.8} at $|n|=l+j$ and orthogonality with
$L_{j}^{(2l-1)}=\sum_{k=0}^{j}L_{k}^{(2l-2)}$
\cite[(13.18.17), (18.9.13), and Table~18.3.1]{DLMF} give
$\int_{\mathbb{R}^{\times}}|W_{n}(\diag(a,1))|^{2}|a|^{-1}d^{\times}a
=4\pi/(2l-1)$, which bounds the weighted integral on $|a|\leq1$.
On $|a|\geq1$, the $L^{2}(d^{\times}a)$ norm bounds the integral in all three cases.

\begin{lemma}\label{lem:discrete-series-coefficient}
Let $H_{R}$ satisfy the hypotheses of Proposition~\ref{prop:direct-polynomial-central-transform}.
For lowest weight $2l$, the transform in Definition~\ref{def:real-dual-transform} is
\begin{equation}\label{eq:direct-central-discrete-period}
 (\mathscr{T} H_{R})^{\mathrm{disc}}(2l)
 =\frac{(-1)^{l+1}}{\pi^{3/2}}
 \int_{-1}^{1}P_{l-1}(y)\Phi_{R}(y)\,dy,
\end{equation}
where $P_{l-1}$ is the Legendre polynomial and, for $-1<y<1$,
\begin{multline}\label{eq:central-discrete-coefficient-definition}
\Phi_{R}(y):=-\frac{\sqrt{\pi}}{1+y}\int_{0}^{\pi/2}
 W_{R}(D\sin\theta)\,d\theta\\
+\frac{1}{\sqrt{\pi}}\int_{0}^{\infty}
 \int_{z^{-1}}^{\infty}\frac{W_{R}(z)dudz}{2u^{2}+(1+y)zu+1-y}.
\end{multline}
Here $D=2\sqrt{2(1-y)}/(1+y)$ and
\begin{align*}
 W_{R}(z):=\int_{\mathbb{R}}H_{R}(t)t\tanh(\pi t)F_{t}(z)\,dt
\end{align*}
with $F_{t}$ from \eqref{eq:central-input-Ft}.
\begin{align*}
 |(\mathscr{T} H_{R})^{\mathrm{disc}}(2l)|
 \ll_{m}\begin{cases}
 1,&l=1,\\
 \min\{1,R^{4}[l(l-1)]^{-2}\},&l\geq2.
 \end{cases}
\end{align*}
\end{lemma}

\begin{proof}
In \eqref{10.7} at $\mathbf{s}=\mathbf{0}$, the term
$n=\sgn(a)(l+j)$ has pairing factor $j!/\Gamma(j+2l)$ by \eqref{10.8}.
For $x=4\pi|a|$ and $|w|<1$, \cite[(18.18.27), index $2l-1$, and
(13.6.9), $(l;2l;2xw/(1-w^{2}))$]{DLMF} give
\begin{multline*}
 Q_{l}(x,w):=x^{l}e^{-x/2}w^{l}
 \sum_{j\geq0}\frac{j!L_{j}^{(2l-1)}(x)w^{j}}{\Gamma(j+2l)}
 \int_{0}^{\infty} u^{l-1}e^{-u/2}L_{j}^{(2l-1)}(u)\,du\\
 =\sqrt{\frac{2\pi xw}{1-w^{2}}}\,I_{l-1/2}\left(\frac{xw}{1-w^{2}}\right)
 \exp\left(-\frac{x(1+w^{2})}{2(1-w^{2})}\right).
\end{multline*}
In \eqref{10.7}, inserting $\rho^{l+j}$ and the Fourier integral for
$J_{\pm2(l+j)}$ gives on each $a$ half-line
\begin{align*}
 \frac{1}{2}\int_{0}^{2\pi}e^{2\pi izb_{1}\sin\theta}
 Q_{l}(4\pi|a|,\rho e^{-2i\theta})\,d\theta,
 \qquad 0<\rho<1.
\end{align*}
For $0<\epsilon<1$, the half-integer Bessel bound gives
\begin{align*}
 |Q_{l}(x,w)|\ll_{l}\min\{1,|xw/(1-w^{2})|^{l}\},\qquad 
\int_{0}^{2\pi}|Q_{l}(x,\rho e^{-2i\theta})|\,d\theta
 \ll_{l,\epsilon}\min\{1,x^{\epsilon}\}.
\end{align*}
In \eqref{eq:central-Whittaker-two-K-kernel}, change
$(x,z)\mapsto(2x/\sqrt{1+1/s},z/(2\pi\sqrt{1+s}))$.
After the angular bound for $Q_{l}$, the absolute $x,z$ integral is
$O_{l,\epsilon}((1+1/s)^{-\epsilon/2}(1+s)^{-1/2})$ at fixed $s$,
which is integrable against $ds/s$.

Integrating the Fourier formula for $J_{2(l+j)}$ twice in $\theta$ gives
\begin{align*}
|J_{2(l+j)}(2\pi zb_{1})|\ll_{l}(1+j)^{-2}(|zb_{1}|+|zb_{1}|^{2}).
\end{align*}
By \eqref{eq:uniform-kirillov-periods} and the second moment of $|f_{v,R}|$,
the $j$th term of \eqref{10.7} after integration in $a,z,b_{1}$ is
$O_{l,R}((1+j)^{-2})$. Thus $\rho\to1^{-}$ removes $\rho^{l+j}$.

In the folded $K_{0}$ integral, the angular factor is
$(-1)^{l+1}P_{l-1}(y)$ by \cite[(10.54.2), $n=l-1$]{DLMF}.
In \eqref{eq:central-input-kernel-pairing}, take $(c,d)=(4\sqrt{1+s},4|\sin\theta|)$.
Including both $a$ half-lines, the coefficient in
\eqref{eq:direct-central-discrete-period} is
$(2\cdot\tfrac12)(8\cdot4\cdot4)/(128\pi^{3/2})=\pi^{-3/2}$.

After \eqref{eq:central-Whittaker-two-K-kernel} and
\eqref{eq:central-input-kernel-pairing}, the folded $s,\theta$ integral
in \eqref{10.7} has measure $ds\,d\theta/s$ on $s>0$, $0<\theta<\pi/2$
and integrand
\begin{align*}
\frac{F_{t}\bigl(((1+s)/\sin^{2}\theta-1)^{1/2}\bigr)}
 {2\sin\theta\sin2\theta}
 \int_{0}^{\infty}
 K_{0}\bigl(x(1+s^{-1})^{1/2}/2\bigr)
 \sin\left(\frac{x(y+\cos2\theta)}{2\sin2\theta}\right)\,dx .
\end{align*}
The changes $s\mapsto w/(1+(1+w)u^{2})$ and
$\theta\mapsto\arcsin\sqrt{(1+u^{2})/(1+(1+w)u^{2})}$ give
$F_{t}(\sqrt w)\,dw/\sqrt w$ times
\begin{align*}
 \int_{0}^{\infty}\left(
 \int_{1/u}^{\infty}\frac{dz}{(z+D_{u})\sqrt{1+w+z^{2}}}
 -\frac12\int_{0}^{D_{u}^{2}}\frac{dh}{(1+w+h)\sqrt{D_{u}^{2}-h}}
 \right)\frac{du}{(1+y)u},
\end{align*}
where $w,u>0$ and $D_{u}=(2u^{2}+1-y)/((1+y)u)$.
In the positive term, $z\geq1/u$ becomes $u\geq1/z$, with
$(1+y)u(z+D_{u})=2u^{2}+(1+y)zu+1-y$.
In the negative term, the two branches of $v=D_{u}$ give
$2\,dv/\sqrt{v^{2}-D^{2}}$ on $v>D$. With $h\mapsto v^{2}\sin^{2}\theta$,
the negative term becomes
\begin{multline*}
 -\frac{2}{1+y}\int_{0}^{\pi/2}\int_{D}^{\infty}
 \frac{v\sin\theta\,dv\,d\theta}
 {(1+w+v^{2}\sin^{2}\theta)\sqrt{v^{2}-D^{2}}}\\
 =-\frac{\pi}{1+y}\int_{0}^{\pi/2}
 \frac{d\theta}{\sqrt{1+w+D^{2}\sin^{2}\theta}}.
\end{multline*}
After $w\mapsto b^{2}$ with $dw/\sqrt{w}=2\,db$ on $b>0$, apply
\eqref{eq:central-hypergeometric-kernel-pairing} with
$c=\sqrt{1+z^{2}}$ for the positive term and
$c=\sqrt{1+D^{2}\sin^{2}\theta}$ for the negative term in
\eqref{eq:central-discrete-coefficient-definition}.

The $e^{-itr}$ term of $g_{t}$ in \eqref{eq:direct-central-input-wave-equation} is
$c_{-}(t)e^{-itr}m_{t}(r)$ with $m_{t}(\infty)=1$ and
\begin{equation}\label{eq:wave-connection-coefficient}
 c_{-}(t)=2^{it}\sqrt{\pi}\cdot \frac{\Gamma(-it)}{\Gamma(1/4-it/2)^{2}}
\end{equation}
by \cite[(15.8.2)]{DLMF}. In $W_{R}(\sinh r)$, shift the $e^{-itr}$ term
to $\Im(t)=-\gamma$, $0<\gamma<m+1$. The product
$H_{R}(t)t\tanh(\pi t)c_{-}(t)$ has no poles between the lines.
On $\Im(t)=-\gamma$ away from poles, Stirling's formula and Cauchy's
estimates for the Jost expansion \cite[\textsection4.1, (4.5)]{DU11} at
$(\zeta,V)=(-t,\operatorname{sech}^{2}r/4)$ give
\begin{align*}
|\partial_{t}^{n}c_{-}(t)|\ll_{n,\gamma}(1+|\Re(t)|)^{-n},\qquad |\partial_{r}^{j}\partial_{t}^{n}m_{t}(r)|
 \ll_{j,n,\gamma}(1+|\Re(t)|)^{-n}.
\end{align*}
By $t\mapsto-t$, \eqref{eq:R-admissible-test} and integration by parts in $t$,
\begin{equation}\label{eq:central-input-wave-localization}
 |\partial_{r}^{j}W_{R}(\sinh r)|
 \ll R^{j+2}(\cosh r)^{-1/2}e^{-\gamma r}(1+Rr)^{-n},
\end{equation}
Integrating \eqref{eq:central-input-wave-localization} at $z\mapsto\sinh r$
with $\gamma>n+1/2$ and decay order greater than $n+1$ gives
\begin{equation}\label{eq:central-input-wave-moments}
 \int_{0}^{\infty} z^{n}|W_{R}(z)|\,dz\ll_{n} R^{1-n}.
\end{equation}
For the even moment of degree $2n$, take $m\ge2n$ in
\eqref{eq:R-admissible-test}. Then
\begin{equation}\label{eq:central-input-even-moments}
 \int_{0}^{\infty} z^{2n}W_{R}(z)\,dz=0.
\end{equation}
The moment is the $t$ average of half of \eqref{eq:central-input-mellin}
at $r=n+1/2$. The crossed half-integer poles cancel against $H_{R}$, and
$1/\Gamma(-n)=0$.

In the second term of \eqref{eq:central-discrete-coefficient-definition},
$2u^{2}+(1+y)zu+1-y\geq2u^{2}+2$ for $u\geq1/z$.
Thus the $j$th $y$ derivative of the $u$ integral is $O_{j}(z+z^{j+1})$,
and \eqref{eq:central-input-wave-moments} gives the $C^{4}[-1,1]$ bound.

In the first term, $D\sin\theta\mapsto z$ gives
\begin{align*}
 -\frac{\sqrt{\pi}}{(1+y)D}\int_{0}^{D}
 W_{R}(z)(1-z^{2}/D^{2})^{-1/2}\,dz .
\end{align*}
For $RD\geq1$, \eqref{eq:central-input-even-moments} converts the Taylor
terms on $z<D/2$ to integrals over $z\geq D/2$.
By \eqref{eq:central-input-wave-moments} and
\eqref{eq:central-input-wave-localization},
the displayed $z$ integral is $O_{n,j}(R(RD)^{-n})$ also after $(D\partial_{D})^{j}$:
$(1+z\partial_{z})^{j}W_{R}$ has the same vanishing moments.
The even Taylor expansion of $W_{R}(D\sin\theta)$ for $RD\leq1$
completes the estimate for $0\leq j\leq4$:
\begin{align*}
 \left|\partial_{y}^{j}\left(
 \frac{\sqrt{\pi}}{1+y}\int_{0}^{\pi/2}W_{R}(D\sin\theta)\,d\theta
 \right)\right|
 \ll_{n,j}
 \begin{cases}
 R^{2+2j}(1+R\sqrt{1-y})^{-n},&0\leq y\leq1,\\
 R^{1-n}(1+y)^{n-j},&-1<y\leq0.
 \end{cases}
\end{align*}
Thus $\Phi_{R}\in C^{4}[-1,1]$. For
$\mathcal{L}_{y}=-\partial_{y}((1-y^{2})\partial_{y})$, the coefficients of
$\partial_{y}^{4},\partial_{y}^{3}$ in $\mathcal{L}_{y}^{2}$ are
$(1-y^{2})^{2},-8y(1-y^{2})$. The change $y\mapsto1-u^{2}/R^{2}$ gives
$\|\Phi_{R}\|_{1}\ll1$, $\|\mathcal{L}_{y}^{2}\Phi_{R}\|_{1}\ll R^{4}$.
By $|P_{l-1}|\leq1$ and $\mathcal{L}_{y}P_{l-1}=l(l-1)P_{l-1}$,
\eqref{eq:direct-central-discrete-period} is $O(1)$ and, for $l\geq2$,
$O(R^{4}[l(l-1)]^{-2})$ by integration by parts.
\end{proof}

\subsection{Bounds for the Residual Terms}\label{app:residual-estimates}

\begin{lemma}\label{lem:fourth-moment-residue-bounds}
For the data of \textsection\ref{sec15}, the original and dual residue
contributions to \eqref{15.11} are bounded, respectively, by
\begin{align*}
N_{F}(\mathfrak{q})^{\varepsilon}N_{F}(\mathfrak{n})^{1/2+\varepsilon},\qquad N_{F}(\mathfrak{q})^{1+\varepsilon}
 N_{F}(\mathfrak{n})^{-1/2+\varepsilon},
\end{align*}
up to constants depending on $F,\mathcal{T},\varepsilon$.
\end{lemma}

\begin{proof}
The original poles on the shift circles of \textsection\ref{sec15}
are distinct and avoid $\pm1/2$. For the original bound in \eqref{14.125},
\eqref{eq:real-weight-evaluated} at $\beta_{v}=\mathbf{L}(\pi_{v})^{-1}$ gives
the real factors, and Lemma~\ref{lem:original-moving-residue-weight} at
$W_{v}^{[r_{v}]}$ gives the finite factors. The product is multiplied by
$\lambda_{\sigma_{\lambda}}(\mathfrak{n})$.

For \eqref{14.127}, take the positive Rankin--Selberg and standard circles
\begin{equation}\label{eq:dual-residue-circles}
 |\lambda-(1/2-s_{1}')|=3\varepsilon/25,\qquad
 |\lambda-(1/2-s_{2}')|=\varepsilon/100
\end{equation}
and their reflections. The Rankin--Selberg radius exceeds
$\max_{i}|\nu_{i}|\leq11\varepsilon/100$. The inequalities
\begin{align*}
|s_{1}'-s_{2}'|=\varepsilon,\qquad |s_{1}'|\geq\varepsilon/8,\qquad |s_{2}'|\geq5\varepsilon/8,\qquad |s_{1}'+s_{2}'|\geq\varepsilon/4
\end{align*}
separate the disks from each other and from $\pm1/2$ and the local
Gamma poles $\pm(1/2+s_{1}'+\nu_{i})$, $\pm(1/2+s_{2}')$.
For small enough $\varepsilon$, both reciprocals $1/\Lambda(1\pm2\lambda,\mathbf{1})$
are holomorphic on both disks $|\lambda\pm1/2|\leq3\varepsilon$,
by the pole at $0$ and nonvanishing at $2$.
At real places, the periods in \eqref{fc10.2} divided by
$L_{v}(1/2+s_{1}',\pi_{v}\times\sigma_{\lambda,v})$ and
$L_{v}(1/2+s_{2}',\sigma_{\lambda,v})$ are holomorphic.
The basis sum converges locally uniformly by angular differentiation
in \eqref{eq:maximal-GL2-Jacquet} and \eqref{eq:real-bessel-norm}.

The superscript $\reg$ in \eqref{eq:finite-residual-estimates} means that
\eqref{eq:residue-weight-interpolation} is applied to the normalized terms
of level $\ell\geq1$ in \eqref{eq:dual-level-superposition} at $v\mid\mathfrak{n}$.
The quotient is unchanged at $v\mid\mathfrak{q}$.
In the sum in \eqref{eq:dual-level-superposition},
the normalized $\ell=0$ term is $1$ by \eqref{fc11.42} and
\eqref{eq:finite-newvector-value}. For $\ell\geq1$,
Lemma~\ref{lem:original-moving-residue-weight} bounds the modified
term, including $q_{v}^{2\ell s_{2}}$, by $q_{v}^{O(\varepsilon\ell)}$.
At $v\mid\mathfrak{q}$, \eqref{eq12.14} and \eqref{eq:finite-newvector-value}
give the bound by $|\lambda_{\sigma_{\lambda,v}}(\mathfrak{p}_{v}^{j})|
\leq(j+1)q_{v}^{j|\Re(\lambda)|}$. Thus on \eqref{eq:dual-residue-circles}
and their reflections,
\begin{equation}\label{eq:finite-residual-estimates}
 \left|\tfrac{\mathcal{H}_{W_{v}}^{\vee,\reg}(\sigma_{\lambda,v},\mathbf{s}^{\vee},\mathbf{1})}
 {L_{v}(1/2+s_{1}',\pi_{v}\times\sigma_{\lambda,v})
 L_{v}(1/2+s_{2}',\sigma_{\lambda,v})}\right|
 =\left\{\begin{array}{@{}l@{\enspace}l@{}}
 O_{\varepsilon}(q_{v}^{-n_{v}(1/2+\Re(s_{2}))+O(\varepsilon n_{v})}),
 &v\mid\mathfrak{n},\\[2pt]
 O((1+r_{v})q_{v}^{(1+O(\varepsilon))r_{v}}),
 &v\mid\mathfrak{q}.
 \end{array}\right.
\end{equation}
The replacement preserves the residues by
Lemma~\ref{lem:original-moving-residue-weight}. On the circles the completed
$L$-factor quotient in \eqref{e4.6} is $O_{F,\varepsilon}(1)$.
In \eqref{eq:residue-cluster}, the product of \eqref{eq:finite-residual-estimates}
over $v\mid\mathfrak{n}\mathfrak{q}$ with the normalized real factors
therefore gives the stated dual bound.
\end{proof}

\subsubsection{Mixed-moment residues}\label{app:mixed-residual-estimates}

In Lemma~\ref{lem:mixed-fixed-residues}, only the standard $L$-factor
for the trivial inducing character gives residues. At $\lambda=1/2-s_{2}'$,
the completed factor is
\begin{equation}\label{eq:mixed-standard-completed-residue}
 \frac{\operatorname*{Res}_{z=1}\Lambda(z,\mathbf{1})}
 {\Lambda(2-2s_{2}',\mathbf{1})}
 \Lambda(s_{1}'+s_{2}',\pi)\Lambda(1+s_{1}'-s_{2}',\pi).
\end{equation}
For the original residue replace $s_{j}'$ by $s_{j}$.
Both factors are $O_{\pi,F}(1)$. At $t=\pm i/2$, $H_{R_{v}},S_{0}$ vanish to order
$\geq2$, $m_{z}$ has pole order $\leq1$, and the Gamma poles of
$L_{v}(1,\sigma_{t},\Ad)$ and $L_{v}(1/2,\sigma_{t})$ cancel in
$1/G_{\pi_{v}}$. This gives
$\mathcal{S}h_{v,R_{v},z}(\pm1/2)=\mathcal{S}h_{v,z}(\pm1/2)=0$
in \eqref{eq:mixed-polynomial-AFE-sieve}.

For the original residue, \eqref{eq:real-weight-evaluated} gives
$\prod_{v\mid\infty}\mathcal{S}h_{v}(1/2-s_{2})$ and
Lemma~\ref{lem:mixed-standard-pole-finite} gives the finite limits.
The bound follows from \eqref{14.36} with the Hecke twist.
For the dual residue, the finite factors at $\lambda=1/2-s_{2}'$
tend to
\begin{align*}
 v\mid\mathfrak{n}&:\quad q_{v}^{-n_{v}/2}g_{\pi}(\mathfrak{p}_{v}^{n_{v}})
 &&\text{by \eqref{eq:mixed-finite-endpoint}},\\
 v\mid\mathfrak{a}&:\quad q_{v}^{r_{v}}
 &&\text{by \eqref{eq12.14} at $(\mathbf{s},\lambda)=(\mathbf{0},\pm1/2)$}.
\end{align*}
Multiply these finite factors by \eqref{eq:mixed-standard-completed-residue}
and the real quotients in \eqref{eq:mixed-real-standard-pole-bound}
to obtain the dual bound in \eqref{eq:mixed-residual-bound}.
At a signed prime of depth one, the moving limit in
Lemma~\ref{lem:mixed-standard-pole-finite} replaces the dual factor $q_{v}$
by $-q_{v}^{-1}\alpha_{v}(\sigma_{1/2,v};1)$, which is also the signed
original factor. Its absolute value is $O(\alpha_{\pi}(\mathfrak{p}_{v})/q_{v})$
by \eqref{eq:old-row-simplified}.

\subsubsection{Mixed-moment residues at dual level one}\label{app:mixed-mass-reciprocity}

\begin{lemma}\label{lem:mixed-level-one-residue-bound}
For the data of \textsection\ref{app:mixed-mass-nonspectral}, the residue
terms in \eqref{8.1} at dual level one are
\begin{align*}
 O_{\pi,F,\varepsilon}\Big(N_{F}(\mathfrak{a})^{1+\varepsilon}
 \prod_{v\mid\infty}R_{v}^{O(1)}\Big).
\end{align*}
\end{lemma}

\begin{proof}
The positive poles are $1/2-s_{2}$ for $\xi=\mathbf{1}$ and
$1/2-s_{1}$ for $\xi=\chi^{-1}$, with $\chi$ a nontrivial quadratic
constituent of $\pi$ and with the inducing shifts before specialization.
For $|s_{1}|=\delta/2$, $|s_{2}|=\delta/4$ and shifts of modulus
$<\delta/1000$, the dual poles lie inside the circles
$|\lambda-(1/2-s_{j}')|=\delta/100$ and their reflections,
as in \eqref{eq:dual-residue-circles}.
On these circles the completed $L$-factor quotient in
$\widetilde{\mathcal{G}}$ is $O_{\pi,F,\delta}(1)$. In the two real periods
of \eqref{fc10.1} at $(\pi_{v}(w_{2})W_{v},\mathbf{s}^{\vee})$,
the compact test data of \textsection\ref{app:mixed-mass-nonspectral}
have Sobolev seminorms $R_{v}^{O(1)}$. The periods divided by their
local $L$-factors are consequently $O(R_{v}^{O(1)})$ on these circles.

At each original Rankin--Selberg pole, the real period is holomorphic
by compact support and Schwartz convolution, so its quotient by the
local Rankin--Selberg $L$-factor vanishes. At the standard pole the terms $l\geq2$ in
\eqref{equ12.11} vanish, and Lemma~\ref{lem:mixed-standard-pole-finite}
gives the $l=1$ value. At finite depth $m$,
the original and dual bounds are
\begin{align*}
 &q_{v}^{(\vartheta_{\pi_{v}}+O(\delta))m}
 &&\text{by Lemma~\ref{lem:original-moving-residue-weight} (first five
 $l_{2}$ coefficients) and \eqref{eq:newform-Schur-bound}},\\
 &(1+m)q_{v}^{(1+O(\varepsilon))m}
 &&\text{by \eqref{eq12.14} at $r_{v}=m$, as in \eqref{eq:finite-residual-estimates}}.
\end{align*}
Multiply these local estimates and sum the finitely many signed
integrals \eqref{eq:residue-cluster} before taking the central constant term.
\end{proof}

\subsection{Local Weights on the Residue Contours}\label{app:polynomial-residual-transform}

\begin{lemma}\label{lem:real-residue-contour-weights}
Take $h_{v}=h_{v,R}$ in \eqref{fc11.3}. For sufficiently large $m$
and $\boldsymbol{\nu}=\mathbf{0}$, uniformly on the circles in
\eqref{eq:dual-residue-circles} and their reflections,
\begin{align*}
 |\pi\mathcal{H}_{W_{v}}^{\vee}
 (\sigma_{\lambda,v},\mathbf{s}^{\vee},\mathbf{1})|
 \ll_{\varepsilon}R^{3+O(\varepsilon)}.
\end{align*}
The other real factors satisfy \eqref{eq:hybrid-real-factors}.
On the $\boldsymbol{\nu}$ circles of \textsection\ref{sec15}, the
$s_{2}$ constant term of \eqref{e14.126} vanishes when taken first.
\end{lemma}

\begin{proof}
We first derive an integral for the dual weight in a region of absolute
convergence.
Write $c=1/2-s_{1}$, $p=1/2-\lambda+s_{1}'$ and
$q_{1,2}=1/2+\lambda\pm s_{2}'$. With $|s_j|\leq 2^{-j}\varepsilon$,
we first take $\lambda$ real with $|\lambda-(1/2-2\varepsilon)|<\varepsilon/8$, so
$\Re(\lambda),\Re(p)>0$ and the Hecke exponents $q_1$ and $1-q_2$
have real parts in $(0,1)$.
The Fourier calculation for \eqref{eq:principal-cosine-kernel},
with Mellin power $q_1$ in the second period, gives the angular weight
$C_0(q_1)|\sin\theta|^{-q_1}|\cos\theta|^{-q_2}$.
The decay \eqref{eq:real-inverse-Abel-decay} gives the moments of $f_v$
needed to sum the $K_v'$-types. The Mellin powers $p$ and $1+2s_1'$
in \eqref{eq:central-Whittaker-two-K-kernel} give an integrable $s$ integral.
After inserting \eqref{eq11.26}, the $z$ power is $c$, with $\Re(c)>0$,
and the singular powers at $\cos(\arctan x)=0$ and
$\sin(\arctan x-\theta)=0$ sum to
$1-2\Re(\lambda)+1/2+\Re(s_2)<1$.

In \eqref{10.3}, the $a_{1},z_{1}$ integral is $4f_v(b_1)$ by
\eqref{eq10.17}. With $\phi=\arctan x-\theta$, the $b_1$ integral
gives $\widehat f_v(z\sin\phi)$, and the Bessel pairing below evaluates
the $z$ integral. The changes $x=\sqrt{1+s^{-1}}\,y$, $s=z^{2}$,
$v=\sqrt{1+z^{2}}\cot\phi$ and $(z,v)=(r\sin\alpha,r\cos\alpha)$ then give
\begin{multline}\label{eq:shifted-residual-radial-transform}
 \pi\mathcal{H}_{W_{v}}^{\vee}(\sigma_{\lambda,v},\mathbf{s}^{\vee},\mathbf{1})
 =\frac{\pi^{s_{1}-s_{2}-1}C_{0}(q_{1})}{4\Gamma(c)}
 \int_{\mathbb{R}}H_{R}(t)t\tanh(\pi t)
 \frac{A_{s_{2}}(t)A_{-s_{1}}(t)^{2}}{A(t)^{4}}\\
 2\int_{0}^{\infty} r^{2c-1}F_{c,t}(r)\mathcal{U}(r^{2})\,dr\,dt,
\end{multline}
where $A(t)$ is defined in \eqref{eq:central-input-Ft} and
\begin{align*}
&A_{z}(t):=\Gamma(1/4+z/2+it/2)\Gamma(1/4+z/2-it/2),\\
&F_{c,t}(r):={}_{2}F_{1}(c/2+it/2,c/2-it/2;c;-r^{2}).
\end{align*}
Thus $F_{1/2,t}=F_{t}$ from \eqref{eq:central-input-Ft}.
For $w>0$, define
\begin{multline*}
\mathcal{U}(w):=(1+w)^{\lambda+s_{1}'}\int_{0}^{\pi}\int_{\mathbb{R}}
B_{i(p-1/2),0}(1,y)(\sin\alpha)^{-1-2s_{1}}(1+w\sin^{2}\alpha)^{-s_{2}}\\
|\cot\alpha\,y-1|^{-q_{1}}
 |(1+w\sin^{2}\alpha)y+w\sin\alpha\cos\alpha|^{-q_{2}}
 \,dy\,d\alpha,
\end{multline*}
with $B_{u,0}$ as in \eqref{eq:central-K-Fourier-Mellin-kernel}.
In \eqref{eq:shifted-residual-radial-transform}, the $K_{it}$ and local
$L$-factor in \eqref{eq11.26} combine with the first $K_{0}$ in
\eqref{eq:central-Whittaker-two-K-kernel} as
\begin{align*}
L_{v}(1/2+s_{2},\sigma_{t})
\int_{0}^{\infty}K_{it}(2\pi b)K_{0}(2\pi b\sqrt{1+r^{2}})b^{c-1}\,db=\frac{A_{s_{2}}(t)A_{-s_{1}}(t)^{2}F_{c,t}(r)}
 {8\pi^{1-s_{1}+s_{2}}\Gamma(c)}.
\end{align*}

With $\cot\alpha\mapsto x$ and $D=1+w$, the $y$ roots are
\begin{align*}
 r_{1}=1/x,\quad r_{2}=-\frac{(D-1)x}{x^{2}+D},\qquad
 r_{1}-r_{2}=\frac{D(1+x^{2})}{x(x^{2}+D)}.
\end{align*}
Split $\mathcal{U}=\mathcal{U}_{0}+\mathcal{U}_{\reg}$ by writing
$B_{i(p-1/2),0}(1,y)=B_{i(p-1/2),0}(1,0)+g(y)$, where
\begin{align*}
 g(y):=\frac{\Gamma(p/2)^{2}}{2\pi^{p}}
 \bigl[{}_{2}F_{1}(p/2,p/2;1/2;-y^{2})-1\bigr].
\end{align*}
In $\mathcal{U}_{0}$, the $y$ integral with $y\mapsto r_{1}+(r_{2}-r_{1})u$
is a beta integral. The $\alpha$ integral then gives
\begin{multline*}
 \mathcal{U}_{0}(w)=\frac{B_{i(p-1/2),0}(1,0)C_{0}(2\lambda)}{C_{0}(q_{1})C_{0}(q_{2})}\cdot
 \frac{\Gamma(p/2)\Gamma(q_{2}/2)}{\Gamma((p+q_{2})/2)}\\
 (1+w)^{p-1/2}{}_{2}F_{1}(p-s_{1},p/2;c;-w).
\end{multline*}
$B_{i(p-1/2),0}(1,0)=\Gamma(p/2)^{2}/(2\pi^{p})$.
On \eqref{eq:dual-residue-circles}, the connection formulas
\cite[(15.8.2), (15.8.4)]{DLMF} at $(p-s_{1},p/2;c)$ give
\begin{align*}
 |\mathcal{U}_{0}(w)|\ll_{\varepsilon}
 \begin{cases}
 |1+w|^{\Re(p)-1/2}+|1+w|^{-\Re(p)/2},& w\to-1,\\
 |w|^{-\Re(c)}+|w|^{-1/2+\Re(p)/2},& |w|\to\infty.
 \end{cases}
\end{align*}

Continue $\mathcal{U}_{\reg}$ from $D>0$ along $D=\rho e^{i\theta}$,
using $x=e^{i\theta/2}u$, $u\in\mathbb{R}$. For $u>0$ and $0<\theta<\pi$,
$\Re(r_{1})>\Re(r_{2})$ and neither root meets $-i[1,\infty)$.
If the roots are on opposite sides of the imaginary axis, then $\rho\geq1$
and at least one is above $-i$: $|\Im(r_{1})|<1$ for $u\geq1$ and
$|\Im(r_{2})|<1$ for $u\leq1$. The path between the roots can therefore
be chosen with length $O(|r_{1}-r_{2}|)$, crossing the imaginary axis
at the height of the higher root if needed. Continue the branches on
the three $y$ paths from $D>0$, and use reflected paths for $u<0$
and for the lower side of the cut.
At each root, subtracting $g(r_j)$ locally leaves an integrable remainder
for $\Re(q_j)<1+1/16$ by the exponent $\Re(1/2-p)>1/16$ at $y=\pm i$.
The constant terms continue by the beta integrals
\cite[(5.12.1), (5.12.3)]{DLMF} after $y=r_1+(r_2-r_1)u$.
With $g(y)=O_{\varepsilon}((1+|y|)^{1/16})$ at infinity
\cite[(15.8.4), (15.8.8) with $m=0$]{DLMF}, this gives, on
\eqref{eq:dual-residue-circles},
\begin{equation}\label{eq:residual-fiber-bound}
 \left|\int
 g(y)(y-r_{1})^{-q_{1}}(y-r_{2})^{-q_{2}}\,dy\right|
 \ll_{\varepsilon} |r_{1}-r_{2}|^{1-\Re(q_{1}+q_{2})}
 (1+|r_{1}|+|r_{2}|)^{1/8}.
\end{equation}

By \eqref{eq:residual-fiber-bound}, the $|x|\leq H$ part of
$\mathcal{U}_{\reg}$ is $O_{\varepsilon}(\rho^{-7/16+O(\varepsilon)})$
for $\rho\geq3$, $H=10\rho$, and $O_{\varepsilon}(\rho^{-1/2-O(\varepsilon)})$
for $\rho\leq1/3$, $H=4$.

For $|x|\geq H$, $g(y)=O_{\varepsilon}(y^{2})$ at $0$ gives an
$O_{\varepsilon}(1)$ contribution to the $y$ integral from
$|y|\gg |r_{1}|+|r_{2}|$. The subtraction at the roots in
\eqref{eq:residual-fiber-bound} bounds the rest by
$O_{\varepsilon}(|r_{1}-r_{2}|^{-2\Re(\lambda)}(|r_{1}|+|r_{2}|)^{2})$.
Since $|r_{1}-r_{2}|\asymp\rho/|x|$ and
$|r_{1}|+|r_{2}|\ll(1+\rho)/|x|$, the $x$ integrand of $\mathcal{U}_{\reg}$
is bounded by
\begin{align*}
 \rho^{\Re(p+2\lambda)-1/2}|x|^{-1-\Re(p+2\lambda)}
 +\rho^{\Re(p)-1/2}(1+\rho)^{2}|x|^{-3-\Re(p)}.
\end{align*}
Integrating for $|x|\geq H$ gives
\begin{equation}\label{eq:residual-constant-tail}
 \left|\left.\mathcal{U}_{\reg}(w)\right|_{|x|\geq H}\right|
 \ll_{\varepsilon}
 \begin{cases}
 \rho^{-1/2},&\rho\geq3,\\
 \rho^{\Re(p)-1/2},&\rho\leq1/3.
 \end{cases}
\end{equation}
The formula for $\mathcal{U}_{0}$ and the subtractions used for
\eqref{eq:residual-fiber-bound} continue $\mathcal{U}(r^{2})$
meromorphically to $|\lambda-1/2|<3\varepsilon$. The same $x$-integral
bounds on $w=r^{2}\geq0$ give $\mathcal U(r^{2})\ll(1+r)^{-3/4}$,
locally uniformly away from the poles. Euler's integral gives
$|F_{c,t}(r)|\ll e^{C|t|}(1+r)^{-\Re(c)}(1+\log(2+r))$ for a fixed $C$.
Together with the decay of $H_R$, this gives locally uniform convergence
in \eqref{eq:shifted-residual-radial-transform}, which therefore holds on
\eqref{eq:dual-residue-circles} by meromorphic continuation.

We take $\theta\to\pi^{-}$ and $\theta\to(-\pi)^{+}$ on the same $u$ line,
keeping the one-sided values on the continued $y$ paths.
For $1/3\leq\rho\leq3$ and $|u|\leq30$,
\eqref{eq:residual-fiber-bound} bounds the $u$ integrand by
\begin{align*}
 |1-u^{2}|^{\Re(p)-1/2}|u|^{-\Re(p+2s_{1})}(1+|u|^{-1})^{1/8},
\end{align*}
using $|1+x^{2}|\geq|1-u^{2}|$. This is integrable even at $\rho=1$,
and the tail estimate preceding \eqref{eq:residual-constant-tail} is uniform
for $|u|>30$ in this range of $\rho$. The two boundary values follow
by dominated convergence on the rotated lines.
Cauchy's formula and the bounds as $\rho\to0^{+}$ and $\rho\to\infty$,
with $\rho=u-1$, give
\begin{equation}\label{eq:residual-kernel-cauchy-norm}
 \begin{gathered}
 \mathcal{U}(w)=\int_{1}^{\infty}\frac{\mu(u)}{u+w}\,du,\qquad
 \mu(u)=\frac{\mathcal{U}(-u-i0)-\mathcal{U}(-u+i0)}{2\pi i},\\[2pt]
 \int_{1}^{\infty}|\mu(u)|u^{\Re(c)/2-1}(1+\log u)\,du\ll_{\varepsilon}1.
 \end{gathered}
\end{equation}

Insert \eqref{eq:residual-kernel-cauchy-norm} into
\eqref{eq:shifted-residual-radial-transform}. The Euler bound gives an
absolute $r$ integral of size $O_{t}(u^{\Re(c)/2-1}(1+\log u))$, so
\eqref{eq:residual-kernel-cauchy-norm} allows integration first in $r$, giving
\begin{align*}
 \int_{0}^{\infty}\frac{r^{2c-1}F_{c,t}(r)}{u+r^{2}}\,dr
 =\frac{\Gamma(c)A_{1+s_{1}}(t)}{2}
 {}_{2}F_{1}(1-c/2+it/2,1-c/2-it/2;1;1-u).
\end{align*}
The change $r\mapsto\sqrt{(1-\xi)/\xi}$ and Pfaff's formula
\cite[(15.8.1)]{DLMF} reduce the $r$ integral to
\cite[(3.22), $z=1-u$, $(a,b,c,d,e)=(1,1-c/2+it/2,1-c/2-it/2,1,1)$]{Koo15}.
For $t\in\mathbb{R}$, \cite[(2.21) and Lemma~2.3, $n=0$]{Koo75}
with $(\alpha,\beta,\lambda)=(0,1-c,t)$ at $\operatorname{arcosh}\sqrt{u}$ give
\begin{align*}
 \left|{}_{2}F_{1}\left(1-\tfrac{c-it}{2},1-\tfrac{c+it}{2};1;1-u\right)\right|
 \ll_{\varepsilon}
 \begin{cases}
 1,&1\leq u\leq2,\\
 u^{\Re(c)/2-1}(1+\log u),&u\geq2.
 \end{cases}
\end{align*}
Integrating against $|\mu(u)|\,du$ and applying \eqref{eq:residual-kernel-cauchy-norm} gives
\begin{multline*}
 \left|2\int_{0}^{\infty}r^{2c-1}F_{c,t}(r)\mathcal{U}(r^{2})\,dr\right|
 \ll_{\varepsilon}|\Gamma(c)A_{1+s_{1}}(t)|
 \int_{1}^{\infty}|\mu(u)|u^{\Re(c)/2-1}(1+\log u)\,du\\
 \ll_{\varepsilon}|\Gamma(c)A_{1+s_{1}}(t)|.
\end{multline*}
In \eqref{eq:shifted-residual-radial-transform}, the factor $|\Gamma(c)|$
cancels. Stirling and \eqref{eq:R-admissible-test} give
\begin{equation}\label{eq:polynomial-shifted-residual-bound}
\begin{aligned}
 |\pi\mathcal{H}_{W_{v}}^{\vee}(\sigma_{\lambda,v},\mathbf{s}^{\vee},\mathbf{1})|
 &\ll_{\varepsilon}\int_{\mathbb{R}}|H_{R}(t)|\,|t|
 \left|\frac{A_{s_{2}}(t)A_{-s_{1}}(t)^{2}A_{1+s_{1}}(t)}{A(t)^{4}}\right|\,dt\\
 &\ll_{\varepsilon}\int_{\mathbb{R}}|H_{R}(t)|(1+|t|)^{2-\Re(s_{1})+\Re(s_{2})}\,dt
 \ll_{\varepsilon}R^{3+O(\varepsilon)}.
\end{aligned}
\end{equation}
The same bound holds on the reflected circles since
$\sigma_{-\lambda,v}\simeq\sigma_{\lambda,v}$.
For $h_{v,R}$, the real factors satisfy
\begin{equation}\label{eq:hybrid-real-factors}
 \renewcommand{\arraystretch}{1.35}
 \begin{array}{c|ccc}
 &\text{Generic (original)}&\text{Degenerate}&\text{Residual}\\ \hline
 \text{at a real place}&O(R^{2+O(\varepsilon)})&O(R^{2+O(\varepsilon)})&O(R^{3+O(\varepsilon)})
 \end{array}
\end{equation}
In \eqref{e14.30} and \eqref{c14.4}, the substitution
$\mathcal{S} h_{v,R}(it)=\pi^{2}H_{R}(t)/(A(t)^{4}\cosh(\pi t))$
gives the first two columns of \eqref{eq:hybrid-real-factors} by Stirling
and \eqref{eq:R-admissible-test}. The crossed residues and evaluations in
\eqref{14.4}, \eqref{14.125} are $O(R^{-m})$.

In \eqref{e14.126}, $|s_{2}|+\max_{i}|\nu_{i}|<|s_{1}|$ keeps the possible
$L$-factor poles $s_{2}=\pm(s_{1}+\nu_{i})$ outside the $s_{2}$ circle.
The $s_{2}$ constant term is the value at zero, where the real factor
$4\widehat{f}_{v}(0)=\mathcal{S}h_{v,R}(1/2+s_{2})$ from
\textsection\ref{sec14.7} and \eqref{11.22} vanishes by the zeros of $H_{R}$.
\end{proof}

\subsection{Real weights for the moment estimates}
\label{app:real-weights-moment-estimates}
For a unitary irreducible representation $\sigma_{v}$ of $\mathrm{PGL}_{2}(\mathbb{R})$,
let $\lambda_{\sigma_{v}}$ be its eigenvalue for the $\mathfrak{sl}_{2}$
Casimir $\Delta_{v}$. Then
$C(\sigma_{v})^{1/2}\asymp1+|\lambda_{\sigma_{v}}|^{1/2}$.

For functions of $k\in K_{v}'$, write
\begin{align*}
\partial_{k}f(k):=\left.\frac{d}{d\theta}
f\bigl(k\left(\begin{smallmatrix}\cos\theta&\sin\theta\\
-\sin\theta&\cos\theta\end{smallmatrix}\right)\bigr)\right|_{\theta=0}.
\end{align*}

\begin{lemma}
\label{lem:uniform-real-dual-decay}
For $|s_{j}|\leq2^{-j}\varepsilon$, $|\nu_{j}|\leq10^{-j}\varepsilon$
with $j=1,2$, and every $\ell>0$,
\begin{equation}\label{eq:uniform-real-dual-decay}
 \mathcal{H}_{W_{v}}^{\vee}(\sigma_{v},\mathbf{s}^{\vee},\mathbf{1})
 \ll_{\ell,h_{v},\pi_{v},\varepsilon,\vartheta_{\infty}}
 C(\sigma_{v})^{-\ell/2},
\end{equation}
uniformly for tempered $\sigma_{v}$ and for complementary series $\sigma_{v}$
with real inducing exponents $\pm\eta_{v}$, $|\eta_{v}|\leq\vartheta_{\infty}$.
\end{lemma}

\begin{proof}
With \eqref{eq10.17}--\eqref{c11.20} and $x=z\diag(a,1)k$, \eqref{10.3} becomes
\begin{align*}
 \mathcal{H}_{W_{v}}^{\vee}(\sigma_{v},\mathbf{s}^{\vee},\mathbf{1})
 =\sum_{n}\overline{\Psi(W_{n};\overline{s_{2}'})}
 \int_{\mathbb{R}^{\times}}\int_{K_{v}'}
 \Phi_{\mathbf{s}}(a,k)W_{n}(\diag(a,1)k)\frac{dk\,d^{\times}a}{|a|},
\end{align*}
where
\begin{align*}
 \Phi_{\mathbf{s}}(a,k):=4|a|^{s_{1}'}\int_{\mathbb{R}^{\times}}
 W_{v}^{\circ}(\diag(za,z,1))|z|^{2s_{1}'}
 \widehat{f}_{v}(zk_{21})\,d^{\times} z.
\end{align*}

By Stirling's formula, taking absolute values in \eqref{10.19} gives
\begin{align*}
 |W_{v}^{\circ}(\diag(za,z,1))|\ll_{c_{1},c_{2},\pi_{v}}
 |z|^{1-c_{1}}|a|^{1-c_{2}},
\end{align*}
also after logarithmic derivatives in $a,z$.
Right differentiation in Iwasawa coordinates and integration in $z$
give, for every $M>0$,
\begin{align*}
 |(1+\Delta_{v}^{2})^{r}(1-\partial_{k}^{2})^{j}\Phi_{\mathbf{s}}(a,k)|
 \ll_{r,j,M,\varepsilon,h_{v},\pi_{v}}
 \begin{cases}
 |a|^{1-2\varepsilon},&0<|a|\leq1,\\
 |a|^{-M},&|a|\geq1.
 \end{cases}
\end{align*}
For $\vartheta_{\infty}<\eta<1/2-2\varepsilon$, the square is integrable
against $|a|^{-1-2\eta}d^{\times}a\,dk$, since the bound at zero gives
$\int_{0}^{1}a^{-4\varepsilon-2\eta}\,da<\infty$.
On $W_{n}(\diag(a,1)k)$, the operators $(1+\Delta_{v}^{2})^{r}$
and $(1-\partial_{k}^{2})^{j}$ act by $(1+|\lambda_{\sigma_{v}}|^{2})^{r}$
and $(1+4n^{2})^{j}$. Moving them to $\Phi_{\mathbf{s}}$ in the
$a,k$ pairing and applying Cauchy--Schwarz with
\eqref{eq:uniform-kirillov-periods} gives
\begin{align*}
 \Big|\int_{\mathbb{R}^{\times}}\int_{K_{v}'}
 \Phi_{\mathbf{s}}(a,k)W_{n}(\diag(a,1)k)|a|^{-1}\,dk\,d^{\times}a\Big|
 \ll_{r,j,\varepsilon,h_{v},\pi_{v}}
 \frac{(1+|n|)^{\vartheta_{\infty}}}{(1+|\lambda_{\sigma_{v}}|^{2})^{r}(1+4n^{2})^{j}}.
\end{align*}
By \eqref{eq:uniform-kirillov-periods} and
\cite[\textsection2.3.3, Proposition~3.2.3]{MV10},
\begin{align*}
 |\Psi(W_{n};\overline{s_{2}'})|
 \ll_{\varepsilon,\vartheta_{\infty}}
 (1+C(\sigma_{v})^{1/2}+|n|)^{m}
\end{align*}
for an $m$ independent of $r,j$. Multiplying the two estimates gives
\begin{align*}
 \left|\mathcal{H}_{W_{v}}^{\vee}(\sigma_{v},\mathbf{s}^{\vee},\mathbf{1})\right|
\ll C(\sigma_{v})^{m/2-2r}
 \sum_{n}(1+|n|)^{m+\vartheta_{\infty}-2j}.
\end{align*}
Taking $2j>m+\vartheta_{\infty}+1$ and $4r>\ell+m$ proves
\eqref{eq:uniform-real-dual-decay}.
\end{proof}
\begin{proof}[Proof of Lemma~\ref{lem:hybrid-afe-multiplier}]
The poles of $m_{z}$ at $t=\pm i(1/2+2\ell)$ are at most simple, so
$H_{R}m_{z}$ still vanishes to order at least four.
Stirling and Cauchy's formula give \eqref{eq:R-admissible-test} at
scale $\sqrt{2}R$, up to the factor
$R^{c}(1+|\Im(z)|)^{O(1)}e^{\pi|\Im(z)|/2}$.
The transform and residue estimates use only these analytic bounds and
the zeros, not positivity. For \eqref{eq:original-weighted-fourth}, sum
the fourth-moment bound against the pointwise bound for $|H_{R}m_{z}|$.
This gives the same factor in Lemma~\ref{lem:hybrid-arch-bookkeeping}.
Multiplying the bounds at the real places
and including $|e^{z^{2}}|$ gives \eqref{eq:afe-multiplier-vertical-bound},
since $\prod_{v\mid\infty}R_{v}^{c}=X^{c/2}$.
\end{proof}

\subsubsection{Real weights for the mixed moment}
For functions $h$ on $G'(\mathbb{R})$, write
\begin{align*}
 \mathscr{L}_{b}(h):=\int_{G'(\mathbb{R})}|h(g)|
 (1+\|g\|+\|g^{-1}\|)^{b}\,dg,
\end{align*}
where $\|\cdot\|$ is the Euclidean operator norm.
\begin{lemma}
\label{lem:mixed-fixed-real-dual-decay}
Fix a spherical unitary $\pi_{v}$ of $\mathrm{PGL}_{3}(\mathbb{R})$
and $0\leq\theta_{0}<1/2$ with
$\vartheta_{\pi_{v}}+\theta_{0}<1/2$. Let $W_{v}$ be the vector defined
from $h_{v}$ by \eqref{fc11.3} at $\mathbf{s}=\mathbf{0}$.
For every $\ell>0$ there is $b=b(\pi_{v},\ell,\theta_{0})$ such that
\begin{equation}\label{eq:mixed-fixed-real-dual-decay}
 |\mathcal{H}_{W_{v}}^{\vee}
 (\sigma_{v},\mathbf{0},\mathbf{1})|
 \ll_{\pi_{v},\ell,\theta_{0}}
 \mathscr{L}_{b}(h_{v})C(\sigma_{v})^{-\ell/2}.
\end{equation}
The estimate includes both principal-series parities, discrete series,
and complementary series with real inducing exponents $\pm\eta_{v}$,
$|\eta_{v}|\leq\theta_{0}$.

For $h_{v,z}$ constructed from $S_{z}$ in
Lemma~\ref{lem:mixed-afe-realization} with strip parameter $l>b+1/2$, we have
\begin{equation}\label{eq:mixed-afe-local-decay}
 |\mathcal{H}_{W_{h_{v,z}}}^{\vee}(\sigma_{v},\mathbf{0},\mathbf{1})|
 \ll_{\pi_{v},\ell,l,R_{0},c}
 (1+|\Im(z)|)^{O(1)} e^{\pi|\Im(z)|/2}
 C(\sigma_{v})^{-\ell/2}.
\end{equation}
\end{lemma}

\begin{proof}
For a smooth Whittaker vector $V$ of $\pi_{v}$ and $k\in K_{v}'$, define
\begin{align*}
 \Phi_{V}(a,k):=\int_{\mathbb{R}^{\times}}
 V(\diag(az,z,1)k)\,d^{\times} z.
\end{align*}
Integration by parts in $k$ gives
\begin{multline*}
 \Psi_{v}(1/2,V,W_{n})
 =\frac{1}{(1+4n^{2})^{2}}\int_{\mathbb{R}^{\times}}\int_{K_{v}'}
 ((1-\partial_{k}^{2})^{2}\Phi_{V})(a,k)\\
 W_{n}(\diag(a,1)k)|a|^{-1}\,dk\,d^{\times}a.
\end{multline*}
For $\theta_{0}<\eta<1/2-\vartheta_{\pi_{v}}-\varepsilon$,
Cauchy--Schwarz and \eqref{eq:uniform-kirillov-periods} give
\begin{align*}
 |\Psi_{v}(1/2,V,W_{n})\overline{\Psi(W_{n};0)}|
 \ll_{\eta,\theta_{0}}\frac{(1+|n|)^{2\theta_{0}}}{(1+4n^{2})^{2}}\cdot  \|(1-\partial_{k}^{2})^{2}\Phi_{V}\|_{L^{2}(|a|^{-1-2\eta}d^{\times}a\,dk)}.
\end{align*}
Summing \eqref{fc10.1} and expanding $(1-\partial_{k}^{2})^{2}$ gives
\begin{equation}\label{eq:mixed-full-K-sobolev-bound}
 |\mathcal{H}_{V}(\sigma_{v},\mathbf{0},\mathbf{1})|
 \ll_{\eta,\theta_{0}}\sum_{j=0}^{4}
 \|\partial_{k}^{j}\Phi_{V}\|_
 {L^{2}(|a|^{-1-2\eta}d^{\times} a\,dk)}.
\end{equation}
Transferring $(1+\Delta_{v}^{2})^{m}$ from the Jacquet vector to $V$
in \eqref{fc10.1} gives
\begin{align*}
 \mathcal{H}_{V}(\sigma_{v},\mathbf{0},\mathbf{1})
 =\frac{\mathcal{H}_{(1+\Delta_{v}^{2})^{m}V}(\sigma_{v},\mathbf{0},\mathbf{1})}
 {(1+|\lambda_{\sigma_{v}}|^{2})^{m}}.
\end{align*}

For $a,z>0$, the Whittaker asymptotics for the fixed representation
$\pi_{v}$ \cite[\textsection 15.2]{Wal92} give, for every fixed element
$U$ of the enveloping algebra and every $n$,
\begin{equation}\label{eq:mixed-fixed-Whittaker-bound}
|(d\pi_{v}(U)W_{v}^{\circ})(\diag(az,z,1)k)|
 \ll_{U,\varepsilon,n}(az)^{1-\vartheta_{\pi_{v}}-\varepsilon}
 (1+a)^{-n}(1+z)^{-n},
\end{equation}
uniformly for compact $k$. Differentiation along $n_{12}$ and $n_{23}$
gives the rapid decay by the factors $a$ and $z$ in the Whittaker character.

For $V=\pi_{v}(w_{2})W_{v}$ in \eqref{fc11.3}, the last one and
two row norms bound the change of the Iwasawa ratios under $g$ by
powers of $1+\|g\|+\|g^{-1}\|$. Fixed derivatives act via
$\operatorname{Ad}((w_{2}g)^{-1})$. Integrating
\eqref{eq:mixed-fixed-Whittaker-bound} first in $z$ and then against
$|h_{v}(g)|\,dg$ gives, for $0\leq j\leq4$ and every $N>0$,
\begin{align*}
 |\partial_{k}^{j}\Phi_{(1+\Delta_{v}^{2})^{m}V}(a,k)|
 \ll_{\pi_{v},m,j,N,\varepsilon}
 \mathscr{L}_{b}(h_{v})\,|a|^{1-\vartheta_{\pi_{v}}-\varepsilon}(1+|a|)^{-N}.
\end{align*}
The weighted $L^{2}$ norms in \eqref{eq:mixed-full-K-sobolev-bound} are
therefore $O_{\pi_{v},m,\eta}(\mathscr{L}_{b}(h_{v}))$.
With $1+|\lambda_{\sigma_{v}}|^{2}\asymp C(\sigma_{v})^{2}$ take $4m\geq\ell$.

On the support of $\alpha$, \eqref{eq:real-inverse-Abel-decay} and
\eqref{eq:mixed-afe-test-seminorms} bound $\mathscr{L}_{b}(h_{v,z})$ by
\begin{align*}
 (1+|\Im(z)|)^{O(1)}e^{\pi|\Im(z)|/2}
 \int_{0}^{\infty}e^{-(l-b-1/2)r}\,dr .
\end{align*}
Applying \eqref{eq:mixed-fixed-real-dual-decay} gives \eqref{eq:mixed-afe-local-decay}.
\end{proof}

\begin{lemma}
\label{lem:mixed-input-Sobolev}
Fix $U$ in the enveloping algebra of $\mathfrak{gl}_{3}(\mathbb{R})$.
Take $h_{R}$ with $\mathcal{S}h_{R}(it)=H_{R}(t)/G_{\pi_{v}}(t)$ and the
fixed central cutoff of \textsection\ref{def:central-test-cutoff}.
Assume \eqref{eq:R-admissible-test} with $m$ sufficiently large in terms
of $\pi_{v},U$. For $W_{h_{R}}$ from \eqref{fc11.3} at $\mathbf{s}=\mathbf{0}$,
\begin{equation}\label{eq:mixed-input-Sobolev}
 \|d\pi_{v}(U)W_{h_{R}}\|_{\pi_{v}}
 +\|d\pi_{v}(U)\pi_{v}(w_{2})W_{h_{R}}\|_{\pi_{v}}
 \ll_{\pi_{v},U}R^{O_{\pi_{v},U}(1)}.
\end{equation}
For fixed $c=\Re(z)>0$, replacing $H_{R}$ by $H_{R}m_{z}$ multiplies the bound by
$R^{c}(1+|\Im(z)|)^{O(1)}e^{\pi|\Im(z)|/2}$.
\end{lemma}

\begin{proof}
The zeros of $H_{R}$ cancel the poles of $1/G_{\pi_{v}}$ in the strip.
The argument for \eqref{eq:real-inverse-Abel-decay} with $a>b+1/2$
then gives $\mathscr{L}_{b}(h_{R})<\infty$, also with $H_{R}m_{z}$.
The adjoint-action estimate in the
proof of Lemma~\ref{lem:mixed-fixed-real-dual-decay} gives convergence of
the compact truncations in the fixed Sobolev norms used below.
For $a,r>0$, put $f_{R}(a,r):=W_{h_{R}}(\diag(ar,r,1))$ and define
\begin{align*}
 \mathcal{M}f_{R}(u_{1},u_{2})
 :=\int_{0}^{\infty}\int_{0}^{\infty}
 f_{R}(a,r)a^{u_{1}}r^{u_{2}}\,\frac{da}{a}\frac{dr}{r}.
\end{align*}

For $\Re(u_{2})$ sufficiently large and $\Re(u_{1})>\Re(u_{2})/2-1/2$,
spherical inversion with \eqref{eq:real-rs-evaluated} at
$(s_{1},\nu)=(u_{2}/2,it)$ and \eqref{10.42} at $\alpha=u_{1}+1/2-u_{2}/2$ gives
\begin{multline}\label{eq:mixed-input-joint-Mellin}
 \mathcal{M} f_{R}(u_{1},u_{2})
 =\frac{\beta_{v}}{16\pi^{2}}
 \left(\int_{\mathbb{R}^{\times}}\alpha(r)|r|^{u_{2}}\,d^{\times}r\right)
 \int_{0}^{\infty}\mathcal{S} h_{R}(it)t\sinh(\pi t)\\
 \prod_{j=1}^{3}\prod_{\epsilon=\pm1}
 \Gamma_{\mathbb{R}}(1/2+u_{2}/2+\nu_{j}+\epsilon it)
 \prod_{\epsilon=\pm1}
 \Gamma_{\mathbb{R}}(u_{1}+1/2-u_{2}/2+\epsilon it)\,dt.
\end{multline}
In \eqref{eq:mixed-input-joint-Mellin}, put
$b_{j}=1/2+u_{2}/2+\nu_{j}$ for $j\leq3$ and $b_{4}=u_{1}+1/2-u_{2}/2$.
For each shift $(j,\ell)$ in \eqref{eq:mixed-input-Lie-Mellin-shifts},
move the Mellin lines to
\begin{align*}
\Re(u_{1})=-1/2+j\pm\epsilon_{1},\qquad \Re(u_{2})=\ell\pm\epsilon_{2},
\end{align*}
with small $\epsilon_{1},\epsilon_{2}>0$ satisfying
\begin{equation}\label{eq:mixed-input-Mellin-lines}
  \Re(b_{j})\notin-2\mathbb{Z}_{\geq0},\quad \Re(b_{i}+b_{j})\notin-2\mathbb{Z}_{\geq0},\quad \Re(b_{k}-b_{4})\notin2\mathbb{Z},
\end{equation}
for $1\leq i,j\leq4$ and $1\leq k\leq3$. By evenness, the
continuation of the $t$ integral in \eqref{eq:mixed-input-joint-Mellin}
is half the real-line integral plus $\pi i$ times the signed residues at
$t=\pm i(b_{j}+2n)$ for $\Re(b_{j})+2n<0$ ($n\geq0$), with signs
$+$ for $+i(b_{j}+2n)$ and $-$ for $-i(b_{j}+2n)$.
Within the strip of $H_{R}$, this gives
\begin{equation}\label{eq:mixed-input-Mellin-bound}
 |\mathcal{M} f_{R}(u_{1},u_{2})|
 \ll_{\pi_{v},m}R^{O_{\pi_{v},m}(1)}
 (1+|\Im(u_{1})|+|\Im(u_{2})|)^{-m}.
\end{equation}
By Stirling and evenness, the exponential factors are bounded by constants times
\begin{align*}
 t\in\mathbb{R}_{\geq0}&:\quad e^{-\pi(|\Im(u_{1}-u_{2}/2)|-t)_{+}/2},\\
 t=i(b_{4}+2n)&:\quad 1,\qquad |\Re(t)|=|\Im(u_{1}-u_{2}/2)|,\\
 t=i(b_{j}+2n),\ j\leq3&:\quad
 e^{-\pi(|\Im(u_{1})|+|\Im(u_{1}-u_{2})|-|\Im(u_{2})|)/4}.
\end{align*}
In the $t$ integral and at the poles, \eqref{eq:R-admissible-test} gives Gaussian decay in
$|\Re(t)|/R$, and the Mellin transform of $\alpha$ in
\eqref{eq:mixed-input-joint-Mellin} decays rapidly in $|\Im(u_{2})|$.
For repeated poles, Cauchy's formula on small fixed-radius circles in
the strip gives the same bound \eqref{eq:mixed-input-Mellin-bound}.
For $H_{R}m_{z}$, the strip estimate in the proof of
Lemma~\ref{lem:hybrid-afe-multiplier} applies on these circles and the
real line, giving the additional factor in \eqref{eq:mixed-input-Sobolev}.

In lower triangular coordinates, right differentiation gives Laurent
polynomial coefficients for the upper-left block and $E_{13},E_{23}$.
The central identities of degrees $1,2,3$ give the same form for
$E_{33},E_{31},E_{32}$. For fixed $U$, there are finitely many
polynomials $P_{j,\ell}$ such that
\begin{equation}\label{eq:mixed-input-Lie-Mellin-shifts}
 (d\pi_{v}(U)W_{h_{R}})(\diag(ar,r,1))
 =\sum_{j,\ell}a^{j}r^{\ell}
 P_{j,\ell}(a\partial_{a},r\partial_{r})f_{R}(a,r).
\end{equation}
The Sobolev convergence passes the smooth-vector bound
\eqref{eq:mixed-fixed-Whittaker-bound}
\cite[\textsection~4.3, Proposition~2]{JS90} to the limit, so the Mellin
integral of the left side converges near $(\Re(u_{1}),\Re(u_{2}))=(-1/2,0)$.
There its transform is the sum of the products
$P_{j,\ell}(-u_{1}-j,-u_{2}-\ell)\mathcal{M} f_{R}(u_{1}+j,u_{2}+\ell)$
by continuation from their common region of convergence.
For $m>1+\max_{j,\ell}\deg P_{j,\ell}$, \eqref{eq:mixed-input-Mellin-bound}
and Mellin--Parseval bound the $L^{2}$ norms for the four weights
$a^{-1\pm2\epsilon_{1}}r^{\pm2\epsilon_{2}}$ by $O(R^{O_{\pi_{v},U}(1)})$,
with the same additional factor for $H_{R}m_{z}$.
The weights sum to at least $a^{-1}$.
Applying the bound to $\operatorname{Ad}(k)U$ for $k\in K_{v}'$ and to
$\operatorname{Ad}(w_{2}^{-1})U$ for the translated vector proves
\eqref{eq:mixed-input-Sobolev}.
\end{proof}

\begin{lemma}
\label{lem:mixed-real-standard-pole}
Let $W_{h}$ be the vector defined by \eqref{fc11.3} at
$\mathbf{s}=\mathbf{0}$ for a fixed self-dual unramified unitary $\pi_{v}$.
With $\mathcal{V}_{v}^{\circ}$ as in Lemma~\ref{lem:mixed-spherical-arch},
the quotient below extends holomorphically to $\lambda=\pm1/2$.
For some fixed $b$,
\begin{equation}\label{eq:mixed-real-standard-pole-bound}
 \left|\lim_{\lambda\to\pm1/2}
 \frac{\pi\mathcal{H}_{W_{h}}^{\vee}(\sigma_{\lambda,v},\mathbf{0},\mathbf{1})}
 {\mathcal{V}_{v}^{\circ}(\sigma_{\lambda,v})}\right|
 \ll_{\pi_{v}}\mathscr{L}_{b}(h).
\end{equation}
\end{lemma}

\begin{proof}
Write $(\nu_{1},\nu_{2},\nu_{3})=(\nu,0,-\nu)$ for the Langlands parameters
of $\pi_{v}$ and set $m=3$ if $\nu=0$ and $m=1$ otherwise.
Let $c_{\pi_{v}}(V)$ be the coefficient of $a(\log a)^{m-1}$ in the
expansion of $\int_{K_{v}'}\Phi_{V}(a,k)\,dk$ as $a\to0^{+}$, with $\Phi_{V}$
as in the proof of Lemma~\ref{lem:mixed-fixed-real-dual-decay}. We show that
\begin{equation}\label{eq:mixed-real-standard-pole-value}
 \lim_{\lambda\to1/2}
 \frac{\pi\mathcal{H}_{W_{h}}^{\vee}(\sigma_{\lambda,v},\mathbf{0},\mathbf{1})}
 {\mathcal{V}_{v}^{\circ}(\sigma_{\lambda,v})}
 =\frac{c_{\pi_{v}}(\pi_{v}(w_{2})W_{h})}{c_{\pi_{v}}(W_{\pi_{v}}^{\circ})}.
\end{equation}

Write the weight-$2n$ vector in the Kirillov model as
\begin{align*}
J_{n,\lambda}(a)=|a|^{1/2-\lambda}\widehat{v}_{n,\lambda}(a),\qquad v_{n,\lambda}(x)=(1+x^{2})^{-1/2-\lambda}e^{2in\arctan x}.
\end{align*}
For $\lambda\in i\mathbb{R}$, Parseval gives
$\|J_{n,\lambda}\|_{\mathrm{Kir}}^{2}=\pi$, and
$\overline{J_{n,\lambda}(a)}=J_{n,-\lambda}(a)$, so the $n$th summand of
\eqref{fc10.1} at $\mathbf{s}=\mathbf{0}$ and $\eta_{v}=\mathbf{1}$ is $\pi^{-1}\Psi_{v}(1/2,V,J_{n,\lambda})\Psi(J_{n,-\lambda};0)$. 

For fixed $\vartheta_{\pi_{v}}<\delta<1$ and $\lambda$ near $1/2$,
$\widehat v_{n,\lambda}(a)=\widehat v_{n,\lambda}(0)+O(|a|^{\delta})$.
The second period $\Psi(J_{n,-\lambda};0)$ is holomorphic near $\lambda=1/2$,
grows at most polynomially in $|n|$, and satisfies $\Psi(J_{0,-1/2};0)=2$.
Integration by parts in $k$ makes the sum in \eqref{fc10.1} locally
uniformly convergent after multiplication by $(1/2-\lambda)^{m}$.
Since $\widehat v_{n,1/2}(0)=\pi\mathbf{1}_{n=0}$, only $n=0$ contributes
to the $(1/2-\lambda)^{-m}$ coefficient, which in the first period is
$2\pi(-1)^{m-1}(m-1)!\,c_{\pi_{v}}(V)$.
Taking the ratio with $V=W_{\pi_{v}}^{\circ}$ gives
\eqref{eq:mixed-real-standard-pole-value} by
Lemma~\ref{lem:mixed-spherical-arch}. Since
$\sigma_{-\lambda,v}\simeq\sigma_{\lambda,v}$, the limit at $\lambda=-1/2$
equals the limit in \eqref{eq:mixed-real-standard-pole-value}.

By \cite[\textsection~4.4, Proposition~3]{JS90}, the coefficient of
$a(\log a)^{m-1}$ in $V(\diag(az,z,1)k)$ is $O_{V}(|z|^{\varepsilon})$
as $z\to0$ for some $\varepsilon>0$ and $O_{V,N}(|z|^{-N})$ as
$|z|\to\infty$ for every $N>0$, uniformly in $k\in K_{v}'$.
Its integral against $|z|^{2s_{1}'}d^{\times}z$ is holomorphic at
$s_{1}'=0$, so the leading pole coefficient is unchanged by this weight.
On $\lambda=1/2-s_{2}'$, Fourier inversion for $\Re(s_{2}')>0$ gives
$\Psi(J_{n,-\lambda};s_{2}')=v_{n,-\lambda}(0)=1$.
Division by $L_{v}(1/2+s_{1}',\pi_{v}\times\sigma_{\lambda,v})$
cancels the first period's pole at $s_{1}'+s_{2}'=0$.
As $s_{2}'\to0$ on this line, only $n=0$ contributes to the leading
coefficient. Its second period cancels in
\eqref{eq:mixed-real-standard-pole-value}.
At $\lambda=1/2$, dividing the two $(s_{1}')^{-m}$ coefficients gives
\begin{align*}
 \left.\frac{\Psi_{v}(1/2+s_{1}',V,J_{0,1/2})}
 {L_{v}(1/2+s_{1}',\pi_{v}\times\sigma_{1/2,v})}\right|_{s_{1}'=0}
 =\frac{2\pi(-1)^{m-1}(m-1)!\,c_{\pi_{v}}(V)}
 {2^{m}L_{v}(1,\pi_{v})\prod_{\nu_{j}\ne0}\Gamma_{\mathbb{R}}(\nu_{j})}.
\end{align*}
Continuity of the normalized period \cite[Theorem~2.3(i)]{Jac09}
bounds $|c_{\pi_{v}}(V)|$ by a fixed Sobolev seminorm.
Combining this with the adjoint action in the proof of
Lemma~\ref{lem:mixed-fixed-real-dual-decay} gives
$|c_{\pi_{v}}(\pi_{v}(w_{2}g)W_{v}^{\circ})|\ll(1+\|g\|+\|g^{-1}\|)^{b}$.
Integrating against $|h(g^{-1})|\,dg$ in \eqref{fc11.3} gives
$|c_{\pi_{v}}(\pi_{v}(w_{2})W_{h})|\ll\mathscr{L}_{b}(h)$.
\end{proof}

\begin{cor}\label{cor:mixed-polynomial-input-decay}
Let $W_{h_{v,R,z}}$ be the Whittaker vector of
Lemma~\ref{lem:mixed-input-Sobolev} for $H_{R}m_{z}$, with $\Re(z)=c$.
For fixed $\ell$ and a sufficiently wide holomorphy strip for $H_{R}$,
the following bound holds for every generic $\sigma_{v}$ in the unitary
automorphic spectrum with trivial central character:
\begin{equation}\label{eq:mixed-polynomial-input-decay}
|\mathcal{H}_{W_{h_{v,R,z}}}^{\vee}(\sigma_{v},\mathbf{0},\mathbf{1})|
 \ll_{\pi_{v},\ell,c}\frac{R^{c+O_{\pi_{v},\ell}(1)}(1+|\Im(z)|)^{O(1)}e^{\pi|\Im(z)|/2}}{C(\sigma_{v})^{\ell/2}}.
\end{equation}
The same bound without $C(\sigma_{v})^{-\ell/2}$ holds for both integrals
in \eqref{7.5} at $\mathbf{s}=\mathbf{0}$ and the limits in
\eqref{eq:mixed-real-standard-pole-bound} for $h=h_{v,R,z}$.
\end{cor}

\begin{proof}
The torus expansion \cite[\textsection~4.3, Proposition~2]{JS90}
gives \eqref{eq:mixed-fixed-Whittaker-bound} for $d\pi_{v}(U)V$ with a
fixed Sobolev seminorm of $V$ as the constant. Integrating in the central
variable of $\Phi_{V}$ bounds the norms in \eqref{eq:mixed-full-K-sobolev-bound}
for $V=(1+\Delta_{v}^{2})^{m}\pi_{v}(w_{2})W_{h_{v,R,z}}$.
The range of $\eta$ in that estimate is nonempty by
Lemma~\ref{lem:mixed-spherical-arch}. With \eqref{eq:mixed-input-Sobolev}
and the Casimir factor $(1+|\lambda_{\sigma_{v}}|^{2})^{-m}$,
$4m\geq\ell$ gives \eqref{eq:mixed-polynomial-input-decay}.

Integrating \eqref{eq:mixed-fixed-Whittaker-bound} in both diagonal
variables bounds the integrals in \eqref{7.5} by fixed Sobolev seminorms.
The proof of Lemma~\ref{lem:mixed-real-standard-pole} bounds the limits
by fixed Sobolev seminorms as well. Equation \eqref{eq:mixed-input-Sobolev}
bounds both sets of seminorms.
\end{proof}

\section{Further Fourier Reciprocity Formulas}\label{app:further-reciprocity}
\subsection{Fourier Reciprocity Formula of Type \RNum{3}}
Let $\varphi\in\pi$, where $\pi$ is generic automorphic with central
character $\omega$. Put $\varphi^{\iota}(g):=\varphi(g^{\iota})$ in
the contragredient realization $\pi^{\iota}$. With respect to
$\overline{\psi}$, its Whittaker coefficient is
\begin{equation}\label{e2.25}
W_{\varphi^{\iota}}(g)=W_{\varphi}(\widetilde{w} g^{\iota}).
\end{equation}
Write $T\varphi:=\pi^{\iota}(\widetilde{w})\varphi^{\iota}$ and put
\begin{align*}
\mathcal{I}_{\mathrm{const}}(\varphi):=
\int_{F\backslash\mathbb{A}_{F}}
\bigl[\varphi(n_{23}(c))-\varphi(n_{12}(c))\bigr]\,dc+\mathcal{F}_{\mathrm{const}}^{\dag}\varphi(I_{3})
-\mathcal{F}_{\mathrm{const}}\varphi(I_{3}).
\end{align*}

\begin{prop}\label{prop2.4}
With $\mathcal{I}_{\mathrm{degen}}$ as in Proposition~\ref{prop2.5}, put
$\mathcal{I}_{\mathrm{degen}}^{\dagger}(\varphi):=\mathcal{I}_{\mathrm{degen}}(T\varphi)$.
Then
\begin{align*}
\mathcal{I}(\mathcal{F}_{\mathrm{gen}}\varphi)
=\mathcal{I}(\mathcal{F}_{\mathrm{gen}}T\varphi)
+\mathcal{I}_{\mathrm{degen}}^{\dag}(\varphi)
-\mathcal{I}_{\mathrm{degen}}(\varphi)
+\mathcal{I}_{\mathrm{const}}(\varphi).
\end{align*}
If the following integral converges absolutely, then
$\mathcal{I}_{\mathrm{degen}}^{\dag}(\varphi)$ equals
\begin{align*}
\sum_{\alpha,\gamma\in F^{\times}}\int_{\mathbb{A}_{F}}
W_{\varphi}^{\mathrm{degen},\dag}
\bigl(w_{2}n_{23}(c)\diag(\alpha\gamma,\gamma,1)\bigr)
\overline{\psi(c)}\,dc.
\end{align*}
\end{prop}

\begin{proof}
By \eqref{e2.31} and $T\varphi(I_{3})=\varphi(I_{3})$,
$\mathcal{I}(\varphi)-\mathcal{I}(T\varphi)
=\mathcal{I}_{0}(T\varphi)-\mathcal{I}_{0}(\varphi)$. The identity
$\widetilde{w} n_{12}(c)^{\iota}\widetilde{w}=n_{23}(-c)$ turns
$\mathcal{I}_{0}(T\varphi)$ into $\int_{F\backslash\mathbb{A}_{F}}
\varphi(n_{23}(c))\,dc$. Expand $\mathcal{I}(\varphi)$ and
$\mathcal{I}(T\varphi)$ by \eqref{2.1}. For the two $\mathcal{F}_{\mathrm{const}}$
terms, \eqref{e2.31} subtracts the same constant term along $N$ at $I_{3}$,
leaving the difference of values
$\mathcal{F}_{\mathrm{const}}^{\dag}\varphi(I_{3})
-\mathcal{F}_{\mathrm{const}}\varphi(I_{3})$.
The sum in $\mathcal{I}$ is unchanged by replacing $\psi$ by
$\overline{\psi}$. Apply \eqref{fc2.37} to $T\varphi$ with that
character. Inverse transpose, rational central invariance and
$(c,\alpha,\gamma)\mapsto(-c,\gamma,\alpha)$ give the stated integral.
\end{proof}

If $\pi$ is cuspidal, the degenerate terms and both
parabolic constant coefficients vanish, but
$\mathcal{I}_{\mathrm{const}}(\varphi)$ need not vanish.

\subsection{Spectral Reciprocity Formula of Type \RNum{3}}\label{sec:type-III-spectral}
\begin{thm}
Suppose $\pi$ is convergent in the sense of Definition~\ref{defn3.2},
and let $\varphi\in\pi$. Put
$\mathbf{s}=(s_{1},s_{2})$ and
$\mathbf{s}^{*}:=(\frac{s_{1}+s_{2}}{2},\frac{3s_{1}-s_{2}}{2})$.
If $\Re(s_{2}-s_{1})\ggg1$ and $\Re(3s_{1}-s_{2})\ggg1$, with constants
allowed to depend on $\varphi$, then
\begin{align*}
I_{\mathrm{spec}}(\mathbf{s},\varphi,\omega',\eta)
=I_{\mathrm{spec}}(\mathbf{s}^{*},T\varphi,\omega\eta,\omega\omega')
+I_{\mathrm{corr}}(\mathbf{s},\varphi,\omega',\eta),
\end{align*}
where
\begin{multline*}
I_{\mathrm{corr}}(\mathbf{s},\varphi,\omega',\eta)
:=\int_{(F^{\times}\backslash\mathbb{A}_{F}^{\times})^{2}}
(\mathcal{I}_{\mathrm{degen}}^{\dag}-\mathcal{I}_{\mathrm{degen}}
+\mathcal{I}_{\mathrm{const}})
\bigl(\pi(\diag(yz,z,1))\varphi\bigr)\\
\omega'(z)\eta(y)|z|^{2s_{1}}|y|^{s_{1}+s_{2}}\,d^{\times} z\,d^{\times} y.
\end{multline*}
The integral defining $I_{\mathrm{corr}}$ converges absolutely in this region.
\end{thm}

\begin{proof}
Choose representatives of $c\bmod F$ and of the norm-one idele
classes in compact sets. Vary $|y|$ and $|z|$ by scaling at infinity.
For $g=n_{12}(c)\diag(yz,z,1)$ as in \eqref{e3.13}, the term indexed by
$\delta$ in the generic sum in \eqref{2.1} has $\mathrm{GL}_{2}$ block
$\delta\left(\begin{smallmatrix}yz&cz\\0&z\end{smallmatrix}\right)$.
Write its local Iwasawa diagonal as $\diag(a_{v}b_{v},b_{v})$.
The bottom row is $z(y\delta_{21},c\delta_{21}+\delta_{22})$.
Its norm is comparable to $|b_{v}|_{v}$, and the determinant gives
$|a_{v}b_{v}^{2}|_{v}=|yz^{2}\det\delta|_{v}$.

At finite places the Whittaker function vanishes unless
$|a_{v}|_{v},|b_{v}|_{v}\ll1$, with upper bound $1$ outside a fixed
finite set. Thus $\delta_{21},\delta_{22}$ and $\det\delta$ lie in
fixed fractional ideals. They determine the coset of $\delta$ in
$N'(F)\backslash G'(F)$. At infinity the same bottom row and determinant give
\begin{align*}
\max\{|\delta_{21}|_{v},|\delta_{22}|_{v}\}
\ll |b_{v}/z|_{v}(1+|y|_{v}^{-1}),\qquad
|\det\delta|_{v}=|a_{v}b_{v}^{2}/(yz^{2})|_{v}.
\end{align*}
By the product formula, the inverse powers of $|a_{v}|_{v},|b_{v}|_{v}$
in the Whittaker estimates are bounded by fixed powers of the
archimedean sizes of the three entries and of $|y|^{\pm1},|yz^{2}|^{\pm1}$.
Since $|\det\delta|_{\infty}\gg1$, the determinant identity yields
$\prod_{v\mid\infty}|a_{v}b_{v}^{2}|_{v}\gg|yz^{2}|$.
The absolute Whittaker sum over $\delta$ is therefore
$O_{\varphi,N}((|y|+|y|^{-1})^{B}|yz^{2}|^{-B}(1+|yz^{2}|)^{-N})$
with $B$ fixed, also for fixed right derivatives of $\varphi$.
By $n_{12}(c)\diag(yz,z,1)=\diag(yz,z,1)n_{12}(c/y)$, integration
by parts in $c$ and summation over the nonzero frequencies in the
fixed fractional ideal contribute a factor
$O_{M}(|y|^{-1}(1+|y|)^{-M})$.
After including the modulus weight in \eqref{e3.13}, the absolute
sum of the nonzero $c$-Fourier coefficients is
\begin{align*}
\ll_{\varphi,M,N}
|y|^{\Re(s_{2})-B}|yz^{2}|^{\Re(s_{1})-B}
(1+|y|)^{-M}(1+|yz^{2}|)^{-N},
\end{align*}
after enlarging $B$ and choosing the decay orders as needed.
Since $d\log|z|=2^{-1}d\log|yz^{2}|$ at fixed $y$, this is
integrable for $\Re(s_{1}),\Re(s_{2})>B$ and large $M,N$.

Splitting the $y$ integral of \eqref{e3.13} into idele classes
sums the nonzero coefficients as $\mathcal{I}(\mathcal{F}_{\mathrm{gen}}
(\pi(\diag(yz,z,1))\varphi))$. For the second term of
Proposition~\ref{prop2.4},
\begin{align*}
T\bigl(\pi(\diag(yz,z,1))\varphi\bigr)
=\omega(yz)\pi^{\iota}(\diag(yz,y,1))T\varphi.
\end{align*}
Interchanging the two idele class variables gives the factor
$(\omega\eta)(z)(\omega\omega')(y)$ and the powers
$2s_{1}^{*}=s_{1}+s_{2}$, $s_{1}^{*}+s_{2}^{*}=2s_{1}$.
Both pairs of Mellin parameters have sufficiently large real parts,
including the shifts in \eqref{eq:moment-character-twists}.
For both vectors the Whittaker bound also holds with $(1+|yz|)^{-N}$:
$\prod_{v\mid\infty}|b_{v}|_{v}\gg|yz|$ if $\delta_{21}\ne0$, and
$\prod_{v\mid\infty}|a_{v}b_{v}|_{v}\gg|yz|$ otherwise, since the
nonzero entry $\delta_{21}$ or $\delta_{11}$ lies in a fixed fractional ideal.
Changing from $|z|$ to $|yz|$ bounds the central Mellin integrals by
$O(|y|^{2B-\Re(s_{1})})$ uniformly on a Siegel set, also after fixed
right derivatives. This gives $L^{2}$ for both integrals.
Proposition~\ref{prop2.4} expresses $I_{\mathrm{corr}}$ as their absolutely
convergent Fourier integrals in \eqref{e3.13} with opposite signs.
Lemma~\ref{lem4.3} gives rapid decay of the first periods, so the
argument of \eqref{3.9} identifies both spectral terms. The rational
diagonals $\diag(1,-1,1)$ and $a(-1)$ reverse the Whittaker characters
in the transformed periods.
\end{proof}

\subsection{Spectral Reciprocity Formula of Type \RNum{2}: Higher Ranks}
Let $n\geq2$ and let $\pi_{1},\pi_{2}$ be generic automorphic
representations of $\mathrm{GL}_{n+1}/F$, $\mathrm{GL}_{n-1}/F$,
with central characters $\omega_{1},\omega_{2}$. Take
$\varphi_{i}\in\pi_{i}$ and a Hecke character $\omega'$.

Let $U_{m}\subset\mathrm{GL}_{m}$ be upper unipotent, set
$[U_{m}]=U_{m}(F)\backslash U_{m}(\mathbb{A}_{F})$ with volume $1$,
and write $\theta_{m}(u)=\psi(\sum_{j=1}^{m-1}u_{j,j+1})$. Define
\begin{align*}
\mathcal{I}_{n}(\varphi_{1}):=
\sum_{\delta\in U_{n-1}(F)\backslash \mathrm{GL}_{n-1}(F)}
\int_{[U_{n}]}\varphi_{1}\bigl(\diag(u,1)\diag(\delta,I_{2})\bigr)
\overline{\theta_{n}(u)}\,du.
\end{align*}
At $I_{n}$, the Fourier expansion of the restriction of $\varphi_{1}$
to $\mathrm{GL}_{n}$ has generic part $\mathcal{I}_{n}(\varphi_{1})$
\cite[Proposition~3.1]{Yan25}.

For $\mathbf{s}=(s_{1},s_{2})\in\mathbb{C}^{2}$, define, whenever the integral
converges absolutely,
\begin{multline}\label{3.35}
I_{n}(\mathbf{s},\varphi_{1},\omega',\varphi_{2}):=\int_{[\mathrm{GL}_{n-1}]}\int_{F^{\times}\backslash\mathbb{A}_{F}^{\times}}\mathcal{I}_{n}(\pi_{1}(\diag(g',I_{2}))\pi_{1}(\diag(zI_{n},1))\varphi_{1})\\
\omega'(z)|z|^{ns_{1}}d^{\times}z\varphi_{2}(g')|\det g'|^{s_{1}+s_{2}}dg'.
\end{multline}

We define $I_{\text{non-gen}}$ by replacing $\mathcal{I}_{n}$ in
\eqref{3.35} with the sum of the non-generic Fourier terms.

\begin{prop}
Let $w_{n}$ be the permutation matrix of $(n,n+1)$ and put
$\mathbf{s}^{\vee}=(\frac{(n-1)s_{2}-s_{1}}{n},\frac{(n+1)s_{1}+s_{2}}{n})$.
Whenever the four integrals converge absolutely,
\begin{multline}\label{f3.37}
I_{n}(\mathbf{s},\varphi_{1},\omega',\varphi_{2})=I_{n}(\mathbf{s}^{\vee},\pi_{1}(w_{n})\varphi_{1},\omega_{1}^{-1}\omega_{2}\omega'^{-1},\varphi_{2})\\
-I_{\text{non-gen}}(\mathbf{s},\varphi_{1},\omega',\varphi_{2})
+I_{\text{non-gen}}(\mathbf{s}^{\vee},\pi_{1}(w_{n})\varphi_{1},\omega_{1}^{-1}\omega_{2}\omega'^{-1},\varphi_{2}),
\end{multline}
\end{prop}

\begin{proof}
The full Fourier expansion of $\varphi_{1}(\diag(zg',z,1))$
in \eqref{3.35} gives $I_n+I_{\text{non-gen}}$.
Substituting $(z,g')\mapsto(z^{-1},zg')$ and using automorphy with
$w_{n}\diag(g',z^{-1},1)=z^{-1}\diag(zg',z,1)w_{n}$ replaces
$\varphi_{1}$ by $\pi_{1}(w_{n})\varphi_{1}$. The new character is
$(\omega_{1}^{-1}\omega_{2}\omega'^{-1})(z)$ and the new exponents are
$n s_{1}^{\vee}=(n-1)s_{2}-s_{1}$ and
$s_{1}^{\vee}+s_{2}^{\vee}=s_{1}+s_{2}$. Splitting the Fourier sum
again gives \eqref{f3.37}.
\end{proof}

When $\varphi_{2}$ is cuspidal, unfolding its Whittaker expansion leaves
only the $N_{n}$-constant term among the non-generic terms. In
$I_{\text{non-gen}}$, replace the Fourier sum in \eqref{3.35} by
$\int_{[N_{n}]}\varphi_{1}(\diag(u,1)\diag(zg',z,1))\,du$,
where $N_{n}$ is the unipotent radical of the $(n-1,1)$ parabolic in
$\mathrm{GL}_{n}$ and $[N_{n}]=N_{n}(F)\backslash N_{n}(\mathbb{A}_{F})$.

If $\pi_{1}$ is also cuspidal, the $z$ Mellin integral of
$\varphi_{1}(\diag(zg,1))$ with the weight in \eqref{3.35}, times
$|\det g|^{s_{1}}$, has central character $\omega'^{-1}$.
Its spectral coefficients are the $\mathrm{GL}_{n+1}\times\mathrm{GL}_{n}$
periods. Subtracting its $N_{n}$-constant term and integrating against
$\varphi_{2}(g')|\det g'|^{s_{2}}$ gives the
$\mathrm{GL}_{n}\times\mathrm{GL}_{n-1}$ periods.
Thus both $I_n$-terms in \eqref{f3.37} are spectral averages of
$\mathrm{GL}_{n+1}\times\mathrm{GL}_{n}$ and
$\mathrm{GL}_{n}\times\mathrm{GL}_{n-1}$ $L$-values for
$\Re(s_{1}),\Re((n-1)s_{2}-s_{1})\ggg1$.
For $n=2$, $\varphi_{2}=\eta$ and cuspidal $\varphi_{1}$, unfolding
\eqref{3.35} gives \eqref{e3.13}, whose spectral expansion is \eqref{3.25}.

For cuspidal $\pi_{1},\pi_{2}$ and $\omega_{1}=\omega_{2}=\omega'=1$,
the period identity for $I_{n}+I_{\text{non-gen}}$ is the one used by
Miao and Jana--Nunes, with parameters $(s_{1}+1/2,s_{2}+1/2)$
\cite[Proposition~3.5]{Mia21}, \cite[Proposition~5.1]{JN24}.

\clearpage
\section*{List of Notation and Transforms}

The entries below explain the main symbols and transforms used in the
paper and point to their definitions.

\begingroup
\newcommand{\notationgroup}[1]{%
 \par\addvspace{.8\baselineskip}%
 \noindent\textbf{#1}\par\nobreak\vspace{.38\baselineskip}\nobreak
}
\newcommand{\notationentry}[2]{%
 \par\noindent
 \makebox[.29\linewidth][l]{#1}%
 \parbox[t]{.71\linewidth}{\raggedright #2\strut}%
 \par\addvspace{.37\baselineskip}%
}

\notationgroup{Fourier expansions and global reciprocity}

\notationentry{$n(b),\ n_{ij}(u)$}
{$n(b)$ is the $\mathrm{GL}_{2}$ unipotent matrix with upper-right entry $b$
defined in \textsection\ref{sec:matrix-conventions}.
For $i\ne j$, $n_{ij}(u)=I_{3}+uE_{ij}$, with $E_{ij}$ the matrix unit
in row $i$, column $j$. Thus $n_{21}(y)$ has $y$ just below the first
diagonal entry. See \textsection\ref{sec1.3}.}

\notationentry{$w',\ w_{1},\ w_{2},\ \widetilde{w}$}
{The Weyl element $w'$ swaps the two coordinates of $G'$.
In $G$, the elements $w_{1}$ and $w_{2}$ swap coordinates $(1,2)$
and $(2,3)$, and $\widetilde{w}=w_{1}w_{2}w_{1}$ is the long Weyl
element. See \textsection\ref{sec:matrix-conventions}.}

\notationentry{$W_{\varphi},\ W_{\varphi}^{\mathrm{degen}},\
W_{\varphi}^{\mathrm{degen},\dagger}$}
{The $\mathrm{GL}_{3}$ Whittaker integrals in and after \eqref{w2.5}
use $\overline{\psi(a+c)}$, $\overline{\psi(c)}$, and
$\overline{\psi(a)}$, respectively, where $a$ and $c$ are the
$(1,2)$ and $(2,3)$ entries. The $\mathrm{GL}_{2}$ integral $W_{\phi}$
uses $\psi$ in \textsection\ref{sec3.2}.}

\notationentry{$g^{\iota},\ \varphi^{\iota},\ \widetilde{W}_{v}$}
{The inverse transpose $g^{\iota}$ is defined in
\textsection\ref{sec:matrix-conventions}. The contragredient form is
$\varphi^{\iota}(g)=\varphi(g^{\iota})$ (before \eqref{2.2}), and
$\widetilde{W}_{v}(g)=W_{v}(\widetilde{w}g^{\iota})$ is the local
Whittaker vector of \textsection\ref{sec3.2}.}

\notationentry{$\mathcal{F}_{\mathrm{gen}}\varphi,\ \mathcal{F}_{\mathrm{gen}}^{\dagger}\varphi$}
{Generic parts of the Fourier expansions for the upper-left and
lower-right copies of $G'$, defined after \eqref{2.1} and \eqref{2.2}.
The constant and degenerate parts use the same $\dagger$ convention.}

\notationentry{$\mathcal{W},\ \mathcal{I},\ \mathcal{I}_{0}$}
{In type~\RNum{1}, $\mathcal{W}$ is the $\psi$-Fourier coefficient
along $n_{12}$ at $I_{3}$, defined in \eqref{f2.9}. The type~\RNum{2} sum
$\mathcal{I}$ of all nonzero $n_{12}$-Fourier modes at $I_{3}$ is defined
before Proposition~\ref{prop2.5}. The constant term $\mathcal{I}_{0}$
is defined before Lemma~\ref{lem2.4}.}

\notationentry{$\Psi,\ \Psi_{v},\ \widetilde{\Psi}$}
{Global periods, their local integrals, and their global meromorphic
continuations. The arguments distinguish $\mathrm{GL}_{3}\times\mathrm{GL}_{2}$
from $\mathrm{GL}_{2}\times\mathrm{GL}_{1}$. See \textsection\ref{sec3.2}.}

\notationentry{$J_{\mathrm{spec}},\ J_{\mathrm{dual}},\ J_{\mathrm{sing}}$}
{The type~\RNum{1} spectral moment \eqref{s3.3}, its
$\mathrm{GL}_{1}$ dual, and the $\beta=0$ term before Mellin inversion.
The last two are defined after \eqref{e3.3}.}

\notationentry{$J_{\mathrm{degen}},\ J_{\mathrm{degen}}^{\dagger}$}
{The two type~\RNum{1} degenerate Whittaker integrals defined after
\eqref{e3.3}. For Eisenstein data, their $\heartsuit$ versions use the
constant term in the inducing parameters \eqref{eq:inducing-constant-term}.}

\notationentry{$I_{\mathrm{spec}},\ I_{\mathrm{gen}},\ I_{\mathrm{degen}}$}
{In type~\RNum{2}, the two-period spectral moment \eqref{3.25}
and the generic and degenerate Whittaker integrals following \eqref{e3.8}.}

\notationentry{$\omega,\ \omega',\ \eta$}
{The central characters are $\omega$ for $\pi$ and $\omega'$ for the
original $\mathrm{GL}_{2}$ spectrum. The Hecke character $\eta$
twists the standard period (before \eqref{3.25}).
In Lemma~\ref{lem:maximal-Steinberg-positive-residue},
$\eta_{v}$ instead twists the Steinberg representation.}

\notationentry{$\mathbf{s}^{\vee},\ \pi(w_{2})\varphi$}
{In type~\RNum{2}, $\mathbf{s}^{\vee}=(s_{1}',s_{2}')
=((s_{2}-s_{1})/2,(3s_{1}+s_{2})/2)$. The dual side right-translates
$\varphi$ by $w_{2}$ and replaces $\omega'$ by
$\overline{\omega}\,\overline{\omega}'\eta$, leaving $\eta$ unchanged.
See Theorem~\ref{thmii}.}

\notationentry{$\mathbf{s}^{*},\ T\varphi,\ I_{\mathrm{corr}}$}
{In \textsection\ref{sec:type-III-spectral},
$\mathbf{s}^{*}=((s_{1}+s_{2})/2,(3s_{1}-s_{2})/2)$ and
$I_{\mathrm{corr}}$ is the correction from the degenerate and constant terms.
Type~\RNum{3} replaces $(\omega',\eta)$ by $(\omega\eta,\omega\omega')$
and uses $T\varphi=\pi^{\iota}(\widetilde{w})\varphi^{\iota}$,
defined just after \eqref{e2.25}.}

\notationentry{$\mathbf{s}^{\vee}$ in higher rank}
{The map for \eqref{f3.37} is defined just above it by
$(((n-1)s_{2}-s_{1})/n,((n+1)s_{1}+s_{2})/n)$.
It reduces to the type~\RNum{2} map at $n=2$.}

\notationgroup{Local weights and transforms}

\notationentry{$\mathcal{H}_{W_{v}}(\sigma_{v},s)$}
{Original type~\RNum{1} weight \eqref{e9.1}. Sum the local
$\mathrm{GL}_{3}\times\mathrm{GL}_{2}$ period times
$\overline{W_{v}'(I_{2})}$ over the Whittaker basis of $\sigma_{v}$.}

\notationentry{$\mathcal{H}_{W_{v}}(\xi_{v},\lambda;\omega_{v}',s)$}
{The double Mellin integral \eqref{eqq9.2} is the local weight for
$J_{\mathrm{dual}}$, indexed by $(\xi_{v},\lambda)$. See
\textsection\ref{sec:type-I-dual-expansion} for the change of Whittaker argument.}

\notationentry{$\mathcal{H}_{W_{v}}^{\sharp}$}
{Finite-place dual weight of type~\RNum{1} after division by
$L_{v}(1/2+2s-\lambda,\widetilde{\pi}_{v}\times
\overline{\xi}_{v}\omega_{v}\omega_{v}')$ and
$L_{v}(1/2+\lambda,\xi_{v})$.
See \textsection\ref{def:normalized-finite-dual-weight}.}

\notationentry{$W_{v}(\cdot;\mathbf{b})$}
{Auxiliary Whittaker vector \eqref{eq9.11-real} with
$\mathbf{b}=(b_{1},b_{2})$. Equation \eqref{fc9.13} recovers the
type~\RNum{1} test vector by integrating its right translates.}

\notationentry{$\mathcal{P}_{v}^{\sharp},\ \mathcal{P}_{v}^{\heartsuit}$}
{The local degenerate integral divided by its Euler factor in
\eqref{eq:local-w2-degenerate}, and its continuation defined before
\eqref{f5.15}. Add $\dagger$ for the opposite term.
These are not the Eisenstein factors $\mathscr{P}_{\chi,v}$ defined
after \eqref{eq:maximal-spectral-density-factors}.}

\notationentry{$\mathcal{H}_{W_{v}}(\sigma_{v},\mathbf{s},\eta_{v})$}
{Original type~\RNum{2} basis sum \eqref{fc10.1}. It pairs the
$\mathrm{GL}_{3}\times\mathrm{GL}_{2}$ period with the conjugate of
the $\mathrm{GL}_{2}\times\mathrm{GL}_{1}$ period at
$(1/2+\overline{s_{2}},\overline{\eta_{v}})$.}

\notationentry{$\mathcal{H}_{W_{v}}^{\vee}(\sigma_{v},\mathbf{s}^{\vee},\eta_{v})$}
{Dual type~\RNum{2} weight. In \eqref{fc10.1} replace $W_{v}$ by
$g\mapsto W_{v}(gw_{2})$ and $\mathbf{s}$ by $\mathbf{s}^{\vee}$,
leaving $\sigma_{v}$ and $\eta_{v}$ unchanged. See \eqref{fc10.2}.}

\notationentry{$\mathcal{S}h_{v},\ \mathcal{T}_{v},\ \mathcal{T}$}
{The spherical transform $\mathcal{T}_{v}=\mathcal{S}h_{v}$, evaluated
at $\nu_{\sigma_{v}}$, is the eigenvalue of $\sigma_{v}(h_{v})$ on
its spherical line. Theorem~\ref{thmC} uses
$\mathcal{T}=(\mathcal{T}_{v})_{v\mid\infty}$. See \eqref{eq10.4}.}

\notationentry{$H_{R},\ h_{v,R},\ W_{v}$}
{For the fifth moment, $H_{R}$ is the example after
\eqref{eq:R-admissible-test} with fourth powers replaced by fifth powers.
The function $h_{v,R}$ is defined by its
spherical transform in \textsection\ref{def:polynomial-spherical-test}.
Taking $h_{v}=h_{v,R}$ in \eqref{fc11.3} gives $W_{v}$ with
$\pi\mathcal{H}_{W_{v}}(\sigma_{t},\mathbf{0},\mathbf{1})=H_{R}(t)$
by \eqref{eq:G-as-L-ratio}.}

\notationentry{$\nu_{\sigma_{v}},\ \sigma_{t},\ G_{\pi_{v}}(t)$}
{The exponents of a real spherical $\sigma_{v}$ are
$\pm\nu_{\sigma_{v}}$ (before \eqref{eq10.4}). The representation
$\sigma_{t}$ has exponents $\pm it$ (before \eqref{eq:G-as-L-ratio}).
The local $L$-ratio $G_{\pi_{v}}(t)$, including $\beta_{v}$, is defined
before Lemma~\ref{lem:mixed-afe-realization}. By
\eqref{eq:real-weight-evaluated}, the normalized spherical weight is
$G_{\pi_{v}}(t)\mathcal{S}h_{v}(it)$.}

\notationentry{$W_{n},\ \mathcal{N}_{\sigma_{v}},\ \mathfrak{B}_{\sigma_{v}}$}
{For real $\sigma_{v}$ of trivial central character, the vectors
$W_{n}$ of rotation weights $2n\in\mathcal{N}_{\sigma_{v}}$
form the orthonormal basis $\mathfrak{B}_{\sigma_{v}}$.
See \textsection\ref{sec10.2.1}.}

\notationentry{$G_{\boldsymbol{\nu}},\ \beta_{v},\ \mathcal{S}_{\sigma_{v}}$}
{$G_{\boldsymbol{\nu}}$ is the Gamma ratio \eqref{eq:real-double-Mellin-Gamma},
and $\beta_{v}$ is the real Whittaker scale in
\eqref{eq:real-jacquet-double-mellin}.
The rotation-weight sum $\mathcal{S}_{\sigma_{v}}$ of
Lemma~\ref{lemm11.3} enters the dual weight in Theorem~\ref{prop10.2}.}

\notationentry{$\mathscr{T}H_{R}$}
{The real weight $\pi\mathcal{H}_{W_{v}}^{\vee}(\sigma_{v},\mathbf{0},\mathbf{1})$
with $W_{v}$ from $h_{v,R}$ and $\sigma_{v}$ of trivial central character.
The notation $(\mathscr{T}_{\epsilon}H_{R})(u)$ uses the inducing
characters $\sgn^{\epsilon}|\cdot|^{\pm iu}$. For the complementary
series $u=i\eta$, where $\eta$ is a real exponent, not a Hecke
character. For lowest weight $2l$ use
$(\mathscr{T}H_{R})^{\mathrm{disc}}(2l)$. See
Definition~\ref{def:real-dual-transform}.}

\notationentry{$f_{v},\ \widehat{f}_{v},\ f_{v,R}$}
{The integral $f_{v}$ of $h_{v}$ is defined in \eqref{eq10.17}.
Its Fourier transform $\widehat{f}_{v}$ is defined in \eqref{c11.20},
with value at zero given by \eqref{11.22}.
Immediately after Definition~\ref{def:real-dual-transform},
$f_{v,R}$ sets $h_{v}=h_{v,R}$ and $s_{2}=0$.}

\notationentry{$\mathcal{K}_{iu,\epsilon}$}
{The cosine kernel \eqref{eq:principal-cosine-kernel} replaces the
sum over $K_{v}'$-types in the principal-series transform.}

\notationentry{$m_{z}(t),\ S_{0},\ S_{z}$}
{The approximate functional equation replaces the tests by
$H_{R}m_{z}$ and $S_{z}=S_{0}m_{z}$, with $S_{0}$ for the bounded
range. The Gamma ratio $m_{z}$ is local, while $e^{z^{2}}$ is inserted
once globally. See
Lemmas~\ref{lem:hybrid-afe-multiplier} and~\ref{lem:mixed-afe-realization}.}

\notationentry{$W_{v}^{[m]},\ W_{v}^{\circ},\ W_{v}'^{\circ}$}
{$W_{v}^{\circ}(\diag(\varpi_{v}^{-2d_{v}},\varpi_{v}^{-d_{v}},1))=1$
in \textsection\ref{sec12}, whereas at real $v$ the scale of
$W_{v}^{\circ}$ is $\beta_{v}$ in \eqref{eq:real-jacquet-double-mellin},
with $\beta_{v}=\mathbf{L}(\pi_{v})^{-1}$ for Theorem~\ref{thmC}
(\textsection\ref{sec11.4}).
Its finite translates give $W_{v}^{[m]}$ in \eqref{c12.2}.
The primed $W_{v}'^{\circ}$ is the unit finite newvector
(\textsection\ref{sec:finite-newvector}), but at real $v$ the proof of
Lemma~\ref{lem:real-rs-evaluation} takes
$W_{v}'^{\circ}(\diag(a,1))=|a|^{1/2}K_{\nu}(2\pi|a|)$.}

\notationentry{$W_{m,v}',\ \xi_{\sigma_{v}},\ \alpha_{\sigma_{v}}$}
{The finite-place oldform basis vector \eqref{c12.4} has coefficients
$\xi_{\sigma_{v}}$ given in \eqref{eq:oldform-basis-coefficients}.
The normalized Hecke coefficient $\alpha_{\sigma_{v}}$ is defined
immediately after that formula.}

\notationentry{$\psi_{v}^{0},\ W_{v}^{0}$}
{The unramified character is $\psi_{v}^{0}(x)=\psi_{v}(\varpi_{v}^{-d_{v}}x)$
in \textsection\ref{def:unramified-additive-character}.
The vector $W_{v}^{0}$ in \textsection\ref{sec14.8.1} uses
\eqref{fc11.41} or \eqref{c12.2} with $m=r_{v}$, according as
$v\nmid\mathfrak{q}$ or $v\mid\mathfrak{q}$, starting from the
spherical $\psi_{v}^{0}$-Whittaker vector with value $1$ at $I_{3}$.}

\notationentry{$f_{\boldsymbol{\nu},v}^{\circ},\ f_{\boldsymbol{\nu},v}$}
{The spherical flat section \eqref{c14.1} and the test section
obtained from it by \eqref{fc14.3} at infinity or
\eqref{c14.5}--\eqref{c14.6} at finite places.}

\notationentry{$\mathcal{J}_{\psi_{v}}f$}
{The Jacquet integral \eqref{f7.11} sends a section $f$ to a
Whittaker vector. It is distinct from the local double Mellin
integral $\mathcal{J}_{v}$ in \eqref{f7.4}. At unramified finite
places the normalized spherical section gives the value
$\mathbf{L}(\pi_{v})^{-1}$ at $I_{3}$.}

\notationentry{$\mathcal{V}_{v}^{\circ},\ \mathcal{W}_{\infty}$}
{The normalized real weight of the spherical $\mathrm{GL}_{3}$ vector,
defined before Lemma~\ref{lem:mixed-spherical-arch}, is unchanged on
the dual side. The product times $\epsilon_{\sigma}D_{F}^{-3/2}$
defines $\mathcal{W}_{\infty}$ in \textsection\ref{def:mixed-real-weight}.}

\notationentry{$\mathcal{W}_{\infty,\mathbf{R}}$}
{The real weight for a box of scales $\mathbf{R}$. For the fifth
moment it is $\epsilon_{\sigma}D_{F}^{-3/2}\prod_{v\mid\infty}
H_{R_{v}}(t_{v})$. In the mixed moment the bounded-range factors
are $G_{\pi_{v}}(t_{v})S_{0}(t_{v})$ instead. See the definitions before
\eqref{eq:original-weighted-fourth} and \eqref{eq:mixed-polynomial-AFE-sieve}.}

\notationgroup{Spectral families and normalizations}

\notationentry{$\mathcal{F}(\mathfrak{q};\omega')$}
{Generic unitary $\mathrm{GL}_{2}$ spectrum of central character $\omega'$
and conductor dividing $\mathfrak{q}$, including the continuous spectrum.
See the definition before Theorem~\ref{thm1.3}.}

\notationentry{$\mathcal{F}^{\circ},\ \mathcal{F}_{\mathrm{cusp}}^{\circ}$}
{The superscript $\circ$ requires spherical Archimedean components.
The subscript $\mathrm{cusp}$ requires cuspidality. Neither imposes
exact finite conductor. See \textsection\ref{sec:explicit-fourth-moment-statement}.}

\notationentry{$\boldsymbol{\nu},\ \boldsymbol{\mu},\ \sigma_{\xi,\lambda}$}
{The inducing parameters on $G$ are $\boldsymbol{\nu}$ for the
minimal parabolic \eqref{2.9} and $\boldsymbol{\mu}$ for the maximal
parabolic \eqref{e2.3}.
The continuous spectrum on $G'$ is indexed by
$\sigma_{\xi,\lambda}=\xi|\cdot|^{\lambda}\boxplus
\overline{\xi}\omega'|\cdot|^{-\lambda}$.
See \textsection\ref{def:spectral-measure}.}

\notationentry{$\mathfrak{B}(\sigma),\ \mathfrak{B}(\xi_{1},\xi_{2})$}
{The global orthonormal bases of a cuspidal $\sigma$ and of the
space of inducing sections for the pair $(\xi_{1},\xi_{2})$,
respectively. Both are defined after \eqref{eq:inducing-global-norm},
which gives the inner product for the inducing sections.}

\notationentry{$d\mu_{\sigma},\ \epsilon_{\sigma}$}
{The spectral measure counts cuspidal representations and uses
$d\lambda/(4\pi i)$ for each inducing character. In the fourth
and mixed moments $\epsilon_{\sigma}$ is $1$ on the cuspidal spectrum
and $2$ on the continuous spectrum. See
\textsection\ref{def:spectral-measure} and \textsection\ref{def:spectral-multiplicity}.}

\notationentry{$L(1,\sigma,\Ad)$}
{In a global spectral integral this denotes $L^{\mathrm{reg}}(1,\sigma,\Ad)$.
On the continuous spectrum at $\lambda=it$ it equals
$(\Res_{z=1}\zeta_{F}(z))|L(1+2it,\xi^{2}\overline{\omega}')|^{2}$.
The local factors are not regularized. See \eqref{eq:regularized-adjoint}.}

\notationentry{$\mathcal{L},\ \mathcal{L}^{\dagger}$}
{Products and ratios of completed $L$-functions in the degenerate terms, with
$\dagger$ for the opposite Fourier expansion. Their local factors
$\mathcal{L}_{v}$ and $\mathcal{L}_{v}^{\dagger}$ are divided out
in the definitions of $\mathcal{P}_{v}^{\sharp}$ and
$\mathcal{P}_{v}^{\dagger,\sharp}$. See the definitions at the start
of \textsection\ref{sec5.5} and \textsection\ref{sec5.7}.}

\notationentry{$\mathbf{L}(\pi_{v}),\ \mathbf{L}_{\fin}(\pi),\ \mathcal{L}_{\fin}$}
{Boldface $\mathbf{L}(\pi_{v})$ is the local Jacquet normalizing factor,
and $\mathbf{L}_{\fin}(\pi)=\prod_{v<\infty}\mathbf{L}(\pi_{v})$.
See \textsection\ref{def:jacquet-euler-factor}.
Calligraphic $\mathcal{L}_{\fin}$ gives the finite fourth-moment factors in
\textsection\ref{def:finite-fourth-moment-factor} (spectral, with adjoint denominator),
\textsection\ref{sec14.8.3} (generic), and
\textsection\ref{sec14.8.2} (degenerate).}

\notationentry{$\lambda_{\sigma}(\mathfrak{n}),\
\widetilde{\lambda}_{\sigma,\mathbf{s}}(\mathfrak{q})$}
{The original twist is the Hecke coefficient
$\lambda_{\sigma}(\mathfrak{n})$ in \eqref{eq12.3}. The dual twist
\eqref{eq14.123} is a product over $v\mid\mathfrak{q}$ of the Hecke
coefficient at $\mathfrak{p}_{v}^{e_{v}(\mathfrak{q})}$ minus the
preceding coefficient times $q_{v}^{-(1+3s_{1}+s_{2})/2}$.}

\notationentry{$\lambda_{\pi_{v}}(\mathfrak{p}_{v}^{i+j},\mathfrak{p}_{v}^{j}),\
s_{i+j,j,0}$}
{The two-index $\mathrm{GL}_{3}$ Hecke coefficient is the Schur
polynomial \eqref{12.6} in the Satake parameters of unramified $\pi_{v}$.
Equation \eqref{12.5} relates it to $W_{v}^{\circ}$.}

\notationentry{$\mathcal{H}_{\sigma,\mathcal{T}},\
\mathcal{H}_{\sigma,\mathcal{T}}^{\vee}$}
{Global spectral weights in Theorem~\ref{thmC}, including the
discriminant and $\epsilon_{\sigma}$ factors. The finite ratios are at
$v\mid\mathfrak{q}$ on the original side and $v\mid\mathfrak{n}$
on the dual side. See \textsection\ref{def:fourth-spectral-weights}.}

\notationentry{$\mathcal{H}_{\mathrm{gen},\mathcal{T}},\
\mathcal{H}_{\mathrm{degen},\mathcal{T}}$}
{Local products in the generic and degenerate fourth-moment terms.
Their finite factors divide out the local Euler factors. The dual
products carry $\vee$. See \textsection\ref{def:fourth-generic-weights}
and \textsection\ref{sec14.8.2}, with $w\in\mathbb{W}
=\{w_{2},w_{1}w_{2},w_{1}w_{2}w_{1}\}$.}

\notationentry{$\widetilde{\mathcal{I}}_{\mathcal{T}},\
\widetilde{\mathcal{I}}_{\mathcal{T}}^{\vee}$}
{The normalized fourth moments defined before Theorem~\ref{thmC}.
The original has level $\mathfrak{q}$, twist $\mathfrak{n}$, and
spherical components at infinity. The dual has level $\mathfrak{n}$
and dual twist at $\mathfrak{q}$, without the spherical restriction.
The subscripts $\mathrm{gen}$, $\mathrm{degen}$ and $\mathrm{res}$ mark
the three correction terms defined after \eqref{fc1.3}.}

\notationentry{$\alpha_{v}(\sigma_{v};m),\ \alpha_{\pi}(\mathfrak{q})$}
{The weight of $W_{v}^{[m]}$ at $\mathbf{s}=\mathbf{0}$ after
division by $L_{v}(1/2,\pi_{v}\times\sigma_{v})L_{v}(1/2,\sigma_{v})$.
At $m=0$ it is $\mathbf{1}_{r_{\sigma_{v}}=0}$. See
\textsection\ref{def:mixed-finite-weight}. The Satake product
$\alpha_{\pi}(\mathfrak{q})$ is defined before Corollary~\ref{thm1.7}
and bounds the sieve coefficients.}

\notationentry{$\mathcal{M}_{\pi},\ \mathcal{M}_{\pi}^{\vee}$}
{The original and dual mixed first moments at $\mathbf{s}=\mathbf{0}$.
The original uses $\mathcal{W}_{\infty}$ and
$\alpha_{v}(\sigma_{v};e_{v}(\mathfrak{a}))$. See
\textsection\ref{def:mixed-finite-weight}.}

\notationentry{$C_{v},\ C_{v}(\pi),\ A_{v}$}
{In the type~\RNum{1} first-moment estimate, $C_{v}$ without an
argument is the chosen Archimedean cutoff for $\sigma_{v}$.
Here $C_{v}(\pi)$ is the analytic conductor of $\pi_{v}$ and
$A_{v}=C_{v}(\pi)/C_{v}$. See Theorem~\ref{thm1.3}.}

\notationentry{$\mathfrak{m},\ Q,\ \mathbf{C}_{\infty}$}
{In Theorem~\ref{thm1.3}, $\mathfrak{m}$ is the finite
conductor of $\pi$, $Q=N_{F}(\mathfrak{q}\cap\mathfrak{m})$,
and $\mathbf{C}_{\infty}=\prod_{v\mid\infty}C_{v}$ is the
product of the chosen cutoffs, not the conductor of a
representation.}

\notationentry{$r_{v},\ n_{v},\ r_{\sigma_{v}}$}
{The type~\RNum{2} test in \textsection\ref{sec11.2.1} uses
$r_{v}=e_{v}(\mathfrak{q})$ and $n_{v}=e_{v}(\mathfrak{n})$.
The type~\RNum{1} first moment in \textsection\ref{sec13} uses
$r_{v}=\max\{e_{v}(\mathfrak{q}),r_{\pi_{v}}\}$. The conductor
exponent $r_{\sigma_{v}}$ uses the convention in
\textsection\ref{sec:automorphic-data-conventions}.}

\notationentry{$R_{v},\ \mathbf{R},\ X$}
{Before \eqref{eq:original-weighted-fourth}, the real dyadic scales
give $\mathbf{R}=(R_{v})_{v\mid\infty}$ and
$X=\prod_{v\mid\infty}R_{v}^{2}$ for $C_{\infty}(\sigma)$, with
$R_{v}=2$ in the bounded and complementary ranges.
In the proof of Lemma~\ref{lem:spherical-projector-dual-mass},
$X$ is instead the variable in the finite Euler factor.}

\notationentry{$\vartheta,\ \vartheta_{\pi_{v}},\ \vartheta_{\infty}$}
{The finite-place Ramanujan exponent $\vartheta$ for $\mathrm{PGL}_{2}$
is fixed before Theorem~\ref{thm1.5}.
The largest absolute real part of the inducing exponents of $\pi_{v}$
is $\vartheta_{\pi_{v}}$ (\textsection\ref{sec:archimedean-analytic-conductors}).
The bound $\vartheta_{\infty}$ on the real complementary-series
exponents $\pm\eta_{v}$ is fixed at the start of \textsection\ref{sec15}.}

\notationgroup{Newform sieves and applications}

\notationentry{$S,\ f_{v}^{\mathrm{AL}}$}
{$S$ is the exceptional set for Euler products in
Definitions~\ref{def:exceptional-places-minimal} and~\ref{def:exceptional-places-maximal}.
In the sieve (Lemma~\ref{lem:polynomial-newform-sieve}), $S$ instead
selects level primes where the section is $q_{v}^{-s_{1}}f_{v}^{\mathrm{AL}}$,
with $f_{v}^{\mathrm{AL}}$ defined before
\eqref{eq:signed-section-Hecke-identity}. At $r_{\sigma_{v}}=1$ this
inserts $q_{v}^{-1/2}\lambda_{\sigma_{v}}(\mathfrak{p}_{v})$
by \eqref{eq:signed-newvector-factor}.}

\notationentry{$c_{\mathfrak{b},S}(\mathfrak{l}),\ b_{t}$}
{In the squarefree sieve, $c_{\mathfrak{b},S}(\mathfrak{l})$ from
\eqref{eq:squarefree-sieve-coefficients} weights level
$\mathfrak{a}/\mathfrak{b}$ with extra twist $\mathfrak{l}$.
The coefficients $b_{t}$ in \eqref{eq:depth-inverse} give the
projector \eqref{eq:depth-projector} for conductor exponents $\geq2$.}

\notationentry{$g_{\pi},\ g_{\pi}^{+}$}
{The twist coefficient in the mixed generic term and dual standard
residue, and its nonnegative majorant. They are defined before and
after \eqref{eq:mixed-generic-coefficient}, respectively.}

\notationentry{$\alpha_{\infty}(\sigma_{\infty})$}
{The real-place weight in the nonvanishing count \eqref{1.5}.
It is built from the $\mathcal{V}_{v}^{\circ}$ and decays rapidly
in $C_{\infty}(\sigma)$. See \textsection\ref{sec:proof-simultaneous-nonvanishing}.}

\notationgroup{Continuation, residues, and central values}

\notationentry{$\mathcal{J},\ \mathcal{J}_{v},\ \widetilde{\mathcal{J}}$}
{The type~\RNum{1} double Mellin integral $\mathcal{J}$ and its
local factor $\mathcal{J}_{v}$ are defined in \textsection\ref{sec7.2}
at \eqref{f7.4}. Its continuation $\widetilde{\mathcal{J}}$ is
defined before \eqref{e4.10}.}

\notationentry{$F_{\xi},\ \widetilde{\mathcal{G}}$}
{The type~\RNum{1} Eisenstein sum $F_{\xi}$ pairs a period with a
Whittaker value at $I_{2}$ (before Proposition~\ref{prop4.1}).
The type~\RNum{2} sum $\widetilde{\mathcal{G}}$ pairs two continued
periods (start of \textsection\ref{sec4.2}).}

\notationentry{$R_{\RNum{1}}^{\pm},\ R_{\mathrm{dual}}^{\pm}$}
{Type~\RNum{1} residues of the $\mathrm{GL}_{2}$ Eisenstein integral
and the dual $\mathrm{GL}_{1}$ integral, respectively.
Their definitions are in Propositions~\ref{prop4.1} and~\ref{prop4.5}.}

\notationentry{$\mathscr{P}_{\chi},\ \mathscr{P}_{\chi,v}$}
{The normalized Eisenstein basis sum for the type~\RNum{1} spectral
residues and its local factor, defined in and just after
\eqref{eq:maximal-spectral-density-factors}. The local factor is $1$
at almost all places.}

\notationentry{$R_{\RNum{2}}^{1,\pm},\ R_{\RNum{2}}^{2,\pm}$}
{Residues from the Rankin--Selberg period at parameter $s_{1}$
and the twisted standard period at parameter $s_{2}$,
respectively. They are defined before Proposition~\ref{prop4.4}.
The dual moment uses $\mathbf{s}^{\vee}$.}

\notationentry{$\mathfrak{X}_{\pi,\omega'}^{\pm},\
\mathfrak{X}_{\pi,\omega'}^{\mathrm{dual}}$}
{Characters indexing possible Rankin--Selberg residues in the
Eisenstein spectrum and in the type~\RNum{1} dual integral,
respectively. These sets are empty for cuspidal $\pi$.
See \textsection\ref{sec4.1} and \eqref{e4.11}.}

\notationentry{$D_{\mathcal{T}},\ D_{\mathcal{T}}^{\vee}$}
{Original and dual fourth-moment densities at
$\sigma_{\lambda}=|\cdot|^{\lambda}\boxplus|\cdot|^{-\lambda}$,
including their respective finite Hecke factors and weights.
See \textsection\ref{def:fourth-residue-densities}.}

\notationentry{$\mathcal{H}_{W_{v}}^{\vee,\reg}$}
{The finite dual weight defined before \eqref{eq:finite-residual-estimates}.
At $v\mid\mathfrak{n}$, \eqref{eq:residue-weight-interpolation} removes
oldform singularities and preserves residues. At $v\mid\mathfrak{q}$
the weight is unchanged.}

\notationentry{$\mathscr{R}_{\mathbf{s},\boldsymbol{\nu}}[D]$}
{The signed half-residues of a meromorphic density $D(\lambda)$
at the Rankin--Selberg poles indexed by $\nu_{i}$ and the two standard poles. A single contour
encloses coincident poles of the same sign. See
\eqref{eq:residue-operator} and \eqref{eq:residue-cluster}.}

\notationentry{$\mathscr{L}_{b}(h),\ \Phi_{V}(a,k)$}
{The norm $\mathscr{L}_{b}(h)$, defined before
Lemma~\ref{lem:mixed-fixed-real-dual-decay}, integrates $|h(g)|$ with
weight $(1+\|g\|+\|g^{-1}\|)^{b}$. For a smooth Whittaker vector $V$
of $\pi_{v}$, the proof defines $\Phi_{V}(a,k)$ by integrating
$V(\diag(az,z,1)k)$ against $d^{\times}z$.}

\notationentry{$c_{\pi_{v}}(V)$}
{For self-dual $\pi_{v}$ as in Lemma~\ref{lem:mixed-real-standard-pole},
the coefficient of $a(\log a)^{m-1}$ in $\int_{K_{v}'}\Phi_{V}(a,k)\,dk$
as $a\to0^{+}$. Here $m$ counts the zero Langlands parameters of $\pi_{v}$.
The definition is just before \eqref{eq:mixed-real-standard-pole-value}.}

\notationentry{$\operatorname{CT}_{0}$}
{Constant term in $(s_{1},s_{2})$ on circles with $|s_{2}|\ll|s_{1}|$.
For a dual term, substitute $\mathbf{s}^{\vee}$ on the same contours.
See \textsection\ref{sec8.2}.}

\notationentry{$J_{*}^{\heartsuit},\ I_{*}^{\heartsuit}$}
{The continued type~\RNum{1} terms are $J_{*}^{\heartsuit}$
(\textsection\ref{sec8.1}), and the continued type~\RNum{2} terms are
$I_{*}^{\heartsuit}$ (\textsection\ref{sec8.2}). When the parameters
are omitted, these sections specify which terms are evaluated at the
center and which use constant terms.}

\endgroup
\clearpage

\bibliographystyle{alpha}
	
\bibliography{LY}

\end{document}